\documentclass[12pt]{article}

\usepackage[colorlinks=true,linkcolor=blue,
  citecolor=blue,
  urlcolor=blue,
  filecolor=blue]{hyperref}

\usepackage[all]{xy}

\usepackage[latin1]{inputenc}
\usepackage{amsfonts}
\usepackage{amsmath}
\usepackage{amssymb}
\usepackage{amscd}
\usepackage{amsthm}
\usepackage{indentfirst}
\usepackage[hmargin=2.98cm
,vmargin=4.0cm
]{geometry}

\usepackage{mathrsfs} 

\usepackage{faktor}

\usepackage{tikz-cd}

\usepackage{epsfig}
\usepackage{enumitem}

\usepackage{booktabs} 

\usepackage[T1]{fontenc}
\usepackage{lmodern} 
\usepackage{fix-cm}   

\usepackage{graphicx}
\usepackage{subcaption}
\usepackage{float}
\usepackage{xcolor}

\newtheorem{thm}{Theorem}[section]
\newtheorem{lemma}[thm]{Lemma}

\newtheorem{prop}[thm]{Proposition}
\newtheorem{coroll}[thm]{Corollary}
\newtheorem{claim}[thm]{Claim}

\theoremstyle{definition}

\newtheorem{defin}[thm]{Definition}
\newtheorem{rem}[thm]{Remark}
\newtheorem{exam}[thm]{Example}

\newcommand{\R}{{\mathbb{R}}}

\newcommand{\Q}{{\mathbb{Q}}}
\newcommand{\T}{{\mathbb{T}}}
\newcommand{\Z}{{\mathbb{Z}}}
\newcommand{\N}{{\mathbb{N}}}
\newcommand{\C}{{\mathbb{C}}}

\newcommand{\CP}{{\mathbb{C}P}}
\newcommand{\RP}{{\mathbb{R}P}}

\newcommand{\cA}{{\mathcal{A}}}
\newcommand{\cB}{{\mathcal{B}}}
\newcommand{\cC}{{\mathcal{C}}}

\newcommand{\cE}{{\mathcal{E}}}
\newcommand{\cF}{{\mathcal{F}}}

\newcommand{\cH}{{\mathcal{H}}}

\newcommand{\cJ}{{\mathcal{J}}}

\newcommand{\cL}{{\mathcal{L}}}
\newcommand{\cM}{{\mathcal{M}}}
\newcommand{\cN}{{\mathcal{N}}}
\newcommand{\cO}{{\mathcal{O}}}
\newcommand{\cP}{{\mathcal{P}}}

\newcommand{\cR}{{\mathcal{R}}}

\newcommand{\cU}{{\mathcal{U}}}

\makeatletter
\newcommand{\hbigl}[1]{\mathopen{\mathpalette\halfbig@{#1}}}
\newcommand{\hbigr}[1]{\mathclose{\mathpalette\halfbig@{#1}}}

\newcommand{\halfbig@}[2]{%
  \vcenter{\hbox{\scalebox{1}[1.5]{$\m@th#1#2$}}}%
}
\makeatother

\newcommand{\fc}{\colon}
\newcommand{\ve}{\varepsilon}

\newcommand{\wh}{\widehat}

\newcommand{\bS}{\mathbf{S}}

\DeclareMathOperator{\RS}{RS}

\DeclareMathOperator{\SL}{SL}
\DeclareMathOperator{\mn}{\mu\nu}
\DeclareMathOperator{\U}{U}

\DeclareMathOperator{\wind}{wind}
\DeclareMathOperator{\Fix}{Fix}
\DeclareMathOperator{\Mat}{Mat}
\DeclareMathOperator{\Id}{Id}
\DeclareMathOperator{\Sp}{Sp}
\DeclareMathOperator{\Mod}{Mod}
\DeclareMathOperator{\Ch}{Ch}

\DeclareMathOperator{\Hom}{Hom}
\DeclareMathOperator{\spn}{span}
\DeclareMathOperator{\sgn}{sgn}

\DeclareMathOperator{\Crit}{Crit}

\DeclareMathOperator{\im}{im}
\DeclareMathOperator{\id}{id}

\DeclareMathOperator{\ter}{ter}

\DeclareMathOperator{\ev}{ev}

\DeclareMathOperator{\ind}{ind}

\DeclareMathOperator{\sign}{sign}
\DeclareMathOperator{\ini}{ini}
\DeclareMathOperator{\tel}{tel}
\DeclareMathOperator{\cone}{cone}

\newcommand{\drctlim}{\underset{\longrightarrow}{\lim}\,}

\DeclareMathOperator{\supp}{supp}

\DeclareMathOperator{\Ham}{Ham}

\DeclareMathOperator{\Cont}{Cont}

\DeclareMathOperator{\val}{val}
\DeclareMathOperator{\Int}{Int}

\DeclareMathOperator{\PD}{PD}

\DeclareMathOperator{\coker}{coker}

\DeclareMathOperator{\res}{res}

\DeclareMathOperator{\0}{\operatorname{\mathbf{0}}}

\begin{document}
	
	\title{Relative symplectic cohomology in complex projective spaces}
	\author{Adi Dickstein\footnote{Partially supported by the Fondation Courtois.}
    \ \ 
     and Yaniv Ganor
     }

	\date{}
	
	\setcounter{tocdepth}{2}
	
	\renewcommand{\labelenumi}{(\roman{enumi})}

	\renewcommand{\labelenumi}{(\roman{enumi})}

	\maketitle
    \begin{abstract}
        Relative symplectic cohomology is an invariant of compact subsets of a closed symplectic manifold, introduced by Varolg\"une\c s \cite{Varolgunes_2021_MV_and_relSH}. There are many examples of computations of this invariant over the Novikov field, but the collection of computed examples over the Novikov ring is still quite limited. One reason for this is that such computations require determining the relevant Floer complexes for Hamiltonians that are not necessarily $C^2$-small Morse functions. In this work, we present a computation of relative symplectic cohomology over the Novikov ring for balls and their complements in $\CP^n$. Our computation relies on explicit descriptions of Floer complexes, in the Morse--Bott setting with cascades, for J-shaped Hamiltonians on $\CP^n$. This allows us to deduce new estimates for the stable displacement energy of the boundaries of balls in $\CP^n$.
    \end{abstract}
    \tableofcontents

   \section{Introduction}
\subsection{Background}
Symplectic topology began as a geometric language for Hamiltonian mechanics, providing a coordinate free framework for the study of the group of Hamiltonian diffeomorphisms. Liouville's theorem states that these diffeomorphisms are volume preserving, and hence a natural question rising from this fact is whether there exist finer symplectic invariants. Also, Darboux's theorem states that every symplectic manifold is locally isomorphic to an open set in $\R^{2n}$ equipped with the standard symplectic form. As a consequence of this, symplectic manifolds admit no local invariants, hence the need for global properties, which are the subject of symplectic topology.

In his seminal paper \cite{Gromov_1985_J_curves} from 1985, Gromov provided an answer to this question using pseudo-holomorphic curves. He proved his famous non-squeezing theorem which later led to the theory of symplectic capacities.

  The next landmark was the development of symplectic homology, pioneered by Floer and Hofer, in the late 1980s and 1990s, and further developed by many other experts throughout the years, see \cite{FH_1994_SH_I_open_sets_Cn,FHW_1994_SH_I,CFH_1995_SH_II_general,CFHW_1996_ApSH_II} for instance. This homology, which is based on Floer theory, serves as a refinement of the notion of capacity and associates a filtered module to a given Liouville domain. This homology theory is endowed with algebraic structures and has functorial properties.

  This theory, however, does not obey a general local-to-global principle. In general, the symplectic homology of a Liouville domain is not reconstructible from that of the elements of an open cover.

In 2018, Varolg\"une\c s introduced a variant of symplectic homology called \textbf{relative symplectic cohomology} \cite{Varolgunes_2021_MV_and_relSH}. Relative symplectic cohomology associates a module over the Novikov ring $\Lambda_{\geq0}$ to any compact subset of a closed symplectic manifold, together with restriction maps for inclusions of compacts. The novelty of this version is that it admits a local-to-global principle: it satisfies a Mayer--Vietoris property for commuting subsets, i.e., zero sets of Poisson-commuting functions.

The definition of relative symplectic cohomology is complicated, combining Floer theory with non-Archimedean algebra; hence, its computation, even in elementary cases, such as the empty set and the manifold itself, is not an obvious task. In this paper, we consider basic interesting examples and develop techniques that provide a full description.

The relative symplectic cohomology can be taken with coefficients in the Novikov field $\Lambda$. In this setting, the local-to-global principle becomes extra useful, since in many examples the relative symplectic cohomology over the field vanishes, and short exact sequences over a field split. As a result of that, several computations over the Novikov field have been carried out.  See \cite{Tonkonog_Varolgunes_Super_rigidity_of_certain_skeleta_using_relative_symplectic_cohomology,Sun_2024_ind-bdd_relSH,DGPZ_2024_Symp_top_and_IVMs,Groman_Varolgunes_2023_locality_of_SH,MSV_2024_heavy_sets_and_SH,Sun_2024_heavy + relSH,MSV_2024_heavy_sets_and_SH} for several prominent results in the Novikov field setting.

Over the Novikov ring $\Lambda_{\geq0}$, however, such a computation requires, first of all, the explicit computation of the Floer complexes associated to Hamiltonians that are not necessarily $C^2$-small, and, in addition, the computation of all continuation maps between these Floer complexes. To the best of our knowledge, the only place where an explicit expression of the cochain complex in such a situation has been displayed is Varolg\"une\c s's thesis \cite{Varolgunes_2018_PhD}, for the case of disks in the sphere $S^2$.

In contrast, the goal of the present work is to obtain a full computation over the Novikov ring.

\subsection{Main results}
The main results of this paper are a full computation of the relative symplectic cohomology with coefficients in the Novikov ring for toric balls and their complements in $\CP^n$. Let us introduce a few notations before formulating the main result.

\begin{itemize}

\item \textbf{Novikov intervals:} 
The Novikov field
$$\Lambda=\left\{\sum_{i=1}^\infty a_i T^{\lambda_i}\,\middle|\,a_i\in \Q,\,\lambda_i\in \R\,\,\text{s.t.}\,\,\lambda_i\to+\infty\right\}$$
admits a natural valuation $\val\fc \Lambda\to\R\cup\{+\infty\}$ defined as follows. Let $x\in\Lambda$. If $x=0$ then $\val(x)=+\infty$; otherwise, $\val(x)=\min\{\lambda_i\,|\,i\in\N,\,a_i\neq0\}$, where $x=\sum_{i=1}^\infty a_i T^{\lambda_i}$, for $a_i\in \Q$, $\lambda_i\in \R$ for every $i\in \N$ with $\lambda_i\to+\infty$.

Using the valuation, we can define the \textbf{interval modules}
$$\Lambda_{\geq r} = \val^{-1}([r,\infty])\,,\qquad \Lambda_{>r}=\val^{-1}((r,\infty])\,,\qquad \Lambda_{(a,b]}=\Lambda_{> a}/\Lambda_{> b}\text{ for some }a< b\,.$$
Note that all of these modules are modules over the \textbf{Novikov ring} $\Lambda_{\geq0}$.

\item \textbf{Relative symplectic cohomology:} 
    Given a closed symplectic manifold $(M,\omega)$ and a compact subset $K \subset M$, the relative symplectic cohomology $SH_M^*(K;\Lambda_{\geq0})$ is a graded module over the Novikov ring $\Lambda_{\geq 0}$, where the grading is taken modulo twice the minimal first Chern number of $M$; see also \cite{DGPZ_2024_Symp_top_and_IVMs} for further explanations about the grading. Moreover, for any two compact subsets $K,K' \subset M$ with $K \subset K'$, there exists a map, called the restriction map, which is a module homomorphism:
    \[
        \mathrm{res}_K^{K'} : SH_M^*(K';\Lambda_{\geq0}) \to SH_M^*(K;\Lambda_{\geq0}).
    \]
    Additionally, Varolg\"une\c s showed that 
    $$
     SH_M^*(M;\Lambda_{\geq0}) = H^*(M;\Z)\otimes_{\Z}\Lambda_{>0}\quad \text{and} \quad SH_M^*(\varnothing;\Lambda_{\geq0}) = 0,
    $$
    where $\Lambda_{>0}$ is \textbf{the maximal ideal} of $\Lambda_{\geq 0}$. 
    He also proved that for any compact subset $K\subset M$ displaceable by a Hamiltonian diffeomorphism, the relative symplectic cohomology $SH_M^*(K;\Lambda_{\geq0})$ is a torsion module; that is,
    $$
        SH_M^*(K;\Lambda_{\geq0}) \otimes_{\Lambda_{\geq 0}} \Lambda = 0,
    $$
    where $\Lambda$ is the Novikov field.

        \item \textbf{Toric balls in $\CP^n$:} For every $\Delta>0$, denote by $K^n(\Delta)$ the $n$-simplex
        $$
        K^n(\Delta)=\left\{(x_1,\ldots,x_n)\in \R^n \,\middle|\, x_1,\ldots,x_n\geq 0 \text{ and } x_1+\cdots+x_n\leq \Delta\right\}.
        $$
        We consider the complex projective space $\CP^n$ equipped with the Fubini--Study symplectic form, normalized to have Gromov width $1$. Throughout the paper, whenever we mention the ``capacity'' of a set, we mean the Gromov width of its interior.
        
        Note that the minimal first Chern number of $\CP^n$ is $n+1$. The standard moment map on $\CP^n$ is the function $\Phi\fc \CP^n\to \R^n$ given by
        $$
        \Phi(z)=\left(\frac{|z_1|^2}{\sum_{j=0}^n |z_j|^2}, \ldots, \frac{|z_n|^2}{\sum_{j=0}^n |z_j|^2}\right),
        $$
        for every $z\in \CP^n$ written in homogeneous coordinates $z=[z_0:\cdots:z_n]$. The image of $\Phi$ is the $n$-simplex $K^n(1)$. The \textbf{standard closed toric ball in $\CP^n$ of capacity $\Delta$}, for some $\Delta\in(0,1)$, is defined as the preimage of the $n$-simplex $K^n(\Delta)$ under $\Phi$, i.e. $\Phi^{-1}(K^n(\Delta))$.

\end{itemize}

The following two theorems constitute the main results of the paper:

    \begin{thm}\label{thm: SH of ball in CP^n}
        Let $B^{}_\Delta$ be the standard closed toric ball in $\CP^n$ of capacity $\Delta$, for some $\Delta\in(0,1)$.  The relative symplectic cohomology $SH^*_{\CP^n}( B^{}_{\Delta}; \Lambda_{\geq 0})$, as a $\Z/2(n+1)\Z$-graded module over $\Lambda_{\geq 0}$, is given explicitly as follows:
        \begin{itemize}
            
            \item If $\Delta<\frac{n}{n+1}$ then
            $$SH^*_{\CP^n}(B^{}_\Delta;\Lambda_{\geq0})=\left\{\begin{array}{cl}
               \bigoplus\limits_{i=0}^\infty \Lambda_{(0,\Delta]}, & * \equiv 0,2,\ldots,2n \mod 2n+2,\\
                0, & *\equiv 1,3,\ldots, 2n+1 \mod 2n+2. 
            \end{array}\right.$$
            \item If $\Delta\geq\frac{n}{n+1}$ then
            $$SH^*_{\CP^n}(B^{}_\Delta;\Lambda_{\geq0})=\left\{\begin{array}{cl}
               \Lambda_{>0}\oplus\bigoplus\limits_{i=0}^\infty \Lambda_{(0,n(1-\Delta)]}, &  * \equiv 0,2,\ldots,2n \mod 2n+2,\\
                0, & *\equiv 1,3,\ldots, 2n+1 \mod 2n+2. 
            \end{array}\right.$$
        \end{itemize}
    \end{thm}

\begin{rem}
 In his thesis \cite{Varolgunes_2018_PhD}, Varolg\"une\c s described a complex computing the relative symplectic cohomology of disks in $S^2$.
\end{rem}
\begin{rem}\phantom{M}
\begin{itemize}
    \item The appearance of the threshold $\frac{n}{n+1}$ is not coincidental. 
Indeed, Entov--Polterovich \cite{EP_2009_rigid_subsets} proved that a standard closed toric ball in $\CP^n$ of capacity $\Delta$ is heavy if and only if 
$\Delta \geq \frac{n}{n+1}$. In this case, the ball contains the Clifford torus of $\CP^n$, which is the fiber of the moment map $\Phi$ over the point $\left(\frac{1}{n+1},\ldots,\frac{1}{n+1}\right) \in \R^n$. Moreover, Mak--Sun--Varolg\"une\c s \cite{MSV_2024_heavy_sets_and_SH} showed that a compact subset is heavy precisely when its relative symplectic cohomology does not vanish over the Novikov field. Equivalently, heaviness is characterized by the non-triviality of the torsion-free part of the relative symplectic cohomology, exactly as we observe in Theorem~\ref{thm: SH of ball in CP^n}.

    \item When $\Delta < \frac{n}{n+1}$, we found that $SH_{\CP^n}(B^{}_\Delta;\Lambda_{\geq0})$ consists of a countable direct sum of Novikov intervals of length $\Delta$, a value which coincides with the Gromov width of the interior of $B$. In a work in progress by the authors, in collaboration with Frol Zapolsky, we prove that for sufficiently small balls, the supremal length of a Novikov interval in the relative symplectic cohomology of a ball, in a general closed symplectic manifold, is at least the Gromov width of its interior.

    \item When $\frac{1}{2}\leq \Delta < \frac{n}{n+1}$, Entov--Polterovich \cite{EP_2009_rigid_subsets} proved that the ball is not heavy and, in fact, stably-displaceable. Moreover they noted that it is still non-displaceable due to a variation on the argument for Gromov's two-ball packing inequality \cite[0.3.B]{Gromov_1985_J_curves}. As we see in the conclusion of Theorem~\ref{thm: SH of ball in CP^n}, the relative symplectic cohomology of the ball does not reflect the non-displaceability property.

\end{itemize}
  
\end{rem}

\begin{rem}
We contrast Theorem~\ref{thm: SH of ball in CP^n}, with the computation of relative symplectic cohomology of a ball in $\C^n$ with cappings (in the sense of \cite{BSV_2022_QH_deformamation of SH}), which follows from \cite[Section 3.3]{Oancea_survey} (see also Section~\ref{ss: CF of J-shaped} here): $$SH_{\C^n}^*(B_\Delta)=\bigoplus_{j\geq1}\cF(-j\Delta,-(j-1)\Delta],$$ where the Conley--Zehnder index of the $j$-th summand is $-2jn$. Here we use the presentation of $SH_{\C^n}^*(B_\Delta)$ as a persistence module, as a countable sum of interval modules with coefficients in a field $\cF$, see \cite{PRSZ_2020_persistence}, for an overview. Note, that all the intervals are of length $\Delta$.
\end{rem}
We turn to state our result of the computation of the relative symplectic cohomology of the complement of the ball. This result involves more cases than for the ball:

\begin{thm}\label{thm: SH of complement of ball in CPn}
Let $B^{}_{\Delta}$ be the standard closed ball of capacity $\Delta \in (0,1)$ in $\CP^n$. The relative symplectic cohomology $SH^*_{\CP^n}(\CP^n \setminus \Int B^{}_{\Delta}; \Lambda_{\geq 0})$, as a $\Z/2(n+1)\Z$-graded module over $\Lambda_{\geq 0}$, is given explicitly as follows:

\begin{enumerate}
    \item If $\Delta \leq \frac{1}{2}$:
    $$SH^*_{\CP^n}(\CP^n \setminus \Int B^{}_{\Delta}; \Lambda_{\geq 0}) = \begin{cases}
        0, & * \text{ is odd}, \\
        \Lambda_{>0} \oplus \bigoplus\limits_{k=0}^\infty\Lambda_{(0, \Delta]}, &   * \text{ is even}.
    \end{cases}$$

    \item If $\frac{j}{j+1} < \Delta \leq \frac{j+1}{j+2}$ for some integer $1 \leq j \leq n-1$:
    $$SH^*_{\CP^n}(\CP^n \setminus \Int B^{}_{\Delta}; \Lambda_{\geq 0}) = \begin{cases}
        0, & * \text{ is odd}, \\
        \Lambda_{(0, (n-i)(1-\Delta)]} \oplus \Lambda_{>0} \oplus \bigoplus\limits_{k=0}^\infty\Lambda_{(0, \Delta]}, & *=2i,\,\,i\in\{n-j,\ldots,n-1\}, \\
        \Lambda_{>0} \oplus \bigoplus\limits_{k=0}^\infty\Lambda_{(0, \Delta]}, & *=2i,\,\,i\in\{n,0,\ldots,n-j-1\}.
    \end{cases}$$

    \item If $\Delta > \frac{n}{n+1}$:
    $$SH^*_{\CP^n}(\CP^n \setminus \Int B^{}_{\Delta}; \Lambda_{\geq 0}) = \begin{cases}
        0, &  * \text{ is odd}, \\
        \Lambda_{(0, (n-i)(1-\Delta)]} \oplus \bigoplus\limits_{k=0}^\infty\Lambda_{(0, n(1-\Delta)]}, &  * = 2i,\,\, i\in\{1,\ldots,n-1\}, \\
        \bigoplus\limits_{i=0}^\infty\Lambda_{(0, n(1-\Delta)]}, &  * \in\{ 0,2n\} .
    \end{cases}$$
\end{enumerate}
\end{thm}

\begin{rem}
In contrast to the relative symplectic cohomology of the ball, Theorem~\ref{thm: SH of complement of ball in CPn} captures a distinction between the cases where the ball $\Int B^{}_{\Delta}$ is displaceable ($\Delta \leq \frac{1}{2}$) and where it is non-displaceable but stably displaceable ($\frac{1}{2} < \Delta \leq \frac{n}{n+1}$). In the case of the ball, the lengths of all finite interval modules, as functions of $\Delta$, achieve a unique maximum at $\Delta = \frac{n}{n+1}$. By contrast, in the relative symplectic cohomology of the complement of the ball, for every $1 \leq i \leq n-1$, there is a direct summand interval module of degree $2i$ whose length, as a function of $\Delta$, is maximized at $\frac{n-i}{n-i+1}$. In particular, on the interval $[\frac{1}{2}, \frac{n}{n+1}]$, it is no longer true that the lengths of all finite interval modules increase monotonically, as they do in the case of the ball.

It would be interesting to find a direct connection between this phenomenon and the displaceability properties of the ball.
\end{rem}

\subsection{Restriction maps}

The relative symplectic cohomology is a functor, it assigns a module to a compact set and a restriction map to an inclusion of compacts. Therefore, it natural to ask for a computation of the restriction maps. Here we do this for two balls.

\begin{thm}\label{thm: res for balls in CP^n}
Let $B^{}_{\Delta'} \subset B^{}_{\Delta}$ be standard closed toric balls in $\CP^n$ with capacities $\Delta'$ and $\Delta$, respectively, where $0 < \Delta' \leq \Delta < 1$. Let $0 \leq j \leq n$. Then the restriction map in degree $2j$ modulo $2(n+1)$,
$$
\res^{B^{}_{\Delta}}_{B^{}_{\Delta'}} \fc SH_{\CP^n}^{2j}(B^{}_{\Delta};\Lambda_{\geq0}) \to SH_{\CP^n}^{2j}(B^{}_{\Delta'};\Lambda_{\geq0}),
$$
is completely determined as follows: 
\begin{enumerate}
    \item If $\Delta' \leq \Delta < \frac{n}{n+1}$, then the map 
    $$\res^{B^{}_{\Delta}}_{B^{}_{\Delta'}} \fc \bigoplus_{i \geq 0} \Lambda_{(0, \Delta]} e_i \to \bigoplus_{i \geq 0} \Lambda_{(0, \Delta']} e'_i$$
    is given by
    $$\res^{B^{}_{\Delta}}_{B^{}_{\Delta'}}(\lambda e_i) = T^{\hbigl(j+i(n+1)\hbigr)(\Delta-\Delta')} \lambda e'_i$$
    for every $i \geq 0$ and $\lambda \in \Lambda_{(0, \Delta]}$, where $e_i$ and $e'_i$ are generators of $\Lambda_{[0,\Delta]}$ and $\Lambda_{[0,\Delta']}$, respectively.

    \item If $\frac{n}{n+1} \leq \Delta' \leq \Delta$, then the map 
    $$\res^{B^{}_{\Delta}}_{B^{}_{\Delta'}} \fc \Lambda_{>0} e \oplus \bigoplus_{i \geq 0} \Lambda_{(0, n(1-\Delta)]} e_i \to \Lambda_{>0} e' \oplus \bigoplus_{i \geq 0} \Lambda_{(0, n(1-\Delta')]} e'_i$$
    is given by
    $$\res^{B^{}_{\Delta}}_{B^{}_{\Delta'}}(\lambda e) = T^{j(\Delta-\Delta')} \lambda e',$$
    for every $\lambda \in \Lambda_{>0}$, where $e$ and $e'$ are generators of copies of $\Lambda_{\geq0}$, and
    $$\res^{B^{}_{\Delta}}_{B^{}_{\Delta'}}(\lambda e_i) = T^{(j+(i+1)(n+1))(\Delta-\Delta')} \lambda e'_i$$
    for every $i \geq 0$ and $\lambda \in \Lambda_{(0, n(1-\Delta)]}$, where $e_i$ and $e'_i$ are generators of $\Lambda_{[0,n(1-\Delta)]}$ and $\Lambda_{[0,n(1-\Delta')]}$, respectively.

    \item If $\Delta' < \frac{n}{n+1} \leq \Delta$, then the map 
    $$\res^{B^{}_{\Delta}}_{B^{}_{\Delta'}} \fc \Lambda_{>0} e \oplus \bigoplus_{i \geq 0} \Lambda_{(0, n(1-\Delta)]} e_i \to \bigoplus_{i \geq 0} \Lambda_{(0, \Delta']} e'_i$$
    is given by
    $$\res^{B^{}_{\Delta}}_{B^{}_{\Delta'}}(\lambda e) = T^{j(\Delta-\Delta')} \sum_{i=0}^\infty (-T^{n(1-\Delta')-\Delta'})^{i} \lambda e'_i,$$
    for every $\lambda \in \Lambda_{>0}$, and
    $$\res^{B^{}_{\Delta}}_{B^{}_{\Delta'}}(\lambda e_i) = T^{\hbigl(j+i(n+1)\hbigr)(\Delta-\Delta')+\Delta-n(1-\Delta)} \lambda e'_i$$
    for every $i \geq 0$ and $\lambda \in \Lambda_{(0, n(1-\Delta)]}$, where $e_i$ and $e'_i$ are generators of $\Lambda_{[0,n(1-\Delta)]}$ and $\Lambda_{[0,\Delta']}$, respectively. Also, $e$ is a generator of $\Lambda_{\geq0}$.
\end{enumerate}  
\end{thm}
For the sake of completeness, we also computed the restriction map from the relative symplectic cohomology of $\CP^n$ itself to that of a closed toric ball.

\begin{prop}\label{prop: res from CPn to a ball}
Let $B^{}_{\Delta}$ be a standard closed toric ball in $\CP^n$ of capacity $\Delta$, for some $\Delta\in(0,1)$. 
 Let $0 \leq j \leq n$. Then the restriction map in degree $2j$ modulo $2(n+1)$,
$$
\res^{\CP^n}_{B^{}_{\Delta}} \fc SH_{\CP^n}^{2j}(\CP^n;\Lambda_{\geq0}) \to SH_{\CP^n}^{2j}(B^{}_{\Delta};\Lambda_{\geq0}),
$$

is completely determined as follows: 
\begin{enumerate}

    \item If $\Delta < \frac{n}{n+1} $, then the map 
    $$\res^{\CP^n}_{B^{}_{\Delta}} \fc \Lambda_{>0} f \to \bigoplus_{i \geq 0} \Lambda_{(0, \Delta]} e_i$$
    is given by 
    $$\res^{\CP^n}_{B^{}_{\Delta}}(\lambda f) = T^{j(1-\Delta)} \sum_{i=0}^\infty (-T^{n(1-\Delta)-\Delta})^{i} \lambda e_i,$$
    for every $\lambda \in \Lambda_{>0}$.

     \item If $\Delta\geq\frac{n}{n+1}$, then the map 
    $$\res^{\CP^n}_{B^{}_{\Delta}} \fc \Lambda_{>0} f \to \Lambda_{>0} e \oplus \bigoplus_{i \geq 0} \Lambda_{(0, n(1-\Delta)]} e_i$$
    is given by 
    $$\res^{\CP^n}_{B^{}_{\Delta}}(\lambda f) = T^{j(1-\Delta)} \lambda e,$$
    for every $\lambda \in \Lambda_{>0}$.

\end{enumerate}

\end{prop}

\begin{prop}\label{prop: res from CPn to the complement a ball}
    Let $B^{}_{\Delta}$ be the standard closed ball of capacity $\Delta \in (0,1)$ in $\CP^n$. Let $0 \leq i \leq n$. Then the restriction map in degree $2i$ modulo $2(n+1)$,
    $$
    \res^{\CP^n}_{\CP^n \setminus \Int B^{}_{\Delta}} \fc SH_{\CP^n}^{2i}(\CP^n;\Lambda_{\geq0}) \to SH_{\CP^n}^{2i}(\CP^n \setminus \Int B^{}_{\Delta};\Lambda_{\geq0}),
    $$
    depends on the value of $\Delta$ and is completely determined as follows, for every $\lambda \in \Lambda_{>0}$, where $f$ is the generator of $SH_{\CP^n}^{2i}(\CP^n;\Lambda_{\geq0}) \cong \Lambda_{>0} f \subseteq \Lambda_{\geq 0} f$:
    
    \begin{enumerate}
        \item Assume $\Delta \leq \frac{1}{2}$. The restriction map 
        $$\res^{\CP^n}_{\CP^n \setminus \Int B^{}_{\Delta}} \fc \Lambda_{>0} f \to \Lambda_{>0} e_1 \oplus \bigoplus_{k=2}^\infty \Lambda_{(0, \Delta]} e_k$$
        is given by
        $$\res^{\CP^n}_{\CP^n \setminus \Int B^{}_{\Delta}}(\lambda f) = \begin{cases}
             
            \lambda e_1, & 0 \leq i \leq n-1,\\
            T^{\Delta} \lambda e_1, & i = n,
        \end{cases}$$
        for every $\lambda \in \Lambda_{>0}$.
        \item Assume $\frac{j}{j+1} < \Delta \leq \frac{j+1}{j+2}$ for some integer $1 \leq j \leq n-1$. 
        \begin{itemize}
            \item If either $i\leq n-1-j$ or $i=n$, then the restriction map 
            $$\res^{\CP^n}_{\CP^n \setminus \Int B^{}_{\Delta}} \fc \Lambda_{>0} f \to  \Lambda_{>0} e_1 \oplus \bigoplus\limits_{k=2}^\infty \Lambda_{(0, \Delta]} e_k$$
            is given by
            $$\res^{\CP^n}_{\CP^n \setminus \Int B^{}_{\Delta}}(\lambda f) = \begin{cases}
                T^{\Delta} \lambda e_1, & i = n, \\
                \lambda e_1, & i \leq n-1,
            \end{cases}$$
          for every $\lambda \in \Lambda_{>0}$.
            \item If $n-j\leq i\leq n-1$, then the restriction map
            $$\res^{\CP^n}_{\CP^n \setminus \Int B^{}_{\Delta}} \fc \Lambda_{>0} f \to \Lambda_{(0, (n-i)(1-\Delta)]} e_1 \oplus \Lambda_{>0} e_2 \oplus \bigoplus\limits_{k=3}^\infty \Lambda_{(0, \Delta]} e_k$$
            is given by
            $$\res^{\CP^n}_{\CP^n \setminus \Int B^{}_{\Delta}}(\lambda f) = [\lambda]_{(n-i)(1-\Delta)} e_1 - T^{\Delta-(n-i)(1-\Delta)} \lambda e_2,$$
              for every $\lambda \in \Lambda_{>0}$.
        \end{itemize}

        \item Assume $\Delta > \frac{n}{n+1}$.
        
        \begin{itemize}
            \item If $i = n$, then the restriction map
            $$\res^{\CP^n}_{\CP^n \setminus \Int B^{}_{\Delta}} \fc \Lambda_{>0} f \to \bigoplus\limits_{k=1}^\infty \Lambda_{(0, n(1-\Delta)]} e_k$$
            is the zero map.
            \item If $i \leq n-1$, then the restriction map
            $$\res^{\CP^n}_{\CP^n \setminus \Int B^{}_{\Delta}} \fc \Lambda_{>0} f \to \Lambda_{(0, (n-i)(1-\Delta)]} e_1 \oplus \bigoplus\limits_{k=2}^\infty \Lambda_{(0, n(1-\Delta)]} e_k$$
            is given by
            $$\res^{\CP^n}_{\CP^n \setminus \Int B^{}_{\Delta}}(\lambda f) = [\lambda]_{(n-i)(1-\Delta)} e_1 - T^{\Delta-(n-i)(1-\Delta)} [\lambda]_{n(1-\Delta)} \sum\limits_{k=2}^\infty (-T^{\Delta-n(1-\Delta)})^{k-2} e_k,$$
              for every $\lambda \in \Lambda_{>0}$.
        \end{itemize}
    \end{enumerate}  
\end{prop}

\begin{rem}\phantom{M}
    \begin{itemize}
        \item The tools we use for these computations can be used for computing the restriction maps between the relative symplectic cohomology of the complements of two standard toric balls. 
        \item It would be very nice to have a computation of the restriction map from the relative symplectic cohomology of a complement of a toric ball to that of a toric ball contained in it.
    \end{itemize}
\end{rem}

\subsection{Stable displaceability and computational results}\label{ss: sd_energy and application}

Relative symplectic cohomology with Novikov ring coefficients provides a more precise characterization of a subset's symplectic properties. In his thesis \cite{Varolgunes_2018_PhD}, Varolg\"une\c s demonstrates that quantitative information can be extracted from this invariant.

Given a module $A$ over the Novikov ring $\Lambda_{\geq0}$, its capacity is defined to be
$$c(A)=\inf\{r>0\,|\,T^r A=0\}.$$
Note that $c$ takes values in $[0,+\infty]$. The module $A$ is called \textbf{torsion} if $A\otimes \Lambda = 0$, in particular, this holds when $c(A)<+\infty$. For example, the capacity of $\Lambda_{>0}$ is $+\infty$, in contrast to the capacity of $\Lambda_{(a,b]}$, for some $a<b$, which is $b-a$.

We caution the reader that the term ``capacity'' is used in two distinct senses: for modules, it refers to the capacity defined above, whereas for subsets, as stated earlier, it refers to the Gromov width of the interior.

In his thesis \cite[Remark 4.2.8]{Varolgunes_2018_PhD}, Varolg\"une\c s proved that given a closed symplectic manifold $(M,\omega)$ and a compact subset $K\subset M$, the displacement energy of $K$ is bounded from below by the capacity of $SH_M^*(K;\Lambda_{\geq0})$.
In particular, if $K$ is displaceable then $SH_M^*(K;\Lambda_{\geq0})$ is torsion, and if $c(SH_M^*(K;\Lambda_{\geq0}))=+\infty$ then $K$ is non-displaceable\footnote{Due to \cite{MSV_2024_heavy_sets_and_SH}, it follows that $c(SH_M^*(K;\Lambda_{\geq0}))=+\infty$ if $K$ is a heavy subset in the sense of \cite{EP_2009_rigid_subsets}.}. 
In addition, Varolg\"une\c s proved that stable displaceability implies that $SH_M^*(K;\Lambda_{\geq0})$ is torsion. We complete the picture as follows:

\begin{thm}\label{thm: c_relSH vs e_st}
    Let $(M,\omega)$ be a closed symplectic manifold and $K\subset M$ be a compact subset. Then the stable displacement energy of $K$ is bounded from below by the capacity of $SH_M^*(K;\Lambda_{\geq0})$.
\end{thm}

The following result is an immediate application of Theorem~\ref{thm: c_relSH vs e_st} and the computation we made of the relative symplectic cohomology of the complement of a ball in $\CP^n$.

\begin{thm}\label{thm: stbly_disp_enrgy_CP^n-B}
    Let $B^{}_{\Delta} \subset \CP^n$ be a standard toric ball of Gromov width $\Delta \in (0,1)$. If $\Int B^{}_{\Delta}$ contains the Clifford torus, i.e., $\Delta > \frac{n}{n+1}$, then the stable displacement energy of $\CP^n \setminus \Int B^{}_{\Delta}$ is at least $n(1-\Delta)$.
\end{thm}

\begin{proof} 
Assume that the capacity $\Delta$ of the open ball $B$ satisfies $\Delta>\frac{n}{n+1}$. By Theorem~\ref{thm: SH of complement of ball in CPn} we know that 
  $$SH^*_{\CP^n}(\CP^n \setminus B^{}_{\Delta}; \Lambda_{\geq 0}) = \begin{cases}
        0, &  * \text{ is odd}, \\
        \Lambda_{(0,  (n-i)(1-\Delta)]} \oplus \bigoplus\limits_{i=0}^\infty\Lambda_{(0,n(1-\Delta)]}, &  * = 2i,\,\, i\in\{1,\ldots,n-1\}, \\
        \bigoplus\limits_{i=0}^\infty\Lambda_{(0, n(1-\Delta)]}, &  * \in\{ 0,2n\},
    \end{cases}$$
    and hence the capacity of $SH^*_{\CP^n}(\CP^n \setminus \Int B^{}_{\Delta})$ is $n(1-\Delta)$. By Theorem~\ref{thm: c_relSH vs e_st} we deduce that the stable displacement energy of $\CP^n \setminus \Int B^{}_{\Delta}$ is at least $n(1-\Delta)$, as required.
\end{proof}

\begin{rem}$\,$
\begin{itemize}
\item The set $\CP^n\setminus \Int B^{}_{\Delta}$ where $\Delta > \frac{n}{n+1}$ is stably displaceable due to \cite[Theorem 1.9]{EP_2009_rigid_subsets}.
    \item To the best of our knowledge, this result is new and improves upon the previously known lower bound of $1-\Delta$. This latter bound is obtained by noting that the Gromov width of $\CP^n \setminus B$ is at least $1-\Delta$; thus, by \cite[Corollary 1.3]{Usher_2010_energy_capacity_inequality}, the stable displacement energy is at least $1-\Delta$. 
\end{itemize}
 
 \end{rem}

\begin{rem}
In joint work with Frol Zapolsky, currently in preparation, we provide further applications of the capacity of the relative symplectic cohomology.
\end{rem}

\subsubsection*{Boundaries of Balls---The Use of Mayer--Vietoris}

Propositions~\ref{prop: res from CPn to a ball} and \ref{prop: res from CPn to the complement a ball} provide a precise formula for the restriction map from the relative symplectic cohomology of $\CP^n$ to that of a ball, and to that of the ball's complement, respectively. Both balls and their complements are preimages of closed subsets under the standard moment map on $\CP^n$. Consequently, these sets commute (see \cite[Theorem 1.3.7]{Varolgunes_2021_MV_and_relSH} and \cite[Remark 1.14]{DGPZ_2024_Symp_top_and_IVMs}). We can therefore consider the corresponding Mayer--Vietoris exact triangle for a toric closed ball and the complement of a toric open ball, provided they form a cover of $\CP^n$. Since the relative symplectic cohomology of these sets vanishes in odd degrees and the restriction maps are known, a purely algebraic computation yields the relative symplectic cohomology of the intersection. This intersection is a shell---specifically, the difference between a closed toric ball and a smaller open toric ball. While we do not carry out the computation in this paper, it is entirely within reach. 

We give a full statement in the case of ball boundaries.

\begin{thm}\label{thm: SH of ball boundary in CPn}
Let $B^{}_{\Delta}$ be the standard closed ball of capacity $\Delta \in (0,1)$ in $\CP^n$, and denote by $\partial B^{}_{\Delta}$ its boundary. The relative symplectic cohomology $SH^*_{\CP^n}(\partial B^{}_{\Delta}; \Lambda_{\geq 0})$, as a $\Z/2(n+1)\Z$-graded module over $\Lambda_{\geq 0}$, is given explicitly as follows:

\begin{enumerate}[labelindent=0pt,
leftmargin=*,
align=left]
    \item If $\Delta < \frac{n}{n+1}$:
    \[\hspace{-1em}SH^*_{\CP^n}(\partial B^{}_{\Delta}; \Lambda_{\geq 0}) = \begin{cases}
        0, & * \text{ is odd}, \\ \bigoplus\limits_{i=0}^\infty\Lambda_{(0, \Delta]}, &   * \text{ is even}.
    \end{cases}\]

    \item If $\Delta = \frac{n}{n+1}$:
    \[\hspace{-1em}SH^*_{\CP^n}(\partial B^{}_{\Delta}; \Lambda_{\geq 0}) = \begin{cases}
        0, & * \text{ is odd}, \\ \Lambda_{>0}\oplus \bigoplus\limits_{i=0}^\infty\Lambda_{(0, n/(n+1)]}, &   * \text{ is even}.
    \end{cases}\]

    \item If $\Delta > \frac{n}{n+1}$:
    \[\hspace{-1em}SH^*_{\CP^n}(\partial B^{}_{\Delta}; \Lambda_{\geq 0}) = \begin{cases}
        0, & * \text{ is odd}, \\  \bigoplus\limits_{i=0}^\infty\Lambda_{(0,n(1-\Delta)]}, &   * \text{ is even}.
    \end{cases}\]

\end{enumerate}
\end{thm}

As a corollary we obtain
\[
c(SH^*_{\CP^n}(\partial B^{}_{\Delta}; \Lambda_{\geq 0})) = 
\begin{cases}
        \Delta, &  \Delta < \frac{n}{n+1}\\
        +\infty, &  \Delta = \frac{n}{n+1}\\
        n(1-\Delta), &  \Delta > \frac{n}{n+1}
\end{cases}.
\]
We thus obtain lower bounds on the stable displacement energy of the toric ball's boundary in $\CP^n$.
\begin{thm}\label{thm: stbly_disp_enrgy_CP^n-BdryB}
    Let $\partial B^{}_{\Delta} \subset \CP^n$ be the boundary of the standard toric ball of Gromov width $\Delta \in (0,1)$. It is stably displaceable unless it contains the Clifford torus\footnote{The set $\partial B^{}_{\Delta} = \partial (\CP^n\setminus B^{}_{\Delta})$, and when $\Delta\neq\frac{n}{n+1}$ either the ball or its complement is stably displaceable due to \cite[Theorem 1.9]{EP_2009_rigid_subsets}.}
    , i.e., unless $\Delta = \frac{n}{n+1}$, and its stable displacement energy is bounded from below by:
    \begin{alignat*}{2}
    &\Delta, &\qquad\qquad &\text{if } \Delta<\frac{n}{n+1}, \text{ and,}\\
    &n(1-\Delta),&\qquad\qquad &\text{if } \Delta>\frac{n}{n+1}
    \end{alignat*}
\end{thm}

Moreover, ideally, inductive use of Mayer--Vietoris allows for further computation of relative symplectic cohomology of more subsets. We demonstrate the result in dimension $2$:
\begin{exam}
    In $S^2 = \CP^1$, consider $k$ smooth closed disks $D_1, \ldots, D_k$ such that their complements are pairwise disjoint. Their intersection is a sphere with $k$ disks removed, denoted by $\Sigma$. Let $\Delta_1, \ldots, \Delta_k$ denote their capacities, and assume that $\Delta_1$ is the smallest. A direct computation using Mayer--Vietoris shows that:
 
 \begin{itemize}
      \item If $\Delta_1 < 1/2$, then $$SH^0_{S^2}(\Sigma; \Lambda_{\geq 0}) =  \bigoplus_{j=0}^\infty\Lambda_{(0, \Delta_1]} \oplus \bigoplus_{i=2}^k \left(\bigoplus_{j=0}^\infty \Lambda_{(0, 1-\Delta_i]}\right).$$
     \item If $\Delta_1 \geq 1/2$, then $$SH^0_{S^2}(\Sigma; \Lambda_{\geq 0}) = \Lambda_{>0} \oplus \bigoplus_{i=1}^k \left(\bigoplus_{j=0}^\infty \Lambda_{(0, 1-\Delta_i]}\right).$$
 \end{itemize}
 
 Similarly, we can compute the relative symplectic cohomology of $\Sigma$ in all other degrees, as well as the restriction maps $\res \fc SH_{S^2}^*(D_i; \Lambda_{\geq 0}) \to SH_{S^2}^*(\Sigma; \Lambda_{\geq 0})$ for every $1 \leq i \leq k$.
 \end{exam}

\subsection{Review--relative symplectic cohomology}

Let $(M,\omega)$ be a closed symplectic manifold. The relative symplectic cohomology is a contravariant functor from the category of compact subsets of $M$ with respect to inclusion to the category of modules over the Novikov ring. Let us describe how it acts on objects. For a full definition see Section~\ref{ss: def of relSH}.
Given a compact subset $K \subset M$, the \textbf{relative symplectic cohomology of $K$ in $M$} is defined as
\begin{equation}\label{eq:def_relSH w lim}
    SH^*_M(K;\Lambda_{\geq0}) = H^*\Big(\wh{\varinjlim\limits_{i}}\, CF^*(H_i) \Big)\,,
\end{equation}
where $H_i \in C^\infty(S^1 \times M)$ is a pointwise increasing sequence of non-degenerate time-periodic Hamiltonians on $M$ such that $H_i|_K < 0$, and such that
$$
\lim_{i\to\infty} H_i(t,x) = \begin{cases}
0, & x \in K, \\
+\infty, & x \notin K,
\end{cases}
$$
for every $(t,x) \in S^1 \times M$. Such a sequence, together with homotopies between consecutive Hamiltonians, is called an \textbf{acceleration datum for $K$}.

Here, $CF^*(H_i)$ denotes the Floer complex
\begin{equation}\label{eqn:defin_Floer_cx}
    CF^*(H_i) = \bigoplus_{x \in \cP^\circ(H_i)} \Lambda_{\geq 0} \cdot x\,,
\end{equation}
where $\cP^\circ(H_i)$ is the set of 1-periodic orbits of $H_i$ that are contractible in $M$. The complexes $CF^*(H_i)$ are connected by Floer continuation maps (See Section~\ref{sss: weighted_CF}). The hat denotes the $T$-adic completion of $\Lambda_{\geq 0}$-modules (see Section~\ref{ss: completion}). 

The grading of a generator $x \in \cP^\circ(H_i)$ is given by a Conley--Zehnder type index which is well-defined modulo $2N_M$, where $N_M$ is the minimal first Chern number of $(M,\omega)$.

\begin{rem}
The cohomology $SH_M^*(K;\Lambda_{\geq0})$ is independent of the specific choice of Hamiltonians $H_i$ and other auxiliary data. To establish this, and more generally to work with $SH_M$, it is often advantageous to use an alternative definition based on homotopy theory; see Section~\ref{ss: def of relSH} and \cite{Varolgunes_2021_MV_and_relSH} for further details.
\end{rem}

\subsection{A brief review of the proof of Theorem \ref{thm: SH of ball in CP^n}}\label{ss: brief review of proof}
We start with the following acceleration data. Consider the standard moment map $\Phi\colon \CP^n \to \R^n$, and let $B_\Delta$ denote the toric ball of capacity $\Delta$, whose image under $\Phi$ is the corner simplex of side length $\Delta$.
We set
$H_\ell(z) = h_\ell\left( \left\langle \Phi(z), \left(1,1,\ldots,1\right)\right\rangle\right)$ for every $z\in \CP^n$, where $h_\ell$ has a graph as depicted in Figure~\ref{fig:fig1} and Figure~\ref{fig:fig2} below.

\begin{figure}[htbp]
    \centering

    \begin{subfigure}{0.44\textwidth}
        \centering
        \includegraphics[width=\linewidth]{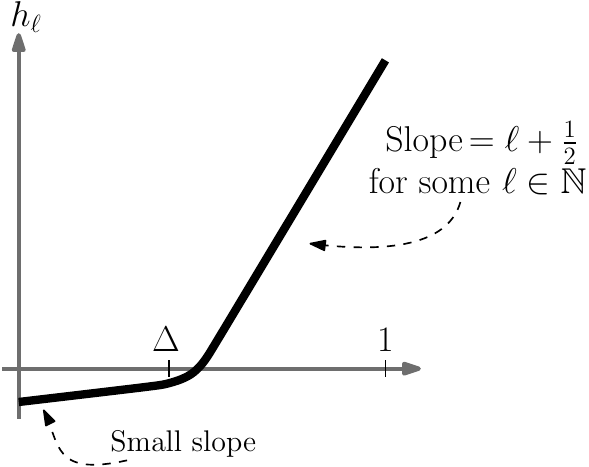}
        \caption{The function $h_\ell$}
        \label{fig:fig1}
    \end{subfigure}
    \hfill
    \begin{subfigure}{0.52\textwidth}
        \centering
        \includegraphics[width=\linewidth]{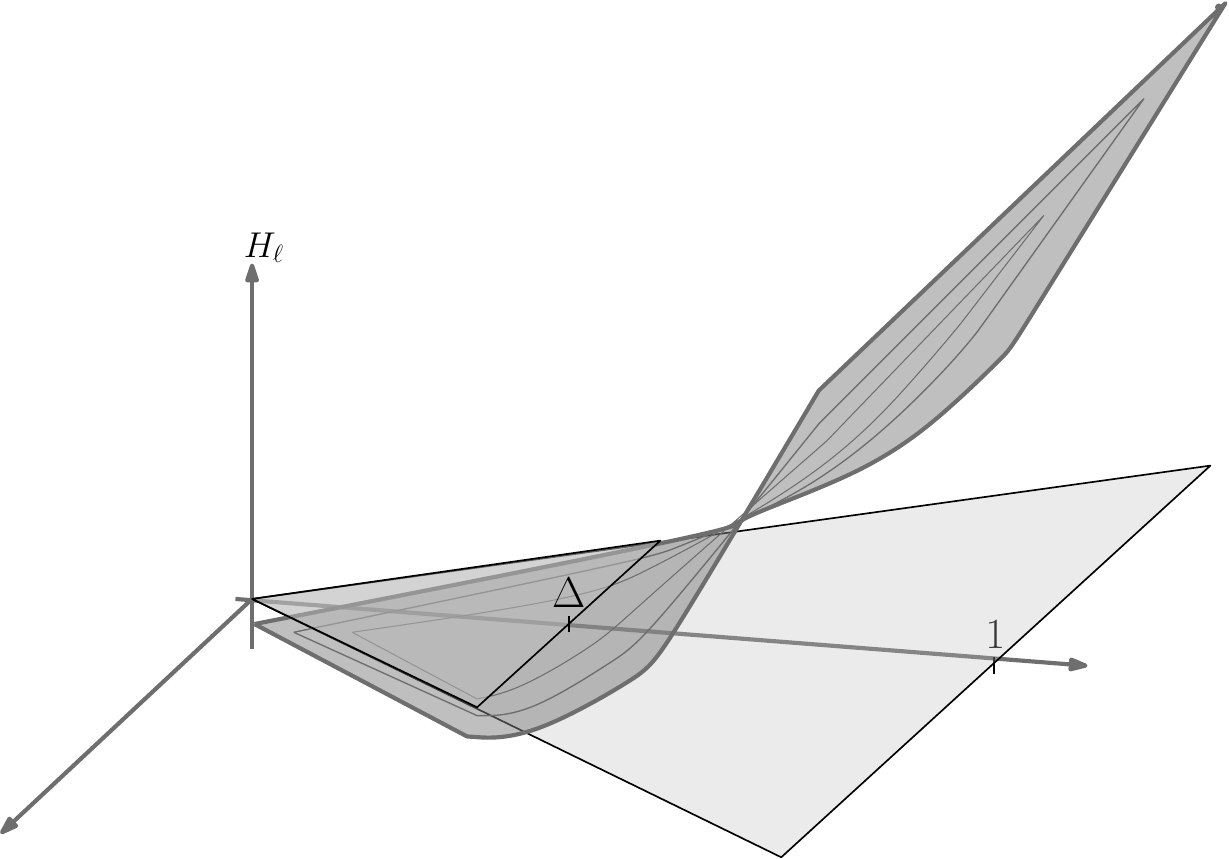}
        \caption{The Hamiltonian $H_\ell$ on the moment picture of $\CP^n$}
        \label{fig:fig2}
    \end{subfigure}

    \label{fig:side_by_side}
\end{figure}

Ideally we would like to construct acceleration data by using the Hamiltonians $\left(H_{n}\right)_{n=0}^\infty$ and choosing homotopies from each $H_n$ to the next, $H_{n+1}$, by a monotone interpolation, (see Section~\ref{sss: relSH} for a precise definition of the notion of an acceleration datum), and for an almost complex structure to use $J_{std}$, the standard complex structure of $\CP^n$. However, to achieve regularity we have to perturb this ideal data: We perturb the almost complex structure $J_{std}$ for the definition of the Floer complexes of $H_n$, and we perturb the homotopies between them for the definition of the continuation maps, which enables us to construct  a direct system of Floer complexes. See the start of Section~\ref{s: computations of CF} for a discussion on the perturbed data.

In the setting of \cite{Varolgunes_2021_MV_and_relSH} that we use, the differentials in Floer complexes we use are weighted by topological energy of the solutions encoded in the exponent of the Novikov parameter. We show that each matrix element of the differential, $\langle dx,y \rangle$, is of the form $AT^E$ where $A\in \Q$, $T$ is the Novikov parameter and $E \in [0,+\infty)$. A similar property holds for a continuation morphism.

This is due to Equation~\eqref{eq: energy}, which shows that in monotone symplectic manifolds, the energy of Floer and continuation solutions is determined by the index difference of the $1$-periodic orbits. Thus, we may also define Floer complexes and continuation maps over $\Z$, equivalently, setting $T=1$, and we compute these first.

Our method of computing the differential involves two stages. First we show that most matrix elements vanish, by proving that the corresponding moduli spaces are empty. Then, we compute the other matrix elements either by identifying them with the matrix elements of a local computation in $\C^n$, or deduce them by an argument using the fact that continuation maps are quasi-isomorphisms.

To prove that the moduli spaces are empty, we first prove that they are empty for the unperturbed acceleration data, then show that emptiness persists under small enough perturbations. See Section~\ref{ss: transversality}.

Our methods require certain symmetries of the Floer data (e.g., a $\U(n)$ symmetry, crucial to our argument\label{location: U(n) symmetry}), which necessitates the use of autonomous Hamiltonians, and thus a Floer theory based on cascades. In this framework, the Hamiltonian is autonomous, and its $1$-periodic orbits form a collection of smooth, finite-dimensional submanifolds of the loop space. By choosing a Morse function on each manifold, we define a chain complex generated by their critical points. The differential in this setting counts ``cascades'': trajectories consisting of Floer cylinders interspersed with gradient segments of the Morse functions, introduced by Frauenfelder, see \cite{Frauenfelder_2004_MB}. See Section~\ref{s: MB in Floer} for a review of Morse--Bott theory in Floer theory, using cascades.

For non-degenerate Hamiltonians, the cascades theory coincides with the standard Floer theory. Consequently, we demonstrate that relative symplectic cohomology may also be computed using Floer complexes defined via cascades, and therefore we may use this model for our computations.

Let us now describe the $1$-periodic orbits generating $CF(H_\ell;\Z)$.
The $1$-periodic orbits of $H_\ell$ come in three types of families.
\begin{enumerate}
    \item The global minimum at $0$,
    \item $\ell$ spheres $S^{2n-1}$ consisting of the points where the slope is integral,
    \item and the global maxima, forming a copy of $\CP^{n-1}$, at the divisor at infinity.
\end{enumerate}
The first family contributes a generator which we denote by $\check{x}_0$, each sphere contributes a pair of generators corresponding to the minimum and maximum of a perfect Morse function, yielding in total $\check{x}_1,\hat{x}_1,\ldots,\check{x}_\ell,\hat{x}_\ell$, where the hats correspond to maxima and the checks to minima. Lastly, the divisor at infinity contributes generators, which we denote by $\check{x}_{\ell+1},\ldots,\check{x}_{\ell + n}$.
The acceleration data and the generators are discussed in Section~\ref{ss: acc. data}.

We consider two symplectic trivializations of $TM$ along the non-constant generators, which we use in index computations. The first, denoted $\tau_B$, is obtained by trivializing $TM$ over $\CP^n \setminus D_\infty$, which is symplectomorphic to a ball. The second, denoted by $\tau_T$ is obtained by identifying $\CP^n \setminus \left(D_1 \cup \ldots \cup D_n \cup D_\infty\right)$ with an open set in $T^*\T^n$, and trivializing the tangent bundle there. These trivializations are discussed in Section~\ref{ss: charts}.

The trivialization $\tau_T$ has the crucial property that the indices of all the non-constant generators are either $0$ or $2n-1$, regardless of winding.

The main tools we use to prove the vanishing of matrix elements are:
\begin{description}
    \item[Index formula:] Given a symplectic trivialization $\tau$ of $TM$ along the generators, and $u$ a flowline with cascades $u_1,\ldots,u_m$, the local dimension of the moduli space at $u$ is computed by the following formula involving index differences and relative first Chern class:
    \begin{equation}\label{eq: index_formula_intro}
         \mu^\tau_{FMB}(q_+;H_\ell) - \mu^\tau_{FMB}(q_-;H_\ell) + 2\sum_{i=1}^m c_1^\tau(u_i).
    \end{equation}
    See Remark~\ref{rem: index formula for Floer and cont flowlines}.
    \item[Intersection formula:] If $\tau$ is $\tau_B$ or $\tau_T$, then the first relative Chern number can be computed in terms of the intersection number with the divisors $D_\infty$ or $D_1\cup \ldots \cup D_n \cup D_\infty$ respectively. This is developed in Section~\ref{ss: c_1 and intersection}
    \item[Positivity of intersections:] If the Hamiltonian vector field is tangent to a complex divisor, then Floer solutions intersect positively with it, see \cite[Lemma 4.3]{Ganatra_Pomerleano_2021_log_PSS}. Moreover when the asymptotes of a Floer solution are fixed points lying on the divisor, one can define asymptotic intersection numbers, for which, under certain assumptions, a lower bound can be obtained. This was proven by Seidel in \cite[Eq. (7.22)]{Seidel_Fukaya_A_infty_Strcs_Assoc_to_Lef_Fib. III} and is slightly refined here in Proposition~\ref{prop: asymp_intersection}. With Floer flowlines with cascades, we encounter a particular type of asymptotic intersections, where two consecutive cascades are asymptotic, one at $\infty$ and the other at $-\infty$, to a divisor consisting of fixed points. We show that such intersections are positive in total, see Lemma~\ref{lemma: estimate for int asymptc intersection}. 
\end{description}
We prove that the following matrix elements vanish, see Theorem~\ref{thm: diff_obst}:
\begin{itemize}
    \item $\langle d\hat{x}_i, \hat{x}_j\rangle=0$ for every $1 \leq i, j \leq \ell$.
    \item $\langle d\check{x}_i, \check{x}_j\rangle=0$ for every $0 \leq i, j \leq \ell + n$.
    \item $\langle d\check{x}_i, \hat{x}_j\rangle=0$ for every $0 \leq i \leq \ell + n$ and $1 \leq j \leq \ell$, except where $i=j$ and $n=1$.
    \item $\langle d\hat{x}_i, \check{x}_j\rangle=0$ for every $1 \leq i \leq \ell$ and $0 \leq j \leq \ell + n$, except where $j \in \{i-1, i+n\}$. 
 
\end{itemize}

To demonstrate our methods, we sketch the proof of the fourth item, in the case $1\le i,j \le \ell$.

First, by the index formula with respect to the trivialization $\tau_B$, substituting the indices computed in Table~\ref{tab:mu-fmb-values}, from Section~\ref{prop: traj and modular indices}, into the index formula~\eqref{eq: index_formula_intro}, we obtain
\[ (n+1) u\cdot D_\infty = n(j+1-i),\]
where $u\cdot D_\infty$ is the total intersection number of $u$ with the divisor $D_\infty$. Since $n$ and $n+1$ are coprime, there exists $k\in \Z$ such that 
\begin{equation}\label{eqn:keqIntro}
    u\cdot D_\infty = kn\qquad\qquad\text{and}\qquad\qquad
    j+1-i = k(n+1).
\end{equation}

The positivity of intersections (see Section~\ref{ss: intersections}), implies 
\[ u\cdot (D_1\cup \ldots \cup D_n \cup D_\infty) \ge  u\cdot D_\infty \ge 0.\]

By the index formula applied with respect to $\tau_T$, combined with the above inequality, we deduce $ u\cdot D_\infty \le n$, and since $0\le u\cdot D_\infty = kn$ we obtain $0\le kn \le n$. Thus either $k=0$ or $k=1$, and Equation \eqref{eqn:keqIntro} implies that $j$ is either $i-1$ or $i+n$, respectively.

As mentioned previously, we define our acceleration data by perturbing $J_{std}$ and the homotopies between each $H_n$ and $H_{n+1}$. The proof that such perturbations preserve the emptiness of moduli spaces, involves a compactness argument, detailed in Section~\ref{ss: transversality}, and bounds on energy and the number of cascades, detailed in Section~\ref{ss: bounds for energy and cascades}.

Thus, the complex of $H_n$ may be depicted as follows. The dots represent the generators, and the dashed lines represent the matrix elements that at this point are yet unknown (the rest are zero).

\begin{figure}[H]
    \centering
    \includegraphics[width=\linewidth]{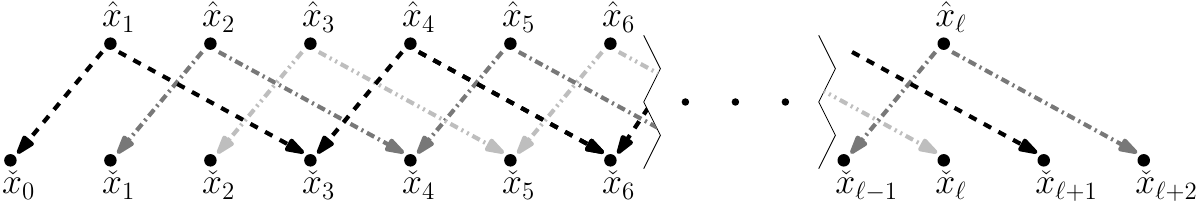}
    \caption{The complex of $H_\ell$ on $\CP^2$. The different dash styles and shades of gray are used as a visual cue to separate the subcomplexes}
    \label{fig:figfullcomplex}
\end{figure}

The complex thus splits as a direct sum of $n+1$ subcomplexes. Since all of these subcomplexes have the same form, it suffices to analyze a single summand; the analysis for the remaining summands is identical.
\begin{figure}[H]
    \centering
    \includegraphics[width=0.95\linewidth]{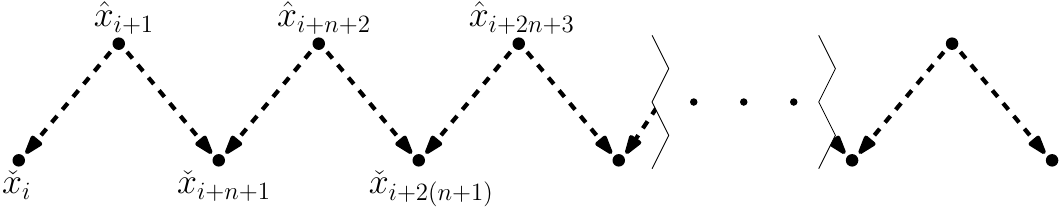}
    \caption{A subcomplex in the direct sum decomposition. These subcomplexes are indexed by $i$, $0\le i \le n$, where each zigzag starts at $\check{x}_i$.}
    \label{fig:figSubcomplexA}
\end{figure}

Those arrows going to the left correspond to matrix elements of the type $\langle d\hat{x}_i, \check{x}_{i-1}\rangle$, for which we show in Theorem~\ref{thm: diff_obst} that the Floer solutions do not intersect $D_\infty$, hence, they equal the same matrix elements, computed in $\C^n$. To that end we employ a ``no escape lemma'' for Floer solutions, similar to the one from \cite{Ritter_TQFT}, which requires adaptations to our setting since the conditions are somewhat different. This is proved in Section~\ref{ss: no_escape}.

The computation of these matrix elements in $\C^n$ was carried out by Oancea in \cite{Oancea_survey} and is reproduced here in the setting of cascades in Section~\ref{ss: CF of J-shaped}. We conclude that they are all equal to either $1$ or $-1$.

Thus, each of the $n+1$ subcomplexes takes the form
\begin{figure}[H]
    \centering
    \includegraphics[width=\linewidth]{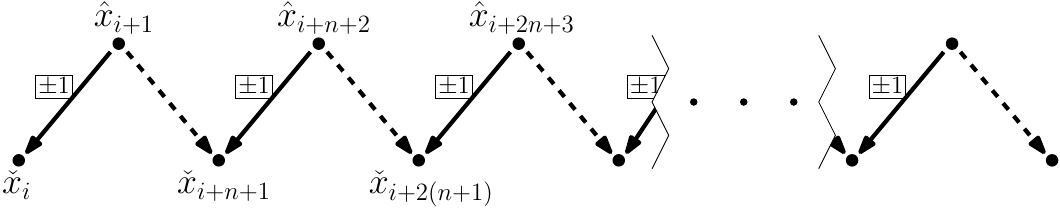}
    \caption{The $i$-th zigzag, together with the information obtained from $\C^n$.}
    \label{fig:figSubcomplexB}
\end{figure}

By a similar combination of index formulas for the trivializations $\tau_B$ and $\tau_T$, and the positivity of intersections, we deduce the vanishing of most matrix elements of the continuation map $\varphi\colon CF(H_\ell;Z)\to CF(H_{\ell'};\Z)$ for $\ell' \ge \ell$. In fact, all matrix elements vanish except for $\langle \varphi\check{x}_i, \check{x}_i\rangle$ and $\langle \varphi\hat{x}_i, \hat{x}_i\rangle$. See Theorem~\ref{thm: continuation_obst}.

\begin{figure}[H]
    \centering
    \includegraphics[width=\linewidth]{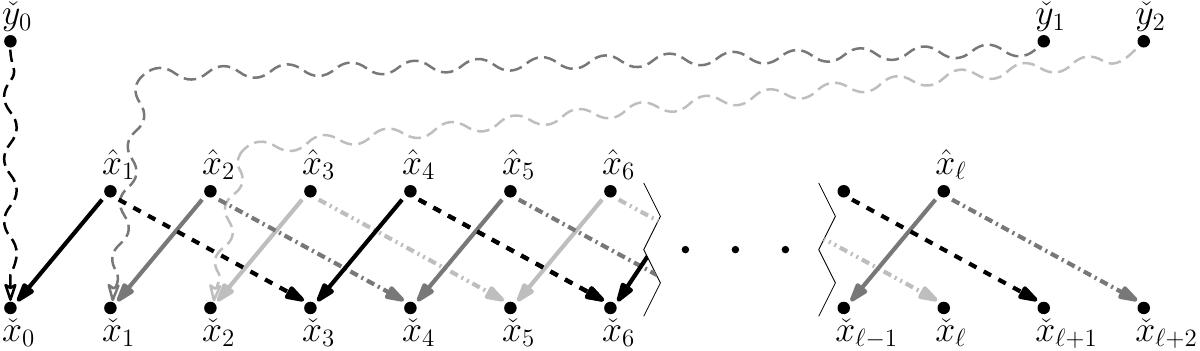}
    \caption{The continuation $\varphi\colon CF(H_0;Z)\to CF(H_\ell;\Z)$ in $\CP^2$. Note the leftmost generator is the generator of $CF(H_0;\Z)$ corresponding to the global minimum, and the pair on the right are the generators corresponding to the global maximum at the $\CP^1$ at infinity.}
    \label{fig:FullComplex1ContFromH0}
\end{figure}

Now, $H_0$ provides exactly $n+1$ generators of different grading modulo $2(n+1)$, and $HF(H_0;\Z)\cong QH(\CP^n;\Z)$, thus the differential in $CF(H_0;\Z)$ vanishes. We now consider the continuation map $\varphi\colon CF(H_0;Z)\to CF(H_\ell;\Z)$. The situation is depicted in Figure~\ref{fig:FullComplex1ContFromH0}.

\begin{figure}[H]
    \centering

    \begin{subfigure}{\textwidth}
        \centering
        \includegraphics[width=\linewidth]{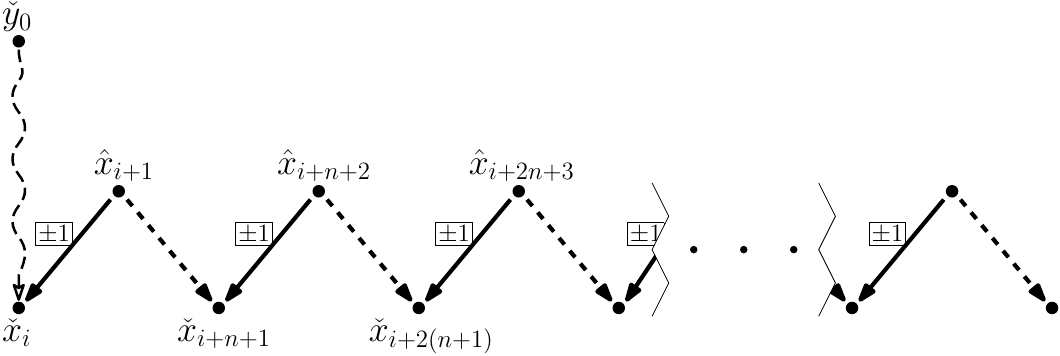}
        \caption{The continuation from $H_0$ forces all dashed arrows to be $\pm1$.}
        \label{fig:SubcomplexBContFromH0}
    \end{subfigure}
    
    \vspace{1em}
    $\Big\Downarrow$
    \vspace{1em}
    
    \begin{subfigure}{\textwidth}
        \centering
        \includegraphics[width=\linewidth]{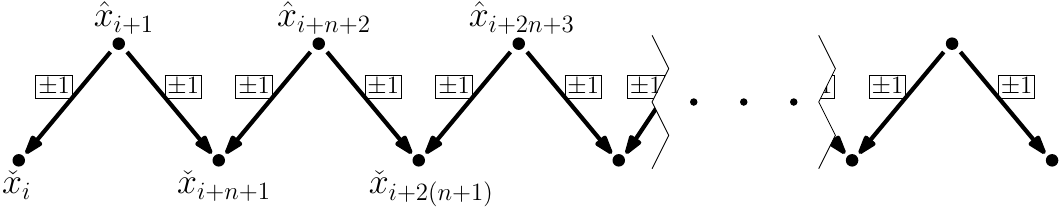}
        \caption{The resulting subcomplex.}
        \label{fig:SubcomplexC}
    \end{subfigure}

    \caption{Differential of Floer subcomplex.}
    \label{fig:ResultingSubcomplex}
\end{figure}

Since $\CP^n$ is a closed symplectic manifold, the continuation map $\varphi\colon CF(H_0;\Z)\to CF(H_\ell;\Z)$ is a quasi-isomorphism, in particular, for all $0\le i \le n$, $\varphi(x_i)$ is closed, non-exact, and in fact, a represents a generator $HF(H_\ell;\Z)$ in degree $2i$ (Figure~\ref{fig:SubcomplexBContFromH0}). As a consequence of Lemma \ref{lemma: zigzag complex}, an algebraic result on the structure of certain zigzag complexes, that utilizes properties of the $\Z$ coefficients, we deduce that $CF(H_\ell;\Z)$ splits as a direct sum of $n+1$ subcomplexes where all the nonzero matrix elements are either 1 or -1, see Figure~\ref{fig:SubcomplexC}.

This argument is detailed in Section~\ref{ss: proof of computations of CF + continuations}.

Next, we return to the complexes over the Novikov ring, using Formula~ \eqref{eq: energy}, relating indices and energy, to compute the topological energy of the solutions.
After a change of basis consisting only of multiplications by $+ 1$, see Proposition~\ref{prop: omission of signs for CF}, we end up with $n+1$ copies of the following complex:

\begin{figure}[H]
    \centering
    \includegraphics[width=\linewidth]{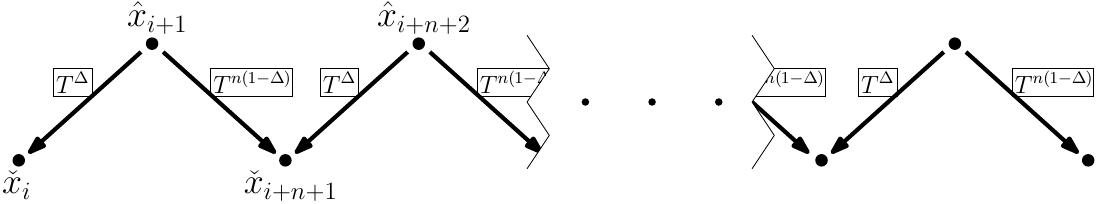}
    \caption{The $i$-th subcomplex with Novikov coefficients weighted by energy.}
    \label{fig:figSubcomplexD}
\end{figure}

The continuation maps $\varphi\colon CF(H_\ell;\Z)\to CF(H_{\ell'};\Z)$ act by $\varphi(\check{x}_i)=\check{x}_i$ and $\varphi(\hat{x}_i)=\hat{x}_i$. Thus in the limit we obtain $n+1$ copies of a zigzag complex, infinite on one side:

\begin{figure}[H]
    \centering
    \includegraphics[width=\linewidth]{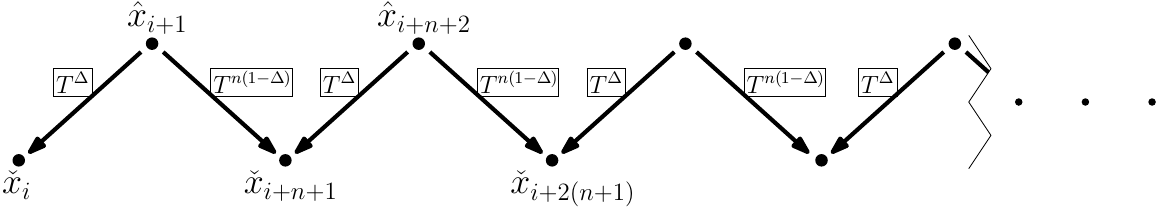}
    \caption{The limit is an infinite zigzag}
    \label{fig:figSubcomplexELimit}
\end{figure}

Now, we take the completion of the complex. All the bottom generators, $(\check{x_i})_{i=0}^\infty$, are closed, and, in fact they generate $\ker d$ as a complete module. In order to compute the homology, namely $\ker d \mathop{/} \im d$, we use the fact that the notion of a Schauder basis applies verbatim in the setting of torsion-free Novikov modules. We thus find a Schauder basis for $\ker d$, $(e_i)_{i\in \N}$, such that $\im d$ is generated by $(T^{\alpha_j} e_j)_{j\in K}$ for nonnegative real $0\le\alpha_j$, where $K\subseteq \N$. Hence the homology is a direct sum of interval modules $\Lambda_{(0,\alpha_i]}$ and $\Lambda_{(0,\infty)}$. We note that different Schauder bases are required, depending on whether the left-going arrows or the right-going arrows carry more energy, or equivalently, whether $\Delta < \frac{n}{n+1}$ or $\Delta \ge \frac{n}{n+1}$, leading to different results in homology.
Schauder bases and their properties are reviewed in Section~\ref{sss: Schauder basis}, and the homology calculation, which completes that proof, is carried out in Section~\ref{ss: proof of main thm}.

As for restriction maps from the relative symplectic cohomology of one ball to another, contained in it, the computation relies on three ideas. The first is that continuation maps $\varphi\colon CF(H_\ell;\Z)\to CF(H'_{\ell'};\Z)$ act by $\varphi(\check{x}_i)=\check{x}'_i$ and $\varphi(\hat{x}_i)=\hat{x}'_i$, where $(H_\ell)_{\ell=0}^\infty$ and $(H'_{\ell'})_{\ell'=0}^\infty$ are acceleration data of the two balls. This is obtained by obstructions on continuation maps and the fact that continuation maps are quasi-isomorphism in closed symplectic manifolds. See Theorem~\ref{thm: continuation maps, over Z}. Second, by using of obstructions on the chain homotopy maps, we prove that the continuation maps involved commute, although the homotopies built into the definition of the restriction maps remain unknown to us. This is shown in Proposition~\ref{prop: commuting of cont maps}.
Third, functoriality in relative symplectic cohomology, utilizes a different model rather than the direct limit, namely, the telescope construction, a model of homotopy colimit. Therefore, we do not obtain a morphism between the complexes defined as completed limits. Thus, we have to compute the restriction maps in cohomology indirectly. Then using an algebraic lemma, see Proposition~\ref{prop: algebraic preparation for computing restriction maps}, and since in our particular case, completion induces a surjective map in homology, see Claim~\ref{claim: surjectivity and completion in CF}, we are able to compute the restriction maps on the homology level using the limit models.

    \subsection{Discussion}
    \subsubsection{More computational examples}

   Expanding the library of explicit computations for the relative symplectic cohomology of various closed symplectic manifolds and compact subsets remains a highly desirable goal. We present here several interesting such cases.

    \subsubsection*{Ellipsoids in $\CP^n$}
    Unlike balls, ellipsoids are not $\U(n)$-invariant, and neither are the critical submanifolds associated with compatible acceleration data. For this reason, we cannot displace the $1$-periodic orbits from the relevant divisors, as we do in the case of the ball, as described on Page~\pageref{location: U(n) symmetry}. This leads to a failure of the positivity of intersections; consequently, we cannot establish the obstructions for Floer and continuation trajectories as we do elsewhere in this paper.
    
    \subsubsection*{Balls in other toric manifolds}

A standard ball in a toric manifold is a symplectic ball that is invariant under the toric action on the manifold, where the induced action coincides with the standard toric action on a ball in $\C^n$. Up to the action of $\SL(n, \Z)$, we can assume that the image of a standard ball in a toric manifold of dimension $2n$ is an $n$-simplex of the form
$$ \{(x_1,\dots,x_n)\in \R^n : \sum_{i=1}^n x_i \leq \Delta \text{ and } x_1,\dots,x_n \geq 0\} $$
for some $\Delta>0$. In this case, there is no obstacle to using the tools we present in this paper. However, the process of finding obstructions for the Floer and continuation trajectories is more complicated, as there are more toric divisors, and hence many more generators with a greater variety of Conley-Zehnder indices.

\subsubsection*{Neighborhoods of complex hyperplanes}

 In Claim~\ref{claim: biholo O(1)->CPn-0}, we show that the complement of a point in $\CP^n$ can be seen as the line bundle $\cO(1)$ over $D_\infty = \CP^{n-1}$. Thus, the complement of a ball in $\CP^n$ can be seen as a disk bundle over the divisor $D_\infty$, and the complement of this disk bundle, as a subset of $\CP^n$, is just a ball. 

This situation is a particular case of the Biran decomposition, see \cite{Biran_2001_Lag_barriers_and_symp_emb}. In $\CP^n$ with the standard Fubini-Study metric, it is claimed there that for every divisor $D \subset \CP^n$, there is an embedding of a disk bundle over $D$ into $\CP^n$ such that the complement of the image is an isotropic skeleton. In the case of $D_\infty\subset \CP^n$, the disk bundle is $\cO(1)$ and the complement is just a point. Another important example is the smooth quadric hypersurface:
$$ D = \{ [z_0 : \cdots : z_n] \in \CP^n : z_0^2 + \cdots + z_n^2 = 0 \}. $$
Here, the disk bundle will be $\cO(2)$ and the skeleton will be a Lagrangian copy of $\RP^n$. In general, we can consider a divisor of the form
$$ D = \{ [z_0 : \cdots : z_n] \in \CP^n \,|\, z_0^k + \cdots + z_n^k = 0 \},$$
for some $k\in \N$. The associated disk bundle will be $\cO(k)$ and the skeleton will be an isotropic CW-complex.

It would be highly interesting to find a systematic way to compute the relative symplectic cohomology of disk bundles over such divisors and of their complements.

After the first version of this paper appeared, Jo--Kang \cite{Jo_Kang_2026_Biran_relSH} posted a computation of an equivalent version of relative symplectic cohomology for circle subbundles arising from Biran decompositions of closed monotone integral symplectic manifolds, for divisors of degree less than the monotonicity constant. The two works have different emphases: Jo--Kang focus on computations in this broader geometric setting and their applications to heaviness, while our computation over the Novikov ring focuses on finite Novikov torsion and the relative symplectic cohomology capacity. Like our approach, theirs uses positivity of intersections with a divisor; it also relies on a neck-stretching argument. We believe that their methods could be extended to compute restriction maps, but such an extension would require addressing the same difficulties that arise in our computations, including the handling of higher chain homotopies.

\subsubsection*{Computations using Maurer--Cartan elements}

In a recent work \cite{BES_2025_Maurer_Cartan_elements_in_SH}, Borman--El-Alami--Sheridan use Maurer--Cartan elements to compute quantum cohomology over a Novikov-type ring, denoted by R in their paper, for Liouville domains obtained by removing divisors $D_1,\ldots,D_k$ from closed monotone symplectic manifolds, under the assumption that
\begin{equation}\label{eq: condition for Maurer_Cartan and SH}
2c_1(TM)=\lambda_1\PD(D_1)+\cdots+\lambda_k \PD(D_k),\qquad\text{where}\,\,\lambda_1,\ldots,\lambda_k\in(0,2]\cap \Q.
\end{equation}

They model the quantum cohomology as Floer homology of Hamiltonians closely related to those appearing in acceleration data, and the authors conjecture a connection between relative symplectic cohomology over the Novikov field of the Liouville skeletons and deformations of symplectic homology. It would be highly interesting to establish a connection between their results and relative symplectic cohomology over the Novikov ring; such a connection would likely lead to many new examples of explicit computations of this invariant.

Let us mention, though, that these techniques do not seem to be capable of recovering the results of the present paper, since the complement of a ball in $\CP^n$ is a not exact and has a concave end rather than being a Liouville domain. Moreover, although balls in $\CP^n$ are Liouville domains and can be obtained by removing a copy of $\CP^{n-1}$ from $\CP^n$ and shrinking along the Liouville flow, Equation~\eqref{eq: condition for Maurer_Cartan and SH} above does not hold in this case, since 
\[2c_1(\CP^n)=2(n+1)\PD(\CP^{n-1})\]
and $2(n+1)>2$.

    \subsubsection{K\"unneth formula and computations}

Given two closed symplectic manifolds $M, M'$ and compact subsets $K \subset M$ and $K' \subset M'$, Varolg\"une\c s constructed in \cite[Section 4.3]{Varolgunes_2018_PhD} a K\"unneth morphism $$SH_M^*(K; \Lambda_{\ge 0}) \otimes SH_{M'}(K'; \Lambda_{\ge 0}) \to SH_{M \times M'}(K \times K'; \Lambda_{\ge 0}).$$
In general, this morphism is not an isomorphism, mainly due to the failure of the tensor product to commute with either homology or the completion. Varolg\"une\c s proved that for $K=M$ and $K'=M'$, this morphism is indeed an isomorphism. 

First of all, it would be very interesting to find exact conditions under which this morphism is either injective or surjective, or to embed it in some sort of a K\"unneth short exact sequence. Second, even if we cannot use a direct combination of the K\"unneth morphism and our computations of relative symplectic cohomology, it seems possible to apply our techniques and compute it directly in simple product cases; for instance, for polydisks in $S^2 \times S^2$.

 \subsubsection{Product structure}
Tonkonog--Varolg\"une\c s defined in \cite{Tonkonog_Varolgunes_Super_rigidity_of_certain_skeleta_using_relative_symplectic_cohomology} a product structure on the relative symplectic cohomology with coefficients in the Novikov field. Later, Dickstein--Ganor--Polterovich--Zapolsky generalized the product structure to coefficients in an arbitrary algebra over the Novikov ring using a construction at the chain level; see \cite[Sections 2, 4]{DGPZ_2024_Symp_top_and_IVMs}. Dickstein provided in his thesis, see \cite[Appendix A]{Dickstein_2025_PhD}, more details on this construction, based on counting perturbed pseudo-holomorphic pairs-of-pants. Also, Abouzaid--Groman--Varolg\"une\c s developed a much more general construction of a chain level product, using the language of framed $E_2$ structures; see \cite{AGV_2024_framed_E2_structures}. 

We believe that our tools can be used for computing the product structure induced by the chain-level construction described in \cite{Dickstein_2025_PhD, DGPZ_2024_Symp_top_and_IVMs}. Moreover, using the positivity of intersections, we might be able to isolate the contribution of pairs-of-pants that do not intersect the divisor $D_\infty = \CP^{n-1}$, and thus apply this computation to balls in $\C^n$.

\subsubsection{IVQMs over $\Lambda_{\geq0}$}

In \cite{DGPZ_2024_Symp_top_and_IVMs}, the authors, in collaboration with Polterovich and Zapolsky, introduce the notion of a quantum cohomology ideal-valued measure, called quantum-IVQM. For a closed symplectic manifold $(M,\omega)$, this measure associates to any compact subset $K$ of $M$ an ideal $\tau(K)$ in the quantum cohomology of $(M,\omega)$ with coefficients in the Novikov field, defined as follows: 
$$\tau(K)=\bigcap_{K\subset U\text{ open}}\ker\left(\res\fc SH_M^*(M;\Lambda)\to SH_M^*(M\setminus U;\Lambda)\right).$$
If $K$ is displaceable, then $\tau(K)=0$. Conversely, if $\tau(K)\neq0$, Mak--Sun--Varolg\"une\c s \cite{MSV_2024_heavy_sets_and_SH} proved that $K$ is a heavy set in the sense of Entov--Polterovich \cite{EP_2009_rigid_subsets}. The other direction was proved under certain assumptions in \cite{DGPZ_2024_Symp_top_and_IVMs} and conjectured in general. 

From the perspective of the current work, it is natural to ask whether it is possible to construct a version of the quantum-IVQM over the Novikov ring. The natural candidate, which we denote by $\tau^{}_{\Lambda_{\geq0}}$, is defined as follows:
$$\tau^{}_{\Lambda_{\geq0}}(K)=\bigcap_{K\subset U\text{ open}}\ker\left(\res\fc SH_M^*(M;\Lambda_{\geq0})\to SH_M^*(M\setminus U;\Lambda_{\geq0})\right).$$
It would be appealing if, instead of vanishing for displaceable sets, this yielded some nonzero modules; however, this is not the case. Indeed, let $K$ be a compact subset of $M$ and assume that $\tau(K)=0$. One can show that there exists an open subset $U$ of $M$ containing $K$ such that
$$\ker\left(\res\fc SH_M^*(M;\Lambda)\to SH_M^*(M\setminus U;\Lambda)\right)=0.$$
Since the Novikov field $\Lambda$ is flat, we deduce that
\begin{align*}
    \ker&\left(\res\fc SH_M^*(M;\Lambda)\to SH_M^*(M\setminus U;\Lambda )\right) \\
    &{\phantom{\res\fc SH_M^*(M;\Lambda)\to}}=\ker\left(\res\fc SH_M^*(M;\Lambda_{\geq0})\to SH_M^*(M\setminus U;\Lambda_{\geq0})\right)\otimes_{\Lambda_{\geq0}}\Lambda.
\end{align*}
This means that 
$$\ker\left(\res\fc SH_M^*(M;\Lambda_{\geq0})\to SH_M^*(M\setminus U;\Lambda_{\geq0})\right)\subset SH_M^*(M;\Lambda_{\geq0})$$
is a torsion submodule of $SH_M^*(M;\Lambda_{\geq0})$ and hence it is $0$. 

Nevertheless, if $\tau(K)$ is non-zero, $\tau^{}_{\Lambda_{\geq0}}$ could carry information about the symplectic measurements of $K$. Propositions~\ref{prop: res from CPn to a ball} and \ref{prop: res from CPn to the complement a ball} allow us to compute $\tau^{}_{\Lambda_{\geq0}}$ for closed toric balls and for the complements of open balls in $\CP^n$. If $B^{}_{\Delta}$ is a closed standard toric ball in $\CP^n$ of capacity $\Delta\in(0,1)$, a direct computation yields:
\[
\tau^{}_{\Lambda_{\geq0}}(B^{}_{\Delta})
=
\begin{cases}
0, & \Delta<\frac{n}{n+1}, \\[0.4em]
\displaystyle
\left(\bigoplus_{i=0}^{n-1}
\Lambda_{\geq (n-i)(1-\Delta)}\right)
\oplus \Lambda_{>0},
& \Delta\geq\frac{n}{n+1},
\end{cases}
\]
where the summand $\Lambda_{\geq (n-i)(1-\Delta)}$ lies in degree $2i$, for
$0\leq i\leq n-1$, and the summand $\Lambda_{>0}$ lies in degree $2n$.

Similarly,
\[
\tau^{}_{\Lambda_{\geq0}}
\left(\CP^n\setminus\Int B^{}_{\Delta}\right)
=
\begin{cases}
\displaystyle
\bigoplus_{j=0}^{n}
\left(
\Lambda_{>0}
\cap
\Lambda_{\geq \Delta-j(1-\Delta)}
\right),
& \Delta\leq\frac{n}{n+1}, \\[0.6em]
0,
& \Delta>\frac{n}{n+1},
\end{cases}
\]
where the $j$-th summand lies in degree $2j$, for $0\leq j\leq n$.
It would be very interesting to find a geometric interpretation for these ``starting points of the bars'' of the quantum-IVQM in these cases. Furthermore, such computations could lead to new applications of symplectic rigidity, using techniques similar to those in \cite{DGPZ_2024_Symp_top_and_IVMs}.

\subsubsection{Categorification}

As we know from \cite{Varolgunes_2021_MV_and_relSH}, the relative symplectic cohomology is a functor from the category of compact subsets of a given closed symplectic manifold to the category of modules over the Novikov ring. It would be highly interesting to generalize this functor to the ``category of pairs of a closed symplectic manifold and a compact subset'', i.e., letting the symplectic manifold itself vary as well. The first problem is that this category is not yet established. While a natural candidate for morphism is a Lagrangian correspondence in the product, the composition of such morphisms remains problematic. Nevertheless, using quilted Floer cohomology, see \cite{Wehrheim_Woodward_2-10_Quilted}, one may try to produce morphisms between relative symplectic cohomology modules of different manifolds, induced by a Lagrangian correspondence in their product. Currently, it seems that such questions are mostly open.

    \subsubsection{Torsion in general symplectic manifolds}
Varolg\"une\c s conjectured in his thesis \cite[Remark 3.3.5]{Varolgunes_2018_PhD} that for every closed symplectic manifold $(M, \omega)$ and every compact subset $K$ with non-empty interior (in the sense of Section~\ref{ss: sd_energy and application}), the capacity of the relative symplectic cohomology of $K$ is positive. In this section, we confirm this conjecture for $\CP^n$. 

Indeed, let $K \subset \CP^n$ be a compact subset with non-empty interior. We can find a closed toric ball $B$ such that its image under a symplectomorphism  is contained in $K$. Let $\Delta$ denote the Gromov width of $B$; we may assume that $\Delta < \frac{n}{n+1}$. By Theorem~\ref{thm: SH of ball in CP^n} and Proposition~\ref{prop: res from CPn to a ball}, the image of the restriction map 
$$\res^{\CP^n}_B \fc SH_{\CP^n}^*(\CP^n; \Lambda_{\geq 0}) \to SH_{\CP^n}^*(B; \Lambda_{\geq 0})$$ 
contains a copy of $\Lambda_{(0, \Delta]}$. In particular, there exists an element $x \in SH_{\CP^n}^*(\CP^n; \Lambda_{\geq 0})$ such that 
$$T^{\Delta/2} \res^{\CP^n}_B(x) = \res^{\CP^n}_B(T^{\Delta/2} x) \neq 0.$$ 
By the functoriality of relative symplectic cohomology, the restriction map factors as $\res^{\CP^n}_B = \res_B^K \circ \res^{\CP^n}_K$. Thus, we deduce that for every $\alpha\in[0,\Delta/2]$
$$T^{\alpha} \res^{\CP^n}_K(x) = \res^{\CP^n}_K(T^{\alpha} x) \neq 0.$$ 
This implies that 
$$c(SH_{\CP^n}(K; \Lambda_{\geq 0})) \geq \frac{\Delta}{2} > 0,$$ 
which confirms Varolg\"une\c s's conjecture for $\CP^n$. Furthermore, a slightly more delicate argument shows that the inequality $c(SH_{\CP^n}(K; \Lambda_{\geq 0})) \geq c_G^T(K)$ holds, where $c_G^T(K)$ denotes the supremal Gromov width of a toric ball which can be embedded by an ambient symplectomophism into the interior of $K$. Note that in $\CP^1$ and $\CP^2$ we know that $c_G^T(K)$ is just the Gromov width of the interior of $K$,\footnote{The case of $\CP^2$ is a conclusion from the fact that $\Ham(\CP^2)$ acts transitively on symplectic balls in $\CP^2$. See \cite{McDuff_1993_symplectic _blow_up}.} in a higher dimension this is an open question.

In a joint work with Frol Zapolsky,  currently in preparation, we establish a similar inequality for any compact subset of an arbitrary closed symplectic manifold. We also apply this result to establish a new symplectic rigidity phenomenon, specifically the non-existence of certain symplectic isotopies subject to specific constraints.

\paragraph*{Organization of the paper.} In Section~\ref{s: MB in Floer}, we survey the adaptation of Floer theory to the setting of Morse--Bott with cascades. We start in Section~\ref{ss: convensions} with our conventions and an overview of the Robbin--Salamon index and classical Floer theory. Section~\ref{ss: MB setting} presents the Morse--Bott condition for Hamiltonians and an adapted index for $1$-periodic orbits. Sections~\ref{ss: flowlines and cascades} and \ref{ss:CF_MB_over_Z} concern the definition of cascades and the construction of the differential and continuation maps. This construction is generalized in Section~\ref{ss:cubes} for induced maps from higher homotopies, using the language of cubes.

Section~\ref{s:Charts and acc. data} concerns the geometry of $\CP^n$. We discuss the construction of $\CP^n$ as a symplectic cut and use it to present radial functions, which constitute the prototype of the Hamiltonians used in this paper. In Section~\ref{ss: charts}, we present three charts of $\CP^n$ that will be extremely helpful for computing Robbin--Salamon indices, estimating asymptotic intersection numbers, and constructing the crucial trivializations needed to prove the main results. Section~\ref{ss: acc. data} presents the Hamiltonians we use for the acceleration data for balls in $\CP^n$, and discusses the generators of their Floer--Morse--Bott complexes, which are collected into three groups.

Section~\ref{s: indices} deals with computing the indices of the generators mentioned in the previous paragraph. In Section~\ref{ss: RS part II}, we discuss some useful properties of the Robbin--Salamon index for paths of symplectic matrices. Then, in Section~\ref{ss:RS index for Ham orbits}, we use those properties alongside the three charts discussed previously to compute the Robbin--Salamon index of those generators. Finally, in Section~\ref{ss: computations of FMB}, we compute the Floer--Morse--Bott indices, which adapt the Robbin--Salamon index to the Morse--Bott with cascades setting.

Sections~\ref{s: obstructions} and \ref{s: computations of CF} contain the Floer-theoretic heart of this paper. In Section~\ref{ss: c_1 and intersection}, we survey the connection between the intersection number of divisors with cylinders and the relative first Chern number of the cylinder with respect to trivializations that satisfy certain conditions related to the divisors themselves. We then prove that the trivializations considered in Section~\ref{s:Charts and acc. data} satisfy those conditions. Section~\ref{ss: intersections} investigates the intersection numbers of flowlines with cascades and divisors, including a discussion on positivity of intersections and Seidel's estimate for the asymptotic intersection number. Combined with the index formula, these results allow us to prove, in Section~\ref{ss: obstructions}, the emptiness of several moduli spaces of Floer, continuation, and homotopy flowlines with cascades.

Next, Section~\ref{s: computations of CF} is dedicated to computing the Floer differential and continuation maps associated with our acceleration data. It begins with a discussion of the no-escape lemma (Section~\ref{ss: no_escape}). In Section~\ref{ss: CF of J-shaped}, we use this lemma to show that we can recover an argument for computing the Floer complex in $\C^n$ and apply it to our acceleration data, thus determining the algebraic counting for certain moduli spaces of Floer solutions. In Section~\ref{ss: proof of computations of CF + continuations}, we apply the algebraic results from Section~\ref{ss: alg preperation} to carry out the computations targeted in Section~\ref{s: computations of CF}.

Section~\ref{s: algebra} constitutes the algebraic preparation required to prove the paper's main results. First, in Section~\ref{ss: Novikov modules}, we discuss modules over the Novikov ring and their completions. We then deal with Schauder bases for completed modules and use them to compute quotients; these will be utilized later when determining the homology at the end of the relative symplectic cohomology computations for balls and their complements. Later, in Section~\ref{ss: lim of 1-ray}, we introduce the concept of a telescope---a model of the homotopy colimit of a sequence of chain complexes upon which the definition of relative symplectic cohomology relies. This involves higher homotopies between the Hamiltonians in two different acceleration data sets. In contrast to the relative symplectic cohomology itself, which can be computed without these homotopies, the restriction maps are usually highly dependent on them. We prove that in our specific case, even though we do not know how to explicitly compute those homotopies, we are still able to compute the restriction maps using an indirect argument.

Section~\ref{s: Def_relSH} provides a brief overview of the original definition of relative symplectic cohomology based on classical Floer theory (Section~\ref{ss: def of relSH}). It then covers the definition of, and equivalence to, a version of relative symplectic cohomology defined using Floer theory in the Morse--Bott with cascades setting.

Section~\ref{s: computations_of_relSH} combines the results from Sections~\ref{s: computations of CF}, \ref{s: algebra}, and \ref{s: Def_relSH} to prove the main results of the paper. Sections~\ref{ss: energy of our cascades} and \ref{ss: direct limits} handle the computation of topological energies and the direct limit of the acceleration data chosen for the ball, presenting it in the form of the Novikov modules discussed in Section~\ref{ss: Novikov modules}. In Section~\ref{ss: proof of main thm}, this allows us to compute the relative symplectic cohomology of balls in $\CP^n$, thereby proving Theorem~\ref{thm: SH of ball in CP^n}. Following this, in Sections~\ref{ss: res maps} and \ref{ss: res maps from CPn to a ball}, we use these ideas together with Proposition~\ref{prop: algebraic preparation for computing restriction maps} from Section~\ref{ss: lim of 1-ray} to compute the restriction maps between balls, as well as from $\CP^n$ to a ball.

Section~\ref{s: Tor and e_st} is the shortest section of the paper. Here, we prove Theorem~\ref{thm: c_relSH vs e_st}, which asserts that the (algebraic) capacity of the relative symplectic cohomology is a lower bound for the stable displacement energy.

The appendix consists of three sections. Appendix~\ref{app: cascades} contains the definitions of Floer and continuation flowlines with cascades, as well as their moduli spaces, including compactifications (Section~\ref{ss: flowlines with cascades and moduli spaces}). We use these to provide upper and lower bounds for the energy of those flowlines (Section~\ref{ss: bounds for energy and cascades}). We then combine this with the compactness results to achieve transversality for data sufficiently close to our ideal data from Section~\ref{ss: acc. data}, while preserving the emptiness of the relevant moduli spaces (see Sections~\ref{ss: transversality} and \ref{ss:tower}).

In Appendix~\ref{app: proof of Seidel's lemma}, we prove Proposition~\ref{prop: asymp_intersection}, which is a slight generalization of Seidel's estimate for the asymptotic intersection number of a Floer cylinder with a divisor (see \cite[Eq. (7.22)]{Seidel_Fukaya_A_infty_Strcs_Assoc_to_Lef_Fib. III}). 

Finally, in Appendix~\ref{app: complement of the ball}, we compute the relative symplectic cohomology of the complement of toric balls in $\CP^n$ and compute several restriction maps. Because the proof is very similar to the case of the ball, we chose to include it in the appendix.

\paragraph*{Acknowledgements.} We would like to express our sincere thanks to Frol Zapolsky for his valuable assistance with the computations of symplectic cohomology and continuation maps for disks in the sphere. We are also grateful to Leonid Polterovich for encouraging us to compute examples of relative symplectic cohomology with coefficients in the Novikov ring, to Umut Varolg\"une\c s for generously sharing his notes on the computation of the symplectic complex associated with a disk in the sphere and for useful discussions, to Alexander Ritter for consulting with us on the equivalence between relative symplectic cohomology defined via Hamiltonian Floer theory and that defined via Morse--Bott theory with cascades, and to Yuan Yao for explanations on the Fredholm analytical details of Floer theory in the Morse-Bott setting. Finally, we thank Octav Cornea, Egor Shelukhin, Dylan Cant, Yoel Groman, Ely Kerman, Sara Tukachinsky and David Keren-Ya'ar for helpful and interesting discussions.

The first author was partially supported by the Milner Foundation and by the Israel Science Foundation (grant no.~1102/20) during his PhD studies at Tel Aviv University. He is currently partially supported by Fondation Courtois during his postdoctoral fellowship at the Centre de Recherches Math\'ematiques, Universit\'e de Montr\'eal. He would also like to thank the Centre de Recherches Math\'ematiques for its warm hospitality and excellent working environment.

\paragraph*{AI Tools Use.} We used LLMs for proofreading, checking mathematical correctness, and simplification of the calculation of the Robin-Salamon index in Section \ref{sss: RS on CPn minus divsors}.

    \section{Morse--Bott approach to Floer theory}\label{s: MB in Floer}

    \subsection{Basic definitions and conventions}\label{ss: convensions}
A \textbf{symplectic manifold} $(M,\omega)$ is a pair consisting of a manifold $M$ and a closed $2$-form $\omega$ that is \textit{non-degenerate}, meaning that for every point $x \in M$, if there exists a vector $X \in T_x M$ such that $\omega_x(X,Y) = 0$ for all $Y \in T_x M$, then necessarily $X = 0$.

An \textbf{almost complex structure} $J$ on $M$ is an automorphism of $TM$ satisfying $J^2=-\Id$. An almost complex structure $J$ is called \textbf{$\omega$-compatible} if $\omega$ is $J$-invariant and the pairing 
$$
\langle\cdot,\cdot\rangle_{\omega,J} = \omega(\cdot, J\cdot)
$$ 
defines a Riemannian metric on $M$.

A \textbf{Hamiltonian} $H$ on $(M,\omega)$ is a smooth function $H\fc S^1\times M  \to \R$, where we identify $S^1$ with $\R/\Z$. If $H$ depends only on $M$, it is called an \textbf{autonomous} or \textbf{time-independent} Hamiltonian; otherwise, it is called a \textbf{time-dependent} Hamiltonian.

The \textbf{symplectic gradient} $X_H$ of $H = (H_t\fc M \to \R)_{t \in S^1}$ is the unique (possibly time-dependent) vector field on $M$ satisfying 
$$
\omega(X_{H_t}, \cdot) = dH_t,
$$
where the differential $d$ is taken with respect to the manifold $M$.
Note that in our conventions we have 
$$X_{H_t}=-J\nabla H_t,$$
where the gradient is computed with respect to the Riemannian metric $\langle\cdot,\cdot\rangle_{\omega,J}$.

The flow of $X_H$ is called the \textbf{Hamiltonian flow} of $H$, denoted by $\varphi_H^t\fc M \to M$, where $t \in \R$ and $\varphi_H^0 = \Id$. The $1$-time map $\varphi_H^1$ is called the \textbf{Hamiltonian diffeomorphism} associated with $H$. Moreover, for every $t\in \R$, the diffeomorphism $\varphi_H^t$ preserves $\omega$, i.e., $\varphi_H^{t*}\omega=\omega$.

A \textbf{$1$-periodic orbit} of $H$ is a smooth map $x\fc S^1 \to M$ satisfying 
$$
x(t) = \varphi_H^t(x(0)) \quad \text{for every } t \in S^1.
$$
Denote by $\cL_0(M) \subset C^\infty(S^1,M)$ the space of contractible smooth loops in $M$. For every $x \in \cL_0(M)$, we can identify $T_x \cL_0(M)$ with the space $\Gamma(S^1,x^*TM)$ of vector fields along $x$. Thus, given an $\omega$-compatible almost complex structure $J$, we can define a Riemannian metric $g_J$ on $\cL_0(M)$ by
$$
g_{J,x}(\eta,\theta) = \int_{S^1} \langle \eta(t),\theta(t)\rangle_J \,dt,
$$
for every $x\in \cL_0(M)$ and $\eta,\theta\in T_x \cL_0(M)=\Gamma(S^1,x^*TM)$.

If $u\fc \R\times S^1\to M$ is a smooth cylinder with asymptotics $x_-,x_+\fc S^1\to M$ at $-\infty,+\infty$, respectively, where $x_-,x_+$ are $1$-periodic orbits of $H$, then its \textbf{topological energy} is defined to be 
\begin{equation}\label{eq: top energy}
  E_{top}(u)=\int_{\R\times S^1} u^*\omega+\int_{S^1} H\circ x_+(t)\,dt-\int_{S^1}H\circ x_-(t)\,dt.  
\end{equation}
In particular, if $v\fc S^2\to M$ is a smooth sphere, then
$$E_{top}(u\#v)=E_{top}(u)+\omega([v]),$$
where $\omega([v])=\int_{S^2} v^*\omega$, and $u\#v$ is the connected sum of $u$ and $v$.

A \textbf{capping} $w$ for a contractible loop $x$ is a smooth map $w\fc D^2 \to M$ with $w|_{\partial D^2} = x$, where $D^2 \subset \CP$ is the closed unit disk with standard orientation and $\partial D^2$ is identified with $S^1$. Two cappings $w,w'$ for a contractible loop $x$ are equivalent if $[w \# \overline{w}'] = 0 \in \pi_2(M)$, where $w \# \overline{w}'\fc S^2 \to M$ is the result of gluing $w$ with $\overline{w}'$, which is equal to $w'$ but with $D^2$ given the opposite orientation.

We denote by $\widetilde{\cL_0}(M)$ the cover of $\cL_0(M)$ which is defined as the collection of pairs $(x,\hat{x})$, where $\hat{x}$ is an equivalence class of cappings of $x$.

Given a Hamiltonian $H\fc M \times S^1 \to \R$, we define its corresponding \textbf{action functional} $\cA_H\fc \widetilde{\cL_0}(M) \to \R$ by 
$$
\cA_H(x,\hat{x}) = \int_{S^1} H(x(t),t)dt + \int_{D^2} w^*\omega,
$$
for every $(x,\hat{x}) \in \widetilde{\cL_0}(M)$, where $w$ is a representative of $\hat{x}$.

Given an $\omega$-compatible almost complex structure $J$ on $M$, a direct computation shows that the gradient of $\cA_H$ with respect to the metric $g_J$ has the form
$$
\nabla \cA_H(x,\hat{x})=-J(\dot{x}-X_H(x)),
$$
for every $(x,\hat{x})\in \widetilde{\cL_0}(M)$. Thus, we deduce that $(x,\hat{x})$ is a critical point of $\cA_H$ if and only if $x$ is a $1$-periodic orbit of $H$. Moreover, the positive gradient flow of $\cA_H$ gives rise to the \textbf{Floer equation} 
\begin{equation}\label{eq: Floer eq}
    \frac{\partial u}{\partial s}(s,t)+J(u(s,t))\left(\frac{\partial u}{\partial t}(s,t)-X_H(u(s,t))\right)=0,
\end{equation}
where $u\fc \R\times S^1\to M$ is a smooth cylinder in $M$. Note that the Floer equation is well-defined over $\cL_0(M)$ and not just over $\widetilde{\cL_0}(M)$.

\subsubsection{Robbin--Salamon index}\label{sss: RS part I}

The Robbin--Salamon index first appeared in \cite{RS_index}. A comprehensive background and collection of useful results concerning it can be found in \cite{Gutt_indexes}. Here we present the essential material required for our purposes. 

\textbf{Reminder:} Let $n\in \N$. The standard complex structure $J_0$ on $\R^{2n}$ is  presented in the standard basis by the block matrix
$$J_0=\left(\begin{array}{cc}
    0 & -I_n \\
    I_n & 0
\end{array}\right),$$
where $I_n\in \Mat(n,\R)$ is the unit matrix for size $n\times n$. Also, the symplectic group, which consists of symplectic matrices, is defined by
$$\Sp(2n)=\left\{ A\in \Mat(2n,\R)\,\middle|\,A^T J_0 A=J_0 \right\}.$$
Also, recall that the \textbf{signature}, $\sign S$ ,of a symmetric matrix $S$ is defined to be the difference between the number of its positive and negative eigenvalues.

Now, let us begin with an axiomatic characterization of the Robbin--Salamon index. 

\begin{thm}[{\cite[Theorem 55]{Gutt_indexes}}]\label{thm: RS_index_characterization}
    There exists a unique function $\mu$ that associates a number from $\frac{1}{2}\Z$ to each continuous path of symplectic matrices that is characterized by the following properties:

\begin{itemize}
  \item \textbf{(Homotopy)} $\mu$ is invariant under homotopies of paths with fixed endpoints.
  
  \item \textbf{(Catenation)} $\mu$ is additive under concatenation of paths: if $\Psi$ is the concatenation of $\Psi_1$ and $\Psi_2$, then $\mu(\Psi) = \mu(\Psi_1) + \mu(\Psi_2)$.
  
  \item \textbf{(Zero)} $\mu$ vanishes on any path $\Psi \fc [a,b] \to \Sp(2n)$ such that $\dim \ker(\Psi(t) - \Id) = k$ is constant for all $t \in [a,b]$.
  
  \item \textbf{(Normalization)} If $S = S^T \in \R^{2n \times 2n}$ is a symmetric matrix with all eigenvalues of absolute value less than $2\pi$, and $\Psi(t) = \exp(J_0 S t)$ for $t \in [0,1]$, then
  $$
  \mu(\Psi) = \frac{1}{2} \sign\, S.
  $$
\end{itemize}

\end{thm}
\begin{defin}\label{def: RS index}
    The function $\mu$ from the previous theorem is called the \textbf{Robbin--Salamon index} and denoted here by $\mu_{RS}$.
\end{defin}

Given a Hamiltonian $H$ on a symplectic manifold $(M,\omega)$, let $x\fc[0,1]\to M$ be a $1$-periodic orbit and fix a symplectic trivialization $\tau=\{\tau_t\}_{t\in[0,1]}$ of $TM$ along $x$, where for every $t\in[0,1]$ the linear map $\tau_t\fc(R^{2n},\omega_{0})\to (T_{x(t)}M,\omega)$ is a symplectic isomorphism and $\tau(0)=\tau(1)$. The linearized flow $d_{x(0)}\varphi_H^t\fc T_{x(0)}M\to T_{x(t)}M$ is converted by $\tau$ into a path of symplectic matrices $\Psi\fc[0,1]\to \Sp(2n)$ given by
$$\Psi(t)=\tau_t^{-1}\circ d_{x(0)}\varphi_H^t\circ\tau_0$$
for every $t\in[0,1]$. The \textbf{Robbin--Salamon index of the $1$-periodic orbit $x$ with respect to the trivialization $\tau$} is then defined by $\mu_{RS}^\tau(x;H):=\mu_{RS}(\Psi)$.

\subsubsection{Classical Floer theory}

Throughout the rest of Section~\ref{s: MB in Floer}, let $(M,\omega)$ be a fixed closed monotone symplectic manifold.

Let $H\in C^\infty(S^1\times M)$ be a Hamiltonian. A $1$-periodic $x\fc S^1\to M$ of $H$ is called \textbf{non-degenerate} if $1$ is not an eigenvalue of the linearization of the return map, i.e. $d_{x(0)}\varphi_H^1\fc T_{x(0)}M\to  T_{x(0)}M$. The Hamiltonian $H$ is called \textbf{non-degenerate} if all of its $1$-periodic orbits are non-degenerate. Since $M$ is closed, any non-degenerate Hamiltonian has a finite collection of $1$-periodic orbits.

 Assume that $H$ is a non-degenerate Hamiltonian, let $J$ be a time-dependent compatible almost complex structure on $(M,\omega)$ and assume that $(H,J)$ are generic in the sense of \cite{FHS_1995_Transversality}. The Floer complex of $H$, generated by the set $\cP^\circ(H)$ of its contractible $1$-periodic orbits, is defined as 
$$CF^*(H;\Z) = \bigoplus_{x\in\cP^\circ(H)}\Z\cdot x.$$ 
This complex is graded over $\Z/2N_M^{}\Z$ by the Robbin--Salamon index, shifted by $n$, where $N_M$ is the minimal first Chern number. The differential $d\fc CF(H;\Z)\to CF(H;\Z)$ is given by
$$dx = \sum_{y\in\cP^\circ(H)}\#\widehat\cM(H;x,y)\,\cdot\,y\,,$$
where $\widehat\cM(H;x,y)$ stands for the moduli space of unparametrized Floer trajectories corresponding to $H$, going from $x$ to $y$, such that every solution has index $1$; $\#\widehat\cM(H;x,y)$ is a signed count. The fact that this sum is well-defined is due to the assumption that $(M,\omega)$ is monotone and closed, thus there is only  a finite number of unparametrized Floer solutions of index $1$.

\subsection{Morse--Bott setting}\label{ss: MB setting}
Let $H$ be an autonomous Hamiltonian on $(M,\omega)$. By the existence and uniqueness theorem for ODEs, the Hamiltonian flow forms a one-parameter subgroup of diffeomorphisms and preserves the level sets of $H$. Moreover, if $x\fc S^1 \to M$ is a $1$-periodic orbit of an autonomous Hamiltonian $H$, then for every $t_0 \in S^1$, the loop $x_{t_0}\fc S^1 \to M$ given by 
$$
x_{t_0}(t) = x(t - t_0), \quad \text{for every } t \in S^1,
$$
is also a $1$-periodic orbit of $H$. Thus, every $1$-periodic orbit of an autonomous Hamiltonian gives rise to an $S^1$-family of $1$-periodic orbits. Such $S^1$-families are examples of critical submanifolds of the action functional. In this section we introduce a Morse--Bott theoretical approach to Floer theory that helps to overcome the fact that the critical points of the action functional are not isolated.

\subsubsection{Morse--Bott functions}

Let $(X,g)$ be a Riemannian manifold, possibly infinite-dimensional. Let $\cA\fc X\to \R$ be a $C^2$-smooth functional. Then the gradient $\nabla \cA\in \Gamma(X, TX)$ and the Hessian $\nabla^2 \cA\in \Gamma(X, TX\otimes T^*X)$ of $\cA$ are well-defined.

\begin{defin}
        A $C^2$ functional $\cA\fc X\to \R$ is \textbf{Morse--Bott} if it has the following two properties:
    \begin{enumerate}
        \item The set of critical points consists of a disjoint union of smooth submanifolds:
        $$\Crit(\cA)=\bigcup_{i} C_i.$$
        \item For all $i$ and $w\in C_i$,
        $$T_w C_i=\ker\left(\nabla^2_w \cA\right)$$
    \end{enumerate}
\end{defin}
\begin{rem}
    If $\cA\fc X\to \R$ is a Morse--Bott function, then the connected components of $\Crit(\cA)$ are called \textbf{critical submanifolds} with respect to $\cA$.
\end{rem}
In the next subsection we discuss this in the context of the action functional of an autonomous Hamiltonian.

\subsubsection{Morse--Bott condition}

Let $H\fc S^1\times M\to \R$ be a Hamiltonian and let $\cA_H\fc \widetilde{\cL_0}(M) \to \R$ be the corresponding action functional. 
\begin{defin}
    We say that $H$ \textbf{satisfies the MB condition} if the action functional $\cA_H\fc \widetilde{\cL_0}(M) \to \R$  is Morse--Bott.
\end{defin}

\begin{rem}
    Note that if $H$ is non-degenerate then it satisfies the \textbf{MB} condition.
\end{rem}

A connected component $C$ of $\Fix(\varphi_H^1)\subseteq M$ satisfies the \textbf{MB} condition if it is a closed manifold, and satisfies $T_p C=\ker(d_p \phi_H^1-\id)$, for every $p\in C$.

The following is a version of \cite[Theorem 23]{Fauck_2016_thesis}:
\begin{thm}
     If every connected component of $\Fix(\varphi_H^1)\subset M$ satisfies the \textbf{MB} condition, then $H$ satisfies the \textbf{MB} condition.
\end{thm}

Denote by $p\fc\widetilde{\cL_0}(M)\to \cL_0(M)$ the covering map $(x,\hat{x})\mapsto x$, and let  $\pi\fc \cL_0(M)\to M$ be the evaluation map $\pi(x)=x(0)$, where we identify $S^1$ with $\R\slash \Z$. Given a connected critical submanifold $C\subset \widetilde{\cL_0}(M)$ for $\cA_H$, the restriction $p|_C\fc C\to \cL_0(M)$ is a smooth embedding. We refer to its image in  $\cL_0(M)$ as a \textbf{critical submanifold with respect to $H$}. We denote by $\bS_H$ the collection of the critical submanifolds with respect to $H$, that is $\bS_H=p(\Crit \cA_H)$.

\begin{rem} Note that the map $\pi\fc \cL_0(M)\to M$ induces a bijection between the connected components of $\Fix(\varphi_H^1)\subset M$ and $\bS_H=p(\Crit(\cA_H))\subset \cL_0(M)$. Moreover, for every connected component $S\subset \Crit(\cA_H)$, the map $\pi\fc S\to M$ is a smooth embedding.

\end{rem}

\subsubsection{Trivializations and indices}

Let $H\fc S^1\times M\to \R$ be a Hamiltonian, and assume that $H$ satisfies the \textbf{MB} condition. Let $S\subset \cL_0(M)$ be a critical submanifold with respect to $H$.

\begin{defin}
    Let $\tau=(\tau_x)_{x\in S}$ be a family of symplectic trivializations of $TM$ along the elements of $S$, i.e. for every $x\in S$, $\tau_x$ is a symplectic trivialization of $x^*TM$. We say that $\tau$ is \textbf{a compatible symplectic trivialization along $S$} if for every $x,y\in S$ there is a path $\gamma\fc[0,1]\to S$ such that the symplectic trivialization $\gamma^*\tau_y$ is homotopic to the trivialization $\tau_x$.\footnote{The pullback $\gamma^*\tau_y$ is defined as follows:
    the symplectic trivialization $\tau_y\fc S^1\times \R^{2n}\to y^*TM$ extends to a symplectic trivialization $[0,1]\times S^1\times \R^{2n}\to \gamma^*TM$ since there is a deformation retraction of $[0,1]\times S^1$ onto $\{1\}\times S^1$. The restriction of this symplectic trivialization to $\{0\}\times S^1$, i.e. the symplectic trivialization $S^1\times \R^{2n}\to x^* TM$ is defined to be $\gamma^*\tau_y$. Note that, up to a homotopy, $\gamma^*\tau_y$ is independent of the extension of $\tau_y$ to $[0,1]\times S^1$. This is because a symplectic trivialization is equivalent to a section of the corresponding principal $\Sp(2n,\R)$ bundle.}
    In this case, we say that $\tau$ is \textbf{obtained from cappings} if there exists $x\in S$ and a capping $\hat{x}$ of $x$ such that $\tau_x$ is homotopic to the trivialization of $x^*TM$ that is induced from $\hat{x}$.
\end{defin}
\begin{rem}\label{rem: compatible trivializations}\phantom{M}
    \begin{enumerate}
    \item If $\tau$ is a compatible symplectic trivialization, then for every $x \in S$ we denote the Robbin--Salamon index of $x$ with respect to the trivialization $\tau_x$ by $\mu_{RS}^\tau(x;H)$.
    \item In this case, for every $x,y \in S$ we have $\mu_{RS}^\tau(x;H) = \mu_{RS}^\tau(y;H)$. Indeed, there is a path $\gamma$ on $S$ connecting $x$ to $y$ such that $\tau_x$ is homotopic to $\gamma^*\tau_y$. Since the Robbin--Salamon index is invariant under homotopies, we deduce that $$\mu_{RS}^{\tau_x}(x;H) = \mu_{RS}^{\gamma^*\tau_y}(x;H) = \mu_{RS}^{\tau_y}(y;H).$$
    \item Given a point $x_0\in S$ and a symplectic trivialization $\tau_0$ of $x_0^*TM$, we can extend it to a compatible symplectic trivialization $\tau$ along $S$ as follows: Choose a family $(\gamma_x)_{x\in S}$ of smooth paths such that for every $x\in S$, $\gamma_x\fc [0,1]\to S$ starts at $x$ and ends at $x_0$, and for  $\gamma_{x_0}$ choose the constant path at $x_0$. Then, define $\tau=(\tau_x)_{x\in S}$ by $\tau_x=\gamma_x^*\tau_0$ for every $x\in S$.
    \item If $\tau$ is obtained from cappings, then for every $y \in S$ there exists a capping $\hat{y}$ such that $\tau_y$ is homotopic to the trivialization of $y^*TM$ induced by $\hat{y}$. This follows immediately from taking the connected sum of the capping $\hat{x}$, which induces $\tau_x$, with a path $\gamma$ for which $\tau_x$ is homotopic to $\gamma^*\tau_y$. 
    \end{enumerate}
\end{rem}

The indices that will be useful in this paper are what we call the Floer--Morse--Bott indices. They are defined as follows:
\begin{defin}\label{def: FMB index}
    Let $\tau=(\tau_x)_{x\in S}$ be a family of symplectic trivializations of $TM$ that is compatible along $S$. Let $h\fc S\to \R$ be a Morse functions and let $p\in S$ be a critical point of $h$. We define \textbf{Floer--Morse--Bott index} of a critical point $p$ of $h$ by
    \[\mu_{FMB}^\tau(p;H,h) = \mu_{RS}^\tau(p;H) + \frac{1}{2}\dim M - \frac{1}{2}\dim S + \ind_h p.\]
Here $\ind_h p$ is the Morse index of $p$ as a critical point of $h$.
\end{defin}
\begin{rem}\label{rem: mu_FMB well def mod 2N_M}
$\,$
\begin{itemize}
    \item When $H$ and $h$ are clear from the context, we abbreviate to $\mu_{FMB}^\tau(p)$.
    \item The Floer--Morse-Bott indices of the same critical point of $h$, with respect to two compatible symplectic trivializations obtained from cappings, differ from each other by a multiple of twice the minimal first Chern number $N_M$ of $(M,\omega)$. For this reason, the Floer--Morse-Bott index is well-defined for all critical points of $h$ as an element of $\Z/2N_M \Z$.
\end{itemize}
   
\end{rem}

    \subsection{Flowlines with cascades}\label{ss: flowlines and cascades}

    In this section we survey basic concepts and sketch the definitions of Floer and continuation flowlines with cascades. For exact definitions see Appendix~\ref{app: cascades}. 

\begin{itemize}
    \item \textbf{Floer flowlines with cascades:} 
Let $H$ be a Hamiltonian satisfying the \textbf{MB} condition, and let $J$ be an almost complex structure compatible with $\omega$, possibly time-dependent. Let $h$ be a Morse function on $\bS_{H}$, which is the union of the critical submanifolds with respect to $H$. Equip $\bS_H$ with a Riemannian metric $g$. Let $q^+,q^-$ be critical points of $h$.

Given $m\in \Z_{\geq0}$, a \textbf{Floer flowline of $(H,J)$ with $m$ cascades connecting $q^-$ to $q^+$} consists of the following information:
\begin{itemize}
    \item $2m+2$ points $q_0^-,q_0^+,\ldots,q_m^-,q_m^+\in \bS_{H} $ such that for every $0\leq i\leq m$, the points $q_i^-,q_i^+$ belong to the same connected component of $\bS_{H}$, where $q_0^-=q^-,\,q_m^+=q^+$. 
    \item For every $1\leq i\leq m-1$, the points $q_i^-$ and $q_i^+$ are connected to each other by a finite gradient flowline of $h$, whereas $q_0^-$ and $q_m^-$ are connected to $q_0^+$ and $q_m^+$, respectively, by a semi-infinite gradient flowline of $h$.

    \item $m$ Floer cascades, i.e. $m$ Floer cylinders $u_1,\ldots,u_m$ such that the cylinder $u_i$ has a negative asymptote at $q_{i-1}^+$ and a positive asymptote at $q_i^-$,   for every $1\leq i\leq m$.
\end{itemize}
The moduli space of the Floer flowlines of $(H,J)$ with $m$ cascades connecting $q^-$ to $q^+$ is denoted by $\cM_m(q^-,q^+;H,J)$, this space admits a free $\R^m$-action by translations of the Floer cylinders constituting the cascades, the quotient is denoted by $\widehat\cM_m(q^-,q^+;H,J)$. When $m=0$ we quotient out the translation of the single gradient flowline. Also, we denote $$\widehat\cM(q^-,q^+;H,J)=\bigcup_{m=0}^\infty \widehat\cM_m(q^-,q^+;H,J), $$
this union can be given the structure of a smooth manifold, using gluing neighborhoods, for families of flowlines with $m$ cascades, degenerating into a flowline with more cascades. See Definitions \ref{def: moduli of unparametrized Floer flowlines of m cascades} and \ref{def: moduli of unparametrized Floer flowlines} for details.
Given a family of symplectic trivializations $\tau$ of $TM$ which is compatible along $\bS_H$, the local dimension of $\widehat\cM(q^-,q^+;H,J)$ at a Floer flowline $u$ with $m$ cascades $u_1,\ldots,u_m$ connecting $q^-$ to $q^+$ is given by
\begin{equation}\label{eq: index formula for Floer}
     \dim_u \widehat\cM(q^-,q^+;H,J)= \mu^\tau_{FMB}(q^+) - \mu^\tau_{FMB}(q^-) + 2\sum_{k=1}^m c_1^\tau(u_k)-1.
\end{equation}
    See Section~\ref{ss: c_1 and intersection} for the definition of the relative first Chern class $c^\tau_1$.

    \item Let $H_-$ and $H_+$ be Hamiltonians satisfying the \textbf{MB} condition, and let $J_-$ and $J_+$ be almost complex structures compatible with $\omega$, possibly time-dependent. Moreover, let $\cH=(H_s)_{s\in \R}$ be a homotopy of Hamiltonians and $\cJ=(J_s)_{s\in \R}$ a homotopy of compatible almost complex structures such that $H_s = H_-$ and $J_s = J_-$ for every $s\leq 0$, and $H_s = H_+$ and $J_s = J_+$ for $s\geq 1$. Let $h_-$ and $h_+$ be two Morse functions on $\bS_{H_-}$ and $\bS_{H_+}$, respectively. Equip $\bS_{H_-},\bS_{H_+}$ with Riemannian metrics $g_-,g_+$, respectively.
     Let $q^+,q^-$ be critical points of $h_+,h_-$, respectively.

Given $m\in \N$, a \textbf{continuation flowline of $(\cH,\cJ)$ with $m$ cascades connecting $q^-$ to $q^+$} consists of the following information:
\begin{itemize}
    \item An integer $0\leq k\leq m-1$.
    \item $2k+2$ points $q_0^-,q_0^+,\ldots,q_{k}^-,q_{k}^+\in \bS_{H_-}$ such that for every $0\leq i\leq k$ the points $q_i^-$ and $q_i^+$ belong to the same connected component of $\bS_{H_-}$, where $q_0^-=q^-$.

    \item For every $1\leq i\leq k$, the points $q_i^-$ and $q_i^+$ are connected to each other by a finite gradient flowline of $h_-$, whereas $q_0^-$ is connected to $q_0^+$ by a semi-infinite gradient flowline of $h_-$.

    \item $k$ Floer cascades $u_1,\ldots,u_k$ with respect to $(H_-,J_-)$, such that the cylinder $u_i$ has a negative asymptote at $q_{i-1}^+$ and a positive asymptote at $q_i^-$,   for every $1\leq i\leq k$.
    
    \item A continuation cascade, i.e. a continuation cylinder with respect to the pair of homotopies $(\cH,\cJ)$ from $(H_-,J_-)$ to $(H_+,J_+)$ connecting $q_k^+$ to $q_{k+1}^-$.
    
    \item $2(m-k)$ points $q_{k+1}^-,q_{k+1}^+,\ldots,q_m^-,q_m^+\in \bS_{H+}$ such that for every $k+1\leq i\leq m$ the points $q_i^-$ and $q_i^+$ belong to the same connected component of $\bS_{H_+}$, where $q_m^+=q^+$.
    \item For every $m-k+1\leq i\leq m-1$, the points $q_i^-$ and $q_i^+$ are connected to each other by a finite gradient flowline of $h_+$, whereas $q_m^-$ is connected to $q_m^+$ by a semi-infinite gradient flowline of $h_+$.

    \item $m-k-1$ Floer cascades $u_{k+1},\ldots,u_m$ with respect to $(H_+,J_+)$, such that the cylinder $u_i$ has a negative asymptote at $q_{i-1}^+$ and a positive asymptote at $q_i^-$, for every $k+1\leq i\leq m$.

\end{itemize}
The moduli space of the continuation flowlines of $(\cH,\cJ)$ with $m$ cascades connecting $q^-$ to $q^+$ is denoted by $\cM_m^{cont}(q^-,q^+;H,J)$; this space admits a free $\R^{m-1}$-action by translations of the Floer cylinders constituting the cascades, the quotient is denoted by $\widehat\cM_m^{cont}(q^-,q^+;H,J)$. Also, we denote $$\widehat\cM^{cont}(q^-,q^+;H,J)=\bigcup_{m=1}^\infty \widehat\cM_m^{cont}(q^-,q^+;H,J), $$
this union can be given the structure of a smooth manifold, using gluing neighborhoods, for families of flowlines with $m$ cascades, degenerating into a flowline with more cascades. See Definitions \ref{def:unparamterizedContWithCascades} and \ref{def:unparamterizedContWithCascades-general} for details.
Given a family of symplectic trivializations $\tau$ of $TM$ which is compatible along $\bS_{H_-}$ and $\bS_{H_+}$, the local dimension of $\widehat\cM^{cont}(q^-,q^+;\cH,\cJ)$ at a continuation flowline $u$ with $m$ cascades $u_1,\ldots,u_m$ connecting $q^-$ to $q^+$ is given by
\begin{equation}\label{eq: index formula for continuation}
    \dim_u \widehat\cM^{cont}(q^-,q^+;\cH,\cJ)=\mu^\tau_{FMB}(q^+;H_+) - \mu^\tau_{FMB}(q^-;H_-) + 2\sum_{k=1}^m c_1^\tau(u_k).
\end{equation}

    \end{itemize}
\begin{rem}\label{rem: index formula for Floer and cont flowlines}Given a regular pair $(H,J)$ where $H$ is non-degenerate, the Fredholm index of any Floer trajectory $u$ with asymptotes $q^-,q^+$ equals the local dimension of $\cM(q^-,q^+;H,J)$. Thus, Equation~\eqref{eq: index formula for Floer} implies the equality
$$\ind u=\mu_{RS}^\tau(q^+;H)-\mu_{RS}^\tau(q^-;H)+2c_1^\tau(u).$$
This equation is traditionally referred to as the \textbf{index formula} for $u$. For this reason, throughout this paper, we call Equation~\eqref{eq: index formula for Floer} the index formula for $u$ when $u$ is a Floer flowline with cascades, and, similarly, Equation~\eqref{eq: index formula for continuation} the index formula for $u$ when $u$ is a continuation flowline with cascades. Therefore, we call the expression 
\[\mu^\tau_{FMB}(q^+) - \mu^\tau_{FMB}(q^-) + 2\sum_{k=1}^m c_1^\tau(u_k)\] the \textbf{index} of $u$.
See also Remark~\ref{rem: index formula for homotopies} below for a more general formulation of the index formula. Also, it is important to note that the left-hand side of the index formula is independent of the compatible symplectic trivialization $\tau$.
\end{rem}

    \subsection{Differential and continuation maps}\label{ss:CF_MB_over_Z}

    In this section we discuss Floer differential and continuation maps by the Morse--Bott approach to Floer theory.

Let $H$ be a Hamiltonian satisfying the \textbf{MB} condition, and let $J$ be an almost complex structure compatible with $\omega$, possibly time-dependent. Let $h$ be a Morse function on $\bS_{H}$, which is the union of the critical submanifolds with respect to $H$. Equip $\bS_H$ with a Riemannian metric $g$. Let us define the associated Floer complex $CF^*(H;\Z)$, assuming regularity for all the moduli spaces involved.
        
        The associated Floer complex $CF^*(H;\Z)$ is freely generated by the set $\Crit(h^{}_H)$ over $\Z$ namely
        \[CF^*(H;\Z) = \bigoplus_{x\in\Crit(h^{}_H)}\Z\cdot x.\]
        This complex is graded over $\Z/2N_M^{}\Z$, where the degree of a generator $x$ is given by the Floer--Morse--Bott index $\mu_{FMB}(x)$, see Definition \ref{def: FMB index}, and $N_M$ is the minimal first Chern number of $M$. 

        The differential $d\colon CF(H;\Z) \to CF(H;\Z)$ is given by

        \[dx = \sum_{y\in\Crit(h^{}_H)}\#\widehat\cM^0(x,y;H,J)\,\cdot\,y\,,\]
       where $\widehat{\cM}^0(x,y;H,J)$ is the union of the zero dimensional components of $\widehat{\cM}(x,y;H,J)$ as defined in Section~\ref{ss: flowlines and cascades}, and $\#\widehat{\cM}^0(x,y;H,J)$ denotes a signed count according to a system of coherent orientations chosen on the moduli spaces, see \cite[Section 4.4]{BO_Morse_Bott_2009} and \cite[Section 8]{Schmaschke_2016_thesis}. The reason that this sum is well-defined is due to the assumption that $(M,\omega)$ is monotone and closed, which, as shown in Section \ref{ss: bounds for energy and cascades}, implies a bound on the number of cascades and the energy of flowlines of a given index. Thus, the zero dimensional moduli spaces are compact and hence, finite.

       Now, let us discuss continuation maps.  Let $H_-$ and $H_+$ be Hamiltonians satisfying the \textbf{MB} condition, and let $J_-$ and $J_+$ be almost complex structures compatible with $\omega$, possibly time-dependent. Moreover, let $\cH=(H_s)_{s\in \R}$ be a homotopy of Hamiltonians and $\cJ=(J_s)_{s\in \R}$ a homotopy of compatible almost complex structures such that $H_s = H_-$ and $J_s = J_-$ for every $s\leq 0$, and $H_s = H_+$ and $J_s = J_+$ for $s\geq 1$. Let $h_-$ and $h_+$ be two Morse functions on $\bS_{H_-}$ and $\bS_{H_+}$, respectively. Equip $\bS_{H_-},\bS_{H_+}$ with Riemannian metrics $g_-,g_+$, respectively. Let us define the associated continuation map $\Psi\fc CF^*(H_-;\Z)\to CF^*(H_+;\Z)$, again, assuming regularity for all the moduli spaces involved.

         The continuation map $\Psi\fc CF(H_-;\Z) \to CF(H_+;\Z)$ is the unique homomorphism of $\Z$-modules that satisfies 
        
        \[\Psi(x) = \sum_{y\in\Crit(h_+)}\#\widehat\cM^{cont,0}(x,y;\cH,\cJ)\,\cdot\,y\,,\]
        for every generator $x\in \Crit(h_-)$ of the Floer complex $CF^*(H_-;\Z)$. Here $\widehat\cM^{cont,0}(x,y;\cH,\cJ)$ stands for the union of the zero dimensional components of $\widehat\cM^{cont}(x,y;\cH,\cJ)$ as defined in Section~\ref{ss: flowlines and cascades}, and $\#\widehat\cM^{cont,0}(x,y;\cH,\cJ)$ denotes a signed count according to a system of coherent orientations chosen on the moduli spaces. The sum is finite for reasons similar to the differential.

    \subsection{Energy of cascades}\label{ss: energy}

The purpose of this section is to express the topological energy of a cascade in terms of action and index with respect to a given symplectic trivialization. Recall that $(M,\omega)$ is assumed to be a monotone closed symplectic manifold, denote by $\kappa>0$ its monotonicity constant, that is, we have $[\omega]=\kappa c_1$ as cohomology classes in $H_{DR}^2(M)$.

Let $H_-$ and $H_+$ be Hamiltonians on $(M,\omega)$
satisfying the \textbf{MB} condition, such that each Hamiltonian is constant along each of its connected critical submanifolds, and let $J_-$ and $J_+$ be almost complex structures compatible with $\omega$, possibly time-dependent. Moreover, let $\cH=(H_s)_{s\in \R}$ be a homotopy of Hamiltonians and $\cJ=(J_s)_{s\in \R}$ a homotopy of compatible almost complex structures such that $H_s = H_-$ and $J_s = J_-$ for every $s\leq 0$, and $H_s = H_+$ and $J_s = J_+$ for $s\geq 1$. Let $h_-$ and $h_+$ be two Morse functions on $\bS_{H_-}$ and $\bS_{H_+}$, respectively. Equip $\bS_{H_-},\bS_{H_+}$ with Riemannian metrics $g_-,g_+$, respectively. Let $\tau$ be a family of symplectic trivializations of $TM$ which is compatible along $\bS_{H_-}$ and $\bS_{H_+}$.

\begin{rem}
     Note that in the case when $H_-=H_+$, $J_-=J_+$ and the homotopies $\cH,\cJ$ are constant then any continuation flowline with cascade is actually a Floer flowline with cascades.
\end{rem}
The topological energy of a smooth cylinder is defined in Equation~\eqref{eq: top energy}. In order to adapt it for flowlines with cascades, we need the following claim:

\begin{claim}\label{claim: area of level set cylinders vanish}
    Let $S\subset \cL_0(M)$ be a critical submanifold with respect to some autonomous Hamiltonian $H$. And, assume that $H$ is constant along its critical submanifolds. Then for every smooth path $\gamma\fc I\to S$, for some interval $I\subset \R$ we have $\int_{S^1\times I}\gamma^*\omega=0$. 
\end{claim}
\begin{proof}
    Write $I=[a,b]$. Note that since every loop $\gamma(s,\cdot)$ is an orbit, then $\frac {\partial}{\partial t} \gamma$ is the symplectic gradient of $H$. Let us now compute the integral
    \[
    \int_{S^1\times I}\gamma^*\omega =  \int_0^1 \int_a^b \omega\left(\frac {\partial}{\partial t}\gamma,\frac {\partial}{\partial s}\gamma\right) dtds = \int_0^1 \int_a^b dH\left(\frac {\partial}{\partial s}\gamma\right) dtds = 0, 
    \]
    but $dH(\partial_s\gamma(s,t))=\partial_s(H\circ\gamma)(s,t)=0$ since $H$ is constant along critical submanifolds and $\gamma(s,t)$ lies on the same critical submanifold for every $(s,t)\in \R\times S^1$.

\end{proof}
As a consequence, it is natural to define the topological energy of a flowline with cascades based just on its ends and on its cascades:
Let $q^-$ and $q^+$ be critical points of $h_-$ and $h_+$, respectively, and let $u$ be a continuation flowline of $(\cH,\cJ)$ connecting $q^-$ to $q^+$. 
The \textbf{topological energy} of $u$ is defined to be
     
\begin{equation}\label{eq: top energy of flowline}
  E_{top}(u)=\int_{S^1} H_+\circ x_+(t)\,dt-\int_{S^1}H_-\circ x_-(t)\,dt+\sum_{j=1}^m \omega([u_j]),  
\end{equation} 
where $\omega([u_j])=\int_{\R\times S^1} u_j^*\omega$.

\begin{rem}
   Let $x, y \fc S^1 \to M$ be loops, and let $\hat{x}, \hat{y}$ be cappings for $x$ and $y$, respectively. Let $\tau$ be a trivialization of $x^*TM$ and $y^* TM$ induced by these cappings. If $u \fc \R \times S^1 \to M$ is a cylinder asymptotic to $x$ at $-\infty$ and to $y$ at $+\infty$, then
  $$c_1^\tau(u) = c_1(\hat{x} \# u \# \hat{y}^{-1}),$$
   where $\hat{y}^{-1}$ denotes the capping $\hat{y}$ with its orientation reversed. See Claim~\ref{claim: c_1 v.s. rel c_1}.
   \end{rem}

The main result of this section is the following:
\begin{prop}
Let $\tau$ be a compatible symplectic trivialization along each connected component of $\bS_{H_-}$ and $\bS_{H_+}$, obtained from cappings. Then the topological energy of a continuation flowline with cascades $u$ with asymptotes $q^-\in \bS_{H_-}, q^+\in\bS_{H_+}$ is given by
\begin{equation}\label{eq: energy}
    E_{top}(u) = \cA_{H_+}(q^+, \hat{q}^+) - \cA_{H_-}(q^-, \hat{q}^-) + \frac{1}{2}\cdot\kappa \left( \dim_u \widehat\cM(q^-, q^+) - \mu_{FMB}^\tau(q^+) + \mu_{FMB}^\tau(q^-)  \right),
\end{equation}
where $\hat{q}^-$ and $\hat{q}^+$ are cappings for $q^-$ and $q^+$, respectively, that induce trivializations homotopic to $\tau_{q^-}$ and $\tau_{q^+}$, respectively.   
\end{prop}
\begin{proof} Let us start by writing $u$ as $\gamma_0\#u_1\#\gamma_1\#\cdots\#\gamma_{m-1}\#u_m\#\gamma_m$, where $u_1,\ldots,u_m$ are the cascades of $u$ with endpoints $q^-=q_0^-,q_0^+,\ldots,q_m^-,q_m^+=q^+$, as described in Section~\ref{ss: flowlines and cascades}, and each $\gamma_k$ is a smooth path on $\bS_{H_-}\cup \bS_{H_+}$ connecting $q_k^-$ to $q_k^+$ such that $\tau_{q_k^-}$ is homotopic to $\gamma_k^*\tau_{q_k^+}$.

    Now, let $v$ be a capping for $\hat{q}^+$ which is equivalent to $\hat{q}^-\#u$, then from Stokes' theorem we get that 
    $$E_{top}(u) = \cA_{H_+}(q^+,v) - \cA_{H_-}(q^-,\hat{q}^-).$$
    Denote by $(\hat{q}^+)^{-1}$ the capping of $q^+$ obtained from $\hat{q}^+$ by reversing the orientation. 
    Thus we get that
        \begin{align*}
            E_{top}(u)&=\cA_{H_+}(q^+,\hat{q}^-\#u)-\cA_{H_-}(q^-,\hat{q}^-)\\
        &=\cA_{H_+}(q^+,\hat{q}^-\#u\# (\hat{q}^+)^{-1}\#\hat{q}^+)-\cA_{H_-}(q^-,\hat{q}^-)\\
            &=\omega(\hat{q}^-\#u\# (\hat{q}^+)^{-1})+\cA_{H_+}(q^+, \hat{q}^+) - \cA_{H_-}(q^-, \hat{q}^-).
        \end{align*}
        Since $(M,\omega)$ is monotone with monotonicity constant $\kappa$ we get that
        $$ E_{top}(u)=\cA_{H_+}(q^+, \hat{q}^+) - \cA_{H_-}(q^-, \hat{q}^-)+\kappa c_1(\hat{q}^-\#u\# (\hat{q}^+)^{-1}).$$
        Now, as $u=\gamma_0\#u_1\#\gamma_1\#\cdots\#\gamma_{m-1}\#u_m\#\gamma_m$, we get that 
        $$c_1(\hat{q}^-\#u\# (\hat{q}^+)^{-1})=c_1^\tau(\hat{q}^-)+c_1^\tau((\hat{q}^+)^{-1})+\sum_{j=0}^m c_1^\tau(\gamma_j)+\sum_{j=1}^m c_1^\tau(u_j).$$
        By assumption, since $\tau$ computed at $q^-$ and $q^+$ is homotopic to the trivializations induced by $\hat{q}^-$ and $\hat{q}^+$, respectively, we deduce that $c_1^\tau(\hat{q}^-)=c_1^\tau((\hat{q}^+)^{-1})=0$. Furthermore, since for every $0\leq k\leq m$ we have that $\tau_{q_k^-}$ is homotopic to $\gamma_k^*\tau_{q_k^+}$, we deduce that $c_1^\tau(\gamma_0)=\cdots=c_1^\tau(\gamma_m)=0$. Consequently,$$c_1(\hat{q}^-\#u\# (\hat{q}^+)^{-1})=\sum_{j=1}^m c_1^\tau(u_j).$$
        Hence, by the local dimension formula~\eqref{eq: index formula for continuation} we get that
        $$ E_{top}(u) = \cA_{H_+}(q^+, \hat{q}^+) - \cA_{H_-}(q^-, \hat{q}^-) + \frac{1}{2}\cdot\kappa \left( \dim_u \widehat\cM(q^-, q^+) - \mu_{FMB}^\tau(q^+) + \mu_{FMB}^\tau(q^-) \right),$$
        as required.
\end{proof}

    \subsection{Cubes}\label{ss:cubes}
        This section deals with algebraic, Hamiltonian and Floer theoretic cubes. We begin with algebraic cubes:
    \subsubsection{Algebraic cubes}\label{sss:algebraic cubes}
        The language of cubes is very convenient in order to work with the algebraic data arising when defining relative symplectic cohomology. Throughout we fix a unital commutative ring $R$ and we work in the category of graded $R$-modules and graded maps, where the grading is over $\Z$ or over $\Z_k$ with $k$ even. Recall that given graded modules $C,D$ with graded components $C^i,D^i$, a module map $f \fc C \to D$ is \textbf{graded of degree $d$} if $f(C^i) \subset D^{i+d}$. We denote the degree of $f$ by $|f|$.
        
        Fix a nonnegative integer $n$. We call a subset of $[0,1]^n$ a \textbf{face} if it is given by setting some of the coordinates to either $0$ or $1$; the rest of the coordinates are referred to as the \textbf{free coordinates of $F$}. Given a face $F$ its \textbf{dimension}, denoted $|F|$, is the number of its free coordinates, while its \textbf{initial vertex} $\ini F$ and \textbf{terminal vertex} $\ter F$, are the points of $F$ closest to or farthest from the origin, respectively (relative to the Euclidean metric); in this case we write $F \fc v \to v'$, where $v = \ini F$, $v' = \ter F$. Note that $F$ is determined by its initial and terminal vertices. We also say that $v,v'$ \textbf{span} $F$. If $F',F''$ are faces, we write $F = F'\cdot F''$ to denote the situation in which
        $$\ini F = \ini F'\,, \quad \ter F' = \ini F''\,,\quad \ter F'' = \ter F\,.$$
        
        \textbf{An (algebraic) $n$-cube $\cC$} is a pair $(\{\cC^v\}_{v\in\{0,1\}^n},\{f^\cC_F\}_{F\subset [0,1]^n\text{ a face}})$, where each $\cC^v$ is a graded module, and for every face $F \subset [0,1]^n$ we have a graded module morphism $f^\cC_F \fc \cC^{\ini F} \to \cC^{\ter F}$ of degree $1-|F|$, subject to the condition that for each face $F$ we have
        \begin{equation}\label{eqn:cube_relation}
        \sum_{F = F'\cdot F''}(-1)^{|F'|}\sgn(F',F'')f^\cC_{F''} f^\cC_{F'} = 0\,,
        \end{equation}
        where $\sgn(F',F'')$ is a sign defined as follows. Any face of $[0,1]^n$ comes equipped with a natural orientation coming from the ordering of its free coordinates. Then $\sgn(F',F'')$ is the intersection index of $F',F''$ inside $F$. For a more explicit description, assume first that we have a linearly ordered finite set $S=\{s_1<\dots<s_{k+l}\}$. Recall that a $(k,l)$-shuffle on $S$ is a permutation $\sigma \in S_S$ such that $\sigma(s_1)<\dots<\sigma(s_k)$ and $\sigma(s_{k+1})<\dots<\sigma(s_{k+l})$. If $S=S'\sqcup  S''$ with $|S'| = k$, $|S''|=l$, then there is a unique $(k,l)$-shuffle $\sigma_{S',S''}$ such that $\sigma(\{s_1,\dots,s_k\})=S'$, $\sigma(\{s_{k+1},\dots,s_{k+l}\})=S''$. Now let $S\subset \{1,\dots,n\}$ be the set of free coordinates of $F$, $S'$ the set of free coordinates of $F'$ and $S''$ the set of free coordinates of $F''$, so that $S = S'\sqcup S''$. Endow $S$ with the order induced from $\{1,\dots,n\}$. Then $\sgn(F',F'') := \sgn \sigma_{S',S''}$.
        
        The above relations mean in particular that for each vertex $v$, $(\cC^v,f^\cC_{\{v\}})$ is a cochain complex, each $1$-dimensional edge yields a cochain map between its vertex modules, each $2$-dimensional face yields a homotopy between the two compositions of the maps running along its perimeter, and so on.
        
        Note that the above sign only depends on the internal ordering of the free coordinates of $F',F''$, which means that for any face $G\subset [0,1]^n$ the pair $(\{\cC^v\}_{v},\{f^\cC_F\}_F)$ where $v$ ranges over the vertices of $G$ while $F$ over the subfaces of $G$, is itself a cube, provided we renumber the coordinates in their natural order in $\{1,\dots,n\}$. We will refer to such a cube as a \textbf{subcube of $\cC$}, or, when the face $G$ is given, \textbf{the subcube obtained by restricting $\cC$ to $G$}, and we denote it by $\cC|_G$.
        
        If $\cA_0,\cA_1$ are $n$-cubes, a \textbf{map from $\cA_0$ to $\cA_1$} is an $(n+1)$-cube $\cC$ such that $\cA_i = \cC|_{[0,1]^n\times\{i\}}$, $i=0,1$. In this case we write $\cA_0\xrightarrow{\cC}\cA_1$. 
        
    \subsubsection{Hamiltonian and Floer theoretic cubes}\label{sss:Ham cubes FMB}
        Fix a natural number $n\in \N$. Define $g\fc [0,1]^n \to \R$ by $g(x_1,\ldots,x_n)=\sum_{i=1}^n \cos(\pi x_i)$         for every $(x_1,\ldots,x_n)\in [0,1]^n$. This is a Morse function whose critical points are the vertices of the cube. Together with the Euclidean metric it forms a Morse-Smale pair, and its negative gradient flow preserves the cube.

        A \textbf{monotone cube of Hamiltonians} is a suitably smooth map $\cH\fc [0,1]^n\to C^\infty(M)$, which is monotone nondecreasing along the negative gradient lines of $g$. Similarly to Varolg\"une\c s' definition, \cite{Varolgunes_2018_PhD,Varolgunes_2021_MV_and_relSH}, we assume $\cH$ is constant on a neighborhood of each vertex. Moreover, since we work in the Floer--Morse--Bott setting, we assume that the Hamiltonians at the vertices satisfy the \textbf{MB} condition, in this case we say that $\cH$ satisfies the \textbf{MB} condition.

        A \textbf{cube of compatible almost complex structures} is a suitably smooth map $\cJ\colon [0,1]^n\to \mathscr{J}(M,\omega)$. we assume $\cJ$ is constant on a neighborhood of each vertex.

      In what follows, we describe Floer theoretic cubes associated to a generic choice of a monotone cube of Hamiltonians $\cH$ and a cube of compatible almost complex structures $\cJ$. For every vertex $v\in[0,1]^n$, choose a Morse function $h_v$ on $\bS_{\cH(v)}$, and denote by $C_v$ the Floer complex $CF^*(\cH(v))$ with the differential given by counting Floer flowlines with cascades, as in Section~\ref{ss:CF_MB_over_Z}.

        For a face $F\subset [0,1]^n$, denote $v=\operatorname{ini}(F)$ and
        $v'=\operatorname{ter}(F)$. Let us define the morphism
        \[
        f_F\fc C^v \to C^{v'}.
        \]
If $F$ is a vertex, then $v=v'$ and $f_F$ is the differential of the complex $C^v=CF^*(\cH(v);\Z)$. Assume that $F$ has positive dimension.

        Given a negative gradient flowline $\gamma$ of $g$ on the cube, $\gamma\colon \R \to [0,1]^n$, $\dot \gamma= -\nabla g \circ \gamma$, we define the following homotopies of Hamiltonians and almost complex structures as follows:
        \[\cH_\gamma=\cH\circ \gamma,\qquad\cJ_\gamma=\cJ\circ \gamma.\]
     
        For every $q^-\in \Crit(h_v)$ and $q^+\in \Crit(h_{v'})$, define the following moduli space
        \[
        \cM(q^-,q^+;F) = \left\{
            (\gamma,u)\,\middle|\,
            \begin{array}{c}
                 
                 \gamma\fc \mathbb{R}\to F,\quad \dfrac{d}{ds}\gamma=-\nabla g\circ\gamma,\\[0.3em]
                 \lim\limits_{s\to -\infty} \gamma(s)=v,\; \lim\limits_{s\to +\infty} \gamma(s)=v',\\[0.3em]
                 u \in \widehat{\cM}^{\Cont}(q^-,q^+;\cH_\gamma,\cJ_\gamma).
            \end{array}
            \right\}.
        \]
        The moduli space $\cM(q^-,q^+;F)$ admits an $\R$ action, where $\R$ acts by additive shifts on $\gamma$, and on a representative of $u$ by additive shift of the cascade that is a continuation solution. We define $\widehat{\cM}(q^-,q^+;F) = \cM(q^-,q^+;F) / \R$ as a quotient of the $\R$ action.

        Given a family $\tau$ of compatible trivializations along $\bS_{\cH(v)}$ for every vertex $v\in\{0,1\}^n$, for every $(\gamma,u) \in \cM(q^-,q^+;F)$, the local dimension of $\cM(q^-,q^+;F)$ at $(\gamma,u)$ is given by 
        \[
        \mu^\tau_{FMB}(q^+;\cH(v')) - \mu^\tau_{FMB}(q^-;\cH(v)) + 2\sum_{k=1}^{m} c_1^\tau(u_k) + \dim F,\]
        where $m$ is the number of cascades in $u$ and $\left(u_k\right)_{k=1}^{m}$ are the cascades of a representative of $u$, and for $\widehat{\cM}(q^-,q^+;F)$ the dimension is locally given by 
        \[
        \mu^\tau_{FMB}(q^+;\cH(v')) - \mu^\tau_{FMB}(q^-;\cH(v)) + 2\sum_{k=1}^{m} c_1^\tau(u_k) + \dim F - 1.\]

        We denote by $\widehat{\cM}^0(q^-,q^+;F)$ the union of the zero dimensional components of $\widehat{\cM}(q^-,q^+;F)$. Note that if $(\gamma,u)$ represent a class in $\widehat{\cM}^0(q^-,q^+;F)$ then 
        \[
        \mu^\tau_{FMB}(q^+;\cH(v')) - \mu^\tau_{FMB}(q^-;\cH(v)) + 2\sum_{k=1}^{m} c_1^\tau(u_k) = 1 - \dim F.\]

\begin{rem}\label{rem: index formula for homotopies} 
Given a monotone cube of Hamiltonians $\cH$ satisfying the \textbf{MB} condition and a cube of almost complex structures $\cJ$ such that $(\cH,\cJ)$ is regular, similarly to Remark~\ref{rem: index formula for Floer and cont flowlines}, for every face $F\in[0,1]^n$ and critical points $q^-\in \Crit(h_{v}), q^+\in \Crit(h_{v'})$, where $v=\ini(F)$ and $v'=\ter(F)$, the \textbf{index formula} for $(\gamma,u)\in \wh\cM(q^-,q^+;F)$ is given by the equality
\begin{equation}
\dim_{(\gamma,u)}\wh\cM(q^-,q^+;F)=\mu^\tau_{FMB}(q^+;\cH(v')) - \mu^\tau_{FMB}(q^-;\cH(v)) + 2\sum_{k=1}^{m} c_1^\tau(u_k) + \dim F-1,
\end{equation}
where $u_1,\ldots,u_m$ are the cascades of $u$. Also, it is important to note that the left-hand side of the index formula is independent of the compatible symplectic trivialization $\tau$.
\end{rem}
        The map $f_F\fc C^v \to C^{v'}$ is defined to be the unique homomorphism of $\Z$-modules that satisfies 
        
        \[f_F(x) = \sum_{y\in\Crit(h_{v'})}\#\widehat{\cM}^0(x,y;F)\,\cdot\,y\,,\]
        for every generator $x\in \Crit(h_v)$ of the Floer complex $CF^*(\cH(v);\Z)$.

        Recall that the face maps that correspond to each vertex are defined to be the Floer differential. Also, note that by unraveling the definition, the face maps that correspond to an edge are the continuation maps.
        
        The cube relations follow from compactifying the moduli spaces $\widehat{\cM}(q^-,q^+;F)$ by broken flowlines corresponding to the possible decompositions of the face $F$. Transversality guarantees that no bubbling occurs in such broken configurations via the classical codimension argument. In a monotone symplectic manifold, sphere bubbling is a codimension-$2$ phenomenon: removing the sphere bubbles from a bubbled broken configuration yields an element of a moduli space whose dimension is at least two less than that of the original configuration. Consequently, such bubbled configurations cannot appear in the relevant $1$-dimensional moduli spaces.

        To achieve transversality, hence regularity of the moduli spaces involved, we use a generic perturbation of the cube of almost complex structures, supported near the vertices, and a generic perturbation of the cube of Hamiltonians, supported away from the vertices.

\section{Geometric setting and acceleration data}\label{s:Charts and acc. data}
Consider the complex projective space, $\CP^n$, equipped with the Fubini--Study
symplectic form. We normalize it by identifying $\CP^n$ with the
symplectic cut of the Euclidean ball $B(1)\subset\C^n$; that is, the boundary of $B(1)$ is collapsed to
$\CP^{n-1}$ along the fibers of the Hopf fibration (see \cite{Lerman_1995_symplectic_cuts}). Here
$$
B(r)=\{\,z\in \C^n : \pi\|z\|^2 \le r\,\}.
$$
We denote by $q\fc B(1)\to \CP^{n}$ the quotient map obtained by collapsing the Hopf fibers on the boundary.

Consider also the map $\bar\rho(z)\colon B(1) \to \left(\C^{n+1}\setminus 0 \right) \big / \C^*$, given in homogeneous coordinates by
$$
\bar\rho(z)=[\sqrt{1/\pi - \|z\|^2} : z],
$$
for every $z\in B(1)$. 
On the boundary $\partial B(1)$, this becomes $\bar\rho(z)=[0 : z]$, hence $\bar\rho$ is constant on Hopf fibers on the boundary. Moreover the map $\bar\rho$ is bijective onto its image on the interior of $B(1)$. Hence, $\bar\rho$ descends to a map from the symplectic cut, $\rho \colon (B(1) \big/\sim) \to \CP^n$, such that $\bar\rho = \rho \circ q$.
\[
\xymatrix@C=4.5em@R=2.5em{
  & \mathbb{C}P^n
    = \bigl(\mathbb{C}^{n+1}\setminus\{0\}\bigr)\big/\mathbb{C}^*
  \\
  B(1) \ar[ru]^{\bar\rho} \ar[rd]_{q}
  &
  \\
  & \mathbb{C}P^n = B(1)\big/{\sim}
    \ar@{-->}_{\rho}[uu]
}
\]
We thus identify these two models of $\CP^n$, see \cite[Page 249]{Lerman_1995_symplectic_cuts}. We discuss the embedding 
$$
\iota = \rho|_{\Int B(1)}\fc \Int B(1)\to \CP^n
$$
in Claim~\ref{claim: Int B(1)->CPn-D_infty} below.

We denote by $D_\infty$ the copy of $\CP^{n-1}$ obtained from the quotient of $\partial B(1)$ by collapsing the fibers of the Hopf fibration; that is, $D_\infty=\rho(\partial B(1))$. Also, for every $1\le j\le n$, we denote by $D_j$ the image of the coordinate hyperplane $\lbrace z_j = 0 \rbrace$ under $\rho$, i.e.
$$
D_j=\rho\!\left(\{(z_1,\ldots,z_n)\in B(1) : z_j=0\right).
$$

The minimal first Chern number of $\CP^{n}$ is known to be $n+1$. Moreover, it is well known that
$$
\PD(c_1)=[D_1]+\cdots+[D_n]+[D_\infty]=(n+1)\cdot[D_\infty].
$$
This fact will be very useful in Section~\ref{s: obstructions}.

$B(1)$ admits a Hamiltonian $S^1$-action, given by $\theta\cdot z=(e^{i\theta}z_1,\ldots,e^{i\theta}z_n)$ for every $\theta\in S^1=\R/2\pi\Z$ and $z=(z_1,\ldots,z_n)\in B(1)\subset \CP^n$. This action descends to $\CP^n$ by the symplectic cut construction. We denote by $\mu\fc \CP^n\to \R$ its moment map, normalized such that it satisfies
$$
\pi\|z\|^2=\mu\circ\rho(z),
$$
for every $z\in B(1)$. Consequently, the image of $\mu$ is the interval $[0,1]$. Note that $\mu$ is invariant under the action of $\U(n)$ on $\CP^n$, induced from the action of $\U(n)$ on $\C^n$ via the symplectic cut.

\begin{rem}\label{rem: dynamics of a radial Ham}
 Note that for every smooth function $f\fc [0,1]\to \R$, the composition 
$$H=f\circ \mu\fc \CP^n\to \R$$ 
defines an autonomous Hamiltonian on $\CP^n$. The divisor $D_\infty$ consists of critical points of $H$, and hence its Hamiltonian flow is constant on $D_\infty$. Over $\Int B(1)\subset \CP^n$ the Hamiltonian vector field is given by
$$
X_H(z)=f'(\mu(z)) X_\mu(z)=-2\pi if'(\mu(z)) z,
$$
for every $z\in \Int B(1)$, and thus the Hamiltonian flow of $H$ is given by 
$$\varphi_H^t(z)=e^{-2\pi i f'(\mu(z)) t} z,$$ for every $z\in \Int B(1)$ and $t\in \R$. Therefore the fixed points of its time-$1$ map decompose into the following three families:  
\begin{enumerate}
    \item $0$, which is the image of the origin in $B(1)\subset\C^n$ under $\rho$;
    \item points in $D_\infty$;  
    \item points of the $(2n-1)$-sphere $\mu^{-1}(\{r\})$ for some $r\in(0,1)$ for which $f'(r)\in \Z$.
\end{enumerate}
Note that $\{0\}\cup D_\infty$ is the set of critical points of $\mu$, and thus the Hamiltonian flow of $H$ is actually constant on this set.

If $x$ is a $1$-periodic orbit of $H$ contained in $\mu^{-1}(r)$ for some $r \in [0,1)$, it admits a capping $\hat{x}$ contained in $\mathrm{Int}\, B(1)$. In this region, $\omega$ is exact with primitive $\lambda$. We can then deduce by Stokes' theorem that
$$\cA_H (x,\hat{x}) = \int_0^1 H(x(t))\,dt - \int_D w^*\omega = f(r) - \int_{S^1} x^* \lambda = f(r) - f'(r)r,$$
where $w$ represents the capping $\hat{x}$. 

Note that constant orbits can be capped by constant capping. For such a capping, their action is simply the value of $H$ at the constant orbit.

\end{rem}

\subsection{Three charts}\label{ss: charts}

As mentioned in the introduction, our strategy for proving the nonexistence of some Floer and continuation flowlines relies on estimates of intersection numbers with divisors, combined with the index formula (see Remark~\ref{rem: index formula for Floer and cont flowlines}). A crucial input for it, is a choice of trivialization and the computation of the Robbin--Salamon indices of the $1$-periodic orbits with respect to the chosen trivialization. In our setting, we consider three different coordinate charts, each over a different subset of $\CP^n$. In this section, we describe these charts, trivializations arising from them, and discuss the geometric properties of the associated subsets.

\subsubsection{The complement of $D_\infty$}\label{sss: CP^n - D_infty}
Since we constructed $\CP^n$ by a symplectic cut of the ball $B(1)$ along its boundary, it follows from \cite{Lerman_1995_symplectic_cuts} that the open ball $\Int B(1)$ is symplectomorphic to the subset
$$
U_0=\{[z_0:\cdots:z_n]\in \CP^n \,|\, z_0\neq 0\}\subset \CP^n,
$$
where $[z_0:\cdots:z_n]$ denotes a point of $\CP^n$ in homogeneous coordinates. It will be useful for us to describe this symplectomorphism in terms of coordinates.

Consider the map $\iota=\rho|_{\Int B(1)}\fc \Int B(1)\to \CP^n$ given by
$$\iota(z)=\left[\sqrt{1/\pi-\|z\|^2}:z\right]=\left[\sqrt{1/\pi-(|z_1|^2+\cdots+|z_n|^2)}:z_1:\cdots:z_n\right],$$

for every $z=(z_1,\ldots,z_n)\in \Int B(1)$. 
\begin{claim}\label{claim: Int B(1)->CPn-D_infty}
    The map $\iota$ is a symplectic embedding and its image is precisely $U_0$.
\end{claim}
\begin{proof}
    First, we show that $\iota$ is well-defined. Indeed For every $z\in \Int B(1)$ we have $\pi\|z\|^2<1$, and hence $1/\pi-\|z\|^2>0$.

Next we show injectivity. Let $u,v\in \Int B(1)$ and assume that $\iota(u)=\iota(v)$. This means that there exists $\lambda\in \C\setminus\{0\}$ such that
$$
\lambda(\sqrt{1/\pi-\|u\|^2},u)=(\sqrt{1/\pi-\|v\|^2},v).
$$
Thus $v=\lambda u$, and we deduce that
$$
\lambda\sqrt{1/\pi-\|u\|^2}=\sqrt{1/\pi-\|v\|^2}=\sqrt{1/\pi-|\lambda|^2\|u\|^2}.
$$

Since the right-hand side is nonnegative, $\lambda\in\R$ and $\lambda\ge0$. Moreover, we have
$$
\lambda^2(1/\pi-\|u\|^2)=1/\pi-\|v\|^2=1/\pi-|\lambda|^2\|u\|^2=1/\pi-\lambda^2\|u\|^2,
$$
and hence $\lambda^2/\pi=1/\pi$, so $\lambda=1$. Therefore $u=v$, as claimed.

We now show that $\iota$ is surjective onto $U_0$. Let $z\in \Int B(1)$. Since $\pi\|z\|^2<1$, we have $\sqrt{1/\pi-\|z\|^2}>0$, and hence $\iota(z)\in U_0$. This shows that $\im \iota\subset U_0$.

We now prove the converse inclusion. Let $w\in U_0$. We can represent $w$ in homogeneous coordinates as $w = [1:w']$ for some $w'\in \C^n$. Now, consider
$$
z=\frac{w'}{\sqrt{\pi(1+\|w'\|^2)}}.
$$
We show that $z\in \Int B(1)$ and that $\iota(z)=w$.

We compute
$$
\pi\|z\|^2=\pi\frac{\|w'\|^2}{\pi(1+\|w'\|^2)}=\frac{\|w'\|^2}{1+\|w'\|^2}<1.
$$

and therefore $z\in \Int B(1)$. Finally, we compute

\begin{multline*}
\iota(z)
=\left[\sqrt{1/\pi-\|z\|^2}:z\right]
=\left[\sqrt{1/\pi-\frac{\|w'\|^2}{\pi(1+\|w'\|^2)}}:
         \frac{w'}{\sqrt{\pi(1+\|w'\|^2)}}\right]\\
=\left[\sqrt{\frac{1}{\pi(1+\|w'\|^2)}}:
          w'\sqrt{\frac{1}{\pi(1+\|w'\|^2)}}\right]
=[1:w']
=w.
\end{multline*}
This shows that $U_0\subset \im \iota$, and hence $\im \iota=U_0$, as claimed.

Now, let us show that $\iota$ is indeed a symplectic embedding. Let 
$$
\lambda_0=\frac{i}{4}\sum_{j=0}^n (w_j\, d\overline{w}_j - \overline{w}_j\,dw_j ),
$$
be the standard Liouville form on $\C^{n+1}$. We define $\alpha=\lambda_0|_{S^{2n+1}(1/\sqrt{\pi})}$ to be its restriction to the sphere of radius $1/\sqrt{\pi}$. The Fubini--Study form on $\CP^n$ (normalized such that the area of a line is 1) satisfies
$$
\pi^*\omega_{FS}\vert^{}_{S^{2n+1}(1/\sqrt{\pi})}=d\alpha.
$$
Define the map $s\fc \Int B(1)\to \C^{n+1}$ given by $s(z)=\left(\sqrt{1/\pi-\|z\|^2},z\right)$ for every $z\in \Int B(1)$. Let us show that $$s^*\alpha=\lambda_0',$$
where
$$
\lambda_0'=\frac{i}{4}\sum_{j=1}^n (z_j\, d\overline{z}_j - \overline{z}_j\,dz_j ),
$$
is the standard Liouville form on $\Int B(1)\subset \C^{n}$. Indeed, since $\sqrt{1/\pi-\|z\|^2}\in \R$ we get
\begin{align*}
 s^*\alpha
 &=\frac{i}{4}\left(\sqrt{1/\pi-\|z\|^2}\, d\left(\sqrt{1/\pi-\|z\|^2}\right)-\sqrt{1/\pi-\|z\|^2}\,d\left(\sqrt{1/\pi-\|z\|^2}\right)+\sum_{j=1}^n (z_j\, d\overline{z}_j - \overline{z}_j\,dz_j )\right)\\ 
 &=\frac{i}{4}\sum_{j=1}^n (z_j\, d\overline{z}_j - \overline{z}_j\,dz_j )
 =\lambda_0'.
\end{align*}
Since $\iota=\pi\circ s$ we deduce that
$$
\iota^*\omega_{FS}=s^*(\pi^*\omega_{FS})=s^*d\alpha=ds^*\alpha=d\lambda_0'=\omega_0.
$$
Therefore $\iota$ is indeed a symplectic embedding.
\end{proof}

As a consequence, we obtain a new expression of $\mu$, in homogeneous coordinates:
\begin{claim}\label{claim: mu is homo coordinates}
    For every $u\in \CP^n$ we have
    $$
    \mu(u)=\frac{|u_1|^2+\cdots+|u_n|^2}{|u_0|^2+|u_1|^2+\cdots+|u_n|^2},
    $$
    where $u$ is represented as $u=[u_0:\cdots:u_n]$ in homogeneous coordinates.
\end{claim}
\begin{proof}
    Let $u\in \CP^n$, represented in homogeneous coordinates as $u=[u_0:\cdots:u_n]$, for some $(u_0,u_1,\ldots,u_n)\in \C^{n+1}\setminus\{0\}$. If $u_0=0$, then $u\in D_\infty$, and thus
    $$
    \mu(u)=1=\frac{|u_1|^2+\cdots+|u_n|^2}{0^2+|u_1|^2+\cdots+|u_n|^2},
    $$
    as required. Otherwise, assume that $u_0\neq0$. Since $u=[u_0:\cdots:u_n]$ and $u_0\neq0$, we can define $u'\in \C^n$ by $u'=(u_1/u_0,\ldots,u_n/u_0)$, and we get that $u=[1:u']$. Now, define
    $$
    z=\frac{u'}{\sqrt{\pi(1+\|u'\|^2)}}\in \C^n
    $$
    and recall that in the proof of Claim~\ref{claim: Int B(1)->CPn-D_infty} we showed that $u=\iota(z)$ and that
    $$
    \pi\|z\|^2=\frac{\|u'\|^2}{1+\|u'\|^2}.
    $$
    Therefore, since $\iota$ is the restriction of $\rho$ onto $\Int B(1)$ we deduce that
    $$
    \mu(u)=\mu(\rho(z))=\pi\|z\|^2=\frac{\|u'\|^2}{1+\|u'\|^2}=\frac{\|u_0\cdot u'\|^2}{|u_0|^2+\|u_0\cdot u'\|^2}=\frac{|u_1|^2+\cdots+|u_n|^2}{|u_0|^2+|u_1|^2+\cdots+|u_n|^2},
    $$
    as required.
\end{proof}

It is also important for us to compute the pullback of the standard complex structure of $\CP^n$ under $\iota$.
\begin{claim}\label{claim: iota^* J}
    Let $J$ be the standard complex structure on $\CP^n$. Then $\iota^*J$ is given by
    $$ (\iota^* J)_z(V) = iV - \pi \text{Im}\langle z, V \rangle z + \pi \frac{\text{Re}\langle z, V \rangle}{1 - \pi\|z\|^2}iz, $$
    for every $z \in \Int B(1)$ and $V \in \C^n$. Here $\langle -,- \rangle$ is the standard Hermitian inner product on $\C^n$.
\end{claim}

\begin{proof}
    Let $J_0$ be the standard complex structure on $\C^n$. Let $s\fc \C^n \to \CP^n$ be the holomorphic embedding given by $s(z)=[1:z]$ for every $z\in \C^n$. Define $F\fc \Int B(1)\to \C^n$ by $F(z)=\frac{1}{\sqrt{1/\pi-\|z\|^2}}\cdot z$, for every $z\in \Int B(1)$. Note that $\iota=s\circ F$, and since $s$ is holomorphic, we get that
    $$\iota^*J=(s\circ F)^*J=F^* s^* J=F^* J_0.$$
    Define the functions $f\fc \Int B(1)\to \R$ and $g\fc \C^n\to \R$ by 
    $$f(z)=\frac{1}{\sqrt{1/\pi-\|z\|^2}},$$ 
    for every $z\in \Int B(1)$, and
    $$g(w)=\frac{1}{\sqrt{\pi(1+\|w\|^2)}},$$
    for every $w\in \C^n$.
    
    Then $F(z)=f(z)\cdot z$ for every $z\in \Int B(1)$. Additionally, the function $G\fc \C^n\to \Int B(1)$ given by $G(w)=g(w)\cdot w$ for every $w\in \C^n$, is the inverse of $F$.

    Hence, for every $z\in \Int B(1)$ and $V\in \C^n$, we have
    $$d_z F(V)=d_z f(V)\cdot z+f(z)V.$$
    Similarly, for every $w\in \C^n$ and $V\in \C^n$, we have 
    $$d_w G(V)=d_w g(V)\cdot w+g(w)V.$$

    A direct computation shows that for $z\in \Int B(1)$ and $V\in \C^n$,
    $$df_z(V) = f(z)^3 \cdot \text{Re}\langle z, V \rangle,$$
    and for $w\in \C^n$ and $V\in \C^n$, we have
    $$dg_w(V) = -\pi g(w)^3\cdot \text{Re}\langle w, V \rangle,$$
    where $\langle \cdot,\cdot \rangle$ is the standard Hermitian inner product on $\C^n$. Additionally, note that for every $z\in \Int B(1)$, we have
    $$g(F(z)) = \frac{1}{\sqrt{\pi(1 + \|f(z)z\|^2)}} = \frac{1}{\sqrt{\pi(1 + \frac{\|z\|^2}{1/\pi - \|z\|^2})}} = \sqrt{1/\pi - \|z\|^2} = \frac{1}{f(z)},$$
    which implies that $g(F(z))f(z)=1$.

    Let $z\in \Int B(1)$ and $V\in \C^n$.
    Then the almost complex structure $\iota^*J$ evaluated at the point $z$ and the vector $V$ is given by
    \begin{align*}
        (\iota^* J)_z(V)&=(F^* J_0)_z(V)\\
        &=d_{F(z)}G\circ (J_0)_{F(z)}\circ d_zF(V)\\
        &=d_{F(z)}G\circ (J_0)_{F(z)}(d_z f(V)\cdot z+f(z)V)\\
        &=d_{F(z)}G\circ (d_z f(V)\cdot iz+f(z)iV)\\
        &=d_{F(z)} g(d_z f(V)\cdot iz+f(z)iV)\cdot F(z)+g(F(z))(d_z f(V)\cdot iz+f(z)iV)\\
        &=d_z f(V)\cdot d_{F(z)} g(iz)\cdot f(z)\cdot z+
        f(z)^2\cdot d_{F(z)} g(iV)\cdot z\\
        &\quad +g(F(z))\cdot d_z f(V)\cdot iz+g(F(z))\cdot f(z)iV.
    \end{align*}

    Let us compute each term in this sum:
    \begin{itemize}
        \item Note that
        $$ d_{F(z)} g(iz)= -\pi g(F(z))^3\cdot \text{Re}\langle F(z), iz \rangle= -\pi g(F(z))^3\cdot \text{Re}(-i\langle f(z)\cdot z, z \rangle)=0,$$
        and hence 
        $$d_z f(V)\cdot d_{F(z)} g(iz)\cdot f(z)\cdot z=0.$$
        \item Since $f(z)\cdot g(F(z))=1$, we get that
        \begin{align*}
            f(z)^2\cdot d_{F(z)} g(iV)\cdot z&= -\pi f(z)^2\cdot g(F(z))^3\cdot \text{Re}\langle F(z), iV \rangle z\\
            &= -\pi g(F(z))\cdot \text{Re}\langle f(z)\cdot z, iV \rangle z\\
            &= -\pi f(z)\cdot g(F(z))\cdot \text{Re}\langle  z, iV \rangle z\\
            &=-\pi \text{Re}\langle z, iV \rangle z.
        \end{align*}
        
        \item and \hfill$\displaystyle\begin{aligned}[t]
            g(F(z))\cdot d_z f(V)\cdot iz&= g(F(z))\cdot f(z)^3 \cdot \text{Re}\langle z, V \rangle \cdot iz\\
            &= f(z)^2 \cdot \text{Re}\langle z, V \rangle \cdot iz.
        \end{aligned}$ \hfill \null
    \end{itemize}
    Therefore:
    $$(\iota^* J)_z(V)=iV-\pi \text{Re}\langle z, iV \rangle z+ f(z)^2 \cdot \text{Re}\langle z, V \rangle \cdot iz.$$
    Hence,
    $$ (\iota^* J)_z(V) = iV - \pi \text{Im}\langle z, V \rangle z + \pi \frac{\text{Re}\langle z, V \rangle}{1 - \pi\|z\|^2}iz,$$
    for every $z \in \Int B(1)$ and $V \in \C^n$.
\end{proof}

\subsubsection{The complement of the origin}\label{sss: CP^n - 0}

In this section we construct a biholomorphism between $\CP^{n}\setminus\{0\}$, where $0=[1:0:\cdots:0]$ is the image of the origin in $B(1)$ under $\rho$, and the hyperplane line bundle $\cO(1)$ over $\CP^{n-1}$. Recall that the hyperplane line bundle $p\fc\cO(1)\to\CP^{n-1}$ is defined by
$$
\cO(1)=\{(w,\lambda)\mid w\in\CP^{n-1},\ \lambda\in\Hom(w,\C)\}.
$$
Here, $w$ is a $1$-dimensional subspace of $\C^n$, and $\lambda$ is a linear functional on that line.

For every $(w,\lambda)\in\cO(1)$ and $z\in w\setminus\{0\}$, the vector $(\lambda(z),z)\in\C\times \C^n=\C^{n+1}$ is nonzero and hence determines a point in $\CP^{n}$. Define a map $\Phi\fc \cO(1)\to \CP^n$ as follows.
For $(w,\lambda)\in \cO(1)$, choose $z\in w\setminus\{0\}\subset \C^n$ and set
$$
\Phi(w,\lambda) = [\lambda(z):z].
$$
The map $\Phi$ is well-defined. Indeed, let $z,z'\in w\setminus\{0\}$.
Since $w$ is a complex line, there exists $c\in \C\setminus\{0\}$ such that
$z'=cz$. Then
$$
[\lambda(z'):z'] = [\lambda(cz):cz]
= [c\lambda(z):cz]
= [\lambda(z):z],
$$
showing that $\Phi(w,\lambda)$ is independent of the choice of $z$.

\begin{claim}\label{claim: biholo O(1)->CPn-0}
The map $\Phi\fc \cO(1) \to\CP^{n}\setminus\{0\}$ is a biholomorphism with respect to the standard complex structures of $\cO(1)$ and $\CP^n\setminus\{0\}\subset\CP^n$.
\end{claim}

\begin{proof}
    Let us check all the required properties:

\medskip

(1) \emph{Injectivity.} Let $(w,\lambda),(w',\lambda')\in \cO(1)$ and suppose $\Phi(w,\lambda)=\Phi(w',\lambda')$. Choosing representatives $z \in w\setminus\{0\}$ and $z' \in w'\setminus\{0\}$, the equality in $\CP^n$ implies there exists a scalar $c \in \C\setminus\{0\}$ such that
$$
(\lambda'(z'), z') = c(\lambda(z), z) \in \C \times \C^n.
$$
Comparing the vector components yields $z' = cz$, which implies $w' = w$. Comparing the scalar components yields $\lambda'(z') = c\lambda(z)$. Substituting $z'=cz$ into the left side gives
$$
\lambda'(z')=\lambda'(cz)=c\lambda'(z).
$$
Equating this with $c\lambda(z)$ implies $\lambda'(z)=\lambda(z)$. Note that $w=w'$ is a $1$-dimensional linear space and $z$ is a nonzero vector, hence it is a basis for $w$. Since the functionals agree on a basis, we conclude that $\lambda'=\lambda$. Thus, $(w,\lambda)=(w',\lambda')$. This shows that $\Phi$ is injective.

\medskip

(2) \emph{Surjectivity onto $\CP^{n}\setminus\{0\}$.} Let $[u_0 : u] \in \CP^n$ where $u_0\in \C$, $u\in \C^n$ and $(u_0,u)\in \C\times \C^n=\C^{n+1}$ is not the zero vector. Note that the excluded point $0$ corresponds to the case $u=0$. For any $u\in \C^n\setminus\{0\}$, let $w = \spn\{u\} \in \CP^{n-1}$. Let $\lambda \in \Hom(w, \C)$ be the unique linear functional satisfying $\lambda(u) = u_0$. Consider $(w,\lambda)\in \cO(1)$. Then
$$
\Phi(w, \lambda) = [\lambda(u) : u] = [u_0 : u].
$$
Thus, the image is exactly $\CP^{n}\setminus\{0\}$.

\medskip

(3) \emph{Holomorphicity.} We work in standard local charts $U_0,\ldots,U_{n-1}\subset \CP^{n-1}$ given by
$$U_j=\{[z_0:\cdots:z_{n-1}]\in \CP^{n-1}\,|\, (z_0,\ldots,z_{n-1})\in \C^n,\; z_j\neq 0\},$$
for every $j\in\{0,\ldots,n-1\}$.

We will prove that for every $j\in\{0,\ldots,n-1\}$ the restriction $\Phi|_{p^{-1}(U_{j})}\fc p^{-1}(U_{j})\to \CP^n$ is holomorphic. For simplicity of notation, we will focus on the case where $j=n-1$. All the other cases are similar. There exists a unique biholomorphism $y\fc U_{n-1}\to \C^{n-1}$ satisfying $w=\spn\{(y(w),1)\}\subset\C^n$ for every $w\in U_{n-1}$\footnote{The map $y$ is given by $y([u':u_{n-1}])=u'/u_{n-1}$ for every $u=[u':u_{n-1}]\in U_{n-1}\subset \CP^{n-1}$ where $u_{n-1}\in \C\setminus\{0\}$ and $u'\in \C^{n-1}$.}. Note that for every $w\in U_{n-1}$ the vector $(y(w), 1)$ is a basis for $w$. Thus, for every $c\in \C$ there exists a unique functional $\lambda_c\in \Hom(w,\C)$ satisfying $\lambda_c(y(w), 1)=c$.
Let $\varphi\fc U_{n-1}\times \C\to p^{-1}(U_{n-1})$ be defined by $\varphi(w,c)=(w,\lambda_c)$ for every $(w,c)\in U_{n-1}\times \C$. Thus $\varphi$ is a holomorphic local trivialization of $\cO(1)$ over $U_{n-1}$. Now, note that
$$
\Phi\circ\varphi(w, c) =\Phi(w,\lambda_c)=[\lambda_c(y(w), 1):y(w): 1]=[c:y(w): 1] \in \CP^n,
$$
for every $(w,c)\in U_{n-1}\times \C$. Since the coordinates here are holomorphic, we deduce that $\Phi|_{p^{-1}(U_{n-1})}$ is holomorphic. By a similar proof, one can show that $\Phi|_{p^{-1}(U_{j})}$ is holomorphic for every $j\in\{0,\ldots,n-1\}$. Since $\cO(1)$ is covered by the charts $p^{-1}(U_0),\ldots,p^{-1}(U_{n-1})$, we deduce that $\Phi$ is holomorphic.

\medskip

(4) \emph{Smoothness of the inverse map.} We have already established that $\Phi$ is a holomorphic bijection onto $\CP^n\setminus\{0\}$. Therefore, according to \cite[Page 19]{Griffiths_Harris_1978} we deduce that $\Phi^{-1}$ is holomorphic as well.
Thus we conclude that $\Phi$ is a biholomorphism, as required.
\end{proof}

Now we endow $\cO(1)$ with a specific symplectic form which we will denote by $\Omega$. Let us construct this form. Let $h$ be the Hermitian metric on the fibers of $\cO(1)$, induced by restricting to each fiber the standard Hermitian metric on $\C^n$. Explicitly, for a pair $(w, \lambda)\in \cO(1)$, for every $z \in w\setminus\{0\}$, the quotient $|\lambda(z)|^2/\|z\|^2$ does not depend on $z$, and thus we can define the squared-norm of $\lambda$ to be given by
$$
\|\lambda\|_h^2= \frac{|\lambda(z)|^2}{\|z\|^2_{\C^n}}.
$$

Since $\cO(1)$ is a holomorphic line bundle, its total space is a complex manifold. This structure induces the decomposition of the exterior derivative $d = \partial + \bar{\partial}$ into the Dolbeault operators. Additionally, since $\Phi$ is a biholomorphism, it respects the Dolbeault operators.

We define a $(1,1)$-form $\Omega$ on the total space of $\cO(1)$ by
\begin{equation}\label{eq: def of Omega}
\Omega := p^*\omega_{FS}^{\CP^{n-1}} + \frac{i}{2\pi}\partial\bar\partial \log\left(1 + \|\lambda\|_h^2\right),
\end{equation}
where $p\fc\cO(1)\to\CP^{n-1}$ is the bundle projection and $\omega_{FS}^{\CP^{n-1}}$ is the normalized Fubini--Study form on the base.

\begin{rem}\phantom{M}
\begin{itemize}
    \item The form $\Omega$ is not the one induced by the inclusion $\cO(1)\hookrightarrow \CP^{n-1}\times \Hom(\C^n,\C)$.
    \item  The form $\Omega$ admits a natural geometric interpretation. The first term is the lift of the symplectic structure from the base. The second term restricts to the normalized Fubini--Study form on each fiber (where the area of the line is 1).
\end{itemize}
\end{rem}

\begin{claim}
The pullback of the Fubini--Study form $\omega_{FS}$ on $\CP^n$ by the map $\Phi$ satisfies
$$
\Phi^*\omega_{FS} = \Omega.
$$
\end{claim}

\begin{proof}
    We compute the pullback in a local trivialization. Let $U_0,\ldots,U_{n-1}\subset\CP^{n-1}$ be the standard affine charts on the base defined by
    $$
    U_{j}=\{[z_0:\cdots:z_{n-1}]\in \CP^{n-1}\,|\,(z_0,\ldots,z_{n-1})\in \C^n,\, z_{j}\neq 0\},
    $$
    for every $j\in\{0,\ldots,n-1\}$. Similarly, let $V_0,\ldots,V_{n}\subset \CP^n$ be the standard affine charts on the target defined by
    $$
    V_k = \{[u_0:\cdots:u_n]\in \CP^n \,|\,(u_0,\ldots,u_n)\in \C^{n+1},\,u_k \neq 0\},
    $$
    for every $k\in\{0,\ldots,n\}$.

    We will perform the computation on the chart $U_{n-1}$. Note that the image of $p^{-1}(U_{n-1})$ under $\Phi$ lies in $V_n$. Recall that on $V_n$, the Fubini--Study form (normalized to area 1) can be expressed as 
    $$\omega_{FS} = \frac{i}{2\pi}\partial\bar\partial K,$$
    where $K\fc V_n \to \R$ is the K\"ahler potential given by
    $$
    K(u) := \log\left(1 + \sum_{j=0}^{n-1} \left|\frac{u_j}{u_n}\right|^2\right),
    $$
    for every $u=[u_0:\cdots:u_n]\in V_n$. Note that since for every $u\in V_n$ we have $u_n\neq0$ where $u_n$ is the $n$-th homogeneous coordinate of $u$, thus $K$ is well-defined.

    Since $\Phi$ is holomorphic (as shown in Claim~\ref{claim: biholo O(1)->CPn-0}), the pullback commutes with the Dolbeault operators:
    $$
    \Phi^*\omega_{FS} = \Phi^*\left(\frac{i}{2\pi}\partial\bar\partial K\right) = \frac{i}{2\pi}\partial\bar\partial (\Phi^* K).
    $$
    We now compute the function $\Phi^* K$ on $p^{-1}(U_{n-1})$. As in the proof of Claim~\ref{claim: biholo O(1)->CPn-0}, let $y\fc U_{n-1}\to \C^{n-1}$ be the unique biholomorphism satisfying $w=\spn\{(y(w),1)\}\subset\C^n$ for every $w\in U_{n-1}$. Recall the holomorphic local trivialization $\varphi\fc U_{n-1}\times \C\to p^{-1}(U_{n-1})$ defined by $\varphi(w,c)=(w,\lambda_c)$, where $\lambda_c\in \Hom(w,\C)$ is the unique functional satisfying $\lambda_c(y(w),1)=c$.
    
    Using this trivialization, the image is 
    $$
    \Phi(\varphi(w,c)) = \Phi(w,\lambda_c) = [\lambda_c(y(w),1) : y(w) : 1] = [c : y(w) : 1].
    $$
    Note that the last coordinate is $1$, which matches our normalization for $K$ on $V_n$. Thus,
    $$
    (\varphi^*\Phi^* K)(w,c) = K([c : y(w) : 1]) = \log\left(|c|^2 + \|y(w)\|^2 + 1\right).
    $$
    Let $\psi\fc U_{n-1} \to \R$ be the potential for the Fubini--Study form on the base $\CP^{n-1}$ defined by $\psi(w) = \log(1 + \|y(w)\|^2)$, such that $\frac{i}{2\pi}\partial\bar\partial \psi = \omega_{FS}^{\CP^{n-1}}$. We can rewrite the pullback of the potential as
    \begin{align*}
    (\varphi^*\Phi^* K)(w,c) &= \log\left( (1 + \|y(w)\|^2) \left(1 + \frac{|c|^2}{1 + \|y(w)\|^2}\right) \right) \\
    &= \psi(w) + \log\left(1 + \frac{|c|^2}{{1 + \|y(w)\|^2}}\right),
    \end{align*}
    for every $(w,c)\in U_{n-1}\times \C$.   The fiber coordinate $c$ over a base point $w\in \CP^{n-1}$ is the value of $\lambda_c$ on $(y(w),1)$, so the fiber norm satisfies 
    $$
    \|\lambda_c\|_h^2 = \frac{|\lambda_c((y(w),1))|^2}{\|(y(w),1)\|^2} = \frac{|c|^2}{1+\|y(w)\|^2}.
    $$ 
    Substituting this back, we get
    $$
    \varphi^*\Phi^* K(w,c) = \psi(w) + \log\left(1 + \|\lambda_c\|_h^2\right),
    $$
    for every $(w,c)\in U_{n-1}\times \C$, and hence
    $$
    \Phi^* K(w,\lambda) = \psi(w) + \log\left(1 + \|\lambda\|_h^2\right),
    $$
   for every $(w,\lambda)\in p^{-1}(U_{n-1})$. Applying $\frac{i}{2\pi}\partial\bar\partial$ to both sides yields
    \[
    \Phi^*\omega_{FS} = \Omega. \qedhere
    \]
\end{proof}

The following computation of a Hamiltonian vector field will be used to calculate Robbin--Salamon indices and to obtain estimates on the asymptotic intersection numbers of Floer solutions with the divisor $D_\infty=\CP^{n-1}$.

\begin{claim}\label{claim: mu in O(1)}
The Hamiltonian $H := \mu \circ \Phi \fc \cO(1) \to \R$ is given by 
$$
H(w, \lambda) = \frac{1}{1 + \|\lambda\|_h^2},
$$
for every $(w, \lambda) \in \cO(1)$. Moreover, the Hamiltonian vector field $X_H$ with respect to the symplectic form $\Omega$ is given by 
$$
X_H(w, \lambda) = \frac{d}{dt} \Big|_{t=0} (w, e^{2\pi it}\lambda)=2\pi i \lambda,
$$
for every $(w, \lambda) \in \cO(1)$.
\end{claim}

\begin{proof}
Let us first verify the expression for the Hamiltonian function. Recall that by Claim~\ref{claim: mu is homo coordinates}, the moment map $\mu \fc \CP^n \to [0,1]$ is given by
$$
\mu(u) = \frac{\|u'\|^2}{\|u'\|^2 + |u_0|^2},
$$
for every $u\in \CP^n$ represented in homogeneous coordinates as $u=[u_0:u']$, where $(u_0,u')\in \C\times \C^n$ is a non-zero vector.

For any $(w, \lambda) \in \cO(1)$, choose a representative $z \in w \setminus \{0\}$. By the definition of the map $\Phi$, we have $\Phi(w, \lambda) = [\lambda(z) : z]$. Substituting this into the expression for $\mu$, we obtain
$$
H(w, \lambda) = \mu([\lambda(z) : z]) = \frac{\|z\|^2}{\|z\|^2 + |\lambda(z)|^2} = \frac{1}{1 + \frac{|\lambda(z)|^2}{\|z\|^2}}.
$$
Recall that the squared-norm of the functional is defined by $\|\lambda\|_h^2 = |\lambda(z)|^2 / \|z\|^2$. Thus, we conclude $H(w, \lambda) = (1 + \|\lambda\|_h^2)^{-1}$, as required.

Next, we determine the Hamiltonian vector field $X_H$. Since $\Phi \fc (\cO(1), \Omega) \to (\CP^n \setminus \{0\}, \omega_{FS})$ is a symplectomorphism, it intertwines the Hamiltonian flows. Specifically, the Hamiltonian vector field $X_H$ satisfies $\Phi_* X_H = X_\mu$, where $X_\mu$ is the Hamiltonian vector field on $\CP^n$ associated with $\mu$.

The flow $\{\phi_\mu^t\}_{t\in \R}$ on $\CP^n$ generated by the moment map $\mu$ is given by 
$$
\phi_\mu^t(u) = [u_0 : e^{-2\pi it}u'],
$$
for every $t\in \R$ and for every $u\in \CP^n$ represented in homogeneous coordinates as $u=[u_0:u']$, where $(u_0,u')\in \C\times \C^n$ is a non-zero vector. Let $\{\psi^t_H\}$ be the flow on $\cO(1)$ generated by $X_H$. By the equivariance of $\Phi$, we must have $\Phi(\psi_H^t(w, \lambda)) = \phi_\mu^t(\Phi(w, \lambda))$ for every $(w,\lambda)\in \cO(1)$ and $t\in \R$. Note that
$$
\phi_\mu^t(\Phi(w, \lambda)) = [ \lambda(z) : e^{-2\pi it}z ],
$$
for every $(w,\lambda)\in \cO(1)$ and $t\in \R$. Since this is an element of $\CP^n$, we can multiply the representative by the non-zero scalar $e^{2\pi it}$ to obtain
$$
[ \lambda(z) : e^{-2\pi it}z ] = [ e^{2\pi it}\lambda(z) : z ] = \Phi(w, e^{2\pi it}\lambda).
$$
Because $\Phi$ is a biholomorphism, and in particular injective, we deduce that $\psi_H^t(w, \lambda) = (w, e^{2\pi it}\lambda)$ for every $(w,\lambda)\in \cO(1)$ and $t\in \R$. The vector field $X_H$ is the infinitesimal generator of this flow, which completes the proof.
\end{proof}

\subsubsection{The complement of the divisors}\label{sss: CPn - all divisors}

Denote by $D$ the union of the divisors $D_1,\ldots,D_n,D_\infty \subset \CP^n$. In this section we explore the complement of $D$. For this task, denote 
$$ V = \left\{(x_1,\ldots,x_n) \in \R^n \, : \, \sum_{j=1}^n x_j < \frac{1}{\pi} \,\,\text{and}\,\, \forall \, j \in \{1,\ldots,n\} \,\, \text{we have} \,\, x_j > 0 \right\}, $$
and define
$$ U = V \times \T^n $$
where $\T^n = (\R/\Z)^n$ is the $n$-torus. We can think about $U = V \times \T^n$ as a Lagrangian product, equipped with the symplectic form $\Omega_0 = \pi \sum_{j=1}^n dx_j \wedge d\theta_j$, where $\theta_1,\ldots, \theta_n \in \R/\Z$ are the angular coordinates on $\T^n$, and $x_1,\ldots,x_n$ are the standard coordinates of $V \subset \R^n$. Define a map $\Psi \fc U \to \CP^n$ by
$$ \Psi(x_1,\ldots,x_n,\theta_1,\ldots,\theta_n) = \left[ \sqrt{\frac{1}{\pi} - \sum_{k=1}^n x_k} : \sqrt{x_1} e^{2\pi i \theta_1} : \dots : \sqrt{x_n} e^{2\pi i \theta_n} \right], $$
for every $(x_1,\ldots,x_n,\theta_1,\ldots,\theta_n) \in U$.

\begin{claim}\label{claim: CPn - all divisors}
    The map $\Psi \fc U \to \CP^n$ is a symplectic embedding and its image is $\CP^n \setminus D$.
\end{claim}

\begin{proof}
    Define a map $F \fc U \to B(1)$ by 
    $$ F(x,\theta) = (\sqrt{x_1} e^{2\pi i \theta_1}, \ldots, \sqrt{x_n} e^{2\pi i \theta_n}), $$
    for every $(x_1,\ldots,x_n, \theta_1,\ldots,\theta_n) \in U$. We claim that $F$ is a symplectic embedding and its image is
    $$ \im F = \Int B(1) \setminus \bigcup_{j=1}^n \{ (z_1,\ldots,z_n) \in \C^n \, \colon \, z_j = 0 \}. $$

    Let us check the injectivity of $F$. Suppose that for $(x, \theta), (x', \theta') \in U$, we have $F(x, \theta) = F(x', \theta')$. This implies that for every $j \in \{1, \dots, n\}$, we have $\sqrt{x_j} e^{2\pi i \theta_j} = \sqrt{x'_j} e^{2\pi i \theta'_j}$. By taking the absolute value of both sides, we obtain $\sqrt{x_j} = \sqrt{x'_j}$, and since $x_j, x'_j > 0$, it follows that $x_j = x'_j$. Therefore $e^{2\pi i \theta_j} = e^{2\pi i \theta'_j}$. In the $n$-toral coordinates $\T^n = (\R/\Z)^n$, this equality implies $\theta_j = \theta'_j$ for every $j$. Thus $(x, \theta) = (x', \theta')$, proving that $F$ is injective.

    Now we show that
    \[ \im F = \Int B(1) \setminus \bigcup_{j=1}^n \{ (z_1,\ldots,z_n) \in \C^n \, \colon \, z_j = 0 \}. \]
    Let $z \in \im F$. There exists $(x, \theta) \in V \times \T^n$ such that $z = F(x, \theta)$. For all $j \in \{1,\ldots,n\}$, we have $z_j = \sqrt{x_j} e^{2\pi i \theta_j}$ and hence $x_j > 0$, implying $z_j \neq 0$. Furthermore, the condition $\sum_{j=1}^n x_j < \frac{1}{\pi}$ implies
    $$ \pi \|z\|^2 = \pi \sum_{j=1}^n |\sqrt{x_j} e^{2\pi i \theta_j}|^2 = \pi \sum_{j=1}^n x_j < 1. $$
    Thus $z \in \Int B(1)$, which shows $\im F \subset \Int B(1) \setminus \bigcup_{j=1}^n \{z_j=0\}$.
    
    Conversely, let $z = (z_1,\ldots,z_n) \in \C^n$ and assume $z_1, \ldots, z_n \neq 0$ and $\pi \|z\|^2 < 1$. For every $j \in \{1,\ldots,n\}$, define $x_j = |z_j|^2$ and $\theta_j = \frac{1}{2\pi} \arg(z_j) \in \R/\Z$. Since $z_j \neq 0$, we have $x_j > 0$. The condition $\pi \|z\|^2 < 1$ implies $\sum x_j = \sum |z_j|^2 < \frac{1}{\pi}$, hence $(x, \theta) \in V \times \T^n$. By construction, $F_j(x, \theta) = |z_j| e^{i \arg z_j} = z_j$, so $z \in \im F$. This confirms the image of $F$.

    The smoothness of $F$ is an immediate consequence of its formula. Let us check that $F$ is a symplectomorphism. Recall that the standard symplectic form on $B(1) \subset \C^n$ is $\omega_0 = \frac{i}{2} \sum_{j=1}^n dz_j \wedge d\bar{z}_j$. For every $j \in \{1,\ldots,n\}$, we have 
    $$ F^* dz_j = d(\sqrt{x_j} e^{2\pi i \theta_j}) = \frac{1}{2\sqrt{x_j}} e^{2\pi i \theta_j} dx_j + 2\pi i \sqrt{x_j} e^{2\pi i \theta_j} d\theta_j, $$
    and
    $$ F^* d\bar{z}_j = d(\sqrt{x_j} e^{-2\pi i \theta_j}) = \frac{1}{2\sqrt{x_j}} e^{-2\pi i \theta_j} dx_j - 2\pi i \sqrt{x_j} e^{-2\pi i \theta_j} d\theta_j. $$
    The wedge product is:
    \begin{align*}
        F^*(dz_j \wedge d\bar{z}_j) &= \left( \frac{e^{2\pi i \theta_j}}{2\sqrt{x_j}} dx_j + 2\pi i \sqrt{x_j} e^{2\pi i \theta_j} d\theta_j \right) \wedge \left( \frac{e^{-2\pi i \theta_j}}{2\sqrt{x_j}} dx_j - 2\pi i \sqrt{x_j} e^{-2\pi i \theta_j} d\theta_j \right) \\
        &= - \left( \frac{2\pi i}{2} + \frac{2\pi i}{2} \right) dx_j \wedge d\theta_j = -2\pi i dx_j \wedge d\theta_j.
    \end{align*}
    Hence, $F^* \omega_0 = \frac{i}{2} \sum_{j=1}^n (-2\pi i dx_j \wedge d\theta_j) = \pi \sum_{j=1}^n dx_j \wedge d\theta_j = \Omega_0$.
    
    By Claim~\ref{claim: Int B(1)->CPn-D_infty}, $\iota \fc \Int B(1) \to \CP^n$ is a symplectic embedding with image $\CP^n \setminus D_\infty$. Since $\Psi = \iota \circ F$, it is a symplectic embedding as a composition of symplectic embeddings. Its image is
    \[ \im \Psi = \iota \left( \Int B(1) \setminus \bigcup_{j=1}^n \{z_j=0\} \right) = \CP^n \setminus (D_\infty \cup D_1 \cup \dots \cup D_n) = \CP^n \setminus D. \qedhere\]
\end{proof}

\begin{rem}\label{rem: mu in CPn - all divisors}
    Note that for every $(x_1,\ldots,x_n,\theta_1,\ldots,\theta_n)\in U$ we have $$\mu\circ\Psi(x_1,\ldots,x_n,\theta_1,\ldots,\theta_n)=\pi(x_1+\cdots+x_n).$$
    Thus for every smooth function $f\fc \R\to \R$ we have $$X_{f\circ\mu\circ\Psi}(x_1,\ldots,x_n,\theta_1,\ldots,\theta_n)=-f'(\pi(x_1+\cdots+x_n))\left(\partial\theta_1+\cdots+\partial\theta_n\right),$$
    for every $(x_1,\ldots,x_n,\theta_1,\ldots,\theta_n)\in U$.
\end{rem}

We will use the map $\Psi$ to define a trivialization of $TM\big\vert_{\CP^n\setminus D}$ in Section \ref{ss: c_1 and intersection}.

\subsection{Acceleration data}\label{ss: acc. data}
In order to prove Theorem~\ref{thm: SH of ball in CP^n} we have to choose acceleration data for the ball of capacity $\Delta\in(0,1)$.

Let $h\fc (0,1)_\Delta\times[0,+\infty)_s\times [0,1]_r\to \R$ be a smooth function, satisfying the following:
\begin{itemize}
    \item  $\frac{\partial h}{\partial \Delta}\leq0$, $\frac{\partial h}{\partial s}\geq0$ and  $\frac{\partial^2}{\partial r^2}h\geq0$.
    \item For every $(\Delta,r)\in(0,1)\times[0,1]$ we have $\lim_{s\to\infty} h(\Delta,s,r)=\left\{\begin{array}{cc}
        0, & r\leq \Delta, \\
        +\infty, & r>\Delta.
    \end{array}\right.$
    
    \item For every $(\Delta,s)\in(0,1)\times[0,+\infty)$ and $r\in [0,\Delta)$ we have $h(\Delta,s,r)<0$, $\frac{\partial}{\partial r} h(\Delta,s,r)\in (0,1)$ and $\frac{\partial^2}{\partial r^2} h(\Delta,s,r)=0$;
    \item There exists a smooth function $\delta\fc(0,1)_\Delta\times[0,+\infty)_s\to \R$ satisfying 
    \begin{itemize}
        \item $\frac{\partial}{\partial s}\delta<0$;
        \item For every $(\Delta,s)\in(0,1)\times[0,+\infty)$ we have $0<\delta(\Delta,s)\leq1-\Delta$;
        \item For every $\Delta\in(0,1)$ we have  $\lim_{s\to\infty}\delta(s)=0$
    \end{itemize}
    such that for every $(\Delta,s)\in(0,1)\times[1,+\infty)$ and $r\in[\Delta+\delta(\Delta,s),1]$ we have $\frac{\partial}{\partial r} h(\Delta,s,r)=s+\frac{1}{2}$, and for every $r\in(\Delta,\Delta+\delta(\Delta,s))$, we have $\frac{\partial^2}{\partial r^2}h>0$. 
     \item For every $(\Delta,r)\in(0,1)\times [0,1]$ we have $h(\Delta,0,r)<0$ and $\frac{\partial}{\partial r} h(\Delta,0,r)\in (0,1)$.
\end{itemize}

Now, for every $\Delta\in(0,1)$ we use the acceleration datum $(H_\ell)_{\ell=0}^\infty$, consisting of autonomous Hamiltonians associated to the ball of capacity $\Delta$ in $\CP^n$, given by 
\begin{equation}\label{eq: acc. data for balls}
    H_\ell(z)=h(\Delta,\ell,\mu(z)),
\end{equation}
for every $z\in \CP^n$ and $\ell\in\Z_{\geq0}.$ The function $h$ allows us to encapsulate in a single function the Hamiltonians defining the acceleration data for the ball of capacity $\Delta$, the homotopies between the Hamiltonians within the acceleration data (by varying $\ell$), as well as the homotopies between acceleration data corresponding to two balls of different capacities, by varying $\Delta$.) This will be crucial for constructing the restriction maps used in the proof of Theorem~\ref{thm: SH of ball in CP^n} and Theorem~\ref{thm: res for balls in CP^n}. 
\begin{rem}\label{rem: H_0}
    Note that $H_0$ is a non-positive $C^2$-small Morse--Bott function.
\end{rem}

Let us describe the generators of the Floer--Morse--Bott complexes for this acceleration datum. Given $\Delta\in(0,1)$ and $\ell\geq0$, the Hamiltonian $H_\ell$ has $\ell+2$ critical submanifolds, which come in three types of families.
\begin{enumerate}
    \item The global minimum at $0$,
    \item $\ell$ spheres, $S^{2n-1}$, consisting of the points where the slope is integral,
    \item and the global maxima, forming a copy of $\CP^{n-1}$ at the divisor at infinity.
\end{enumerate}
By choosing a perfect Morse function on each of these critical submanifolds we get the following generators: The first family contributes a single generator which we denote by $\check{x}_0^\ell$; each sphere contributes a pair of generators corresponding to the minimum and maximum of a perfect Morse function, yielding in total $\check{x}_1^\ell,\hat{x}_1^\ell,\ldots,\check{x}_\ell^\ell,\hat{x}_\ell^\ell$, where the hats correspond to maxima and the checks to minima; lastly, the divisor at infinity contributes $n$ generators, which we denote by $\check{x}_{\ell+1}^\ell,\ldots,\check{x}_{\ell + n}^\ell$.

In Section~\ref{s: indices} we compute the Floer--Morse--Bott indices of these generators with respect to two different trivializations. Later in Section~\ref{s: obstructions} and Section~\ref{s: computations of CF} we compute the differential of this complex and continuation maps between such complexes.

\section{Indices}\label{s: indices}
This section is dedicated to computing the Robbin--Salamon and Floer--Morse--Bott indices of the $1$-periodic orbits for our acceleration data from Section~\ref{ss: acc. data}. Section~\ref{ss: RS part II} continues the discussion of the Robbin--Salamon index for paths of symplectic matrices initiated in Section~\ref{sss: RS part I}, and lists its basic properties. Section~\ref{ss:RS index for Ham orbits} focuses on the Robbin--Salamon index of Hamiltonian orbits and provides its computation in several important cases for the Hamiltonians used to compute the relative symplectic cohomology of balls in $\CP^n$, defined in Equation~\eqref{eq: acc. data for balls}. Finally, in Section~\ref{ss: computations of FMB}, we compute the Floer--Morse--Bott indices of these $1$-periodic orbits.

\subsection{Properties of the Robbin--Salamon index}\label{ss: RS part II} In this section, we continue the discussion from Section~\ref{sss: RS part I} regarding the Robbin-Salamon index for paths of symplectic matrices (see Definition~\ref{def: RS index}).

\begin{rem}
    The Robbin--Salamon index satisfies several further properties that will be useful below.
    \begin{itemize}
        \item \textbf{(Naturality)} For every pair of paths $\Phi,\Psi\fc[0,1]\to\Sp(2n)$, we have $\mu_{RS}(\Phi\Psi\Phi^{-1})=\mu_{RS}(\Psi)$;
        \item \textbf{(Product)} Given two paths $\Psi_1\fc[0,1]\to \Sp(2n_1)$ and $\Psi_2\fc[0,1]\to \Sp(2n_2)$ we have $\mu_{RS}(\Psi_1\oplus\Psi_2)=\mu_{RS}(\Psi_1)+\mu_{RS}(\Psi_2)$;
        \item \textbf{(Shear)} Given a symmetric matrix $B$ of size $n\times n$, the Robbin-Salamon index of the symplectic shear $\Psi\fc[0,1]\to \Sp(2n)$, which has the form 
        $$\Psi_t=\begin{pmatrix}
            I&-tB\\0&I
        \end{pmatrix}$$
        for every $t\in[0,1]$, is given by $\mu_{RS}(\Psi)=\frac{1}{2}\sign(B)$,
    \end{itemize}
    See, for instance, \cite[Lemma 56]{Gutt_indexes}.
\end{rem}

Given a smooth path $\Psi \fc [0,1] \to \Sp(2n)$ of symplectic matrices, it is a straightforward verification to show that the path of matrices $S \fc [0,1] \to \Mat(2n,\R)$ defined by
\[
\dot{\Psi}_t = J_0 S_t \Psi_t, \qquad t\in[0,1],
\]
is a path of symmetric matrices. 
Based on this observation, one can provide a computational description of the Robbin--Salamon index. 
The following definition and proposition are dedicated to it; see \cite[pages 6--7]{RS_index}, \cite[Section 3.1]{Oancea_survey} and \cite[Proposition 51]{Gutt_indexes}.

\begin{defin}
   Let $\Psi \fc [0,1] \to \Sp(2n)$ be a smooth path of symplectic matrices.
   \begin{itemize}
       \item A point $t \in [0,1]$ is called a \emph{crossing} if $\det(\Psi_t - \Id) = 0$.
       \item For $t \in [0,1]$, the \emph{crossing form} $\Gamma(\Psi, t)$ is defined as the quadratic form which is the restriction of $S_t$ to $\ker(\Psi_t - \Id)$, where $S_t$ is the symmetric matrix satisfying $\dot{\Psi}_t = J_0 S_t \Psi_t$. Note that in $\ker(\Psi_t - \Id)$, we have $\dot{\Psi}_t = J_0 S_t \Psi_t = J_0 S_t$, and thus \[
       \Gamma(\Psi, t)(v) = \langle S_t v, v \rangle = \omega(S_t v, J_0v) = \omega(-J_0\dot{\Psi}_t v, J_0v) = \omega(v, \dot{\Psi}_tv) = \langle v,-J_0\dot{\Psi}_tv\rangle.
       \]
       \item A crossing $t_0$ is called \emph{regular} if the crossing form $\Gamma(\Psi, t_0)$ is nondegenerate, i.e. $0$ is not an eigenvalue of $\Gamma(\Psi, t_0)$.
   \end{itemize}
\end{defin}

\begin{prop}\label{prop: RS crossing formula}
    Let $\Psi\fc[0,1]\to \Sp(2n)$ be a smooth path of symplectic matrices. If $\Psi$ has only regular crossings then its Robbin--Salamon index satisfies
    \[
    \mu_{RS}(\Psi) = \frac{1}{2} \sign\, \Gamma(\Psi,0) + \sum_{\substack{t \in (0,1) \\ t \text{ crossing}}} \sign\, \Gamma(\Psi,t) + \frac{1}{2} \sign\, \Gamma(\Psi,1). \tag*{\qedsymbol}
    \]
\end{prop}
The following corollary of Proposition~\ref{prop: RS crossing formula} is well known and will be useful for our computations. See \cite[Section 6]{Gutt_indexes} for a similar result.
\begin{coroll}\label{coroll: RS of rotation}
Let $k\in\Z$ and let $\Phi\fc[0,1]\to\Sp(2n)$ be the path given by $\Phi(t)=\exp(\pi k t\, J_0)$ for every $t\in[0,1]$. Then $\mu_{RS}(\Phi)=kn$.
\end{coroll}

The following technical lemma will be used for a computation in Section~\ref{sss: RS on CPn minus divsors}. 
Since we did not find it in the literature, we include a proof for the sake of completeness.

\begin{lemma}\label{lemma: RS of T}
	Let $\Psi\fc[0,1]\to\Sp(2n)$ be a continuous path of symplectic matrices. Then $$\mu_{RS}(\Psi^T)=-\mu_{RS}(\Psi).$$
\end{lemma}

\begin{proof}
	It is enough to prove the lemma for regular paths.
	Let $\Psi\fc[0,1]\to\Sp(2n)$ be a regular path. Since for every $t\in[0,1]$ the matrix $\Psi(t)$ is symplectic, we have $\Psi(t)^TJ_0\Psi(t)=J_0$ where $J_0$ is the standard complex structure on $\R^{2n}$. Therefore for every $t\in[0,1]$ we have
	$$\Psi(t)^TJ_0=J_0\Psi(t)^{-1},\,\,\,\,\, \dot{\Psi}(t)^T J_0\Psi(t)=-\Psi(t)^T J_0\dot{\Psi}(t).$$
	
	Let $t_0\in[0,1]$ be a crossing and let $v\in\ker(I-\Psi(t_0))$. Since $\Psi(t_0)v=v$ we get that $\Psi(t_0)^{-1}v=v$, and thus
	$$\Psi(t_0)^TJ_0v=J_0\Psi(t_0)^{-1}v=J_0v,$$
    which means $(I-\Psi(t_0)^T)J_0v=0$. Therefore $J_0|_{\ker(I-\Psi(t_0))}\fc \ker(I-\Psi(t_0))\to \ker(I-\Psi(t_0)^T)$ is an isomorphism.
	Moreover, for every $v\in\ker(I-\Psi(t_0))$ we have

    \begin{align*}
		\Gamma(\Psi^T,t_0)(J_0v)&=\langle J_0v,-J_0\dot{\Psi}(t_0)^TJ_0v\rangle=\langle v,-\dot{\Psi}(t_0)^TJ_0v\rangle=\langle v,-\dot{\Psi}(t_0)^TJ_0 \Psi(t_0) v\rangle \\&= \langle v,\Psi(t_0)^TJ_0 \dot{\Psi}(t_0) v\rangle = \langle \Psi(t_0)v,J_0 \dot{\Psi}(t_0) v\rangle = \langle v, J_0\dot{\Psi}(t_0)v\rangle \\&=-\Gamma(\Psi,t_0)(v).
	\end{align*}

	Therefore $\mu_{RS}(\Psi^T)=-\mu_{RS}(\Psi)$.
\end{proof}

\subsection{Robbin--Salamon index for Hamiltonian orbits}\label{ss:RS index for Ham orbits}
In this section we compute the Robbin--Salamon index of the $1$-periodic orbits of the Hamiltonians used in the computation of the relative symplectic cohomology of balls in $\CP^n$, as defined in Equation~\eqref{eq: acc. data for balls}. The computation will be carried out with respect to three different choices of trivializations. These calculations will play an essential role in the analysis of Floer differentials and continuation maps.
Throughout this section we fix $\Delta \in (0,1)$ and $\ell \in \Z_{\geq0}$, define $h_\ell\fc[0,1]\to\R$ by $h_\ell(r)=h(\Delta,\ell,r)$ for every $r\in[0,1]$ and focus on the Hamiltonian $H_\ell$ on $\CP^n$ given by
$H_\ell(z) = h_\ell(\mu(z))=h(\Delta,\ell,\mu(z))$, for every $z \in \CP^n$, as in Equation~\eqref{eq: acc. data for balls}.

\subsubsection{The Robbin--Salamon index for orbits in $\CP^n\setminus D_\infty$}\label{sss:RS in Cn}

According to Claim~\ref{claim: Int B(1)->CPn-D_infty}, the complement of $D_\infty=\CP^{n-1}$ in $\CP^n$ is symplectomorphic to the open ball $\Int B(1)\subset\C^n$ equipped with the standard symplectic form $\omega_0$. 

The computations we present here are similar to those in \cite[Section 3.3]{Oancea_survey}, and they will be carried out with respect to the standard trivialization $\tau$ of the tangent bundle of $\Int B(1)\subset\C^n$.

By Remark~\ref{rem: dynamics of a radial Ham} the Hamiltonian vector field of $H_\ell=h_\ell\circ \mu$ over $\Int B(1)$ is given by
$$X_{H_\ell}(z)=-2\pi  h'_\ell(\pi\|z\|^2)J_0z,$$ 
$z\in \Int B(1) $, and its Hamiltonian flow is given by
$$
\varphi_{H_\ell}^t(z)=e^{-2\pi h'_\ell(\pi\|z\|^2)J_0 t}z,
$$
for every $t\in\R$ and $z\in \Int B(1) $, where $J_0$ is the standard complex structure on $\R^{2n}$. Therefore, the linearization of the Hamiltonian flow is
$$
d_z\varphi_{H_\ell}^t(Y)
=
e^{-2\pi h'_\ell(\pi\|z\|^2)J_0t}Y
-
4\pi^2 t\,h''_\ell(\pi\|z\|^2)\langle z,Y\rangle\,
e^{-2\pi  h'_\ell(\mu(z))J_0t}J_0z,
$$
for every $t\in \R$, $z\in \Int B(1)$ and $Y\in T_z\R^{2n}$. In particular, along the constant orbit at $\{0\}$ we have $h''_\ell(0)\langle 0,Y\rangle=0$, and therefore
$$
d_0\varphi_{H_\ell}^t(Y)=e^{-2\pi  h'_\ell(0)J_0t}Y,
$$
for every $t\in\R$ and $Y\in T_0\R^{2n}$. Since $|2\pi h'_\ell(0)|<2\pi$ we deduce that $d_0\varphi_{H_\ell}^t$ has a unique crossing at $t=0$, where $d_0\varphi_{H_\ell}^0=\Id$ and $\ker(d_0\varphi_{H_\ell}^0-\Id_{2n})=\R^{2n}$. A direct computation shows that the intersection form at $t=0$ equals $-2\pi h'_\ell(0)\cdot \Id_{2n}$, whose signature is $-2n$. Since this crossing is regular, by Proposition~\ref{prop: RS crossing formula} we get that
$$
\mu_{RS}^\tau(x;H_\ell)=-n,
$$
where $x\fc[0,1]\to \Int B(1)$ is the constant loop $x(t)=0$ for every $t\in[0,1]$.

According to our computation, a non-constant $1$-periodic orbit $x\fc [0,1]\to \Int B(1)$ of $H_\ell$ has the form $x(t)=e^{-2\pi ktJ_0}z$ for every $t\in[0,1]$, for some $1\leq k\leq \ell$, and for some $z\in \Int B(1)$ satisfying $h'_\ell(\pi\|z\|^2)=k$.

Take $1\leq k\leq \ell$ and denote
$$
S_k=\{\,z\in \Int B(1)\,|\, h'_\ell(\pi\|z\|^2)=k\,\}.
$$
Thus $S_k$ is a sphere.  
Take $z\in S_k$, denote $\xi_z=T_z S_k\cap J_0 T_z S_k$, and note that $\xi_z$ is a symplectic linear subspace of $T_z\R^{2n}$ of codimension $2$.
Consider the symplectic trivialization of $T\R^{2n}$ given by
$$
T_w\R^{2n}
=
\R\langle (2\pi \|z\|^2)^{-1}z\rangle
\oplus
\R\langle  2\pi J_0z\rangle
\oplus
\xi_z,
$$
for every $w\in\R^{2n}$.

Note that for every $t\in[0,1]$ we have
\begin{multline*}
    d_z\varphi_{H_\ell}^t\big((2\pi \|z\|^2)^{-1}z\big)=\\=e^{-2\pi kJ_0t}\big((2\pi \|z\|^2)^{-1}z\big)- 4\pi^2 t\,h''_\ell(\pi\|z\|^2)\left\langle z,\big((2\pi \|z\|^2)^{-1}z\big)\right\rangle\,
e^{-2\pi  kJ_0t}J_0z=\\
=e^{-2\pi kJ_0t}\left((2\pi \|z\|^2)^{-1}z
    -  t\,h''_\ell(\pi\|z\|^2)\,
2\pi J_0z\right).
\end{multline*}
Additionally, for every $Y\in (\R\langle z\rangle)^\bot$ we have
$$
d_z\varphi_{H_\ell}^t(Y)=e^{-2\pi k tJ_0}Y.
$$

Thus we can present $d_z\varphi_{H_\ell}^t$ in our trivialization as $[d_z\varphi_{H_\ell}^t]=\Psi(t)\cdot\chi(t)$, where $\Psi,\chi\fc[0,1]\to\Sp(2n)$ are given by $\Psi(t)=e^{-2\pi k tJ}$ and
$$
\chi(t)=
\begin{pmatrix}
    \begin{pmatrix}
        1&0\\-t\,h''_\ell&1
    \end{pmatrix}&\\&\Id_{2n-2}
\end{pmatrix},$$
for every $t\in[0,1]$.

Let $K\fc[0,1]\times[0,1]\to\Sp(2n)$ be the following homotopy, that connects, with fixed endpoints, the path $\Psi(t)\chi(t)$ to the concatenation of $\Psi(t)$ and $\Psi(1)\chi(t)$, given by
$$
K(s,t)=
\begin{cases}
\Psi\!\left(\dfrac{2t}{s+1}\right)\,\chi(st), & t\le \dfrac{s+1}{2},\\[0.8em]
\Psi(1)\,\chi\big((s+2)t-(s+1)\big), & t\ge \dfrac{s+1}{2},
\end{cases}
$$
for every $(s,t)\in[0,1]\times[0,1]$, as in \cite{CFHW_1996_ApSH_II}.

Since the Robbin--Salamon index is homotopy invariant and additive under concatenation, we deduce that
$$
\mu_{RS}([d_z\varphi_{H_\ell}^t])
=
\mu_{RS}(\Psi)+\mu_{RS}(\Psi(1)\chi(\cdot)).
$$
Additionally, since $\Psi(t)=e^{-2\pi k tJ_0}$ for every $t\in[0,1]$, by Corollary~\ref{coroll: RS of rotation}, we know that $\mu_{RS}(\Psi)=-2nk$. Moreover, since $\Psi(1)=\Id$, we get from the product, zero, and shear properties of the Robbin--Salamon index, together with Lemma~\ref{lemma: RS of T}, that
\begin{align*}
 \mu_{RS}(\Psi(1)\chi(\cdot))&=\mu_{RS}(\chi)  \\
 &=\mu_{RS}\left(\begin{pmatrix}
        1&0\\-t\,h''_\ell&1
    \end{pmatrix}_{t\in[0,1]}\right)+\mu_{RS}(\Id_{2n-2})\\
    &=-\mu_{RS}\left(\begin{pmatrix}
        1&-t\,h''_\ell\\0&1
    \end{pmatrix}_{t\in[0,1]}\right)+0\\
    &=-\frac{1}{2}\sign{h''_\ell(\|z\|^2)}\\
    &=-\frac{1}{2}.
\end{align*}
Therefore, we obtain that
$$
\mu_{RS}^\tau(x;H_\ell)
=
\mu_{RS}([d_z\varphi_{H_\ell}^t])
=
-(2nk+\tfrac{1}{2}).
$$

\subsubsection{RS-index of orbits on $D_\infty = \CP^{n-1}$}\label{sss: RS on CPn-1}

According to Claim~\ref{claim: mu in O(1)} we can identify a neighborhood of $D_\infty \subset \CP^n$ with a neighborhood $U$ of the zero section of $\cO(1)$ over $D_\infty=\CP^{n-1}$. In these coordinates, $H_\ell$ takes the form
$$H_\ell(w,\lambda) = \frac{\ell+\frac{1}{2}}{1+\|\lambda\|_h^2},$$
and the Hamiltonian flow is given by
$$\varphi_{H_\ell}^t(w,\lambda) = (w, e^{2\pi i(\ell+1/2)t}\lambda),$$
for every $(w,\lambda) \in U$ and $t \in \R$. 

Over the zero section, the tangent bundle $T\cO(1)|_{\CP^{n-1}}$ decomposes into the sum of the tangent bundle of the zero section and the vertical bundle:
$$T\cO(1)|_{\CP^{n-1}} = T\CP^{n-1} \oplus \cO(1)|_{\CP^{n-1}}.$$
As a consequence of Remark~\ref{rem: dynamics of a radial Ham}, each point of $D_\infty \subset \CP^n$ is a constant $1$-periodic orbit for $H_\ell$. The linearization of the Hamiltonian flow along such a constant orbit $w_0 \in D_\infty$ can be expressed as
$$d_{w_0}\varphi_{H_\ell}^t(X,Y) = (X, e^{2\pi i (\ell+1/2)t}Y),$$
for every $(X,Y) \in T_{w_0} \CP^{n-1} \oplus \cO(1)_{w_0}$.

By the product and the zero axioms of the Robbin--Salamon index, the index of this path of symplectic matrices reduces to the index of the path restricted to the normal bundle:
$$\Psi(t) = D\varphi^t|_{\cO(1)_{w_0}} = e^{2\pi(\ell+1/2)it} \cdot \id.$$
Note that $\frac{d}{dt}\Psi(t) = 2\pi(\ell+1/2)i e^{2\pi(\ell+1/2)it} \cdot \id$. The crossing form is given by
$$\omega\left( Y, \dot{\Psi}(t) Y \right) = \left\langle Y, 2\pi(\ell+1/2) e^{2\pi(\ell+1/2)it} Y \right\rangle.$$
A time $t \in [0,1]$ is a crossing if and only if $(\ell+1/2)t \in \Z$, which occurs when
$$t \in \left\{0, \frac{1}{\ell+\frac{1}{2}}, \dots, \frac{\ell}{\ell+\frac{1}{2}}\right\}.$$
Since the crossing form is positive definite at each crossing (acting on the $2$-dimensional real fiber), its signature is $2$. Thus, the Robbin--Salamon index is
$$\mu_{RS}(\Psi) = \frac{1}{2}\text{sign}(\Psi(0)) + \sum_{j=1}^\ell \text{sign}(\Psi(t_j)) = \frac{1}{2}(2) + \ell(2) = 2\ell+1.$$

Here the computation was made, essentially, with respect to the trivialization induced by the choice of the constant disk at the constant orbit $(w,\lambda)$.

\subsubsection{RS-index of orbits on $\CP^n\setminus\left(D_1\cup\ldots\cup D_n\cup D_\infty\right)$}\label{sss: RS on CPn minus divsors}

This section focuses on computing the Robbin--Salamon index for $1$-periodic orbits of $H_\epsilon$ that do not intersect the divisors $D_1,\ldots,D_n,D_\infty$. Recall that Claim~\ref{claim: CPn - all divisors} asserts the existence of a symplectomorphism $\Psi\fc (V\times \T^n,\Omega_0)\to (\CP^n\setminus (D_1\cup\ldots\cup D_n\cup D_\infty),\omega_{FS})$, where
$$ V = \left\{(x_1,\ldots,x_n) \in \R^n \, : \, \sum_{j=1}^n x_j < \frac{1}{\pi} \,\,\text{and}\,\, \forall \, j \in \{1,\ldots,n\} \,\, \text{we have} \,\, x_j > 0 \right\}, $$
and $\T^n = (\R/\Z)^n$ is the $n$-torus equipped with angular coordinates $(\theta_1,\ldots,\theta_n)$, and symplectic form $\Omega_0=\pi\sum_{j=1}^n dx_j\wedge d\theta_j$.

Since $\T^n$ is a Lie group, its tangent bundle admits a natural trivialization. This gives rise to a trivialization $\tau$ of $T(V\times \T^n)$, since $V$ is an open subset of $\R^n$. The Robbin--Salamon index in this section will be computed with respect to this trivialization.

Recall, by Remark~\ref{rem: mu in CPn - all divisors}, that the Hamiltonian $\Psi^*H_\ell$ is given by
$$ \Psi^*H_\ell(x,\theta_1,\ldots,\theta_n) = h_\ell(\pi\|x\|_1), $$
and the Hamiltonian vector field is given by
$$ X_{\Psi^*H_\ell}(x,\theta_1,\ldots,\theta_n) = -h'_\ell(\pi\|x\|_1) \sum_{j=1}^n \partial\theta_j, $$
for every $(x,\theta_1,\ldots,\theta_n) \in V \times \T^n$, where $\|x\|_1 = |x_1| + \cdots + |x_n|$ for every $x = (x_1,\ldots,x_n) \in \R^n$.
Therefore, the flow of $X_{\Psi^*H_\ell}$ is a translation in the $\theta$-variables:
$$ \varphi_{\Psi^*H_\ell}^t(x,\theta_1,\ldots,\theta_n) = (x, \theta_1 - th'_\ell(\pi\|x\|_1), \ldots, \theta_n - th'_\ell(\pi\|x\|_1)), $$
for every $(x,\theta_1,\ldots,\theta_n) \in V \times \T^n$ and $t \in \R$.

Thus, a point $(x_0,\theta_0) \in V \times \T^n$ is a fixed point of the time-$1$ map $\varphi_{\Psi^* H_\ell}^1$ if and only if $h'_\ell(\pi\|x_0\|_1) \in \Z$.

Let $(x_0,\theta_0) \in V \times \T^n$ be a fixed point of the time-$1$ map $\varphi_{\Psi^* H_\ell}^1$, and let $x$ be its $1$-periodic orbit. With respect to the trivialization $\tau$, the linearization  of $\{d_{(x_0,\theta_0)}\varphi_{\Psi^* H_\ell}^t\}_{t\in[0,1]}$ along $x$ is represented by the path $\Psi \fc [0,1] \to \Sp(2n)$ given by
$$
\Psi(t) = \begin{pmatrix}
I_n & 0 \\
- t a J_n & I_n
\end{pmatrix},
$$
for $t \in [0,1]$, where $J_n \in \R^{n \times n}$ is the all-ones matrix and $a = \pi h''_\ell(\pi\|x_0\|_1)$. Thus, by Lemma~\ref{lemma: RS of T} and the Shear property of the Robbin--Salamon index, we obtain
$$ \mu_{RS}(\Psi) = -\mu_{RS}(\Psi^T) = -\frac{1}{2}\sign(aJ_n) = -\frac{1}{2}\sign(a). $$
Here we used the fact that $na$ is the only non-zero eigenvalue of $aJ_n$ with multiplicity $1$; hence $\sign(aJ_n) =\sign(na) = \sign(a)$. According to the definition of $h_\ell$, we have $h''_\ell(\pi\|x_0\|_1) > 0$ (see Section~\ref{ss: acc. data}), and thus
$$ \mu_{RS}^\tau(x; \Psi^*H_\ell) = -\frac{1}{2}. $$

\subsection{Computations of Floer--Morse--Bott indices}\label{ss: computations of FMB} 
Recall that the Floer--Morse--Bott index of a critical point $p$ of a Morse function $h$ on a critical submanifold $S$, with respect to an autonomous Hamiltonian $H$ satisfying that \textbf{MB} condition and a symplectic trivialization $\tau$ of $TM$ along $p$, is given by
$$ \mu_{FMB}^\tau(p;H) = \mu_{RS}^\tau(p;H) + \frac{1}{2}\dim M - \frac{1}{2}\dim S + \ind_h p $$
where $\ind_h p$ is the Morse index of $p$ as a critical point of $h$.

\begin{rem}\label{rem: FMB defined modulo} 
Note that with respect to trivializations induced by capping disks of contractible $1$-periodic orbits, the Floer--Morse--Bott index, gives a well-defined number modulo $2N$, independent of the trivialization, where $N$ is the minimal first Chern number of $(M,\omega)$. We denote this index by $\mu_{FMB}$. The modulo $2N$ ambiguity arises because changing the capping by attaching a sphere results in the index being shifted by twice the value of the first Chern class evaluated on that sphere.
\end{rem}

We now return to the Hamiltonians $H_\ell$ for $\ell \in \Z_{\geq 0}$ and $\Delta \in (0,1)$. We go over $1$-periodic orbits corresponding to orbits which are the critical points of Morse functions chosen on the families of orbits for which we computed the Robbin--Salamon index, and calculate the Floer--Morse--Bott index.
  \subsubsection*{$\mathbf{\CP^n\setminus D_\infty}$:}
    In $\CP^n \setminus D_\infty$ there are two types of $1$-periodic orbits for $H_\ell$: the minimum of $H_\ell$ and the non-constant periodic orbits. Their Robbin--Salamon indices, with respect to the trivialization $\tau_B$ over $\CP^n \setminus D_\infty$, discussed in Section~\ref{sss: CP^n - D_infty}, were computed in Section~\ref{sss:RS in Cn}:
    
    \begin{itemize}
        \item The minimum $x$ of $H_\ell$ is a constant and isolated $1$-periodic orbit; thus, the critical submanifold $S$ containing it is a singleton. For this reason, $\dim S=0$ and $\ind_h x = 0$, where $h \fc S \to \R$ is a Morse function. Thus,
        $$ \mu_{FMB}^\tau(x;H_\ell) = -n + \frac{1}{2} \cdot 2n - \frac{1}{2} \cdot 0 + 0 = 0. $$

        \item If $x$ is a non-constant orbit of $H_\ell$ lying in the $(2n-1)$-dimensional sphere 
        $$ S_k = \{z \in \CP^n \colon h'_\ell(\mu(z)) = k\} $$
        for some $1 \leq k \leq \ell$, we showed that $\mu_{RS}^\tau(x;H_\ell) = -(2nk + \frac{1}{2})$. Therefore,
        $$ \mu_{FMB}^\tau(x;H_\ell) = -(2nk + \frac{1}{2}) + \frac{1}{2} \cdot 2n - \frac{1}{2}(2n - 1) + \ind_h x = -2nk + \ind_h x, $$
        where $h$ is a perfect Morse function on the sphere $S_k$. Such a function has exactly two critical points (the minimum and the maximum) with Morse indices $0$ and $2n-1$, respectively. Thus, $\ind_h x \in \{0, 2n-1\}$.

        Thus, modulo $2(n+1)$, which is twice the minimal first Chern number, we get that 
        $$\mu_{FMB}(\hat{x}_k^\ell;H_\ell)=-2nk+2n-1=-2(n+1)(k-1)+2k-3\equiv 2k-3 \mod 2(n+1),$$
        for every $1\leq k\leq \ell$, and 
        $$\mu_{FMB}(\check{x}_k^\ell;H_\ell)=-2nk=-2(n+1)k+2k\equiv 2k \mod 2(n+1),$$
        for every $0\leq k\leq \ell$.
    \end{itemize}

    \subsubsection*{$\mathbf{\CP^n\setminus \lbrace 0 \rbrace}$:}
     In $\CP^n \setminus \{0\}$ we considered the 
     orbits of $H_\ell$ that are constant and lie in the divisor $D_\infty = \CP^{n-1}$.  Their Robbin--Salamon indices with respect to the constant capping were computed in Section~\ref{sss: RS on CPn-1}, using the coordinates we computed in Section~\ref{sss: CP^n - 0}.
     
     If $h$ is a perfect Morse function on $\CP^{n-1}$ and $w$ is a critical point of $h$, we proved that $\mu_{RS}^\tau(w;H_\ell) = 2\ell + 1$. Therefore,
    $$ \mu_{FMB}^\tau(w;H_\ell) = 2\ell + 1 + \frac{1}{2} \cdot 2n - \frac{1}{2}(2n - 2) + \ind_h w = 2\ell + 2 + \ind_h w, $$
    where $\ind_h w \in \{0, 2, \dots, 2n-2\}$.

In particular, all the possibilities for $\mu^\tau_{FMB}(w;H_\ell)$ are even numbers, thus $\mu_{FMB}(w;H_\ell)$ is even modulo $2(n+1)$.

\subsubsection*{$\mathbf{\CP^n \setminus (D_1 \cup \dots \cup D_n \cup D_\infty)}$:}  
In $\CP^n \setminus (D_1 \cup \dots \cup D_n \cup D_\infty)$, the $1$-periodic orbits of $H_\ell$ lie on the spheres $S_1, \dots, S_\ell$ as before, and the perfect Morse function $h$ can be chosen such that the critical points correspond to orbits that are disjoint from the divisors. 
In Section~\ref{sss: RS on CPn minus divsors}, we computed Robbin--Salamon indices with respect to the trivialization $\tau_T$, discussed in Section~\ref{sss: CPn - all divisors}. Any such orbit $x$ has $\mu_{RS}^\tau(x;H_\ell) = -\frac{1}{2}$. Therefore,
    $$ \mu_{FMB}^\tau(x;H_\ell) = -\frac{1}{2} + \frac{1}{2} \cdot 2n - \frac{1}{2}(2n - 1) + \ind_h x = \ind_h x, $$
    where $h$ is a perfect Morse function on the $(2n-1)$-sphere, yielding $\ind_h x \in \{0,2n-1\}$ 

    \section{Obstructions}\label{s: obstructions}

\subsection{Relative first Chern number and intersection number}\label{ss: c_1 and intersection}

\subsubsection*{Relative first Chern number for complex bundles}
Let $\Sigma$ be a closed surface, with $\Gamma\subset \Sigma$ a finite set of points. An open neighborhood $U$ of $\Gamma$ is called \emph{disk-like} if its closure $\overline{U}$ has a finite number of connected components, each of them is diffeomorphic to a smooth closed disk and contains a single point from $\Gamma$.
Let $E\to \Sigma\setminus \Gamma$ be a complex vector bundle, of complex rank $n$. Since $U \setminus \Gamma$ has a homotopy type of a finite collection of circles, $E\vert^{}_{U\setminus\Gamma}$ is trivial, hence $\det E\vert^{}_{U\setminus\Gamma} := \Lambda^n E\vert^{}_{U\setminus\Gamma}$ is trivial. Let $\tau$ be some complex trivialization of $\det E\vert^{}_{U\setminus\Gamma} := \Lambda^n E\vert^{}_{U\setminus\Gamma}$, thus $\tau \fc U\setminus \Gamma \times \C \to \det E\vert^{}_{U\setminus \Gamma}$ is an isomorphism of complex vector bundles. The relative first Chern number $c^\tau_1(E)$ of $E$ with respect to the trivialization $\tau$ is defined to be the algebraic signed count of zeroes of a generic section of $\det E\vert^{}_{\Sigma\setminus\Gamma}$, extending the locally constant section $\tau(\cdot,1)$ of $E\vert^{}_{U\setminus \Gamma}$. 

One can define an equivalence relation on the set of pairs
$(U, \tau)$ of disk-like neighborhoods of $\Gamma$ and trivializations $\tau$ of $E\vert^{}_{U\setminus \Gamma}$ by defining
$(U, \tau) \sim (U^\prime, \tau^\prime )$ provided that there exists a disk-like neighborhood $V \subseteq U \cap U^\prime$ such that 
that $\tau\vert^{}_{V\setminus \Gamma}$ is homotopic to $\tau^\prime \vert^{}_{V\setminus \Gamma}$. 
The relative first Chern number $c^\tau_1(E)$ 
only depends on the equivalence class of the pair $(U, \tau)$.

Note that a trivialization of the complex vector bundle $E$ induces a trivialization of $\det E$ by the $n$-th exterior product map, where $n$ is the rank of $E$. Therefore we can define a relative first Chern class with respect to a trivialization of $E$ itself.

\subsubsection*{Relative first Chern number for complex real bundles}
Let $E$ be a real vector bundle over a base $B$ of rank $2n$ for some $n\in \N$. A complex structure $J$ on $E$ is an automorphism of $E$ satisfying $J^2_x=-\id$ on the fiber $E_x$ for every $x\in B$. We call the pair $(E,J)$ a \textbf{complex real bundle}, to distinguish it from a complex bundle where the fibers are vector spaces over $\C$. A \textbf{complex real trivialization} of $(E,J)$ is a trivialization $\tau\fc B\times \R^{2n}\to E$ satisfying $\tau_*J_0=J$, where $J_0$ is the standard complex structure on $\R^{2n}$.

Next, we show how a rank $2n$ complex real bundle $(E,J)$ induces a rank $n$ complex vector bundle $E_\C$.
Consider the complexification $E^\C := E\otimes_\R\C$. This is a complex vector bundle of rank $2n$. Extending $J$ linearly over $\C$, $J$ acts on $E^\C$, and satisfies $J^2=-\id$. Thus, its eigenvalues are $i$ and $-i$, and $E^\C$ splits as the sum of two eigenbundles $E^\C = E^{1,0} \oplus E^{0,1}$, where $E^{1,0}=\ker (J-i)$ is the eigenbundle for $i$, and $E^{0,1}=\ker (J+i)$ is the eigenbundle for $-i$. We consider $E^{1,0}$ as the complex bundle associated to the complex real bundle $(E,J)$.

A complex real trivialization of $(E,J)$, given in terms of a $J$-adapted frame $e_1,\ldots,e_n,$ $Je_1$, $\ldots,Je_n$, induces a complex trivialization of $E^{1,0}$ by the frame $\frac{1}{2}(e_1 - iJe_1),\ldots, \frac{1}{2}(e_n - iJe_n)$. This is the projection of $e_1,\ldots,e_n$ on $E^{1,0}$ along $E^{0,1}$ with respect to the decomposition $E^{\C}=E^{1,0}\oplus E^{0,1}$.

Thus, given a complex real trivialization of $(E,J)$, the relative first Chern number is defined with respect to $\det E^{1,0}$ with its induced trivialization.

\subsubsection*{Relative first Chern number for symplectic bundles}

A \textbf{unitary bundle} $(E,J,\omega)$ is a triplet where $E$ is a $2n$ rank bundle over a base $B$, $(E,J)$ is a complex real bundle, $(E,\omega)$ is a symplectic bundle and $J$  is compatible with $\omega$, that is, $\omega(\cdot,J\cdot)$ is a Riemannian metric on $E$. A \textbf{unitary trivialization} of $(E,J,\omega)$ is a trivialization $\tau\fc B\times \R^{2n}\to E$ of $E$ satisfying $\tau^*J=J_0$ and $\tau^*\omega=\omega_0$.

Our next goal is to define a relative first Chern number for symplectic vector bundles. To achieve this, we should prove the following two results:
\begin{lemma}\label{lemma: symp triv homotopic to a unitary triv}
    Let $(E,J,\omega)$ be a unitary bundle over $B$. For every symplectic trivialization $\tau_\omega$ of $(E,\omega)$, there is a unitary trivialization $\tau$ of $(E,J,\omega)$, homotopic to $\tau_\omega$ through a path of symplectic trivializations.
\end{lemma}

\begin{proof}
    Let $\tau_\omega$ be a symplectic trivialization of $(E,\omega)$. The trivialization $\tau_\omega$ defines the complex real bundle $(\R^{2n}\times B,\tau_\omega^*J,\omega_0)$ over $B$. By the proof of \cite[Lemma 2.5.5]{McDuff_Salamon_intro_2016} for every $x\in B$ the endomorphism $P_x=-J_0\tau_\omega^*J_x$ of $\R^{2n}$ is symmetric positive definite and symplectic. Since $P_x$ is symmetric positive definite, it has a unique symmetric positive square root, denoted $A_x$, and $A_x$ depends smoothly on $x$. Moreover, for every $x\in B$, by \cite[Lemma 2.2.3]{McDuff_Salamon_intro_2016}, the endomorphism $A_x$ is symplectic. Define a trivialization $\tau$ of $E$ by $\tau_x=(\tau_\omega)_x\circ A_x^{-1}$ for every $x\in B$. Note that for every $x\in B$ we have
    \begin{align*}
        \tau_x^*J&=\tau_x^{-1}\circ J\circ \tau_x\\
        &=((\tau_\omega)_x\circ A_x^{-1})^{-1}\circ J\circ (\tau_\omega)_x\circ A_x^{-1}\\
        &=A_x\circ (\tau_\omega)_x^{-1}\circ J\circ (\tau_\omega)_x\circ A_x^{-1}\\
        &=A_x\circ (\tau_\omega^* J)_x\circ A_x^{-1}\\
        &=A_xJ_0P_xA_x^{-1}\\
        &=A_x J_0 A_x^2A_x^{-1}\\
        &=A_xJ_0A_x\\
        &=(A_x)^TJ_0A_x\\
        &=J_0.
    \end{align*}

    Again as a conclusion of \cite[Lemma 2.2.3]{McDuff_Salamon_intro_2016}, the space of positive definite symmetric symplectic matrices is contractible, therefore we deduce that $\tau$ and $\tau_\omega$ are homotopic through a path of symplectic trivializations.
    
\end{proof}

\begin{lemma}\label{lemma: two U triv that symp htpy are U htpy}
    Let $(E,J,\omega)$ be a unitary bundle over $B$ and let $\tau_0,\tau_1$ be a pair of unitary trivializations of $(E,J,\omega)$. If they are homotopic to one another through a path of symplectic trivializations then they are homotopic to one another through a path of unitary trivializations.
\end{lemma}
\begin{proof}
    Assume that there is a homotopy $(\tau_s)_{[0,1]}$ between $\tau_0$ and $\tau_1$ consisting of symplectic trivializations. Thus, $(\tau_0^{-1}\circ \tau_s)_{s\in [0,1]}$ is a homotopy of symplectic trivialization of the trivial bundle $B\times (\R^{2n},\omega_0)\to B$. Moreover, this homotopy connects the identity trivelization with $\tau_0^{-1}\circ \tau_1$, which is unitary. Since the inclusion of $\U(n)$ into $\Sp(2n)$ is a homotopy equivalence \cite[Lemma 2.2.4]{McDuff_Salamon_intro_2016}, the homotopy $(\tau^{-1}_0\circ \tau_s)_{s\in [0,1]}$ is homotopic with fixed endpoints to a homotopy $(\tau'_s)_{s\in [0,1]}$ of unitary trivializations on the trivial unitary bundle  $B\times (\R^{2n},J_0,\omega_0)\to B$. Now consider the homotopy $(\tilde\tau_s)_{s\in [0,1]}$ of the unitary bundle $(E,J,\omega)$ given by $\tilde\tau_s=\tau_0\circ\tau'_s$ for every $s\in[0,1]$. This homotopy consists of unitary trivializations and connects $\tilde\tau_0=\tau_0\circ\id=\tau_0$ and $$\tilde\tau_1=\tau_0\circ \tau'_1=\tau_0\circ\tau_0^{-1}\circ\tau_1=\tau_1,$$
    as required.
\end{proof}

As before, let $\Sigma$ be a closed surface, with $\Gamma\subset \Sigma$ a finite set of points and let $U$ be  a disk-like neighborhood of $\Gamma$. 
Let $(E,\omega)$ be a symplectic vector bundle of real rank $2n$ over $\Sigma$, and let $\tau_\omega\fc 
(U\setminus\Gamma)\times (\R^{2n},\omega_0)\to (E|_{U\setminus\Gamma},\omega)$ be a symplectic trivialization of $(E|_{U\setminus\Gamma},\omega)$ along $U\setminus \Gamma$. We now explain that $\tau_\omega$ determines a well-defined relative first Chern number.

First, given an $\omega$-compatible complex structure $J$ of $E$, by Lemma~\ref{lemma: symp triv homotopic to a unitary triv} the trivialization $\tau_\omega$ is homotopic by symplectic trivializations to a unitary trivialization $\tau$ of $(E,J,\omega)$. Define the relative first Chern number $c^{\tau_\omega}_1(E,\omega)$ of $(E,\omega)$ with respect to the trivialization $\tau_\omega$ to be the relative first Chern number $c^{\tau}_1(E,J)$ of $(E,J)$ with respect to the trivialization $\tau$. 

We claim that $c^{\tau_\omega}_1(E,\omega)$ is well-defined. First of all, if there is another unitary trivialization $\tilde\tau$ of $(E,J,\omega)$ that is homotopic to $\tau_\omega$ via symplectic trivializations, then $\tau$ and $\tilde\tau$ are homotopic via symplectic trivializations. Also, by Lemma~\ref{lemma: two U triv that symp htpy are U htpy} there is a homotopy via unitary trivializations that connects $\tau$ to $\tilde{\tau}$, and hence $c^{\tau}_1(E,J)=c^{\tilde\tau}_1(E,J)$. This means that given $J$, relative first Chern number $c^{\tau_\omega}_1(E,\omega)$ does not depend on the trivialization $\tau$.

Now, let $J'$ be another $\omega$-compatible complex structure of $E$. As before, there is a unitary trivialization $\tau'$ of $(E,J',\omega)$ that is homotopic to $\tau_\omega$ via symplectic trivializations. Since the space of $\omega$-compatible complex structures on $(E,\omega)$ is contractible, there is a homotopy $(J_s)_{s\in [0,1]}$ of $\omega$-compatible complex structures from $J$ to $J'$. Since relative Chern number is  homotopy invariant, we conclude that $c^{\tau}_1(E,J)=c^{\tau'}_1(E,J')$, and hence $c^{\tau_\omega}_1(E,\omega)$ is well-defined.

Indeed, let $\pi\fc\Sigma\times[0,1]\to \Sigma$ be the projection, and consider the pullback bundle $\pi^*E$ equipped with the complex structure $\cJ=(J_s)_{s\in[0,1]}$, which is compatible with $\pi^*\omega$. Note that $\pi^*\tau_\omega$ is a symplectic trivialization of $(\pi^*E,\pi^*\omega)$ over $U\setminus \Gamma$, and hence there is a unitary trivialization $\tilde\tau$ of $(\pi^*E,\cJ,\pi^*\omega)$ that is homotopic to $\pi^*\tau_\omega$ via symplectic trivializations. Let $\sigma$ be a section of $\det(\pi^*E,\cJ)$ which equals $1$ over $U\setminus\Gamma$ with respect to $\tilde\tau$. For $i\in\{0,1\}$, let $\tau_i=\tilde\tau|_{(U\setminus\Gamma)\times\{i\}}$ be the restricted trivialization over $(U\setminus\Gamma)\times\{i\}$, and let $Z_i$ be the intersection of $\sigma|_{(U\setminus\Gamma)\times\{i\}}$ with the zero section of $\det(E,J_i)$. Since $Z_0$ and $Z_1$ are oriented $0$-dimensional manifolds and cobordant to each other, their signed count is the same, which means that $c_1^{\tau_0}(E,J_0)=c_1^{\tau_1}(E,J_1)$. Since the relative first Chern number with respect to a symplectic trivialization of a complex bundle is well-defined, we get that
$$c_1^\tau(E,J)=c_1^{\tau_0}(E,J)=c_1^{\tau_1}(E,J')=c_1^{\tau'}(E,J'),$$
as required.

\subsubsection*{Relative first Chern number and divisors}

Let $(M,\omega)$ be a symplectic manifold. Given a smooth symplectic divisor $D$, we denote by $\cO(D)$ the complex line bundle associated to $D$. This complex line bundle is defined up to isomorphism of complex line bundles as a line bundle over $M$, admitting a section vanishing transversally precisely at $D$.  In particular $c_1(\cO(D)) = \PD(D)$. We briefly recall a construction below:

 Since $D$ is a closed symplectic submanifold, there exists a compatible almost complex structure $J$ on $M$ such that $D$ is $J$-holomorphic. Let $\cN_D := TM\vert^{}_D / TD$ be the normal bundle of $D$, the quotient construction equips it with an almost complex structure $\tilde J$ descending from $J$. There exists a tubular neighborhood $U$ of $D$, which can be identified with a neighborhood of the zero section, $V\subset\cN_D$. Denote by $\Phi\fc U \to V$ the identification, and by $\pi \fc\cN_D \to D$ the bundle projection. Consider the tautological section of $\pi^*\cN_D$ whose base is $\cN_D$ and for every $v=(p,q)\in \cN_D$ where $\pi(v)=q$, the fiber $(\pi^*\cN_D)_v$ of $\pi^*\cN_D$ over $v$ equals the fiber $(\cN_D)_q$ of $\cN_D$ over $q$. Now, consider the tautological section $s\fc V\to \pi^*\cN_D$ given by $s(v)=v\in (\cN_D)_{\pi(v)}$ for every $v\in V$. This section vanishes transversely on $D$ and away from $D$ it is nonzero. Finally, consider the pullback line bundle $\Phi^*(\pi^*\cN_D)$ over $U\subset M$, and its section $\Phi^* s$. Since $\Phi^*s$ does not vanish along $U\setminus D$, it defines a trivialization of $\Phi^*(\pi^*\cN_D)|_{U\setminus D}$ over $U\setminus D$. Consider the trivial line bundle $(M\setminus D)\times \C\to M\setminus D$ over $M\setminus D$ and glue it with the line bundle $\Phi^*(\pi^*\cN_D)$ using the trivialization over $U\setminus D$, which is determined by $\Phi^*s$, to obtain a line bundle over $M$ that coincides with $\Phi^*(\pi^*\cN_D)$ over $U$ and admits a section that vanishes exactly along $D$ and is transversal to $D$. The obtained line bundle is $\cO(D)$.

Let $m\geq 1$ and let $D_1,\ldots,D_m$ be pairwise transversely intersecting smooth divisors. Denote $D=D_1\cup\ldots\cup D_m$. We define $\cO(D)$ by $\cO(D_1)\otimes \ldots \otimes \cO(D_m)$.

Assume that $D$ is Poincar\'e dual to $\frac{1}{N}c_1(TM)$, namely that $c_1(TM) = N\cdot \PD(D)$, for some $N\in \N$. Since complex line bundles are classified up to bundle isomorphism by their first Chern classes, there is a bundle isomorphism $\det TM \cong \cO(D)^{\otimes N}$. As seen by their construction, the bundles $\cO(D_i)$ admit a section vanishing transversely on $D_i$, for every $1\leq i\leq m$; we call such a section a canonical section $\sigma_{D_i} \fc M \to \cO(D_i)$. The tensor product $\sigma_D=\sigma_{D_1}\otimes\cdots\otimes\sigma_{D_m}\fc M\to \cO(D)$ is called a canonical section of $\cO(D)$. The identification $\det TM \cong \cO(D)^{\otimes N}$ defines a section of $\det TM$ which we denote by $\sigma_{N\cdot D}$ which is equal to $\sigma_D^{\otimes N}$, under this identification. 

Note that the section $\sigma_{N\cdot D}$ vanishes to order $N$ along $D$.

The section $\sigma_{N\cdot D}$ is an example of a $D$-canonical section: A section $\sigma$ of $\det TM$ is called \textbf{canonical with respect to $D$}, or simply $D$-canonical, if $\sigma^{-1}(0)=D$ and there exists $N\in \N$ such that for every $x\in D$ there exists a neighborhood $U_x$ such that $\sigma= (g_1\cdot\ldots\cdot g_m)^N\cdot f\cdot e$ where:
\begin{enumerate}
    \item $e$ is a local section of $\det TM$, nowhere vanishing in $U_x$;
    \item $f \fc U_x \to \C$ is a nonzero function;
    \item for every $1\leq i\leq m$, the function $g_i \fc U_x \to \C$ has simple zeros along $D_i\cap U_x$. 
\end{enumerate}

Let $V\subset M$ be an open subset, let $J$ be a compatible almost complex on $M$ and let $\tau$ be a symplectic trivialization of $\det TM$ over $V$. Let $u\fc \Sigma\to M$ be a smooth map from a Riemann surface with punctures $\Sigma$ such that all the asymptotes of $u$ are contained in $V$. Then the \textbf{relative first Chern number} of $u$ with respect to $\tau$ is denoted $c_1^\tau(u)$ and is defined to be $c_1^{u^*\tau}(u^*\det TM|_V)$.

We now formulate and prove a relation between relative Chern numbers and intersection numbers, when the trivialization is taken to be via such canonical sections.

\begin{lemma}\label{lemma: c_1 wrt D canonical section and intersection number}
    Let $(M,\omega)$ be a closed symplectic manifold, let $J$ be a compatible almost complex structure and let $D$ be a symplectic divisor (possibly a union of transversely intersecting smooth divisors). Assume that
    \begin{itemize}
        \item There exists $N\in \N$ such that $c_1(TM)=N\cdot\PD(D)$;
        \item There is a trivialization $\tau\fc (M\setminus D)\times \C\to \det TM$ of $\det TM$ along $M\setminus D$;
        \item There exists a $D$-canonical section $s\fc M\to \det TM$ such that $\tau(\cdot,1)=s$.
    \end{itemize}
    
    Let $u\fc \Sigma\to M$ be a smooth map from a Riemann surface with punctures $\Sigma$ such that all the asymptotes of $u$ are contained in $M\setminus D$. Then the relative first Chern number of $u$ equals to $N$ times the intersection number of $u$ with $D$, i.e. $c_1^\tau(u)=N(u\cdot D)$.
\end{lemma}
\begin{proof}
Without loss of generality, we may assume that $u$ is transverse to $D$. If $D$ is a union of smooth divisors intersecting transversely, then the intersection of any two or more components has complex codimension at least $2$. Hence, after a perturbation, $u$ may be assumed to intersect only the smooth locus of $D$. This preserves $c_1^\tau$, as transversality can be achieved by a homotopy of $u$ supported near the intersection locus, and $c_1^\tau$ is invariant under such homotopies.

    The trivialization $u^*\tau$ of $(u^*\det TM)|_{u^{-1}(M\setminus D)}$ is determined by the section $\sigma = u^*s$, which is a section of $u^*\det TM$. We may compute $c_1^\tau(u)$ by a signed count with multiplicities of the zeroes of this section $\sigma$. 

    Note that for every $x\in \Sigma$, $\sigma(x)=0$ if and only if $u(x)\in D$. Since $s$ vanishes to order $N$ along $D$, and since $u$ is transversal to $D$, every zero of $\sigma$ is of order $N$. This shows that the relative first Chern number $c^\tau_1(u)$ equals $N$ times the intersection number of $u$ with $D$, i.e. $c^\tau_1(u) = N(u\cdot D)$, as required.

\end{proof}

\subsubsection*{Two useful trivializations}

Recall that in \ref{sss: CPn - all divisors} we have shown a symplectomorphism $\Psi:U \to \CP^n\setminus D$, where $U=V\times \T^n$, for $V$ the interior of the moment simplex, and $D=D_1\cup\ldots\cup D_n\cup D_\infty$ is the union of the toric divisors.

In $U=V\times \T^n$, let us denote by $x_1,\ldots, x_n$ the coordinates in $V$ and by $\theta_1,\ldots,\theta_n$ the cyclic coordinates on $\T^n$.
We define the trivialization of $TU$ given by the frame $(\partial_{x_1},\ldots\partial_{x_n},\partial_{\theta_1},\ldots\partial_{\theta_n})$. It is a symplectic trivializtion with respect to $\Omega_0$, and compatible with the $n$-torus action.
We then use the map $\Psi$ to define a trivializtion of $TM\big\vert_{\CP^n\setminus D}$ by pushing the frame from $U$, namely by the frame $(\Psi_*\partial_{x_1},\ldots \Psi_*\partial_{x_n},\Psi_*\partial_{\theta_1},\ldots \Psi_*\partial_{\theta_n})$, this is a symplectic frame with respect to $\omega_{FS}$ due to $\Psi$ being a symplectomorphism. Let us denote by $\tau_T$ (T stands for the \emph{torus}) the trivialization of $\det TM\vert^{}_{\CP^n\setminus D}$ induced by this frame.\label{p: def tau_T}
\begin{prop}
    Let $u\colon \Sigma\to M$ be a smooth map from a Riemann surface with punctures $\Sigma$ such that all the asymptotes of $u$ are contained in $\CP^n\setminus D$. Then for $\tau = \tau_T$, \[c_1^{\tau}(u)=u\cdot D.\]
\end{prop}
\begin{proof}
By Lemmata \ref{lemma: symp triv homotopic to a unitary triv} and \ref{lemma: two U triv that symp htpy are U htpy}, $c_1^\tau(u)$ is independent of the compatible almost complex structure on $TM$, thus let us first choose a convenient compatible almost complex structure $\tilde {J}_0$ on $\CP^n$ in order to compute $c_1^\tau(u)$.

On $TU$ consider the frame $(\partial_{x_1},\ldots\partial_{x_n},\partial_{\theta_1},\ldots\partial_{\theta_n})$ and define an almost complex structure $J_0$ on $TU$ to be the unique almost complex structure on $TU$ that satisfies $J_0\partial_{x_j}=\partial_{\theta_j}$ for all $1\le j \le n$. Since this frame is also symplectic with respect to $\Omega_0$, we deduce that $J_0$ is compatible with $\Omega_0$. Define the almost complex structure $\tilde{J}_0$ on $\CP^n$ as follows: In a compact neighborhood of the asymptotes of $u$ put $\tilde{J}_0=\Psi_* J_0 \Psi_*^{-1}$, and then extend it arbitrarily to all $\CP^n$ to a compatible almost complex structure. Thus, $c_1^\tau(u)$ is defined as $ c_1^\tau(\det u^*(T\CP^n,\widetilde{J}_0))$.

On $\CP^n\setminus D$ the induced trivialization by $\tau$ of the complex bundle $(TM,\tilde{J}_0)$ is determined by the complex The frame $(\Psi_*\partial_{x_1},\ldots \Psi_*\partial_{x_n},\widetilde{J}_0\Psi_*\partial_{x_1},\ldots \widetilde{J}_0\Psi_*\partial_{x_n})$.  Let us use homotopy of vector bundles and trivialization in order to make computations of relative first Chern number more convenient.

We start with a homotopy $(R_t)_{t\in[0,\pi/2]}$ of trivializations, given by rotation in the fibers of $TU$. This homotopy satisfies $R_t(\partial_{x_j}) = \cos(t)\partial_{x_j} + \sin(t)J_0 \partial_{x_j}$, for $0\le t \le \pi/2$  and all $1\le j \le n$. Precomposing it with $\Psi$ yields a homotopy of complex frames from 
$(\Psi_*\partial_{x_1},\ldots \Psi_*\partial_{x_n},\widetilde{J}_0\Psi_*\partial_{x_1},\ldots \widetilde{J}_0\Psi_*\partial_{x_n})$ 
to $(\Psi_*\partial_{\theta_1},\ldots \Psi_*\partial_{\theta_n},\widetilde{J}_0\Psi_*\partial_{\theta_1},\ldots\widetilde{J}_0\Psi_*\partial_{\theta_n})$.
Next, since both $\widetilde{J}_0$ and $J_{std}$ are $\omega_{FS}$ compatible, they are connected by a homotopy of compatible complex structures, $(\widetilde{J}_s)_{s\in[0,1]}$, such that $\widetilde{J}_1 = J_{std}$.
Thus we get a homotopy of real complex vector bundles $((T\CP^n,\widetilde{J}_s))_{s\in[0,1]}$ such that for every $s\in [0,1]$ the vector bundle $(T\CP^n,\widetilde{J}_s)$ admits a trivialization $\tilde\tau_s$ along $\CP^n\setminus D$ determined by the frame $(\Psi_*\partial_{\theta_1},\ldots \Psi_*\partial_{\theta_n},\widetilde{J}_s\Psi_*\partial_{\theta_1},\ldots\widetilde{J}_s\Psi_*\partial_{\theta_n})$.

For every $s\in [0,1]$ consider the line bundle $\mathcal{L}_s := \det(T\CP^n,\widetilde{J}_s)$, over $\CP^n$, and the trivialization $\tau_s$  over $\CP^n \setminus D$ which is induced by the trivialization $\tilde\tau_s$ of the vector bundle $(T\CP^n,\widetilde{J}_s)$, namely given by the section:
\[(\Psi_*\partial_{\theta_1}-i\widetilde{J}_s\Psi_*\partial_{\theta_1})\wedge\ldots \wedge(\Psi_*\partial_{\theta_n}-i\widetilde{J}_s\Psi_*\partial_{\theta_n}).\]
Thus, by homotopy invariance of the relative first Chern number, $c^{\tau_s}_1(u^*\mathcal{L}_s)$, does not depend on $s$, in particular $c^{\tau_0}_1(u^*\mathcal{L}_0) = c^{\tau_1}_1(u^*\mathcal{L}_1)$, and since  $c^{\tau}_1(u^*\mathcal{L}_0) =  c^{\tau_0}_1(u^*\mathcal{L}_0)$, we get $c^{\tau}_1(u^*\mathcal{L}_0) =  c^{\tau_1}_1(u^*\mathcal{L}_1)$. We now aim to show that the section defining the trivialization $\tau_1$ of $\mathcal{L}_1 = (T\CP^n,J_{std})$, is a $D$-canonical section, hence the conclusion of Lemma \ref{lemma: c_1 wrt D canonical section and intersection number} holds.

Let us compute $\Psi_* \partial_{\theta_j}$.
$\Psi$ lifts to a map $\tilde{\Psi}\colon U \to \C^{n+1}\setminus\{0\}$, 
\[ \tilde\Psi(x,\theta) = \left( \sqrt{\frac{1}{\pi} - \sum_{k=1}^n x_k} , \sqrt{x_1} e^{2\pi i \theta_1} , \dots , \sqrt{x_n} e^{2\pi i \theta_n} \right), \]
for every $(x,\theta)=(x_1,\ldots,x_n,\theta_1,\ldots,\theta_n)\in U$, where $\Psi = p\circ\tilde{\Psi}$, with $p\colon \C^{n+1}\setminus\{0\} \to \CP^n$ the Hopf projection. By writing $\tilde\Psi=(z_0,\ldots,z_n)$ we get that for every $1\leq j\leq n$ we have $\frac{\partial z_0}{\partial \theta_j}=0$ and for every $1\leq j,k\leq n$ we have
$$\frac{\partial z_k}{\partial \theta_j} (x,\theta)= \frac{\partial}{\partial \theta_j}\left(\sqrt{x_k} e^{2\pi i \theta_k}\right) = \delta_{kj}2\pi i\sqrt{x_k} e^{2\pi i \theta_k} = \delta_{kj}2\pi i z_k(x,\theta),$$

for every $(x,\theta)\in U$, where $\delta_{kj}$ is Kronecker's delta. Similarly  for every $1\leq j\leq n$ we have $\frac{\partial \overline z_0}{\partial \theta_j}=0$ and for every $1\leq j,k\leq n$ we have
$$\frac{\partial \overline z_k}{\partial \theta_j} (x,\theta)= \frac{\partial}{\partial \theta_j}\left(\sqrt{x_k} e^{-2\pi i \theta_k}\right) = -\delta_{kj}2\pi i\sqrt{x_k} e^{-2\pi i \theta_k} = -\delta_{kj}2\pi i \overline z_k(x,\theta),$$
for every $(x,\theta)\in U$

And therefore for every $1\leq j\leq n$ we get
\[\tilde{\Psi}_*\partial_{\theta_j} = 2\pi i z_j\partial_{z_j} - 2\pi i \bar{z}_j\partial_{\bar{z}_j},\]
and hence, since $J_{std} \partial_{z_j}=i \partial_{z_j}$, and $J_{std} \partial_{\bar{z}_j}=-i\partial_{\bar{z}_j}$, we get that
\[\tilde{\Psi}_*\partial_{\theta_j}-iJ_{std}\tilde{\Psi}_*\partial_{\theta_j} = 4\pi i z_j\partial_{z_j}. \]

Thus the trivialization $\tau_1$ is given by the section
\[4\pi i \cdot p_*\left(z_1\partial_{z_1}\right) \wedge\ldots \wedge 4\pi i \cdot p_*\left(z_n\partial_{z_n}\right),\]
which smoothly extends onto $\CP^n\setminus D_\infty$ and has simple zeroes along the intersection of the toric divisors $D_1,\ldots, D_n$ With $\CP^n\setminus D_\infty$.

To see that this section also extends regularly to the divisor at infinity $D_\infty$ and has a simple zeros along it, we compute it in affine charts meeting $D_\infty$. For every $0\leq j\leq n$ consider the affine chart $$U_j = \{[z_0:\cdots:z_n]\in \CP^n\,|\,z_j \ne 0\}.$$
Note that $U_0=\CP^n\setminus D_\infty$ and $D_\infty\subset \bigcup_{j=1}^n U_j$, thus it is enough to prove that for every $1\leq j\leq n$, the section extends regularly to $U_j$ and has a simple zero along $D_\infty\cap U_j$. For simplicity, let us prove this for $j=n$.

The chart $U_n$ is equipped with the coordinates $w_0,\ldots,w_{n-1}$ given by $w_j = \frac{z_j}{z_n}$ for $0 \le j \le n-1$. In this chart, the intersection $D_\infty \cap U_n$ is given by $\{w\in U\,|\,w_0 = 0\}$. On the overlap $U_0 \cap U_n$, the coordinate transition is given by $z_j = \frac{w_j}{w_0}$ for $1 \le j \le n-1$ and $z_n = \frac{1}{w_0}$. By the chain rule, the vector fields transform as follows: for $1 \le j \le n-1$, we have$$p_*\left(z_j \frac{\partial}{\partial z_j}\right) = w_j \frac{\partial}{\partial w_j},$$and for $j=n$, we have
$$p_*\left(z_n \frac{\partial}{\partial z_n}\right) = -w_0 \frac{\partial}{\partial w_0} - \sum_{k=1}^{n-1} w_k \frac{\partial}{\partial w_k}.$$
Substituting these into the wedge product yields:$$\bigwedge_{j=1}^n p_*\left(z_j \frac{\partial}{\partial z_j}\right) = \left( \bigwedge_{j=1}^{n-1} w_j \frac{\partial}{\partial w_j} \right) \wedge \left( -w_0 \frac{\partial}{\partial w_0} - \sum_{k=1}^{n-1} w_k \frac{\partial}{\partial w_k} \right).$$
Since all terms in the sum repeat one of the vector fields $\frac{\partial}{\partial w_1},\ldots,\frac{\partial}{\partial w_{n-1}}$, their respective wedge products vanish. Therefore only one term contributes, and after reordering to the standard orientation, we obtain:
$$\bigwedge_{j=1}^n p_*\left(z_j \frac{\partial}{\partial z_j}\right) = (-1)^n w_0 w_1 \cdot \ldots \cdot w_{n-1} \frac{\partial}{\partial w_0} \wedge \frac{\partial}{\partial w_1} \wedge \ldots \wedge \frac{\partial}{\partial w_{n-1}}.$$Therefore, the section defining $\tau_1$ becomes
$$(4\pi i)^n (-1)^n w_0 w_1 \cdot \ldots \cdot w_{n-1} \frac{\partial}{\partial w_0} \wedge \frac{\partial}{\partial w_1} \wedge \ldots \wedge \frac{\partial}{\partial w_{n-1}}.$$
Since $\frac{\partial}{\partial w_0} \wedge \ldots \wedge \frac{\partial}{\partial w_{n-1}}$ is a non-vanishing local section of $\det T\CP^n$ over $U_n$, this explicit expression shows that the section extends regularly across $D_\infty\cap U_n$. Moreover, the presence of the variable $w_0$ implies that it vanishes exactly to first order along $D_\infty\cap U_n= \{w\in U\,|\,w_0 = 0\}$. Combining this with the behavior along the toric divisors, we conclude that the section defining $\tau_1$ is $D$-canonical, which completes the proof.
\end{proof}

In Section \ref{sss: CP^n - D_infty} we construct a symplectomorphism $\iota:B(1)\to \CP^n\setminus D_\infty$. Since $B(1)$ is contractible, this induces a trivialization of $T\CP^n\vert^{}_{\CP^n\setminus D_\infty}$. We call this trivialization $\tau_B$ (B stands for the \emph{ball}.)
\begin{prop}
    Let $\tau = \tau_B$,
    let $u\colon \Sigma\to M$ be a smooth map from a Riemann surface with punctures $\Sigma$ such that all the asymptotes of $u$ are contained in $\CP^n\setminus D_\infty$. Then $c_1^{\tau_B}(u)=(n+1)\cdot(u\cdot D_\infty)$.
\end{prop}
\begin{proof}
    $\CP^n\setminus D_\infty$ is contractible, hence any two trivializations of $\det T\CP^n\vert^{}_{\CP^n\setminus D_\infty}$ are homotopic, in particular homotopic to a trivialization given by a $D_\infty$-canonical section. Since $PD([D_\infty]) = \frac {1}{n+1} c_1(T\CP^n)$ the conclusion follows.
\end{proof}

\subsection{Positivity and estimates for intersections}\label{ss: intersections}

The purpose of this section is to explore sums of the form
$$\sum_{k=1}^m c_1^\tau(u_k)+\sum_{\alpha\in A}c_1(B_\alpha),$$
where $u_1, \dots, u_m$ are the cascades of a bubbled flowline with cascades and $\{B_\alpha\}_{\alpha\in A}$ are its bubbles, that is,  $\{B_\alpha\}_{\alpha\in A}$ are holomorphic spheres, see Definition~\ref{def:bubbled floer flowline with cascades}. The main tool we use for this exploration is a comparison with the sum of the intersection numbers
$$\sum_{k=1}^m u_k\cdot D +\sum_{\alpha\in A} B_\alpha\cdot D,$$
where $D$ is a divisor Poincar\'e dual to a multiple of the first Chern class. Note that the standard definition of the intersection number does not hold in cases where a cascade has an asymptotes on $D$; we will therefore generalize it in order to encompass this case.

The first step will be dealing with the removal of punctures:
\begin{prop}\label{prop: c_1 and removing of punctures}
    Let $\Sigma$ be a punctured Riemann surface and let $u\fc \Sigma\to M$ be a smooth map, asymptotic at the punctures to loops. Assume that at the puncture $p$, the map $u$ is asymptotic to a constant loop. Let $\tau$ be a trivialization of $u^*\det TM$ over the asymptotic loops such that at $p$ the trivialization $\tau$ is induced by the constant capping. Let $\tilde{\Sigma}$ be the surface obtained from $\Sigma$ by filling in the puncture $p$ and let $\tilde{u}\fc \tilde{\Sigma}\to M$ be the continuous extension of $u$ to $\tilde{\Sigma}$. Then $c_1^\tau(u)=c_1^\tau(\tilde{u})$.
\end{prop}

\begin{proof}
     $c_1^\tau(u)$ equals the signed count of zeros of a section of $u^* \det TM$ that is transversal to the zero section and equal to a non-zero constant near the asymptotics with respect to the trivialization $\tau$. Since $\tau$ equals at $p$ to the trivialization induced from the constant capping, we can extend such a section to $\tilde{\Sigma}$ by requiring it to be constant over the capping. This extension preserves the signed count of zeros because it introduces no new zeros and remains transversal to the zero section. Consequently, this count equals the signed count of zeros of a generic section of $\tilde{u}^* \det TM$, which is, by definition, $c_1^\tau(\tilde{u})$.
\end{proof}

\subsubsection{Positivity of intersections}\label{sss: pos of int}

The first type of intersections with which we deal is intersections of a Floer solution with a divisor. Let us start with a definition:

Let $X$ be a finite dimensional smooth manifold and let $Y\subset X$ be a submanifold of codimension $2$. Let $\Sigma$ be a Riemann surface and let $f \fc \Sigma \to X$ be a smooth map. Let $z\in \Sigma$ and suppose that $f(z) \in Y$ is an \emph{isolated intersection} of $f$ with $Y$. This means that there exists an open disk $D_0 \subset \Sigma$ around $z$ and a contractible open (in relative topology) neighborhood $U \subset Y$ around $f(z)$ such that $f^{-1}(U) \cap D_0 = \{z\}$. Then the \emph{local intersection number}
$$
(f \cdot Y)_z := (f|_{D_0}) \cdot B
$$
is defined as the signed count of intersections obtained by a small perturbation of $f$ which makes it transversal to $Y$, and is independent of the choice of such perturbation.

More generally, suppose that $\Sigma$ and $Y$ are compact and connected with (possibly empty) boundaries such that $f^{-1}(Y) \cap \partial \Sigma = f^{-1}(\partial Y) \cap \Sigma = \varnothing$. Then we have a well-defined \emph{intersection number} $f \cdot Y$ of $f$ with $Y$.

Now, the following result is known as \textbf{positivity of intersections} for Floer solutions and symplectic divisors, it appears in the literature in several versions, here we provide a variation on the formulation due to Ganatra--Pomerleano:

\begin{prop}[{\cite[Lemma 4.3]{Ganatra_Pomerleano_2021_log_PSS}}]\label{prop: pos_of_interior_intersection}
   Let $(M,\omega)$ be a symplectic manifold and let $D\subset M$ be a smooth divisor. Let $J$ be a compatible almost complex structure on $(M,\omega)$ and assume that $D$ is a complex hypersurface with respect to $J$. Let $\beta$ be a $1$-form on a connected Riemann surface $(\Sigma, j)$, and let $H \fc \Sigma \times M \to \R$ be a $\Sigma$-parametrized Hamiltonian on $M$ satisfying $X_H(z,x) \in T_x D$ for every $z \in \Sigma$ and $x \in D$. Let $u \fc \Sigma \to M$ be a Floer solution for $H$ and $\beta$, i.e., a solution of the equation
    \begin{equation}\label{eq: Floer eq on a surface}
         (du - X_H \otimes \beta)^{0,1} = 0.
    \end{equation}
    If $D$ is connected and there exists an open subset $V\subset \Sigma$ with compact closure such that $u(\Sigma\setminus V)\cap D=\varnothing$, then $u^{-1}(D)$ is a finite set, and for each $z \in \im u\cap D$, the local intersection number $(u \cdot D)_z$ is at least $1$. \qed
\end{prop}

\begin{rem}
    The positivity of intersections established in Proposition~\ref{prop: pos_of_interior_intersection} applies to both the cascades part and the bubbles part of a bubbled flowline with cascades. The latter are holomorphic, satisfying the condition $(dB_\alpha)^{0,1}=0$.
\end{rem}

\subsubsection{Asymptotic intersection number}

The second type of intersections with which we deal, are intersections of the asymptote of a Floer solution with a divisor.

Let $\Sigma$ be a Riemann surface with a puncture $p$, and let $X$ be a finite dimensional smooth manifold containing a submanifold $Y$ of codimension $2$. Suppose $u\fc \Sigma \to X$ is a smooth map that is asymptotic at $p$ to a point $y \in Y$, and assume that $\im (u) \cap Y = \varnothing$. The \textbf{asymptotic intersection number} of $u$ with $Y$, denoted $u * Y$, is defined to be the intersection number of the capped map $u \# c$ with $Y$ at $y$, where $c$ is the constant capping at $y$. The orientation of $c$ is taken to be positive if the puncture at $y$ is negative, and negative if the puncture is positive.

Let $R\in\R$ and let us denote by $S_+(R), S_-(R)$ the half-infinite cylinders: 
$$S_+(R)=[R,+\infty)\times S^1\qquad\text{and}\qquad S_-(R)=(-\infty,R]\times S^1.$$

Let $(M^{2n},\omega)$ be a symplectic manifold and let $J$ be a compatible almost complex structure on $(M,\omega)$. Assume that $(M,J)$ is biholomorphic to a product $(W,J')\times (B,i)$ where $(W,J')$ is a $(2n-2)$-dimensional almost complex manifold and $B$ is an open disk in $\C$ centered at the origin equipped with the standard complex structure. 

Let $H\fc M\to \R$ be a Hamiltonian and assume that its Hamiltonian vector field $X_H$ has the form
$$X_H(w,z)=V(w)+2\pi \alpha i z,$$
where $V$ is a smooth vector field on $W$ and $\alpha\in \R\setminus\Z$ is a non-integer real number.

The following result provides a slight generalization of a formula due to Seidel \cite[Eq. (7.22)]{Seidel_Fukaya_A_infty_Strcs_Assoc_to_Lef_Fib. III}.

\begin{prop}\label{prop: asymp_intersection}
    Let $p\in W$ be a zero of $V$, let $S$ be a half-infinite cylinder and let $u\fc S\to M$ be a Floer solution with respect to $H$ and $J$ that converges to $(p,0)$ and satisfies $\im u \cap (W\times\{0\})=\varnothing$.
    \begin{itemize}
        \item If $S=S_+(R)$, for some $R\in \R$, then $u* (W\times\{0\})\geq-\lfloor\alpha\rfloor$.
        \item If $S=S_-(R)$, for some $R\in \R$, then $u* (W\times\{0\})\geq 1+\lfloor\alpha\rfloor$. \qed
    \end{itemize}
\end{prop}

For completeness, and to ensure compatibility with our sign conventions, the proofs of Proposition~\ref{prop: pos_of_interior_intersection} and Proposition~\ref{prop: asymp_intersection} are provided in Appendix~\ref{app: proof of Seidel's lemma}.

Let us present some examples that will be useful later:

\begin{exam}\label{exam: computations using Seidel's lemma}
    Let $\Delta \in (0,1)$ and $\ell \in \Z_{\geq 0}$. Let $S$ be a half-infinite cylinder and let $u \fc S \to \CP^n$ be a Floer solution with respect to $H_\ell$ and the standard complex structure $J$ on $\CP^n$. Assume that $\im u \cap (D_1 \cup \dots \cup D_n \cup D_\infty) = \varnothing$.
    
    \begin{itemize}
        \item Assume that the asymptote of $u$ lies on $D_\infty = \CP^{n-1}$. According to Claim~\ref{claim: mu in O(1)} and the definition of the Hamiltonian $H_\ell$, we can identify a neighborhood of $D_\infty \subset \CP^n$ with a neighborhood $U$ of the zero section of the line bundle $\cO(1) \to \CP^{n-1}$. In these coordinates, $H_\ell$ takes the form
        $$H_\ell(w,\lambda) = \frac{\ell + \frac{1}{2}}{1 + \|\lambda\|_h^2},$$
        and the Hamiltonian vector field is given by
        $$X_{H_\ell}(w,\lambda) = 2\pi i (\ell + 1/2) \lambda,$$
        for every $(w,\lambda) \in U$. 

        Since $\cO(1)$ is a holomorphic bundle, there is a neighborhood of the asymptote of $u$ which is biholomorphic to a product $W \times B$, where $W$ is an open set in $\CP^{n-1}$ and $B \subset \C$ is a small disk, both equipped with the standard complex structure. Moreover, we can choose this biholomorphism such that $X_{H_\ell}$ takes the form
        $$X_{H_\ell}(w,z) = 2\pi i (\ell + 1/2) z,$$
        for every $(w,z) \in W \times B$.

        Thus, by Proposition~\ref{prop: asymp_intersection} and the fact that $\im u \cap D_\infty = \varnothing$, we obtain:
        \begin{itemize}
        \item 
        If $u$ has an asymptote at $+\infty$ on $D_\infty$, then
        $$u \cdot D_\infty = u \cdot (W \times \{0\}) \geq -\ell.$$
        \item
        Similarly, if $u$ has an asymptote at $-\infty$ on $D_\infty$, then
        $$u \cdot D_\infty = u \cdot (W \times \{0\}) \geq \ell + 1.$$
        \end{itemize}
        \item Assume that the asymptote of $u$ is the constant loop $\check{x}_{0}^\ell$ located at $0 \in \Int B(1) \subset \CP^n$. Let $B_1, \dots, B_n \subset \C$ be small standard disks centered at the origin, such that $P = B_1 \times \dots \times B_n \subset B^{}_{\Delta} \subset \CP^n$. By definition, on $P$ the Hamiltonian $H_\ell$ takes the form
        $$H_\ell(z_1, \dots, z_n) = \pi \alpha_\ell (|z_1|^2 + \dots + |z_n|^2),$$
        for some constant $\alpha_\ell \in (0,1)$. Thus, the Hamiltonian vector field is given by
        $$X_{H_\ell}(z_1, \dots, z_n) = \sum_{j=1}^n (-2\pi i \alpha_\ell z_j).$$

        Let $1 \leq j \leq n$. By Proposition~\ref{prop: asymp_intersection} and the fact that $\im u \cap D_j = \varnothing$:
        \begin{itemize}
            \item If $u$ asymptotic at $+\infty$ to $0$, then
            $$u \cdot D_j = u \cdot (B_j \times \{0\}) \geq -\lfloor -\alpha_\ell \rfloor = -(-1) = 1.$$
            \item
            Similarly, if $u$ asymptotic at $-\infty$ to $0$, then
            $$u \cdot D_j = u \cdot (B_j \times \{0\}) \geq \lfloor -\alpha_\ell \rfloor + 1 = -1 + 1 = 0. $$
        \end{itemize}
    \end{itemize}
\end{exam}

\subsubsection{Internal asymptotic intersections}

The notion of internal asymptotic intersection introduced here is particular for broken flowlines with cascades. Let $H$ be an autonomous Hamiltonian on a closed symplectic manifold $(M,\omega)$ satisfying the \textbf{MB} condition. A critical submanifold with respect to $H$ will be called \textbf{mixed} if it contains both constant and non-constant $1$-periodic orbits. We assume that $H$ has no mixed critical submanifolds. We will identify critical submanifolds consisting of constant orbits with submanifolds of $M$ by sending each constant loop to its image.

Let $q_+$ and $q_-$ be critical points of a Morse function defined on the critical submanifolds of $H$. Suppose that  $u$ is a broken bubbled flowline with cascades connecting $q_-$ to $q_+$, with cascades $u_1, \dots, u_m \fc \R \times S^1 \to M$.

Note that for every $1 \leq k \leq m-1$, the asymptote of the cylinder $u_k$ at $+\infty$ is a constant orbit if and only if the asymptote of cylinder $u_{k+1}$ at $-\infty$ is a constant orbit. 

Let $1 \leq k \leq m-1$, and assume that the asymptote of cylinder $u_k$ at $+\infty$ is a constant orbit belonging to a critical submanifold $S$ with respect to $H$. In such case, we say that \textbf{$u$ has an internal asymptotic intersection along $S$}. In the case where $S$, considered as a submanifold of $M$, is contained in a smooth divisor $D$, we define the \textbf{internal asymptotic intersection number of $u$ with $D$ at $k$} as
\begin{equation*}
    u_k\mathord{\raisebox{-0.8ex}{$\big\vert$}}_{[0,\infty)\times S^1} * D \,\,+\,\, u_{k+1}\mathord{\raisebox{-0.8ex}{$\big\vert$}}_{[(-\infty,0]\times S^1} * D.  
\end{equation*}

Under certain conditions we can provide a positivity of intersections result for internal asymptotic intersections. This result is an immediate conclusion from Proposition~\ref{prop: asymp_intersection}.

\begin{lemma}\label{lemma: estimate for int asymptc intersection}
    Let $(M,\omega)$ be a closed symplectic manifold and $J$ a compatible almost complex structure. Let $H$ be an autonomous Hamiltonian on $M$ satisfying the \textbf{MB} condition. Let $u$ be a broken bubbled flowline with cascades denoted as $u_1, \dots, u_m \fc \R \times S^1 \to M$. Let $S$ be a critical submanifold with respect to $H$ consisting of constant $1$-periodic orbits, and let $D$ be a divisor containing $S$. 

    For $1 \leq k \leq m-1$, assume that the asymptotes $q_k^+$ and $q_{k+1}^-$ of the cylinders $u_k$ and $u_{k+1}$ at $+\infty$ and $-\infty$, respectively, are elements of $S$. Assume further that there exist $\alpha, \beta \in \R$ and open neighborhoods $U, U'$ near $q_k^+, q_{k+1}^-$ , respectively, such that 
    \begin{itemize}
        \item $U$ is biholomorphic to $(U\cap D,J|_D)\times (B,i)$ where $B\subset \C$ is a small ball concentrated at the origin and equipped with the standard complex structure.
        \item the Hamiltonian vector field of $H$ takes the form
        $$ X_H(w,z) = V(w) + 2\pi i \alpha z $$
        in $U$, for every $w\in U\cap D$ and $z\in B$, where $V\fc U\cap D\to TD$ is a smooth vector field.

        \item $U'$ is biholomorphic to $(U'\cap D,J|_D)\times (B',i)$ where $B'\subset \C$ is a small ball concentrated at the origin and equipped with the standard complex structure.
        \item the Hamiltonian vector field of $H$ takes the form
        $$ X_H(w,z) = V'(w) + 2\pi i \beta z $$
        in $U'$, for every $w\in U'\cap D$ and $z\in B'$, where $V'\fc U'\cap D\to TD$ is a smooth vector field.
    \end{itemize}
    Then the internal asymptotic intersection number of $u$ with $D$ at $k$ is greater than or equal to
    $$
    -\lfloor \alpha \rfloor + \lfloor \beta \rfloor + 1.$$
\end{lemma}
\begin{rem}
    Under the conditions of Lemma~\ref{lemma: estimate for int asymptc intersection},  if $\alpha\leq \beta$ then $u$ enjoys positivity of internal asymptotic intersections with $D$, at $k$.
\end{rem}
Note that Example~\ref{exam: computations using Seidel's lemma} provides concrete examples for Lemma~\ref{lemma: estimate for int asymptc intersection}:
\begin{coroll}\label{coroll: positivity of asymptotic intersections for acc. data}
    Let $\ell\in \Z_{\geq0}$ and $\Delta\in (0,1)$. Let $u$ be a bubbled flowline with cascades for $H_\ell\fc \CP^n\to \R$.
    \begin{itemize}
        \item If $u$ has an internal asymptotic intersection along the critical submanifold $\{0\}$, then for every $1\leq j\leq n$ the internal asymptotic intersection of $u$ with the divisor $D_j$ is at least $1$.

        \item If $u$ has an internal asymptotic intersection along the critical submanifold $D_\infty$, then the internal asymptotic intersection of $u$ with the divisor $D_\infty$ is at least $1$.
    \end{itemize}
\end{coroll}

\subsubsection{Intersection number of flowlines with cascades and divisors} \label{sss: Intersection and cascades}

Let $H$ be an autonomous Hamiltonian on a closed symplectic manifold $(M,\omega)$ satisfying the \textbf{MB} condition. We assume that $H$ has no mixed critical submanifolds. Let $u$ be a broken bubbled flowline with cascades $u_1, \dots, u_m \fc \R \times S^1 \to M$ and bubbles $\{B_\alpha\}_{\alpha\in A}$. Let $D$ be a transverse crossing divisor. We define the \textbf{intersection number of $u$ with $D$}, denoted by $u \cdot D$, to be the sum of all intersection numbers and asymptotic intersection numbers of $u_1, \dots, u_m$ with $D$ plus the sum of the intersection numbers of the bubbles with $D$:
$$u \cdot D = \sum_{k=1}^m (u_k \cdot D + u_k * D)+\sum_{\alpha\in A} B_\alpha\cdot D.$$

\begin{rem}
    Note that this definition assumes that for every $1 \leq k \leq m$, the intersection $\im u_k \cap D$ is finite. As well as, for every $\alpha\in A$ we should have $B_\alpha\not\subset D$.
\end{rem}

Let us formulate the connection between the relative first Chern number of the cascades of a flowline with cascades and their intersection numbers:

\begin{thm}\label{thm: c1 vs intersection number}
    Let $H$ be an autonomous Hamiltonian on a closed symplectic manifold $(M,\omega)$ satisfying the \textbf{MB} condition. Assume that $H$ has no mixed critical submanifolds. Let $u$ be a broken bubbled flowline with cascades $u_1, \dots, u_m \fc \R \times S^1 \to M$. Let $D \subset M$ be a transverse crossing divisor, Poincar\'e dual to $\frac{1}{N}c_1(TM)$ for some $N \in \N$. Assume that all the non-constant asymptotes of $u$ are contained in $M\setminus D$. Let $\tau$ be a trivialization of $\det TM$ restricted to $M \setminus D$, admitting a $D$-canonical section $s\fc M\to \det TM$ such that $\tau(\cdot,1)=s$. Then the sum of the relative first Chern numbers of the cascades of $u$ plus the sum for the first Chern class of the bubbles of $u$ equals $N$ times the intersection number of $u$ with $D$:
    $$\sum_{k=1}^m c_1^\tau(u_k)+\sum_{\alpha\in A}c_1(B_\alpha) = N(u \cdot D).$$
\end{thm}

\begin{proof}
    For every $1 \leq k \leq m$, if $u_k$ is asymptotic to a constant $1$-periodic orbit at a puncture of $\Sigma = \R \times S^1$, we fill that puncture and continuously extend $u_k$. We denote the resulting map by $\tilde{u}_k$. Thus, for every $1 \leq k \leq m$, we have $c_1^\tau(u_k) = c_1^\tau(\tilde{u}_k)$. Since $\tilde{u}_k$ has no asymptotic intersections with $D$, we deduce that 
    $$c_1^\tau(\tilde{u}_k) = N(\tilde{u}_k \cdot D) = N(u_k \cdot D + u_k * D).$$
    It follows that
    \begin{align*}
    \sum_{k=1}^m c_1^\tau(u_k) +\sum_{\alpha\in A}c_1(B_\alpha)&= \sum_{k=1}^m c_1^\tau(\tilde{u}_k) +\sum_{\alpha\in A}c_1(B_\alpha) \\
    &= \sum_{k=1}^m N(u_k \cdot D + u_k * D) +\sum_{\alpha\in A}N(B_\alpha\cdot D)\\
    &= N(u \cdot D).\qedhere
    \end{align*}
\end{proof}

\begin{rem}\label{rem: positivity of intersections for our acc. data}
    Let $\ell \in \Z_{\geq 0}$ and $\Delta \in (0,1)$. Let $u$ be a bubbled flowline with cascades for $H_\ell \fc \CP^n \to \R$. Proposition~\ref{prop: pos_of_interior_intersection} and Corollary~\ref{coroll: positivity of asymptotic intersections for acc. data} imply that $u$ satisfies positivity of intersections with respect to the divisors $D_1, \ldots, D_n, D_\infty$, if the asymptotes of $u_1$ and $u_m$ at $-\infty$ and $+\infty$ respectively, are not on these divisors. This ensures that any internal or asymptotic internal intersection with these divisors results in an intersection number of at least $1$. In particular, we deduce that
    $$\sum_{k=1}^m c_1^\tau(u_k) + \sum_{\alpha\in A} c_1(B_\alpha)=(n+1)u\cdot D_\infty\geq0.$$
\end{rem}

\subsection{The obstructions}\label{ss: obstructions}

The purpose of this section is to prove nonexistence of some Floer, continuation and homotopy solutions for the acceleration data $(H_\ell)_{\ell\geq0}$ we defined in Section~\ref{ss: acc. data}. We use the formula for dimensions of moduli spaces, together with intersection formula for relative first Chern number, computed in two trivializations.

We start with a technical claim, similar to Proposition~\ref{prop: c_1 and removing of punctures}, regarding relative and non-relative first Chern class. 
\begin{claim}\label{claim: c_1 v.s. rel c_1}
   Let $x, y \fc S^1 \to M$ be loops, and let $\hat{x}, \hat{y} \fc \bar{D}^2 \to M$ be cappings for $x$ and $y$, respectively. Let $\tau$ be a trivialization of $x^* \det TM$ and $y^* \det TM$ induced by these cappings. If $u \fc \R \times S^1 \to M$ is a cylinder asymptotic to $x$ at $-\infty$ and to $y$ at $+\infty$, then
  $$c_1^\tau(u) = c_1(\hat{x} \# u \# \hat{y}^{-1}),$$
   where $\hat{y}^{-1}$ denotes the capping $\hat{y}$ with its orientation reversed.
\end{claim}
\begin{proof}
    First, $c_1^\tau(u)$ is defined as the signed count of zeros of a section of $u^* \det TM$ that is transverse to the zero section and equals $1$ near the asymptotes, considered in the trivialization $\tau$. Since $\tau$ is induced by the cappings, we can extend such a section to the closed sphere $\hat{x} \# u \# \hat{y}^{-1}$ by setting it to be constant and equal to $1$ over the cappings. This extension preserves the signed count of zeros because it introduces no new zeros and remains transversal to the zero section. Consequently, this count equals $c_1(\hat{x} \# u \# \hat{y}^{-1})$, by definition. 
\end{proof}

The following result is an immediate corollary of Claim~\ref{claim: c_1 v.s. rel c_1} and will serve as a convenient shortcut in the proofs of the subsequent results about obstructions on certain flowlines with cascades.

\begin{prop}\label{prop: traj and modular indices}
     Let $(M,\omega)$ be a closed symplectic manifold and let $H$ be an autonomous Hamiltonian satisfying the \textbf{MB} condition. Let $q_-,q_+$ be critical points of a Morse function defined on the critical submanifolds associated to $\mathcal{A}_H$. Suppose $u$ is a broken bubbled flowline with cascades connecting $q_-$ to $q_+$, with cascades $u_1, \dots, u_m$ and a finite collection $\{B_\alpha\}_{\alpha\in A}$ of bubbles. Let $\tau$ be a trivialization of $\det TM$ along the asymptotes of $u_1,\ldots,u_m$, that is obtained from cappings of the $1$-periodic orbits. Denote
    $$ k=\mu^\tau_{FMB}(q_+;H_\ell) - \mu^\tau_{FMB}(q_-;H_\ell) + 2\sum_{k=1}^m c_1^\tau(u_k) + 2\sum_{\alpha\in A} c_1(B_\alpha). $$
    Then 
    $$\mu_{FMB}(q_+)-\mu_{FMB}(q_-)\equiv k \mod 2N,$$
    where $N$ is the minimal first Chern number of $(M,\omega)$.
\end{prop}

  Theorems~\ref{thm: diff_obst}, \ref{thm: continuation_obst} and \ref{thm: htpy_obst} below focus on the non-existence of certain Floer, continuation, and chain homotopy flowlines with cascades, respectively, associated with our acceleration data. The main tool in their proofs is the positivity of intersection discussed in the previous section, together with computations in two different trivializations. To ensure that the conditions for positivity of intersections hold, we should verify that for a given broken bubbled flowline with cascades $u$ with bubbles $\{B_\alpha\}_{\alpha\in A}$ and $m$ cascades $u_1, \ldots, u_m$, with each $u_k$ having asymptotes $q_-^k$ and $q_+^k$ at $-\infty$ and $+\infty$ respectively, \emph{each of the asymptotes is either a constant, or does not intersect the divisors $D_1, \ldots, D_n$ for every $1 \leq k \leq m$}.
  
  Given Hamiltonians $H_\ell$ and $H'_{\ell'}$ from our acceleration data associated with balls of capacity $\Delta$ and $\Delta'$ respectively, for some $\ell, \ell' \in \Z_{\geq 0}$, we note that any element of the unitary group $\U(n)$ induces a symplectomorphism of $\CP^n$ that preserves the divisor $D_\infty$, the Hamiltonians $H_\ell$ and $H'_{\ell'}$, and all the Floer--Morse--Bott indices of their $1$-periodic orbits. We thus deduce that there exists a unitary matrix $U$ such that the transformed divisors $U(D_1) \cup \dots \cup U(D_n)$ are preserved by the Hamiltonian flow, do not intersect any of the non-constant $1$-periodic orbits in the collection $\{q_-^1, q_+^1, \dots, q_-^m, q_+^m\}$, and do not contain the images of $u_1, \ldots, u_m$ or $\{B_\alpha\}_{\alpha \in A}$. Moreover, $J_0$ is also invariant under the action of $U$, thus the divisors $U(D_1) \cup \dots \cup U(D_n)$ are complex as well. Therefore, we can relabel $U(D_1), \ldots, U(D_n)$ as $D_1, \ldots, D_n$, allowing us to assume that $D_1, \ldots, D_n$ satisfy the required properties in the proofs of each of these theorems.

We recall the definition of the Hamiltonians: given $0 < \Delta < 1$ and $\ell \in \Z_{\geq 0}$, define $H_\ell \fc \CP^n \to \R$ by $H_\ell(z) = h(\Delta, \ell, \mu(z))$ for every $z \in \CP^n$. Throughout this section and under the assumption above, we utilize two distinct trivializations of the tangent bundle of $\CP^n$ along the $1$-periodic orbits of $H_\ell$ and $H'_{\ell'}$, denoted by $\tau_B$ and $\tau_T$. 

The trivialization $\tau_B$ (where the subscript $B$ refers to the \textbf{ball}) is defined as follows: for the constant generators $$\check{x}_0^\ell, \check{x}_{\ell+1}^\ell, \ldots, \check{x}_{\ell+n}^\ell,\check{y}_0^{\ell'}, \check{y}_{\ell'+1}^{\ell'}, \ldots, \check{y}_{\ell'+n}^{\ell'},$$
the trivialization $\tau_B$ is obtained from the constant cappings; for the non-constant generators $$\check{x}_1^\ell, \hat{x}_1^\ell, \ldots, \check{x}_\ell^\ell, \hat{x}_\ell^\ell,\check{y}_1^{\ell'}, \hat{y}_1^{\ell'}, \ldots, \check{y}_{\ell'}^{\ell'}, \hat{y}_{\ell'}^{\ell'},$$ the trivialization $\tau_B$ is the one induced by cappings in $\CP^n \setminus D_\infty$, as discussed in Section~\ref{sss: CP^n - D_infty} and the end of Section~\ref{ss: c_1 and intersection}. 

The trivialization $\tau_T$ (where the subscript $T$ refers to the \textbf{torus}) is defined as follows: for the constant generators $$\check{x}_0^\ell, \check{x}_{\ell+1}^\ell, \ldots, \check{x}_{\ell+n}^\ell,\check{y}_0^{\ell'}, \check{y}_{\ell'+1}^{\ell'}, \ldots, \check{y}_{\ell'+n}^{\ell'},$$
the trivialization $\tau_T$ is obtained from the constant cappings; for the non-constant generators $$\check{x}_1^\ell, \hat{x}_1^\ell, \ldots, \check{x}_\ell^\ell, \hat{x}_\ell^\ell,\check{y}_1^{\ell'}, \hat{y}_1^{\ell'}, \ldots, \check{y}_{\ell'}^{\ell'}, \hat{y}_{\ell'}^{\ell'},$$
the trivialization $\tau_T$ is given by trivialzing $T\CP^n$ over $\CP^n \setminus (D_1 \cup \ldots \cup D_n \cup D_\infty)$ is discussed in Section~\ref{sss: CPn - all divisors} and the end of Section~\ref{ss: c_1 and intersection}.

The Floer--Morse--Bott indices of the generators with respect to these trivializations are summarized in the table below:

\begin{table}[ht]
\centering
\renewcommand{\arraystretch}{1.45}
\setlength{\tabcolsep}{12pt}
\begin{tabular}{c c c c c}
\toprule
&
\begin{tabular}{c}
$\mathrm{const}$ \\[2pt]
$\check{x}_0^\ell$\\[-1pt]
\vphantom{$1 \leq j \leq \ell$}
\end{tabular}
&
\begin{tabular}{c}
$\mathrm{non\ const}$ \\[2pt]
$\check{x}_j^\ell$ \\[-1pt]
$1 \leq j \leq \ell$
\end{tabular}
&
\begin{tabular}{c}
$\mathrm{non\ const}$ \\[2pt]
$\hat{x}_j^\ell$ \\[-1pt]
$1 \leq j \leq \ell$
\end{tabular}
&
\begin{tabular}{c}
$\mathrm{const}$ \\[2pt]
$\hat{x}_j^\ell$ \\[-1pt]
$\ell + 1 \leq j \leq \ell+n$
\end{tabular}
\\
\midrule
$\mu_{\mathrm{FMB}}^{\tau_T}$
&
$0$
&
$0$
&
$2n-1$
&
$2j$
\\
$\mu_{\mathrm{FMB}}^{\tau_B}$
&
$0$
&
$-2nj$
&
$-2nj + 2n - 1$
&
$2j$
\\
\bottomrule
\end{tabular}
\caption{Values of $\mu_{\mathrm{FMB}}$ for the indicated generators.}
\label{tab:mu-fmb-values}
\end{table}
\noindent The value $2j$ on the rightmost column is obtained by $2\ell + 2 + 2(j-\ell-1) = 2j$. \bigskip

The following proposition concerns cascades and bubbles contained in the divisor $D_\infty$. We prove that for our acceleration data, all of these components are constant. Fix $\ell,\ell'\in \Z_{\geq0}$ with $\ell\leq \ell'$ and $0<\Delta'\leq \Delta<1$. To use this result for  concatenation of bubbled broken flowlines with cascades, we consider a finite sequence of Hamiltonians $(F_i)_{i=0}^N$ on $\CP^n$, given by $F_i=H_{a_i}^{\Delta_{i}}=h(\Delta_{i},a_i,\mu(\cdot))$ for every $0\leq i\leq N$. Here, the indices form a finite sequence $\ell=a_0\leq \cdots\leq a_N=\ell'$ in $\Z_{\geq0}$, with corresponding parameters $\Delta'=\Delta_{N}\leq\cdots\leq\Delta_{0}=\Delta$. We then consider $N$ homotopies, one from $F_i$ to $F_{i+1}$ for each $i\in\{0,\ldots,N-1\}$, defined via the parameters of $h$, and we use $J_0$ for the almost complex strucutre. We denote $F_-=F_0$ and $F_+=F_N$.

\begin{prop}\label{prop: contained components ball}
Let $u$ be a concatenation of broken bubbled flowline with
cascades with respect to the Hamiltonians $(F_i)_{i=0}^N$, with first negative asymptote at a generator $q_-$ for $F_-$ and last positive asymptote at a generator $q_+$ for
$F_+$. Denote the collection of the cascades of $u$ by $\cU$ and the collection of the
sphere bubbles by $\cB$. Let $\tau$ be any symplectic trivialization of $\det TM$ along the asymptotes of the cascades of $u$. If 
    $$ \mu^\tau_{FMB}(q_+;F_+) - \mu^\tau_{FMB}(q_-;F_-) + 2\sum_{v\in\cU} c_1^\tau(v) + 2\sum_{B\in\cB} c_1(B) \;\in\;\{-1,0, 1\}, $$
then every cascade or bubble of $u$ that is contained in $D_\infty$ is constant.

\end{prop}

\begin{proof} First, denote
\[\delta=\mu^\tau_{FMB}(q_+;F_+) - \mu^\tau_{FMB}(q_-;F_-) + 2\sum_{v\in\cU} c_1^\tau(v) + 2\sum_{B\in\cB} c_1(B), \]
and assume that $\delta\in\{-1,0,1\}$. We will show that every cascade or bubble of $u$ that is contained in $D_\infty$ is constant.

Second, note that $\delta$ is independent of the trivialization $\tau$, thus we can assume that $\tau=\tau_T$, where  $\tau_T$ is defined on page~\pageref{p: def tau_T}.

We denote by $\cU_\infty$ and $\cB_\infty$ the subcollections of $\cU$ and $\cB$, respectively, whose
images are contained in $D_\infty$.
Also, set
\[
\begin{split}
C^{\mathrm{tot}}=\sum_{v\in\cU}c_1^{\tau_T}(v)
    +\sum_{B\in\cB}c_1(B),\qquad\qquad 
C^{\mathrm{rem}}=\sum_{v\in\cU\setminus\cU_\infty}c_1^{\tau_T}(v)
+\sum_{B\in\cB\setminus\cB_\infty}c_1(B).
\end{split}
\]
and
\[
R:=C^{\mathrm{tot}}-C^{\mathrm{rem}}
=\sum_{v\in\cU_\infty}c^{\tau_T}_1( v)+\sum_{B\in\cB_\infty}c_1(B)
.
\]
Now, note that the Hamiltonian vector field vanishes on $D_\infty$. Thus, every cascade $v$
contained in $D_\infty$ extends to a $J_0$-holomorphic sphere $\bar v$, by removal of singularities. 
Since for every $v\in \cU_\infty$, over the asymptotes of $v$, the trivialization $\tau$ is obtained from the constant cappings, thus we deduce that $c_1^\tau(v)=c_1(\bar v)$. Since the minimal first Chern number of $\CP^n$ is $n+1$, we deduce that $R\in (n+1) \Z$. Also, since for every $v\in \cU_\infty$, the sphere $\bar v$ is pseudo-holomorphic, we deduce together with monotonicity that $c_1(\bar v)=(n+1)\omega(\bar v)\geq0$. Similarly, for every $B\in \cB_\infty$ we have $c_1(B)\geq0$. Therefore, $R\in (n+1)\Z_{\geq0}$. We shall show that $R<n+1$, and thus deduce that $R=0$.

Let us separate into several cases:
\begin{itemize}
\item Assume $\cU_\infty=\cU$ or all the cascades of $u$ are constant at the origin. In these cases the first Chern numbers of all the elements of $\cU\cup\cB$  are nonnegative, thus we can write
$C_{\tau_T}^{\mathrm{tot}}=(n+1)d$ for some $d\geq0$. Also, in these cases $q_-,q_+$ are constant and are located either at the origin or contained in $D_\infty$, thus their Floer--Morse-Bott indices are
even, and hence $\delta=0$.

If $q_-,q_+$ are located at the origin, their Floer--Morse--Bott indices are zero, and hence the fact $\delta=0$ implies that $d=0$

Otherwise, assume that $q_-,q_+$ are on $D_\infty$, they have the form $q_-=\check{x}_i^\ell$ and $q_+=\check{y}_j^{\ell'}$ for some $\ell+1\leq i\leq \ell+n$ and $\ell'+1\leq j\leq \ell'+n$, and we computed that $\mu_{FMB}^{\tau}(q_-;F_-)=2i$ and $\mu_{FMB}^{\tau}(q_+;F_+)=2j$. Thus, since $\delta=0$, we get that \[(n+1)d=i-j\leq\ell+n-(\ell'+1)\leq n-1,\]
which implies that $d=0$.

Thus in both cases we have $d=0$, and in particular $R=0$. Since each term defining $R$ is nonnegative and is positive for a
non-constant component, It follows that all the components are constant.

\item Assume that $\cU\setminus \cU_\infty\neq \varnothing$ and not all the cascades are constant and located at the origin.

As explained before, we can assume without loss of generality that the divisors $D_1, \ldots, D_n$ do not intersect any of the non-constant asymptotes of the cascades of $u$, and do not contain their images, as well as the its bubbles.

Fix $0\leq i\leq N$ and a generator $q$ of $F_i$. For now, for brevity, denote $a=a_i$. Let denote $\mathcal D(q)$ the collection of the divisors from $D_1,\ldots,D_n,D_\infty$ the intersect $q$. For every $D\in\mathcal D(q)$, write the normal Hamiltonian vector
field as $2\pi i\beta_D z$, as in Example~\ref{exam: computations using Seidel's lemma}, and define
\[
e_-(q)=\sum_{D\in\mathcal D(q)}(1+\lfloor\beta_D\rfloor),
\qquad
e_+(q)=-\sum_{D\in\mathcal D(q)}\lfloor\beta_D\rfloor.
\]
By Proposition~\ref{prop: asymp_intersection}, these are lower bounds for
the sums of asymptotic intersections at a negative and a positive end,
respectively, of all the elements $v\in \cU\setminus\cU_\infty$, which are not contained in $D_\infty$. Note that $\beta_{D_\infty}=a+\frac{1}{2}$, and for every $j\in\{1,\ldots,n\}$ we have $\beta_{D_j}\in(-1,0)$. As usual, denote by $\check{x}_0^{a},\ldots,\check{x}_{a+n}^{a},\hat{x}_1^{a},\ldots,\hat{x}_{a}^{a}$ the generators associated to $F$ with a suitable choice of Morse functions on the critical submanifolds. Writing $\mu_T(q)=\mu_{FMB}^{\tau_T}(q;F_i)$, we obtain:
\begin{center}
\small
\setlength{\tabcolsep}{3pt}
\renewcommand{\arraystretch}{1.3}
\begin{tabular}{@{}c@{\quad}c@{\quad}c@{\quad}c@{\quad}c@{\quad}c@{\quad}c@{}}
\toprule
$q$ & Location or range & $\mu_T(q)$ & $e_-$ & $e_+$
    & $\mu_T(q)-2e_-$ & $\mu_T(q)+2e_+$\\
\midrule
$\check{x}_{0}^{a}$ & $\text{origin}$ & $0$ & $0$ & $n$ & $0$ & $2n$\\
$\check{x}_{j}^{a}$ & $1\leq j\leq a$ & $0$ & $0$ & $0$ & $0$ & $0$\\
$\hat{x}_{j}^{a}$ & $1\leq j\leq a$ & $2n-1$ & $0$ & $0$ & $2n-1$ & $2n-1$\\
$\check{x}_{j}^{a}$ & $a+1\leq j\leq a+n$ & $2j$ & $a+1$ & $-a$ & $2(j-a-1)$ & $2(j-a)$\\
\bottomrule
\end{tabular}
\end{center}
In particular,
\[
0\leq\mu_T(q)-2e_-(q)\leq2n-1,
\qquad 0\leq\mu_T(q)+2e_+(q)\leq2n.\tag{$\ast$}
\]

Omit the cascades from $\cU_\infty$, the bubbles from $\cB_\infty$, and constant cascades at the origin. By the assumption, at least one cascade remains. We denote by $u_1,\ldots,u_m$ the cascades from $\cU\setminus\cU_\infty$ that are not the constant ones at the origin, in their order of appearance in $u$.

Recall that $F_i=H^{\Delta_i}_{a_i}$ for every $i\in\{0,\ldots, N\}$, and for every $0\leq i\leq j\leq N$ it holds that $a_i \le a_{j}$.
Given any two consecutive remaining
cascades,
The adjacent asymptotes contribute a nonnegative intersection, indeed, they contribute either at least $-a_i+(a_j+1)\geq1$, for some $0\leq i\leq j \leq N$, to the asymptotic intersection with $D_\infty$, or at least $n$ to the asymptotic intersection with the origin. The same estimates apply to consecutive cascades between consecutive broken flowlines. Therefore, we have 

\begin{equation*}
    u_k\mathord{\raisebox{-0.8ex}{$\big\vert$}}_{[0,+\infty)\times S^1} * D \,\,+\,\, u_{k+1}\mathord{\raisebox{-0.8ex}{$\big\vert$}}_{(-\infty,0]\times S^1} * D \;\ge\; 0,  
\end{equation*}
for every $1\leq k\leq m-1$, where $D=D_1 \cup \ldots D_n \cup D_\infty$.  Therefore, together with positivity of intersections, and Theorem~\ref{thm: c1 vs intersection number} we have 
\begin{equation*}
    C^{\mathrm{rem}}\;\geq\; 
    \left(u_1\mathord{\raisebox{-0.8ex}{$\big\vert$}}_{(-\infty,0]\times S^1} * D \,\,+\,\, u_{m}\mathord{\raisebox{-0.8ex}{$\big\vert$}}_{[0,\infty)\times S^1} * D \right).
\end{equation*}

Assume that $u_1$ is either a Floer solution for $F_i$ or a continuation trajectory from $F_i$ to $F_{i+1}$, and $u_m$ is either a Floer solution of $F_j$ or a continuation trajectory from $F_{j-1}$ to $F_{j}$, for some $0\leq i\leq j\leq N$.

Denote by $q$ the negative asymptote of $u_1$ and by $q'$ the positive asymptote of $u_m$. Since $a_0\leq\ldots\leq a_N$, from the table above we deduce that
\[e_-(q_-;F_-) \;\le\; e_-(q;F_i) \;\le\; u_1\mathord{\raisebox{-0.8ex}{$\big\vert$}}_{(-\infty,0]\times S^1} * D,\]
and 
\[ e_+(q_+;F_+) \;\le\; e_+(q';F_j) \;\le\; u_{m}\mathord{\raisebox{-0.8ex}{$\big\vert$}}_{[0,\infty)\times S^1} * D. \]

Therefore, we get
\[
C^{\mathrm{rem}}\geq e_-(q_-;F_-)+e_+(q_+;F_+).
\]
Consequently, from $(\ast)$ we deduce that
\[
\begin{split}
\delta&=\mu_T(q_+) - \mu_T(q_-) + 2C^{\mathrm{tot}}\\
&=\mu_T(q_+) - \mu_T(q_-) + 2C^{\mathrm{rem}}+2R\\
&\geq\mu_T(q_+) - \mu_T(q_-) + 2\left( e_-(q_-;F_-)+e_+(q_+;F_+)\right)+2R\\
&=
\bigl(\mu_T(q_+)+2e_+(q_+;F_+)\bigr)
-\bigl(\mu_T(q_-)-2e_-(q_-;F_-)\bigr)+2R\\
&\geq-(2n-1)+2R.
\end{split}
\]
Thus
\[2R\leq\delta+2n-1\leq 1+2n-1=2n.\]
Hence $0\leq R\leq n<n+1$, and since $R\in (n+1)\Z_{\geq0}$ we have $R=0$. Thus, as in the previous case, since each term defining $R$ is nonnegative and is positive for a
non-constant component, It follows that all the components are constant.

\end{itemize}

\end{proof}

The first list of obstructions we present here is for the Floer differential:

\begin{thm}\label{thm: diff_obst}
    Let $q_+$ and $q_-$ be critical points of the Morse function defined on the critical submanifolds with respect to $H_\ell$

    Suppose $u$ is a broken bubbled Floer flowline with cascades connecting $q_-$ to $q_+$, with a list of cascades $u_1, \dots, u_m$ and a finite collection $\{B_\alpha\}_{\alpha\in A}$ of bubbles. Let $\tau$ be any symplectic trivialization of $\det TM$ along the asymptotes of $u_1,\ldots,u_m$. If 
    $$ \mu^\tau_{FMB}(q_+;H_\ell) - \mu^\tau_{FMB}(q_-;H_\ell) + 2\sum_{k=1}^m c_1^\tau(u_k) + 2\sum_{\alpha\in A} c_1(B_\alpha) = 1, $$
    then the pair $(q_-, q_+)$ cannot be any of the following:
    \begin{itemize}
        \item $(\hat{x}^\ell_i, \hat{x}^\ell_j)$ for every $1 \leq i, j \leq \ell$.
        \item $(\check{x}^\ell_i, \check{x}^\ell_j)$ for every $0 \leq i, j \leq \ell + n$.
        \item $(\check{x}^\ell_i, \hat{x}^\ell_j)$ for every $0 \leq i \leq \ell + n$ and $1 \leq j \leq \ell$, except where $i=j$ and $n=1$, in which case $u$ does not intersect $D_\infty$.
        \item $(\hat{x}^\ell_i, \check{x}^\ell_j)$ for every $1 \leq i \leq \ell$ and $0 \leq j \leq \ell + n$, except where $j \in \{i-1, i+n\}$. Moreover, if $j=i-1$ then $u$ does not intersect $D_\infty$.
    \end{itemize}
\end{thm}

Remark~\ref{rem: index formula for Floer and cont flowlines} allows us to replace the trivialization $\tau$ from the formulation of Theorem~\ref{thm: diff_obst} with $\tau_B$ or $\tau_T$ as needed during the proof. We do the same in the proofs of Theorem~\ref{thm: continuation_obst} and Theorem~\ref{thm: htpy_obst} as well.

\begin{proof}[Proof of Theorem \ref{thm: diff_obst}] 
As explained before, we can assume without loss of generality that the divisors $D_1, \ldots, D_n$ do not intersect any of the non-constant asymptotes of $u_1, \ldots, u_m$, and do not contain the images of $u_1, \ldots, u_m$ and $\{B_\alpha\}_{\alpha \in A}$. Moreover, Proposition \ref{prop: contained components ball} guarantees that all cascades and bubbles contained in $D_\infty$ are constant and hence can be safely ignored in the proof.

We deal with each case separately:

    \medskip
    \noindent\textbf{Case 1: }\emph{
    Let $1\leq i,j\leq \ell$ and assume that $(q_-,q_+)=(\hat{x}^\ell_i,\hat{x}^\ell_j)$.}
    \medskip

    As computed in Section~\ref{ss: computations of FMB}, both $\mu_{FMB}(\hat{x}_i^\ell;H_\ell)$ and $\mu_{FMB}(\hat{x}_j^\ell;H_\ell)$ are odd modulo $2(n+1)$. In particular, the difference 
    $$\mu_{FMB}(\hat{x}_j^\ell;H_\ell) - \mu_{FMB}(\hat{x}_i^\ell;H_\ell) $$
    is even, and therefore not equal to $1$ modulo $2(n+1)$. Proposition~\ref{prop: traj and modular indices} then leads to a contradiction.

    \medskip
    \noindent\textbf{Case 2: }\emph{Let $0\leq i,j\leq \ell+n$ and assume that $(q_-,q_+)=(\check{x}^\ell_i,\check{x}^\ell_j)$.}
    \medskip

    As computed in Section~\ref{ss: computations of FMB}, both $\mu_{FMB}(\check{x}_i^\ell;H_\ell)$ and $\mu_{FMB}(\check{x}_j^\ell;H_\ell)$ are even modulo $2(n+1)$. In particular, the difference 
   $$\mu_{FMB}(\hat{x}_j^\ell;H_\ell) - \mu_{FMB}(\hat{x}_i^\ell;H_\ell) $$
    is even, and therefore not equal to $1$ modulo $2(n+1)$. Proposition~\ref{prop: traj and modular indices} then leads to a contradiction.

    \medskip
    \noindent\textbf{Case 3: }\emph{Let $0\leq i\leq \ell+n$ and $1\leq j\leq \ell$, and assume that $(q_-,q_+)=(\check{x}^\ell_i,\hat{x}^\ell_j)$. } 
    \medskip
    
    We separate into two sub cases:
    \begin{enumerate}
        \item Assume that $0 \leq i \leq \ell$. Let us use the trivialization $\tau_T$. 

With respect to this trivialization, one has that $\mu_{FMB}^{\tau_T}(\hat{x}_j^\ell;H_\ell) = 2n-1$ and $\mu_{FMB}^{\tau_T}(\check{x}_i^\ell;H_\ell) = 0$. Thus, by the positivity of intersections from Remark~\ref{rem: positivity of intersections for our acc. data} and Example~\ref{exam: computations using Seidel's lemma}, we deduce that 
$$ \sum_{k=1}^m c_1^\tau(u_k) + \sum_{\alpha\in A} c_1(B_\alpha) \geq 0. $$
Note that in the case $i=0$, the fact that $c_1^\tau(u_1)=u_1 \cdot (D_1\cup \ldots\cup D_n \cup D_\infty)\geq0$ is obtained from Proposition~\ref{prop: asymp_intersection}. It follows that
$$ 1 = \mu^\tau_{FMB}(q_+;H_\ell) - \mu^\tau_{FMB}(q_-;H_\ell) + 2\sum_{k=1}^m c_1^\tau(u_k) + 2\sum_{\alpha\in A} c_1(B_\alpha) \geq 2n-1. $$
For $n \in \N$, the only possible case in which this inequality holds is $n=1$. Moreover, since in this case the inequality is actually an equality, we deduce from Theorem~\ref{thm: c1 vs intersection number} and from the fact that $\PD(D_1\cup D_\infty)=c_1$ that 
$$ u \cdot (D_1 \cup D_\infty) = \sum_{k=1}^m c_1^\tau(u_k) + \sum_{\alpha\in A} c_1(B_\alpha) = 0. $$
When $n=1$, the divisors $D_1$ and $D_\infty$ are simply the poles of $S^2 \cong \CP^1$.

Positivity of intersections gives
$u\cdot D_\infty=0$. Next, we use the trivialization $\tau_B$.
We obtain
\begin{align*}
     1 &= \mu^\tau_{FMB}(q_+;H_\ell) - \mu^\tau_{FMB}(q_-;H_\ell) + 2\sum_{k=1}^m c_1^\tau(u_k) + 2\sum_{\alpha\in A} c_1(B_\alpha) \\
     &=(-2j+1)-(-2i)+2(n+1)u\cdot D_\infty \\
     &=1+2(i-j),
\end{align*}
and hence $i=j$.

\item Assume that $\ell+1 \leq i \leq \ell+n$. Let us use the trivialization $\tau_B$. 

By the positivity of intersections from Remark~\ref{rem: positivity of intersections for our acc. data} and by Example~\ref{exam: computations using Seidel's lemma}, we deduce:
$$ \sum_{k=1}^m c_1^\tau(u_k) + \sum_{\alpha\in A} c_1(B_\alpha) \geq c_1^\tau(u) \geq (n+1) (u\cdot D_\infty) \geq (n+1)(\ell+1). $$

Also, with respect to $\tau_B$, one has $\mu_{FMB}^{\tau_B}(\hat{x}_j^\ell;H_\ell) = -2nj+2n-1$ and $\mu_{FMB}^{\tau_B}=2i$.

It follows that
\begin{align*}
     1 &= \mu^\tau_{FMB}(q_+;H_\ell) - \mu^\tau_{FMB}(q_-;H_\ell) + 2\sum_{k=1}^m c_1^\tau(u_k) + 2\sum_{\alpha\in A} c_1(B_\alpha) \\
     &\geq  -2nj+2n-1 - 2i + 2(n+1)(\ell+1) \\
     &\geq -2n\ell+2n-1 - 2(\ell+n) + 2(n+1)(\ell+1) \\
     &= 2n+1 \\
     &\geq 3,
\end{align*}
which is a contradiction.
    \end{enumerate}

    \medskip
    \noindent\textbf{Case 4: }\emph{Let $1\leq i\leq \ell$ and $0\leq j\leq \ell+n$, and assume that $(q_-,q_+)=(\hat{x}_i^\ell,\check{x}_j^\ell)$.}
    \medskip
    
    We separate into two cases:
    \begin{enumerate}
        \item Assume that $0 \leq j \leq \ell$. 
        In this case we have to handle first the cascades contained in $D_\infty$.
        
        Let us first consider the trivialization $\tau_B$.

Theorem~\ref{thm: c1 vs intersection number} implies that
$$ \sum_{k=1}^m c_1^\tau(u_k) + \sum_{\alpha\in A} c_1(B_\alpha) = (n+1) u \cdot D_\infty. $$

Also, with respect to this trivialization, we have $\mu_{FMB}^\tau(\check{x}_j^\ell;H_\ell) = -2nj$ and $\mu_{FMB}^\tau(\hat{x}_i^\ell;H_\ell) = -2ni+2n-1$.
It follows that
\begin{align*}
1 &= \mu^\tau_{FMB}(q_+;H_\ell) - \mu^\tau_{FMB}(q_-;H_\ell) + 2\sum_{k=1}^m c_1^\tau(u_k) + 2\sum_{\alpha\in A} c_1(B_\alpha) \\
&= -2nj + 2ni - 2n + 1 + 2(n+1)u \cdot D_\infty \\
&= -2n(j+1-i) + 1 + 2(n+1)u \cdot D_\infty,
\end{align*}
yielding the relation
$$ u \cdot D_\infty = \frac{n(j+1-i)}{n+1}. $$
 $\gcd(n,n+1)=1$ which implies that $j+1-i$ must be a multiple of $n+1$. Writing $j+1-i = k(n+1)$ for some integer $k$, it follows that $u \cdot D_\infty = nk$. By the positivity of intersections (Remark~\ref{rem: positivity of intersections for our acc. data}), we have $u \cdot D_\infty \geq 0$, hence $k\in \Z_{\geq 0}$.

Next we move to the trivialization $\tau_T$. 
Theorem~\ref{thm: c1 vs intersection number} and Remark~\ref{rem: positivity of intersections for our acc. data} imply that
$$ \sum_{k=1}^m c_1^\tau(u_k) + \sum_{\alpha\in A} c_1(B_\alpha) = u \cdot (D_1 \cup \ldots \cup D_n \cup D_\infty) \geq u \cdot D_\infty. $$
Moreover, with respect to this trivialization, we have $\mu_{FMB}^\tau(\check{x}_j^\ell;H_\ell) = 0$ and $\mu_{FMB}^\tau(\hat{x}_i^\ell;H_\ell) = 2n-1$. Thus, it follows that
\begin{align*}
1 &= \mu^\tau_{FMB}(q_+;H_\ell) - \mu^\tau_{FMB}(q_-;H_\ell) + 2\sum_{k=1}^m c_1^\tau(u_k) + 2\sum_{\alpha\in A} c_1(B_\alpha) \\
&= 0 - (2n-1) + 2u \cdot (D_1 \cup \ldots \cup D_n \cup D_\infty) \\
&\geq -2n + 1 + 2u \cdot D_\infty,
\end{align*}
and thus
$$ nk = u \cdot D_\infty \leq n. $$
Since $k \in \Z_{\geq 0}$ and $nk \leq n$, we deduce that either $k=0$ or $k=1$. If $k=1$ we get that $j+1-i=n+1$, namely, $j=i+n$. Otherwise, $k=0$ and thus $j+1-i=0$ which means that $j=i-1$. Additionally, in this case $u\cdot D_\infty =k\cdot n=0$, and from the positivity of intersections, Remark~\ref{rem: positivity of intersections for our acc. data} we deduce that $\im u \cap D_\infty=\varnothing$.

\item Assume that $\ell+1 \leq j \leq \ell+n$. Let us first consider the trivialization $\tau_B$.
Theorem~\ref{thm: c1 vs intersection number} implies that
$$ \sum_{k=1}^m c_1^\tau(u_k) + \sum_{\alpha\in A} c_1(B_\alpha) = (n+1) u \cdot D_\infty. $$
Moreover, with respect to this trivialization, we have $\mu_{FMB}^\tau(\check{x}_j^\ell;H_\ell) =2j$ and $\mu_{FMB}^\tau(\hat{x}_i^\ell;H_\ell) = -2ni+2n-1$. Thus, it follows that
\begin{align*}
1 &= \mu^\tau_{FMB}(q_+;H_\ell) - \mu^\tau_{FMB}(q_-;H_\ell) + 2\sum_{k=1}^m c_1^\tau(u_k) + 2\sum_{\alpha\in A} c_1(B_\alpha) \\
&= 2j + 2ni - 2n + 1 + 2(n+1)u \cdot D_\infty ,
\end{align*}
yielding the relation
$$ u \cdot D_\infty = -\frac{j+n(i-1)}{n+1}. $$

Next we move to the trivialization $\tau_T$. 

Theorem~\ref{thm: c1 vs intersection number} and Remark~\ref{rem: positivity of intersections for our acc. data} imply that
$$ \sum_{k=1}^m c_1^\tau(u_k) + \sum_{\alpha\in A} c_1(B_\alpha) = u \cdot (D_1 \cup \ldots \cup D_n \cup D_\infty) \geq u \cdot D_\infty. $$
Also, with respect to this trivialization, we have $\mu_{FMB}^\tau(\check{x}_j^\ell;H_\ell) = 2j$ and $\mu_{FMB}^\tau(\hat{x}_i^\ell;H_\ell) = 2n-1$. Thus, it follows that
\begin{align*}
1 &= \mu^\tau_{FMB}(q_+;H_\ell) - \mu^\tau_{FMB}(q_-;H_\ell) + 2\sum_{k=1}^m c_1^\tau(u_k) + 2\sum_{\alpha\in A} c_1(B_\alpha) \\
&= 2j - (2n-1) + 2u \cdot (D_1 \cup \ldots \cup D_n \cup D_\infty) \\
&\geq 2j-2n + 1 + 2u \cdot D_\infty\\
&= 2j-2n + 1 -2\cdot\frac{j+n(i-1)}{n+1},
\end{align*}
from which it follows that $i \geq j-n$.

Since $u \cdot D_\infty = -(j+n(i-1))/(n+1)$ is an integer, we deduce that
$$ j+n(i-1) \equiv 0 \pmod{n+1}. $$
Since $n \equiv -1 \pmod{n+1}$, we get $j-i+1 \equiv 0 \pmod{n+1}$. Consequently, there exists $k \in \Z$ such that $j-i+1 = k(n+1)$. Based on the ranges $j-n \leq i \leq \ell$ and $\ell+1 \leq j \leq \ell+n$, it follows that
$$ 2 \leq j-i+1 \leq n+1, $$
which proves that $j-i+1 = n+1$, and hence $j=i+n$. \qedhere
\end{enumerate}

\end{proof}

Let $0<\Delta'\leq\Delta<1$. Consider additionally the corresponding acceleration data $(H'_\ell)_{\ell\geq0}$ for the ball $B^{}_{\Delta'}=\mu^{-1}([0,\Delta'])$, given by $H'_\ell(z)=h(\Delta',\ell,\mu(z))$, for every $z\in \CP^n$ and $\ell\in \Z_{\geq0}$. For every $\ell\in \Z_{\geq0}$ denote by
$\check{y}_0^\ell,\ldots,\check{y}_{\ell+n}^\ell,\hat{y}_1^\ell,\ldots,\hat{y}_\ell^\ell$ the generators for $CF(H'_\ell)$, as in Theorem~\ref{thm: CF(H_l;Z)}.

Let $\ell,\ell'\in \Z_{\geq0}$ and assume that $\ell\leq \ell'$. Let $\chi\fc \R\to \R$ be a smooth non-increasing function, satisfying $\chi(s)=\Delta$ for every $s<0$ and $\chi(s)=\Delta'$ for every $s>1$. Similarly, let $\lambda\fc \R\to \R$ be a smooth non-decreasing function, satisfying $\lambda(s)=\ell$ for every $s<0$ and $\lambda(s)=\ell'$ for every $s>1$. 

Now, define a monotone homotopy $K\fc \R\times \CP^n\to \R$ form $H_\ell$ to $H'_{\ell'}$ by $K(s,z)=h(\chi(s),\lambda(s),\mu(z))$ for every $(s,z)\in \R\times \CP^n$.

\begin{thm}\label{thm: continuation_obst}
 
    Let $q_-$ and $q_+$ be critical points of the Morse functions associated with the critical submanifolds with respect to $H_\ell$ and $H'_{\ell'}$, respectively.
    
    Suppose $u$ is a broken bubbled continuation flowlines with cascades connecting $q_-$ to $q_+$, with a list of cascades $u_1, \dots, u_m$ and a finite collection $\{B_\alpha\}_{\alpha\in A}$ of bubbles. Let $\tau$ be a symplectic trivialization of $\det TM$ along the asymptotes of $u_1,\ldots,u_m$. If 
    $$ \mu^\tau_{FMB}(q_+;H'_{\ell'}) - \mu^\tau_{FMB}(q_-;H_\ell) + 2\sum_{k=1}^m c_1^\tau(u_k) + 2\sum_{\alpha\in A} c_1(B_\alpha) = 0, $$
    then the pair $(q_-, q_+)$ cannot be any of the following:
    \begin{itemize}
        \item $(\check{x}^\ell_i, \hat{y}^{\ell'}_j)$ for every $0 \leq i \leq \ell + n$ and $1 \leq j \leq \ell'$.
        \item $(\hat{x}^\ell_i, \check{y}^{\ell'}_j)$ for every $1 \leq i \leq \ell$ and $0 \leq j \leq \ell' + n$.
        \item $(\hat{x}^\ell_i, \hat{y}^{\ell'}_j)$ for every $1 \leq i \leq \ell$ and $1\leq j\leq \ell'$, except where $i=j$.
        \item $(\check{x}^\ell_i, \check{y}^{\ell'}_j)$  for every $0 \leq i \leq \ell+n$ and $0\leq j\leq \ell'+n$, except where $i=j$.
    \end{itemize}
\end{thm}

\begin{proof}
    As explained before, we can assume without loss of generality that the divisors $D_1, \ldots, D_n$ do not intersect any of the non-constant asymptotes of $u_1, \ldots, u_m$, and do not contain the images of $u_1, \ldots, u_m$ and $\{B_\alpha\}_{\alpha \in A}$. Moreover, Proposition \ref{prop: contained components ball} guarantees that all cascades and bubbles contained in $D_\infty$ are constant and hence can be safely ignored in the proof.

We deal with each case separately:

    \medskip
    \noindent\textbf{Case 1: }\emph{Let $0 \leq i \leq \ell + n$ and $1 \leq j \leq \ell'$, and assume that $(q_-,q_+)=(\check{x}^\ell_i, \hat{y}^{\ell'}_j)$}.
    \medskip
    
    As computed in Section~\ref{ss: computations of FMB},  $\mu_{FMB}(\hat{y}^{\ell'}_j;H'_{\ell'})$ is odd and $\mu_{FMB}(\check{x}^\ell_i;H_\ell)$ is even modulo $2(n+1)$. In particular, the difference 
    $$\mu_{FMB}(\hat{y}^{\ell'}_j;H'_{\ell'}) - \mu_{FMB}(\check{x}^\ell_i;H_\ell) $$
    is odd, and therefore not equal to $0$ modulo $2(n+1)$. Proposition~\ref{prop: traj and modular indices} then leads to a contradiction.

    \medskip
    \noindent\textbf{Case 2: }\emph{Let $1 \leq i \leq \ell$ and $0 \leq j \leq \ell'+n$, and assume that $(q_-,q_+)=(\hat{x}^\ell_i, \check{y}^{\ell'}_j)$}.
    \medskip
    
    As computed in Section~\ref{ss: computations of FMB},  $\mu_{FMB}(\check{y}^{\ell'}_j;H'_{\ell'})$ is even and $\mu_{FMB}(\hat{x}^\ell_i;H_\ell)$ is odd modulo $2(n+1)$. In particular, the difference 
    $$\mu_{FMB}(\check{y}^{\ell'}_j;H'_{\ell'}) - \mu_{FMB}(\hat{x}^\ell_i;H_\ell) $$
    is odd, and therefore not equal to $0$ modulo $2(n+1)$. Proposition~\ref{prop: traj and modular indices} then leads to a contradiction.

    \medskip
    \noindent\textbf{Case 3: }\emph{Let $1\leq i\leq \ell$ and $1\leq j\leq \ell'$, and assume that $(q_-,q_+)=(\hat{x}^\ell_i,\hat{y}^{\ell'}_j)$}.
    \medskip
    
    Let us use the trivialization $\tau_T$. First, with respect to this trivialization, we have that $\mu_{FMB}^\tau(\hat{y}_j^{\ell'};H'_{\ell'}) =\mu_{FMB}^\tau(\hat{x}_i^\ell;H_\ell) = 2n-1$. Second, by the positivity of intersections of Remark~\ref{rem: positivity of intersections for our acc. data}, we deduce that 
$$ \sum_{k=1}^m c_1^\tau(u_k) + \sum_{\alpha\in A} c_1(B_\alpha) \geq 0. $$
Hence, It follows that
\begin{align*}
    0 &= \mu^\tau_{FMB}(q_+;H'_{\ell'}) - \mu^\tau_{FMB}(q_-;H_\ell) + 2\sum_{k=1}^m c_1^\tau(u_k) + 2\sum_{\alpha\in A} c_1(B_\alpha)\\
    &=2n-1 - (2n-1) + 2\sum_{k=1}^m c_1^\tau(u_k) + 2\sum_{\alpha\in A} c_1(B_\alpha)\\
    &= 2\sum_{k=1}^m c_1^\tau(u_k) + 2\sum_{\alpha\in A} c_1(B_\alpha) \geq0,
\end{align*}
which yields, using Theorem~\ref{thm: c1 vs intersection number} and Remark~\ref{rem: positivity of intersections for our acc. data} that
   \[ 0\leq u\cdot D_\infty\leq u\cdot(D_1\cup\ldots\cup D_n\cup D_\infty)=\sum_{k=1}^m c_1^\tau(u_k) + \sum_{\alpha\in A} c_1(B_\alpha)\\
    =0,\]
    and thus, we deduce that $u\cdot D_\infty=0$.

     Next we move to the trivialization $\tau_B$. With respect to this trivialization, we find that $\mu_{FMB}^\tau(\hat{y}_j^{\ell'};H'_{\ell'}) = -2nj+2n-1$ and $\mu_{FMB}^\tau(\hat{x}_i^\ell;H_\ell) = -2ni+2n-1$. From Theorem~\ref{thm: c1 vs intersection number} we deduce that
$$ \sum_{k=1}^m c_1^\tau(u_k) + \sum_{\alpha\in A} c_1(B_\alpha)=(n+1) u\cdot D_\infty=0. $$
It follows that
\begin{align*}
     0 &= \mu^\tau_{FMB}(q_+;H'_{\ell'}) - \mu^\tau_{FMB}(q_-;H_{\ell}) + 2\sum_{k=1}^m c_1^\tau(u_k) + 2\sum_{\alpha\in A} c_1(B_\alpha) \\
     &= -2nj+2n-1-(-2ni+2n-1) \\
     &=-2n(j-i),
\end{align*}
and hence $i=j$.

    \medskip
    \noindent\textbf{Case 4: }\emph{Let $0\leq i\leq \ell+n$ and $0\leq j\leq \ell'+n$, and assume that $(q_-,q_+)=(\check{x}^\ell_i,\check{y}^{\ell'}_j)$}.\medskip

    We have eight cases to check:
    \begin{enumerate}
    \item Assume that $i=0$ and $1\leq j\leq \ell'$. Let us use the trivialization $\tau_T$. 

    With respect to this trivialization, we have that $\mu_{FMB}^\tau(\check{y}_j^{\ell'};H'_{\ell'}) =\mu_{FMB}^\tau(\check{x}_i^\ell;H_\ell) = 0$. By Theorem~\ref{thm: c1 vs intersection number} and Example~\ref{exam: computations using Seidel's lemma}, we deduce that 
$$ \sum_{k=1}^m c_1^\tau(u_k) + \sum_{\alpha\in A} c_1(B_\alpha) =u\cdot(D_1\cup\ldots\cup D_n\cup D_\infty)\geq u_1\cdot(D_1\cup\ldots\cup D_n)\geq 0. $$

It follows that
\begin{align*}
    0 &= \mu^\tau_{FMB}(q_+;H'_{\ell'}) - \mu^\tau_{FMB}(q_-;H_\ell) + 2\sum_{k=1}^m c_1^\tau(u_k) + 2\sum_{\alpha\in A} c_1(B_\alpha)\\
    &= 2\sum_{k=1}^m c_1^\tau(u_k) + 2\sum_{\alpha\in A} c_1(B_\alpha)\geq0,
\end{align*}
which yields, using Theorem~\ref{thm: c1 vs intersection number} and Remark~\ref{rem: positivity of intersections for our acc. data} that
   $$ 0\leq u\cdot D_\infty\leq u\cdot(D_1\cup\ldots\cup D_n\cup D_\infty)=\sum_{k=1}^m c_1^\tau(u_k) + \sum_{\alpha\in A} c_1(B_\alpha)\\
    =0,$$ 
    thus we deduce that $u\cdot D_\infty=0$.

   Next we move to the trivialization $\tau_B$. With respect to this trivialization, we have that $\mu_{FMB}^\tau(\check{y}_j^{\ell'};H'_{\ell'}) = -2nj$ and $\mu_{FMB}^\tau(\check{x}_i^\ell;H_\ell) = 0$. From Theorem~\ref{thm: c1 vs intersection number} we deduce that
$$ \sum_{k=1}^m c_1^\tau(u_k) + \sum_{\alpha\in A} c_1(B_\alpha)=(n+1) u\cdot D_\infty=0. $$
It follows that
\begin{align*}
     0 &= \mu^\tau_{FMB}(q_+;H_\ell) - \mu^\tau_{FMB}(q_-;H'_{\ell'}) + 2\sum_{k=1}^m c_1^\tau(u_k) + 2\sum_{\alpha\in A} c_1(B_\alpha) \\
     &=-2nj,
\end{align*}
and hence $j=0$ which is a contradiction.

 \item Assume that $i=0$ and $\ell'+1\leq j\leq \ell'+n$. Let us use the trivialization $\tau_T$. 

    With respect to this trivialization, we have that $\mu_{FMB}^\tau(\check{y}_j^{\ell'};H'_{\ell'}) =2j$ and $\mu_{FMB}^\tau(\check{x}_i^\ell;H_\ell) = 0$. By Theorem~\ref{thm: c1 vs intersection number} and Example~\ref{exam: computations using Seidel's lemma}, we deduce that 
$$ \sum_{k=1}^m c_1^\tau(u_k) + \sum_{\alpha\in A} c_1(B_\alpha) =u\cdot(D_1\cup\ldots\cup D_n\cup D_\infty)\geq u_1\cdot(D_1\cup\ldots\cup D_n)+u_m\cdot D_\infty\geq -\ell'. $$

It follows that
\begin{align*}
    0 &= \mu^\tau_{FMB}(q_+;H'_{\ell'}) - \mu^\tau_{FMB}(q_-;H_\ell) + 2\sum_{k=1}^m c_1^\tau(u_k) + 2\sum_{\alpha\in A} c_1(B_\alpha)\\
    &= 2j-0+2\sum_{k=1}^m c_1^\tau(u_k) + 2\sum_{\alpha\in A} c_1(B_\alpha)\\
    &\geq2j-2\ell'.
\end{align*}
Thus we deduce that $j\leq \ell'$, contradicting the assumption $\ell'+1\leq j$.

        \item Assume that $1\leq i\leq \ell$ and $j=0$. Let us use the trivialization $\tau_T$. 

    With respect to this trivialization, we have that $\mu_{FMB}^\tau(\check{y}_j^{\ell'};H'_{\ell'}) =\mu_{FMB}^\tau(\check{x}_i^\ell;H_\ell) = 0$. By Theorem~\ref{thm: c1 vs intersection number} and Example~\ref{exam: computations using Seidel's lemma}, we deduce that 
$$ \sum_{k=1}^m c_1^\tau(u_k) + \sum_{\alpha\in A} c_1(B_\alpha) =u\cdot(D_1\cup\ldots\cup D_n\cup D_\infty)\geq u_m\cdot(D_1\cup\ldots\cup D_n)\geq n\cdot1=n. $$

It follows that
\begin{align*}
    0 &= \mu^\tau_{FMB}(q_+;H'_{\ell'}) - \mu^\tau_{FMB}(q_-;H_\ell) + 2\sum_{k=1}^m c_1^\tau(u_k) + 2\sum_{\alpha\in A} c_1(B_\alpha)\\
    &= 2\sum_{k=1}^m c_1^\tau(u_k) + 2\sum_{\alpha\in A} c_1(B_\alpha)\\
    &\geq 2n,
\end{align*}
a contradiction.

  \item Assume that $1\leq i\leq \ell$ and $1\leq j\leq \ell'$. Let us use the trivialization $\tau_T$. 

    With respect to this trivialization, we have that $\mu_{FMB}^\tau(\check{y}_j^{\ell'};H'_{\ell'}) =\mu_{FMB}^\tau(\check{x}_i^\ell;H_\ell) = 0$. By the positivity of intersections from Remark~\ref{rem: positivity of intersections for our acc. data}, we deduce that 
$$ \sum_{k=1}^m c_1^\tau(u_k) + \sum_{\alpha\in A} c_1(B_\alpha) \geq 0. $$

It follows that
\begin{align*}
    0 &= \mu^\tau_{FMB}(q_+;H'_{\ell'}) - \mu^\tau_{FMB}(q_-;H_\ell) + 2\sum_{k=1}^m c_1^\tau(u_k) + 2\sum_{\alpha\in A} c_1(B_\alpha)\\
    &= 2\sum_{k=1}^m c_1^\tau(u_k) + 2\sum_{\alpha\in A} c_1(B_\alpha) \geq0,
\end{align*}
which yields, using Theorem~\ref{thm: c1 vs intersection number} and Remark~\ref{rem: positivity of intersections for our acc. data} that
   $$ 0\leq u\cdot D_\infty\leq u\cdot(D_1\cup\ldots\cup D_n\cup D_\infty)=\sum_{k=1}^m c_1^\tau(u_k) + \sum_{\alpha\in A} c_1(B_\alpha)\\
    =0,$$ 
    thus we deduce that $u\cdot D_\infty=0$.

Next we move to the trivialization $\tau_B$. With respect to this trivialization, we have that $\mu_{FMB}^\tau(\check{y}_j^{\ell'};H'_{\ell'}) = -2nj$ and $\mu_{FMB}^\tau(\check{x}_i^\ell;H_\ell) = -2ni$. From Theorem~\ref{thm: c1 vs intersection number} we deduce that
$$ \sum_{k=1}^m c_1^\tau(u_k) + \sum_{\alpha\in A} c_1(B_\alpha)=(n+1) u\cdot D_\infty=0. $$
It follows that
\begin{align*}
     0 &= \mu^\tau_{FMB}(q_+;H'_{\ell'}) - \mu^\tau_{FMB}(q_-;H_{\ell}) + 2\sum_{k=1}^m c_1^\tau(u_k) + 2\sum_{\alpha\in A} c_1(B_\alpha) \\
     &=-2n(j-i),
\end{align*}
and hence $i=j$.

  \item Assume that $1\leq i\leq \ell$ and $\ell'+1\leq j\leq \ell'+n$. Let us use the trivialization $\tau_T$. 

    With respect to this trivialization, we have that $\mu_{FMB}^\tau(\check{y}_j^{\ell'};H'_{\ell'})=2j$ and $\mu_{FMB}^\tau(\check{x}_i^\ell;H_\ell) = 0$. By Theorem~\ref{thm: c1 vs intersection number} and Example~\ref{exam: computations using Seidel's lemma}, we deduce that  
$$ \sum_{k=1}^m c_1^\tau(u_k) + \sum_{\alpha\in A} c_1(B_\alpha) \geq u_m\cdot D_\infty\geq -\ell'. $$

It follows that
\begin{align*}
    0 &= \mu^\tau_{FMB}(q_+;H'_{\ell'}) - \mu^\tau_{FMB}(q_-;H_\ell) + 2\sum_{k=1}^m c_1^\tau(u_k) + 2\sum_{\alpha\in A} c_1(B_\alpha)\\
    &= 2j+2\sum_{k=1}^m c_1^\tau(u_k) + 2\sum_{\alpha\in A} c_1(B_\alpha)\\
    &\geq2j-2\ell',
\end{align*}
hence $\ell'+1\le j \le \ell'$ which is a contradiction.

    \item Assume that $\ell+1\leq i\leq \ell+n$ and $j=0$. Let us use the trivialization $\tau_T$. 

    With respect to this trivialization, we find that $\mu_{FMB}^\tau(\check{y}_j^{\ell'};H'_{\ell'})=0$ and $\mu_{FMB}^\tau(\check{x}_i^\ell;H_\ell) = 2i$. By Theorem~\ref{thm: c1 vs intersection number} and Example~\ref{exam: computations using Seidel's lemma}, we deduce that 
    \begin{align*}
        \sum_{k=1}^m c_1^\tau(u_k) + \sum_{\alpha\in A} c_1(B_\alpha) &=u\cdot(D_1\cup\ldots\cup D_n\cup D_\infty)\geq u_1\cdot D_\infty+u_m\cdot(D_1\cup\ldots\cup D_n)\\
        &\geq \ell+1+n\cdot1 = \ell+1+n. 
    \end{align*}

It follows that
\begin{align*}
    0 &= \mu^\tau_{FMB}(q_+;H'_{\ell'}) - \mu^\tau_{FMB}(q_-;H_\ell) + 2\sum_{k=1}^m c_1^\tau(u_k) + 2\sum_{\alpha\in A} c_1(B_\alpha)\\
    &=0-2i+ 2\sum_{k=1}^m c_1^\tau(u_k) + 2\sum_{\alpha\in A} c_1(B_\alpha)\\
    &\geq -2i+2(\ell+1+n).
\end{align*}
thus we deduce that
$$\ell+n+1\leq i\leq \ell+n$$
and this is a contradiction.

  \item Assume that $\ell+1\leq i\leq \ell+n$ and $1\leq j\leq \ell'$. Let us use the trivialization $\tau_B$. 

    With respect to this trivialization, we have that $\mu_{FMB}^\tau(\check{y}_j^{\ell'};H'_{\ell'})=-2nj$ and $\mu_{FMB}^\tau(\check{x}_i^\ell;H_\ell) = 2i$. By Theorem~\ref{thm: c1 vs intersection number} we know that
$$ \sum_{k=1}^m c_1^\tau(u_k) + \sum_{\alpha\in A} c_1(B_\alpha) = (n+1)u\cdot D_\infty. $$

It follows that
\begin{align*}
    0 &= \mu^\tau_{FMB}(q_+;H'_{\ell'}) - \mu^\tau_{FMB}(q_-;H_\ell) + 2\sum_{k=1}^m c_1^\tau(u_k) + 2\sum_{\alpha\in A} c_1(B_\alpha)\\
    &=-2nj-2i+2(n+1)u\cdot D_\infty
\end{align*}
which yields that
$$ u\cdot D_\infty=\frac{nj+i}{n+1},$$
is an integer.
Since $n\equiv -1 \pmod{n+1}$ and $nj+i\equiv0 \pmod{n+1}$ we deduce that 
$$-j+i\equiv0 \pmod{n+1},$$
thus there exists $r\in \Z$ such that $i=j+(n+1)r$. Thus we can write $u\cdot D_\infty=j+r$. Additionally, since $\ell+1\leq i\leq \ell+n$ we get that
$$\ell+1-j\leq (n+1)r\leq \ell+n-j.\qquad(\star)$$

Additionally, by Example~\ref{exam: computations using Seidel's lemma} and Remark~\ref{rem: positivity of intersections for our acc. data} we deduce that
$$u\cdot D_\infty\geq u_1* D_\infty\geq\ell+1,$$
thus 
$j+r\geq \ell+1.$

 Next we move to the trivialization $\tau_T$. With respect to this trivialization, we have that $\mu_{FMB}^\tau(\check{y}_j^{\ell'};H'_{\ell'}) = 0$ and $\mu_{FMB}^\tau(\check{x}_i^\ell;H_\ell) = 2i$. By Theorem~\ref{thm: c1 vs intersection number} and Example~\ref{exam: computations using Seidel's lemma}, we deduce that  
$$ \sum_{k=1}^m c_1^\tau(u_k) + \sum_{\alpha\in A} c_1(B_\alpha) \geq u\cdot D_\infty=j+r. $$

It follows that
\begin{align*}
    0 &= \mu^\tau_{FMB}(q_+;H'_{\ell'}) - \mu^\tau_{FMB}(q_-;H_\ell) + 2\sum_{k=1}^m c_1^\tau(u_k) + 2\sum_{\alpha\in A} c_1(B_\alpha)\\
    &=0-2i+2u\cdot (D_1\cup \ldots\cup D_n\cup D_\infty)\\
    &\geq-2i+2u\cdot D_\infty\\
    &=-2i+2(j+r),
\end{align*}
Thus $j+r\leq i=j+(n+1)r$, which implies that $r\leq nr+r$, so $0\leq nr$ which means that $r\geq0$.

If $j>\ell+n$, then by $(\star)$ we get that $(n+1)r<0$ which means that $r\leq -1$ and this is a contradiction. Thus $j\leq \ell+n$. 

If $j<\ell+1$, then $j\leq \ell$, so $0\leq \ell-j$. Since $j+r\geq \ell+1$ we deduce that $r\geq \ell-j+1>0$, thus by $(\star)$ we get that
$$(n+1)(\ell-j+1)\leq (n+1)r\leq \ell+n-j.$$
Since $(n+1)(\ell-j+1)=n\ell-nj+1+\ell+n-j$ we deduce that
$$n(\ell-j)+1\leq 0,$$
but \[ n(\ell-j)+1 \ge 0 + 1 = 1,\]
hence $1\le 0 $, a contradiction.

Thus we get that $j\geq \ell+1$. This shows that
$$\ell+1\leq i,j\leq \ell+n,$$
and we know that $i=j+(n+1)r$, which implies that $r=0$ and therefore $i=j$.

\item Assume that $\ell+1\leq i\leq \ell+n$ and $\ell'+1\leq j\leq \ell'+n$. Let us use the trivialization $\tau_B$. 

    With respect to this trivialization, we have that $\mu_{FMB}^\tau(\check{y}_j^{\ell'};H'_{\ell'})=2j$ and $\mu_{FMB}^\tau(\check{x}_i^\ell;H_\ell) = 2i$. By Theorem~\ref{thm: c1 vs intersection number} we know that
$$ \sum_{k=1}^m c_1^\tau(u_k) + \sum_{\alpha\in A} c_1(B_\alpha) = (n+1)u\cdot D_\infty. $$

It follows that
\begin{align*}
    0 &= \mu^\tau_{FMB}(q_+;H'_{\ell'}) - \mu^\tau_{FMB}(q_-;H_\ell) + 2\sum_{k=1}^m c_1^\tau(u_k) + 2\sum_{\alpha\in A} c_1(B_\alpha)\\
    &=2(j-i)+2(n+1)u\cdot D_\infty
\end{align*}
which yields that
$$ u\cdot D_\infty=\frac{i-j}{n+1},$$
is an integer. Thus we can write $i=j+(n+1)r$, where $r=u\cdot D_\infty$.

 Next we move to the trivialization $\tau_T$. With respect to this trivialization, we still have that $\mu_{FMB}^\tau(\check{y}_j^{\ell'};H'_{\ell'}) = 2j$ and $\mu_{FMB}^\tau(\check{x}_i^\ell;H_\ell) = 2i$. By Theorem~\ref{thm: c1 vs intersection number} and Example~\ref{exam: computations using Seidel's lemma}, we deduce that  
$$ \sum_{k=1}^m c_1^\tau(u_k) + \sum_{\alpha\in A} c_1(B_\alpha) \geq u\cdot D_\infty=r. $$
It follows that
\begin{align*}
    0 &= \mu^\tau_{FMB}(q_+;H'_{\ell'}) - \mu^\tau_{FMB}(q_-;H_\ell) + 2\sum_{k=1}^m c_1^\tau(u_k) + 2\sum_{\alpha\in A} c_1(B_\alpha)\\
    &=2j-2i+2u\cdot (D_1\cup \ldots\cup D_n\cup D_\infty)\\
    &\geq2(j-i)+2u\cdot D_\infty\\
    &=2(j-i)+2r,
\end{align*}
therefore
$$r\leq i-j=(n+1)r=nr+r,$$
which means that $nr\geq 0$ and hence $r\geq 0$. 

If $r\geq1$ then since $i\leq \ell+n$ and $j\geq\ell'+1$ we get that
$$\ell+n\geq i=j+(n+1)r\geq j+n+1\geq\ell'+1+n+1=\ell'+n+2,$$
thus $\ell\geq\ell'+2>\ell'\geq\ell$, and this is a contradiction. Thus we deduce that $r=0$ and hence $i=j$. \qedhere
\end{enumerate}
\end{proof}

The next theorem deals with obstructions to solutions related to the analysis of chain homotopy maps.

\begin{thm}\label{thm: htpy_obst}
    Let $0 < \Delta'\leq\Delta < 1$ and $\ell,\ell' \in \Z_{\geq 0}$ with $\ell\leq \ell'$. Define $H_\ell,H'_{\ell'} \fc \CP^n \to \R$ by $H_\ell(z) = h(\Delta, \ell, \mu(z))$ and $H'_{\ell'}(z) = h(\Delta', \ell', \mu(z))$ for every $z \in \CP^n$. Let $q_-$ and $q_+$ be critical points of the Morse functions associated with the critical submanifolds with respect to $H_\ell$ and $H'_{\ell'}$, respectively. 

    Suppose either that 
    \begin{enumerate}
        \item $K$ is some homotopy connecting $h(\Delta,\ell,\mu(z))$ to $h(\Delta',\ell',\mu(z))$ through autonomous $\U(n)$-invariant Hamiltonians, and there exists $u$, a broken bubbled continuation flowline with cascades connecting $q_-$ to $q_+$, with cascades $u_1, \dots, u_m$ and a finite collection $\{B_\alpha\}_{\alpha\in A}$ of bubbles, such that 
    $$ \mu^\tau_{FMB}(q_+;H'_{\ell'}) - \mu^\tau_{FMB}(q_-;H_\ell) + 2\sum_{k=1}^m c_1^\tau(u_k) + 2\sum_{\alpha\in A} c_1(B_\alpha) = -1, $$
    where $\tau$ is a symplectic trivialization of $\det TM$ along the asymptotes of $u_1,\ldots,u_m$,
    \item or that $K_1$ is a homotopy connecting $H_- := h(\Delta,\ell,\mu(z))$ to $H_0 := h(\Delta,\ell',\mu(z))$, and $K_2$ connects $H_0 = h(\Delta,\ell',\mu(z))$ to $H_+ := h(\Delta',\ell',\mu(z))$, both through autonomous $\U(n)$-invariant Hamiltonians, and $u$, $v$ are broken bubbled continuation flowline with cascades for $K_1$ and $K_2$ respectively, with cascades $u_1, \dots, u_{m'}$ and a finite collection $\{B_\xi\}_{\xi\in X}$ of bubbles for $u$ and cascades $v_1, \dots, v_{m''}$ and a finite collection $\{B_\zeta\}_{\zeta\in Z}$ of bubbles, for $v$, such that $u$ connects $q_-$ to some $q_0$, and $v$ connects $q_0$ to $q_+$, and such that, denoting 
    \begin{align*}
    \mu^\tau_{FMB}(q_0;H_0) - \mu^\tau_{FMB}(q_-;H_-) + 2\sum_{k=1}^{m'} c_1^\tau(u_k) + 2\sum_{\xi\in X} c_1(B_\xi) &= i_u, \text{ and}\\
   \mu^\tau_{FMB}(q_+;H_+) - \mu^\tau_{FMB}(q_0;H_0) + 2\sum_{k=1}^{m''} c_1^\tau(v_k) + 2\sum_{\zeta\in Z} c_1(B_\zeta) &= i_v,
    \end{align*}
    with respect to $\tau$, a symplectic trivialization of $\det TM$ along the asymptotes of $u_1,\ldots,u_{m'}$, $v_1,\ldots,v_{m''}$, and such that $i_u + i_v = -1$,
    \item or that $K_1$ is a homotopy connecting $H_-:=h(\Delta,\ell,\mu(z))$ to $H_0:=h(\Delta',\ell,\mu(z))$, and $K_2$ connects $H_0 = h(\Delta',\ell,\mu(z))$ to $H_+ := h(\Delta',\ell',\mu(z))$, both through autonomous $\U(n)$-invariant Hamiltonians, and $u$, $v$ are broken bubbled continuation flowline with cascades for $K_1$ and $K_2$ respectively, with cascades $u_1, \dots, u_{m'}$ and a finite collection $\{B_\xi\}_{\xi\in X}$ of bubbles for $u$ and cascades $v_1, \dots, v_{m''}$ and a finite collection $\{B_\zeta\}_{\zeta\in Z}$ of bubbles, for $v$, such that $u$ connects $q_-$ to some $q_0$, and $v$ connects $q_0$ to $q_+$, and such that, denoting 
    \begin{align*}
    \mu^\tau_{FMB}(q_0;H_0) - \mu^\tau_{FMB}(q_-;H_-) + 2\sum_{k=1}^{m'} c_1^\tau(u_k) + 2\sum_{\xi\in X} c_1(B_\xi) &= i_u, \text{ and}\\
   \mu^\tau_{FMB}(q_+;H_+) - \mu^\tau_{FMB}(q_0;H_0) + 2\sum_{k=1}^{m''} c_1^\tau(v_k) + 2\sum_{\zeta\in Z} c_1(B_\zeta) &= i_v,
    \end{align*}
    with respect to $\tau$, a symplectic trivialization of $\det TM$ along the asymptotes of $u_1,\ldots,u_{m'}$, $v_1,\ldots,v_{m''}$, and such that $i_u + i_v = -1$.
    \end{enumerate}
    Then, the pair $(q_-, q_+)$ cannot be any of the following:
    \begin{itemize}
        \item $(\check{x}^\ell_i, \check{y}^{\ell'}_j)$ for every $0 \leq i \leq \ell + n$ and $0 \leq j \leq \ell'+n$.
        \item $(\hat{x}^\ell_i, \hat{y}^{\ell'}_j)$ for every $1 \leq i \leq \ell$ and $1 \leq j \leq \ell' $.
        \item $(\check{x}^\ell_i, \hat{y}^{\ell'}_j)$ for every $0 \leq i \leq \ell+n$ and $1\leq j\leq \ell'$.
        \item $(\hat{x}^\ell_i, \check{y}^{\ell'}_j)$ for every
        $1\leq i\leq\ell$ and $0\leq j\leq\ell'+n$,
    except where $j=i+n-1$.
    \end{itemize}
\end{thm}

\begin{proof}
First, we reduce cases (ii) and (iii) to case (i) by relabling all the cascades appearing both in $u$ and $v$ as $u_1,\ldots,u_m$ and all the bubbles appearing in both as $(B_\alpha)_{a\in A}$ for some indexing set $A$. Then, as explained before, we can assume without loss of generality that the divisors $D_1, \ldots, D_n$ do not intersect any of the non-constant asymptotes of $u_1, \ldots, u_m$, and do not contain the images of $u_1, \ldots, u_m$ and $\{B_\alpha\}_{\alpha \in A}$. Moreover, Proposition \ref{prop: contained components ball} guarantees that all cascades and bubbles contained in $D_\infty$ are constant and hence can be safely ignored in the proof.

We treat each case separately:

    \medskip
    \noindent\textbf{Case 1: }Let $0\leq i\leq \ell+n$ and $0\leq j\leq \ell'+n$, and assume that $(q_-,q_+)=(\check{x}^\ell_i,\check{y}^{\ell'}_j)$.
    \medskip

    As computed in Section~\ref{ss: computations of FMB}, both $\mu_{FMB}(\check{x}_i^\ell;H_\ell)$ and $\mu_{FMB}(\check{y}_j^{\ell'};H'_{\ell'})$ are even modulo $2(n+1)$. In particular, the difference 
   $$\mu_{FMB}(\check{y}_j^{\ell'};H'_{\ell'}) - \mu_{FMB}(\hat{x}_i^\ell;H_\ell) $$
    is even, and therefore not equal to $-1$ modulo $2(n+1)$. Proposition~\ref{prop: traj and modular indices} then leads to a contradiction.

    \medskip
    \noindent\textbf{Case 2: } Let $1\leq i\leq \ell$ and $1\leq j\leq \ell'$, and assume that $(q_-,q_+)=(\hat{x}^\ell_i,\hat{y}^{\ell'}_j)$.
    \medskip
    
    As computed in Section~\ref{ss: computations of FMB}, both $\mu_{FMB}(\hat{x}_i^\ell;H_\ell)$ and $\mu_{FMB}(\hat{y}_j^{\ell'};H'_{\ell'})$ are odd modulo $2(n+1)$. In particular, the difference 
    $$\mu_{FMB}(\hat{y}_j^{\ell'};H'_{\ell'}) - \mu_{FMB}(\hat{x}_i^\ell;H_\ell) $$
    is even, and therefore not equal to $-1$ modulo $2(n+1)$. Proposition~\ref{prop: traj and modular indices} then leads to a contradiction.

    \medskip
    \noindent\textbf{Case 3: }Let $0\leq i\leq \ell+n$ and $1\leq j\leq \ell'$, and assume that $(q_-,q_+)=(\check{x}^\ell_i,\hat{y}^{\ell'}_j)$.
    \medskip

    Let us split into two cases:
    \begin{enumerate}
        \item Assume that $0\leq i\leq \ell$. Let us use the trivialization $\tau_T$. 

    With respect to this trivialization, we have $\mu_{FMB}^\tau(\hat{y}_j^{\ell'};H'_{\ell'})=2n-1$ and $\mu_{FMB}^\tau(\check{x}_i^\ell;H_\ell) = 0$. By the positivity of intersections from Remark~\ref{rem: positivity of intersections for our acc. data} and Example~\ref{exam: computations using Seidel's lemma} we get that
    $$ \sum_{k=1}^m c_1^\tau(u_k) + \sum_{\alpha\in A} c_1(B_\alpha) = u\cdot (D_1\cup\ldots D_n\cup D_\infty)\geq0. $$

    It follows that
    \begin{align*}
        -1 &= \mu^\tau_{FMB}(q_+;H'_{\ell'}) - \mu^\tau_{FMB}(q_-;H_\ell) + 2\sum_{k=1}^m c_1^\tau(u_k) + 2\sum_{\alpha\in A} c_1(B_\alpha)\\
      &=2n-1-0+2u\cdot (D_1\cup\ldots     D_n\cup D_\infty)\\
      &\geq 2n-1 \geq 1,
    \end{align*}
   
and this is a contradiction.

 \item Assume that $\ell+1\leq i\leq \ell+n$. Let us use the trivialization $\tau_T$. 

    With respect to this trivialization, we have  $\mu_{FMB}^\tau(\hat{y}_j^{\ell'};H'_{\ell'})=2n-1$ and $\mu_{FMB}^\tau(\check{x}_i^\ell;H_\ell) = 2i$. By the positivity of intersections from Remark~\ref{rem: positivity of intersections for our acc. data} and Example~\ref{exam: computations using Seidel's lemma} we get that
$$ \sum_{k=1}^m c_1^\tau(u_k) + \sum_{\alpha\in A} c_1(B_\alpha) = u\cdot (D_1\cup\ldots D_n\cup D_\infty)\geq u\cdot D_\infty \geq \ell+1. $$

It follows that
\begin{align*}
    -1 &= \mu^\tau_{FMB}(q_+;H'_{\ell'}) - \mu^\tau_{FMB}(q_-;H_\ell) + 2\sum_{k=1}^m c_1^\tau(u_k) + 2\sum_{\alpha\in A} c_1(B_\alpha)\\
    &=2n-1-2i+2u\cdot (D_1\cup\ldots D_n\cup D_\infty)\\
    &\geq 2n-1-2i+2(\ell+1)\\
    &=2(\ell+n-i)+1,
\end{align*}
and since $i\leq \ell+n$, we deduce that $-1\geq2(\ell+n-i)+ 1\geq1$, which is a contradiction.

    \end{enumerate}
    \medskip
    \noindent\textbf{Case 4: }Let $1\leq i\leq \ell$ and $0\leq j\leq \ell'+n$, and assume that $(q_-,q_+)=(\hat{x}^\ell_i,\check{y}^{\ell'}_j)$.
    \medskip
    
    There are three cases that we have to check:
    \begin{enumerate}
        \item Assume that $j=0$. Let us use the trivialization $\tau_T$. 

    With respect to this trivialization, we have  $\mu_{FMB}^\tau(\check{y}_j^{\ell'};H'_{\ell'})=0$ and $\mu_{FMB}^\tau(\hat{x}_i^\ell;H_\ell) = 2n-1$. By the positivity of intersections from Remark~\ref{rem: positivity of intersections for our acc. data} and Example~\ref{exam: computations using Seidel's lemma} we get that
$$ \sum_{k=1}^m c_1^\tau(u_k) + \sum_{\alpha\in A} c_1(B_\alpha) = u\cdot (D_1\cup\ldots D_n\cup D_\infty)\geq u\cdot (D_1\cup\ldots\cup D_n) \geq n\cdot 1=n. $$

It follows that
\begin{align*}
    -1 &= \mu^\tau_{FMB}(q_+;H'_{\ell'}) - \mu^\tau_{FMB}(q_-;H_\ell) + 2\sum_{k=1}^m c_1^\tau(u_k) + 2\sum_{\alpha\in A} c_1(B_\alpha)\\
    &=0-2n+1+2u\cdot (D_1\cup\ldots D_n\cup D_\infty)\\
    &\geq -2n+1+2n=1,
\end{align*}
which is a contradiction.

    \item Assume that $1\leq j\leq \ell'$. Let us use the trivialization $\tau_T$.

With respect to this trivialization, we have that $\mu_{FMB}^\tau(\check{y}_j^{\ell'};H'_{\ell'})=0$ and $\mu_{FMB}^\tau(\hat{x}_i^\ell;H_\ell)=2n-1$. By Theorem~\ref{thm: c1 vs intersection number}, we get that
$$
\sum_{k=1}^m c_1^\tau(u_k)+\sum_{\alpha\in A}c_1(B_\alpha)
=u\cdot(D_1\cup\ldots\cup D_n\cup D_\infty).
$$
It follows that
\begin{align*}
    -1 &= \mu^\tau_{FMB}(q_+;H'_{\ell'})
    -\mu^\tau_{FMB}(q_-;H_\ell)
    +2\sum_{k=1}^m c_1^\tau(u_k)
    +2\sum_{\alpha\in A}c_1(B_\alpha)\\
    &=0-2n+1+2u\cdot(D_1\cup\ldots\cup D_n\cup D_\infty),
\end{align*}
which means that
$$
u\cdot(D_1\cup\ldots\cup D_n\cup D_\infty)=n-1.
$$
By the positivity of intersections from Remark~\ref{rem: positivity of intersections for our acc. data}, we deduce that
$$
0\leq u\cdot D_\infty
\leq u\cdot(D_1\cup\ldots\cup D_n\cup D_\infty)=n-1.
$$

Next we move to the trivialization $\tau_B$.

With respect to this trivialization, we have that $\mu_{FMB}^\tau(\check{y}_j^{\ell'};H'_{\ell'})=-2nj$ and $\mu_{FMB}^\tau(\hat{x}_i^\ell;H_\ell)=-2ni+2n-1$. By Theorem~\ref{thm: c1 vs intersection number}, we deduce that
$$
\sum_{k=1}^m c_1^\tau(u_k)+\sum_{\alpha\in A}c_1(B_\alpha)
=(n+1)u\cdot D_\infty.
$$
It follows that
\begin{align*}
    -1 &= \mu^\tau_{FMB}(q_+;H'_{\ell'})
    -\mu^\tau_{FMB}(q_-;H_\ell)
    +2\sum_{k=1}^m c_1^\tau(u_k)
    +2\sum_{\alpha\in A}c_1(B_\alpha)\\
    &=-2nj+2ni-2n+1+2(n+1)u\cdot D_\infty,
\end{align*}
which means that
$$
u\cdot D_\infty=\frac{n(j-i+1)-1}{n+1},
$$
and hence
\begin{align*}
    0&\equiv n(j-i+1)-1 \pmod{n+1}\\
     &\equiv -(j-i+1)-1 \pmod{n+1}\\
     &\equiv -j+i-2 \pmod{n+1}.
\end{align*}
This implies that there exists $r\in\Z$ such that $j-i+2=(n+1)r$. Hence
\[
    u\cdot D_\infty
    =\frac{n(j-i+1)-1}{n+1}
    =\frac{n((n+1)r-1)-1}{n+1}
    =nr-1.
\]
Since we previously proved that $0\leq u\cdot D_\infty\leq n-1$, we get that
$$
0\leq nr-1\leq n-1,
$$
and therefore $1\leq nr\leq n$. Since $r\in\Z$ and $n\geq1$, we deduce that $r=1$, which implies that
$$
j-i+2=n+1,
$$
and hence $j=i+n-1$, as required.

\item Assume that $\ell'+1\leq j\leq \ell'+n$. Let us use the trivialization $\tau_B$. 

    With respect to this trivialization, we still have that $\mu_{FMB}^\tau(\check{y}_j^{\ell'};H'_{\ell'}) = 2j$ and $\mu_{FMB}^\tau(\hat{x}_i^\ell;H_\ell) = -2ni+2n-1$.
It follows that
\begin{align*}
    -1 &= \mu^\tau_{FMB}(q_+;H'_{\ell'}) - \mu^\tau_{FMB}(q_-;H_\ell) + 2\sum_{k=1}^m c_1^\tau(u_k) + 2\sum_{\alpha\in A} c_1(B_\alpha)\\
    &=2j+2ni-2n+1+2(n+1) u\cdot D_\infty,
\end{align*}
which means that
$$u\cdot D_\infty=\frac{-j-ni+n-1}{n+1},$$
and hence
\begin{align*}
    0&\equiv -j-ni+n-1 \pmod{n+1}\\
    &\equiv -j-(n+1-1)i+n+1-2 \pmod{n+1}\\
    &\equiv -j+i-2 \pmod{n+1}.
\end{align*}
This implies that there exists $r\in \Z$ such that $i=j+2+(n+1)r$. Hence
\begin{align*}
    u\cdot D_\infty&=\frac{-j-ni+n-1}{n+1}\\
    &=\frac{-j-n(j+2+(n+1)r)+n-1}{n+1}\\
    &=\frac{-j-nj-2n-n(n+1)r+n-1}{n+1}\\
    &=\frac{-j(n+1)-n(n+1)r-n-1}{n+1}\\
    &=-j-nr-1.
\end{align*}
 By Example~\ref{exam: computations using Seidel's lemma} we know that
    $$u\cdot D_\infty\geq-\ell',$$
and hence 
$$j+nr+1\leq \ell'$$
therefore
$$j\leq \ell'-nr-1.$$
By assumption, in this case, $\ell'+1\leq j$, thus $\ell' +1 \le\ell'-nr-1$, which implies that
$nr\leq -2$, so either $r\leq -2$ or $r=-1$ and $n\geq 2$. In both cases $r\leq -1$.

Next we move to the trivialization $\tau_T$. 

With respect to this trivialization, we still have that $\mu_{FMB}^\tau(\check{y}_j^{\ell'};H'_{\ell'}) = 2j$ and $\mu_{FMB}^\tau(\hat{x}_i^\ell;H_\ell) = 2n-1$. By Remark~\ref{rem: positivity of intersections for our acc. data} we deduce that  
$$ \sum_{k=1}^m c_1^\tau(u_k) + \sum_{\alpha\in A} c_1(B_\alpha) =u\cdot(D_1\cup \ldots\cup D_n\cup D_\infty)\geq u\cdot D_\infty=-j-nr-1. $$

It follows that
\begin{align*}
    -1 &= \mu^\tau_{FMB}(q_+;H'_{\ell'}) - \mu^\tau_{FMB}(q_-;H_\ell) + 2\sum_{k=1}^m c_1^\tau(u_k) + 2\sum_{\alpha\in A} c_1(B_\alpha)\\
    &=2j-2n+1+2u\cdot (D_1\cup \ldots\cup D_n\cup D_\infty)\\
    &\geq2j-2n+1-2(j+nr+1)\\
    &=-2n-1-2nr
\end{align*}
and thus $nr\geq-n$, and hence $r\geq -1$. Since we previously proved that $r\leq -1$, we deduce that $r=-1$, which implies that $i=j+2-(n+1)=j-n+1$, as required.\qedhere
\end{enumerate}

\end{proof}

In fact, even though Theorem~\ref{thm: htpy_obst} does not imply that the chain homotopy map vanishes for every $n\in \N$, its combination with Theorem~\ref{thm: diff_obst} still implies that the above diagram commutes.

\begin{prop}\label{prop: commuting of cont maps}
For every $n\in \N$, $0<\Delta'\leq \Delta<1$ and $\ell\in \Z_{\geq 0}$ the diagram

$$\xymatrix{
CF(H_\ell) \ar[r] \ar[d] & CF(H_{\ell+1}) \ar[d] \\
CF(H'_\ell) \ar[r]       & CF(H'_{\ell+1})
}$$
commutes, where all the arrows are continuation maps.
    
\end{prop}
\begin{proof}
   The commutativity of the diagram is equivalent to the equality
   \[d'H+Hd=0,\]
   where $H$ is the chain homotopy map induced by the $2$-face of the cube:
   $$H:CF^*(H_\ell)\to CF^{*-1}(H'_{\ell+1}),$$
   and $d,d'$ are the differentials on $CF(H_\ell),CF(H'_{\ell+1})$, respectively.

   To prove this equality, we evaluate $d'H+Hd$ on the generators $\check{x}_0^\ell,\ldots,\check{x}_{\ell+n}^\ell,\hat{x}_1^\ell,\ldots,\hat{x}_\ell^\ell$ of $CF(H_\ell)$ and show it vanishes.

   Theorem~\ref{thm: diff_obst} implies that for every $0\leq i\leq \ell+n$ we have $d\check{x}_i^\ell=0$ and by Theorem~\ref{thm: htpy_obst} we know that $H(\check{x}_i^\ell)=0$. Thus, for every $0\leq i\leq \ell+n$ we get that 
   $$d'H(\check{x}_i^\ell)+Hd(\check{x}_i^\ell)=d'(0)+H(0)=0.$$
   Now, let $1\leq i\leq \ell$. By Theorem~\ref{thm: diff_obst} we deduce that there exist $\lambda_1,\lambda_2\in \Lambda_{\geq0}$ such that $d\hat{x}_i^\ell=\lambda_1\check{x}_{i-1}^\ell+\lambda_2\check{x}_{i+n}^\ell$ and by Theorem~\ref{thm: htpy_obst} we deduce that there exist $\lambda_3,\lambda_4\in \Lambda_{\geq0}$, possibly $0$, such that $H(\hat{x}_i^\ell)=\lambda_3 \check{y}_i^{\ell+1}+\lambda_4 \check{y}_{i+n-1}^{\ell+1}$. Since Theorem~\ref{thm: diff_obst} implies that $d'\check{y}_i^{\ell+1}=d'\check{y}_{i+n-1}^{\ell+1}=0$, and Theorem~\ref{thm: htpy_obst} implies that $H(\check{x}_{i-1}^\ell)=H(\check{x}_{i+n}^\ell)=0$ we get that
   $$d' H(\hat{x}_i^\ell)+Hd(\hat{x}_i^\ell)=d'(\lambda_3 \check{y}_i^{\ell+1}+\lambda_4 \check{y}_{i+n-1}^{\ell+1})+H(\lambda_1\check{x}_{i-1}^\ell+\lambda_2\check{x}_{i+n}^\ell)=0,$$
   as required.
\end{proof}

    \section{Computations of Floer complexes}\label{s: computations of CF}

    The purpose of this section is to compute the Floer cochain complexes over $\Z$ associated with the Hamiltonians in the acceleration data introduced in Section~\ref{ss: acc. data}, together with the continuation maps between them. \label{p: transversality paragraph}To achieve transversality for the Floer differentials, we perturb the almost complex structures, keeping them equal to $J_{std}$ on a fixed shell outside $B_\Delta$, so that the no-escape lemma remains applicable. For the continuation maps, we perturb the Hamiltonian homotopies and the associated families of almost complex structures, keeping them linear monotone increasing on a on a fixed shell outside $B_\Delta$, together with keeping there the almost complex structure $J_{std}$, so that the no-escape lemma remains applicable. Sections~\ref{ss: transversality} and~\ref{ss:tower} explain how to choose these perturbations coherently and sufficiently close to the original data for the obstruction statements of Section~\ref{s: obstructions} to persist.
    
    \textbf{Note:}  To streamline the text and make it easier to follow, Throughout this section, references to the obstruction theorems of Section~\ref{s: obstructions} also refer to their extensions to the chosen regular perturbations, as justified in Sections~\ref{ss: transversality} and~\ref{ss:tower}.

Let $\Delta \in (0,1)$. For every $\ell \in \Z_{\geq 0}$, we will compute the Floer cochain complex of the Hamiltonian $H_\ell \fc \CP^n \to \R$ which is given by $H_\ell(z) = h(\Delta, \ell, \mu(z))$ for every $z \in \CP^n$.

Fix $\ell \in \Z_{\geq 0}$. For every $1 \leq i \leq \ell$, let $r_{\ell,i}^\Delta \in [0,1)$ denote the number satisfying $\frac{\partial}{\partial r}|_{r=r_{\ell,i}^\Delta} h(\Delta, \ell, r) = i$. also, put $r_{\ell,0}=0$. For every $0\leq i\leq \ell$ let us denote $h_{\ell,i}^\Delta=h(\Delta,\ell,r_{\ell,i}^\Delta)$, additionally, for every $1 \leq i \leq \ell$, we define the $(2n-1)$-spheres
$$S_{\ell,i} = \{z \in \CP^n : h'_\ell(\mu(z)) = i\},$$
and equip each with a perfect Morse function. Furthermore, we equip the divisor $D_\infty = \CP^{n-1}$ with a perfect Morse function. 

Let $\check{x}_0^\ell$ denote the constant orbit of $H_\ell$ at $0 \in \Int B(1) \subset \CP^n$. For every $1 \leq i \leq \ell$, let $\check{x}^\ell_i$ and $\hat{x}^\ell_i$ denote the minimum and maximum points on $S_{\ell,i}$, respectively. Finally, for every $\ell+1 \leq i \leq \ell + n$, let $\check{x}_i^\ell$ denote the critical point on $D_\infty$ of Morse index $2(i-\ell-1)$.
\begin{thm}\label{thm: CF(H_l;Z)}
    For every $\ell\in \Z_{\geq0}$ there exist signs $A_{\ell,1},\ldots,A_{\ell,\ell},B_{\ell,1},\ldots,B_{\ell,\ell}\in \{-1,1\}$,

    such that the differential $d\fc CF(H_\ell;\Z)\to CF(H_\ell;\Z)$ satisfies
    \newlength{\condwidth}
    \settowidth{\condwidth}{For every $0\leq i\leq \ell+n:$}
    \begin{itemize}
        \item \makebox[\condwidth][l]{For every $1\leq i\leq \ell:$}\qquad $d \hat{x}^\ell_i=A_{\ell,i}\check{x}^\ell_{i-1}+B_{\ell,i}\check{x}^\ell_{i+n}$.
     
    \item \makebox[\condwidth][l]{For every $0\leq i\leq \ell+n:$}\qquad $d\check{x}^\ell_i=0$.
    \end{itemize}
  
\end{thm}

\begin{figure}[hbtp]
    \centering
    \includegraphics[width=\linewidth]{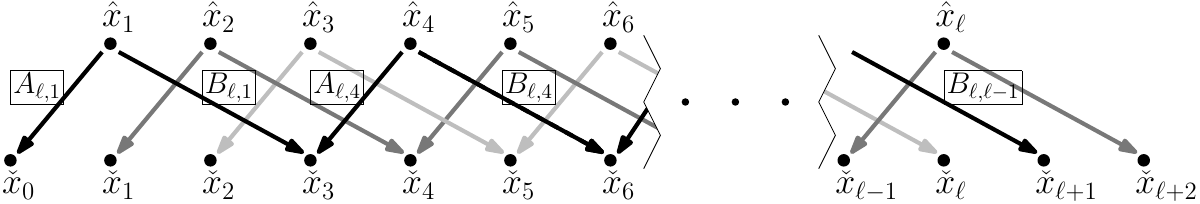}
    \caption{The complex of $H_\ell$ on $\CP^2$. The different shades of gray are used as a visual cue to separate the subcomplexes. One of them is labeled with the signs $A$ and $B$.}
    \label{fig:figfullcomplexproof}
\end{figure}

Employing similar techniques we also compute the continuation maps between the Floer cochain complexes of the acceleration data of two different balls: Let $0<\Delta'\leq\Delta<1$. Consider additionally the corresponding acceleration data $(H'_\ell)_{\ell\geq0}$ for the ball $B^{}_{\Delta'}=\mu^{-1}([0,\Delta'])$, given by $H'_\ell(z)=h(\Delta',\ell,\mu(z))$, for every $z\in \CP^n$ and $\ell\in \Z_{\geq0}$. For every $\ell\in \Z_{\geq0}$ denote by
$\check{y}_0^\ell,\ldots,\check{y}_{\ell+n}^\ell,\hat{y}_1^\ell,\ldots,\hat{y}_\ell^\ell$ the generators for $CF(H'_\ell)$, as in Theorem~\ref{thm: CF(H_l;Z)}.

Let $\ell,\ell'\in \Z_{\geq0}$ and assume that $\ell\leq \ell'$. Let $\chi\fc \R\to \R$ be a smooth non-increasing function, satisfying $\chi(s)=\Delta$ for every $s<0$ and $\chi(s)=\Delta'$ for every $s>1$. Similarly, let $\lambda\fc \R\to \R$ be a smooth non-decreasing function, satisfying $\lambda(s)=\ell$ for every $s<0$ and $\lambda(s)=\ell'$ for every $s>1$. 

Now, define a monotone homotopy $K\fc \R\times \CP^n\to \R$ from $H_\ell$ to $H'_{\ell'}$ by $K(s,z)=h(\chi(s),\lambda(s),\mu(z))$ for every $(s,z)\in \R\times \CP^n$. The homotopy $K$ defines a continuation map $\Phi\fc CF(H_\ell;\Z)\to CF(H'_{\ell'};\Z)$.

\begin{thm}\label{thm: continuation maps, over Z}
    For every $0<\Delta'\leq\Delta<1$ and $\ell,\ell'\in \Z_{\geq0}$ with $\ell\leq \ell'$ there exist signs
    $$\check{C}_{0},\ldots,\check{C}_{\ell+n},\hat{C}_{1},\ldots,\hat{C}_{\ell}\in \{-1,1\}$$
    such that the continuation map $\Phi\fc CF(H_\ell;\Z)\to CF(H'_{\ell'};\Z)$ satisfies
    
    \newlength{\condwidthB}
    \settowidth{\condwidthB}{For every $0\leq i\leq \ell+n:$}
   
    \begin{itemize}
        \item \makebox[\condwidthB][l]{For every $0\leq i\leq \ell+n:$}\qquad  $\Phi(\check{x}_i^\ell)=\check{C}_{i} \check{y}_i^{\ell'}$.

         \item \makebox[\condwidthB][l]{For every $1\leq i\leq \ell:$}\qquad $\Phi(\hat{x}_i^\ell)=\hat{C}_{i} \hat{y}_i^{\ell'}$.
    \end{itemize}
\end{thm}

Our approach for the proof of these results combines the obstructions on Floer solutions and continuation maps we discussed in Section~\ref{s: obstructions} with well known ideas from the computation of the Floer complex of $J$-shaped Hamiltonians on $\C^n$, such as those \cite[Section 3]{Oancea_survey}. To apply these computations for our case in $\CP^n$, we first discuss in Section~\ref{ss: no_escape} the no escape lemma, and in Section~\ref{ss: CF of J-shaped} we apply it and adapt the ideas from \cite[Section 3]{Oancea_survey} to $\CP^n$ to compute some matrix elements of the Floer differential. The proof of Theorem~\ref{thm: CF(H_l;Z)} and Theorem~\ref{thm: continuation maps, over Z} appear in Section~\ref{ss: proof of computations of CF + continuations}, using algebraic properties of the cochain complexes and continuation maps to deduce the missing matrix elements.

\subsection{No escape lemma and Floer theory on open balls}\label{ss: no_escape}

In this section, we survey a construction of a robust Floer theory for Hamiltonians on $\Int B(1)$ that are linear with respect to the Liouville coordinate in the complement of a compact set and are equipped with a special almost complex structure. This is a variation on the classical Floer theory in completed Liouville domains. The construction is based on what is known as the ``no escape lemma.'' We present a version of this lemma adapted from \cite{Ritter_TQFT}, which allows us to show that the algebraic count of Floer trajectories between generators is well-defined, provided that these trajectories do not escape from a certain compact subset. Consequently, this framework enables us to establish a well-defined Floer theory for such Hamiltonians.

We start with the formulation of the no escape lemma. Let $(V, \omega)$ be an exact, not necessarily compact, symplectic manifold with boundary $\partial V$, and let $\theta$ be a primitive $1$-form for $\omega$. We assume that $\partial V$ is of negative contact type: the Liouville field $Y$, defined by $\iota_Y \omega = \theta$, points strictly inwards. We denote the flow of $Y$ by $(\varphi_Y^t)_t$, which is well-defined for sufficiently small $t$. In a sufficiently small neighborhood $U$ of $\partial V$, we define the Liouville coordinate $R$ such that $\partial V = \{z \in V \,:
\, R(z) = R_0\}$ for some $R_0 \in \R$. Specifically, for every $z \in U$, let $r(z) \in \R$ be the unique value such that $\varphi_Y^{-r(z)}(z) \in \partial V$. The Liouville coordinate $R$ at $z$ is then defined as $R(z) = R_0 e^{r(z)}$.

Let $J$ be a compatible almost complex structure on $(V,\omega)$. Define the Reeb vector field by $\cR = JY$.

 We consider the cascade $\Sigma=\R\times S^1$ equipped with the complex coordinate $s+it$, for $s\in \R$ and $t\in S^1$. $\Sigma$ is equipped with the closed $1$-form $\beta=dt$. Denote by $j$ the standard complex structure on $\Sigma$.

 We present here a version of Lemma D.3. from \cite{Ritter_TQFT}.

 \begin{lemma}\label{lemma: no escape}
 Let $S$ be a compact submanifold with boundary of $\Sigma$ of codimension $0$, with smooth boundary. Let $H\fc S\times V\to \R$ be a domain-dependent Hamiltonian and assume that there exists a neighborhood $U$ of $\partial V$ and smooth functions $m\fc \R\to \R\setminus\{0\}$ and $b\fc \R\to \R$ such that for every $((s,t),z)\in S\times V$ we have $H(s,t,z)=m(s)R(z)+b(s)$. Additionally, we assume that $m'(s)\geq0$ for every $s\in \R$. Let $u\fc S\to V$ be a solution of the Floer equation $(du-X_{H}\otimes \beta)^{0,1}=0$ satisfying $u(\partial S)\subset \partial V$. If there exists a $c>0$ such that along $\partial V$ we have $J^*\theta=cdR$, then $u(S)\subset \partial V$ and $du = X \otimes \beta$.
 \end{lemma}

\begin{proof}
Without loss of generality, we assume that $b(s)=0$ for every $s\in \R$, since this does not change $X_H$, therefore it does not change the Floer equation. By the assumption, along $\partial V$ we have 
$$\iota_{X_H}\omega=dH=m(s)dR=\frac{m(s)}{c}J^*\theta=\frac{m(s)}{c}\theta(J\cdot)=\frac{m(s)}{c}\omega(Y,J\cdot).$$
Since $J$ is $\omega$-compatible we get that $$\omega(Y,J\cdot)=\omega(JY,J^2\cdot)=-\omega(JY,\cdot)=-\omega(\cR,\cdot).$$
Thus
$$\omega(X_H,\cdot)=\omega\left(-\frac{1}{c}m(s)\cR,\cdot\right),$$
which implies that $X_{H}=-\frac{1}{c}m(s)\cR=-\frac{1}{c}m(s)JY$.

Let us calculate $\theta(X_H)$ on $\partial V$:
$$\theta(X_H)=-\frac{m(s)}{c}\theta(JY)=-\frac{m(s)}{c}J^*\theta(Y)=-\frac{m(s)}{c}cdR(Y)=-m(s)R=-H.$$
We will now show that $E(u) \le 0$, therefore $E(u) = 0$, hence $du = X \otimes \beta$, so $du$ belongs in the span of $X \in \spn\{\cR\} \subset T\partial V$. Thus, $u(S) \subset \partial V$ as required. Define $\tilde{U}\fc S\to S\times V$ by $\tilde{u}(s,t)=(s,t,u(s,t))$ for every $(s,t)\in S$.

\begin{align*}
E(u) &= \int_S u^*d\theta + \tilde{u}^*(dH) \wedge \beta \quad (\text{Energy identity,\, } \omega = d\theta) \\
&= \int_S u^*d\theta + (d(\tilde{u}^*H)-\tilde{u}^*\partial_sHds-\tilde{u}^*\partial_tHdt)\wedge \beta \quad  \\
&\leq \int_S u^*d\theta + d(\tilde{u}^*H)\wedge \beta \quad  (dt\wedge \beta=0 \text{ and }\partial_s H\geq0)\\
&= \int_S d(\tilde{u}^*\theta) + d(u^*H\beta)-\tilde{u}^*Hd\beta \\
&= \int_S d(\tilde{u}^*\theta) + d(u^*H\beta) \quad (d\beta =0) \\
&= \int_{\partial S} u^*\theta + (\tilde{u}^*H)\beta \quad (\text{Stokes' Theorem}) \\
&= \int_{\partial S} u^*\theta - \theta(X_H)\beta \quad (\text{on } u(\partial S) \subset \partial V : H =-\theta(X_H)) \\
&= \int_{\partial S} \theta(du - X _H\otimes \beta) \\
&= \int_{\partial S} -\theta J(du - X_H \otimes \beta)j \quad (\text{since } (du - X_H \otimes \beta)^{0,1} = 0) \\
&= \int_{\partial S} -cdR(du - X_H \otimes \beta)j \quad (J^*\theta=cdR \text{ along } u(\partial S) \subset \partial V) \\
&= \int_{\partial S} -cdR(du)j \quad (dR=\frac{1}{m(s)}dH \text{ vanishes on } X_H \text{ on } u(\partial S) \subset \partial V)
\end{align*}

Let $\hat{n} =$ outward normal direction along $\partial S \subset S$. Then $(\hat{n}, j\hat{n})$ is an oriented frame, so $\partial S$ is oriented by $j\hat{n}$. Now $dR(du)j(j\hat{n}) = d(R \circ u) \cdot (-\hat{n}) \ge 0$ since in the inward direction $-\hat{n}$, $R \circ u$ can only increase since $\partial V$ minimizes $R$. So $E(u) \le 0$. 
\end{proof}

Recall that Claim~\ref{claim: Int B(1)->CPn-D_infty} asserts that the map $\iota\fc (\Int B(1),\omega_0)\to(\CP^n,\omega_{FS})$, given by
$$\iota(z)=[\sqrt{1/\pi - \|z\|^2} : z],$$
for every $z\in \Int B(1)$, where $B(1)=\{z\in \C^n\,|\,\pi\|z\|^2\leq 1\}$, is a symplectic embedding and its image is $\CP^n\setminus D_\infty$. 

The $1$-form $\theta$ on $\Int B(1)\subset \R^{2n}$ given by
$$\theta=\frac{1}{2}\sum_{j=1}^n (x_j d y_j-y_j dx_j)=\frac{1}{2}\text{Im}\sum_{j=1}^n \bar{z}_j dz_j,$$
is a primitive of $\omega_0=\iota^*\omega_{FS}$. The vector field $Y$ on $\Int B(1)$, given by $Y(z)=\frac{1}{2}z$ for every $z\in \Int B(1)$, is a Liouville vector field with respect to the primitive, i.e. $\omega_0(Y,\cdot)=\theta$. The flow $(\varphi_Y^t)_t$ is locally defined and is given by $\varphi_Y^t(z)=e^{t/2}z$ for every $z\in \Int B(1)$, and $t\in(-\infty,-\log(\pi\|z\|^2))$.

Given $R_0\in(0,1)$, the Liouville coordinate $R\fc \Int B(1)\setminus\{0\}\to \R$ with respect to the sphere $S_{R_0}=\{z\in\C^n\,|\,\pi\|z\|^2=R_0\}$ is defined as follows: for every $z\in \Int B(1)\setminus\{0\}$ denote by $r(z)\in \R$ the value for which $\varphi_Y^{-r(z)}(z)\in S_{R_0}$, and then take $R(z)=R_0 e^{r(z)}$. Note that 
$$R_0=\pi\|\varphi_Y^{-r(z)}(z)\|^2=\pi\|e^{-r(z)/2}z\|^2=\pi\|z\|^2 e^{-r(z)},$$
which means that
$$R(z)=R_0e^{r(z)}=R_0\frac{\pi\|z\|^2}{R_0}=\pi\|z\|^2.$$
Thus $R$ smoothly extends to all of $\Int B(1)$, and satisfies
$$dR=d(\pi\|z\|^2)=\pi d(z\cdot \bar z)=\pi(zd\bar{z}+\bar{z}dz)=2\pi \text{Re}(\bar{z}dz)=2\pi \text{Re}(\sum_{j=1}^n \bar{z}_j dz_j).$$
Therefore,
$$dR=2\pi \text{Re}(\langle z,\cdot)\rangle.$$

By Claim~\ref{claim: iota^* J} if $J$ is the standard complex structure on $\CP^n$, then $\iota^*J$ is given by
    $$ (\iota^* J)_z(V) = iV - \pi \text{Im}(\langle z, V \rangle) z + \pi \frac{\text{Re}(\langle z, V \rangle)}{1 - \pi\|z\|^2}iz, $$
    for every $z \in \Int B(1)$ and $V \in \C^n$, where $\langle \cdot,\cdot \rangle$ is the standard Hermitian inner product on $\C^n$.

Note that for every $z\in \Int B(1)$ and $V\in \C^n$, if we decompose $V$ as $V=U+az+biz$ for $a,b\in \R$ and $U\in \C^n$ with $\langle U,z\rangle=0$ we get that

\begin{align*}
    (\iota^* J)^*\theta(V)&=\theta(\iota^*J V)\\
    &=\theta\left(  iV - \pi \text{Im}(\langle z, V \rangle) z + \pi \frac{\text{Re}(\langle z, V \rangle)}{1 - \pi\|z\|^2}iz\right)\\
     &=\theta\left(  iV - \pi \text{Im}(\langle z, U+az+biz \rangle) z + \pi \frac{\text{Re}(\langle z, U+az+biz \rangle)}{1 - \pi\|z\|^2}iz\right)\\
     &=\theta\left(  i(U+az+biz) - \pi \text{Im}(a\|z\|^2+bi\|z\|^2) z + \pi \frac{\text{Re}(a\|z\|^2+bi\|z\|^2)}{1 - \pi\|z\|^2}iz\right)\\
     &=\theta\left(  iU+aiz-bz - \pi b\|z\|^2 z + \pi \frac{a\|z\|^2}{1 - \pi\|z\|^2}iz\right)\\
      &=\theta\left(  iU-b(1+ \pi \|z\|^2) z +  \frac{a}{1 - \pi\|z\|^2}iz\right)\\
      &=\frac{1}{2}\text{Im}\left(\left(\sum_{j=1}^n \bar{z}_j dz_j\right)\left(  iU-b(1+ \pi \|z\|^2) z +  \frac{a}{1 - \pi\|z\|^2}iz\right)\right)\\
      &=\frac{1}{2}\text{Im}\left(i\langle U,z\rangle-b(1+ \pi \|z\|^2) \|z\|^2 +  \frac{a}{1 - \pi\|z\|^2}i\|z\|^2\right)\\
      &=\frac{1}{2}\cdot\frac{a\|z\|^2}{1 - \pi\|z\|^2},
\end{align*}
and since $a\|z\|^2=\text{Re}(\langle z,V\rangle)$, we conclude that
$$(\iota^* J)^*\theta(V)=\frac{1}{2}\cdot\frac{\text{Re}(\langle z,V\rangle)}{1 - \pi\|z\|^2},$$
and thus
$$(\iota^* J)^*\theta=\frac{1}{2}\cdot\frac{\text{Re}(\langle z,\cdot\rangle)}{1 - R}=\frac{1}{4\pi(1-R)}dR.$$

This shows that for every $R_0 \in (0,1)$, we can apply Lemma~\ref{lemma: no escape} to the exact symplectic manifold $V_{R_0} = R^{-1}([R_0,1)) = \{z \in \Int B(1) : R_0 \leq \pi\|z\|^2 < 1\}$. This allows us to define the Floer complex for \textit{admissible} Hamiltonians on $\Int B(1)$.

Let $H$ be a Hamiltonian on $\Int B(1)$ that is either nondegenerate or autonomous and satisfies the \textbf{MB} condition. We say that $H$ is \textbf{admissible} if there exist $R_0 \in (0,1)$, $m \in \R \setminus \Z$ and $b\in \R$ such that for every $z \in \Int B(1)$, $R(z) \geq R_0$ implies $H(z) = mR(z)+b$. In this case the constant $m$ is called the \textbf{slope} of $H$. The assumption $m \notin \Z$ guarantees that all $1$-periodic orbits of $H$ are contained in$$\Int B(1) \setminus V_{R_0} = R^{-1}([0,R_0)) = \{z \in \Int B(1) : \pi\|z\|^2 < R_0\}.$$
By Lemma~\ref{lemma: no escape}, all Floer trajectories connecting the $1$-periodic orbits of $H$ must be contained in \[\overline{\Int B(1) \setminus V_{R_0}} = R^{-1}([0,R_0]).\] Therefore, by standard compactness techniques, the Floer complex $CF(H;\Z)$ of $H$ is well-defined. The homology of $CF(H;\Z)$ is called the Floer homology of $H$ and it is denoted by $HF(H;\Z)$.

    Let $H,H'$ be admissible Hamiltonians on $\Int B(1)$, and assume that $H \leq H'$. A homotopy $(H_s)_{s \in \R}$ from $H$ to $H'$ is called \textbf{admissible} if it satisfies the following conditions:
\begin{enumerate}
    \item For every $s \leq 0$, we have $H_s = H$, and for every $s \geq 1$, we have $H_s = H'$.
    \item $\partial_s H_s \geq 0$.
    \item There exists a non-decreasing functions $m \fc \R \to \R \setminus \{0\}$, $b\fc \R\to \R$ and a constant $R_0 \in (0,1)$ such that for every $s \in \R$ and $z \in \Int B(1)$, the condition $R(z) \geq R_0$ implies $H_s(z) = m(s) R(z)+b(s)$.
\end{enumerate}

Note that given two admissible Hamiltonians $H,H'$ on $\Int B(1)$, if $H \leq H'$ and the slope of $H$ is less than or equal to the slope of $H'$, then there exists an admissible homotopy between them.

By Lemma~\ref{lemma: no escape}, given two admissible Hamiltonians $H,H'$ with $H \leq H'$ and an admissible homotopy between them, all continuation trajectories connecting the $1$-periodic orbits of $H$ to those of $H'$ must be contained in a compact set. Therefore, the continuation map $CF(H;\Z) \to CF(H';\Z)$ is well-defined and the induced map in homology is independent of the choice of admissible homotopy.

Additionally, given three admissible Hamiltonians $H \leq H' \leq H''$ such that the slope of $H$ is less than or equal to that of $H'$, which in turn is less than or equal to that of $H''$, the continuation map $CF(H;\Z) \to CF(H'';\Z)$ is chain homotopic to the composition of the continuation maps 
$$CF(H;\Z) \to CF(H';\Z) \to CF(H'';\Z),$$
hence they induce the same map on the homology level: $HF(H;\Z) \to HF(H'';\Z)$.

Note that for every admissible Hamiltonian $H$ and a constant $C > 0$, the Hamiltonian $H+C$ is admissible and satisfies $CF(H+C;\Z) = CF(H;\Z)$. Since $H$ and $H+C$ have the same slope, there exists an admissible homotopy from $H$ to $H+C$, and the induced continuation map is the identity.

Let us formulate the following conclusion of the above discussion as a claim:

\begin{claim}\label{claim: qi for Ham with the same slope}
    Let $H,H'$ be admissible Hamiltonians on $\Int B(1)$ with $H \leq H'$. If $H$ and $H'$ have the same slope, then the continuation map $CF(H;\Z) \to CF(H';\Z)$ is a quasi-isomorphism.
\end{claim}

\begin{proof}
    Since $\Int B(1)$ is bounded in $\C^n$ and $H,H'$ are admissible, there exists $C > 0$ such that $H' \leq H+C$. Since there exist admissible homotopies from $H$ to $H'$ and from $H'$ to $H+C$, the identity map $\id \fc HF(H;\Z) \to HF(H+C;\Z) = HF(H;\Z)$ is equal to the composition:
    $$HF(H;\Z) \to HF(H';\Z) \to HF(H+C;\Z) = HF(H;\Z).$$
    This implies that the continuation map $HF(H;\Z) \to HF(H';\Z)$ is injective.
    Now we note that  $H+C\le H'+C$, hence, applying the same argument with $H' \le H+C \le H' + C$ yields that the map above is also surjective, hence an isomorphism; in other words, the continuation map at the chain level $CF(H;\Z) \to CF(H';\Z)$ is a quasi-isomorphism, as required.
\end{proof}

\subsection{Computation of Floer complexes of admissible Hamiltonians}\label{ss: CF of J-shaped}

In this section, we adapt a computation from \cite[Section 3]{Oancea_survey} concerning the Floer complex of ``$J$-shaped'' Hamiltonians on $\C^n$ with the standard complex structure, to the setting of $(\Int B(1), \iota^*J)\cong(\CP^n\setminus D_\infty,J)$.

Let $\ell \in \Z_{\geq 0}$, and define $\tilde{H}_\ell \fc \Int B(1) \to \R$ by 
$$ \tilde{H}_\ell(z) = H_\ell \circ \iota(z) = h(\Delta, \ell, \pi\|z\|^2), $$
for every $z \in \Int B(1)$, where $H_\ell$ is the $\ell$-th function in the acceleration data defined in Section~\ref{ss: acc. data}. Moreover, the Hamiltonian $\tilde{H}_\ell$ is autonomous and satisfies the \textbf{MB} condition. Its critical submanifolds are $\{0\}$ and the spheres 
$$ S_j = \{z \in \Int B(1) : h_\ell'(\pi\|z\|^2) = j\}, $$
for every $1 \leq j \leq \ell$, where $h_\ell \fc [0,1] \to \R$ is given by $h_\ell(r) = h(\Delta, \ell, r)$ for every $r \in [0,1]$.

We denote the generators of $CF(\tilde{H}_\ell; \Z)$ by $\check{x}_0^\ell, \hat{x}_1^\ell, \dots, \hat{x}_\ell^\ell, \check{x}_\ell^\ell$, where for every $1 \leq j \leq \ell$ the generator $\check{x}_j^\ell$ corresponds to the minimum of a perfect Morse function on $S_j$, and $\hat{x}_j^\ell$ to the maximum of that function.

The main result of this section is the following proposition:

\begin{prop}\label{prop: CF of J-shaped in C^n}
    There exist signs $A_{\ell,1}, \dots, A_{\ell,\ell} \in \{-1,1\}$ such that the differential $d \fc CF(\tilde{H}_\ell; \Z) \to CF(\tilde{H}_\ell; \Z)$ satisfies:
    \begin{itemize}
        \item $d\hat{x}_i^\ell = A_{\ell,i} \check{x}_{i-1}^\ell$, for every $1 \leq i \leq \ell$.
        \item $d\check{x}_i^\ell = 0$, for every $0 \leq i \leq \ell$.
    \end{itemize}
\end{prop}

\begin{proof}
    We begin by computing the Floer complex of an auxiliary Hamiltonian. Let $\ell \in \Z_{\geq 0}$ and define $F_\ell \fc \Int B(1) \to \R$ by 
    $$ F_\ell(z) = \left(\ell + \frac{1}{2}\right)R(z) = \pi\left(\ell + \frac{1}{2}\right)\|z\|^2 $$
    for every $z \in \Int B(1)$. Then $F_\ell$ has only one $1$-periodic orbit $\gamma$, a constant orbit located at the origin. The Robbin-Salamon index of $\gamma$ with respect to the constant trivialization is $-n(2\ell+1)$; hence,
    $$ \mu_{FMB}^\tau(\gamma;F_\ell) = \mu_{RS}^\tau(\gamma;F_\ell) + \frac{1}{2}\dim \R^{2n} - \frac{1}{2}\dim \{0\} + \ind_h \gamma = -n(2\ell+1) + n = -2n\ell, $$
    where $h$ is a Morse function on $\{0\}$, which must be constant. Thus, the Floer homology of $F_\ell$ is
    $$ HF^*(F_\ell;\Z) = \begin{cases}
        \Z \cdot [\gamma], & * = -2n\ell, \\
        0, & * \neq -2n\ell.
    \end{cases} $$

    Returning to the Hamiltonian $\tilde{H}_\ell$, let $\tau$ denote the trivialization of the $1$-periodic orbits of $\tilde{H}_\ell$ obtained from cappings. From Section~\ref{ss: computations of FMB}, we know that $\mu_{FMB}^\tau(\check{x}_0^\ell)=0$, and for every $1 \leq i \leq \ell$:
    $$ \mu_{FMB}^\tau(\check{x}_i^\ell) = -2ni \quad \text{and} \quad \mu_{FMB}^\tau(\hat{x}_i^\ell) = -2ni + 2n - 1 = -2n(i-1) - 1. $$

    In the case $n=1$, since the differential increases the FMB-index by $+1$, we deduce that for every $1 \leq i \leq \ell$ there exists $C_i \in \Z$ such that $d\check{x}_i^\ell = C_i \hat{x}_i^\ell$. Fix $1\le i \le \ell$.
    Note that $\C^n$ is an exact symplectic manifold, hence action is well defined on all asymptotic orbits of the cascades. Since Floer solutions are action increasing, and since the action of  $\hat{x}^\ell_i$ equals the action of $\check{x}^\ell_i$, all cascades are of zero energy, hence constant. Therefore, there are no cascades in the flowline. Thus $C_i=0$ by the usual count of Morse differential in the Morse homology of $S^1$.

    Consequently, for any $n \in \N$, since the differential increases the FMB-index by $+1$, there exist $A_{\ell,1}, \dots, A_{\ell,\ell} \in \Z$ such that the differential $d \fc CF(\tilde{H}_\ell;\Z) \to CF(\tilde{H}_\ell;\Z)$ satisfies:
    \begin{itemize}
        \item $d\hat{x}_i^\ell = A_{\ell,i} \check{x}_{i-1}^\ell$, for every $1 \leq i \leq \ell$.
        \item $d\check{x}_i^\ell = 0$, for every $0 \leq i \leq \ell$.
    \end{itemize}
    To complete the proof, it remains to show that $A_{\ell,i} \in \{-1, 1\}$ for every $1 \leq i \leq \ell$.

    Note that $F_\ell$ and $\tilde{H}_\ell$ are admissible Hamiltonians with $\tilde{H}_\ell \leq F_\ell$, and they share the same slope at infinity. Thus, by Claim~\ref{claim: qi for Ham with the same slope}, we have an isomorphism $HF^*(\tilde{H}_\ell;\Z) \cong HF^*(F_\ell;\Z)$.

    Let $1 \leq i \leq \ell$. Since $HF^{-2n(i-1)}(\tilde{H}_\ell;\Z) \cong HF^{-2n(i-1)}(F_\ell;\Z) = 0$, and since all the elements of $CF^{-2n(i-1)}(\tilde{H}_\ell;\Z)$ are closed, we deduce that the restriction of the differential
    $$ d|_{CF^{-2n(i-1)-1}(\tilde{H}_\ell;\Z)} \fc CF^{-2n(i-1)-1}(\tilde{H}_\ell;\Z) \to CF^{-2n(i-1)}(\tilde{H}_\ell;\Z) $$
    is surjective. Since $CF^{-2n(i-1)-1}(\tilde{H}_\ell;\Z) = \Z \cdot \hat{x}_{i}^\ell$ and $CF^{-2n(i-1)}(\tilde{H}_\ell;\Z) = \Z \cdot \check{x}_{i-1}^\ell$, it follows that $A_{\ell,i} \in \{-1, 1\}$ as required.
\end{proof}

\subsection{Algebraic preparation}\label{ss: alg preperation}

In this section, we examine a specific class of cochain complexes over $\Z$ and deduce the structure of their differentials, as well as cochain maps between them. These complexes arise naturally in our computations of the Floer cochain complexes in $\CP^n$ associated with the Hamiltonians under consideration.

\begin{figure}[H]
    \centering
    \includegraphics[width=\linewidth]{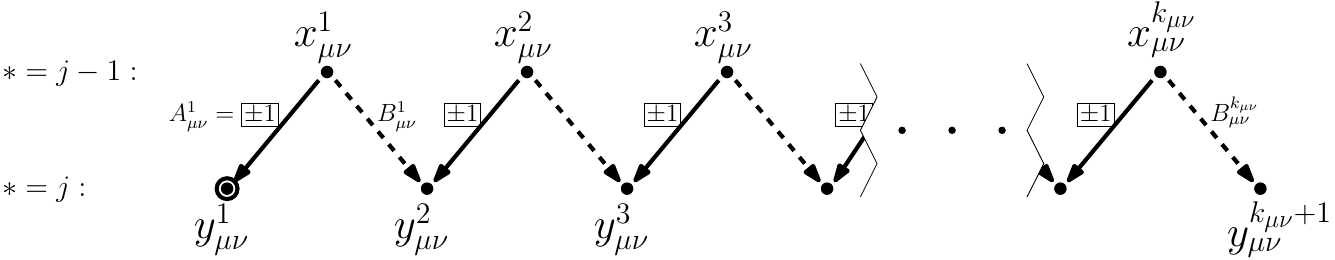}
    \caption{}
    \label{fig:Alg-Prep-Complex}
    \end{figure}

Let $j\in \Z$. For every $\mu,\nu\in\{0,1\}$ take $k_{\mn}\in \N$ such that
$$k_{00}\leq \min\{k_{10},k_{01}\}\leq \max\{k_{10},k_{01}\}\leq k_{11}.$$
For every $\mu,\nu\in\{0,1\}$ define a cochain complex $(C_{\mn}^*,d_{\mn})$ as follows:
\begin{itemize}
    \item For every $\ell\in\Z\setminus\{j-1,j\}$ we have $C_{\mn}^\ell=0$;

    \item $C_{\mn}^{j-1}$ is a free $\Z$-module with $k_{\mn}$ generators, denoted by $x_{\mn}^1,\ldots,x_{\mn}^{k_{\mn}}$;
    \item $C_{\mn}^{j}$ is a free  $\Z$-module with $k_{\mn}+1$ generators, denoted by $y_{\mn}^1,\ldots,y_{\mn}^{k_{\mn}+1}$;
    \item For every $1\leq i\leq k_{\mn}$ there are $A_{\mn}^i,B_{\mn}^i\in \Z$ such that $d_{\mn} x_{\mn}^i=A_{\mn}^i y_{\mn}^i+B_{\mn}^i y_{\mn}^{i+1}$;
    \item $H^j(C_{\mn}^*,d_{\mn})=\Z\cdot[y_{\mn}^1]$;
    \item $|A_{\mn}^1|=\cdots=|A_{\mn}^{k_{\mn}}|=1$.
\end{itemize}
Such a complex is depicted in Figure \ref{fig:Alg-Prep-Complex}.

We start with the following lemma.
\begin{lemma}\label{lemma: zigzag complex}
    For every $\mu,\nu\in\{0,1\}$ we have $|B_{\mn}^1|=\cdots=|B_{\mn}^{k_{\mn}}|=1$.
\end{lemma}
\begin{proof}
 
For every $1\leq i\leq k_{\mn}$ we have $(A_{\mn}^i)^2=1$, and $d_{\mu\nu}x^i_{\mu\nu}=0$. 
Hence, 
\begin{align*}
    [y_{\mn}^i]&=A_{\mn}^i\cdot[A_{\mn}^i y_{\mn}^i]\\
    &=A_{\mn}^i\cdot([A_{\mn}^i y_{\mn}^i+B_{\mn}^i y_{\mn}^{i+1}]-[B_{\mn}^i y_{\mn}^{i+1}])\\
    &=A_{\mn}^i\cdot([d_{\mn} x_{\mn}^i]-B_{\mn}^i[ y_{\mn}^{i+1}])\\
    &=-A_{\mn}^i B_{\mn}^i[ y_{\mn}^{i+1}].
\end{align*}
In particular, we deduce that for every $1\leq i\leq k_{\mn}$ we have $B_{\mn}^i\in \Z\setminus\{0\}$, and hence $|B_{\mn}^i|\geq1$.
By inductively applying the last equality, we deduce that 
$$[y_{\mn}^1]=(-1)^{k_{\mn}}\left(\prod_{i=1}^{k_{\mn}} A_{\mn}^i B_{\mn}^i\right)[y_{\mn}^{k_{\mn}+1}].$$
Since $H^j(C_{\mn}^*,d_{\mn})=\Z\cdot[y_{\mn}^1]$ we get that
$$[y_{\mn}^{k_{\mn}+1}]\in H^j(C_{\mn}^*,d_{\mn})=\Z\cdot[y_{\mn}^1]=\Z\cdot (-1)^{k_{\mn}}\left(\prod_{i=1}^{k_{\mn}} A_{\mn}^i B_{\mn}^i\right)[y_{\mn}^{k_{\mn}+1}],$$
and this shows that 
$$1=\left|(-1)^{k_{\mn}}\left(\prod_{i=1}^{k_{\mn}} A_{\mn}^i B_{\mn}^i\right)\right|=\prod_{i=1}^{k_{\mn}} |B_{\mn}^i|.$$
Thus, since $|B_{\mn}^1|,\ldots,|B_{\mn}^{k_{\mn}}|\geq1$ we deduce that in fact $|B_{\mn}^1|=\cdots=|B_{\mn}^{k_{\mn}}|=1$.   
\end{proof}

The next lemma shows that by a change of basis, we can make all the signs be positive.
\begin{lemma}\label{lemma: positive signs for zigzag}
    For every $\mu,\nu\in\{0,1\}$ there exist signs $\epsilon_{\mn}^1,\ldots,\epsilon_{\mn}^{k_{\mn}},\delta_{\mn}^1,\ldots,\delta_{\mn}^{k_{\mn}+1}\in\{-1,1\}$ such that, denoting
    \begin{align*}
        &\bar{x}_{\mn}^i=\epsilon_{\mn}^i x_{\mn}^i \ \  \text{for every}\ \  1\leq i\leq k_{\mn}\ \ \text{and}, \\  &\bar{y}_{\mn}^i=\delta_{\mn}^i y_{\mn}^i\ \  \text{for every}\ \  1\leq i\leq k_{\mn}+1,
    \end{align*} 
    it holds that for every $1\leq i\leq k_{\mn}$ we have $d_{\mn}\bar{x}_{\mn}^i=\bar{y}_{\mn}^i+\bar{y}_{\mn}^{i+1}$.
\end{lemma}

\begin{proof}
    
    Let us fix $\mu,\nu\in\{0,1\}$. First, set $\delta_{\mn}^1=1$. Now, for every $1\leq i\leq k_{\mn}$, construct iteratively:
    $$\epsilon_{\mn}^i=A_{\mn}^i\delta_{\mn}^i,\qquad \delta_{\mn}^{i+1}=B_{\mn}^i\epsilon_{\mn}^i.$$
    Then for every $1\leq i\leq k_{\mn}$ we have
    \begin{align*}
        d_{\mn}\bar{x}_{\mn}^i&=d_{\mn}(\epsilon_{\mn}^i x_{\mn}^i)\\
        &=\epsilon_{\mn}^id_{\mn}x_{\mn}^i\\
        &=\epsilon_{\mn}^i(A_{\mn}^i y_{\mn}^i+B_{\mn}^i y_{\mn}^{i+1})\\
        &=(A_{\mn}^i)^2\delta_{\mn}^i y_{\mn}^i+\delta_{\mn}^{i+1}y_{\mn}^{i+1}\\
        &=\bar{y}_{\mn}^i+\bar{y}_{\mn}^{i+1},
    \end{align*}
    as required.
\end{proof}

Let $\nu\in\{0,1\}$ and let $\Phi_\nu\fc C_{0\nu}\to C_{1\nu}$ be a cochain map. Assume that there are integers $C_{\nu}^1,\ldots,C_{\nu}^{k_{0\nu}},D_{\nu}^1,\ldots,D_{\nu}^{k_{0\nu}+1}\in \Z$ such that $\Phi_\nu(x_{0\nu}^i)=C_\nu^i x_{1\nu}^i$ for every $1\leq i\leq k_{0\nu}$ and $\Phi_\nu(y_{0\nu}^i)=D_\nu^i y_{1\nu}^i$ for every $1\leq i\leq k_{0\nu}+1$. 
\begin{lemma}\label{lemma: ABCD=1}
    For every $\nu\in\{0,1\}$ and $1\leq i\leq k_{0\nu}$ we have
    $$C_\nu^i A_{1\nu}^i=A_{0\nu}^i D_\nu^i,\qquad\text{and}\qquad C_\nu^iB_{1\nu}^i=B_{0\nu}^i D_\nu^{i+1}.$$
\end{lemma}
\begin{proof}
    Let $\nu\in\{0,1\}$ and $1\leq i\leq k_{0\nu}$. We know that
    $$d_{1\nu} (\Phi_\nu x_{0\nu}^i)=d_{1\nu} (C_\nu^i x_{1\nu}^i)=C_\nu^i d_{1\nu} ( x_{1\nu}^i)=C_\nu^i ( A_{1\nu}^i y_{1\nu}^i+B_{1\nu}^i y_{1\nu}^{i+1})=C_\nu^i A_{1\nu}^i y_{1\nu}^i+C_\nu^iB_{1\nu}^i y_{1\nu}^{i+1},$$
    and
    $$\Phi_\nu(d_{0\nu}x_{0\nu}^i)=\Phi_\nu(A_{0\nu}^i y_{0\nu}^i+B_{0\nu}^i y_{0\nu}^{i+1})=A_{0\nu}^i \Phi_\nu(y_{0\nu}^i)+B_{0\nu}^i \Phi_\nu(y_{0\nu}^{i+1})=A_{0\nu}^i D_\nu^iy_{1\nu}^i+B_{0\nu}^i D_\nu^{i+1}y_{1\nu}^{i+1}.$$
    Thus, since $\Phi_\nu$ is a cochain map, i.e., $\Phi_\nu\circ d_{0\nu}=d_{1\nu}\circ \Phi_\nu$, we deduce that 
    $$C_\nu^i A_{1\nu}^i y_{1\nu}^i+C_\nu^i B_{1\nu}^i y_{1\nu}^{i+1}=d_{1\nu} (\Phi_\nu x_{0\nu}^i)=\Phi_\nu(d_{0\nu}x_{0\nu}^i)=A_{0\nu}^i D_\nu^i y_{1\nu}^i+B_{0\nu}^i D_\nu^{i+1} y_{1\nu}^{i+1}.$$
    Since the elements $y_{1\nu}^i,y_{1\nu}^{i+1}$ are linearly independent, we deduce that 
    $$C_\nu^i A_{1\nu}^i=A_{0\nu}^i D_\nu^i,\qquad\text{and}\qquad C_\nu^iB_{1\nu}^i=B_{0\nu}^i D_\nu^{i+1},$$
    as required.
\end{proof}
Now, assume additionally that $\Phi_\nu\fc (C_{0\nu}^*,d_{0\nu})\to (C_{1\nu}^*,d_{1\nu})$ is a quasi-isomorphism. Then we obtain the following:

\begin{lemma}\label{lemma: zigzag and morphisms}
    For every $\nu\in\{0,1\}$ we have $|C_{\nu}^1|=\cdots=|C_{\nu}^{k_{0\nu}}|=|D_{\nu}^1|=\cdots=|D_{\nu}^{k_{0\nu}+1}|=1$.
\end{lemma}

\begin{proof}
    Let $\nu\in\{0,1\}$. Since $\Phi_\nu$ is a cochain map, it induces a morphism $(\Phi_\nu)_*\fc H^j(C_{0\nu}^*,d_{0\nu})\to H^j(C_{1\nu}^*,d_{1\nu})$, and this morphism satisfies
    $$(\Phi_\nu)_*([y_{0\nu}^1])=[\Phi_\nu(y_{0\nu}^1)]=[D_\nu^1 y_{1\nu}^1]=D_\nu^1[y_{1\nu}^1].$$
    Since $H^j(C_{0\nu}^*,d_{0\nu})=\Z\cdot[y_{0\nu}^1]$ and $H^j(C_{1\nu}^*,d_{1\nu})=\Z\cdot[y_{1\nu}^1]$, and $\Phi_\nu$ is a quasi-isomorphism, it follows that $(\Phi_\nu)_*$ is an isomorphism, and hence 
    $$[y_{1\nu}^1]\in H^j(C_{1\nu}^*,d_{1\nu})=\im (\Phi_\nu)_*=\Z\cdot D_\nu^1[y_{1\nu}^1],$$
    and this implies that $|D_\nu^1|=1$.

    Let us now prove that for every $1\leq i\leq k_{0\nu}$, if $|D_\nu^i|=1$, then $|C_\nu^i|=|D_\nu^{i+1}|=1$, which, by induction, will complete the proof.
    Indeed, let $1\leq i\leq k_{0\nu}$ and assume that $|D_\nu^i|=1$. Since $|A_{0\nu}^i|=|A_{1\nu}^i|=|B_{0\nu}^i|=|B_{1\nu}^i|=1$, by Lemma~\ref{lemma: ABCD=1} we get that
    $$|C_\nu^i|=|C_\nu^iA_{1\nu}^i|=|A_{0\nu}^i D_\nu^i|=1,$$
    and
    $$|D_\nu^{i+1}|=|B_{0\nu}^i D_\nu^{i+1}|=|C_\nu^i B_{1\nu}^i|=1,$$
    as required.
\end{proof}

Next we show that by a change of basis, we can arrange for both the signs in the differential and the continuation map to be positive.
\begin{lemma}\label{lemma: positive signs for zigzag and maps}
    Let $\nu\in\{0,1\}$, and assume that 
$A_{0\nu}^1=\cdots=A_{0\nu}^{k_{0\nu}}=B_{0\nu}^1=\cdots=B_{0\nu}^{k_{0\nu}}=1$.  Then there exist signs $\epsilon_{1\nu}^1,\ldots,\epsilon_{1\nu}^{k_{1\nu}},\delta_{1\nu}^1,\ldots,\delta_{1\nu}^{k_{1\nu}+1}\in\{-1,1\}$ such that, by denoting $\bar{x}_{1\nu}^i=\epsilon_{1\nu}^i x_{1\nu}^i$ for every $1\leq i\leq k_{1\nu}$ and $\bar{y}_{1\nu}^i=\delta_{1\nu}^i y_{1\nu}^i$ for every $1\leq i\leq k_{1\nu}+1$, the following hold:
    \begin{itemize}
        \item For every $1\leq i\leq k_{1\nu}$, we have $d_{1\nu}\bar{x}_{1\nu}^i=\bar{y}_{1\nu}^i+\bar{y}_{1\nu}^{i+1}$;
        \item For every $1\leq i\leq k_{0\nu}$, we have $\Phi_\nu(x_{0\nu}^i)=\bar{x}_{1\nu}^i$;
        \item For every $1\leq i\leq k_{0\nu}+1$, we have $\Phi_\nu(y_{0\nu}^i)=\bar{y}_{1\nu}^i$.
    
    \end{itemize}
\end{lemma}

\begin{proof}

By applying Lemma~\ref{lemma: positive signs for zigzag} to the cochain complex $(C_{1\nu}^*,d_{1\nu})$, we deduce that there exist signs $\tilde\epsilon_{1\nu}^1,\ldots,\tilde\epsilon_{1\nu}^{k_{1\nu}},\tilde\delta_{1\nu}^1,\ldots,\tilde\delta_{1\nu}^{k_{1\nu}+1}\in\{-1,1\}$ such that if we denote $\tilde{x}_{1\nu}^i=\tilde\epsilon_{1\nu}^i x_{1\nu}^i$ for every $1\leq i\leq k_{1\nu}$ and $\tilde{y}_{1\nu}^i=\tilde\delta_{1\nu}^i y_{1\nu}^i$ for every $1\leq i\leq k_{1\nu}+1$, then for every $1\leq i\leq k_{1\nu}$ we have $d_{1\nu}\tilde{x}_{1\nu}^i=\tilde{y}_{1\nu}^i+\tilde{y}_{1\nu}^{i+1}$.

For every $1\leq i\leq k_{0\nu}$ let $\tilde C_\nu^i=\tilde\epsilon_{1\nu}^i C_\nu^i$; then we have 
$$\Phi_\nu(x_{0\nu}^i)=C_\nu^i x_{1\nu}^i=  C_\nu^i \tilde\epsilon_{1\nu}^i \tilde x_{1\nu}^i=\tilde C_\nu^i \tilde x_{1\nu}^i.$$

Similarly, for every $1\leq i\leq k_{0\nu}+1$ let $\tilde D_\nu^i=\tilde\delta_{1\nu}^i D_\nu^i$; then we have 
$$\Phi_\nu(y_{0\nu}^i)=D_\nu^i y_{1\nu}^i=  D_\nu^i \tilde\delta_{1\nu}^i \tilde y_{1\nu}^i=\tilde D_\nu^i \tilde y_{1\nu}^i.$$

By applying Lemma~\ref{lemma: ABCD=1} with respect to this new basis, we deduce that for every $1\leq i\leq k_{0\nu}$ we have $\tilde D_\nu^i=\tilde C_\nu^i=\tilde D_\nu^{i+1}$, and hence
$$\tilde C_\nu^1=\cdots=\tilde C_\nu^{k_{0\nu}}=\tilde D_\nu^1=\cdots=\tilde D_\nu^{k_{0\nu}+1}\in\{-1,1\}.$$

Consider the signs that are given by $\epsilon_{1\nu}^i=\tilde C_\nu^1\tilde\epsilon_{1\nu}^i$ for every $1\leq i\leq k_{1\nu}$, and $\delta_{1\nu}^i=\tilde C_\nu^1\tilde\delta_{1\nu}^i$ for every $1\leq i\leq k_{1\nu}+1$. For every $1\leq i\leq k_{1\nu}$ denote $\bar x_{1\nu}^i=\epsilon_{1\nu}^i x_{1\nu}^i$ and for every $1\leq i\leq k_{1\nu}+1$ denote $\bar y_{1\nu}^i=\delta_{1\nu}^i y_{1\nu}^i$. Then:

\begin{itemize}
    \item For every $1\leq i\leq k_{1\nu}$, we have 
    \begin{align*}
        d_{1\nu}\bar{x}_{1\nu}^i&=d_{1\nu}(\epsilon_{1\nu}^i x_{1\nu}^i)\\
        &=d_{1\nu}(\tilde C_\nu^1 \tilde \epsilon_{1\nu}^i x_{1\nu}^i)\\
        &= \tilde C_\nu^1 d_{1\nu}(\tilde x_{1\nu}^i)\\
        &=\tilde C_\nu^1(\tilde{y}_{1\nu}^i+\tilde{y}_{1\nu}^{i+1})\\
        &=\tilde C_\nu^1\tilde\delta_{1\nu}^i {y}_{1\nu}^i+\tilde C_\nu^1 \tilde\delta_{1\nu}^{i+1}  {y}_{1\nu}^{i+1}\\
        &=\delta_{1\nu}^i {y}_{1\nu}^i+\delta_{1\nu}^{i+1}  {y}_{1\nu}^{i+1}\\
        &=\bar y_{1\nu}^i+\bar y_{1\nu}^{i+1}.
    \end{align*}
    \item For every $1\leq i\leq k_{0\nu}$ we have 
    $$\Phi_\nu(x_{0\nu}^i)=\tilde C_\nu^i \tilde x_{1\nu}^i=\tilde C_\nu^1 \tilde x_{1\nu}^i= \tilde C_\nu^1 \tilde\epsilon_{1\nu}^i x_{1\nu}^i=\epsilon_{1\nu}^i x_{1\nu}^i=\bar x_{1\nu}^i.$$
    \item For every $1\leq i\leq k_{0\nu}+1$ we have 
    $$\Phi_\nu(y_{0\nu}^i)=\tilde D_\nu^i \tilde y_{1\nu}^i=\tilde C_\nu^1 \tilde y_{1\nu}^i= \tilde C_\nu^1 \tilde\delta_{1\nu}^i y_{1\nu}^i=\delta_{1\nu}^i y_{1\nu}^i=\bar y_{1\nu}^i.$$
\end{itemize}

This completes the proof.
\end{proof}

Similarly to the previous case,
let $\mu\in\{0,1\}$ and let $\Psi_\mu\fc C_{\mu0}\to C_{\mu1}$ be a cochain map. Assume that there are integers $E_{\mu}^1,\ldots,E_{\mu}^{k_{\mu0}},F_{\mu}^1,\ldots,F_{\mu}^{k_{\mu0}+1}\in \Z$ such that $\Psi_\mu(x_{\mu0}^i)=E_\mu^i x_{\mu1}^i$ for every $1\leq i\leq k_{\mu0}$ and $\Psi_\mu(y_{\mu0}^i)=F_\mu^i y_{\mu1}^i$ for every $1\leq i\leq k_{\mu0}+1$. Additionally, assume that $\Psi_\mu$ is a quasi-isomorphism.
\begin{lemma}\label{lemma: zigzag and morphisms 2}
    For every $\mu\in\{0,1\}$ we have $|E_{\mu}^1|=\ldots=|E_{\mu}^{k_{\mu0}}|=|F_{\mu}^1|=\ldots=|F_{\mu}^{k_{\mu0}+1}|=1$. \qed
\end{lemma}

\begin{lemma}\label{lemma: positive signs for zigzag and maps 2}
    Let $\mu\in\{0,1\}$, and assume that
$A_{\mu0}^1=\cdots=A_{\mu0}^{k_{\mu0}}=B_{\mu0}^1=\cdots=B_{\mu0}^{k_{\mu0}}=1$.  Then there exist signs $\epsilon_{\mu1}^1,\ldots,\epsilon_{\mu1}^{k_{\mu1}},\delta_{\mu1}^1,\ldots,\delta_{\mu1}^{k_{\mu1}+1}\in\{-1,1\}$ such that, by denoting $\bar{x}_{\mu1}^i=\epsilon_{\mu1}^i x_{\mu1}^i$ for every $1\leq i\leq k_{\mu1}$ and $\bar{y}_{\mu1}^i=\delta_{\mu1}^i y_{\mu1}^i$ for every $1\leq i\leq k_{\mu1}+1$, the following hold:
    \begin{itemize}
        \item For every $1\leq i\leq k_{\mu1}$, we have $d_{\mu1}\bar{x}_{\mu1}^i=\bar{y}_{\mu1}^i+\bar{y}_{\mu1}^{i+1}$;
        \item For every $1\leq i\leq k_{\mu0}$, we have $\Psi_\mu(x_{\mu0}^i)=\bar{x}_{\mu1}^i$;
        \item For every $1\leq i\leq k_{\mu0}+1$, we have $\Psi_\mu(y_{\mu0}^i)=\bar{y}_{\mu1}^i$. \qed
    \end{itemize}
\end{lemma}

The proofs of Lemma~\ref{lemma: zigzag and morphisms 2} are Lemma~\ref{lemma: positive signs for zigzag and maps 2} are identical to those of Lemma~\ref{lemma: zigzag and morphisms} and Lemma~\ref{lemma: positive signs for zigzag and maps}, respectively, and are therefore omitted.

Finally, if, moreover, the following square commutes
\begin{equation}\label{eq: comm square alg lemma}
\vcenter{\xymatrix{
    (C_{00}^*,d_{00}) \ar[r]^{\Phi_0} \ar[d]_{\Psi_0} & (C_{10}^*,d_{10}) \ar[d]^{\Psi_1} \\
    (C_{01}^*,d_{01}) \ar[r]_{\Phi_1} & (C_{11}^*,d_{11})
}}
\end{equation}
i.e., $\Phi_1\circ\Psi_0=\Psi_1\circ\Phi_0$, then we get the following result:

\begin{lemma}\label{lemma: positive signs for zigzag and maps in commuting case}
Assuming the commutativity of the square \eqref{eq: comm square alg lemma}, assume moreover that the signs in the diffrerentials in all vertices except for the bottom right, have been set to $1$, namely, that
$$A_{\mu0}^1=\cdots=A_{\mu0}^{k_{\mu0}}=B_{\mu0}^1=\cdots=B_{\mu0}^{k_{\mu0}}=1,$$
for all $\mu\in\{0,1\}$, and that
$$A_{01}^1=\cdots=A_{01}^{k_{01}}=B_{01}^1=\cdots=B_{01}^{k_{01}}=1.$$
Moreover, assume that the vertical maps have had their signs set to $1$, namely, that 
$$C_0^1=\cdots=C_0^{k_{00}}=D_0^1=\cdots=D_0^{k_{00}+1}=1,$$
and
$$E_0^1=\cdots=E_0^{k_{00}}=F_0^1=\cdots=F_0^{k_{00}+1}=1.$$

Then, by a change of basis, we may set the signs of the horizontal maps to $1$ as well. That is, there exist signs $\epsilon_{11}^1,\ldots,\epsilon_{11}^{k_{11}},\delta_{11}^1,\ldots,\delta_{11}^{k_{11}+1}\in\{-1,1\}$ such that, by denoting $\bar{x}_{11}^i=\epsilon_{11}^i x_{11}^i$ for every $1\leq i\leq k_{11}$ and $\bar{y}_{11}^i=\delta_{11}^i y_{11}^i$ for every $1\leq i\leq k_{11}+1$, the following hold:
    \begin{itemize}
        \item For every $1\leq i\leq k_{11}$, we have $d_{11}\bar{x}_{11}^i=\bar{y}_{11}^i+\bar{y}_{11}^{i+1}$;
        \item For every $1\leq i\leq k_{01}$, we have $\Phi_1(x_{01}^i)=\bar{x}_{11}^i$;
        \item For every $1\leq i\leq k_{01}+1$, we have $\Phi_1(y_{01}^i)=\bar{y}_{11}^i$;
        \item For every $1\leq i\leq k_{10}$, we have $\Psi_1(x_{10}^i)=\bar{x}_{11}^i$;
         \item For every $1\leq i\leq k_{10}+1$, we have $\Psi_1(y_{10}^i)=\bar{y}_{11}^i$.
    \end{itemize}
\end{lemma}

\begin{proof}

Consider the signs $\epsilon_{11}^1,\ldots,\epsilon_{11}^{k_{11}},\delta_{11}^1,\ldots,\delta_{11}^{k_{11}+1}\in\{-1,1\},$ obtained by applying Lemma~\ref{lemma: positive signs for zigzag and maps} to the cochain complexes $(C_{01}^*,d_{01})$, $(C_{11}^*,d_{11})$ and the quasi-isomorphism $\Phi_1$ between them. Thus, if we denote $\bar{x}_{11}^i=\epsilon_{11}^i x_{11}^i$ for every $1\leq i\leq k_{11}$ and  $\bar{y}_{11}^i=\delta_{11}^i y_{11}^i$ for every $1\leq i\leq k_{11}+1$, then we have 
\begin{itemize}
    \item For every $1\leq i\leq k_{11}$, we have $d_{11}\bar{x}_{11}^i=\bar{y}_{11}^i+\bar{y}_{11}^{i+1}$;
    \item For every $1\leq i\leq k_{01}$, we have $\Phi_1(x_{01}^i)=\bar{x}_{11}^i$;
    \item For every $1\leq i\leq k_{01}+1$, we have $\Phi_1(y_{01}^i)=\bar{y}_{11}^i$.
    \end{itemize}

First, we prove that there exists a sign $S\in\{-1,1\}$ such that for every $1\leq i\leq k_{10}$ we have $\Psi_1(x_{10}^i)=S \bar{x}_{11}^i$, and for every $1\leq i\leq k_{10}+1$ we have $\Psi_1(y_{10}^i)=S \bar{y}_{11}^i$.

Indeed, by the definitions of our maps and basis elements, there exist signs $$\bar{E}_1^1,\ldots,\bar{E}_1^{k_{10}},\bar{F}_1^1,\ldots,\bar{F}_1^{k_{10}+1}\in\{-1,1\}$$ such that for every $1\leq i\leq k_{10}$ we have 
$\Psi_1(x_{10}^i)=\bar{E}_1^i \bar{x}_{11}^i$ 
and for every $1\leq i\leq k_{10}+1$ we have $\Psi_1(y_{10}^i)=\bar{F}_1^i \bar{y}_{11}^i$. Since $\Psi_1\fc C_{10}^*\to C_{11}^*$ is a cochain map, we have for every $1\leq i\leq k_{10}$:
\begin{align*}
    \bar{F}_1^i \bar{y}_{11}^i + \bar{F}_1^{i+1} \bar{y}_{11}^{i+1} &= \Psi_1(y_{10}^i + y_{10}^{i+1}) \\
    &= \Psi_1(d_{10}x_{10}^i) \\
    &= d_{11}(\Psi_1 x_{10}^i) \\
    &= d_{11}(\bar{E}_1^i \bar{x}_{11}^i) \\
    &= \bar{E}_1^i(\bar{y}_{11}^i + \bar{y}_{11}^{i+1}) \\
    &= \bar{E}_1^i \bar{y}_{11}^i + \bar{E}_1^i \bar{y}_{11}^{i+1}.
\end{align*}
For every $1\leq i\leq k_{10}$, by the linear independence of $\bar{y}_{11}^i$ and $\bar{y}_{11}^{i+1}$, it follows that $\bar{F}_1^i = \bar{E}_1^i = \bar{F}_1^{i+1}$. Thus, all these signs are identical; we denote this common sign by $S$.

Second, we show that $S=1$. Indeed, by the assumption, the square commutes, i.e., $\Phi_1\circ\Psi_0=\Psi_1\circ\Phi_0$, therefore
$$S\bar{x}_{11}^1 = \Psi_1(x_{10}^1) = \Psi_1(\Phi_0(x_{00}^1)) = \Phi_1(\Psi_0(x_{00}^1)) = \Phi_1(x_{01}^1) = \bar{x}_{11}^1.$$
Since $\bar{x}_{11}^i\neq0$, we deduce that $S=1$, and this completes the proof.
\end{proof}

    \subsection{Proof of Theorem~\ref{thm: CF(H_l;Z)} and Theorem~\ref{thm: continuation maps, over Z}}\label{ss: proof of computations of CF + continuations}

    We start with a remark, followed by an auxiliary lemma:
    \begin{rem}\label{rem: conc from thm: diff_obst}
    Let $n\in \N$ and $\ell\in \Z_{\geq0}$. By Theorem~\ref{thm: diff_obst}, there are integers
$$A_{\ell,1}, \dots, A_{\ell,\ell}, B_{\ell,1}, \dots, B_{\ell,\ell} \in \Z,$$
such that the differential $d \fc CF(H_\ell; \Z) \to CF(H_\ell; \Z)$ satisfies:
\begin{itemize}
    \item For every $1 \leq i \leq \ell$, we have $d \hat{x}^\ell_i = A_{\ell,i} \check{x}^\ell_{i-1} + B_{\ell,i} \check{x}^\ell_{i+n}$.
    \item For every $0 \leq i \leq \ell+n$, we have $d \check{x}^\ell_i = 0$.
\end{itemize}

This is indeed an almost immediate consequence of Theorem~\ref{thm: diff_obst}, except for the case where $n=1$. In that case, we have $d\check{x}_i^\ell=C_i\hat{x}_i^\ell$ for every $1\leq i\leq \ell$, where $C_1,\ldots,C_\ell\in \Z$. 

Fix $1 \le i \le \ell$.
By Theorem~\ref{thm: diff_obst} every flowline with cascades connecting $\check{x}^\ell_i$ to $\hat{x}^\ell_i$ does not intersect $D_\infty$, and by the no escape lemma (Lemma \ref{lemma: no escape}), it is thus contained in a ball, and hence the action functional is well defined on all asymptotic orbits of the cascades. Since Floer solutions are action increasing, and since the action of  $\hat{x}^\ell_i$ equals the action of $\check{x}^\ell_i$, all cascades are of zero energy, hence constant. Therefore, there are no cascades in the flowline. Thus $C_i=0$ by the usual count of Morse differential in the Morse homology of $S^1$, and we deduce that $C_1=\cdots=C_\ell=0$.
    \end{rem}

    \begin{lemma}\label{lemma: HF(H_l)}
    For every $\ell\in \Z_{\geq0}$, the Floer homology $HF^*(H_\ell;\Z)$ of $H_\ell$ is non-zero only for $*\in\{0,2,\ldots,2n\}$. Moreover, for every $0\leq j\leq n$, the Floer homology group $HF^{2j}(H_\ell;\Z)$ is a copy of $\Z$ which is generated by $[\check{x}_j^\ell]$, i.e.
    $$HF^*(H_\ell;\Z)=\left\{\begin{array}{ll}
   \Z\cdot[\check{x}_j^\ell],  & *=2j,\,0\leq j\leq n, \\
    0, & \text{otherwise}. 
\end{array}\right.$$
 
    \end{lemma}
    \begin{proof}
        We start with the case $\ell=0$. The Hamiltonian $H_0$ is Morse--Bott and has no non-constant $1$-periodic orbits. Thus, the generators of $CF(H_0;\Z)$, which are $\check{x}_0^0, \check{x}_1^0, \ldots, \check{x}_n^0$, are critical points of $H_0$, each of a different $\mu_{FMB}$ index. Comparing this with the quantum cohomology of $\CP^n$, we deduce that the differential in $CF(H_0;\Z)$ vanishes.
$$HF^*(H_0;\Z)=\left\{\begin{array}{ll}
   \Z\cdot[\check{x}_j^0],  & *=2j,\,0\leq j\leq n, \\
    0, & \text{otherwise}. 
\end{array}\right.$$
Let $F \fc CF(H_0;\Z) \to CF(H_\ell;\Z)$ be the continuation map. By Theorem~\ref{thm: continuation_obst}, there are integers $k_0,\ldots,k_n\in \Z$ such that $F(\check{x}^0_j) = k_j \check{x}^{\ell}_{j}$ for every $0 \leq j \leq n$. Let us compute the continuation map $F \fc CF(H_0;\Z) \to CF(H_\ell;\Z)$ and the differential of $CF(H_\ell;\Z)$.

Let $0\leq j\leq n$. Since 
$$F_*([\check{x}_j^0]) = [k_j\check{x}_j^\ell] = k_j[\check{x}_j^\ell],$$
and $[\check{x}_j^0]$ is a generator for $HF^{2j}(H_0;\Z)$, we deduce that 
$ F_*({HF^{2j}(H_0;\Z)})=\Z\cdot k_j[\check{x}_j^\ell]$. On the other hand, since $\CP^n$ is a closed symplectic manifold, we know that the continuation map $F$ is a quasi-isomorphism; see \cite[Theorem 4]{Floer_1989_monotone_mfd}. Therefore, the induced map 
$$F_*|_{HF^{2j}(H_0;\Z)} \fc HF^{2j}(H_0;\Z) \to HF^{2j}(H_\ell;\Z)$$
is an isomorphism, and hence
$$HF^{2j}(H_\ell;\Z) = \im F_*|_{HF^{2j}(H_0;\Z)} = \Z\cdot k_j[\check{x}_j^\ell].$$
Thus, $[\check{x}_j^\ell] \in HF^{2j}(H_\ell;\Z) = \Z \cdot k_j[\check{x}_j^\ell]$, which implies that $k_j \in \{-1, 1\}$, and this shows that $HF^{2j}(H_\ell;\Z) = \Z \cdot [\check{x}_j^\ell]$.

    \end{proof}

    \begin{proof}[Proof of Theorem~\ref{thm: CF(H_l;Z)}]

Let $\ell\in \Z_{\geq0}$, and let us compute the differential of $CF(H_\ell;\Z)$. By Remark~\ref{rem: conc from thm: diff_obst}, the differential $d \fc CF(H_\ell; \Z) \to CF(H_\ell; \Z)$ satisfies
\begin{itemize}
    \item $d \hat{x}^\ell_i = A_{\ell,i} \check{x}^\ell_{i-1} + B_{\ell,i} \check{x}^\ell_{i+n}$, for every $1 \leq i \leq \ell$, 
    \item $d \check{x}^\ell_i = 0$ for every $0 \leq i \leq \ell+n$,
\end{itemize}
for some $A_{\ell,1},\ldots,A_{\ell,\ell},B_{\ell,1},\ldots,B_{\ell,\ell}\in \Z$. As explained on Page~\pageref{p: transversality paragraph}, the conclusion of Theorem~\ref{thm: diff_obst} holds, thus, for every $1\leq i\leq \ell$ we know that if $u$ is a bubbled Floer flowline with cascades that connects $\hat{x}_i^{\ell}$ to $\check{x}_{i-1}^\ell$ then $u$ does not intersect $D_\infty$. Therefore the flowlines that contribute to the differential of $\hat{x}_i^\ell$  in $CF(H_\ell;\Z)$ are the same as those that contribute to the differential of $\hat{x}_i^\ell$ in $CF(\tilde{H}_\ell;\Z)$, where $\tilde{H}_\ell=H_\ell|_{\CP^n\setminus D_\infty}$. Thus by Proposition~\ref{prop: CF of J-shaped in C^n} we deduce that $A_{\ell,1},\ldots, A_{\ell,\ell}\in \{-1,1\}$.

Note that the Floer cochain complex $CF(H_\ell;\Z)$ decomposes into $n+1$ subcomplexes. For each $0 \leq j \leq n$, the submodule 
$$C_{\ell,j}=CF^{2j-1}(H_\ell;\Z)\oplus CF^{2j}(H_\ell;\Z),$$
generated by
$$
\check{x}_{j}^\ell,\ \hat{x}_{j+1}^\ell,\ \check{x}_{j+n+1}^\ell,\ \hat{x}_{j+n+2}^\ell,\ \ldots,\ 
\hat{x}_{j+(k_{\ell,j}-1)(n+1)+1}^\ell,\ \check{x}_{j+k_{\ell,j}(n+1)}^\ell,
$$
where $k_{\ell,j} = \left\lceil \frac{\ell - j}{n+1} \right\rceil$, is a subcomplex of $CF(H_\ell;\Z)$. Additionally, by Lemma~\ref{lemma: HF(H_l)} we obtain that for every $0\leq j\leq n$ we have 
$$H^*(C_{\ell,j})=\left\{\begin{array}{lc}
    HF^{2j}(H_\ell;\Z), & *=2j, \\
     0, & *\neq 2j, 
\end{array}\right.=\left\{\begin{array}{lc}
    \Z\cdot[\check{x}_j^\ell], & *=2j, \\
     0, & *\neq 2j. 
\end{array}\right.$$

 Thus Lemma~\ref{lemma: zigzag complex}, applied for each of the cochain subcomplexes $C_{\ell,0},\ldots,C_{\ell,n}$, implies that for every $1\leq i\leq \ell$ we have $B_{\ell,i}\in\{-1,1\}$, which completes the proof.

    \end{proof}

\begin{proof}[Proof of Theorem~\ref{thm: continuation maps, over Z}]
     As we mentioned in the proof of Theorem~\ref{thm: CF(H_l;Z)}, each of the Floer cochain complexes $CF(H_\ell;\Z)$, $CF(H'_{\ell'};\Z)$ decomposes into $n+1$ subcomplexes. For each $0 \leq j \leq n$, the submodules 
$$C_{\ell,j}=CF^{2j-1}(H_\ell;\Z)\oplus CF^{2j}(H_\ell;\Z),\qquad\text{and}\qquad C'_{\ell',j}=CF^{2j-1}(H'_{\ell'};\Z)\oplus CF^{2j}(H'_{\ell'};\Z),$$
generated by
$$
\check{x}_{j}^\ell,\ \hat{x}_{j+1}^\ell,\ \check{x}_{j+n+1}^\ell,\ \hat{x}_{j+n+2}^\ell,\ \ldots,\ 
\hat{x}_{j+(k_{\ell,j}-1)(n+1)+1}^\ell,\ \check{x}_{j+k_{\ell,j}(n+1)}^\ell,
$$ and
$$
\check{y}_{j}^{\ell'},\ \hat{y}_{j+1}^{\ell'},\ \check{y}_{j+n+1}^{\ell'},\ \hat{y}_{j+n+2}^{\ell'},\ \ldots,\ 
\hat{y}_{j+(k_{{\ell'},j}-1)(n+1)+1}^{\ell'},\ \check{y}_{j+k_{{\ell'},j}(n+1)}^{\ell'},
$$
respectively, where $k_{\ell,j} = \left\lceil \frac{\ell - j}{n+1} \right\rceil$ and $k_{\ell',j} = \left\lceil \frac{\ell' - j}{n+1} \right\rceil$, are subcomplexes of $CF(H_\ell;\Z)$ and $CF(H'_{\ell'};\Z)$, respectively. Additionally, by Lemma~\ref{lemma: HF(H_l)} we obtain that for every $0\leq j\leq n$ we have 
$$H^*(C_{\ell,j})=\left\{\begin{array}{lc}
    HF^{2j}(H_\ell;\Z), & *=2j, \\
     0, & *\neq 2j,
\end{array}\right.=\left\{\begin{array}{lc}
    \Z\cdot[\check{x}_j^\ell], & *=2j, \\
     0, & *\neq 2j,
\end{array}\right.$$
and 
$$H^*(C'_{\ell',j})=\left\{\begin{array}{lc}
    HF^{2j}(H'_{\ell'};\Z), & *=2j, \\
     0, & *\neq 2j,
\end{array}\right.=\left\{
\begin{array}{lc}
    \Z\cdot[\check{y}_j^{\ell'}], & *=2j, \\
     0, & *\neq 2j.
\end{array}\right.$$
By Theorem~\ref{thm: continuation_obst} there are $\check{C}_{0},\ldots,\check{C}_{\ell+n},\hat{C}_{1},\ldots,\hat{C}_{\ell}\in \Z$ such that $\Phi\check{x}_i^\ell=\check{C}_i \check{y}_i^{\ell'}$ for every $0\leq i\leq \ell+n$ and $\Phi\hat{x}_i^\ell=\hat{C}_i \hat{y}_i^{\ell'}$ for every $1\leq i\leq \ell$.

Since $\CP^n$ is a closed symplectic manifold, the continuation map $\Phi\fc CF(H_\ell;\Z)\to CF(H'_{\ell'};\Z)$ is a quasi-isomorphism. Thus Theorem~\ref{thm: CF(H_l;Z)} and  Lemma~\ref{lemma: zigzag and morphisms}, applied for each of the pairs $(C_{\ell,0},C'_{\ell,0}),\ldots,(C_{\ell,n},C'_{\ell,n})$ of subcomplexes the pair of cochain complexes $(CF(H_\ell;\Z) ,CF(H'_{\ell'};\Z))$, imply that 
$$|\check{C}_0|=\cdots=|\check{C}_{\ell+n}|=|\hat{C}_1|=\cdots=|\hat{C}_\ell|=1,$$
as required.

\end{proof}

    \section{Algebra}\label{s: algebra}

    \subsection{Novikov modules}\label{ss: Novikov modules}

For every ring $\cR$, let us define 
$$\Lambda^\cR = \left\{ \sum_{i=1}^\infty a_i T^{\lambda_i} : a_i \in \cR, \, \lambda_i \in \R \text{ s.t. } \lambda_i \to +\infty \right\},$$
the Novikov module over $\cR$, which is also an algebra over $\cR$. 
The Novikov module admits a natural valuation $\val \fc \Lambda^\cR \to \R \cup \{+\infty\}$ defined as follows: for $x \in \Lambda^\cR$, if $x = 0$ then $\val(x) = +\infty$; otherwise, $\val(x) = \min \{ \lambda_i : a_i \neq 0 \}$ where $x = \sum a_i T^{\lambda_i}$.

Using this valuation, we define the \textbf{interval modules}:
$$\Lambda^\cR_{\geq r} = \val^{-1}([r, \infty]), \qquad \Lambda^\cR_{>r} = \val^{-1}((r, \infty]), \qquad \Lambda^\cR_{[a,b)} = \Lambda^\cR_{\geq a} / \Lambda^\cR_{\geq b} \text{ for } a < b.$$
Note that all of these are modules over $\Lambda^\cR_{\geq 0}$.

When $\cR = \mathbb{Q}$, we denote $\Lambda = \Lambda^\mathbb{Q}$ and $\Lambda_{\geq 0} = \Lambda^\mathbb{Q}_{\geq 0}$ and refer to them as the Novikov field and Novikov ring, respectively. A module over $\Lambda_{\geq 0}$ is called a \textbf{Novikov module}.

   \subsubsection{Completion of $\Lambda_{\geq 0}$-modules}\label{ss: completion}
	
Completion of $\Lambda_{\geq 0}$-modules plays an essential role in Varolg\"une\c s's definition of symplectic cohomology; therefore, we dedicate this section to describing it and to recollection of important facts about it.

For every $0 \leq r \leq r'$, we have $\Lambda_{\geq r'} \subset \Lambda_{\geq r}$, which provides a natural morphism
$$\Lambda_{[0,r')} = \Lambda_{\geq 0} / \Lambda_{\geq r'} \to \Lambda_{\geq 0} / \Lambda_{\geq r} =  \Lambda_{[0,r)}.$$
Therefore, $(\Lambda_{[0,r)})_{r>0}$ is an inverse system.

\begin{defin}
	Let $A$ be a Novikov module. The \textbf{completion} of $A$, denoted by $\widehat{A}$, is defined as
	$$\widehat{A} = \varprojlim_r A \otimes_{\Lambda_{\geq 0}} \Lambda_{\geq 0} / \Lambda_{\geq r},$$
	where $r$ runs over $(0, +\infty)$.
\end{defin}

\begin{rem}\label{rem: properties of completion}\phantom{M}
	\begin{itemize}
		\item Completion is a functor $\widehat{\phantom{x}} \fc \Mod(\Lambda_{\geq 0}) \to \Mod(\Lambda_{\geq 0})$.

        \item As a functor, the completion of an isomorphism is again an isomorphism.
        
		\item The completion functor automatically extends to a functor $\widehat{\phantom{x}} \fc \Ch(\Lambda_{\geq 0}) \to \Ch(\Lambda_{\geq 0})$. Specifically, if $(C, d)$ is a cochain complex over $\Lambda_{\geq 0}$, then the completion $(\widehat{C}, \widehat{d})$ is obtained by applying the completion functor to the underlying module and to the map $d \fc C \to C$. We also extend the completion functor to the categories of graded $\Lambda_{\geq 0}$-modules and graded cochain complexes by degree-wise completion.

		\item Given a Novikov module $A$, the maps $A \to A \otimes_{\Lambda_{\geq 0}} \Lambda_{\geq 0} / \Lambda_{\geq r}$ given by $x \mapsto x \otimes 1$ for every $r > 0$ induce a natural map $A \to \widehat{A}$. If this map is an isomorphism, $A$ is called \textbf{complete}. The interval modules are typical examples of complete modules. On the other hand, $\Lambda$ is not complete since $\widehat{\Lambda} = 0$ (see \cite[Page 599]{Varolgunes_2021_MV_and_relSH}).

        \item Given two Novikov modules $A,B$ and a morphism $f\fc A\to B$, the following  diagram commutes:

        $$
        \xymatrix@R=2pc@C=4pc{
            A \ar[r] \ar[d]_f & \widehat{A} \ar[d]^{\hat{f}} \\
            B \ar[r] & \widehat{B}
        }$$
        where the right arrow is the completion of $f$ and the horizontal arrows are those described in the previous item.

		\item When $A$ is free, this description becomes simpler. Choose a basis $\{v_i\}_{i \in I}$. Then $\widehat{A}$ is isomorphic to
		$$\left\{\sum_{i \in I} \beta_i v_i : \forall i \in I, \, \beta_i \in \Lambda_{\geq 0}, \text{ and } \forall R \geq 0, \, \#\{i \in I : \val(\beta_i) < R\} < \infty \right\}.$$

		\item For every $r > 0$, we have a natural isomorphism $A \otimes_{\Lambda_{\geq 0}} \Lambda_{\geq 0} / \Lambda_{\geq r} \to A / T^r A$. Thus, the completion can be expressed as
		$$\widehat{A} = \varprojlim A / T^r A.$$
	\end{itemize}
\end{rem}

\begin{defin}
	Let $A$ be a Novikov module and let $(a_i)_{i \in \N}$ be a sequence of elements in $A$.
	\begin{itemize}
		\item Let $a \in A$. The sequence $(a_i)_{i \in \N}$ \textbf{converges} to $a$ if for every $r > 0$ there exists $N \in \N$ such that for every $i > N$ we have $a_i - a \in T^r A$.
		\item The sequence $(a_i)_{i \in \N}$ is a \textbf{Cauchy sequence} if for every $r > 0$ there exists $N \in \N$ such that for every $i, j > N$ we have $a_i - a_j \in T^r A$.
	\end{itemize}
\end{defin}

By a standard argument, any convergent sequence is also a Cauchy sequence. Additionally, as mentioned in \cite[Section 2.3]{Varolgunes_2021_MV_and_relSH}, the completion of a Novikov module is isomorphic to the module of equivalence classes of Cauchy sequences, where two Cauchy sequences are equivalent if their difference converges to $0$. Note that the completeness of a module is equivalent to all Cauchy sequences being convergent. The following result is well known; we include its proof for completeness.

\begin{claim}\label{claim: series in a Novikov module 0} 
    Let $A$ be a complete Novikov module and let $(\lambda_i)_{i \in \N}$ be a sequence in $A$. If $\lim_{i\to\infty} \lambda_i = 0$, then the series $\sum_{i=1}^\infty \lambda_i$ converges.
\end{claim}

\begin{proof}
    We will show that the sequence $(\sum_{i=1}^n \lambda_i)_{n \in \N}$ is a Cauchy sequence. 
    Let $r > 0$. Since $\lim_{i\to\infty} \lambda_i = 0$, there exists $i_0 \in \N$ such that for every $i > i_0$, we have $\lambda_i \in T^r A$. Let $n > m > i_0$. We have 
    $$\sum_{i=1}^n \lambda_i - \sum_{i=1}^m \lambda_i = \sum_{i=m+1}^n \lambda_i \in T^r A.$$ 
    Therefore, the sequence $(\sum_{i=1}^n \lambda_i)_{n \in \N}$ is a Cauchy sequence. Since $A$ is complete, we deduce that the sequence $(\sum_{i=1}^n \lambda_i)_{n \in \N}$ converges; i.e., the series $\sum_{i=1}^\infty \lambda_i$ converges.
\end{proof}

An immediate corollary of this result is:

\begin{coroll}\label{coroll: series in a Novikov module}
    Let $A$ be a complete Novikov module. Let $(\lambda_i)_{i \in \N}$ be a sequence of elements in $A$ and let $(\alpha_i)_{i \in \N}$ be a sequence of elements in $\Lambda_{\geq 0}$. If $\val(\alpha_i) \to +\infty$, then the series $\sum_i \alpha_i \lambda_i$ converges.
\end{coroll}

\begin{proof}
    Let $r > 0$. Then there exists $N \in \N$ such that for every $i > N$, we have $\val(\alpha_i) > r$ and hence $\alpha_i \in T^r \Lambda_{\geq 0}$; therefore, there exists $\beta_i \in \Lambda_{\geq 0}$ such that $\alpha_i = T^r \beta_i$. Thus, for every $i > N$, we have $\alpha_i \lambda_i = T^r \beta_i \lambda_i \in T^r A$. Hence $\lim_{i \to \infty} \alpha_i \lambda_i = 0$, so the series $\sum_{i=1}^\infty \alpha_i \lambda_i$ converges.
\end{proof}

\begin{prop}\label{prop: rearrangement}
    Let $A$ be a complete Novikov module. Let $(w_{i,j})_{i,j \in \N}$ be a double sequence in $A$. Suppose that for every $r > 0$, there exists a finite set $S_r \subset \N \times \N$ such that $w_{i,j} \in T^r A$ for all $(i,j) \in (\N \times \N) \setminus S_r$. Then:
    \begin{enumerate}
        \item For each $i \in \N$, the series $\sum_{j=1}^\infty w_{i,j}$ converges, and for each $j \in \N$, the series $\sum_{i=1}^\infty w_{i,j}$ converges .
        \item The iterated sums are equal:
        $$\sum_{i=1}^\infty \left( \sum_{j=1}^\infty w_{i,j} \right) = \sum_{j=1}^\infty \left( \sum_{i=1}^\infty w_{i,j} \right).$$
    \end{enumerate}
\end{prop}

\begin{proof}\phantom{M}

\begin{itemize}
    \item Fix $i \in \N$ and $r > 0$. By assumption, there exists a finite set $S_r \subset \N \times \N$ such that $w_{i,j} \in T^r A$ for all $(i,j) \in (\N \times \N) \setminus S_r$. Let
    $$N = \max\{n \in \N : \exists m \in \N \text{ such that } (n,m) \in S_r \text{ or } (m,n) \in S_r\}.$$
    Thus, for every $j > N$, we have $w_{i,j} \in T^r A$. Hence $\lim_{j\to\infty} w_{i,j} = 0$, and since $A$ is complete, by Claim~\ref{claim: series in a Novikov module 0}, we deduce that the series $\sum_{j=1}^\infty w_{i,j}$ converges. The same argument holds for the series $\sum_{i=1}^\infty w_{i,j}$ for every $j \in \N$.

    \item For every $i \in \N$, denote $W_i = \sum_{j=1}^\infty w_{i,j}$. Similarly, for every $j \in \N$, denote $W^j = \sum_{i=1}^\infty w_{i,j}$. First, we show that the series $\sum_{i=1}^\infty W_i$ converges.
    
    Let $r > 0$. There exists a finite set $S_r \subset \N \times \N$ such that $w_{i,j} \in T^r A$ for all $(i,j) \in (\N \times \N) \setminus S_r$. Let $N$ be defined as above. For every $i > N$ and every $j \in \N$, we have $w_{i,j} \in T^r A$, and therefore $W_i \in T^r A$. Hence $\lim_{i\to\infty} W_i = 0$, so by Claim~\ref{claim: series in a Novikov module 0}, the series $\sum_{i=1}^\infty W_i$ converges; denote its limit by $W$.
    
    We now show that the series $\sum_{j=1}^\infty W^j$ converges to $W$. Let $r > 0$. As before, there exists a finite set $S_r \subset \N \times \N$ such that $w_{i,j} \in T^r A$ for all $(i,j) \in (\N \times \N) \setminus S_r$, with $N$ defined as the maximum index in $S_r$. 
    
    As shown, for every $i > N$, we have $W_i \in T^r A$. Thus, for every $i_0 > N$:
    $$W - \sum_{i=1}^{i_0} W_i = \sum_{i=i_0 + 1}^\infty W_i \in T^r A.$$
    Note that for every $i_0 > N$:
    \begin{align*}
    \sum_{i=1}^{i_0} W_i &= \sum_{i=1}^N W_i + \sum_{i=N+1}^{i_0} W_i \\
    &= \sum_{i=1}^N \left( \sum_{j=1}^\infty w_{i,j} \right) + \sum_{i=N+1}^{i_0} W_i \\
    &= \sum_{i=1}^N \left( \left( \sum_{j=1}^N w_{i,j} \right) + \left( \sum_{j=N+1}^\infty w_{i,j} \right) \right) + \sum_{i=N+1}^{i_0} W_i \\
    &= \sum_{i=1}^N \sum_{j=1}^N w_{i,j} + \sum_{i=1}^N \left( \sum_{j=N+1}^\infty w_{i,j} \right) + \sum_{i=N+1}^{i_0} W_i.
    \end{align*}
    Since $w_{i,j} \in T^r A$ for $i > N$ or $j > N$, the last two terms in the summation belong to $T^r A$. Therefore:
    $$\sum_{i=1}^{i_0} W_i - \sum_{i=1}^N \sum_{j=1}^N w_{i,j} \in T^r A.$$
    Similarly, for every $j_0 > N$, we have:
    $$\sum_{j=1}^{j_0} W^j - \sum_{i=1}^N \sum_{j=1}^N w_{i,j} \in T^r A.$$
    In particular, for every $i_0, j_0 > N$, we have:
    $$\sum_{j=1}^{j_0} W^j - \sum_{i=1}^{i_0} W_i = \left( \sum_{j=1}^{j_0} W^j - \sum_{i,j=1}^N w_{i,j} \right) - \left( \sum_{i=1}^{i_0} W_i - \sum_{i,j=1}^N w_{i,j} \right) \in T^r A.$$
    It follows that for every $j_0 > N$:
    $$\sum_{j=1}^{j_0} W^j - W = \left( \sum_{j=1}^{j_0} W^j - \sum_{i=1}^{i_0} W_i \right) + \left( \sum_{i=1}^{i_0} W_i - W \right) \in T^r A,$$
    completing the proof.
\end{itemize}
\end{proof}

\begin{claim}\label{claim: countable linearity}
	Let $A$ be a complete Novikov module and let $L \fc A \to A$ be a $\Lambda_{\geq 0}$-linear endomorphism of $A$. Then for every sequence $(\lambda_i)_{i \in \N}$ in $A$ and every sequence $(\alpha_i)_{i \in \N}$ in $\Lambda_{\geq 0}$ with $\val(\alpha_i) \to +\infty$, we have $L(\sum_i \alpha_i \lambda_i) = \sum_i \alpha_i L(\lambda_i)$.
\end{claim}

\begin{proof}
	Since $\val(\alpha_i) \to +\infty$, Claim~\ref{coroll: series in a Novikov module} implies that the series $\sum_i \alpha_i \lambda_i$ and $\sum_i \alpha_i L(\lambda_i)$ converge. Let $r > 0$. There exists $N \in \N$ such that for every $i > N$ we have $\val(\alpha_i) > r$.
	
	Then $\sum_{i > N} \alpha_i \lambda_i \in T^r A$, and we can define $x_r = T^{-r} \sum_{i > N} \alpha_i \lambda_i \in A$. The $\Lambda_{\geq 0}$-linearity of $L$ implies that
	$$L\left(\sum_i \alpha_i \lambda_i\right) = L\left( T^r x_r + \sum_{i \leq N} \alpha_i \lambda_i \right) = T^r L(x_r) + \sum_{i \leq N} \alpha_i L(\lambda_i).$$
	Therefore,
	$$L\left(\sum_i \alpha_i \lambda_i\right) - \sum_{i \leq N} \alpha_i L(\lambda_i) \in T^r A.$$
	Hence, $L(\sum_i \alpha_i \lambda_i)$ is the limit of the series $\sum_i \alpha_i L(\lambda_i)$, as required.
\end{proof}

\subsubsection{Endomorphisms and Schauder bases for $\widehat{\Lambda_{\geq 0}^\infty}$}\label{sss: Schauder basis}

For a Novikov module $A$, let $A^\infty$ be the direct sum of countably infinitely many copies of $A$: $A^\infty = \bigoplus_{i \in \N} A$.

The module $\Lambda_{\geq 0}^\infty$ can be identified with the space of sequences having only finitely many nonzero terms:
$$\left\{ (a_1, a_2, a_3, \dots) : a_i \in \Lambda_{\geq 0}, \, \exists N \in \N \text{ s.t. } a_n = 0 \text{ for every } n \geq N \right\}.$$
Let us denote by $e_i \in \Lambda_{\geq 0}^\infty$ the $i$-th element of the standard basis; that is, $e_i$ has $1$ in the $i$-th coordinate and zero everywhere else. Then, every $a \in \Lambda_{\geq 0}^\infty$ is uniquely expressible as $a = \sum_{i \in \N} a_i e_i$, where the $a_i$ are the coordinates.

\begin{coroll}\label{coroll: description of completion}\phantom{M}
    \begin{itemize}
        \item The module $\widehat{\Lambda_{\geq 0}^\infty}$ can be identified with 
        $$\left\{ \sum_{i \in \N} a_i e_i : a_i \in \Lambda_{\geq 0}, \, \val(a_i) \to +\infty \right\}.$$
        \item $\Lambda_{\geq 0}^\infty$ is the submodule of $\widehat{\Lambda_{\geq 0}^\infty}$ consisting of sequences which have only finitely many nonzero terms. 
        \item The canonical map $c\fc \Lambda_{\geq0}^\infty\to\widehat{\Lambda_{\geq0}^\infty}$ equals the identity map on $\Lambda_{\geq 0}^\infty$ if we consider $\Lambda_{\geq 0}^\infty$ as the submodule of $\widehat{\Lambda_{\geq 0}^\infty}$ consisting of sequences which have only finitely many nonzero terms.
    \end{itemize}
\end{coroll}

Since $\Lambda_{\geq 0}^\infty$ is a free $\Lambda_{\geq 0}$-module with basis $(e_i)_{i \in \N}$, an endomorphism $L$ of $\Lambda_{\geq 0}^\infty$ is determined by the sequence $(L(e_i))_{i \in \N}$ of elements in $\Lambda_{\geq 0}^\infty$. In fact, the resulting map from the collection of endomorphisms to the collection of sequences is a bijection. The module $\widehat{\Lambda_{\geq 0}^\infty}$ is not a free module; however, its endomorphisms are described in the same way.
\begin{claim}\phantom{M}
\begin{enumerate}
    \item 
    The map
    $$\cE \fc \left\{\begin{array}{c} \text{Endomorphisms} \\ \text{of } \widehat{\Lambda_{\geq 0}^\infty} \end{array}\right\} \longrightarrow \left\{\begin{array}{c} \text{Sequences of} \\ \widehat{\Lambda_{\geq 0}^\infty} \text{ elements} \end{array}\right\}$$
    given by $L \mapsto (L(e_i))_{i \in \N}$ is a bijection.
    \item For every sequence $(b_i)_{i\in \N}$ of $\widehat{\Lambda_{\geq 0}^\infty}$ elements and for every sequence $(\alpha_i)_{i\in \N}$ of $\Lambda_{\geq 0}$ elements with $\val(\alpha_i) \to \infty$, we have
    $$\cE^{-1}\left((b_i)_{i\in \N}\right)\left(\sum_{i=1}^\infty \alpha_i e_i\right) = \sum_{i=1}^\infty \alpha_i b_i.$$
\end{enumerate}
\end{claim}

\begin{proof}\phantom{M}
\begin{enumerate}
    \item
        \textbf{Let us show that $\cE$ is injective:} Let $K, L$ be two endomorphisms of $\widehat{\Lambda_{\geq 0}^\infty}$ and assume that $K(e_i) = L(e_i)$ for every $i \in \N$. Let $x \in \widehat{\Lambda_{\geq 0}^\infty}$. Corollary~\ref{coroll: description of completion} implies that there exists a sequence $(a_i)_{i \in \N}$ of elements of $\Lambda_{\geq 0}$ with $\val(a_i) \to +\infty$ such that $x = \sum_{i \in \N} a_i e_i$. By Claim~\ref{claim: countable linearity}, we have:
        $$K(x) = K\left(\sum_{i \in \N} a_i e_i\right) = \sum_{i \in \N} a_i K(e_i) = \sum_{i \in \N} a_i L(e_i) = L\left(\sum_{i \in \N} a_i e_i\right) = L(x).$$
        This shows that $K = L$, and hence $\cE$ is injective.
        
        \textbf{Let us show that $\cE$ is surjective:} Let $(\lambda_i)_{i \in \N}$ be a sequence of elements in $\widehat{\Lambda_{\geq 0}^\infty}$. Consider the function $L \fc \widehat{\Lambda_{\geq 0}^\infty} \to \widehat{\Lambda_{\geq 0}^\infty}$ defined by $L(x) = \sum_{i \in \N} a_i \lambda_i$ for every $x \in \widehat{\Lambda_{\geq 0}^\infty}$, where $x = \sum_{i \in \N} a_i e_i$. The map $L$ is well-defined according to Corollary~\ref{coroll: series in a Novikov module} and is $\Lambda_{\geq 0}$-linear; thus, $L$ is an endomorphism of $\widehat{\Lambda_{\geq 0}^\infty}$. The image of $L$ under $\cE$ is the sequence $(\lambda_i)_{i \in \N}$ by construction. Thus, $\cE$ is surjective.
    \item Let $(b_i)_{i \in \N}$ be a sequence in $\widehat{\Lambda_{\geq 0}^\infty}$. By the surjectivity proof above, the endomorphism $L = \cE^{-1}((b_i)_{i \in \N})$ is defined exactly as the unique linear map such that $L(e_i) = b_i$ for every $i \in \N$. Given a sequence $(\alpha_i)_{i \in \N}$ in $\Lambda_{\geq 0}$ with $\val(\alpha_i) \to \infty$, the element $x = \sum_{i=1}^\infty \alpha_i e_i$ is a well-defined element of $\widehat{\Lambda_{\geq 0}^\infty}$. By the definition of $L$ used to establish surjectivity, we have:
    $$ \cE^{-1}\left((b_i)_{i \in \N}\right)\left(\sum_{i=1}^\infty \alpha_i e_i\right) = L\left(\sum_{i=1}^\infty \alpha_i e_i\right) = \sum_{i=1}^\infty \alpha_i L(e_i) = \sum_{i=1}^\infty \alpha_i b_i, $$
    where the second equality follows from the countable linearity (Claim~\ref{claim: countable linearity}) of endomorphisms on $\widehat{\Lambda_{\geq 0}^\infty}$.
\end{enumerate}
\end{proof}
  
Schauder bases are usually defined for normed vector spaces; here, we show that the definition extends to $\widehat{\Lambda_{\geq 0}^\infty}$. 

\begin{defin}
	A sequence $(b_i)_{i=1}^\infty$ of elements in $\widehat{\Lambda_{\geq 0}^\infty}$ is called a \textbf{Schauder basis} if for every $a \in \widehat{\Lambda_{\geq 0}^\infty}$ there is a unique sequence $(a_i)_{i=1}^\infty$ of elements in $\Lambda_{\geq 0}$ with $\val(a_i) \to +\infty$ such that $a = \sum_{i=1}^\infty a_i b_i$.
\end{defin}

\begin{lemma}\label{lemma: conditions for Schauder basis}
    Let $(b_i)_{i=1}^\infty$ be a sequence of elements in $\widehat{\Lambda_{\geq 0}^\infty}$. The following statements are equivalent:
    \begin{enumerate}
        \item $(b_i)_{i=1}^\infty$ is a Schauder basis.
        \item The following two conditions hold:
        \begin{itemize}
            \item For every sequence $(a_i)_{i \in \N}$ in $\Lambda_{\geq 0}$, if $\val(a_i) \to \infty$ and $\sum_i a_i b_i = 0$, then $a_i = 0$ for every $i \in \N$.
            \item For every $j \in \N$, there is a sequence $(\alpha_i^j)_{i \in \N}$ in $\Lambda_{\geq 0}$ such that $\val(\alpha_i^j) \to \infty$ and $e_j = \sum_i \alpha_i^j b_i$.
        \end{itemize}
        \item The endomorphism $L = \cE^{-1}((b_i)_{i=1}^\infty)$ is an isomorphism.
    \end{enumerate} 
\end{lemma}

\begin{proof}
    The direction $(i) \Rightarrow (ii)$ is immediate from the definition of a Schauder basis.
    
    $(ii) \Rightarrow (iii)$: Let $L = \cE^{-1}((b_i)_{i=1}^\infty)$. We show that $L$ is injective and surjective. 
    
    \textbf{Injectivity:} Suppose $x \in \ker L$. Write $x = \sum_i a_i e_i$ with $a_i \in \Lambda_{\geq 0}$ and $\val(a_i) \to \infty$. Then $L(x) = \sum_i a_i b_i = 0$ by countable linearity. By the first condition of (2), $a_i = 0$ for all $i$, so $x = 0$.
    
    \textbf{Surjectivity:} Let $y = \sum_{j \in \N} a_j e_j \in \widehat{\Lambda_{\geq 0}^\infty}$ with $\val(a_j) \xrightarrow{j \to \infty} \infty$. By the second condition of (2), for each $j \in \N$, there is a sequence $(\alpha_i^j)_{i \in \N}$ such that $\val(\alpha_i^j) \xrightarrow{i \to \infty} \infty$ and $e_j = \sum_i \alpha_i^j b_i$.
    
    First, for every $i \in \N$, since $\val(a_j) \xrightarrow{j \to \infty} \infty$, the series $\beta_i = \sum_{j=1}^\infty a_j \alpha_i^j$ converges by Corollary~\ref{coroll: series in a Novikov module}. Let us show that $\val(\beta_i) \to \infty$.

    Let $r > 0$. 
    \begin{itemize}
        \item Since $\val(a_j) \to \infty$, there exists $J \in \N$ such that for every $j > J$, there exists $c_j \in \Lambda_{\geq 0}$ such that $a_j = T^r c_j$. Therefore, for every $j > J$ and $i \in \N$, we have
        $$\alpha_i^j a_j = T^r \alpha_i^j c_j \in T^r \Lambda_{\geq 0}.$$

        \item For each $j \in \{1, \dots, J\}$, we have $\val(\alpha_i^j) \to \infty$ as $i \to \infty$. Thus, there exists $N \in \N$ such that for all $i > N$ and all $j \leq J$, there exists $\gamma_i^j \in \Lambda_{\geq 0}$ such that $\alpha_i^j = T^r \gamma_i^j$. Therefore, for every $i > N$ and $j \leq J$, we have
        $$\alpha_i^j a_j = T^r \gamma_i^j a_j \in T^r \Lambda_{\geq 0}.$$
    \end{itemize}
    Therefore, for every $(i,j) \in (\N \times \N) \setminus (\{1, \dots, N\} \times \{1, \dots, J\})$, we have $a_j \alpha_i^j \in T^r \Lambda_{\geq 0}$.

    In particular, for every $i > N$ and every $j \in \N$, we have $a_j \alpha_i^j \in T^r \Lambda_{\geq 0}$, and hence $\beta_i = \sum_{j=1}^\infty a_j \alpha_i^j \in T^r \Lambda_{\geq 0}$. This shows that $\val(\beta_i) \to \infty$.

    Hence, by Corollary~\ref{coroll: series in a Novikov module}, the series $\sum_{i=1}^\infty \beta_i e_i$ converges. Denote its limit by $x \in \widehat{\Lambda_{\geq 0}^\infty}$. We will prove that $y = L(x)$, which concludes that $L$ is surjective.
    
    For every $(i,j) \in (\N \times \N) \setminus (\{1, \dots, N\} \times \{1, \dots, J\})$, we have $a_j \alpha_i^j \in T^r \Lambda_{\geq 0}$, therefore $\alpha_i^j a_j b_i \in T^r \widehat{\Lambda_{\geq 0}^\infty}$. By Proposition~\ref{prop: rearrangement}, we deduce that:
    $$y = \sum_{j=1}^\infty a_j e_j = \sum_{j=1}^\infty a_j \left( \sum_{i=1}^\infty \alpha_i^j b_i \right) = \sum_{j=1}^\infty \sum_{i=1}^\infty a_j \alpha_i^j b_i = \sum_{i=1}^\infty \left( \sum_{j=1}^\infty a_j \alpha_i^j \right) b_i = \sum_{i=1}^\infty \beta_i b_i = L(x),$$
    as required.
    
    $(iii) \Rightarrow (i)$: Assume $L$ is an isomorphism. For any $a \in \widehat{\Lambda_{\geq 0}^\infty}$, there exists a unique $x = \sum a_i e_i$ with $\val(a_i) \to \infty$ such that $L(x) = a$. Then $a = \sum a_i L(e_i) = \sum a_i b_i$. The uniqueness and existence of $(a_i)$ follow from the bijectivity of $L$.
\end{proof}

We present below examples of Schauder bases which will turn out to be very useful in the proofs of Theorem~\ref{thm: SH of ball in CP^n} and Theorem~\ref{thm: res for balls in CP^n}.

\begin{exam}\label{exam: Schauder bases}\phantom{M}
    \begin{enumerate}
        \item Any basis of $\Lambda_{\geq 0}^\infty$, viewed as a sequence in $\widehat{\Lambda_{\geq 0}^\infty}$, is a Schauder basis. Indeed, a basis of $\Lambda_{\geq 0}^\infty$ induces an isomorphism of $\Lambda_{\geq 0}^\infty$, mapping the $i$-th element of the standard basis to the $i$-th basis element. Consequently, it induces an isomorphism on $\widehat{\Lambda_{\geq 0}^\infty}$, which corresponds to a Schauder basis of $\widehat{\Lambda_{\geq 0}^\infty}$ by Lemma~\ref{lemma: conditions for Schauder basis}. A particularly important example is given by the basis:
        $$(e_1, e_1 + e_2, e_2 + e_3, \dots).$$

        \item This can be generalized as follows. For $\alpha \geq 0$, let $\lambda_1 = e_1$, and for $i \geq 2$, let $\lambda_i = T^{\alpha} e_{i-1} + e_i$. The first example is the special case where $\alpha = 0$. We claim that $(\lambda_i)_{i \in \N}$ is a basis of $\Lambda_{\geq 0}^\infty$. Linear independence is immediate. Furthermore, for each $i \in \N$, we have:
        $$e_i = \sum_{j=0}^{i-1} (-T^\alpha)^j \lambda_{i-j} = \lambda_i - T^\alpha \lambda_{i-1} + T^{2\alpha} \lambda_{i-2} + \dots + (-T^\alpha)^{i-1} \lambda_1.$$
        This shows that the $(\lambda_i)_{i \in \N}$ span $\Lambda_{\geq 0}^\infty$ and thus constitute a basis. In particular, they also form a Schauder basis of $\widehat{\Lambda_{\geq 0}^\infty}$.

        \item Let $\alpha > 0$. The sequence $(\epsilon_j)_{j \in \N}$ given by $\epsilon_j = e_j + T^\alpha e_{j+1}$ is a Schauder basis. Indeed, let $(\alpha_j)_{j \in \N}$ be a sequence of elements in $\Lambda_{\geq 0}$ such that $\val(\alpha_j) \to +\infty$ and $\sum_j \alpha_j \epsilon_j = 0$. Then:
        $$0 = \sum_{j=1}^\infty \alpha_j (e_j + T^\alpha e_{j+1}) = \alpha_1 e_1 + \sum_{j=2}^\infty (\alpha_j + \alpha_{j-1} T^\alpha) e_j.$$
        Since $(e_j)_{j \in \N}$ is a Schauder basis and $\val(\alpha_j + \alpha_{j-1} T^\alpha) \to +\infty$, the uniqueness property implies $\alpha_1 = 0$ and $\alpha_j + \alpha_{j-1} T^\alpha = 0$ for each $j \geq 2$. Proceeding inductively, we see that $\alpha_j = 0$ for all $j \in \N$. 
        
        Furthermore, for every $i \in \N$, we have:
        $$e_i = \sum_{j=i}^\infty (-T^\alpha)^{j-i} \epsilon_j.$$
        Since $\alpha > 0$, the valuations of the coefficients $(-T^\alpha)^{j-i}$ tend to infinity as $j \to \infty$. By Lemma~\ref{lemma: conditions for Schauder basis}, $(\epsilon_i)_{i \in \N}$ is indeed a Schauder basis. Note, however, that no $e_i$ can be expressed as a finite linear combination of the $\epsilon_j$'s. Thus, despite belonging to $\Lambda_{\geq 0}^\infty$ and being linearly independent, $(\epsilon_i)_{i \in \N}$ does not form a basis of $\Lambda_{\geq 0}^\infty$.
    \end{enumerate}
\end{exam}

Let us turn our attention to elements of strictly positive valuation. Note that since a linear endomorphism of $\widehat{\Lambda_{\geq0}^\infty}$ is given
by a matrix with entries in $\Lambda_{\geq0}$, the multiplication by such a matrix preserves $\widehat{\Lambda_{>0}^\infty}$. It follows that an automorphism of $\widehat{\Lambda_{\geq0}^\infty}$ restricts to an automorphism of $\widehat{\Lambda_{>0}^\infty}$. Therefore a Schauder basis $(b_i)_{i\in\N}$ for $\widehat{\Lambda_{\geq0}^\infty}$ determines such an automorphism of $\widehat{\Lambda_{>0}^\infty}$ via $a = (a_1,a_2,\ldots)\mapsto\sum_i a_ib_i$. Similar reasoning shows that this automorphism further restricts to an automorphism of $\widehat{\Lambda_{>\alpha}^\infty}$ for any $\alpha>0$.

Given $\alpha, \beta > 0$, there is a unique endomorphism $\tilde{D} \fc \widehat{\Lambda_{\geq 0}^\infty} \to \widehat{\Lambda_{\geq 0}^\infty}$ that satisfies 
$$\tilde{D}(e_i) = T^\alpha e_i + T^\beta e_{i+1}$$
for every $i \in \N$. Denote by $D$ the restriction of $\tilde{D}$ to $\widehat{\Lambda_{>0}^\infty}$. Thus, $D$ satisfies
\begin{equation}\label{eq: diff of drct CF(Ball)}
    D(\lambda e_i) = T^\alpha \lambda e_i + T^\beta \lambda e_{i+1}
\end{equation}
for every $i \in \N$ and $\lambda\in\Lambda_{>0}$. Our computation of the relative symplectic cohomology of the ball, and the proof of Theorem~\ref{thm: SH of ball in CP^n}, is based on the following computation of the cokernel of $D$. We denote by $q\fc \widehat{\Lambda_{>0}^\infty}\to\coker D$ the natural quotient map associated to $D$.

\begin{prop}\label{prop: coker of D} $D$ is injective. In addition:
    \begin{itemize}
        \item If $\alpha < \beta$, then  $\coker D$ is isomorphic to $\Lambda_{(0,\alpha]}^\infty$. Moreover, the natural quotient map $q\fc \widehat{\Lambda_{>0}^\infty}\to \coker D$ satisfies
        $$q(\lambda \epsilon_i)=\overline{\lambda} \epsilon_i,$$
        for every $i\in \N$ and $\lambda\in \Lambda_{>0}$, where $\overline\lambda$ is the image of $\lambda$ under the quotient $\Lambda_{>0}\to\Lambda_{(0,\alpha]}=\Lambda_{>0}/\Lambda_{>\alpha}$, and $(\epsilon_i)_{i=1}^\infty$ is the Schauder basis of $\Lambda_{\geq0}^\infty$ given by $\epsilon_i=e_i+T^{\beta-\alpha}e_{i+1}$ for every $i\in \N$.

        \item If $\alpha \geq \beta$, then $\coker D$ is isomorphic to  $\Lambda_{>0} \oplus \Lambda_{(0,\beta]}^\infty$. Moreover, the natural quotient map $q\fc \widehat{\Lambda_{>0}^\infty}\to \coker D$ satisfies
        $q(\lambda \epsilon_1)=\lambda \epsilon_1$ for every $\lambda\in\Lambda_{>0}$ and 
    $$q(\lambda \epsilon_i)=\overline{\lambda} \epsilon_i,$$
for every $i\geq2$ and $\lambda\in \Lambda_{>0}$, where $\overline\lambda$ is the image of $\lambda$ under the quotient $\Lambda_{>0}\to\Lambda_{(0,\beta]}=\Lambda_{>0}/\Lambda_{>\beta}$, and $(\epsilon_i)_{i=1}^\infty$ is the Schauder basis of $\Lambda_{\geq0}^\infty$ given by $\epsilon_1=e_1$ and  $\epsilon_i=T^{\alpha-\beta}e_{i-1}+e_i$ for every $i\geq 2$.
    
    \end{itemize}
    
\end{prop}

The following results constitute a preparation for the proof of Proposition~\ref{prop: coker of D}.

\begin{lemma}\label{lemma: quotient of completes modules 1}
    For every $\alpha > 0$, $\widehat{\Lambda_{>0}^\infty} / \widehat{\Lambda_{>\alpha}^\infty}$ is isomorphic to  $\Lambda_{(0,\alpha]}^\infty$.
\end{lemma}

\begin{proof}
    Let $\alpha > 0$.
    Let us show that $\Lambda_{>0}^\infty + \widehat{\Lambda_{>\alpha}^\infty} = \widehat{\Lambda_{>0}^\infty}$. The inclusion $\Lambda_{>0}^\infty + \widehat{\Lambda_{>\alpha}^\infty} \subset \widehat{\Lambda_{>0}^\infty}$ is obvious; let us prove the opposite inclusion.
    Let $a \in \widehat{\Lambda_{>0}^\infty}$. There is a sequence $(a_i)_{i=1}^\infty$ of elements in $\Lambda_{>0}$ with $\val(a_i) \to +\infty$ such that $a = \sum_{i \geq 1} a_i e_i$. There exists $j_0 \in \N$ such that for every $i > j_0$, we have $\val(a_i) > \alpha$. Denote $b = \sum_{1 \leq i \leq j_0} a_i e_i$ and $c = \sum_{i > j_0} a_i e_i$. Then $b \in \Lambda_{>0}^\infty$, $c \in \widehat{\Lambda_{>\alpha}^\infty}$, and we have $a = b + c$.
    
    Additionally, note that the elements of $\Lambda_{>0}^\infty$ which also belong to $\widehat{\Lambda_{>\alpha}^\infty}$ are exactly those for which the valuations of all coordinates are $>\alpha$. This is precisely the submodule $\Lambda_{>\alpha}^\infty$. 
    
    Thus, the second isomorphism theorem implies that
    $$\widehat{\Lambda_{>0}^\infty} / \widehat{\Lambda_{>\alpha}^\infty} = (\Lambda_{>0}^\infty + \widehat{\Lambda_{>\alpha}^\infty}) / \widehat{\Lambda_{>\alpha}^\infty} \simeq \Lambda_{>0}^\infty / (\Lambda_{>0}^\infty \cap \widehat{\Lambda_{>\alpha}^\infty}) = \Lambda_{>0}^\infty / \Lambda_{>\alpha}^\infty.$$
    Since the quotient commutes with the direct sum, we obtain
    \[\widehat{\Lambda_{>0}^\infty} / \widehat{\Lambda_{>\alpha}^\infty} \simeq \Lambda_{>0}^\infty / \Lambda_{>\alpha}^\infty \simeq (\Lambda_{>0} / \Lambda_{>\alpha})^\infty = \Lambda_{(0,\alpha]}^\infty.\qedhere\]
\end{proof}
\begin{lemma}\label{lemma: quotient of completes modules 2}
    For every $\beta > 0$, the quotient $\widehat{\Lambda_{>0}^\infty} / \widehat{(0 \oplus \Lambda_{>\beta}^\infty)}$ is isomorphic to $\Lambda_{>0} \oplus \Lambda_{(0,\beta]}^\infty$, where in the denominator, the zero is contained in the first $\Lambda_{>0}$ factor of the infinite direct sum:
    $$0 \oplus \Lambda_{>\beta}^\infty = 0 \oplus \Lambda_{>\beta} \oplus \Lambda_{>\beta} \oplus \dots \;\subset\; \Lambda_{>0} \oplus \Lambda_{>0} \oplus \Lambda_{>0} \oplus \dots = \Lambda_{>0}^\infty.$$
\end{lemma}

\begin{proof}
    The module $\Lambda_{>0}$ is complete. Since completion commutes with finite direct sums, we have:
    $$\widehat{\Lambda_{>0}^\infty}=\widehat{\Lambda_{>0}\oplus \Lambda_{>0}^\infty}=\widehat{(\Lambda_{>0}\oplus0)\oplus(0\oplus \Lambda_{>0}^\infty)}=\widehat{\Lambda_{>0}\oplus0}\oplus \widehat{0 \oplus\Lambda_{>0}^\infty}=\widehat{\Lambda_{>0}}\oplus \widehat{\Lambda_{>0}^\infty}=\Lambda_{>0}\oplus \widehat{\Lambda_{>0}^\infty}.$$
    Similarly,
    $$\widehat{0 \oplus \Lambda_{>\beta}^\infty} = 0 \oplus \widehat{\Lambda_{>\beta}^\infty}.$$
    Thus, Lemma~\ref{lemma: quotient of completes modules 1} implies that
    \[\widehat{\Lambda_{>0}^\infty} / \widehat{0 \oplus \Lambda_{>\beta}^\infty} = (\Lambda_{>0} \oplus \widehat{\Lambda_{>0}^\infty}) / (0 \oplus \widehat{\Lambda_{>\beta}^\infty}) \simeq (\Lambda_{>0} / 0) \oplus (\widehat{\Lambda_{>0}^\infty} / \widehat{\Lambda_{>\beta}^\infty}) \simeq \Lambda_{>0} \oplus \Lambda_{(0,\beta]}^\infty. \qedhere\]
\end{proof}\begin{proof}[Proof of Proposition~\ref{prop: coker of D}]\phantom{M}
\begin{itemize}
    \item Assume that $\alpha < \beta$. In this case, we use the Schauder basis $(\epsilon_i)_{i \in \N}$ given by $\epsilon_i = e_i + T^{\beta-\alpha} e_{i+1}$ for every $i \in \N$; see Example~\ref{exam: Schauder bases}. For every $i\in \N$ and $\lambda\in \Lambda_{>0}$ we have
    $$D(\lambda e_i)=T^\alpha \lambda e_i+T^\beta \lambda e_{i+1}=T^\alpha \left(\lambda (e_i+T^{\beta-\alpha} e_{i+1})\right)=T^\alpha \lambda \epsilon_i.$$
    Therefore 
    $$\im D=\widehat{\bigoplus_{i=1}^\infty \Lambda_{>\alpha}\epsilon_i},$$
    and hence it follows from Lemma~\ref{lemma: quotient of completes modules 1} that
    $$\coker D=\widehat{\Lambda_{>0}^\infty}/\im D=\widehat{\bigoplus_{i=1}^\infty \Lambda_{>0}\epsilon_i}\Big/\widehat{\bigoplus_{i=1}^\infty \Lambda_{>\alpha}\epsilon_i}=\bigoplus_{i=1}^\infty \Lambda_{(0,\alpha]}\epsilon_i=\Lambda_{(0,\alpha]}^\infty.$$
    Moreover, the quotient map $q\fc \widehat{\bigoplus_{i=1}^\infty\Lambda_{> 0} \epsilon_i}\to\bigoplus_{i=1}^\infty \Lambda_{(0,\alpha]}\epsilon_i$  satisfies
    $$q(\lambda \epsilon_i)=\overline{\lambda} \epsilon_i,$$
for every $i\in \N$ and $\lambda\in \Lambda_{>0}$, where $\overline\lambda$ is the image of $\lambda$ under the quotient $\Lambda_{>0}\to\Lambda_{(0,\alpha]}=\Lambda_{>0}/\Lambda_{>\alpha}$.

    \item Assume that $\alpha \geq \beta$. In this case, we use the Schauder basis $(\epsilon_i)_{i \in \N}$ given by $\epsilon_1 = e_1$ and $\epsilon_i =  T^{\alpha-\beta} e_{i-1}+e_i $ for every $i \geq 2$; see Example~\ref{exam: Schauder bases}. As before,  for every $i\in \N$ and $\lambda\in \Lambda_{>0}$ we have
    $$D(\lambda e_i)=T^\alpha \lambda e_i+T^\beta \lambda e_{i+1}=T^\beta \left(\lambda (T^{\alpha-\beta} e_i+e_{i+1})\right)=T^\beta \lambda \epsilon_{i+1}.$$
    Therefore 
    $$\im D=\widehat{\bigoplus_{i=2}^\infty \Lambda_{>\beta}\epsilon_i},$$
    and hence it follows from Lemma~\ref{lemma: quotient of completes modules 2} that
    $$\coker D=\widehat{\Lambda_{>0}^\infty}/\im D=\widehat{\bigoplus_{i=1}^\infty \Lambda_{>0}\epsilon_i}\Big/\widehat{\bigoplus_{i=2}^\infty \Lambda_{>\beta}\epsilon_i}=\Lambda_{>0}\epsilon_1\oplus\bigoplus_{i=2}^\infty \Lambda_{(0,\beta]}\epsilon_i=\Lambda_{>0}\oplus\Lambda_{(0,\beta]}^\infty.$$
    Moreover, the quotient map $q\fc \widehat{\bigoplus_{i=1}^\infty\Lambda_{> 0}\epsilon_i}\to\Lambda_{>0}\epsilon_1\oplus\bigoplus_{i=2}^\infty \Lambda_{(0,\beta]}\epsilon_i$  satisfies that $q(\lambda \epsilon_1)=\lambda \epsilon_1$ for every $\lambda\in\Lambda_{>0}$ and 
    $$q(\lambda \epsilon_i)=\overline{\lambda} \epsilon_i,$$
for every $i\geq2$ and $\lambda\in \Lambda_{>0}$, where $\overline\lambda$ is the image of $\lambda$ under the quotient $\Lambda_{>0}\to\Lambda_{(0,\beta]}=\Lambda_{>0}/\Lambda_{>\beta}$.   \qedhere
 
\end{itemize} 
\end{proof}

The computation of the restriction maps associated with the relative symplectic cohomology of balls involves several independent calculations. The first of these is Proposition~\ref{prop: surjectivity and completion}, presented below.

Let $c \fc \Lambda_{>0}^\infty \to \widehat{\Lambda_{>0}^\infty}$ be the canonical map to the completion. Let us define an uncompleted morphism parallelling $D$ from Equation \eqref{eq: diff of drct CF(Ball)}. There is a unique endomorphism $\tilde{d} \fc \Lambda_{\geq 0}^\infty \to \Lambda_{\geq 0}^\infty$ that satisfies 
$$\tilde{d}(e_i) = T^\alpha e_i + T^\beta e_{i+1}$$
for every $i \in \N$. Denote by $d$ the restriction of $\tilde{d}$ to $\Lambda_{>0}^\infty$. Thus $d$ satisfies
$$d(\lambda e_i) = T^\alpha \lambda e_i + T^\beta \lambda e_{i+1}$$ 
for every $i \in \N$ and $\lambda \in \Lambda_{>0}$. 

Since completion is a functor and $D$ is the completion of $d$, the following diagram commutes:
$$
\xymatrix@R=2pc@C=3pc{
  \Lambda_{>0}^\infty \ar[r]^d \ar[d]_{c} & \Lambda_{>0}^\infty \ar[d]^{c} \\
  \widehat{\Lambda_{>0}^\infty} \ar[r]^{D} & \widehat{\Lambda_{>0}^\infty} 
}
$$
Consequently, the universal property of the cokernel implies the existence of a unique morphism $\bar{c} \fc \coker d \to \coker D$ making the following diagram commute:

$$
\xymatrix@R=2pc@C=3pc{
  \Lambda_{>0}^\infty \ar[r]^d \ar[d]_{c} & \Lambda_{>0}^\infty \ar[d]_{c} \ar[r]^{p} & \coker d \ar[d]_{\bar{c}} \\
  \widehat{\Lambda_{>0}^\infty} \ar[r]^{D} & \widehat{\Lambda_{>0}^\infty} \ar[r]^{q} & \coker D
}
$$
where $p$ and $q$ are the canonical quotient maps.

\begin{prop}\label{prop: surjectivity and completion}
  The map $\bar{c} \fc \coker d \to \coker D$ is surjective.
\end{prop}
\begin{proof}
    There are two cases to check.
    \begin{itemize}
        \item Assume that $\alpha<\beta$. Let $x\in \coker D$, and let us show that there exists $y\in \coker d$ such that $\bar{c}(y)=x$.
        
        As in the proof of Proposition~\ref{prop: coker of D}, consider the Schauder basis $(\epsilon_i)_{i \in \N}$ given by $\epsilon_i = e_i + T^{\beta-\alpha} e_{i+1}$ for every $i \in \N$. The proof of Proposition~\ref{prop: coker of D} shows that $\coker D$ is generated by $(\epsilon_i)_{i\in \N}$ in the sense that
        $$\coker D=\bigoplus_{i\geq1} (\Lambda_{(0,\alpha]}\epsilon_i).$$
        Since $x\in \coker D$, there exists $N\in \N$ and $\lambda_1,\ldots,\lambda_N\in \Lambda_{(0,\alpha]}$ such that $x=\sum_{i=1}^N \lambda _i \epsilon_i$.  Since $\Lambda_{(0,\alpha]}=\Lambda_{>0}/\Lambda_{>\alpha}$, there exist $\eta_1,\ldots,\eta_N\in \Lambda_{>0}$ with finite valuation, such that the quotient map $\Lambda_{>0}\to \Lambda_{>0}/\Lambda_{>\alpha}=\Lambda_{(0,\alpha]}$ sends $\eta_i$ to $\lambda_i$ for every $1\leq i\leq N$. Since for every $1\leq i\leq N$ we have $\eta_i\epsilon_i\in \Lambda_{>0}^\infty$, we deduce that $z=\sum_{i=1}^N \eta_i \epsilon_i$ is an element of $\Lambda_{>0}^\infty\subset \widehat{\Lambda_{>0}^\infty}$. Consider $y=p(z)$. Note that $x=q(z)$ and $c(z)=z$. Since $q\circ c=\bar{c}\circ p$ we deduce that 
        $$\bar{c}(y)=\bar{c}\circ p(z)=q\circ c(z)=q(z)=x.$$
        This shows that $\bar{c}$ surjective.
        
        \item Assume that $\alpha \geq \beta$.  Let $x\in \coker D$, and let us show that there exists $y\in \coker d$ such that $\bar{c}(y)=x$.
        
        As in the proof of Proposition~\ref{prop: coker of D}, consider the Schauder basis $(\epsilon_i)_{i \in \N}$ given by $\epsilon_1=e_1$ and $\epsilon_i = T^{\alpha-\beta}e_{i-1} + e_{i}$ for every $i \geq 2$. The proof of Proposition~\ref{prop: coker of D} shows that $\coker D$ is generated by $(\epsilon_i)_{i\in \N}$ in the sense that
        $$\coker D=\Lambda_{>0}\epsilon_1 \oplus\bigoplus_{i\geq2} (\Lambda_{(0,\beta]}\epsilon_i).$$
        Since $x\in \coker D$, there exists $N\in \N$, $\lambda_1\in \Lambda_{>0}$ and $\lambda_2,\ldots,\lambda_N\in \Lambda_{(0,\beta]}$ such that $x=\sum_{i=1}^N \lambda _i \epsilon_i$.  There exist $\eta_2,\ldots,\eta_N\in \Lambda_{>0}$ with finite valuation, such that the quotient map $\Lambda_{>0}\to \Lambda_{>0}/\Lambda_{>\beta}=\Lambda_{(0,\beta]}$ sends $\eta_i$ to $\lambda_i$ for every $2\leq i\leq N$. Since for every $2\leq i\leq N$ we have $\eta_i\epsilon_i\in \Lambda_{>0}^\infty$, and $\lambda_1\epsilon_1\in\Lambda_{>0}$, we deduce that $z=\lambda_1 \epsilon_1 +\sum_{i=2}^N \eta_i \epsilon_i$ is an element of $\Lambda_{>0}^\infty\subset \widehat{\Lambda_{>0}^\infty}$. Consider $y=p(z)$. Note that $x=q(z)$ and $c(z)=z$. Since $q\circ c=\bar{c}\circ p$ we deduce that 
        $$\bar{c}(y)=\bar{c}\circ p(z)=q\circ c(z)=q(z)=x.$$
        This shows that $\bar{c}$ surjective. \qedhere
    \end{itemize}
\end{proof}

For every $a>0$ and $\eta\in \Lambda_{>0}$ we denote by $[\eta]_a$ the image of $\eta$ under the quotient map $\Lambda_{>0}\to \Lambda_{>0}/\Lambda_{>a}=\Lambda_{(0,a]}$. Also, for every $a>0$ there is a unique section $s_a\fc\Lambda_{(0,a]}\to \Lambda_{>0}$ so that for every $\lambda\in \Lambda_{(0,a]}\setminus\{0\}$, the of exponents of $T$ appearing in $s_a(\lambda)$ with a nonzero coefficient is supported in $(0,a]$. Note that this is not a Novikov module morphism, but rather, a group homomorpism. Additionally, for every $a,b>0$ we can define the group homomorphism $[\cdot]_{a,b}\fc\Lambda_{(0,a]}\to \Lambda_{(0,b]}$ by $[\lambda]_{a,b}=[s_a(\lambda)]_b$ for every $\lambda\in \Lambda_{(0,a]}$. Note that for every $a,b>0$ and $\gamma\geq0$, the map $T^\gamma\cdot[\cdot]_{a,b}$ is a Novikov module morphism if and only if $\gamma\geq b-a$.

The following, models the algebraic computation occurring when restricting from $B_{\Delta}$ to $B_{\Delta'}$ for $\Delta' \le \Delta$. Hence, the intended use case of what follows, is with $\alpha' = \Delta' \;\le\; \Delta = \alpha$, and $\beta' = n(1-\Delta') \;\ge\; n(1-\Delta) = \beta$. This is applied in the proof of Theorem \ref{thm: res for balls in CP^n} in Section \ref{ss: res maps} with $\alpha,\alpha',\beta,\beta'$ corresponding to topological energies.

Fix $n\in \N$ and $0\leq j\leq n$. Given $\alpha', \beta' > 0$ satisfying $\alpha'\leq \alpha$, $\beta'\geq\beta$ and
\begin{equation}
    n(\alpha-\alpha')=\beta'-\beta,
\end{equation}
consider the unique endomorphism $\tilde{D}' \fc \widehat{\Lambda_{\geq 0}^\infty} \to \widehat{\Lambda_{\geq 0}^\infty}$ satisfying 
$$\tilde{D}'(e_i) = T^{\alpha'} e_i + T^{\beta'} e_{i+1}$$
for every $i \in \N$. Let $D'$ be the restriction of $\tilde{D}'$ to $\widehat{\Lambda_{>0}^\infty}$. 

Additionally, consider the unique endomorphisms $\tilde\Psi_1,\tilde\Psi_2 \fc \widehat{\Lambda_{\geq 0}^\infty} \to \widehat{\Lambda_{\geq 0}^\infty}$ satisfying
$$\tilde{\Psi}_1(e_i) = T^{\hbigl(j+(i-1)(n+1)+1\hbigr)(\alpha-\alpha')} e_i,\qquad \tilde{\Psi}_2(e_i) = T^{\hbigl(j+(i-1)(n+1)\hbigr)(\alpha-\alpha')} e_i,$$
for every $i \in \N$. Let $\Psi_1,\Psi_2$ denote the restrictions of $\tilde\Psi_1,\tilde\Psi_2$ to $\widehat{\Lambda_{>0}^\infty}$, which satisfy
\begin{equation}\label{eq: cont between drct CF(Ball)}
    \Psi_1(\lambda e_i) = T^{\hbigl(j+(i-1)(n+1)+1\hbigr)(\alpha-\alpha')} \lambda e_i,\qquad \Psi_2(\lambda e_i) = T^{\hbigl(j+(i-1)(n+1)\hbigr)(\alpha-\alpha')} \lambda e_i
\end{equation}
for every $i \in \N$ and $\lambda\in\Lambda_{>0}$. 

The computation of the restriction maps between the relative symplectic cohomology of the balls is based on the induced map between the cokernels of $D$ and $D'$. First, we observe that the following diagram commutes:
\begin{equation}\label{diag: D' Psi1=D Psi2}
    \vcenter{
    \xymatrix@R=2pc@C=3pc{
    \widehat{\Lambda_{>0}^\infty} \ar[r]^D \ar[d]_{\Psi_1} &  \widehat{\Lambda_{>0}^\infty}  \ar[d]^{\Psi_2} \\
    \widehat{\Lambda_{>0}^\infty}  \ar[r]^{D'} &  \widehat{\Lambda_{>0}^\infty} 
}}
\end{equation}
Indeed, for every $i \in \N$ and $\lambda\in\Lambda_{>0}$, we have
\begin{align*}
    \Psi_2 \circ D(\lambda e_i) &= \Psi_2\left(T^\alpha\lambda e_i + T^\beta \lambda e_{i+1}\right) \\
    &= T^{\hbigl(j+(i-1)(n+1)\hbigr)(\alpha-\alpha')} T^\alpha \lambda e_i + T^{\hbigl(j+i(n+1)\hbigr)(\alpha-\alpha')} T^\beta \lambda e_{i+1} \\
    &= T^{\hbigl(j+(i-1)(n+1)\hbigr)(\alpha-\alpha')+\alpha}\lambda e_i + T^{\hbigl(j+i(n+1)\hbigr)(\alpha-\alpha')+\beta} \lambda e_{i+1}.
\end{align*}
Conversely, applying $D' \circ \Psi_1$ yields
\begin{align*}
    D' \circ \Psi_1(\lambda e_i) &= D'\left(T^{\hbigl(j+(i-1)(n+1)+1\hbigr)(\alpha-\alpha')}\lambda e_i\right) \\
    &= T^{\alpha'}T^{\hbigl(j+(i-1)(n+1)+1\hbigr)(\alpha-\alpha')}\lambda e_i + T^{\beta'}T^{\hbigl(j+(i-1)(n+1)+1\hbigr)(\alpha-\alpha')}\lambda e_{i+1} \\
    &= T^{\hbigl(j+(i-1)(n+1)\hbigr)(\alpha-\alpha')+\alpha-\alpha'+\alpha'}\lambda e_i + T^{\hbigl(j+i(n+1)\hbigr)(\alpha-\alpha')-n(\alpha-\alpha')+\beta'}\lambda e_{i+1}.
\end{align*}
Using the relation $n(\alpha-\alpha') = \beta'-\beta$, the second term becomes
$$ T^{\hbigl(j+i(n+1)\hbigr)(\alpha-\alpha') - (\beta'-\beta) + \beta'} \lambda e_{i+1} = T^{\hbigl(j+i(n+1)\hbigr)(\alpha-\alpha') + \beta} \lambda e_{i+1}. $$
This confirms that $D' \circ \Psi_1(\lambda e_i) = \Psi_2 \circ D(\lambda e_i)$, completing the proof that of the commutation of the diagram from Equation~\eqref{diag: D' Psi1=D Psi2}..
The commutativity of the diagram in Equation~\eqref{diag: D' Psi1=D Psi2} and the universal property of the cokernel imply the existence of a unique morphism $\overline{\Psi} \fc \coker D \to \coker D'$ such that the following diagram commutes:

\begin{equation}\label{diag: alg res maps}
  \vcenter{
  \xymatrix@R=2pc@C=3pc{
    \widehat{\Lambda_{>0}^\infty} \ar[r]^D \ar[d]_{\Psi_1} & \widehat{\Lambda_{>0}^\infty} \ar[d]_{\Psi_2} \ar[r]^{q} & \coker D \ar[d]_{\overline{\Psi}} \\
    \widehat{\Lambda_{>0}^\infty} \ar[r]^{D'} & \widehat{\Lambda_{>0}^\infty} \ar[r]^{q'} & \coker D'
  }}
\end{equation}
where $q,q'$ are the quotient maps. Finally, we compute $\overline{\Psi}$.
\begin{prop}\label{prop: alg restriction maps}
    The map $\overline{\Psi} \fc \coker D \to \coker D'$ is given as follows:
    \begin{enumerate}
        \item If $\alpha<\beta$ then $\alpha'<\beta'$ and the map
        $$\overline{\Psi} \fc \bigoplus_{i\geq1}\Lambda_{(0,\alpha]}\epsilon_i\to \bigoplus_{i\geq1}\Lambda_{(0,\alpha']}\epsilon'_i$$
        satisfies
        $$\overline{\Psi}(\lambda \epsilon_i)=T^{\hbigl(j+(i-1)(n+1)\hbigr)(\alpha-\alpha')}[\lambda ]_{\alpha,\alpha'}\epsilon_i,$$
for every $i\in \N$ and $\lambda\in\Lambda_{(0,\alpha]}$.
         \item If $\alpha'\geq\beta'$ then  $\alpha\geq \beta$ and the map
        $$\overline{\Psi} \fc \Lambda_{>0}\epsilon_0\oplus \bigoplus_{i\geq1}\Lambda_{(0,\beta]}\epsilon_i\to \Lambda_{>0}\epsilon'_0\oplus \bigoplus_{i\geq1}\Lambda_{(0,\beta']}\epsilon'_i$$
        satisfies
        $$\overline{\Psi}(\lambda \epsilon_0)=T^{j(\alpha-\alpha')}\lambda \epsilon'_0,$$
        for every $\lambda\in \Lambda_{>0}$, where $\epsilon_0$ is a generator of the $\Lambda_{>0}$ factor in $\coker D$, and
        $$\overline{\Psi}(\lambda \epsilon_i)=T^{\hbigl(j+i(n+1)\hbigr)(\alpha-\alpha')}[\lambda ]_{\beta,\beta'}\epsilon'_i,$$
for every $i\in \N$ and $\lambda\in\Lambda_{(0,\beta]}$.

        \item If $\alpha'<\beta'$ and $\alpha\geq \beta$ then the map
        $$\overline{\Psi} \fc \Lambda_{>0}\epsilon_0\oplus \bigoplus_{i\geq1}\Lambda_{(0,\beta]}\epsilon_i\to \bigoplus_{i\geq1}\Lambda_{(0,\alpha']}\epsilon'_i$$
        satisfies
        $$\overline{\Psi}(\lambda \epsilon_0)=T^{j(\alpha-\alpha')}\sum_{i=1}^\infty (-T^{\beta'-\alpha'})^{i-1}  [\lambda  ]_{\alpha'}\epsilon'_i\in\bigoplus_{i\geq1}\Lambda_{(0,\alpha']}\epsilon'_i,$$
        for every $\lambda\in \Lambda_{>0}$, and
        $$\overline{\Psi}(\lambda \epsilon_i)= T^{\hbigl(j+(i-1)(n+1)\hbigr)(\alpha-\alpha')+\alpha-\beta} [\lambda ]_{\beta,\alpha'}\epsilon'_i,$$
for every $i\in \N$ and $\lambda\in\Lambda_{(0,\beta]}$, where $\epsilon_i$ is a generator of the $i$-th summand in the direct sum $\bigoplus_{i\geq1}\Lambda_{(0,\beta]}\subset\coker D$.

    \end{enumerate}
\end{prop}

\begin{proof}\phantom{M}
    \begin{enumerate}
        \item Assume that $\alpha<\beta$. Then $\alpha'\leq \alpha<\beta\leq \beta'$, so $\alpha'<\beta'$. In this case, we use the Schauder bases $(\epsilon_i)_{i \in \N}$, $(\epsilon'_i)_{i \in \N}$ given by $\epsilon_i = e_i + T^{\beta-\alpha} e_{i+1}$ and $\epsilon_i = e_i + T^{\beta'-\alpha'} e_{i+1}$ for every $i \in \N$; see Example~\ref{exam: Schauder bases}. 

        By Proposition~\ref{prop: coker of D}, the quotient map $q\fc \widehat{\Lambda_{>0}^\infty}\to\coker D$ satisfies $q(\lambda \epsilon_i)=[\lambda]_\alpha\epsilon_i$, 
        for every $i\in \N$ and $\lambda\in \Lambda_{>0}$.
        Similarly, the quotient map $q'\fc \widehat{\Lambda_{>0}^\infty}\to\coker D'$ satisfies $q'(\lambda \epsilon'_i)=[\lambda]_{\alpha'}\epsilon'_i$, for every $i\in \N$ and $\lambda\in \Lambda_{>0}$.
        
Let $i\in \N$ and $\lambda\in\Lambda_{(0,\alpha]}$, let us compute $\overline{\Psi}(\lambda \epsilon_i)$. Denote $\eta=s_\alpha(\lambda)$, and note that $q(\eta\epsilon_i)=\lambda \epsilon_i$. Since $\overline{\Psi}\circ q=q'\circ\Psi_2$ we get that 
\begin{align*}
    \overline{\Psi}(\lambda \epsilon_i)&=\overline{\Psi}\circ q(\eta \epsilon_i)\\
    &=q'\circ\Psi_2(\eta \epsilon_i)\\
    &=q'\circ\Psi_2\left(\eta (e_i+T^{\beta-\alpha}e_{i+1})\right)\\
    &=q' (T^{\hbigl(j+(i-1)(n+1)\hbigr)(\alpha-\alpha')} \eta e_i + T^{\hbigl(j+i(n+1)\hbigr)(\alpha-\alpha')} T^{\beta-\alpha} \eta e_{i+1})\\
 &=T^{\hbigl(j+(i-1)(n+1)\hbigr)(\alpha-\alpha')} q' ( \eta e_i + T^{(n+1)(\alpha-\alpha')+\beta-\alpha} \eta e_{i+1}).
\end{align*}
Using the relation $n(\alpha-\alpha') = \beta'-\beta$ we get that 
$$(n+1)(\alpha-\alpha')+\beta-\alpha=\beta'-\beta+\alpha-\alpha'+\beta-\alpha=\beta'-\alpha',$$
therefore
\begin{align*}
 \overline{\Psi}(\lambda \epsilon_i)&=T^{\hbigl(j+(i-1)(n+1)\hbigr)(\alpha-\alpha')} q' ( \eta e_i + T^{\beta'-\alpha'} \eta e_{i+1})\\
 &=T^{\hbigl(j+(i-1)(n+1)\hbigr)(\alpha-\alpha')} q' ( \eta\epsilon'_i)\\
 &=T^{\hbigl(j+(i-1)(n+1)\hbigr)(\alpha-\alpha')}  [\eta]_{\alpha'}\epsilon'_i\\
 &=T^{\hbigl(j+(i-1)(n+1)\hbigr)(\alpha-\alpha')}  [\lambda]_{\alpha,\alpha'}\epsilon'_i,
\end{align*}
as required.

       \item Assume that $\alpha'\geq\beta'$. Then $\beta\leq \beta'\leq\alpha'\leq \alpha$, so $\alpha\geq\beta$. In this case, we use the Schauder bases $(\epsilon_i)_{i \in \N}$, $(\epsilon'_i)_{i \in \N}$ given by $\epsilon_1=\epsilon'_1=e_1$ and $\epsilon_i = T^{\alpha -\beta}e_{i-1} + e_{i}$ and $\epsilon'_i = T^{\alpha' -\beta'}e_{i-1} + e_{i}$ for every $i \geq2$; see Example~\ref{exam: Schauder bases}.

   By Proposition~\ref{prop: coker of D}, the quotient map $q\fc \widehat{\Lambda_{>0}^\infty}\to\coker D$ satisfies $q(\lambda \epsilon_1)=\lambda\epsilon_1$ for every $\lambda\in \Lambda_{>0}$ and $q(\lambda \epsilon_i)=[\lambda]_\beta\epsilon_i$, for every $i \geq2$ and $\lambda\in \Lambda_{>0}$. Similarly, the quotient map $q'\fc \widehat{\Lambda_{>0}^\infty}\to\coker D'$ satisfies $q'(\lambda \epsilon'_1)=\lambda\epsilon'_1$ for every $\lambda\in \Lambda_{>0}$ and $q'(\lambda \epsilon'_i)=[\lambda]_{\beta'}\epsilon'_i$, for every $i \geq2$ and $\lambda\in \Lambda_{>0}$.

Let $i\in \N$. If $i=1$, then for every $\lambda\in \Lambda_{>0}$ we have 

\begin{align*}
 \overline{\Psi}(\lambda \epsilon_1)&= \overline{\Psi}\circ q(\lambda \epsilon_1)\\
 &=q'\circ \Psi_2(\lambda \epsilon_1)\\
 &=q'\circ \Psi_2(\lambda e_1)\\
 &=q' (T^{j(\alpha-\alpha')} \lambda e_1)\\
 &=T^{j(\alpha-\alpha')} q' ( \lambda e_1)\\
  &=T^{j(\alpha-\alpha')} q' ( \lambda \epsilon'_1)\\
&=T^{j(\alpha-\alpha')} \lambda \epsilon'_1.
\end{align*}

Otherwise, $i\geq 2$. Let $\lambda\in \Lambda_{(0,\beta]}$,  denote $\eta=s_\beta(\lambda)$ and note that $q(\eta\epsilon_i)=\lambda\epsilon_i$. Thus we get
\begin{align*}
 \overline{\Psi}(\lambda \epsilon_i)&= \overline{\Psi}\circ q(\eta \epsilon_i)\\
 &=q'\circ \Psi_2(\eta \epsilon_i)\\
 &=q'\circ \Psi_2(\eta (T^{\alpha-\beta}e_{i-1} +  e_{i}))\\
 &=q' (T^{(j+(i-2)(n+1))(\alpha-\alpha')} T^{\alpha-\beta} \eta e_{i-1} + T^{\hbigl(j+(i-1)(n+1)\hbigr)(\alpha-\alpha')}  \eta e_{i}).
\end{align*}
Using the relation $n(\alpha-\alpha') = \beta'-\beta$ we get that 
$$\alpha-\beta=\alpha-\alpha'+\beta'-\beta+\alpha'-\beta'=(n+1)(\alpha-\alpha')+\alpha'-\beta'.$$

therefore
\begin{align*}
 \overline{\Psi}(\lambda \epsilon_i)&=q' (T^{(j+(i-2)(n+1))(\alpha-\alpha')} T^{\alpha-\beta} \eta e_{i-1} + T^{\hbigl(j+(i-1)(n+1)\hbigr)(\alpha-\alpha')}  \eta e_{i})\\
 &=q' (T^{(j+(i-2)(n+1))(\alpha-\alpha')+(n+1)(\alpha-\alpha')+\alpha'-\beta'} \eta e_{i-1} + T^{\hbigl(j+(i-1)(n+1)\hbigr)(\alpha-\alpha')}  \eta e_{i})\\
 &= T^{\hbigl(j+(i-1)(n+1)\hbigr)(\alpha-\alpha')} q' (\eta(T^{\alpha'-\beta'}  e_{i-1} +  e_{i}))\\
 &= T^{\hbigl(j+(i-1)(n+1)\hbigr)(\alpha-\alpha')} q' ( \eta \epsilon'_{i})\\
  &= T^{\hbigl(j+(i-1)(n+1)\hbigr)(\alpha-\alpha')}  [\eta]_{\beta'} \epsilon'_{i}\\
  &= T^{\hbigl(j+(i-1)(n+1)\hbigr)(\alpha-\alpha')}  [\lambda]_{\beta,\beta'} \epsilon'_{i}.
\end{align*}
as required.

       \item Assume that $\alpha'<\beta'$ and $\alpha\geq\beta$. In this case, we use the Schauder basis $(\epsilon_i)_{i \in \N}$ given by $\epsilon_1=e_1$ and $\epsilon_i = T^{\alpha -\beta}e_{i-1} + e_{i}$ every $i \geq2$, and the Schauder basis $\epsilon'_i =e_i + T^{\beta'-\alpha'} e_{i+1}$ for every $i\in \N$.

        By Proposition~\ref{prop: coker of D}, the quotient map $q\fc \widehat{\Lambda_{>0}^\infty}\to\coker D$ satisfies $q(\lambda \epsilon_1)=\lambda\epsilon_1$ for every $\lambda\in \Lambda_{>0}$ and $q(\lambda \epsilon_i)=[\lambda]_\beta\epsilon_i$, for every $i \geq2$ and $\lambda\in \Lambda_{>0}$. Similarly, the quotient map $q'\fc \widehat{\Lambda_{>0}^\infty}\to\coker D'$ satisfies $q(\lambda \epsilon'_i)=[\lambda]_{\alpha'}\epsilon'_i$, for every $i\in \N$ and $\lambda\in \Lambda_{>0}$.
       
Let $i\in \N$. If $i=1$, then for every $\lambda\in \Lambda_{>0}$ we have 

\begin{align*}
 \overline{\Psi}(\lambda \epsilon_1)&= \overline{\Psi}\circ q(\lambda \epsilon_1)\\
 &=q'\circ \Psi_2(\lambda \epsilon_1)\\
 &=q'\circ \Psi_2(\lambda e_1)\\
 &=q' (T^{j(\alpha-\alpha')} \lambda e_1)\\
 &=T^{j(\alpha-\alpha')} q' ( \lambda e_1)\\
&=T^{j(\alpha-\alpha')} q' \left(\lambda\sum_{k=1}^\infty(-T^{\beta'-\alpha'})^{k-1}\epsilon'_k\right) \\
&=T^{j(\alpha-\alpha')} \sum_{k=1}^\infty(-T^{\beta'-\alpha'})^{k-1} q' \left(\lambda\epsilon'_k\right) \\
&=T^{j(\alpha-\alpha')} \sum_{k=1}^\infty(-T^{\beta'-\alpha'})^{k-1} [\lambda]_{\alpha'}\epsilon'_k .
\end{align*}

Otherwise, $i\geq 2$. Let $\lambda\in \Lambda_{(0,\beta]}$,  denote $\eta=s_\beta(\lambda)$ and note that $q(\eta\epsilon_i)=\lambda\epsilon_i$.. Thus we get

       \begin{align*}
 \overline{\Psi}(\lambda \epsilon_i)&=\overline{\Psi}\circ q(\eta \epsilon_i)\\
 &=q'\circ \Psi_2(\eta \epsilon_i)\\
 &=q'\circ \Psi_2(\eta (T^{\alpha-\beta}e_{i-1} +  e_{i}))\\
 &=q' (T^{(j+(i-2)(n+1))(\alpha-\alpha')} T^{\alpha-\beta} \eta e_{i-1} + T^{\hbigl(j+(i-1)(n+1)\hbigr)(\alpha-\alpha')}  \eta e_{i})\\
 &=q' (T^{(j+(i-2)(n+1))(\alpha-\alpha')+\alpha-\beta} \eta e_{i-1} + T^{\hbigl(j+(i-1)(n+1)\hbigr)(\alpha-\alpha')}  \eta e_{i}).
\end{align*}
Using the relation $n(\alpha-\alpha') = \beta'-\beta$ we get that 
$$(n+1)(\alpha-\alpha')-(\alpha-\beta)=\alpha-\alpha'+\beta'-\beta-\alpha+\beta=\beta' -\alpha'>0.$$

therefore
\begin{align*}
 \overline{\Psi}(\eta \epsilon_i)&=q' (T^{(j+(i-2)(n+1))(\alpha-\alpha')+\alpha-\beta} \eta e_{i-1} + T^{\hbigl(j+(i-1)(n+1)\hbigr)(\alpha-\alpha')}  \eta e_{i})\\
 &= T^{(j+(i-2)(n+1))(\alpha-\alpha')+\alpha-\beta} q' ( \eta e_{i-1} + T^{(n+1)(\alpha-\alpha')-(\alpha-\beta)}  \eta e_{i})\\
 &= T^{(j+(i-2)(n+1))(\alpha-\alpha')+\alpha-\beta} q' ( \eta e_{i-1} + T^{\beta'-\alpha'}  \eta e_{i})\\
 &= T^{(j+(i-2)(n+1))(\alpha-\alpha')+\alpha-\beta} q' ( \eta \epsilon'_{i-1})\\
 &= T^{(j+(i-2)(n+1))(\alpha-\alpha')+\alpha-\beta} [\eta]_{\alpha'}\epsilon'_{i-1},\\
 &= T^{(j+(i-2)(n+1))(\alpha-\alpha')+\alpha-\beta} [\lambda]_{\beta,\alpha'}\epsilon'_{i-1}.
\end{align*}
as required.\qedhere

    \end{enumerate}
\end{proof}

Let $j \in \Z_{\geq0}$. Another computation, which is important to us for computing the restriction maps from $\CP^n$ to $B_\Delta$, is an explicit form of the map $\overline{\Psi}\fc \Lambda_{>0}\to \coker D$, which is induced from the map $\Psi\fc \Lambda_{>0}\to \widehat{\Lambda_{>0}^\infty}$ that satisfies
$$\Psi(\lambda e) = T^{j(1-\alpha)} \lambda e_1,$$
for every $\lambda\in \Lambda_{>0}$, where $e$ is the unit of $\Lambda_{\geq0}$ and $e_1$ is the first element in the standard basis of $\coker D$, as it is described in Proposition~\ref{prop: coker of D}, via the following commutative diagram:

\begin{equation}\label{diag: alg res maps simplified}
\vcenter{
  \xymatrix@R=2pc@C=3pc{
    0 \ar[r] \ar[d] & \Lambda_{>0} \ar[d]_{\Psi} \ar[r]^{\id} & \Lambda_{>0} \ar[d]_{\overline{\Psi}} \\
    \widehat{\Lambda_{>0}^\infty} \ar[r]^{D} & \widehat{\Lambda_{>0}^\infty} \ar[r]^{q} & \coker D
  }}
\end{equation}
where $q\fc \widehat{\Lambda_{>0}^\infty}\to \coker D$ is the quotient map.

\begin{prop}\label{prop: alg restriction maps simplified}
    The map $\overline{\Psi} \fc \Lambda_{>0} \to \coker D$ can be explicitly described as follows:
    \begin{enumerate}
       
        \item If $\alpha<\beta$, then the map
        $$\overline{\Psi} \fc \Lambda_{>0}e\to \bigoplus_{i\geq1}\Lambda_{(0,\alpha]}\epsilon_i$$
        satisfies, 
        $$\overline{\Psi}(\lambda e)=T^{j(1-\alpha)} [\lambda]_\alpha \cdot\sum_{i=1}^\infty(-T^{\beta-\alpha})^{i-1}\epsilon_i,$$
        for every $\lambda\in \Lambda_{>0}$.
        
         \item If $\alpha\geq\beta$, then the map
        $$\overline{\Psi} \fc \Lambda_{>0}e\to \Lambda_{>0}\epsilon_1\oplus \bigoplus_{i\geq2}\Lambda_{(0,\beta]}\epsilon_i$$
        satisfies,
        $$\overline{\Psi}(\lambda e)=T^{j(1-\alpha)} \lambda\epsilon_1,$$
         for every $\lambda\in \Lambda_{>0}$.
    \end{enumerate}
\end{prop}
\begin{proof}
    The commutativity of the diagram in Equation~\eqref{diag: alg res maps simplified} dictates that $\overline{\Psi} \circ \id = q \circ \Psi$. Therefore, for every $\lambda \in \Lambda_{>0}$, we have:
    $$\overline{\Psi}(\lambda e) = q(\Psi(\lambda e)) = q(T^{j(1-\alpha)} \lambda e_1) = T^{j(1-\alpha)} q(\lambda e_1).$$

    \begin{enumerate}
        \item Assume that $\alpha<\beta$. In this case, we use the Schauder basis $(\epsilon_i)_{i \in \N}$ given by $\epsilon_i = e_i + T^{\beta-\alpha} e_{i+1}$ for every $i \in \N$; see Example~\ref{exam: Schauder bases}. By Proposition~\ref{prop: coker of D}, the quotient map $q\fc \widehat{\Lambda_{>0}^\infty}\to\coker D$ satisfies $q(\lambda \epsilon_i)=[\lambda]_\alpha\epsilon_i$, 
        for every $i\in \N$ and $\lambda\in \Lambda_{>0}$. Therefore
        \begin{align*}
            \overline{\Psi}(\lambda e) &= T^{j(1-\alpha)} q(\lambda e_1)\\
            &= T^{j(1-\alpha)} q\left(\lambda \sum_{i=1}^\infty(-T^{\beta-\alpha})^{i-1}\epsilon_i\right)\\
            &=T^{j(1-\alpha)} [\lambda]_\alpha \cdot\sum_{i=1}^\infty(-T^{\beta-\alpha})^{i-1}\epsilon_i,
        \end{align*}
        for every $\lambda\in\Lambda_{>0}$.
        
        \item Assume $\alpha\geq\beta$. Let us use the Schauder basis $(\epsilon_i)_{i \in \N}$ given by $\epsilon_1=e_1$ and $\epsilon_i = T^{\alpha -\beta}e_{i-1} + e_{i}$ every $i \geq2$.

        By Proposition~\ref{prop: coker of D}, the quotient map $q\fc \widehat{\Lambda_{>0}^\infty}\to\coker D$ satisfies $q(\lambda \epsilon_1)=\lambda\epsilon_1$ for every $\lambda\in \Lambda_{>0}$.  Therefore, we get that
        $$\overline{\Psi}(\lambda e) = T^{j(1-\alpha)} q(\lambda e_1)= T^{j(1-\alpha)} q(\lambda \epsilon_1) =T^{j(1-\alpha)} \lambda\epsilon_1,$$
        for every $\lambda\in\Lambda_{>0}$, as required.
    \end{enumerate}
\end{proof}

\subsection{Limits of 1-rays}\label{ss: lim of 1-ray}

In this section, we discuss the notion of the telescope of a $1$-ray, its connection to the direct limit, and their completions. The main result of this section is Proposition~\ref{prop: algebraic preparation for computing restriction maps}, which will be crucial for the proof of Theorem~\ref{thm: res for balls in CP^n}.

\begin{defin}\label{def:1-ray of clx}
	A \textbf{$1$-ray of cochain complexes} $\cC$ is an infinite diagram
	$$
	\cC\fc
	C_1 \overset{f_1}{\longrightarrow}
	C_2 \overset{f_2}{\longrightarrow}
	C_3 \longrightarrow \cdots,
	$$
	where $(C_i)_{i \geq 1}$ are cochain complexes over $\Lambda_{\geq 0}$, and for every
	$i \geq 1$ the map $f_i \fc C_i \to C_{i+1}$ is a cochain map.
\end{defin}

	Recall that given modules $C_i$, $i \in \N$, and module homomorphisms
	$f_i \fc C_i \to C_{i+1}$, the direct limit $\drctlim_i C_i$ can be realized
	as the cokernel of the map
	$$
	\id - f \fc \bigoplus_i C_i \to \bigoplus_i C_i,
	$$
	where $f \fc \bigoplus_i C_i \to \bigoplus_i C_i$ is defined by
	$$
	f(x_1,x_2,\ldots) = (0,f_1(x_1),f_2(x_2),\ldots)
	$$
	for every $(x_1,x_2,\ldots) \in \bigoplus_i C_i$.

	In particular, every $1$-ray of cochain complexes $\cC$ admits a direct limit in the category of cochain complexes, which we denote by $\drctlim \cC$.

Also, recall that given two cochain complexes $(A,d_A)$ and $(B,d_B)$ and a cochain map
$\psi\fc A \to B$, the cone of $\psi$, denoted $\cone(\psi)$, is the cochain complex defined by
$$
\cone(\psi)
=
\left(
B \oplus A[1],
\begin{pmatrix}
d_B & \psi \\
0 & -d_A
\end{pmatrix}
\right).
$$
One of the most important properties of the mapping cone is that $f$ is a quasi-isomorphism if and only if $\cone(f)$ is acyclic, i.e., $H^*(\cone(f))=0$ (see \cite[Corollary 10.41]{Rotman_intro_Homo_alg} for instance).

In homotopical algebra, the cone of a morphism is regarded as a homotopical model for its cokernel. The viewpoint of the cone as a homotopical model for the cokernel leads to the definition of a homotopical model for the direct limit, called the \emph{telescope}.

\begin{defin}
	The \textbf{telescope} of a $1$-ray of cochain complexes
	$$
	\cC\fc
	C_1 \overset{f_1}{\longrightarrow}
	C_2 \overset{f_2}{\longrightarrow}
	C_3 \longrightarrow \cdots
	$$
	is denoted by $\tel \cC$ and is defined as the cone of the map
	$$
	\id - f \fc \bigoplus_i C_i \to \bigoplus_i C_i,
	$$
	where $f \fc \bigoplus_i C_i \to \bigoplus_i C_i$ is given by
	$$
	f(x_1,x_2,\ldots) = (0,f_1(x_1),f_2(x_2),\ldots)
	$$
	for every $(x_1,x_2,\ldots) \in \bigoplus_i C_i$.
\end{defin}

Let us present two important results from \cite{Varolgunes_2021_MV_and_relSH} regarding telescopes and directs limits.

\begin{lemma}[\cite{Varolgunes_2021_MV_and_relSH}]\label{lemma: tel->colim}
Let $\cC$ an $1$-ray.
\begin{itemize}
    \item (Lemma 2.2.1) There is a quasi-isomorphism 
    \begin{equation}\label{eq: tel-drctlim}
        \tel(\cC)\to \drctlim \cC.
    \end{equation}
    \item (Lemma 2.3.7) The completion of the map from Equation~\eqref{eq: tel-drctlim} 
    $$\widehat{\tel}\,(\cC)\to\widehat{\drctlim}\cC,$$
    is a quasi-isomorphism.
\end{itemize}
\end{lemma}

Next, we discuss maps between $1$-rays, which are crucial for the definition of restriction maps in relative symplectic cohomology. A map $G\fc \cC\to \cC'$ between two $1$-rays 
$$\cC\fc C_1\overset{f_1}{\longrightarrow} C_2\overset{f_3}{\longrightarrow} C_3\longrightarrow\cdots,$$
$$\cC'\fc C_1'\overset{f_1'}{\longrightarrow} C_2'\overset{f_3'}{\longrightarrow} C_3'\longrightarrow\cdots,$$
is a pair of sequences $(g_i\fc C_i\to C_i')_{i\in \N}$, $(h_i\fc C_i\to C_{i+1}')_{i\in \N}$, such that for every $i\in\N$ we have $$g_{i+1}f_i-f_i'g_i=d_{C_{i+1}'}h_i+h_i d_{C_i}.$$ 
We will present such a map by the following diagram
\begin{equation}\label{diagram: 2-ray}
\vcenter{
\xymatrix@R=3pc@C=3pc{
	C_1 \ar[r]^{f_1} \ar[d]^{g_1} \ar@{=>}[rd]^{h_1} & C_2 \ar[r]^{f_2} \ar[d]^{g_2} \ar@{=>}[rd]^{h_2} & C_3 \ar[r] \ar[d]^{g_3} \ar@{=>}[rd] & \cdots \\
	C_1' \ar[r]_{f_1'} & C_2' \ar[r]_{f_2'} & C_3' \ar[r] & \cdots
}}
\end{equation}
 A map $G\fc \cC\to \cC'$ induces a cochain map, $\tel G\fc \tel(\cC)\to\tel(\cC')$, that is given by
 $$\tel G\left(\left(x_i,y_i\right)_{i=1}^
 \infty\right)=\left((g_1(x_1),g_1(y_1)),\left(g_i(x_i)-h_{i-1}(y_{i-1}),g_i(y_i)\right)_{i=2}^\infty\right),$$
as depicted in the following diagram.

$$
\xymatrix@R=2pc@C=3pc{
	C_1 \ar@(ur,ul)[] \ar[ddd]^-{g_1} & & C_2 \ar@(ur,ul)[] \ar[ddd]^-{g_2} & & C_3 \ar@(ur,ul)[] \ar[ddd]^-{g_3} & & \cdots \\
	& C_1[1] \ar[lu]_-{\id} \ar[ru]^-{-f_1} \ar@(dr,dl)[] \ar[ddd]^-{g_1} \ar[rdd]^-{-h_1} & & C_2[1] \ar[lu]_-{\id} \ar[ru]^-{-f_2} \ar@(dr,dl)[] \ar[ddd]^-{g_2} \ar[rdd]^-{-h_2} & & C_3[1] \ar[lu]_-{\id} \ar@{-->}[ru] \ar@(dr,dl)[] \ar[ddd]^-{g_3} \ar@{-->}[rdd] & \\
	& & & & & & \\
	C_1' \ar@(ur,ul)[] & & C_2' \ar@(ur,ul)[] & & C_3' \ar@(ur,ul)[] & & \cdots \\
	& C_1'[1] \ar[lu]_-{\id} \ar[ru]^-{-f_1'} \ar@(dr,dl)[] & & C_2'[1] \ar[lu]_-{\id} \ar[ru]^-{-f_2'} \ar@(dr,dl)[] & & C_3'[1] \ar[lu]_-{\id} \ar@{-->}[ru] \ar@(dr,dl)[] &
}
$$

Another important result from \cite{Varolgunes_2021_MV_and_relSH} is summarized below.

\begin{prop}[{\cite[Lemma 2.2.2]{Varolgunes_2021_MV_and_relSH}}]\label{prop: 2-ray, tel, colim}
Let $\cC, \cC'$ be $1$-rays of cochain complexes, and let $G \fc \cC \to \cC'$ be a map of $1$-rays. 
\begin{itemize}
    \item The map $G$ induces a commutative diagram:
    \[
    \xymatrix{
    H(C_1) \ar[r] \ar[d] & H(C_2) \ar[r] \ar[d] & H(C_3) \ar[r] \ar[d] & \cdots \\
    H(C'_1) \ar[r] & H(C'_2) \ar[r] & H(C'_3) \ar[r] & \cdots
    }
    \]
    \item $G$ induces a canonical morphism $\drctlim H(C_n) \to \drctlim H(C'_n)$.
    \item The following diagram commutes:
    $$
    \xymatrix@R=2pc@C=2pc{
        H(\tel\,\cC) \ar[r] \ar[d] & H(\tel\,\cC') \ar[d] \\
        \drctlim H(C_n) \ar[r] & \drctlim H(C'_n)
    }$$
    where the vertical arrows are isomorphisms, the top arrow is induced by $\tel(G)$, and the bottom arrow is induced by the morphism from the previous item.
\end{itemize}
\end{prop}

The main result of this section is the following.

\begin{prop}\label{prop: algebraic preparation for computing restriction maps}
Let 
$$ \cC \fc C_1 \xrightarrow{f_1} C_2 \xrightarrow{f_2} C_3 \to \cdots $$
$$ \cC' \fc C'_1 \xrightarrow{f'_1} C'_2 \xrightarrow{f'_2} C'_3 \to \cdots $$
be two $1$-rays, and let $G \fc \cC \to \cC'$ be a morphism given by the sequences $(g_i \fc C_i \to C'_i)_{i \in \N}$ and $(h_i \fc C_i \to C'_{i+1})_{i \in \N}$. Assume that for every $i \in \N$ we have $g_{i+1} f_i - f'_i g_i = 0$, that is, the $1$-skeleton of the diagram in \eqref{diagram: 2-ray} commutes. Then:

\begin{enumerate}
    \item There exists an induced morphism between the direct limits:
    $$ \drctlim \cC \to \drctlim \cC'.$$

    \item The following diagram commutes:
    \begin{equation}\label{diagram: 2D, lim, tel}
    \vcenter{
    \xymatrix@R=2pc@C=4pc{
        H(\tel\, \cC) \ar[r] \ar[d] & H(\tel\, \cC') \ar[d]\\ 
        H(\drctlim \cC) \ar[r] & H(\drctlim \cC')
    }}
    \end{equation}
    where the vertical arrows are isomorphisms.

    \item If the morphism $\drctlim \cC \to \widehat{\drctlim} \cC$ induces a surjection in homology, then the diagram
    \begin{equation}\label{diagram: 2D, lim, tel and completion}
    \vcenter{
    \xymatrix@R=2pc@C=4pc{
        H(\widehat{\tel}\, \cC) \ar[r] \ar[d] & H(\widehat{\tel}\, \cC') \ar[d]\\ 
        H(\widehat{\drctlim} \cC) \ar[r] & H(\widehat{\drctlim} \cC')
    }}
    \end{equation}
    commutes and its vertical arrows are isomorphisms.
\end{enumerate}
\end{prop}

 \begin{proof}\phantom{M}
\begin{enumerate}
    \item The assumption that $g_{i+1}f_i - f'_i g_i = 0$ for every $i \in \N$, implies that the collection $\{g_i\}$ defines a morphism of direct systems between $\langle C_i \rangle_{i \in \N}$ and $\langle C'_i \rangle_{i \in \N}$. Consequently, there is an induced map between the direct limits (see, e.g., \cite[pp. 246--247]{Rotman_intro_Homo_alg}).
        
    \item Consider the following diagram:
    \begin{equation}\label{diagram: proof_step_2}
    \vcenter{
    \xymatrix@R=2pc@C=4pc{
        H(\tel \cC) \ar[r] \ar[d] & H(\tel \cC') \ar[d] \\
        \drctlim H(C_n) \ar[r] \ar[d] & \drctlim H(C'_n) \ar[d] \\
        H(\drctlim C_n) \ar[r] & H(\drctlim C'_n)
    }}
    \end{equation}
    where we identify $\drctlim C_n = \drctlim \cC$. The upper square commutes by Proposition~\ref{prop: 2-ray, tel, colim}. The bottom square commutes because the homology functor commutes with direct limits \cite[Page 58]{Weibel_1994_hom_alg}. Since both squares commute, the outer rectangle commutes, establishing the result for Equation~\eqref{diagram: 2D, lim, tel}.

    \item Consider the cube diagram:
    $$ \xymatrix@R=2pc@C=2pc{
        H(\tel \cC) \ar[rr]^F \ar[rd]^{c_1} \ar[dd]_{g_1} & & H(\tel \cC') \ar[rd]^{c_2} \ar[dd]^(.3){g_2} & \\
        & H(\widehat{\tel}\,\cC) \ar[rr]^(.4){\hat{F}} \ar[dd]^(.3){\hat{g}_1} & & H(\widehat{\tel}\,\cC') \ar[dd]^{\hat{g}_2} \\
        H(\drctlim \cC) \ar[rr]^(.7)f \ar[rd]^{c_3} & & H(\drctlim \cC') \ar[rd]^{c_4} & \\
        & H(\widehat{\drctlim} \cC) \ar[rr]_{\hat{f}} & & H(\widehat{\drctlim} \cC')
    } $$
    Here, the back square commutes by the previous item. The maps $c_1, c_2,c_3, c_4$ are induced by the naturality of the completion maps. The vertical arrows $g_1, g_2, \hat{g}_1, \hat{g}_2$ are the isomorphisms in homology induced by Lemma~\ref{lemma: tel->colim}. The remaining side, top, and bottom squares commute by the properties of completion (Remark~\ref{rem: properties of completion}). 

    We aim to show the front square commutes, i.e., $\hat{g}_2 \circ \hat{F} = \hat{f} \circ \hat{g}_1$. Let $x \in H(\widehat{\tel}\,\cC)$. By assumption, the map $c_3$ is surjective. Thus, there exists $y \in H(\drctlim\cC)$ such that $\hat{g}_1(x) = c_3(y)$. Since the map $g_1$ is invertible, we may denote $z = g_1^{-1}(y)$. The left square commutes, so $\hat{g}_1 \circ c_1 = c_3 \circ g_1$. Because $\hat{g}_1$ is invertible, we deduce that:
    $$c_1(z) = \hat{g}_1^{-1} \circ c_3 \circ g_1(z) = \hat{g}_1^{-1} \circ c_3(y) = \hat{g}_1^{-1} \circ \hat{g}_1(x) = x.$$

    Therefore:
    $$\hat{g}_2 \circ \hat{F}(x) = \hat{g}_2 \circ \hat{F} \circ c_1(z).$$
    Since the top square commutes, we have $\hat{F} \circ c_1 = c_2 \circ F$, which implies:
    $$\hat{g}_2 \circ \hat{F}(x) = \hat{g}_2 \circ c_2 \circ F(z).$$
    The right square commutes, meaning $\hat{g}_2 \circ c_2 = c_4 \circ g_2$; thus:
    $$\hat{g}_2 \circ \hat{F}(x) = c_4 \circ g_2 \circ F(z).$$
    The back square commutes as well, i.e., $g_2 \circ F = f \circ g_1$. Therefore:
    $$\hat{g}_2 \circ \hat{F}(x) = c_4 \circ f \circ g_1(z).$$
    Finally, the bottom square commutes, namely $c_4 \circ f = \hat{f} \circ c_3$, hence:
    $$\hat{g}_2 \circ \hat{F}(x) = \hat{f} \circ c_3 \circ g_1(z).$$
    By definition, $z = g_1^{-1}(y)$ and $c_3(y) = \hat{g}_1(x)$. We conclude that:
    $$\hat{g}_2 \circ \hat{F}(x) = \hat{f} \circ c_3 \circ g_1(g_1^{-1}(y)) = \hat{f} \circ c_3(y) = \hat{f} \circ \hat{g}_1(x),$$
    as required. \qedhere
\end{enumerate}
\end{proof}

\begin{rem}
    If, in addition to the assumptions of Proposition~\ref{prop: algebraic preparation for computing restriction maps}, we assume that the chain homotopy maps $(h_i)$ are zero, then the diagram
    \begin{equation*}
    \xymatrix@R=2pc@C=3pc{
        \tel\, \cC \ar[r] \ar[d] & \tel\, \cC' \ar[d]\\ 
        \drctlim \cC \ar[r] & \drctlim \cC'
    }
    \end{equation*}
    commutes at the level of cochain complexes. Consequently, one can deduce the same conclusions as in Proposition~\ref{prop: algebraic preparation for computing restriction maps}.
\end{rem}

    \section{Relative symplectic cohomology}\label{s: Def_relSH}

    This section is dedicated to relative symplectic cohomology. We survey its definition, an analogous definition based on the theory of Floer--Morse--Bott with cascades, and finally, we explain the equivalence between these definitions.
    
\subsection{Definition of the relative symplectic cohomology}\label{ss: def of relSH}
This section is dedicated to the definition of relative symplectic cohomology and its associated restriction maps; for further details, see also \cite{Varolgunes_2021_MV_and_relSH} and \cite[Section 4]{DGPZ_2024_Symp_top_and_IVMs}. Throughout this section we fix a closed symplectic manifold $(M,\omega)$.

\subsubsection{Weighted Floer complexes}\label{sss: weighted_CF}

We start by specifying the Floer complexes we will be using. If $H \in C^\infty(M\times S^1)$ is a non-degenerate Hamiltonian, its Floer complex, generated by the set $\cP^\circ(H)$ of its contractible $1$-periodic orbits, is defined as 
$$CF^*(H) = \bigoplus_{x\in\cP^\circ(H)}\Lambda_{\geq0}\cdot x.$$ 
This complex is graded over $\Z/2N_M^{}\Z$ by the Conley--Zehnder index shifted by $n$, where $N_M^{}$ is the minimal first Chern number. The differential is determined by its matrix elements, given by
$$\langle dx,y\rangle = \sum_{A \in \pi_2(M,x,y)}\#\wh\cM(H;x,y;A)\,T^{E_{top}(A)}\,,$$
where $\pi_2(M,x,y)$ is the set of homotopy classes of continuous maps of the cylinder $\R\times S^1$ to $M$ that are asymptotic at $-\infty$ to $x$ and at $+\infty$ to $y$, $\wh\cM(H;x,y;A)$ stands for the moduli space of unparametrized Floer flowlines with cascades corresponding to $H$, going from $x$ to $y$, and representing the class $A$, such that every Floer solution has index $1$; $\#\wh\cM(H;x,y;A)$ is a suitable virtual count as in Pardon \cite{Pardon_2016_alg_vir_counting,Pardon_2019_contact} or just a signed count if $M$ is assumed to be semipositive\footnote{Bai and Xu \cite{Bai_Xu_2026_integral_Floer} recently introduced an integral Hamiltonian Floer theory on general symplectic manifolds. It is highly likely that this construction can serve as another Floer-theoretic framework for relative symplectic cohomology, analogous to the use of Pardon's framework in Varolg\"une\c s' original construction.}. Finally $E_{top}(A) = \int_{S^1}\big(H_t(y(t))-H_t(x(t))\big)\,dt+\langle \omega,A\rangle$ is the topological energy of solutions in class $A$ (see Section~\ref{s: MB in Floer}). Continuation maps between such complexes, corresponding to monotone nondecreasing homotopies of Hamiltonians, are defined in a similar manner, with weights given by powers of $T$ whose exponents are the topological energies.

\subsubsection{Relative symplectic cohomology}\label{sss: relSH}

We start with the following definitions.

\begin{defin}\label{def: acc. data}\phantom{M}
\begin{itemize}
  \item An \textbf{acceleration datum} is a pointwise increasing sequence $\cH = (H_i)_{i=1}^\infty \subset C^\infty(S^1 \times M)$ of nondegenerate Hamiltonians, equipped with monotone homotopies between each successive Hamiltonian in the sequence.
  \item Given two acceleration data $\cH = (H_i)_i$ and $\cH' = (H_i')_i$, we write $\cH \preceq \cH'$ if $H_i \leq H_i'$ for all $i$.
  \item The \textbf{Floer $1$-ray} of an acceleration datum $\cH = (H_i)_i$, denoted by $CF(\cH)$, is the $1$-ray of cochain complexes:
  $$CF^*(H_1) \xrightarrow{\phi_1} CF^*(H_2) \xrightarrow{\phi_2} CF^*(H_3) \to \dots$$
  where the connecting morphisms $(\phi_i)_i$ are the continuation maps determined by the homotopies.
\end{itemize}
\end{defin}

\begin{defin}\label{def: acc. data for set}
  Let $K \subset M$ be compact. We say that an acceleration datum $\cH = (H_i)_i$ is an \textbf{acceleration datum for $K$} if for every $(t,x)\in S^1\times M$ we have
  $$\lim_{i\to\infty}H_i(t,x)=\left\{\begin{array}{cc}
      0, & x\in K, \\
     +\infty,  & x\notin K.
  \end{array}\right.$$
\end{defin}

\begin{rem}
    As a consequence of Dini's theorem, an acceleration datum is an acceleration datum for $K$ if and only if it is a cofinal sequence in $$C^\infty_{K \subset M} = \{H \in C^\infty( S^1\times M) \,|\, H|_{S^1\times K} < 0\}$$ relative to the usual order on functions.\footnote{This is a variant of the notation used in \cite{Varolgunes_2021_MV_and_relSH}.} See \cite[Lemma 3.3.1]{Varolgunes_2021_MV_and_relSH}.
\end{rem}

\begin{defin}(Varolg\"une\c s \cite[Section 3.3.2]{Varolgunes_2021_MV_and_relSH})
  Let $K \subset M$ be compact and let $\cH = (H_i)_i$ be an acceleration datum for $K$. Define the following complex:
  $$SC^*(\cH) := \wh\tel \, CF(\cH).$$
  The \textbf{relative symplectic cohomology of $K$ inside $M$} is defined as
  $$SH_M^*(K) := H\big(SC^*(\cH)\big).$$
\end{defin}

Given two acceleration data $\cH$ and $\cH'$ satisfying $\cH \preceq \cH'$, Varolg\"une\c s constructs in \cite{Varolgunes_2021_MV_and_relSH}  a map $G \fc CF(\cH) \to CF(\cH')$ between the corresponding Floer $1$-rays, using a Floer-theoretic construction together with properties of the telescope.
We call the maps constructed in this way \textbf{Floer-theoretic maps}.

\begin{defin}(Varolg\"une\c s \cite[Section 3.3.2]{Varolgunes_2021_MV_and_relSH})
  Let $K' \subset K \subset M$ be compact subsets, let $\cH$ and $\cH'$ be acceleration data for $K$ and $K'$, respectively, such that $\cH \preceq \cH'$, and let $G \fc CF(\cH) \to CF(\cH')$ be any Floer-theoretic map. The \textbf{restriction map} $\res^K_{K'} \fc SH_M^*(K) \to SH_M^*(K')$ is the map induced on homology by the map $\wh\tel\, G \fc SC^*(\cH) \to SC^*(\cH')$.
\end{defin}

\begin{rem}
  Varolg\"une\c s proves in \cite[Proposition 3.3.2]{Varolgunes_2021_MV_and_relSH} that both the relative symplectic cohomology and its associated restriction maps are well-defined. Specifically, any two acceleration data associated with a compact set define the same cohomology up to a canonical isomorphism, and similarly for the restriction maps.
\end{rem}

A practical way to compute relative symplectic cohomology in concrete examples is to use the direct limit instead of a telescope. Indeed, given an acceleration datum $\cH$, Lemma~\ref{lemma: tel->colim} asserts there exists a quasi-isomorphism:
$$SC^*(\cH) = \wh\tel\,CF(\cH) \to \wh\drctlim CF(\cH).$$
Thus, if $\cH$ is an acceleration datum for a compact subset $K \subset M$, then
$$SH^*_M(K) = H\left(\wh\drctlim CF(\cH)\right).$$
Thus, we base our proof of Theorem~\ref{thm: SH of ball in CP^n} on the computation of a direct limit of Floer complexes. The computation of the restriction maps is more delicate and is discussed in detail in Section~\ref{ss: res maps}.

    \subsection{Relative SH in Floer--Morse--Bott setting}\label{ss: relSH using MB}
    
The proof of Theorem~\ref{thm: SH of ball in CP^n} is based on a computation of the relative symplectic cohomology in the setting of Morse--Bott with cascades (the Floer--Morse--Bott setting). Let us briefly introduce this setting here. We also explain why the resulting relative symplectic cohomology coincides with the original version surveyed in the previous section. Let us start with the relevant modification of the Floer complex.

        \subsubsection{Weighted Floer--Morse--Bott complexes}\label{sss: weighted_CF_MB}
        
    As in Section~\ref{sss: weighted_CF}, we define the weighted Floer complexes, but this time, in the Floer--Morse--Bott setting.

    Throughout this subsection, we assume that the Hamiltonians defining the Floer complexes are constant along each of their connected critical submanifolds. This assumption applies to the Hamiltonians at the vertices of the cubes and to every Hamiltonian in an acceleration datum.
    
        Let $H \in C^\infty(S^1\times M)$ be a Hamiltonian satisfying the \textbf{MB} condition, let $h$ be a Morse function on $\bS_H$, let $g$ be a Riemannian metric on  $\bS_H$, and let $J$ be a compatible almost complex structure on $M$.
        
        Assuming regularity for all the moduli spaces involved, define the Floer complex associated to the pair $(H,J)$ as the free $\Lambda_{\geq0}$-module generated by the set $\Crit(h^{}_H)$, namely
        \[CF^*(H) = \bigoplus_{x\in\Crit(h^{}_H)}\Lambda_{\geq0}\cdot x.\]
        This complex is graded over $\Z/2N_M^{}\Z$ by the Floer--Morse--Bott index $\mu_{FMB}$, where $N_M$ is the minimal first Chern number of $M$, see Definition~\ref{def: FMB index} and  Remark~\ref{rem: mu_FMB well def mod 2N_M}. 
 The differential is determined by its matrix elements, given by
$$\langle dx,y\rangle = \sum_{u \in \widehat{\cM}^0(x,y;H,J)}\sgn(u)\cdot\,T^{E_{top}(u)}\,,$$
  where $\widehat{\cM}^0(x,y;H,J)$ is the union of the zero dimensional components of $\widehat{\cM}(x,y;H,J)$, and $\sgn(u)$ is a sign determined by a system of coherent orientations chosen on the moduli spaces. The exponent $E_{top}(u)$ is the topological energy of the Floer flowline with cascades $u$, see Equation~\eqref{eq: top energy of flowline} for the definition.   Continuation maps are defined similarly.

\subsubsection{Relative symplectic cohomology}\label{sss: relSH with MB}
Given a monotone Hamiltonian cube $\cH$ satisfying the \textbf{MB} condition and a cube of compatible almost complex structures $\cJ$, if the pair $(\cH,\cJ)$ is regular, one can apply the construction from Section~\ref{sss:Ham cubes FMB} to obtain a Floer--Morse--Bott-theoretic cube, using weighted Floer complexes, see Section~\ref{sss: weighted_CF_MB}. The vertices of this cube are the Floer complexes of the vertices of $(\cH,\cJ)$, and the face maps are defined by counting homotopy flowlines with cascades of the corresponding local dimension. Following Definitions~\ref{def: acc. data} and~\ref{def: acc. data for set}, acceleration data---both in general and for a compact set---must satisfy the \textbf{MB} condition here, and Floer $1$-rays are defined by considering weighted Floer--Morse--Bott complexes.

Similarly to Varolg\"une\c s \cite[Section 3.3.2]{Varolgunes_2021_MV_and_relSH}, given a compact subset $K \subset M$ and an acceleration datum $\cH = (H_i)_i$ for $K$ that satisfies the \textbf{MB} condition, the corresponding Floer--Morse--Bott complex is defined as
$$SC^*_{FMB}(\cH) := \wh\tel \, CF(\cH).$$
Thus, the \textbf{Floer--Morse--Bott relative symplectic cohomology of $K$ inside $M$} is
$$SH_{M,FMB}^*(K) := H\big(SC_{FMB}^*(\cH)\big).$$

First, Varolg\"une\c s's filling result for partial monotone Hamiltonian cubes, \cite[Proposition 3.2.18]{Varolgunes_2021_MV_and_relSH} (see also \cite[p. 558, 594]{Varolgunes_2021_MV_and_relSH}), directly translates to the Floer--Morse--Bott setting. Consequently, since Floer--Morse--Bott theory with cascades extends classical Floer theory, given two acceleration data $\cH, \cH'$ for a compact subset $K$ such that $\cH'$ consists only of non-degenerate Hamiltonians, the proof of \cite[Proposition 3.3.3]{Varolgunes_2021_MV_and_relSH} by Varolg\"une\c s is entirely adaptable. Thus, we conclude that $SH_{M,FMB}^*(K)$ is canonically isomorphic to the original relative symplectic cohomology $SH_{M}^*(K)$.

Let us remark that restriction maps are defined similarly and coincide with those from the original definition.

\section{Computations of relative symplectic cohomology}\label{s: computations_of_relSH}

 In this section, we apply our computations for Floer complexes from Section~\ref{s: computations of CF} to compute the relative symplectic cohomology of toric balls in $\CP^n$ and the restriction maps between them, thereby proving Theorem~\ref{thm: SH of ball in CP^n} and Theorem~\ref{thm: res for balls in CP^n}.

  \subsection{Energy of Floer trajectories}\label{ss: energy of our cascades}

Let $\Delta \in (0,1)$. Let us describe the complexes and connecting maps of the $1$-ray
$$CF(H_0) \to CF(H_1) \to CF(H_2) \to \cdots$$
where for every $\ell \in \Z_{\geq 0}$, the Hamiltonian $H_\ell \fc \CP^n \to \R$ is given by $H_\ell(z) = h(\Delta, \ell, \mu(z))$ for every $z \in \CP^n$.

Theorem~\ref{thm: CF(H_l;Z)} and Theorem~\ref{thm: continuation maps, over Z} describe the Floer complexes of the Hamiltonians $(H_i)_i$ over $\Z$. To describe them with coefficients in the Novikov ring, in the sense of Section~\ref{sss: weighted_CF}, we must compute the energy of the Floer and continuation flowlines with cascades that appear in the differential and continuation maps.

Fix $\ell \in \Z_{\geq 0}$. For every $1 \leq i \leq \ell$, let $r_{\ell,i}^\Delta \in [0,1)$ denote the number satisfying $\frac{\partial}{\partial r}|_{r=r_{\ell,i}^\Delta} h(\Delta, \ell, r) = i$. For every $1\leq i\leq \ell$ let us denote $h_{\ell,i}^\Delta=h(\Delta,\ell,r_{\ell,i}^\Delta)$, and set $h_{\ell,0}^\Delta=h(\Delta,\ell,0).$ In addition, for every $1 \leq i \leq \ell$, define the $(2n-1)$-spheres
$$S_{\ell,i} = \{z \in \CP^n : h'_\ell(\mu(z)) = i\},$$
each equipped with a perfect Morse function. Furthermore, we equip the divisor $D_\infty = \CP^{n-1}$ with a perfect Morse function. 

As before, let $\check{x}_0^\ell$ denote the constant orbit of $H_\ell$ at $0 \in \Int B(1) \subset \CP^n$. For every $1 \leq i \leq \ell$, let $\check{x}^\ell_i$ and $\hat{x}^\ell_i$ denote the minimum and maximum on $S_{\ell,i}$, respectively. Finally, for every $\ell+1 \leq i \leq \ell + n$, let $\check{x}_i^\ell$ denote the critical point on $D_\infty$ with Morse index $2(i-\ell-1)$.
\begin{thm}\label{thm: CF(H_l) + continuation maps}
    For every $\ell\in \Z_{\geq0}$ there exist signs $$A_{\ell,1},\ldots,A_{\ell,\ell},B_{\ell,1},\ldots,B_{\ell,\ell},\check{C}_{\ell,0},\ldots,\check{C}_{\ell,\ell+n+1},\hat{C}_{\ell,1},\ldots,\hat{C}_{\ell,\ell}\in \{-1,1\},$$

    such that for every $\ell\in \Z_{\geq0}$, the differential $d\fc CF(H_\ell)\to CF(H_\ell)$ satisfies
    \begin{itemize}
        \item  $$d \hat{x}^\ell_i=A_{\ell,i}T^{E(\hat{x}_{i}^\ell,\check{x}_{i-1}^\ell)}\check{x}^\ell_{i-1}+B_{\ell,i}T^{E(\hat{x}_{i}^\ell,\check{x}_{i+n}^\ell)}\check{x}^\ell_{i+n},$$
    where $$E(\hat{x}_{i}^\ell,\check{x}_{i-1}^\ell)=h_{\ell,i-1}^\Delta-h_{\ell,i}^\Delta+i(r_{\ell,i}^\Delta-r_{\ell,i-1}^\Delta)+r_{\ell,i-1}^\Delta$$ and $$E(\hat{x}_{i}^\ell,\check{x}_{i+n}^\ell)=\left\{\begin{array}{ll}
        h_{\ell,i+n}^\Delta-h_{\ell,i}^\Delta+i(r_{\ell,i}^\Delta-r_{\ell,i+n}^\Delta) +n(1-r_{\ell,i+n}^\Delta), & i\leq \ell-n, \\
        h(\Delta,\ell,1)-h_{\ell,i}^\Delta+ir_{\ell,i}^\Delta-i, & i\geq\ell-n+1,
    \end{array}\right.$$
    for every $1\leq i\leq \ell$, 
    \item $d\check{x}^\ell_i=0$ for every $0\leq i\leq \ell+n$,
    \end{itemize}
   and the continuation map $\Phi_\ell\fc CF(H_\ell)\to CF(H_{\ell+1})$ satisfies 
   
    \begin{itemize}
        \item  $$\Phi_\ell \check{x}^\ell_i=\check{C}_{\ell,i}T^{E(\check{x}_i^\ell,\check{x}_i^{\ell+1})}\check{x}^{\ell+1}_{i},$$
    where 
    $$E(\check{x}_i^\ell,\check{x}_i^{\ell+1})=\left\{\begin{array}{ll}
       h_{\ell+1,i}^\Delta-h_{\ell,i}^\Delta-i(r_{\ell+1,i}^\Delta-r_{\ell,i}^\Delta), & i\leq \ell, \\
       h_{\ell+1,\ell+1}^\Delta+(\ell+1)(1-r_{\ell+1,\ell+1}^\Delta)-h(\Delta,\ell,1)  & i=\ell+1,\\
       h(\Delta,\ell+1,1)-h(\Delta,\ell,1),& i\geq\ell+2,
    \end{array}\right.$$
    for every $0\leq i\leq \ell+n$,

    \item $$\Phi_\ell \hat{x}^\ell_i=\hat{C}_{\ell,i}T^{E(\hat{x}_i^\ell,\hat{x}_i^{\ell+1})}\hat{x}^{\ell+1}_{i},$$
     where 
    $$E(\hat{x}_i^\ell,\hat{x}_i^{\ell+1})=h_{\ell+1,i}^\Delta-h_{\ell,i}^\Delta-i(r_{\ell+1,i}^\Delta-r_{\ell,i}^\Delta)$$
    for every $1\leq i\leq \ell$.

    \end{itemize}
\end{thm}

As described in Section~\ref{ss: energy}, the topological energy of a Floer or continuation flowline $u$ connecting the $1$-periodic orbits $p$ and $q$ is given by
\begin{equation}\label{eq: Energy determined by index difference}
    \begin{aligned}
        &\begin{aligned}
        E_{top}(u) = \cA_H(q, \hat{q}) &- \cA_H(p, \hat{p}) \\ &+ \frac{1}{2(n+1)} \left( \dim_u \widehat{\cM}_m(p, q) - \mu_{FMB}^\tau(q) + \mu_{FMB}^\tau(p) + 1 \right),
        \end{aligned}\\
    &\text{for a Floer flowline, and}\\
    &\begin{aligned}
     E_{top}(u) = \cA_H(q, \hat{q}) &- \cA_H(p, \hat{p}) \\ &+ \frac{1}{2(n+1)} \left( \dim_u \widehat{\cM}_m(p, q) - \mu_{FMB}^\tau(q) + \mu_{FMB}^\tau(p) \right),
     \end{aligned}\\
     &\text{for a continuation flowline.}
     \end{aligned}
\end{equation}
where $\hat{p}$ and $\hat{q}$ are cappings for $p$ and $q$, respectively. The trivialization $\tau$ is induced by these cappings, and $\cA_H$ denotes the action functional associated with the Hamiltonian in the case of the differential, or the homotopy in the case of continuation maps.

\begin{proof}[Proof of Theorem~\ref{thm: CF(H_l) + continuation maps}]
   Let $\ell \in \Z_{\geq 0}$. We begin by computing the topological energy of flowlines with cascades that contribute to the differential. In this case, the dimension $\dim_u \cM_m(p, q)$ equals $m$ for every $p, q$, and $u$. Thus, the topological energy of such a flowline with cascades is given by
$$ E_{top}(u) = \cA_H(q, \hat{q}) - \cA_H(p, \hat{p}) + \frac{1}{2(n+1)} \left( \mu_{FMB}^\tau(p) - \mu_{FMB}^\tau(q) + 1 \right). $$
Note that for constant orbits, we can choose the cappings to be constant. Consequently, their action equals the value of the Hamiltonian itself. Thus,
$$ \cA_{H_\ell}(\check{x}_0^\ell) = h(\Delta, \ell, 0)=h_{\ell,0}^\Delta, \qquad \text{and} \qquad \cA_{H_\ell}(\check{x}_{\ell+1}^\ell) = \cdots = \cA_{H_\ell}(\check{x}_{\ell+n}^\ell) = h(\Delta, \ell, 1). $$
For the orbits $\hat{x}_1^\ell, \check{x}_1^\ell, \ldots, \hat{x}_\ell^\ell, \check{x}_\ell^\ell$, we choose cappings contained in $\CP^n \setminus D_\infty$, which is an exact symplectic manifold symplectomorphic to $\Int B(1)$. Let $\lambda = \frac{1}{2} \sum_{i=1}^n (x_i dy_i - y_i dx_i)$ be a primitive of $\omega_0$ on $\Int B(1)$. Then, for every $1 \leq i \leq \ell$, we find that
$$ \cA_{H_\ell}(\hat{x}_i^\ell) = \cA_{H_\ell}(\check{x}_i^\ell) = \int_{S^1} H_\ell \circ \check{x}_i^\ell(t) \, dt + \int_{S^1} (\check{x}_i^\ell)^* \lambda = h_{\ell,i}^\Delta - i r_{\ell,i}^\Delta. $$
Additionally, in Section~\ref{ss: computations of FMB}, we found that the Floer--Morse--Bott indices for these $1$-periodic orbits with respect to the trivialization $\tau$ induced by these cappings are as follows:
\begin{itemize}
    \item $\mu_{FMB}^\tau(\check{x}_0^\ell) = 0$.
    \item For every $1 \leq i \leq \ell$, $\mu_{FMB}^\tau(\check{x}_i^\ell) = -2ni$.
    \item For every $1 \leq i \leq \ell$, $\mu_{FMB}^\tau(\hat{x}_i^\ell) = -2ni + 2n - 1$.
    \item For every $\ell+1 \leq i \leq \ell+n$, $\mu_{FMB}^\tau(\check{x}_i^\ell) = 2i$.
\end{itemize}

The only Floer flowlines with cascades that contribute to the differential are those connecting $\hat{x}_i^\ell$ to $\check{x}_{i-1}^\ell$ and $\check{x}_{i+n}^\ell$ for $1 \leq i \leq \ell$. Thus, for every $1 \leq i \leq \ell$, we have
$$ \mu_{FMB}^\tau(\hat{x}_{i}^\ell) -\mu_{FMB}^\tau(\check{x}^\ell_{i-1})= -2ni + 2n - 1-(-2n(i-1)) = -1,$$
and hence
\begin{align*}
E(\hat{x}_i^\ell, \check{x}_{i-1}^\ell) &= \cA_{H_\ell}(\check{x}_{i-1}^\ell) - \cA_{H_\ell}(\hat{x}_i^\ell) \\
&= h_{\ell,i-1}^\Delta - (i-1)r_{\ell,i-1}^\Delta - h_{\ell,i}^\Delta + ir_{\ell,i}^\Delta \\
&= h_{\ell,i-1}^\Delta - h_{\ell,i}^\Delta + i(r_{\ell,i}^\Delta - r_{\ell,i-1}^\Delta) + r_{\ell,i-1}^\Delta.
\end{align*}

Now consider the second case for $1 \leq i \leq \ell$. If $i \leq \ell - n$, we find
$$\mu_{FMB}^\tau(\hat{x}_{i}^\ell) - \mu_{FMB}^\tau(\check{x}^\ell_{i+n}) = -2ni + 2n - 1 - (-2n(i+n)) = 2n(n+1) - 1,$$
which implies 
\begin{align*}
E(\hat{x}_i^\ell, \check{x}_{i+n}^\ell) &= \cA_{H_\ell}(\check{x}_{i+n}^\ell) - \cA_{H_\ell}(\hat{x}_i^\ell) + \frac{1}{2(n+1)} \left( \mu_{FMB}^\tau(\hat{x}_i^\ell) - \mu_{FMB}^\tau(\check{x}_{i+n}^\ell) + 1 \right) \\
&= h_{\ell,i+n}^\Delta - (i+n)r_{\ell,i+n}^\Delta - h_{\ell,i}^\Delta + ir_{\ell,i}^\Delta + n \\
&= h_{\ell,i+n}^\Delta - h_{\ell,i}^\Delta + i(r_{\ell,i}^\Delta - r_{\ell,i+n}^\Delta) + n(1 - r_{\ell,i+n}^\Delta).
\end{align*}

On the other hand, if $i \geq \ell - n + 1$, we have
$$\mu_{FMB}^\tau(\hat{x}_{i}^\ell) - \mu_{FMB}^\tau(\check{x}^\ell_{i+n}) = -2ni + 2n - 1 - 2(i+n) = -2i(n+1) - 1,$$
yielding 
\begin{align*}
E(\hat{x}_i^\ell, \check{x}_{i+n}^\ell) &= \cA_{H_\ell}(\check{x}_{i+n}^\ell) - \cA_{H_\ell}(\hat{x}_i^\ell) + \frac{1}{2(n+1)} \left( \mu_{FMB}^\tau(\hat{x}_i^\ell) - \mu_{FMB}^\tau(\check{x}_{i+n}^\ell) + 1 \right) \\
&= h(\Delta, \ell, 1) - h_{\ell,i}^\Delta + ir_{\ell,i}^\Delta - i.
\end{align*}

We now compute the topological energy of flowlines with cascades that contribute to the continuation map $\Phi_\ell \fc CF(H_\ell) \to CF(H_{\ell+1})$. In this case, the dimension $\dim_u \cM_m(p, q) $ equals $m-1$ for every $p, q$, and $u$. Thus, the topological energy of such a flowline with cascades is given by
$$ E_{top}(u) = \cA_{H_{\ell+1}}(q) - \cA_{H_\ell}(p) + \frac{1}{2(n+1)} \left( \mu_{FMB}^\tau(p) - \mu_{FMB}^\tau(q) \right). $$
As before, for constant orbits of $H_{\ell+1}$, we can choose the cappings to be constant. Consequently, their action equals the value of the Hamiltonian itself. Thus, $\cA_{H_{\ell+1}}(\check{x}_0^{\ell+1}) = h(\Delta, \ell+1, 0) = h_{\ell+1,0}^\Delta$ and
$$ \cA_{H_{\ell+1}}(\check{x}_{\ell+2}^{\ell+1}) = \cdots = \cA_{H_{\ell+1}}(\check{x}_{\ell+n+1}^{\ell+1}) = h(\Delta, \ell+1, 1). $$
For the orbits $\hat{x}_1^{\ell+1}, \check{x}_1^{\ell+1}, \ldots, \hat{x}_{\ell+1}^{\ell+1}, \check{x}_{\ell+1}^{\ell+1}$, we again choose cappings contained in $\CP^n \setminus D_\infty$. Thus, for every $1 \leq i \leq \ell+1$, we find
$$ \cA_{H_{\ell+1}}(\hat{x}_i^{\ell+1}) = \cA_{H_{\ell+1}}(\check{x}_i^{\ell+1}) = \int_{S^1} H_{\ell+1} \circ \check{x}_i^{\ell+1}(t) \, dt + \int_{S^1} (\check{x}_i^{\ell+1})^* \lambda = h_{\ell+1,i}^\Delta - i r_{\ell+1,i}^\Delta. $$
Additionally, the Floer--Morse--Bott indices for these $1$-periodic orbits with respect to the trivialization $\tau$ are as follows:
\begin{itemize}
    \item $\mu_{FMB}^\tau(\check{x}_0^{\ell+1}) = 0$.
    \item For $1 \leq i \leq \ell+1$, $\mu_{FMB}^\tau(\check{x}_i^{\ell+1}) = -2ni$.
    \item For $1 \leq i \leq \ell+1$, $\mu_{FMB}^\tau(\hat{x}_i^{\ell+1}) = -2ni + 2n - 1$.
    \item For $\ell+2 \leq i \leq \ell+n+1$, $\mu_{FMB}^\tau(\check{x}_i^{\ell+1}) = 2i$.
\end{itemize}

Given $1 \leq i \leq \ell$, note that $\mu_{FMB}^\tau(\hat{x}_i^\ell) = \mu_{FMB}^\tau(\hat{x}_i^{\ell+1})$, and hence
\begin{align*}
    E(\hat{x}_i^\ell, \hat{x}_i^{\ell+1}) &= \cA_{H_{\ell+1}}(\hat{x}_i^{\ell+1}) - \cA_{H_\ell}(\hat{x}_i^\ell) \\
    &= h_{\ell+1,i}^\Delta - i r_{\ell+1,i}^\Delta - h_{\ell,i}^\Delta + i r_{\ell,i}^\Delta \\
    &= h_{\ell+1,i}^\Delta - h_{\ell,i}^\Delta - i(r_{\ell+1,i}^\Delta - r_{\ell,i}^\Delta).
\end{align*}

Now, let $0 \leq i \leq \ell+n$.
\begin{itemize}
    \item If $i \leq \ell$, then $\mu_{FMB}^\tau(\check{x}_i^\ell) = \mu_{FMB}^\tau(\check{x}_i^{\ell+1})$, yielding
    \begin{align*}
        E(\check{x}_i^\ell, \check{x}_i^{\ell+1}) &= \cA_{H_{\ell+1}}(\check{x}_i^{\ell+1}) - \cA_{H_\ell}(\check{x}_i^\ell) \\
           &= h_{\ell+1,i}^\Delta - i r_{\ell+1,i}^\Delta - h_{\ell,i}^\Delta + i r_{\ell,i}^\Delta \\
    &= h_{\ell+1,i}^\Delta - h_{\ell,i}^\Delta - i(r_{\ell+1,i}^\Delta - r_{\ell,i}^\Delta).
    \end{align*}

    \item If $i = \ell+1$, then
    $$\mu_{FMB}^\tau(\check{x}_{\ell+1}^\ell) - \mu_{FMB}^\tau(\check{x}_{\ell+1}^{\ell+1}) = 2\ell + 2 - (-2n(\ell+1)) = 2(n+1)(\ell+1).$$
    Therefore,
    \begin{align*}
        E(\check{x}_{\ell+1}^\ell, \check{x}_{\ell+1}^{\ell+1}) &= \cA_{H_{\ell+1}}(\check{x}_{\ell+1}^{\ell+1}) - \cA_{H_\ell}(\check{x}_{\ell+1}^\ell) + \frac{1}{2(n+1)} \cdot 2(n+1)(\ell+1) \\
        &= h_{\ell+1,\ell+1}^\Delta - (\ell+1)r_{\ell+1,\ell+1}^\Delta - h(\Delta, \ell, 1) + (\ell+1) \\
        &= h_{\ell+1,\ell+1}^\Delta + (\ell+1)(1 - r_{\ell+1,\ell+1}^\Delta) - h(\Delta, \ell, 1).
    \end{align*}

    \item If $i \geq \ell+2$, then $\mu_{FMB}^\tau(\check{x}_i^\ell) = \mu_{FMB}^\tau(\check{x}_i^{\ell+1})$, and hence
    $$ E(\check{x}_i^\ell, \check{x}_i^{\ell+1}) = \cA_{H_{\ell+1}}(\check{x}_i^{\ell+1}) - \cA_{H_\ell}(\check{x}_i^\ell) = h(\Delta, \ell+1, 1) - h(\Delta, \ell, 1). $$
\end{itemize}
This completes the proof.
\end{proof}

\subsection{Direct limit of Floer complexes}\label{ss: direct limits}

Throughout this section, we fix $\Delta \in (0,1)$. As usual, for every $\ell \in \Z_{\geq0}$, define $h_\ell \fc [0,1] \to \R$ by $h_\ell(r) = h(\Delta, \ell, r)$ for every $r \in [0,1]$ and focus on the Hamiltonians $H_\ell$ on $\CP^n$ given by 
$H_\ell(z) = h_\ell(\mu(z)) = h(\Delta, \ell, \mu(z))$, for every $z \in \CP^n$, as in Equation~\eqref{eq: acc. data for balls}.

The purpose of this section is to compute the direct limit of the 1-ray of the cochain Floer complexes
$$CF(H_0) \to CF(H_1) \to CF(H_2) \to \cdots$$
where the connecting maps are the continuation maps. The description of the Floer complexes and the continuation maps can be found in Theorem~\ref{thm: CF(H_l) + continuation maps}.

Define a cochain complex $(C, d)$ as follows: the module $C$ is defined as
$$C = \Lambda_{>0}\check{x}_0 \oplus \bigoplus_{j\in \N} \left( \Lambda_{>j\Delta}\check{x}_j \oplus \Lambda_{>j\Delta}\hat{x}_j \right),$$
and the differential $d$ satisfies $d(\lambda\check{x}_j) = 0$ for every $j \in \Z_{\geq0}$ and $\lambda \in \Lambda_{>j\Delta}$, and
$$d(\lambda\hat{x}_j) = \lambda\check{x}_{j-1} + T^n\lambda\check{x}_{j+n}$$
for every $j \in \N$ and $\lambda \in \Lambda_{>j\Delta}$.

\begin{thm}\label{thm: drct lim of CF(H_l)}
    The direct limit of the 1-ray
    $$CF(H_0) \to CF(H_1) \to CF(H_2) \to \cdots$$
    is isomorphic to $(C, d)$. Additionally, the canonical maps $(f_\ell\fc CF(H_\ell)\to C)_{\ell\geq0}$ satisfy:
     \begin{itemize}
        \item for every $0 \leq j \leq \ell$ we have $f_\ell(\check{x}_j^\ell) =  T^{-(h_{\ell,j}^\Delta-jr_{\ell,j}^\Delta)}\check{x}_j$, 
        \item for every $\ell+1 \leq j \leq \ell+n$ we have $f_\ell(\check{x}_j^\ell) = T^{j-h(\Delta,\ell,1)}\check{x}_j$,
        \item for every $1 \leq j \leq \ell$ we have $f_\ell(\hat{x}_j^\ell) =  T^{-(h_{\ell,j}^\Delta-jr_{\ell,j}^\Delta)}\hat{x}_j$,
    \end{itemize}
    for every $\ell\in \Z_{\geq0}$.
\end{thm}

First, let us show that we can ignore the signs in the differentials and continuation maps associated with the Floer complexes described in Theorem~\ref{thm: CF(H_l) + continuation maps}. Specifically, we will prove that one can change the bases of the complexes $CF(H_\ell)$ so that all signs appearing in the formulas for these differentials and continuation maps are equal to $1$. This result is a consequence of Lemma~\ref{lemma: positive signs for zigzag} and Lemma~\ref{lemma: positive signs for zigzag and maps}.

\begin{prop}\label{prop: omission of signs for CF}
    For every $\ell\in \Z_{\geq0}$ there is a basis
    $$\hat{x}_1^\ell,\ldots,\hat{x}_\ell^\ell,\check{x}_0^\ell,\ldots,\check{x}_{\ell+n}^\ell$$ 
    for $CF(H_\ell)$ such that the differential $d\fc CF(H_\ell)\to CF(H_\ell)$ satisfies
    \begin{itemize}
        \item  $$d \hat{x}^\ell_i=T^{E(\hat{x}_{i}^\ell,\check{x}_{i-1}^\ell)}\check{x}^\ell_{i-1}+T^{E(\hat{x}_{i}^\ell,\check{x}_{i+n}^\ell)}\check{x}^\ell_{i+n},$$
    where $$E(\hat{x}_{i}^\ell,\check{x}_{i-1}^\ell)=h_{\ell,i-1}^\Delta-h_{\ell,i}^\Delta+i(r_{\ell,i}^\Delta-r_{\ell,i-1}^\Delta)+r_{\ell,i-1}^\Delta$$ and $$E(\hat{x}_{i}^\ell,\check{x}_{i+n}^\ell)=\left\{\begin{array}{ll}
        h_{\ell,i+n}^\Delta-h_{\ell,i}^\Delta+i(r_{\ell,i}^\Delta-r_{\ell,i+n}^\Delta) +n(1-r_{\ell,i+n}^\Delta), & i\leq \ell-n, \\
        h(\Delta,\ell,1)-h_{\ell,i}^\Delta+ir_{\ell,i}^\Delta-i, & i\geq\ell-n+1,
    \end{array}\right.$$
    for every $1\leq i\leq \ell$, 
    \item $d\check{x}^\ell_i=0$ for every $0\leq i\leq \ell+n$,
    \end{itemize}
   and the continuation map $\Phi_\ell\fc CF(H_\ell)\to CF(H_{\ell+1})$ satisfies 
   
    \begin{itemize}
        \item  $$\Phi_\ell \check{x}^\ell_i=T^{E(\check{x}_i^\ell,\check{x}_i^{\ell+1})}\check{x}^{\ell+1}_{i},$$
    where 
    $$E(\check{x}_i^\ell,\check{x}_i^{\ell+1})=\left\{\begin{array}{ll}
       h_{\ell+1,i}^\Delta-h_{\ell,i}^\Delta-i(r_{\ell+1,i}^\Delta-r_{\ell,i}^\Delta), & i\leq \ell, \\
       h_{\ell+1,\ell+1}^\Delta+(\ell+1)(1-r_{\ell+1,\ell+1}^\Delta)-h(\Delta,\ell,1)  & i=\ell+1,\\
       h(\Delta,\ell+1,1)-h(\Delta,\ell,1),& i\geq\ell+2,
    \end{array}\right.$$
    for every $0\leq i\leq \ell+n$,

    \item $$\Phi_\ell \hat{x}^\ell_i=T^{E(\hat{x}_i^\ell,\hat{x}_i^{\ell+1})}\hat{x}^{\ell+1}_{i},$$
     where 
    $$E(\hat{x}_i^\ell,\hat{x}_i^{\ell+1})=h_{\ell+1,i}^\Delta-h_{\ell,i}^\Delta-i(r_{\ell+1,i}^\Delta-r_{\ell,i}^\Delta)$$
    for every $1\leq i\leq \ell$.
    \end{itemize}
\end{prop}

\begin{proof}
For simplicity, throughout this proof we omit the powers of the formal variable $T$ in the Novikov ring, as they do not affect the signs; thus, we translate the discussion to the case of cochain complexes with $\Z$-coefficients. In our proof, we will show that by changing the signs of our given basis for $CF(H_\ell;\Z)$ for every $\ell\in \Z_{\geq0}$, we can achieve the desired result of this proposition.

By Theorem~\ref{thm: CF(H_l) + continuation maps}, for every $\ell \in \Z_{\geq 0}$, the cochain complex $CF(H_\ell)$ decomposes into $n+1$ subcomplexes. For each $0 \leq j \leq n$, denote by $C_{\ell,j}$ the subcomplex generated by
$$
\check{x}_{j}^\ell,\ \hat{x}_{j+1}^\ell,\ \check{x}_{j+n+1}^\ell,\ \hat{x}_{j+n+2}^\ell,\ \ldots,\ 
\hat{x}_{j+(k_{\ell,j}-1)(n+1)+1}^\ell,\ \check{x}_{j+k_{\ell,j}(n+1)}^\ell,
$$
where $k_{\ell,j} = \left\lceil \frac{\ell - j}{n+1} \right\rceil$. 

Fix $0 \leq j \leq n$. Define $C_{-1,j}=0$, let $\Phi_{-1}\fc C_{-1,j}\to C_{0,j}$ be the zero map, and denote $k_{-1,j}=-1$. Let us prove by induction that for every $\ell\in \Z_{\geq0}$, we can change the signs of the basis elements of the subcomplex $C_{\ell,j}$ such that:
\begin{itemize}
    \item For every $0\leq i\leq k_{\ell,j}-1$ we have $d\hat{x}_{j+i(n+1)+1}^\ell=\check{x}_{j+i(n+1)}^\ell+\check{x}_{j+(i+1)(n+1)}^\ell$;
    \item For every $0\leq i\leq k_{\ell-1,j}-1$ we have $\Phi_{\ell-1}(\hat{x}_{j+i(n+1)+1}^{\ell-1})= \hat{x}_{j+i(n+1)+1}^{\ell}$;
    \item For every $0\leq i\leq k_{\ell-1,j}$ we have $\Phi_{\ell-1}(\check{x}_{j+i(n+1)}^{\ell-1})= \check{x}_{j+i(n+1)}^{\ell}$.
\end{itemize} 

\textbf{The base case of the induction:} For $\ell=0$, we have $k_{0,j}=0$, therefore all the conditions of the induction statement are satisfied in the trivial sense.
 
\textbf{The step case of the induction:} Let $\ell\in \N$. Assume that for every $0\leq \ell'\leq \ell$, we can change the signs of the basis elements of the subcomplex $C_{\ell',j}$ such that:
\begin{itemize}
    \item For every $0\leq i\leq k_{\ell',j}-1$ we have $d\hat{x}_{j+i(n+1)+1}^{\ell'}=\check{x}_{j+i(n+1)}^{\ell'}+\check{x}_{j+(i+1)(n+1)}^{\ell'}$;
    \item For every $0\leq i\leq k_{\ell'-1,j}-1$ we have $\Phi_{\ell'-1}(\hat{x}_{j+i(n+1)+1}^{\ell'-1})= \hat{x}_{j+i(n+1)+1}^{\ell'}$;
    \item For every $0\leq i\leq k_{\ell'-1,j}$ we have $\Phi_{\ell'-1}(\check{x}_{j+i(n+1)}^{\ell'-1})= \check{x}_{j+i(n+1)}^{\ell'}$.
\end{itemize} 
 
If $\ell\leq j$, then $k_{\ell,j}=0$, and therefore the first two conditions of the induction statement are satisfied in the trivial sense. By Theorem~\ref{thm: continuation maps, over Z} we know that there is $\epsilon\in\{-1,1\}$ such that $\Phi(\check{x}_j^{\ell-1})=\epsilon\check{x}_j^\ell$, thus replacing $\check{x}_j^\ell$ by $\epsilon\check{x}_j^\ell$ completes the proof in this case. 

If $\ell=j+1$, then the second condition is trivial. Moreover, we have
$$C_{\ell-1,j}=\Z\cdot \check{x}_j^{\ell-1}\qquad\text{and}\qquad C_{\ell,j}=\Z\cdot \check{x}_j^{\ell}\oplus \Z\cdot \hat{x}_{j+1}^{\ell}\oplus \Z\cdot \check{x}_{j+n+1}^{\ell},$$
and by Theorem~\ref{thm: CF(H_l;Z)} and Theorem~\ref{thm: continuation maps, over Z} there are $\alpha,\beta,\gamma\in\{-1,1\}$ such that $\Phi(\check{x}_j^{\ell-1})=\alpha\check{x}_j^\ell$ and $d\hat{x}_{j+1}^\ell=\beta\check{x}_j^\ell+\gamma\check{x}_{j+n+1}^\ell$. Then replace the basis $\check{x}_j^\ell, \hat{x}_{j+1}^\ell, \check{x}_{j+n+1}^\ell$ by the basis $\alpha\check{x}_j^\ell, \alpha\beta \hat{x}_{j+1}^\ell, \alpha\beta\gamma\check{x}_{j+n+1}^\ell$. This change of basis works since
$$d(\alpha\beta\hat{x}_{j+1}^\ell)=\alpha\beta^2\check{x}_j^\ell+\alpha\beta\gamma\check{x}_{j+n+1}^\ell=\alpha\check{x}_j^\ell+\alpha\beta\gamma\check{x}_{j+n+1}^\ell.$$
Thus, the first and the third conditions hold as well.

Otherwise, assume that $\ell\geq j+2$. This case is fully covered by Lemma~\ref{lemma: positive signs for zigzag and maps}, and this completes the induction step and the whole proof.

\end{proof}

Our proof of Theorem~\ref{thm: drct lim of CF(H_l)} requires two more technical results as preparation.

\begin{lemma}\label{lemma: drct lim of 1-ray}
	Let $(\alpha_i)_{i=1}^\infty$ be a bounded strictly increasing sequence. For every $i \in \N$, define the homomorphism $T_i \fc \Lambda_{\geq 0} \to \Lambda_{\geq 0}$ by $T_i(x) = T^{\alpha_{i+1}-\alpha_i}x$ for every $x \in \Lambda_{\geq 0}$. Denote by $A$ the limit of the sequence $(\alpha_i)_{i=1}^\infty$.
	Then the direct limit of the 1-ray
	$$\Lambda_{\geq 0} \overset{T_1}{\longrightarrow} \Lambda_{\geq 0} \overset{T_2}{\longrightarrow} \Lambda_{\geq 0} \overset{T_3}{\longrightarrow} \cdots$$
	is isomorphic to $\Lambda_{>-A}$. 
	
	Moreover, the canonical maps $f_i \fc \Lambda_{\geq 0} \to \Lambda_{>-A}$ given by $f_i(x) = T^{-\alpha_i}x$ for every $x \in \Lambda_{\geq 0}$ and $i \in \N$ satisfy $f_j \circ T_i^j = f_i$ for every $j > i \geq 1$, where $T_i^j = T_{j-1} \circ \cdots \circ T_i$.
\end{lemma}

\begin{proof}
	For every $j > i \geq 1$ and $x \in \Lambda_{\geq 0}$, we have
	$$f_j \circ T_i^j(x) = f_j(T^{\alpha_j-\alpha_i}x) = T^{-\alpha_j} T^{\alpha_j-\alpha_i} x = T^{-\alpha_i} x = f_i(x).$$
	This shows that $\im f_i \subset \im f_j$. Since the maps $(f_i)_{i \in \N}$ are injective, we conclude that $\bigcup_{i \in \N} \im f_i \subset \Lambda$ is a direct limit of the 1-ray (see \cite[Example 5.32]{Rotman_intro_Homo_alg}). Since $\alpha_i \nearrow A$, it follows that $\bigcup_{i \in \N} \im f_i = \Lambda_{>-A}$.
\end{proof}

\begin{lemma}\label{lemma: drct lim in mod cat is drct lim in ch cat}
	Let
	$$A_1 \overset{\phi_1}{\longrightarrow} A_2 \overset{\phi_2}{\longrightarrow} A_3 \overset{\phi_3}{\longrightarrow} \cdots$$
	be a 1-ray of $\Lambda_{\geq 0}$-cochain complexes. Assume that this 1-ray has a direct limit in the category of $\Lambda_{\geq 0}$-modules, denoted by $A$. Then $A$ admits a differential $d$ such that $(A, d)$ is a direct limit of that 1-ray in the category of cochain complexes.
\end{lemma}

\begin{proof}
	For every $1 \leq i < j$, denote $\phi_i^j = \phi_{j-1} \circ \cdots \circ \phi_i$, and $\phi_i^i = \id$ for every $i \in \N$.
	
	Since $A$ is a direct limit of the 1-ray as a module, there are canonical maps $f_i \fc A_i \to A$ for every $i \in \N$ such that $f_i = f_j \circ \phi_i^j$ for every $1 \leq i < j$.
	
	Let us define a map $d \fc A \to A$ as follows: for any $x \in A$, there exists $i \in \N$ and $x_i \in A_i$ such that $f_i(x_i) = x$. We define $dx = f_i(d_i(x_i))$, where $d_i \fc A_i \to A_i$ is the differential of the complex $A_i$.  
	
	To show that $d$ is well-defined, let $x \in A$ and assume there are $i, j \in \N$ and $x_i \in A_i, x_j \in A_j$ such that $f_i(x_i) = f_j(x_j) = x$. By the properties of the direct limit, there exists $k \geq i, j$ such that $\phi_i^k(x_i) = \phi_j^k(x_j)$. Since $\phi_i^k$ and $\phi_j^k$ are cochain maps, we deduce:
	$$f_i(d_i x_i) = f_k(\phi_i^k(d_i x_i)) = f_k(d_k \phi_i^k x_i) = f_k(d_k \phi_j^k x_j) = f_k(\phi_j^k(d_j x_j)) = f_j(d_j x_j).$$
	Thus, $d$ is well-defined. The fact that $d^2 = 0$ is straightforward since $d_i^2 = 0$.
	
	Finally, let us show that $(A, d)$ satisfies the universal property in the category of cochain complexes. Let $(Y, d_Y)$ be a cochain complex and let $g_i \fc A_i \to Y$ be cochain maps satisfying $g_i = g_j \circ \phi_i^j$ for every $1 \leq i < j$. By the universal property of $A$ in the category of modules, there exists a unique module morphism $g \fc A \to Y$ such that $g_i = g \circ f_i$ for every $i \in \N$. To see that $g$ is a cochain map, let $x \in A$ with $x = f_i(x_i)$. Then:
	$$d_Y(g(x)) = d_Y(g(f_i(x_i))) = d_Y(g_i(x_i)) = g_i(d_i(x_i)) = g(f_i(d_i x_i)) = g(dx),$$
	as required.
\end{proof}
Finally, we can present the proof of Theorem~\ref{thm: drct lim of CF(H_l)}.

\begin{proof}[Proof of Theorem~\ref{thm: drct lim of CF(H_l)}]
	
	Let us compute the direct limit of the 1-ray
    
    $$CF(H_0) \xrightarrow{\Phi_0} CF(H_1) \xrightarrow{\Phi_1} CF(H_2) \longrightarrow \cdots.$$
	
	Let $j \in \Z_{\geq 0}$. The sequence $(h_{\ell,j}^\Delta-jr_{\ell,j}^\Delta)_{\ell=j}^\infty$ is strictly increasing and converges to $-j\Delta$. By Lemma~\ref{lemma: drct lim of 1-ray}, the module $\Lambda_{>j\Delta}$ is a direct limit for the 1-ray
	$$\Lambda_{\geq 0}\check{x}_j^j \xrightarrow{\Phi_j} \Lambda_{\geq 0}\check{x}_j^{j+1} \xrightarrow{\Phi_{j+1}} \Lambda_{\geq 0}\check{x}_j^{j+2} \xrightarrow{\Phi_{j+2}} \cdots.$$
	
	If $j \geq 1$, the same argument shows that $\Lambda_{>j\Delta}$ is a direct limit for the 1-ray
	$$\Lambda_{\geq 0}\hat{x}_j^j \xrightarrow{\Phi_j} \Lambda_{\geq 0}\hat{x}_j^{j+1} \xrightarrow{\Phi_{j+1}} \Lambda_{\geq 0}\hat{x}_j^{j+2} \xrightarrow{\Phi_{j+2}} \cdots.$$
	
	Since direct limits and direct sums commute, we conclude that the $\Lambda_{\geq 0}$-module
	$$C = \Lambda_{>0}\check{x}_0 \oplus \bigoplus_{j=1}^\infty \left( \Lambda_{>j\Delta} \check{x}_j \oplus \Lambda_{>j\Delta} \hat{x}_j \right)$$ 
	is a direct limit in the category of modules of the 1-ray
			$$CF(H_0) \xrightarrow{\Phi_0} CF(H_1)\xrightarrow{\Phi_{1}} CF(H_2)\xrightarrow{\Phi_{2}} \cdots.$$

	For every $\ell \in \Z_{\geq 0}$, let $f_\ell \fc CF(H_\ell) \to C$ denote the $\Lambda_{\geq 0}$-linear homomorphisms guaranteed by the universal property of the direct limit. By Lemma~\ref{lemma: drct lim of 1-ray}, we deduce that for every $\ell \in \Z_{\geq 0}$:
    \begin{itemize}
        \item for every $0 \leq j \leq \ell$ we have $f_\ell(\check{x}_j^\ell) =  T^{-(h_{\ell,j}^\Delta-jr_{\ell,j}^\Delta)}\check{x}_j$, 
        \item for every $\ell+1 \leq j \leq \ell+n$ we have $f_\ell(\check{x}_j^\ell) = T^{j-h(\Delta,\ell,1)}\check{x}_j$,
        \item for every $1 \leq j \leq \ell$ we have $f_\ell(\hat{x}_j^\ell) =  T^{-(h_{\ell,j}^\Delta-jr_{\ell,j}^\Delta)}\hat{x}_j$.
    \end{itemize}
    Additionally, these maps satisfy $f_i = f_j \circ \Phi_i^j$ for every $j > i \geq 0$, where $\Phi_i^j = \Phi_{j-1} \circ \cdots \circ \Phi_i$.
	Lemma~\ref{lemma: drct lim in mod cat is drct lim in ch cat} implies that there exists a differential $d$ on $C$ such that $(C, d)$ is a direct limit in the category of cochain complexes of the 1-ray, and $f_j \fc CF(H_\ell) \to (C, d)$ is a cochain map for every $j \in \Z_{\geq 0}$. 
	
	Now, let us compute the differential $d$:
	\begin{itemize}
		\item Let $j \in \Z_{\geq 0}$ and $\lambda \in \Lambda_{>j\Delta}$. There exists $\ell \geq j$ such that $\lambda T^{h_{\ell,j}^\Delta-jr_{\ell,j}^\Delta} \in \Lambda_{\geq 0}$ (note that while $h_{\ell,j}^\Delta-jr_{\ell,j}^\Delta$ is negative, the valuation of $\lambda$ is greater than $j\Delta$), which implies $f_\ell(T^{h_{\ell,j}^\Delta-jr_{\ell,j}^\Delta}\lambda \check{x}_j^\ell) = \lambda\check{x}_j$. Since $d_\ell\check{x}_j^\ell = 0$ and $f_\ell$ is a cochain map, we have:
		$$d(\lambda\check{x}_j) = d(f_\ell(T^{h_{\ell,j}^\Delta-jr_{\ell,j}^\Delta}\lambda \check{x}_j^\ell)) = f_\ell(d_\ell(T^{h_{\ell,j}^\Delta-jr_{\ell,j}^\Delta}\lambda \check{x}_j^\ell)) = 0.$$
		
		\item Let $j \in \N$ and $\lambda \in \Lambda_{>j\Delta}$. As before, there exists $\ell \geq j+n$ such that $\lambda T^{h_{\ell,j}^\Delta-jr_{\ell,j}^\Delta} \in \Lambda_{\geq 0}$, so $f_\ell(T^{h_{\ell,j}^\Delta-jr_{\ell,j}^\Delta}\lambda \hat{x}_j^\ell) = \lambda\hat{x}_j$. Since 
		$$d_\ell\hat{x}^\ell_j=T^{E(\hat{x}_{j}^\ell,\check{x}_{j-1}^\ell)}\check{x}^\ell_{j-1}+T^{E(\hat{x}_{j}^\ell,\check{x}_{j+n}^\ell)}\check{x}^\ell_{j+n},$$
    where $$E(\hat{x}_{j}^\ell,\check{x}_{j-1}^\ell)=h_{\ell,j-1}^\Delta-h_{\ell,j}^\Delta+j(r_{\ell,j}^\Delta-r_{\ell,j-1}^\Delta)+r_{\ell,j-1}^\Delta$$ and $$E(\hat{x}_{j}^\ell,\check{x}_{j+n}^\ell)=h_{\ell,j+n}^\Delta-h_{\ell,j}^\Delta+j(r_{\ell,j}^\Delta-r_{\ell,j+n}^\Delta) +n(1-r_{\ell,j+n}^\Delta).$$
		and $f_\ell$ is a cochain map, we compute:
		\begin{align*}
			d(\lambda\hat{x}_j) &= d(f_\ell(\lambda T^{h_{\ell,j}^\Delta-jr_{\ell,j}^\Delta} \hat{x}_j^\ell)) \\
			&= \lambda T^{h_{\ell,j}^\Delta-jr_{\ell,j}^\Delta} f_\ell(d_\ell \hat{x}_j^\ell) \\
			&= \lambda T^{h_{\ell,j}^\Delta-jr_{\ell,j}^\Delta} f_\ell\left(T^{E(\hat{x}_{j}^\ell,\check{x}_{j-1}^\ell)}\check{x}^\ell_{j-1}+T^{E(\hat{x}_{j}^\ell,\check{x}_{j+n}^\ell)}\check{x}^\ell_{j+n}\right) 
		\end{align*}
        Since 
        $$ \lambda T^{h_{\ell,j}^\Delta-jr_{\ell,j}^\Delta} f_\ell\left(T^{E(\hat{x}_{j}^\ell,\check{x}_{j-1}^\ell)}\check{x}^\ell_{j-1}\right)= \lambda T^{h_{\ell,j}^\Delta-jr_{\ell,j}^\Delta} T^{-(h_{\ell,j-1}^\Delta-(j-1)r_{\ell,j-1}^\Delta)}T^{E(\hat{x}_{j}^\ell,\check{x}_{j-1}^\ell)}\check{x}_{j-1}$$

        and
        $$ h_{\ell,j}^\Delta-jr_{\ell,j}^\Delta-(h_{\ell,j-1}^\Delta-(j-1)r_{\ell,j-1}^\Delta)+E(\hat{x}_{j}^\ell,\check{x}_{j-1}^\ell)=0,$$
        we deduce that

    $$ \lambda T^{h_{\ell,j}^\Delta-jr_{\ell,j}^\Delta} f_\ell\left(T^{E(\hat{x}_{j}^\ell,\check{x}_{j-1}^\ell)}\check{x}^\ell_{j-1}\right)= \lambda \check{x}_{j-1}$$

Similarly, since
$$\lambda T^{h_{\ell,j}^\Delta-jr_{\ell,j}^\Delta} f_\ell\left(T^{E(\hat{x}_{j}^\ell,\check{x}_{j+n}^\ell)}\check{x}^\ell_{j+n}\right)=\lambda T^{h_{\ell,j}^\Delta-jr_{\ell,j}^\Delta} T^{-(h_{\ell,j+n}^\Delta-(j+n)r_{\ell,j+n}^\Delta)}T^{E(\hat{x}_{j}^\ell,\check{x}_{j+n}^\ell)}\check{x}_{j+n}$$

        and
        $$h_{\ell,j}^\Delta-jr_{\ell,j}^\Delta-(h_{\ell,j+n}^\Delta-(j+n)r_{\ell,j+n}^\Delta)+E(\hat{x}_{j}^\ell,\check{x}_{j+n}^\ell)=n,$$
        
we deduce that

$$\lambda T^{h_{\ell,j}^\Delta-jr_{\ell,j}^\Delta} f_\ell\left(T^{E(\hat{x}_{j}^\ell,\check{x}_{j+n}^\ell)}\check{x}^\ell_{j+n}\right)=\lambda T^{n}\check{x}_{j+n}.$$
This shows that
$$d(\lambda\hat{x}_j)=\lambda \check{x}_{j-1}+\lambda T^{n}\check{x}_{j+n}$$
		as required.
	\end{itemize}
\end{proof}

\subsection{Proof of Theorem~\ref{thm: SH of ball in CP^n}}\label{ss: proof of main thm}
Fix $\Delta \in (0,1)$. Theorem~\ref{thm: drct lim of CF(H_l)} asserts that the direct limit of the 1-ray
$$CF(H_0) \to CF(H_1) \to CF(H_2) \to \cdots$$
is isomorphic to the cochain complex $(C, d)$ defined by
$$C = \Lambda_{>0}\check{x}_0 \oplus \bigoplus_{j \in \N} \left( \Lambda_{>j\Delta}\check{x}_j \oplus \Lambda_{>j\Delta}\hat{x}_j \right),$$
where the differential $d$ satisfies $d(\lambda\check{x}_j) = 0$ for every $j \in \Z_{\geq 0}$ and $\lambda \in \Lambda_{>j\Delta}$, and
$$d(\lambda\hat{x}_j) = \lambda\check{x}_{j-1} + T^n \lambda\check{x}_{j+n}$$
for every $j \in \N$ and $\lambda \in \Lambda_{>j\Delta}$.

In this section, we compute the homology of the completion of $(C, d)$. We start with a new presentation of $(C, d)$. Define a cochain complex $(\tilde{C}, \tilde{d})$ by
$$\tilde{C} = \Lambda_{>0}\check{x}_0 \oplus \bigoplus_{j \in \N} \left( \Lambda_{>0}\check{x}_j \oplus \Lambda_{>0}\hat{x}_j \right),$$
where the differential $\tilde{d}$ satisfies $\tilde{d}(\lambda\check{x}_j) = 0$ for every $j \in \Z_{\geq 0}$ and $\lambda \in \Lambda_{>0}$, and
$$\tilde{d}(\lambda\hat{x}_j) = T^\Delta \lambda\check{x}_{j-1} + T^{n(1-\Delta)} \lambda\check{x}_{j+n}$$
for every $j \in \N$ and $\lambda \in \Lambda_{>0}$. Now, define a map $K \fc \tilde{C} \to C$ by 
$$K = (\id, T^\Delta \id, T^\Delta \id, T^{2\Delta} \id, T^{2\Delta} \id, T^{3\Delta} \id, T^{3\Delta} \id, \dots)$$

\begin{claim}\label{claim: new presentation of drct lim of CF}
    The map $K$ is a cochain isomorphism.
\end{claim}

\begin{proof}
    The fact that $K$ is an isomorphism of modules is straightforward. We check that $K$ is a cochain map.
    For every $j \in \Z_{\geq 0}$ and $\lambda \in \Lambda_{>0}$, we have:
    $$K \circ \tilde{d}(\lambda\check{x}_j) = K(0) = 0 = d(T^{j\Delta}\lambda\check{x}_j) = d \circ K(\lambda\check{x}_j).$$
    Additionally, for every $j \in \N$ and $\lambda \in \Lambda_{>0}$, we have:
    \begin{align*}
        K \circ \tilde{d}(\lambda\hat{x}_j) &= K(T^\Delta \lambda\check{x}_{j-1} + T^{n(1-\Delta)} \lambda\check{x}_{j+n}) \\
        &= T^{(j-1)\Delta} T^\Delta \lambda\check{x}_{j-1} + T^{(j+n)\Delta} T^{n(1-\Delta)} \lambda\check{x}_{j+n} \\
        &= T^{j\Delta} \lambda\check{x}_{j-1} + T^{j\Delta} T^n \lambda\check{x}_{j+n} \\
        &= T^{j\Delta} (\lambda\check{x}_{j-1} + T^n \lambda\check{x}_{j+n}) \\
        &= d(T^{j\Delta} \lambda\hat{x}_j) \\
        &= d \circ K(\lambda\hat{x}_j).
    \end{align*}
    Therefore, $d \circ K = K \circ \tilde{d}$, and we conclude that $K$ is a cochain isomorphism.
\end{proof}
        
        Note that $\tilde{C}$ decomposes into $n+1$ subcomplexes $C_0, \dots, C_{n}$ defined by
$$C_j = \bigoplus\limits_{i=0}^\infty \left( \Lambda_{>0} \check{x}_{j+i(n+1)} \oplus \Lambda_{>0} \hat{x}_{j+i(n+1)+1} \right),$$
for every $0 \leq j \leq n$. Since,  in Section~\ref{s: indices}, we computed the Floer--Morse--Bott indices of the generators, see Table~\ref{tab:mu-fmb-values}, we know that $C_j^* \neq 0$ if and only if $* \equiv 2j \pmod{2(n+1)}$ or $* \equiv 2j-1 \pmod{2(n+1)}$, for every $0 \leq j \leq n$. Since completion is defined using the inverse limit, which commutes with finite direct sums, we deduce that to compute $H(\widehat{\tilde{C}})$, it is enough to compute the homology of the completion of each of these subcomplexes separately. 

\begin{claim}\label{claim: H(hat{C_j})}
For every $0 \leq j \leq n$:
\begin{itemize}
    \item If $\Delta < \frac{n}{n+1}$, then
    $$H^*(\widehat{C_j}) = \left\{\begin{array}{ll}
    \bigoplus_{i \geq 0} \Lambda_{(0,\Delta]}, & * \equiv 2j \pmod{2(n+1)}, \\
    0, & \text{otherwise}.
\end{array}\right.$$
    \item If $\Delta \geq \frac{n}{n+1}$, then
    $$H^*(\widehat{C_j}) = \left\{\begin{array}{ll}
    \Lambda_{>0} \oplus \bigoplus_{i \geq 0} \Lambda_{(0, n(1-\Delta)]}, & * \equiv 2j \pmod{2(n+1)}, \\
    0, & \text{otherwise}.
\end{array}\right.$$
\end{itemize}
\end{claim}

\begin{proof} 
Let $0 \leq j \leq n$. The subcomplex $(C_j, d)$ is nonzero only at degrees $2j-1$ and $2j$ modulo $2(n+1)$. We can identify $C_j^{2j-1}$ and $C_j^{2j}$ with two copies of $\Lambda_{>0}^\infty$ via the identifications:
$$\lambda\hat{x}_{j+i(n+1)+1} \mapsto \lambda e_{i}, \qquad\qquad \lambda\check{x}_{j+i(n+1)} \mapsto \lambda e_{i},$$
for every $i \geq 0$ and $\lambda\in \Lambda_{>0}$, where $(e_i)_{i \geq 0}$ is the standard basis of $\Lambda_{\geq 0}^\infty$. This identification is preserved under completion. According to Theorem~\ref{thm: drct lim of CF(H_l)}, the differential $\widehat{d} \fc \widehat{C_j} \to \widehat{C_j}$, whose restriction is nonzero only on $\widehat{C_j^{2j-1}}$, can be identified with the map $D \fc \widehat{\Lambda_{>0}^\infty} \to \widehat{\Lambda_{>0}^\infty}$ given by 
$$D(\lambda e_i) = T^\Delta \lambda e_i + T^{n(1-\Delta)} \lambda e_{i+1}$$
for every $i \geq 0$ and $\lambda\in \Lambda_{>0}$. This $D$ satisfies Equation~\eqref{eq: diff of drct CF(Ball)} with $\alpha = \Delta$ and $\beta = n(1-\Delta)$. 

By Proposition~\ref{prop: coker of D}, $D$ is injective, which implies $H^{2j-1}(\widehat{C_j}) = \ker D = 0$. Furthermore, the cohomology at degree $2j$ is given by $H^{2j}(\widehat{C_j}) = \coker D$. The condition $\Delta < \frac{n}{n+1}$ is equivalent to $\alpha < \beta$. Thus, applying the cases from Proposition~\ref{prop: coker of D} yields the forms guaranteed in the statement of Claim~\ref{claim: H(hat{C_j})}.
\end{proof}

\begin{proof}[Proof of Theorem~\ref{thm: SH of ball in CP^n}]
As we explained in Section~\ref{ss: relSH using MB} we can compute the relative symplectic cohomology using an acceleration datum of Hamiltonians that satisfy the \textbf{MB} condition. Therefore we have:
$$SH_{\CP^n}(B; \Lambda_{\geq 0}) = H\left( \widehat{\mathrm{tel}} \, CF(H_\ell) \right).$$
Applying Lemma~\ref{lemma: tel->colim}, we deduce that:
$$SH_{\CP^n}(B; \Lambda_{\geq 0}) = H\left( \widehat{\mathrm{tel}} \, CF(H_\ell) \right) \cong H\left( \widehat{\drctlim} \, CF(H_\ell) \right).$$
Since $ \widehat{\drctlim} \, CF(H_\ell)$ splits as a cochain complex as the direct sum $\bigoplus_{j=0}^n \wh C_j$, we get that
$$SH_{\CP^n}(B; \Lambda_{\geq 0}) =\bigoplus_{j=0}^n H(\wh C_j). $$
The homology $H(\wh C_j)$ has already been computed in Claim~\ref{claim: H(hat{C_j})}, for every $0\leq j\leq n$, which completes the proof.
\end{proof}

\subsection{Restriction maps, Proof of Theorem \ref{thm: res for balls in CP^n}}\label{ss: res maps}

Let $0<\Delta'\leq\Delta<1$. Consider the corresponding acceleration data $(H_\ell)_{\ell\geq0}$ and $(H'_\ell)_{\ell\geq0}$ for the balls $B^{}_{\Delta}=\mu^{-1}([0,\Delta])$ and $B^{}_{\Delta'}=\mu^{-1}([0,\Delta'])$, respectively, given by
$$H_\ell(z)=h(\Delta,\ell,\mu(z)),\qquad H'_\ell(z)=h(\Delta',\ell,\mu(z)),$$
for every $z\in \CP^n$ and $\ell\in \Z_{\geq0}$.

For every $\ell\in \Z_{\geq0}$ denote by $\check{x}_0^\ell,\ldots,\check{x}_{\ell+n}^\ell,\hat{x}_1^\ell,\ldots,\hat{x}_\ell^\ell$ the generators for $CF(H_\ell)$, as in Theorem~\ref{thm: CF(H_l) + continuation maps}, and by
$\check{y}_0^\ell,\ldots,\check{y}_{\ell+n}^\ell,\hat{y}_1^\ell,\ldots,\hat{y}_\ell^\ell$ the generators for $CF(H'_\ell)$, again as in Theorem~\ref{thm: CF(H_l) + continuation maps}.

By Proposition~\ref{prop: omission of signs for CF}, for every $\ell\in \Z_{\geq0}$, the bases 
$$\check{x}_0^\ell,\ldots,\check{x}_{\ell+n}^\ell,\hat{x}_1^\ell,\ldots,\hat{x}_\ell^\ell \qquad\text{and}\qquad \check{y}_0^\ell,\ldots,\check{y}_{\ell+n}^\ell,\hat{y}_1^\ell,\ldots,\hat{y}_\ell^\ell$$ 
for $CF(H_\ell)$ and $CF(H'_\ell)$, respectively, can be chosen such that the differentials and continuation maps in the $1$-rays
$$CF(H_0)\xrightarrow{\Phi_0}CF(H_1)\xrightarrow{\Phi_1}CF(H_2)\to\cdots$$
and 
$$CF(H'_0)\xrightarrow{\Phi'_0}CF(H'_1)\xrightarrow{\Phi'_1}CF(H'_2)\to\cdots$$
are represented by matrices with entries in $\{0, 1\}$ after setting $T=1$, namely, where we ignore the powers of the formal variable $T$ of the Novikov ring. We  refer to such bases as \textbf{normalized bases}.

Let $\chi\fc \R\to \R$ be a smooth non-increasing function, satisfying $\chi(s)=\Delta$ for every $s<0$ and $\chi(s)=\Delta'$ for every $s>1$.  Then for every $\ell\in \Z_{\geq0}$ one can define a monotone homotopy $K_\ell\fc \R\times \CP^n\to \R$ from $H_\ell$ to $H'_\ell$ by $K_\ell(s,z)=h(\chi(s),\ell,\mu(z))$ for every $(s,z)\in \R\times \CP^n$. Thus for every $\ell\in \Z_{\geq0}$ the homotopy $K_\ell$ defines a continuation map $\Psi_\ell\fc CF(H_\ell)\to CF(H'_\ell)$.
Similarly to Theorem~\ref{thm: CF(H_l) + continuation maps}, to describe the continuation maps with coefficients in the Novikov ring, in the sense of Section~\ref{sss: weighted_CF}, we must compute the energy of the Floer and continuation flowlines with cascades that appear in the differential and continuation maps.

\begin{thm}\label{thm: continuation maps for restrictions}
    For every $\ell\in \Z_{\geq0}$ there exist signs
    $$\check{D}_0,\ldots,\check{D}_{\ell+n},\hat{D}_1,\ldots,\hat{D}_{\ell}\in \{-1,1\}$$
    such that the continuation map $\Psi_\ell\fc CF(H_\ell)\to CF(H'_\ell)$ satisfies
    \begin{itemize}
    \item $\Phi'_\ell\circ \Psi_\ell=\Psi_{\ell+1}\circ\Phi_\ell$ where $\Phi_\ell\fc CF(H_\ell)\to CF(H_{\ell+1})$ and $\Phi'_\ell\fc CF(H'_\ell)\to CF(H'_{\ell+1})$ are the continuation maps.
    \item For every $1\leq i\leq \ell$ we have
        $$\Psi(\hat{x}_i^\ell)=\hat{D}_i T^{E(\hat{x}_i^\ell,\hat{y}_i^\ell)}\hat{y}_i^\ell,$$
        where
        $$E(\hat{x}_i^\ell,\hat{y}_i^\ell)= h_{\ell,i}^{\Delta'}-h_{\ell,i}^{\Delta}-i(r_{\ell,i}^{\Delta'}-r_{\ell,i}^{\Delta}).$$
        \item For every $0\leq i\leq \ell+n$ we have
        $$\Psi(\check{x}_i^\ell)=\check{D}_i T^{E(\check{x}_i^\ell,\check{y}_i^\ell)}\check{y}_i^\ell,$$
        where
        $$E(\check{x}_i^\ell,\check{y}_i^\ell)=\left\{\begin{array}{cc}
            h_{\ell,i}^{\Delta'}-h_{\ell,i}^{\Delta}-i(r_{\ell,i}^{\Delta'}-r_{\ell,i}^{\Delta}), & i\leq \ell, \\
            h(\Delta',\ell,1)-h(\Delta,\ell,1), & i\geq\ell+1, 
        \end{array}\right.$$

    \end{itemize}
\end{thm}
\begin{proof}
   Let $\ell \in \Z_{\geq 0}$. The fact that $\Phi'_\ell\circ \Psi_\ell=\Psi_{\ell+1}\circ\Phi_\ell$ is an immediate corollary of Proposition~\ref{prop: commuting of cont maps}.

By Theorem~\ref{thm: continuation maps, over Z}, each of the coefficients $\check{D}_0,\ldots,\check{D}_{\ell+n},\hat{D}_1,\ldots,\hat{D}_{\ell}\in \Z$, whose existence is guaranteed by Theorem~\ref{thm: CF(H_l) + continuation maps}, is either $1$ or $-1$.
   
   The rest of the proof is similar to that of Theorem~\ref{thm: CF(H_l) + continuation maps}. Since we are interested in continuation flowlines with cascades, the dimension $\dim_u \wh\cM(p, q)$ equals $0$ for every $p, q$, and $u$. Thus, the topological energy of such a flowline with cascades is given by
$$ E_{top}(u) = \cA_H(q, \hat{q}) - \cA_H(p, \hat{p}) + \frac{1}{2(n+1)} \left( \mu_{FMB}^\tau(p) - \mu_{FMB}^\tau(q) \right). $$
As in the proof of Theorem~\ref{thm: CF(H_l) + continuation maps}, for constant orbits of $H_\ell$ and $H'_\ell$, choose the cappings to be constant, thus their action equals the value of the Hamiltonian itself, that is
$$ \cA_{H_\ell}(\check{x}_0^\ell) = h(\Delta, \ell, 0)=h_{\ell,0}^\Delta, \qquad \text{and} \qquad \cA_{H_\ell}(\check{x}_{\ell+1}^\ell) = \cdots = \cA_{H_\ell}(\check{x}_{\ell+n}^\ell) = h(\Delta, \ell, 1), $$
and 
$$ \cA_{H'_\ell}(\check{y}_0^\ell) = h(\Delta', \ell, 0)=h_{\ell,0}^{\Delta'}, \qquad \text{and} \qquad \cA_{H'_\ell}(\check{y}_{\ell+1}^\ell) = \cdots = \cA_{H'_\ell}(\check{y}_{\ell+n}^\ell) = h(\Delta', \ell, 1). $$
For the orbits $\hat{x}_1^\ell, \check{x}_1^\ell, \ldots, \hat{x}_\ell^\ell, \check{x}_\ell^\ell,\hat{y}_1^\ell, \check{y}_1^\ell, \ldots, \hat{y}_\ell^\ell, \check{y}_\ell^\ell$, we choose cappings contained in $\CP^n \setminus D_\infty$, which is an exact symplectic manifold symplectomorphic to $\Int B(1)$. For every $1 \leq i \leq \ell$, we find that
$$ \cA_{H_\ell}(\hat{x}_i^\ell) = \cA_{H_\ell}(\check{x}_i^\ell)  = h_{\ell,i}^\Delta - i r_{\ell,i}^\Delta,\qquad\text{and}\qquad  \cA_{H'_\ell}(\hat{y}_i^\ell) = \cA_{H'_\ell}(\check{y}_i^\ell)  = h_{\ell,i}^{\Delta'} - i r_{\ell,i}^{\Delta'}. $$
Additionally, in Section~\ref{ss: computations of FMB}, we found that the Floer--Morse--Bott indices for these $1$-periodic orbits of $H_\ell$ with respect to the trivialization $\tau$ induced by these cappings do not depend on $\Delta$, therefore for every $1\leq i\leq \ell$ we have $\mu_{FMB}^\tau(\hat{x}_{i}^\ell)=\mu_{FMB}^\tau(\hat{y}_{i}^\ell)$ and for every $0\leq i\leq \ell+n$ we have $\mu_{FMB}^\tau(\check{x}_{i}^\ell)=\mu_{FMB}^\tau(\check{y}_{i}^\ell)$. 

Now, given $1 \leq i \leq \ell$, we know that
\begin{align*}
    E(\hat{x}_i^\ell, \hat{y}_i^{\ell}) &= \cA_{H'_{\ell}}(\hat{y}_i^{\ell}) - \cA_{H_\ell}(\hat{x}_i^\ell) \\
    &= h_{\ell,i}^{\Delta'} - i r_{\ell,i}^{\Delta'} - h_{\ell,i}^\Delta + i r_{\ell,i}^\Delta \\
    &= h_{\ell,i}^{\Delta'} - h_{\ell,i}^\Delta - i(r_{\ell,i}^{\Delta'} - r_{\ell,i}^\Delta).
\end{align*}

Next, let $0 \leq i \leq \ell+n$.
\begin{itemize}
    \item If $i \leq \ell$, we have
    \begin{align*}
        E(\check{x}_i^\ell, \check{y}_i^{\ell}) &= \cA_{H'_{\ell}}(\check{y}_i^{\ell}) - \cA_{H_\ell}(\check{x}_i^\ell) \\
        &= h_{\ell,i}^{\Delta'} - i r_{\ell,i}^{\Delta'} - h_{\ell,i}^\Delta + i r_{\ell,i}^\Delta \\
    &= h_{\ell,i}^{\Delta'} - h_{\ell,i}^\Delta - i(r_{\ell,i}^{\Delta'} - r_{\ell,i}^\Delta).
    \end{align*}

    \item If $i \geq \ell+1$, we have
    $$ E(\check{x}_i^\ell, \check{y}_i^{\ell}) = \cA_{H'_{\ell}}(\check{y}_i^{\ell}) - \cA_{H_\ell}(\check{x}_i^\ell) = h(\Delta', \ell, 1) - h(\Delta, \ell, 1). $$
\end{itemize}
This completes the proof.

\end{proof}

The following result is a continuation of Proposition~\ref{prop: omission of signs for CF} and demonstrates that the signs from Theorem~\ref{thm: continuation maps for restrictions} can be normalized to $1$ via a suitable change of bases.

\begin{prop}\label{prop: omission of signs for CF->CF'}
     For every $\ell\in \Z_{\geq0}$, there are normalized bases
$$\check{x}_0^\ell,\ldots,\check{x}_{\ell+n}^\ell,\hat{x}_1^\ell,\ldots,\hat{x}_\ell^\ell \qquad\text{and}\qquad \check{y}_0^\ell,\ldots,\check{y}_{\ell+n}^\ell,\hat{y}_1^\ell,\ldots,\hat{y}_\ell^\ell$$ 
for $CF(H_\ell)$ and $CF(H'_\ell)$, respectively, such that the continuation map $\Psi_\ell\fc CF(H_\ell)\to CF(H'_\ell)$ satisfies
    \begin{itemize}
  
        \item For every $0\leq i\leq \ell+n$ we have
        $$\Psi(\check{x}_i^\ell)=T^{E(\check{x}_i^\ell,\check{y}_i^\ell)}\check{y}_i^\ell,$$
        where
        $$E(\check{x}_i^\ell,\check{y}_i^\ell)=\left\{\begin{array}{cc}
            h_{\ell,i}^{\Delta'}-h_{\ell,i}^{\Delta}-i(r_{\ell,i}^{\Delta'}-r_{\ell,i}^{\Delta}), & i\leq \ell, \\
            h(\Delta',\ell,1)-h(\Delta,\ell,1), & i\geq\ell+1, 
        \end{array}\right.$$
        
         \item For every $1\leq i\leq \ell$ we have
        $$\Psi(\hat{x}_i^\ell)= T^{E(\hat{x}_i^\ell,\hat{y}_i^\ell)}\hat{y}_i^\ell,$$
        where
        $$E(\hat{x}_i^\ell,\hat{y}_i^\ell)= h_{\ell,i}^{\Delta'}-h_{\ell,i}^{\Delta}-i(r_{\ell,i}^{\Delta'}-r_{\ell,i}^{\Delta}).$$
    \end{itemize}
\end{prop}

\begin{proof}

For every $\ell\in \Z_{\geq0}$, let
$$\check{x}_0^\ell,\ldots,\check{x}_{\ell+n}^\ell,\hat{x}_1^\ell,\ldots,\hat{x}_\ell^\ell$$ 
be a normalized basis for $CF(H_\ell)$, which is guaranteed by Proposition~\ref{prop: omission of signs for CF}.

As in the proof of Proposition~\ref{prop: omission of signs for CF}, for simplicity, throughout this proof we omit the powers of the formal variable $T$ in the Novikov ring, as they do not affect the signs. 
For every $\ell \in \Z_{\geq 0}$, the cochain complexes $CF(H_\ell;\Z)$ and $CF(H'_\ell;\Z)$ decompose into $n+1$ subcomplexes. For each $0 \leq j \leq n$, denote by $C_{\ell,j}$ and $C'_{\ell,j}$ the subcomplexes generated by
$$
\check{x}_{j}^\ell,\ \hat{x}_{j+1}^\ell,\ \check{x}_{j+n+1}^\ell,\ \hat{x}_{j+n+2}^\ell,\ \ldots,\ 
\hat{x}_{j+(k_{\ell,j}-1)(n+1)+1}^\ell,\ \check{x}_{j+k_{\ell,j}(n+1)}^\ell,
$$ 
and by
$$
\check{y}_{j}^\ell,\ \hat{y}_{j+1}^\ell,\ \check{y}_{j+n+1}^\ell,\ \hat{y}_{j+n+2}^\ell,\ \ldots,\ 
\hat{y}_{j+(k_{\ell,j}-1)(n+1)+1}^\ell,\ \check{y}_{j+k_{\ell,j}(n+1)}^\ell,
$$
respectively, where $k_{\ell,j} = \left\lceil \frac{\ell - j}{n+1} \right\rceil$.

Fix $0\leq j\leq n$ for the remainder of this proof. We will prove by induction that by changing the signs of the normalized bases of the subcomplexes $(C'_{\ell,j})_{\ell\geq0}$, we can ensure that for every $\ell \geq 0$, all entries of the representative matrix of the continuation map $\Psi_\ell\fc C_{\ell,j}\to C'_{\ell,j}$, with respect to these newly signed generators, are either $0$ or $1$.

\textbf{The base case of the induction:} For $\ell=0$, we have $C_{0,j}=\Z\cdot\check{x}_j^0$, $C'_{0,j}=\Z\cdot\check{y}_j^0$, and $\Psi_0(\check{x}_j^0)=\alpha \check{y}_j^0$ for some $\alpha\in\{-1,1\}$. Replacing $\check{y}_j^0$ by $\alpha \check{y}_j^0$ completes this case.

\textbf{The step case of the induction:} Let $\ell\in \Z_{\geq0}$ and assume that by changing the signs of the normalized bases of the subcomplexes $(C'_{\ell',j})_{0\leq \ell'\leq \ell}$, we can ensure that all entries of the representative matrix of the continuation map $\Psi_{\ell'}\fc C_{\ell',j}\to C'_{\ell',j}$, with respect to these newly signed generators, are either $0$ or $1$, for every $0\leq \ell'\leq \ell$.

By Proposition~\ref{thm: continuation maps for restrictions}, we have $\Phi'_\ell\circ\Psi_\ell=\Psi_{\ell+1}\circ\Phi_\ell$; therefore, Lemma~\ref{lemma: positive signs for zigzag and maps in commuting case} fully completes the step case, and thus the proof itself.

\end{proof}

Consider the cochain complexes $(C, d)$ and $(C', d')$, defined as follows: the modules $C$ and $C'$ are given by
$$C = \Lambda_{>0}\check{x}_0 \oplus \bigoplus_{j\in \N} \left( \Lambda_{>j\Delta}\check{x}_j \oplus \Lambda_{>j\Delta}\hat{x}_j \right),\qquad
C' = \Lambda_{>0}\check{y}_0 \oplus \bigoplus_{j\in \N} \left( \Lambda_{>j\Delta}\check{y}_j \oplus \Lambda_{>j\Delta}\hat{y}_j \right),$$
and the differentials $d$ and $d'$ satisfy $d(\lambda\check{x}_j) = 0$ and $d'(\lambda\check{y}_j) = 0$ for every $j \in \Z_{\geq0}$ and $\lambda \in \Lambda_{>j\Delta}$, while
$$d(\lambda\hat{x}_j) = \lambda\check{x}_{j-1} + T^n\lambda\check{x}_{j+n},\qquad d'(\lambda\hat{y}_j) = \lambda\check{y}_{j-1} + T^n\lambda\check{y}_{j+n}$$
for every $j \in \N$ and $\lambda \in \Lambda_{>j\Delta}$. By Theorem~\ref{thm: drct lim of CF(H_l)}, the direct limits of the $1$-rays
$$CF(H_0) \to CF(H_1) \to CF(H_2) \to \cdots$$
and 
$$CF(H'_0) \to CF(H'_1) \to CF(H'_2) \to \cdots$$
are isomorphic to $(C, d)$ and $(C', d')$, respectively.

Let $\Psi \fc C \to C'$ denote the direct limit of the continuation maps $(\Psi_\ell \fc CF(H_\ell) \to CF(H'_\ell))_{\ell\geq 0}$. Note that $\Psi$ is well-defined because Theorem~\ref{thm: continuation maps for restrictions} asserts that $\Psi_{\ell+1} \circ \Phi_\ell = \Phi'_\ell \circ \Psi_\ell$, which ensures that the maps $(\Psi_\ell)_{\ell\geq 0}$ define a morphism of directed systems. The following proposition describes this map.

\begin{prop}
    The map 
    $$\Psi \fc \Lambda_{>0}\check{x}_0 \oplus \bigoplus_{j\in \mathbb{N}} \left( \Lambda_{>j\Delta}\check{x}_j \oplus \Lambda_{>j\Delta}\hat{x}_j \right) \to \Lambda_{>0}\check{y}_0 \oplus \bigoplus_{j\in \mathbb{N}} \left( \Lambda_{>j\Delta'}\check{y}_j \oplus \Lambda_{>j\Delta'}\hat{y}_j \right)$$
    is the natural inclusion map. Specifically, for every $\lambda \in \Lambda_{>0}$, we have $\Psi(\lambda\check{x}_0) = \lambda\check{y}_0$, and for every $j \in \mathbb{N}$ and $\lambda \in \Lambda_{>j\Delta}$, we have $\Psi(\lambda\check{x}_j) = \lambda\check{y}_j$ and $\Psi(\lambda\hat{x}_j) = \lambda\hat{y}_j$.
\end{prop}

\begin{proof}
By the universal property of the direct limit, for every $\ell \in \Z_{\geq 0}$, the following diagram commutes:
$$
\xymatrix@R=2pc@C=3pc{
CF(H_\ell) \ar[r]^{f_\ell} \ar[d]_{\Psi_\ell} & C \ar[d]^{\Psi} \\
CF(H'_\ell) \ar[r]_{f'_\ell} & C'
}
$$
where $f_\ell$ and $f'_\ell$ are the canonical maps of the direct limits $C$ and $C'$, respectively. 

Let $j \in \Z_{\geq 0}$ and $\lambda \in \Lambda_{>j\Delta}$. From Theorem~\ref{thm: drct lim of CF(H_l)} and the properties of direct limits, for any $\ell \geq j$, the generator $\lambda \check{x}_j \in C$ can be represented as:
$$\lambda \check{x}_j = f_\ell \left( \lambda T^{h_{\ell,j}^\Delta - j r_{\ell,j}^\Delta} \check{x}_j^\ell \right).$$
Using the commutativity of the diagram, we have:
$$\Psi(\lambda \check{x}_j) = \Psi \circ f_\ell \left( \lambda T^{h_{\ell,j}^\Delta - j r_{\ell,j}^\Delta} \check{x}_j^\ell \right) = f'_\ell \circ \Psi_\ell \left( \lambda T^{h_{\ell,j}^\Delta - j r_{\ell,j}^\Delta} \check{x}_j^\ell \right).$$
By Proposition~\ref{prop: omission of signs for CF->CF'}, we have $\Psi_\ell(\check{x}_j^\ell) = T^{E(\check{x}_j^\ell, \check{y}_j^\ell)} \check{y}_j^\ell$. Thus:
$$\Psi(\lambda \check{x}_j) = \lambda T^{h_{\ell,j}^\Delta - j r_{\ell,j}^\Delta} f'_\ell \left( T^{E(\check{x}_j^\ell, \check{y}_j^\ell)} \check{y}_j^\ell \right).$$
Applying the formula for $f'_\ell \fc CF(H'_\ell) \to C'$ from Theorem~\ref{thm: drct lim of CF(H_l)}, we obtain:
$$\Psi(\lambda \check{x}_j) = \lambda T^{(h_{\ell,j}^\Delta - j r_{\ell,j}^\Delta) + E(\check{x}_j^\ell, \check{y}_j^\ell) - (h_{\ell,j}^{\Delta'} - j r_{\ell,j}^{\Delta'})} \check{y}_j.$$
By substituting the formula for $E(\check{x}_j^\ell, \check{y}_j^\ell) = h_{\ell,j}^{\Delta'} - h_{\ell,j}^\Delta - j(r_{\ell,j}^{\Delta'} - r_{\ell,j}^\Delta)$, we observe that the exponent is:
$$(h_{\ell,j}^\Delta - j r_{\ell,j}^\Delta) + (h_{\ell,j}^{\Delta'} - h_{\ell,j}^\Delta - j r_{\ell,j}^{\Delta'} + j r_{\ell,j}^\Delta) - (h_{\ell,j}^{\Delta'} - j r_{\ell,j}^{\Delta'}) = 0.$$
Therefore, $\Psi(\lambda \check{x}_j) = \lambda T^0 \check{y}_j = \lambda \check{y}_j$.

For the generators $\hat{x}_j$ where $j \in \N$ and $\lambda \in \Lambda_{>j\Delta}$, the same calculation applies. For $\ell \geq j$, we have $\lambda \hat{x}_j = f_\ell \left( \lambda T^{h_{\ell,j}^\Delta - j r_{\ell,j}^\Delta} \hat{x}_j^\ell \right)$. Since the energy term $E(\hat{x}_j^\ell, \hat{y}_j^\ell)$ is identical to that of the $\check{x}$ generators, we conclude:
$$\Psi(\lambda \hat{x}_j) = \lambda \hat{y}_j.$$
This completes the proof.
\end{proof}

As in the proof of Theorem~\ref{thm: drct lim of CF(H_l)}, define cochain complexes $(\tilde{C},\tilde{d})$ and $(\tilde{C}',\tilde{d}')$ over $\Lambda_{\geq 0}$ as follows:
$$\tilde{C}=\Lambda_{>0}\check{x}_0\oplus\bigoplus_{i=1}^\infty \left(\Lambda_{>0} \check{x}_i\oplus \Lambda_{>0} \hat{x}_i\right),$$
$$\tilde{C}'=\Lambda_{>0}\check{y}_0\oplus\bigoplus_{i=1}^\infty \left(\Lambda_{>0} \check{y}_i\oplus \Lambda_{>0} \hat{y}_i\right),$$
where for every $i\geq 0$ and $\lambda\in\Lambda_{>0}$ we have $\tilde{d}(\lambda \check{x}_i)=\tilde{d}'(\lambda \check{y}_i)=0$. Also, for every $i\geq 1$ and $\lambda\in\Lambda_{>0}$, the differentials satisfy
$$\tilde{d}(\lambda \hat{x}_i)=T^\Delta \lambda \check{x}_{i-1}+T^{n(1-\Delta)} \lambda \check{x}_{i+n},$$
$$\tilde{d}'(\lambda\hat{y}_i)=T^{\Delta'} \lambda\check{y}_{i-1}+T^{n(1-\Delta')}\lambda\check{y}_{i+n}.$$
By Claim~\ref{claim: new presentation of drct lim of CF}, the cochain complex $(\tilde{C},\tilde{d})$ is isomorphic to $(C,d)$ via the isomorphism $K\fc \tilde{C}\to C$ given by $K=\id\bigoplus_{i=1}^\infty T^{i\Delta}(\id\oplus \id)$. Similarly, $(\tilde{C}',\tilde{d}')$ is isomorphic to $(C',d')$ via the isomorphism $K'\fc \tilde{C}'\to C'$ given by $K'=\id\bigoplus_{i=1}^\infty T^{i\Delta'}(\id\oplus \id)$.
Define the cochain map $\tilde{\Psi} \fc \tilde{C} \to \tilde{C}'$ by $\tilde{\Psi} = (K')^{-1} \circ \Psi \circ K$. A straightforward computation shows that $\tilde{\Psi}$ satisfies 
$$\tilde{\Psi}(\lambda \check{x}_i) = T^{i(\Delta - \Delta')} \lambda \check{y}_i$$ 
for every $i \in \Z_{\geq 0}$ and $\lambda \in \Lambda_{> 0}$. Similarly, for every $i \in \N$ and $\lambda \in \Lambda_{> 0}$, we have 
$$\tilde{\Psi}(\lambda \hat{x}_i) = T^{i(\Delta - \Delta')} \lambda \hat{y}_i.$$

As seen in the proof of Theorem~\ref{thm: SH of ball in CP^n}, and as it is mentioned in Section~\ref{ss: relSH using MB}, we can compute the relative symplectic cohomology using an acceleration datum of Hamiltonians that satisfy the \textbf{MB} condition, therefore we have $SH_{\CP^n}(B;\Lambda_{\geq0}) \cong H(\widehat{\tilde{C}})$ and $SH_{\CP^n}(B';\Lambda_{\geq0}) \cong H(\widehat{\tilde{C}'})$. Our next objective is to compute the induced map $(\widehat{\Psi})_* \fc SH_{\CP^n}(B;\Lambda_{\geq0}) \to SH_{\CP^n}(B';\Lambda_{\geq0})$.

\begin{prop}\label{prop: drctlim restriction maps}
    Let $0 \leq j \leq n$. The morphism 
    $$H(\widehat{\Psi})_{2j} \fc SH_{\CP^n}^{2j}(B;\Lambda_{\geq0}) \to SH_{\CP^n}^{2j}(B';\Lambda_{\geq0})$$
    is determined by the following: 
    \begin{enumerate}
        \item If $\Delta' \leq \Delta < \frac{n}{n+1}$, then the map 
        $$H(\widehat{\Psi})_{2j} \fc \bigoplus_{i \geq 0} \Lambda_{(0, \Delta]} e_i \to \bigoplus_{i \geq 0} \Lambda_{(0, \Delta']} e'_i$$
        satisfies
        $$H(\widehat{\Psi})_{2j}(\lambda e_i) = T^{\hbigl(j+i(n+1)\hbigr)(\Delta-\Delta')} \lambda e'_i$$
        for every $i \geq 0$ and $\lambda \in \Lambda_{(0, \Delta]}$.

        \item If $\frac{n}{n+1} \leq \Delta' \leq \Delta$, then the map 
        $$H(\widehat{\Psi})_{2j} \fc \Lambda_{>0} e \oplus \bigoplus_{i \geq 0} \Lambda_{(0, n(1-\Delta)]} e_i \to \Lambda_{>0} e' \oplus \bigoplus_{i \geq 0} \Lambda_{(0, n(1-\Delta')]} e'_i$$
        satisfies 
        $$H(\widehat{\Psi})_{2j}(\lambda e) = T^{j(\Delta-\Delta')} \lambda e',$$
        for every $\lambda \in \Lambda_{>0}$, and
        $$H(\widehat{\Psi})_{2j}(\lambda e_i) = T^{(j+(i+1)(n+1))(\Delta-\Delta')} \lambda e'_i$$
        for every $i \geq 0$ and $\lambda \in \Lambda_{(0, n(1-\Delta)]}$.

        \item If $\Delta' < \frac{n}{n+1} \leq \Delta$, then the map 
        $$H(\widehat{\Psi})_{2j} \fc \Lambda_{>0} e \oplus \bigoplus_{i \geq 0} \Lambda_{(0, n(1-\Delta)]} e_i \to \bigoplus_{i \geq 0} \Lambda_{(0, \Delta']} e'_i$$
        satisfies 
        $$H(\widehat{\Psi})_{2j}(\lambda e) = T^{j(\Delta-\Delta')} \sum_{i=0}^\infty (-T^{n(1-\Delta')-\Delta'})^{i} \lambda e'_i,$$
        for every $\lambda \in \Lambda_{>0}$, and
        $$H(\widehat{\Psi})_{2j}(\lambda e_i) = T^{\hbigl(j+i(n+1)\hbigr)(\Delta-\Delta')+\Delta-n(1-\Delta)} \lambda e'_i$$
        for every $i \geq 0$ and $\lambda \in \Lambda_{(0, n(1-\Delta)]}$.
    \end{enumerate}  
\end{prop}

\begin{proof}
Consider the subcomplexes
$$C_j = \bigoplus_{i \geq 0} \left( \Lambda_{>0} \check{x}_{j+i(n+1)} \oplus \Lambda_{>0} \hat{x}_{j+i(n+1)+1} \right), \quad C'_j = \bigoplus_{i \geq 0} \left( \Lambda_{>0} \check{y}_{j+i(n+1)} \oplus \Lambda_{>0} \hat{y}_{j+i(n+1)+1} \right),$$
of $(\tilde{C}, \tilde{d})$ and $(\tilde{C}', \tilde{d}')$, respectively. Denote by $d, d'$ the restrictions of $\tilde{d}, \tilde{d}'$ to $C_j, C'_j$. We observe that $(C_j, d)$ and $(C'_j, d')$ are concentrated in degrees $2j-1$ and $2j$ modulo $2(n+1)$. We identify these components with copies of $\Lambda_{>0}^\infty$ via the maps:
$$\lambda \hat{x}_{j+i(n+1)+1} \mapsto \lambda e_{i+1}, \quad \lambda \check{x}_{j+i(n+1)} \mapsto \lambda e_{i+1}, \quad \lambda \hat{y}_{j+i(n+1)+1} \mapsto \lambda e'_{i+1}, \quad \lambda \check{y}_{j+i(n+1)} \mapsto \lambda e'_{i+1},$$
for every $i \geq 0$ and $\lambda \in \Lambda_{>0}$. Under these identifications, the maps $\tilde{\Psi}_{2j-1}$ and $\tilde{\Psi}_{2j}$ take the form:
$$\tilde{\Psi}_{2j-1}(\lambda e_i) = T^{\hbigl(j+(i-1)(n+1)+1\hbigr)(\Delta-\Delta')} \lambda e_i, \quad \tilde{\Psi}_{2j}(\lambda e_i) = T^{\hbigl(j+(i-1)(n+1)\hbigr)(\Delta-\Delta')} \lambda e_i,$$
for every $i \in \N$.

Let $\alpha = \Delta, \beta = n(1-\Delta), \alpha' = \Delta'$, and $\beta' = n(1-\Delta')$. The following diagram, obtained from Diagram~\eqref{diag: alg res maps}, commutes:
$$ \xymatrix@R=2pc@C=3pc{
    \widehat{C_j^{2j-1}} \ar[r]^{\widehat{d}} \ar[d]_{\widehat{\tilde{\Psi}}_{2j-1}} & \widehat{C_j^{2j}} \ar[d]_{\widehat{\tilde{\Psi}}_{2j}} \ar[r] & SH_{\CP^n}^{2j}(B;\Lambda_{\geq0}) \ar[d]_{H(\widehat{\Psi})_{2j}} \\
    \widehat{(C'_j)^{2j-1}} \ar[r]^{\widehat{d'}} & \widehat{(C'_j)^{2j}} \ar[r] & SH_{\CP^n}^{2j}(B';\Lambda_{\geq0})
  }$$
Here we identify $\Psi_1 = \widehat{\tilde{\Psi}}_{2j-1}$, $\Psi_2 = \widehat{\tilde{\Psi}}_{2j}$, $D = \widehat{d}$, $D' = \widehat{d'}$, and $\overline{\Psi} = H(\widehat{\Psi})_{2j}$. Since 
$$\beta' - \beta = n(1-\Delta') - n(1-\Delta) = n(\Delta - \Delta') = n(\alpha - \alpha'),$$
and noting that $\alpha < \beta$ if and only if $\Delta < \frac{n}{n+1}$, the conditions of Proposition~\ref{prop: alg restriction maps} are satisfied. Shifting the index $i$ by $1$ in Proposition~\ref{prop: alg restriction maps} yields the required formulas for $H(\widehat{\Psi})_{2j}$, completing the proof.
\end{proof}

The last step before we present the proof of Theorem~\ref{thm: res for balls in CP^n} is the following claim:

\begin{claim}\label{claim: surjectivity and completion in CF}
    Let $c \fc \drctlim CF(H_\ell) \to \widehat{\drctlim} CF(H_\ell)$ be the canonical completion map. Then the induced map $H(c) \fc H(\drctlim CF(H_\ell)) \to H(\widehat{\drctlim} CF(H_\ell))$ is surjective.
\end{claim}

\begin{proof}
    The direct limit $\drctlim CF(H_\ell)$ is isomorphic to $(\tilde{C}, \tilde{d})$. Let $0 \leq j \leq n$. As in the proof of Proposition~\ref{prop: drctlim restriction maps}, consider the subcomplex
    $$C_j = \bigoplus_{i \geq 0} \left( \Lambda_{>0} \check{x}_{j+i(n+1)} \oplus \Lambda_{>0} \hat{x}_{j+i(n+1)+1} \right)$$
    of $(\tilde{C}, \tilde{d})$, and denote by $d$ the restriction of $\tilde{d}$ to $C_j$. We observe that $(C_j, d)$ is concentrated in degrees $2j-1$ and $2j$ modulo $2(n+1)$. We identify these components with copies of $\Lambda_{>0}^\infty$ via the maps:
    $$\lambda \hat{x}_{j+i(n+1)+1} \mapsto \lambda e_{i+1}, \qquad \lambda \check{x}_{j+i(n+1)} \mapsto \lambda e_{i+1},$$
    for every $i \geq 0$ and $\lambda \in \Lambda_{>0}$. 
    
    By denoting $D = \widehat{d}$, we deduce from Proposition~\ref{prop: surjectivity and completion} that the induced map $H(c)_j \fc H(C_j) \to H(\widehat{C_j})$ is surjective. Consequently, the map 
    $$H(c) \fc H(\drctlim CF(H_\ell)) \to H(\widehat{\drctlim} CF(H_\ell))$$ 
    is surjective, as required.
\end{proof}
\begin{proof}[Proof of Theorem~\ref{thm: res for balls in CP^n}]
By Theorem~\ref{thm: continuation maps for restrictions}, we have the commutativity $\Phi'_\ell \circ \Psi_\ell = \Psi_{\ell+1} \circ \Phi_\ell$, where the continuation maps $\Phi_\ell, \Phi'_\ell, \Psi_\ell$, and $\Psi_{\ell+1}$ appear in the following diagram:

\begin{equation*}
\xymatrix@R=2pc@C=3pc{
    CF(H_\ell) \ar[r]^{\Phi_\ell} \ar[d]^{\Psi_\ell} & CF(H_{\ell+1}) \ar[d]^{\Psi_{\ell+1}} \\ 
    CF(H'_\ell) \ar[r]^{\Phi'_\ell} & CF(H'_{\ell+1})
}
\end{equation*}
for every $\ell \in \Z_{\geq 0}$. According to Claim~\ref{claim: surjectivity and completion in CF}, the canonical map 
$$H(\drctlim CF(H_\ell)) \to H(\widehat{\drctlim} CF(H_\ell))$$
is surjective. Thus, by Proposition~\ref{prop: algebraic preparation for computing restriction maps}, the following diagram commutes:

\begin{equation*}
\xymatrix@R=2pc@C=3pc{
    H(\widehat{\tel}\, CF(H_\ell)) \ar[r] \ar[d]^{\cong} & H(\widehat{\tel}\, CF(H'_\ell)) \ar[d]^{\cong} \\ 
    H(\widehat{\drctlim} CF(H_\ell)) \ar[r] & H(\widehat{\drctlim} CF(H'_\ell))
}
\end{equation*}

Furthermore, the vertical arrows in this diagram are isomorphisms. The top horizontal arrow represents the restriction map $\res \fc SH_{\CP^n}(B;\Lambda_{\geq 0}) \to SH_{\CP^n}(B';\Lambda_{\geq 0})$. Consequently, the commutativity of the diagram and the fact that the vertical arrows are isomorphisms allow us to identify the restriction map with the bottom horizontal arrow. The explicit formula for this bottom arrow is provided in Proposition~\ref{prop: drctlim restriction maps}, which completes the proof.
\end{proof}

\subsection{Proof of Proposition~\ref{prop: res from CPn to a ball}}\label{ss: res maps from CPn to a ball}

The proof of Proposition~\ref{prop: res from CPn to a ball} is similar to and easier than that of Theorem~\ref{thm: res for balls in CP^n}. Let $0<\Delta<1$. Consider the corresponding acceleration data $(H_\ell)_{\ell\geq0}$ for the ball $B^{}_{\Delta}=\mu^{-1}([0,\Delta])$, given by $H_\ell(z)=h(\Delta,\ell,\mu(z))$ for every $z\in \CP^n$ and $\ell\in \Z_{\geq0}$. For every $\ell\in \Z_{\geq 0}$, we have $H_0\leq H_\ell$ and $H_0\leq 0$. Thus, since 
$$\lim\limits_{\ell\to \infty}H_\ell(z)=\left\{\begin{array}{ll}
   0, & z\in B,\\
    +\infty, &  z\notin B,
\end{array}\right.$$ 
for every $z\in \CP^n$, set $\varepsilon_0=1$, and there exists a non-increasing sequence $(\varepsilon_\ell)_{\ell\geq1}$ such that $\lim\limits_{\ell\to \infty}\varepsilon_\ell=0$ and for every $\ell\in \Z_{\geq1}$ we have $\varepsilon_\ell\cdot H_0\leq H_\ell$. For every $\ell\in \Z_{\geq1}$, we denote $H'_\ell=\varepsilon_\ell\cdot H_0$. Note that $(H'_\ell)_{\ell\geq0}$ is an acceleration datum for $\CP^n$ consisting of $C^2$-small Morse--Bott functions. For every $\ell\in \Z_{\geq0}$, the Floer complex $CF(H'_\ell)$ of $H'_\ell$ has $n+1$ generators, denoted $\check{z}_0^\ell,\ldots,\check{z}_n^\ell$. The differential is $0$, and the continuation map $\Phi'_\ell\fc CF(H'_\ell)\to CF(H'_{\ell+1})$ satisfies 
\begin{equation}\label{eq: cont for H of CPn}
    \Phi'_\ell(\check{z}_i^\ell)=T^{H'_{\ell+1}(\check{z}_i^{\ell+1})-H'_{\ell}(\check{z}_i^{\ell})}\check{z}_i^{\ell+1}
\end{equation} for every $0\leq i\leq n$, see \cite[Page 599]{Varolgunes_2021_MV_and_relSH} for the Morse case.

For every $\ell\in \Z_{\geq0}$ denote by $\check{x}_0^\ell,\ldots,\check{x}_{\ell+n}^\ell,\hat{x}_1^\ell,\ldots,\hat{x}_\ell^\ell$ the generators for $CF(H_\ell)$. By Proposition~\ref{prop: omission of signs for CF}, for every $\ell\in \Z_{\geq0}$, the bases $\check{x}_0^\ell,\ldots,\check{x}_{\ell+n}^\ell,\hat{x}_1^\ell,\ldots,\hat{x}_\ell^\ell $ for $CF(H_\ell)$ can be chosen such that the differentials and continuation maps in the $1$-ray
$$CF(H_0)\xrightarrow{\Phi_0}CF(H_1)\xrightarrow{\Phi_1}CF(H_2)\to\cdots$$

are represented by matrices with entries in $\{0, 1\}$, where we ignore the powers of the formal variable $T$ of the Novikov ring. We  refer to such bases as \textbf{normalized bases}.

Let $\chi\fc \R\to \R$ be a smooth non-decreasing function, satisfying $\chi(s)=0$ for every $s<0$ and $\chi(s)=1$ for every $s>1$.  Then for every $\ell\in \Z_{\geq0}$ one can define a monotone homotopy $K_\ell\fc \R\times \CP^n\to \R$ from $H'_\ell$ to $H_\ell$ by 
\[K_\ell(s,z)
=(1-\chi(s))\varepsilon_\ell H_0(z)+\chi(s)H_\ell(z)\]
for every $(s,z)\in \R\times \CP^n$. Thus for every $\ell\in \Z_{\geq0}$ the homotopy $K_\ell$ defines a continuation map $\Psi_\ell\fc CF(H'_\ell)\to CF(H_\ell)$.
Similarly to Theorem~\ref{thm: CF(H_l) + continuation maps}, to describe the continuation maps with coefficients in the Novikov ring, in the sense of Section~\ref{sss: weighted_CF}, we must compute the energy of the Floer and continuation flowlines with cascades that appear in the differential and continuation maps.

\begin{prop}\label{prop: continuation maps for restrictions for CP^n->B}
    For every $\ell\in \Z_{\geq0}$ there exist signs
    $$\check{D}_0,\ldots,\check{D}_{n}\in \{-1,1\}$$
    such that the continuation map $\Psi_\ell\fc CF(H'_\ell)\to CF(H_\ell)$ satisfies
    \begin{itemize}
    \item $\Phi_\ell\circ \Psi_\ell=\Psi_{\ell+1}\circ\Phi'_\ell$ where $\Phi_\ell\fc CF(H_\ell)\to CF(H_{\ell+1})$ and $\Phi'_\ell\fc CF(H'_\ell)\to CF(H'_{\ell+1})$ are the continuation maps.
   
        \item For every $0\leq i\leq n$ we have
        $$\Psi(\check{z}_i^\ell)=\check{D}_i T^{E(\check{z}_i^\ell,\check{x}_i^\ell)}\check{x}_i^\ell,$$
        where
        $$E(\check{z}_i^\ell,\check{x}_i^\ell)=\left\{\begin{array}{ll}
        h_{\ell,0}^\Delta-\varepsilon_\ell h_{0,0}^\Delta,& i=0,\\
            h_{\ell,i}^\Delta-\varepsilon_\ell\cdot h(\Delta,0,1)+i(1-r_{\ell,i}^\Delta), & 1\leq i\leq \ell, \\
            h(\Delta,\ell,1)-\varepsilon_\ell\cdot h(\Delta,0,1), & i\geq\ell+1, 
        \end{array}\right.$$

    \end{itemize}
\end{prop}
\begin{proof}
   Let $\ell \in \Z_{\geq 0}$. Since $H'_\ell$ and $H_0$ obtained one from the other just by multiplication of a small constant, the fact that $\Phi_\ell\circ \Psi_\ell=\Psi_{\ell+1}\circ\Phi'_\ell$ is an immediate corollary of Proposition~\ref{prop: commuting of cont maps}. From the same reason, Theorem~\ref{thm: continuation maps, over Z} implies that there exist $\check{D}_0,\ldots,\check{D}_n\in \{-1,1\}$ and $E(\check{z}_0^\ell,\check{x}_0^\ell),\ldots,E(\check{z}_n^\ell,\check{x}_n^\ell)>0$ such that
 $$\Psi(\check{z}_i^\ell)=\check{D}_i T^{E(\check{z}_i^\ell,\check{x}_i^\ell)}\check{x}_i^\ell,$$
 for every $0\leq i\leq n$.

   The rest of the proof is similar to that of Theorem~\ref{thm: CF(H_l) + continuation maps}. Since we are interested in continuation flowlines with cascades, the dimension $\dim_u \wh\cM_m(p, q)$ equals $0$ for every $p, q$, and $u$. Thus, the topological energy of such a flowline with cascades is given by
$$ E_{top}(u) = \cA_{H_\ell}(q, \hat{q}) - \cA_{H'_\ell}(p, \hat{p}) + \frac{1}{2(n+1)} \left( \mu_{FMB}^\tau(p) - \mu_{FMB}^\tau(q) \right). $$
As in the proof of Theorem~\ref{thm: CF(H_l) + continuation maps}, for constant orbits of $H_\ell$ and $H'_\ell$, choose the cappings to be constant, thus their action equals the value of the Hamiltonian itself, that is
$$ \cA_{H_\ell}(\check{x}_0^\ell) = h(\Delta, \ell, 0)=h_{\ell,0}^\Delta, \qquad \text{and} \qquad \cA_{H_\ell}(\check{x}_{\ell+1}^\ell) = \cdots = \cA_{H_\ell}(\check{x}_{\ell+n}^\ell) = h(\Delta, \ell, 1), $$
and 
$$ \cA_{H'_\ell}(\check{z}_0^\ell) = \varepsilon_\ell h(\Delta, 0, 0)=\varepsilon_\ell h_{0,0}^{\Delta}, \qquad \text{and} \qquad \cA_{H'_\ell}(\check{z}_{1}^\ell) = \cdots = \cA_{H'_\ell}(\check{z}_{n}^\ell) = \varepsilon_\ell\cdot h(\Delta, 0, 1). $$
For the orbits $\hat{x}_1^\ell, \check{x}_1^\ell, \ldots, \hat{x}_\ell^\ell, \check{x}_\ell^\ell$, we choose cappings contained in $\CP^n \setminus D_\infty$, which is an exact symplectic manifold symplectomorphic to $\Int B(1)$. For every $1 \leq i \leq \ell$, we find that
$$ \cA_{H_\ell}(\hat{x}_i^\ell) = \cA_{H_\ell}(\check{x}_i^\ell)  = h_{\ell,i}^\Delta - i r_{\ell,i}^\Delta. $$
Additionally, in Section~\ref{ss: computations of FMB}, we found that the Floer--Morse--Bott indices for these $1$-periodic orbits of $H_\ell$ with respect to the trivialization $\tau$ induced by these cappings satisfy:
\begin{itemize}
    \item For every $0\leq i\leq n$ we have $\mu_{FMB}^\tau(\check{z}_{i}^\ell)=2i$;
    \item $\mu_{FMB}^\tau(\check{x}_{0}^\ell)=0$;
    \item For every $1\leq i\leq \ell$ we have $\mu_{FMB}^\tau(\check{x}_{i}^\ell)=-2ni$;
    \item For every $\ell+1\leq i\leq \ell+n$ we have $\mu_{FMB}^\tau(\check{x}_{i}^\ell)=2i$.
\end{itemize}

Now, let $0 \leq i \leq n$.
\begin{itemize}
\item If $i=0$, then $E(\check{z}_0^\ell, \check{x}_0^{\ell}) = \cA_{H_{\ell}}(\check{x}_0^{\ell}) - \cA_{H'_\ell}(\check{z}_0^\ell) = h_{\ell,0}^\Delta-\varepsilon_\ell h_{0,0}^\Delta.$

    \item Assume that $1\leq i\leq n$.
    \begin{itemize}
        \item If $i\leq \ell$ then
        $$\mu_{FMB}^\tau(\check{z}_i^\ell)-\mu_{FMB}^\tau(\check{x}_i^\ell)=2i-(-2in)=2i(n+1),$$
        and hence
        \begin{align*}
             E(\check{z}_i^\ell, \check{x}_i^{\ell}) &= \cA_{H_{\ell}}(\check{x}_i^{\ell}) - \cA_{H'_\ell}(\check{z}_i^\ell)  + \frac{1}{2(n+1)} \left( \mu_{FMB}^\tau(\check{z}_i^\ell)-\mu_{FMB}^\tau(\check{x}_i^\ell)\right)\\
             &=h_{\ell,i}^\Delta-ir_{\ell,i}^\Delta-\varepsilon_\ell\cdot h(\Delta,0,1)+\frac{2i(n+1)}{2(n+1)}\\
             &=h_{\ell,i}^\Delta-\varepsilon_\ell\cdot h(\Delta,0,1)+i(1-r_{\ell,i}^\Delta).
        \end{align*}

        \item If $i\geq\ell+1$, then $\mu_{FMB}^\tau(\check{z}_i^\ell)=\mu_{FMB}^\tau(\check{x}_i^\ell)=2i$, and hence
        \begin{align*}
             E(\check{z}_i^\ell, \check{x}_i^{\ell}) &= \cA_{H_{\ell}}(\check{x}_i^{\ell}) - \cA_{H'_\ell}(\check{z}_i^\ell) \\
             &=h(\Delta,\ell,1)-\varepsilon_\ell\cdot h(\Delta,0,1).
        \end{align*}

    \end{itemize}
\end{itemize}
This completes the proof.

\end{proof}

The following result is an analogue of Proposition~\ref{prop: omission of signs for CF->CF'} and demonstrates that the signs from Theorem~\ref{prop: continuation maps for restrictions for CP^n->B} can be normalized to $1$ via a suitable change of bases.

\begin{prop}\label{prop: omission of signs for CF'->CF}
     For every $\ell\in \Z_{\geq0}$, there is a normalized basis, $\check{x}_0^\ell,\ldots,\check{x}_{\ell+n}^\ell,\hat{x}_1^\ell,\ldots,\hat{x}_\ell^\ell $ 
for $CF(H_\ell)$ such that the continuation map $\Psi_\ell\fc CF(H'_\ell)\to CF(H_\ell)$ satisfies
        $$\Psi(\check{z}_i^\ell)=T^{E(\check{z}_i^\ell,\check{x}_i^\ell)}\check{x}_i^\ell,$$
        where
        $$E(\check{z}_i^\ell,\check{x}_i^\ell)=\left\{\begin{array}{ll}
        h_{\ell,0}^\Delta-\varepsilon_\ell h_{0,0}^\Delta,& i=0,\\
            h_{\ell,i}^\Delta-\varepsilon_\ell\cdot h(\Delta,0,1)+i(1-r_{\ell,i}^\Delta), & 1\leq i\leq \ell, \\
            h(\Delta,\ell,1)-\varepsilon_\ell\cdot h(\Delta,0,1), & i\geq\ell+1, 
        \end{array}\right.,$$
       for every $0\leq i\leq n$.
\end{prop}

\begin{proof} Let $\ell\in \Z_{\geq0}$. The desired result is an immediate corollary of Proposition~\ref{prop: omission of signs for CF} applied to $CF(H_\ell)$, combined with the fact that the differential on $CF(H'_\ell)$ is zero and that $CF(H'_\ell)$ has at most one generator in each degree. Thus, we are able to change the signs of the generators of $CF(H'_\ell)$ to deduce the desired result.
\end{proof}

A straitforward use of the combination of Lemma~\ref{lemma: drct lim of 1-ray} and Equation~\eqref{eq: cont for H of CPn} implies the direct limit of the $1$-ray
$$CF(H'_0) \to CF(H'_1) \to CF(H'_2) \to \cdots$$
is isomorphic to the cochain complex $(C',d')$ given by $$C'=\bigoplus_{j=0}^n\Lambda_{>0}\check{z}_j$$ where the differential $d'$ is zero, and for every $\ell\in \Z_{\geq0}$ the canonical map $f'_\ell\fc CF(H'_\ell)\to C'$ satisfies $f'_\ell(\check{z}_j^\ell)=T^{-H'_\ell(\check{z}_j^\ell)}\check{z}_j$ for every $0\leq j\leq n$.

Consider the cochain complex $(C, d)$ defined as follows: the module $C$ is given by
$$C = \Lambda_{>0}\check{x}_0 \oplus \bigoplus_{j\in \N} \left( \Lambda_{>j\Delta}\check{x}_j \oplus \Lambda_{>j\Delta}\hat{x}_j \right),$$
and the differential $d$ satisfies $d(\lambda\check{x}_j) = 0$ for every $j \in \Z_{\geq0}$ and $\lambda \in \Lambda_{>j\Delta}$, while $d(\lambda\hat{x}_j) = \lambda\check{x}_{j-1} + T^n\lambda\check{x}_{j+n}$ for every $j \in \N$ and $\lambda \in \Lambda_{>j\Delta}$. By Theorem~\ref{thm: drct lim of CF(H_l)}, the direct limits of the $1$-ray
$$CF(H_0) \to CF(H_1) \to CF(H_2) \to \cdots$$
is isomorphic to $(C, d)$.

Let $\Psi \fc C \to C'$ denote the direct limit of the continuation maps $(\Psi_\ell \fc CF(H'_\ell) \to CF(H_\ell))_{\ell\geq 0}$. Note that $\Psi$ is well-defined because Proposition~\ref{prop: continuation maps for restrictions for CP^n->B} asserts that $\Psi_{\ell+1} \circ \Phi'_\ell = \Phi_\ell \circ \Psi_\ell$, which ensures that the maps $(\Psi_\ell)_{\ell\geq 0}$ define a morphism of directed systems. The following proposition describes this map.

\begin{prop}
    The map 
    $$\Psi \fc \bigoplus_{j=0}^n\Lambda_{>0}\check{z}_j\to\Lambda_{>0}\check{x}_0 \oplus \bigoplus_{j\in \mathbb{N}} \left( \Lambda_{>j\Delta}\check{x}_j \oplus \Lambda_{>j\Delta}\hat{x}_j \right)$$
    satisfies $\Psi(\lambda \check{z}_j) = T^j \lambda \check{x}_j$ for every $0\leq j\leq n$ and $\lambda\in \Lambda_{>0}$.
\end{prop}

\begin{proof}
By the universal property of the direct limit, for every $\ell \in \Z_{\geq 0}$, the following diagram commutes:
$$
\xymatrix@R=2pc@C=3pc{
CF(H'_\ell) \ar[r]^{f'_\ell} \ar[d]_{\Psi_\ell} & C' \ar[d]^{\Psi} \\
CF(H_\ell) \ar[r]_{f_\ell} & C
}
$$
where $f_\ell$ and $f'_\ell$ are the canonical maps of the direct limits $C$ and $C'$, respectively. 

Let $0\leq j\leq n$ and $\lambda \in \Lambda_{>0}$. From the properties of direct limits, for any $\ell \geq j$, the generator $\lambda \check{z}_j \in C'$ can be represented as:
$$\lambda \check{z}_j = f'_\ell \left( \lambda T^{ H'_\ell(\check{z}_j^\ell)} \check{x}_j^\ell \right).$$
Using the commutativity of the diagram, we have:
$$\Psi(\lambda \check{z}_j) = \Psi \circ f'_\ell \left( \lambda T^{ H'_\ell(\check{z}_j^\ell)} \check{z}_j^\ell \right) = f_\ell \circ \Psi_\ell \left( \lambda T^{ H'_\ell(\check{z}_j^\ell)} \check{z}_j^\ell \right).$$
By Proposition~\ref{prop: omission of signs for CF'->CF}, we have $\Psi_\ell(\check{z}_j^\ell) = T^{E(\check{z}_j^\ell, \check{x}_j^\ell)} \check{x}_j^\ell$. Thus:
$$\Psi(\lambda \check{z}_j) = \lambda T^{H'_\ell(\check{z}_j^\ell)} f_\ell \left( T^{E(\check{z}_j^\ell, \check{x}_j^\ell)} \check{x}_j^\ell \right).$$
Applying the formula for $f_\ell \fc CF(H_\ell) \to C$ from Theorem~\ref{thm: drct lim of CF(H_l)}, we obtain:
$$\Psi(\lambda \check{z}_j) = \lambda T^{H'_\ell(\check{z}_j^\ell) + E(\check{z}_j^\ell, \check{x}_j^\ell) - (h_{\ell,j}^{\Delta} - j r_{\ell,j}^{\Delta})} \check{x}_j.$$
Let us show that 
$$H'_\ell(\check{z}_j^\ell) + E(\check{z}_j^\ell, \check{x}_j^\ell) - (h_{\ell,j}^{\Delta} - j r_{\ell,j}^{\Delta})=j.$$
Indeed, if $j=0$, then $H'_\ell(\check{z}_0^\ell)=\varepsilon_\ell\cdot  h_{0,0}^\Delta$ and $E(\check{z}_0^\ell, \check{x}_0^\ell)=h_{\ell,0}^\Delta-\varepsilon_\ell h_{0,0}^\Delta$ which imply that
    $$H'_\ell(\check{z}_0^\ell) + E(\check{z}_0^\ell, \check{x}_0^\ell) - (h_{\ell,0}^{\Delta} - 0\cdot r_{\ell,0}^{\Delta})=\varepsilon_\ell h_{0,0}^\Delta + h_{\ell,0}^\Delta-\varepsilon_\ell h_{0,0}^\Delta - h_{\ell,0}^{\Delta}=0.$$
   Otherwise, $j\geq 1$, and hence $H'_\ell(\check{z}_0^\ell)=\varepsilon_\ell\cdot  h(\Delta,0,1)$ and $E(\check{z}_j^\ell, \check{x}_j^\ell)=h_{\ell,j}^\Delta-\varepsilon_\ell\cdot h(\Delta,0,1)+j(1-r_{\ell,j}^\Delta)$ which imply that
\begin{align*}
    H'_\ell(\check{z}_j^\ell) + E(\check{z}_j^\ell, \check{x}_j^\ell) - (h_{\ell,j}^{\Delta} - j r_{\ell,j}^{\Delta})&=\varepsilon_\ell \cdot h(\Delta,0,1)+h_{\ell,j}^\Delta-\varepsilon_\ell\cdot h(\Delta,0,1)\\
    &+j(1-r_{\ell,j}^\Delta)- (h_{\ell,j}^{\Delta} - j r_{\ell,j}^{\Delta})\\
    &=j.
\end{align*}
 Therefore, $\Psi(\lambda \check{z}_0) = T^j \lambda \check{x}_0$. This completes the proof.
\end{proof}

As in the proof of Theorem~\ref{thm: drct lim of CF(H_l)}, define cochain complex $(\tilde{C},\tilde{d})$  over $\Lambda_{\geq 0}$ as follows:
$$\tilde{C}=\Lambda_{>0}\check{x}_0\oplus\bigoplus_{i=1}^\infty \left(\Lambda_{>0} \check{x}_i\oplus \Lambda_{>0} \hat{x}_i\right),$$
where for every $i\geq 0$ and $\lambda\in\Lambda_{>0}$ we have $\tilde{d}(\lambda \check{x}_i)=0$, also, for every $i\geq 1$ and $\lambda\in\Lambda_{>0}$, the differential satisfies
$$\tilde{d}(\lambda \hat{x}_i)=T^\Delta \lambda \check{x}_{i-1}+T^{n(1-\Delta)} \lambda \check{x}_{i+n}.$$
By Claim~\ref{claim: new presentation of drct lim of CF}, the cochain complex $(\tilde{C},\tilde{d})$ is isomorphic to $(C,d)$ via the isomorphism $K\fc \tilde{C}\to C$ given by $K=\id\bigoplus_{i=1}^\infty T^{i\Delta}(\id\oplus \id)$. 
Define the cochain map $\tilde{\Psi} \fc C' \to \tilde{C}$ by $\tilde{\Psi} = (K)^{-1} \circ \Psi$. A straightforward computation shows that $\tilde{\Psi}$ satisfies 
$$\tilde{\Psi}(\lambda \check{z}_i) = T^{j(1 - \Delta)} \lambda \check{z}_i$$ 
for every $0\leq i\leq n$ and $\lambda \in \Lambda_{> 0}$.

As seen in the proof of Theorem~\ref{thm: SH of ball in CP^n}, and as it is mentioned in Section~\ref{ss: relSH using MB}, we can compute the relative symplectic cohomology using an acceleration datum of Hamiltonians that satisfy the \textbf{MB} condition, therefore we have $SH_{\CP^n}(B;\Lambda_{\geq0}) \cong H(\widehat{\tilde{C}})$ and $SH_{\CP^n}(\CP^n;\Lambda_{\geq0}) \cong H(\widehat{C'})$. Our next objective is to compute the induced map $(\widehat{\Psi})_* \fc SH_{\CP^n}(\CP^n;\Lambda_{\geq0}) \to SH_{\CP^n}(B;\Lambda_{\geq0})$.

\begin{prop}\label{prop: drctlim restriction maps, CPn - ball}
    Let $0 \leq j \leq n$. The morphism 
    $$H(\widehat{\Psi})_{2j} \fc SH_{\CP^n}^{2j}(\CP^n;\Lambda_{\geq0}) \to SH_{\CP^n}^{2j}(B;\Lambda_{\geq0})$$
    is given by following: 
    \begin{enumerate}

        \item If $\Delta < \frac{n}{n+1} $, then the map 
        $$H(\widehat{\Psi})_{2j} \fc \Lambda_{>0} e  \to \bigoplus_{i \geq 0} \Lambda_{(0, \Delta]} e_i$$
        satisfies 
        $$H(\widehat{\Psi})_{2j}(\lambda e) = T^{j(1-\Delta)} \lambda \cdot\sum_{i=0}^\infty(-T^{n(1-\Delta)-\Delta})^i e_i,$$
        for every $\lambda \in \Lambda_{>0}$.

          \item If $ \Delta\geq \frac{n}{n+1}$, then the map 
        $$H(\widehat{\Psi})_{2j} \fc \Lambda_{>0} e\to \Lambda_{>0} e_0 \oplus \bigoplus_{i \geq 1} \Lambda_{(0, n(1-\Delta)]} e_i$$
        satisfies 
        $$H(\widehat{\Psi})_{2j}(\lambda e) = T^{j(1-\Delta)} \lambda e_0,$$
        for every $\lambda \in \Lambda_{>0}$.

    \end{enumerate}  
\end{prop}

\begin{proof}
Consider the subcomplexes
$$C_j = \bigoplus_{i \geq 0} \left( \Lambda_{>0} \check{x}_{j+i(n+1)} \oplus \Lambda_{>0} \hat{x}_{j+i(n+1)+1} \right), \quad C'_j = \Lambda_{>0} \check{z}_j,$$
of $(\tilde{C}, \tilde{d})$ and $(C', d')$, respectively. Denote by $d$ the restriction of $\tilde{d}$ to $C_j$. We observe that $(C_j, d)$ is concentrated in degrees $2j-1$ and $2j$ modulo $2(n+1)$. We identify these components with copies of $\Lambda_{>0}^\infty$ via the maps:
$$\lambda \hat{x}_{j+i(n+1)+1} \mapsto \lambda e_{i+1}, \quad \lambda \check{x}_{j+i(n+1)} \mapsto \lambda e_{i+1},$$
for every $i \geq 0$ and $\lambda \in \Lambda_{>0}$, and replace the notation $\check{z}_j$ by $e'_1$. Under these identifications, the map $\tilde{\Psi}_{2j-1}$ is $0$ and the map $\tilde{\Psi}_{2j}$ satisfies $\tilde{\Psi}_{2j}(\lambda e'_1) = T^{j(1-\Delta)} \lambda e_1$.

Let $\alpha = \Delta, \beta = n(1-\Delta)$. The following diagram, obtained from Diagram~\eqref{diag: alg res maps simplified}, commutes:

$$ \xymatrix@R=2pc@C=3pc{
    0 \ar[r] \ar[d] & \widehat{(C'_j)^{2j}} \ar[d]_{\widehat{\tilde{\Psi}}_{2j}} \ar[r] & SH_{\CP^n}^{2j}(\CP^n;\Lambda_{\geq0}) \ar[d]_{H(\widehat{\Psi})_{2j}} \\
    \widehat{(C_j)^{2j-1}} \ar[r]^{\widehat{d}} & \widehat{(C_j)^{2j}} \ar[r] & SH_{\CP^n}^{2j}(B;\Lambda_{\geq0})
  }$$

and noting that $\alpha < \beta$ if and only if $\Delta < \frac{n}{n+1}$. Thus the conditions of Proposition~\ref{prop: alg restriction maps simplified} are satisfied. Shifting the index $i$ by $1$ in Proposition~\ref{prop: alg restriction maps simplified} yields the required formulas for $H(\widehat{\Psi})_{2j}$, completing the proof.
\end{proof}

\begin{proof}[Proof of Proposition~\ref{prop: res from CPn to a ball}]
By Theorem~\ref{prop: continuation maps for restrictions for CP^n->B}, we have the commutativity $\Phi_\ell \circ \Psi_\ell = \Psi_{\ell+1} \circ \Phi'_\ell$, where the continuation maps $\Phi_\ell, \Phi'_\ell, \Psi_\ell$, and $\Psi_{\ell+1}$ appear in the following diagram:

\begin{equation*}
\xymatrix@R=2pc@C=3pc{
    CF(H'_\ell) \ar[r]^{\Phi'_\ell} \ar[d]^{\Psi_\ell} & CF(H'_{\ell+1}) \ar[d]^{\Psi_{\ell+1}} \\ 
    CF(H_\ell) \ar[r]^{\Phi_\ell} & CF(H_{\ell+1})
}
\end{equation*}
for every $\ell \in \Z_{\geq 0}$. Since $\drctlim CF(H'_\ell)=\Lambda_{>0}^{n+1}$ is a complete module,we deduce that the canonical map 
$$H(\drctlim CF(H'_\ell)) \to H(\widehat{\drctlim} CF(H'_\ell))$$
is an isomorphism (even in the cochain level), and in particular, surjective. Thus, by Proposition~\ref{prop: algebraic preparation for computing restriction maps}, the following diagram commutes:

\begin{equation*}
\xymatrix@R=2pc@C=3pc{
    H(\widehat{\tel}\, CF(H'_\ell)) \ar[r] \ar[d]^{\cong} & H(\widehat{\tel}\, CF(H_\ell)) \ar[d]^{\cong} \\ 
    H(\widehat{\drctlim} CF(H'_\ell)) \ar[r] & H(\widehat{\drctlim} CF(H_\ell))
}
\end{equation*}

Furthermore, the vertical arrows in this diagram are isomorphisms. The top horizontal arrow represents the restriction map $\res \fc SH_{\CP^n}(\CP^n;\Lambda_{\geq 0}) \to SH_{\CP^n}(B;\Lambda_{\geq 0})$. Consequently, the commutativity of the diagram and the fact that the vertical arrows are isomorphisms allow us to identify the restriction map with the bottom horizontal arrow. The explicit formula for this bottom arrow is provided in Proposition~\ref{prop: drctlim restriction maps}, which completes the proof.
\end{proof}

\section{Torsion and stable displacement energy}\label{s: Tor and e_st}
Given a subset $X$ of a symplectic manifold $(M,\omega)$, we denote by $e(X)$ its \textbf{displacement energy}. The \textbf{stable displacement energy of $X$}, denoted $e_{st}(X)$, is defined to be the displacement energy of the stabilization of $X$, that is, of $X \times S^1$ as a subset of $(M \times T^* S^1, \omega \oplus \omega_0)$, so that $e_{st}(X) = e(X \times S^1)$. Trivially, for every subset $X$ we have $e_{st}(X) \leq e(X)$.

Given a module $A$ over the Novikov ring $\Lambda_{\geq0}$, its capacity is defined to be
$$c(A)=\inf\{r>0\,|\,T^r A=0\}.$$
Note that $c$ takes values in $[0,+\infty]$. The module $A$ is called \textbf{torsion} if $A\otimes \Lambda = 0$, in particular, this holds when $c(A)<+\infty$.

In his thesis \cite[Remark 4.2.8]{Varolgunes_2018_PhD}, Varolg\"une\c s proved that given a closed symplectic manifold $(M,\omega)$ and a compact subset $K\subset M$, the displacement energy of $K$ is bounded from below by the capacity of $SH_M^*(K;\Lambda_{\geq0})$, that is, $$c(SH_M^*(K;\Lambda_{\geq0}))\leq e(K).$$
In particular, if $K$ is displaceable then $SH_M^*(K;\Lambda_{\geq0})$ is torsion, and if $c(SH_M^*(K;\Lambda_{\geq0}))=+\infty$ then $K$ is non-displaceable\footnote{Due to \cite{MSV_2024_heavy_sets_and_SH}, it is possible to prove that $c(SH_M^*(K;\Lambda_{\geq0}))=+\infty$ if $K$ is a heavy subset in the sense of \cite{EP_2009_rigid_subsets}}. The main result of this section is the following generalization of Varolg\"une\c s's result, denoted as Theorem \ref{thm: c_relSH vs e_st} in the introduction:

\begin{thm}
    Let $(M,\omega)$ be a closed symplectic manifold and $K\subset M$ be a compact subset. Then the stable displacement energy of $K$ is bounded from below by the capacity of $SH_M^*(K;\Lambda_{\geq0})$, that is, $c(SH_M^*(K;\Lambda_{\geq0}))\leq e_{st}(K)$.
\end{thm}
\begin{proof}[Proof of Theorem~\ref{thm: c_relSH vs e_st}]
 If $K$ is stably non-displaceable, its stable displacement energy is infinity, and the statement is trivially true. Assume that $K$ is stably displaceable; this implies that $K \times S^1$ is displaced by a Hamiltonian diffeomorphism of $M \times T^* S^1$. Since this diffeomorphism is Hamiltonian, we can replace it with a compactly supported Hamiltonian diffeomorphism $(\varphi_t)_{t \in [0,1]}$ that displaces $K \times S^1$, that is, $\varphi_0 = \id$, $\varphi_1(K \times S^1) \cap (K \times S^1) = \varnothing$ and $\supp \varphi = \overline{\bigcup_{t \in [0,1]} \supp \varphi_t}$ is contained in a compact subset of $M \times T^* S^1$. 

Thus, we can symplectically embed an open neighborhood of the support of $\varphi$ into $M \times \T^{2}$, where $\T^{2}$ is the standard torus $(S^1\times S^1, A dx\wedge dy)$ with large enough area $A$, such that the projection of $K\times S^1\subset M\times \T^2$ to $\T^2=S^1\times S^1$ has the form $S^1 \times \{0\}$ for some constant point $0\in S^1=\R/\Z$.

 We will prove that a copy of $SH_M^*(K;\Lambda_{\geq0})\otimes \Lambda_{>0}$ appears as a submodule of $SH_{M\times \T^{2}}(K\times S^1;\Lambda_{\geq0})$. Thus, from Varolg\"une\c s's result \cite[Remark 4.2.8]{Varolgunes_2018_PhD} and the fact that the capacity of modules over the Novikov ring is monotone with respect to inclusions, it follows that
\begin{align*}
    c(SH_{M}(K;\Lambda_{\geq0}))&=c(SH_{M}(K;\Lambda_{\geq0})\otimes \Lambda_{>0}) \\&\leq c(SH_{M\times \T^{2}}(K\times S^1;\Lambda_{\geq0})) \\&\leq e(K\times S^1) \\&\le e_{st}(K),
\end{align*}
and this completes the proof.

Let $(H_\ell)_{\ell\in \N}$ be a nondegenerate acceleration datum for $K\subset M$ and let $(J_\ell)_{\ell\in \N}$ be almost complex structures on $M$ such that $(H_\ell,J_\ell)_{\ell\in \N}$ is regular for every $\ell\in \N$. Also, for every $\ell\in \N$ define the function $\tilde{h}_\ell\fc \T^2\to \R$ by $\tilde{h}_\ell(x,y)=\ell(1-\cos(2\pi y))-\frac{1}{\ell}$ for every $(x,y)\in \T^2=\R^2/\Z^2$. For every $\ell\in \N$, the Hamiltonian $\tilde{h}_\ell$ has two $S^1$-families of $1$-periodic orbits, located at $S^1\times\{0\}$ and $S^1\times \{\frac{1}{2}\}$. Also, note that for every $\ell\in \N$ we have $\tilde{h}_\ell < \tilde{h}_{\ell+1}$. 

For every $\ell\in \N$, let $h_\ell$ be a time-dependent nondegenerate Hamiltonian on $\T^2$ obtained from $\tilde{h}_\ell$ by a small perturbation, as in \cite{CFHW_1996_ApSH_II}. Thus $h_\ell$ has exactly $4$ generators: $\check{x}_\ell$ and $\hat{x}_\ell$ correspond to the minimum and the maximum, respectively, of a Morse function on $S^1\times\{0\}$, and $\check{y}_\ell$ and $\hat{y}_\ell$ are the minimum and the maximum, respectively, of a Morse function on $S^1\times\{\frac{1}{2}\}$. Moreover, we can make the perturbation small enough such that for every $\ell\in \N$ we have $h_\ell < h_{\ell+1}$ and additionally $h_\ell(x,0) < 0$ for every $x \in S^1$. Therefore the sequence $(h_\ell)_{\ell\in \N}$ is an acceleration datum for $S^1\times\{0\} \subset \T^2$.

For every $\ell\in \N$, choose an almost complex structure $j_\ell$ on $\T^2$ such that $(h_\ell,j_\ell)$ is regular. Note that the Floer cohomology of $h_\ell$ with $\Q$-coefficients is isomorphic to the singular cohomology $H^*(\T^2;\Q)$ of $\T^2$, which is $4$-dimensional. Since the Floer complex of $h_\ell$ is $4$-dimensional, we deduce that the differential of the complex is zero.

Since the first Chern class of $\T^2$ is zero, we deduce that the Robbin--Salamon index of any $1$-periodic orbit is well defined in $\frac{1}{2}\Z$. Since $S^1 \times \{0\}$ is the minimum critical submanifold of $\tilde{h}_\ell$, for every $\ell \in \N$, we deduce that the Robbin--Salamon index of any point on $S^1 \times \{0\}$ with respect to $\tilde{h}_\ell$ and the trivialization induced from the constant capping is $-\frac{1}{2}$, and thus from \cite[Lemma 2.2]{CFHW_1996_ApSH_II} we deduce that the Robbin--Salamon index of $\check{x}_\ell$ is $-1$ and that of $\hat{x}_\ell$ is $0$. Similarly, the Robbin--Salamon index of any point on $S^1 \times \{\frac{1}{2}\}$ with respect to $\tilde{h}_\ell$ and the trivialization induced from the constant capping is $\frac{1}{2}$, and thus from \cite[Lemma 2.2]{CFHW_1996_ApSH_II} we deduce that the Robbin--Salamon index of $\check{y}_\ell$ is $0$ and that of $\hat{y}_\ell$ is $1$. In particular, the continuation map $ CF(h_\ell;\Z) \to CF(h_{\ell+1};\Z)$ sends $\check{x}_\ell$ to $\pm\check{x}_{\ell+1}$; after changing signs of generators, we may assume this sign is $1$.

Consider the acceleration datum $(\tilde{H}_\ell)_{\ell\in \N}$ over $M\times \T^2$, corresponding to the compact subset $K\times (S^1\times\{0\})$ and given by $$\tilde{H}_\ell(t,p,(x,y))=H_\ell(t,p)+h_\ell(t,x,y),$$ for every
$(t,p,(x,y))\in S^1\times M\times \T^2$. For every $\ell\in \N$, denote by $\tilde{J}_\ell$ the almost complex structure $J_\ell\oplus j_\ell$ on $M\times \T^2$, and note that $(\tilde{H}_\ell,\tilde{J}_\ell)$ is regular. In this case, for every $\ell\in \N$, the moduli spaces of $(\tilde{H}_\ell,\tilde{J}_\ell)$ are the products of the moduli spaces of $(H_\ell,J_\ell)$ and $(h_\ell,j_\ell)$, therefore
$$CF(\tilde{H}_{\ell};\Lambda_{\geq0})=CF(H_\ell;\Lambda_{\geq0})\otimes CF(h_\ell;\Lambda_{\geq0})$$

Since the differential of $(h_\ell,j_\ell)$ is zero, we deduce that 
$$d_{\tilde{H}_\ell}(\gamma\otimes\check{x}_\ell)=d_{H_\ell}\gamma \otimes \check{x}_\ell,$$
for every $\gamma\in CF(H_\ell;\Lambda_{\geq0})$, where $d_{\tilde{H}_\ell}$ and $d_{H_\ell}$ are the differentials of $CF(\tilde{H}_\ell;\Lambda_{\geq0})$ and $CF(H_\ell;\Lambda_{\geq0})$, respectively. Similarly, for every $\ell\in \N$, the continuation map $CF(\tilde{H}_\ell;\Lambda_{\geq0})\to CF(\tilde{H}_{\ell+1};\Lambda_{\geq0})$ sends each element of the form $\gamma\otimes \check{x}_\ell$, for some $\gamma\in CF(H_\ell;\Lambda_{\geq0})$, to 
$$\Phi_\ell(\gamma)\otimes T^{h_{\ell+1}(\check{x}_{\ell+1})-h_\ell(\check{x}_\ell)}\check{x}_{\ell+1}$$
where $\Phi_\ell\fc CF(H_\ell;\Lambda_{\geq0})\to CF(H_{\ell+1};\Lambda_{\geq0})$ is the continuation map.

Additionally, since continuation maps preserve the index, we deduce that $\drctlim CF(h_\ell;\Lambda_{\geq0})$, as a cochain complex, has a direct summand of the form $\Lambda_{>0}\check{x}=\drctlim\Lambda_{\geq0}\check{x}_\ell$, where the maps in the direct system are induced by the continuation maps.
This shows that $(\drctlim CF(H_\ell;\Lambda_{\geq0}))\otimes \Lambda_{>0}\check{x}$ is a direct summand of $\drctlim CF(\tilde{H}_\ell;\Lambda_{\geq0})$.

This shows that
$$(\drctlim CF(H_\ell;\Lambda_{\geq0}))\widehat{\otimes} \Lambda_{>0}\check{x}=(\widehat{\drctlim }CF(H_\ell;\Lambda_{\geq0}))\otimes(\widehat{\Lambda_{>0}\check{x}})=(\widehat{\drctlim }CF(H_\ell;\Lambda_{\geq0}))\otimes \Lambda_{>0}\check{x}$$
is a direct summand of 
$$SC(K\times(S^1\times\{0\}))=\widehat{\drctlim}CF(\tilde{H}_\ell;\Lambda_{\geq0}).$$
Hence $SH_M^*(K;\Lambda_{\geq0})\otimes \Lambda_{>0}\check{x}$ is a direct summand, and in particular a submodule, of $SH_{M\times \T^2}(K\times(S^1\times\{0\});\Lambda_{\geq0})$, as promised.
\end{proof}

    \appendix
   
\section{Floer theory with cascades}\label{app: cascades}

\subsection{Flowlines with cascades and moduli spaces}\label{ss: flowlines with cascades and moduli spaces}

Throughout this section, let us fix some variables: Let $(M,\omega)$ be a fixed closed symplectic manifold. Let $H,H_+,H_-$ be Hamiltonians over $M$, satisfying the \textbf{MB} condition. Let $\mathbf{S}_H,\mathbf{S}_{H_+},\mathbf{S}_{H_-}$ denote the union of all the images of the covering map $p\fc\widetilde{\cL_0}(M)\to \cL_0(M)$ of all the critical submanifolds of $\mathcal{A}_H,\mathcal{A}_{H_+},\mathcal{A}_{H_-}$, respectively. Thus, each of $\mathbf{S}_H,\mathbf{S}_{H_+},\mathbf{S}_{H_-}$ is a finite disjoint union of smooth closed, finite-dimensional manifolds of varying dimensions.
   Let $h,h_+,h_-$ be Morse functions on $\mathbf{S}_H,\mathbf{S}_{H_+},\mathbf{S}_{H_-}$, respectively. 
   Let $J,J_+,J_-$ be compatible almost complex structures, possibly time dependent. Also, let $g,g_+,g_-$ be Riemannian metrics on $\mathbf{S}_H,\mathbf{S}_{H_+},\mathbf{S}_{H_-}$, respectively. Moreover, let $\cH=(H_s)_{s\in \R}$ be a homotopy of Hamiltonians and $\cJ=(J_s)_{s\in \R}$ a homotopy of compatible almost complex structures such that $H_s = H_-$ and $J_s = J_-$ for every $s\leq 0$, and $H_s = H_+$ and $J_s = J_+$ for $s\geq 1$.

   Let $g$ be a Riemannian metric  on $\mathbf{S}_H$, chosen generically so that the Floer--Morse--Bott package below is well-defined from the point of view of transversality.

Given two disjoint submanifolds $S^-, S^+ \subset \cL_0(M)$, we denote by $\pi_2(S^-, S^+)$ the collection of all homotopy classes of cylinders in $M$ relative to $S^-$ and $S^+$; that is, homotopy classes of smooth maps $u \fc \R \times S^1 \to M$ with a negative asymptote in $S^-$ and a positive asymptote in $S^+$.

   \begin{defin}
      
    Let $S^-$ and $S^+$ be two connected components of $\bS_{H_-}$ and $\bS_{H_+}$, respectively. Let us denote by $\cM^{cont}(S^-, S^+; \cH,\cJ)$ the collection of all continuation solutions connecting a $1$-periodic orbit in $S^-$ to a $1$-periodic orbit in $S^+$. Similarly, for every $A \in \pi_2(S^-, S^+)$, we denote by $\cM^{cont}(S^-, S^+; \cH,\cJ; A)$ the moduli space of all continuation solutions representing the homotopy class $A$.

   \end{defin}

\begin{rem}
\begin{itemize}
    \item In the case where $H_-=H_+$ and $\cH,\cJ$ are constant homotopies, the moduli space $\cM^{cont}(S^-,S^+;\cH,\cJ)$ consists of Floer solutions, for this reason we denote it by $\cM(S^-,S^+;H,J)$, where $J=J_-=J_+$, possibly time-dependent, and $H=H_-=H_+$. Additionally,  $\cM(S^-,S^+;H,J)$ admits a natural $\R$-action.
    
    \item  
In cases where $\cH$ and $\cJ$ are clear from the context, we use the shortened notation $\cM(S^-, S^+) = \cM^{cont}(S^-, S^+; \cH,\cJ)$, and for every $A \in \pi_2(S^-, S^+)$, we write $\cM^{A,cont}(S^-, S^+) = \cM(S^-, S^+; \cH,\cJ; A)$.
\end{itemize}    

\end{rem}

    \begin{prop}\label{prop: ind formula for Continuation sol}
        
        Let $S^-$ and $S^+$ be critical submanifolds, components of $\bS_{H_-}$ and $\bS_{H_+}$ respectively, and let $A\in\pi_2(S^-,S^+)$ be a homotopy class.
        Moreover, let $\tau_-$ and $\tau_+$ be families of compatible trivializations of $\det TM$ along the orbits of $S^-$ and $S^+$, respectively, and denote $\tau=\tau_+\cup\tau_-$.
        
        Then, for a generic choice of a homotopy of Hamiltonians $\cH=(H_s)_{s\in \R}$, such that $H_s = H_-$ for $s\leq0$ and $H_s = H_+$ for $s\geq1$, the moduli space $\cM(S^-,S^+;\cH,\cJ;A)$ is a smooth manifold and its dimension is given by
        \[\dim \cM^{A,cont}(S^-,S^+)=\RS^\tau(S^+)-\RS^\tau(S^-)+\frac{1}{2}\left(\dim S^- +\dim S^+\right)+2c_1^\tau(A),\]
        where $c^\tau_1$ is the relative first Chern class\footnote{See \cite[Definition 5.1]{Wendl_16_Lec_SFT}} with respect to $\tau$.
    \end{prop}

    For this moment, let us focus on Floer trajectories.

   \begin{defin}\label{def:floer flowline with cascades}
   
   Let $S^-,S^+\subset \bS_H$ be critical submanifolds, and let  $q_-,q_+\in\Crit(h)$ be a pair of critical points of $h$, lying on $S^-$ and $S^+$ respectively.
       \begin{enumerate}
        \item A \textbf{Floer flowline with zero cascades} from $q_-$ to $q_+$, $\mathbf u = (\gamma)$ is a non-constant gradient flowline of $h$, i.e. a smooth map $\gamma \fc \R \to \mathbf{S}$, such that 
           $$  \dot\gamma=\nabla_gh\circ \gamma,\qquad\text{and}\qquad
            \lim_{t\to\pm\infty}\gamma(t) = q_\pm.$$
            Note that in this case $S^- = S^+$.
        \item A \textbf{Floer flowline with} $\mathbf{1}$ \textbf{cascade} from $q_-$ to $q_+$ consists of data
            \[\mathbf{u}=\left( \left(S^-,S^+\right),\left(u\right),\left(\gamma_1,\gamma_2\right)\right),\]
            where
            \begin{enumerate}
                \item $u\in \cM\left(S^-,S^{+}\right)$, a Floer solution, non-constant in the $\R$ variable.
                
                \item $\gamma_1\fc (-\infty,0] \to \bS_H$ and $\gamma_{2}\fc [0,+\infty) \to \bS_H$, are all gradient flowlines of $h$, namely $\dot{\gamma}_1 = \nabla_gh\circ\gamma_1$ and $\dot{\gamma}_2 = \nabla_gh\circ\gamma_2$, such that:
                \begin{enumerate}
                    \item $q_- =\lim\limits_{t\to-\infty} \gamma_1(t)$, and $\gamma_1(0) = \lim\limits_{s\to-\infty} u(s,\cdot)$.
                    \item $\lim\limits_{s\to+\infty} u(s,\cdot) = \gamma_{2}(0)$, and $\lim\limits_{t\to+\infty} \gamma_{2}(t) = q_+$.
                \end{enumerate}
            \end{enumerate}

                    \item Let $m\in\N$ with $m\geq2$. A \textbf{Floer flowline with} $\mathbf{m}$ \textbf{cascades} from $q_-$ to $q_+$ consists of data
            \[\mathbf{u}=\left( \left(S^k\right)_{k=1}^{m+1},\left(u_k\right)_{k=1}^m,\left(t_k\right)_{k=2}^{m},\left(\gamma_k\right)_{k=1}^{m+1}\right),\]
            where
            \begin{enumerate}
                \item each $S^k$ is a connected component of $\bS_H$,
                \item $S^1 = S^-$ and $S^+ = S^{m+1}$.
                \item $u_k \in \cM\left(S^k,S^{k+1}\right)$, a Floer solution, non-constant in the $\R$ variable.
                \item $t_2,\ldots,t_m\in (0,+\infty)$ are positive numbers.
                \item $\gamma_1\fc (-\infty,0] \to \bS_H$, $\gamma_k\fc [0,t_k] \to \bS_H$, for $2\le k \le m$ and $\gamma_{m+1}\fc [0,+\infty) \to \bS_H$, are all gradient flowlines of $h$, namely $\dot{\gamma}_k = \nabla_gh\circ\gamma_k$ for $1\le k \le m+1$, such that:
                \begin{enumerate}
                    \item $q_- =\lim\limits_{t\to-\infty} \gamma_1(t)$, and $\gamma_1(0) = \lim\limits_{s\to-\infty} u_1(s,\cdot)$.
                    \item For all $2\le k \le m$, $\lim\limits_{s\to\infty} u_{k-1}(s,\cdot) = \gamma_k(0)$ and $\gamma_k(t_k)=\lim\limits_{s\to-\infty} u_{k}(s,\cdot)$.
                    \item $\lim\limits_{s\to+\infty} u_{m}(s,\cdot) = \gamma_{m+1}(0)$, and $\lim\limits_{t\to+\infty} \gamma_{m+1}(t) = q_+$.
                \end{enumerate}
            \end{enumerate}
            
       \end{enumerate}
   \end{defin}
   \begin{rem}
       Note that by the existence and uniqueness of solutions to ODEs, the data of $\left(t_k\right)_{k=2}^{m},\left(\gamma_k\right)_{k=1}^{m+1}$, in the last definition, is determined by $\left(u_k\right)_{k=1}^{m}$; We choose to make them explicit for the geometric clarity and ease of dealing with Floer-Gromov limits.

   \end{rem}
    \begin{defin}

    Let $q_-,q_+\in\Crit(h)$ be a pair of critical points of $h$.
     For every $m\in \Z_{\geq0}$ we denote \textbf{the moduli space of Floer flowlines with $m$ cascades connecting $q_-$ to $q_+$} by $\cM_m(q_-,q_+;H,J)$. Moreover,  we denote the \textbf{space of all Floer flowlines with cascades connecting $q_-$ to $q_+$} by $\cM(q_-,q_+) := \bigcup_{m=0}^\infty \cM_m(q_-,q_+)$.

   \end{defin}
   \begin{rem}
   \phantom{ }
   \begin{itemize}
       \item As before, we abbreviate to $\cM_m(q_-,q_+)$ when $H$ and $J$ are clear from the context.
    \item 
       Note that $\R^m$ acts freely on $\cM_m(q_-,q_+)$ by additive shifts in the $\R$-coordinate of the Floer solutions.
   \end{itemize}
       
   \end{rem}

    \begin{defin}
        For every pair of critical submanifolds $S_-$ and $S_+$, denote the evaluation maps $\ev^-\fc\cM(S^-,S^+) \to S^-$,  and $\ev^+\fc \cM(S^-,S^+) \to  S^+$ defined by 
        $$ \ev^-(u) =\lim\limits_{s\to-\infty}u_1(s,\cdot),\qquad\text{and}\qquad
            \ev^+(u) = \lim\limits_{s\to+\infty}u_m(s,\cdot),$$
            for every $u\in  \cM(S^-,S^+)$.
    \end{defin}
    \begin{prop}\label{prop: structure_and_dimension_of_moduli_space_of_floer_w_cascascades}

        For every $t\in \R$ denote by $\varphi_{h,g}^t$ the gradient flow of $h$ on $\bS_H$ at time $t$. For every connected critical submanifold $S$, define $\varphi_{S}\fc S\times \R\to S$ by $\varphi(q,t)=\varphi_{h,g}^t(q)$ for every $(q,t)\in S\times \R$.
        
         Let $q_-,q_+\in\Crit(h)$ be a pair of critical points of $h$, lying on the critical submanifolds $S^-,S^+\subset\bS_H$, respectively. Denote by $W^u(q^-)\subset S^-$ and $W^s(q^+)\subset S^+$ the unstable and stable manifolds of $q^-$ and $q^+$ respectively.
         
        \begin{itemize}
            \item Let $m\in \Z_{\geq0}$. The moduli space $\cM_m(q_-,q_+)$ can be constructed as follows: 
            \begin{itemize}
                \item If $m=0$ then $S^-=S^+$ and $\cM_0(q_-,q_+)=W^u(q_-)\cap W^s(q_+)$ if $q_-\neq q_+$, and $\cM_0(q_-,q_+)=\varnothing$ if $q_- = q_+$.
                \item If $m=1$, then $\cM_1(q_-,q_+)$ is given by the following fiber product:
                $$W^u(q_-)\times_{(\id,\ev^-)}\cM(S^-,S^+)\times_{(\ev^+,\id)}W^s(q_+).$$
                \item If $m\geq 2$, then $\cM_m(q_-,q_+)$ is given by the union over all the critical submanifolds $S^2,\ldots,S^m\subset\bS_H$ of the fiber products of the following type:
                
            \end{itemize}
          \end{itemize}
        \begin{multline*}            W^u(q^-)\times_{\left(\id,\ev^-\right)}\cM(S^-,S^2)\times_{\left(\ev^+,\,\varphi_{S^2}\circ(\ev^-,-\id)\right)}\left(\cM(S^2,S^3)\times [0,+\infty)\right)\times_{\left(\ev^+,\,\varphi_{S^3}\circ(\ev^-,-\id)\right)} \ldots \\\ldots \times_{\left(\ev^+,\,\varphi_{S^{m-1}}\circ(\ev^-,-\id)\right)} \left(\cM(S^{m-1},S^{m})\times [0,+\infty)\right)\times_{\left(\ev^+,\,\varphi_{S^m}\circ(\ev^-,-\id)\right)} \cM(S^{m},S^+)\times_{\left(\ev^+,\id\right)}W^s(q_+),
        \end{multline*}
           
         \begin{itemize}
        \item  A generic choice of $J$ makes all the evaluation maps submersions, hence the fiber product is a smooth manifold with corners, as a fiber product of smooth manifolds with boundary. 

    \item  Let $m\in \N$ and let $u\in \cM_m(q_-,q_+)$. For every choice $\tau$ of families of compatible trivializations
        of $\det TM$ along the orbits in the connected components of $\bS_H$, the local dimension of $\cM_m(q_-,q_+)$ near $u$ is given by 
    
        $$
        \mu^\tau_{FMB}(q_+) - \mu^\tau_{FMB}(q_-) + 2\sum_{k=1}^m c_1^\tau(u_k) + m-1,$$
        where  $\left(u_k\right)_{k=1}^m$ are the cascades of $u$. \qed
        \end{itemize}
    \end{prop}

    \begin{defin}  
    
    Let $q_-\in\Crit(h_-)$ and $q_+\in\Crit(h_+)$ be two critical points lying on two critical submanifolds $S^-\subseteq \bS_{H_-}$ and $S^+\subseteq\bS_{H_+}$ respectively.
       \begin{enumerate}

        \item A \textbf{continuation flowline with cascades of type $(0;0)$} from $q_-$ to $q_+$
         consists of data
            \[\left(\left( \left(S^-\right),\left(\gamma_1^-\right)\right),(w),\left( \left(S^+\right),\left(\gamma_1^+\right)\right)\right)\]
            where
            \begin{enumerate}
               
                \item $w \in \cM^{cont}\left(S^-,S^+;\cH,\cJ\right)$, a continuation solution.
                \item $\gamma^-_1\fc (-\infty,0] \to S^-\subset \bS_{H_-}$ is a gradient flowline of $h_-$ and $\gamma^+_1\fc [0,+\infty) \to S^+\subset \bS_{H_+}$ is a gradient flowline of $h_+$.  They are required to satisfy:
                \begin{enumerate}
                    \item $q_- =\lim\limits_{t\to-\infty} \gamma^-_1(t)$, and $\gamma^-_1(0) = \lim\limits_{s\to-\infty} w(s,\cdot)$.
                    \item 
                    $\lim\limits_{s\to+\infty}w(s,\cdot)=\gamma^+_1(0)$, and $\lim\limits_{t\to+\infty} \gamma^+_{1}(t) = q_+$.
                \end{enumerate}
            \end{enumerate}

             \item Let $m_1\in \N$. A \textbf{continuation flowline with cascades of type $(m_1;0)$} from $q_-$ to $q_+$
         consists of data
            \[\left(\left( \left(S_-^k\right)_{k=1}^{m_1+1},\left(u^-_k\right)_{k=1}^{m_1},\left(t^-_k\right)_{k=2}^{m_1+1},\left(\gamma^-_k\right)_{k=1}^{m_1+1}\right),(w),\left( \left(S^+\right),\left(\gamma^+_1\right)\right)\right)\]
            where
            \begin{enumerate}
                \item For every $1\leq k \leq m_1+1 $, $S_-^k$ is a connected component of $\bS_{H_-}$, moreover, $S_-^1 = S^-$.
                \item For every $1\leq k\leq m_1$, $u^-_k \in \cM\left(S_-^k,S_-^{k+1};H_-,J_-\right)$ is a Floer solution.
                \item $w \in \cM^{cont}\left(S^{m_1+1}_-,S^+;\cH,\cJ\right)$, a continuation solution.
                \item $\gamma^-_1\fc (-\infty,0] \to S^-\subset\bS_{H_-}$, $\gamma^-_k\fc [0,t^-_k] \to \bS_{H_-}$, for $2\le k \le m_1+1$, are all gradient flowlines of $h_-$. Similarly, $\gamma^+_1\fc [0,+\infty) \to S^+\subset\bS_{H_+}$, is a gradient flowline of $h_+$.
                They are required to satisfy:
                \begin{enumerate}
                    \item $q_- =\lim\limits_{t\to-\infty} \gamma^-_1(t)$, and $\gamma^-_1(0) = \lim\limits_{s\to-\infty} u^-_1(s,\cdot)$.
                    \item For all $2\le k \le m_1$, $\lim\limits_{s\to+\infty} u^-_{k-1}(s,\cdot) = \gamma^-_k(0)$ and $\gamma^-_k(t^-_k)=\lim\limits_{s\to-\infty} u^-_{k}(s,\cdot)$.
                    \item $\lim\limits_{s\to+\infty}u^-_{m_1}(s,\cdot)=\gamma^-_{m_1+1}(0)$ and $\gamma^-_{m_1+1}(t^-_{m_1+1})= \lim\limits_{s\to-\infty}w(s,\cdot).$
                    \item 
                    $\lim\limits_{s\to+\infty}w(s,\cdot)=\gamma^+_1(0)$, and $\lim\limits_{t\to+\infty} \gamma^+_{1}(t) = q_+$.
                \end{enumerate}
            \end{enumerate}

        \item Let $m_2\in \N$. A \textbf{continuation flowline with cascades of type $(0;m_2)$} from $q_-$ to $q_+$
         consists of data
            \[\left(\left( \left(S^-\right),\left(\gamma_1^-\right)\right),(w),\left( \left(S_+^k\right)_{k=1}^{m_2+1},\left(u^+_k\right)_{k=1}^{m_2},\left(t^+_k\right)_{k=1}^{m_2},\left(\gamma^+_k\right)_{k=1}^{m_2+1}\right)\right)\]
            where
            \begin{enumerate}
                \item For every $1\leq k\leq m_2+1$, $S_+^k$ is a connected component of $\bS_{H_+}$, moreover, $S^+=S_+^{m_2+1}$.
                \item For every $1\leq k\leq m_2$, $u^+_k \in \cM\left(S_+^k,S_+^{k+1};H_+,J_+\right)$ is a Floer solution.
                \item $w \in \cM^{cont}\left(S^{-},S_+^{1};\cH,\cJ\right)$, a continuation solution.
                \item $\gamma^-_1\fc (-\infty,0] \to S^-\subset \bS_{H_-}$ is a gradient flowline of $h_-$. Additionally, $\gamma^+_k\fc[0,t^+_k] \to \bS_{H_+}$, for $1\le k \le m_2$, and $\gamma^+_{m_2+1}\fc [0,+\infty) \to S^+\subset\bS_{H_+}$, are all gradient flowlines of $h_+$.  They are required to satisfy:
                \begin{enumerate}
                    \item $q_- =\lim\limits_{t\to-\infty} \gamma^-_1(t)$, and $\gamma^-_1(0) = \lim\limits_{s\to-\infty} w(s,\cdot)$.
                    \item 
                    $\lim\limits_{s\to+\infty}w(s,\cdot)=\gamma^+_1(0)$, and $\gamma^+_1(t_1)=\lim\limits_{s\to-\infty}u^+_1(s,\cdot)$.
                    \item For all $2\le k \le m_2$, $\lim\limits_{s\to+\infty} u^+_{k-1}(s,\cdot) = \gamma^+_k(0)$ and $\gamma^+_k(t^+_k)=\lim\limits_{s\to-\infty} u^+_{k}(s,\cdot)$.
                    
                    \item$\lim\limits_{s\to+\infty} u_{m_2}(s,\cdot) = \gamma^+_{m_2+1}(0)$, and $\lim\limits_{t\to+\infty} \gamma^+_{m_2+1}(t) = q_+$.
                \end{enumerate}
            \end{enumerate}

             \item Let $m_1,m_2\in \N$. A \textbf{continuation flowline with cascades of type $(m_1;m_2)$} from $q_-$ to $q_+$
         consists of data
            \[\left(\left( \left(S_-^k\right)_{k=1}^{m_1+1},\left(u^-_k\right)_{k=1}^{m_1},\left(t^-_k\right)_{k=2}^{m_1+1},\left(\gamma^-_k\right)_{k=1}^{m_1+1}\right),(w),\left( \left(S_+^k\right)_{k=1}^{m_2+1},\left(u^+_k\right)_{k=1}^{m_2},\left(t^+_k\right)_{k=1}^{m_2},\left(\gamma^+_k\right)_{k=1}^{m_2+1}  \right)\right)\]
            where
            \begin{enumerate}
                \item For every $1\leq k\leq m_1+1$,  $S_-^k$ is a connected component of $\bS_{H_-}$, moreover $S_-^1=S^-$.
                \item For every $1\leq k\leq m_2+1$, $S_+^k$ is a connected component of $\bS_{H_+}$, moreover,  $S^+=S_+^{m_2+1}$.
                \item For every $1\leq k\leq m_1$, $u^-_k \in \cM\left(S_-^k,S_-^{k+1};H_-,J_-\right)$ is a Floer solution.
                \item For every $1\leq k\leq m_2$,  $u^+_k \in \cM\left(S_+^k,S_+^{k+1};H_+,J_+\right)$ is Floer solution.
                \item $w \in \cM^{cont}\left(S^{m_1+1}_-,S_+^{1};\cH,\cJ\right)$, a continuation solution.
                \item $\gamma^-_1\fc (-\infty,0] \to S^-\subset\bS_{H_-}$, $\gamma^-_k\fc [0,t^-_k] \to \bS_{H_-}$, for $2\le k \le m_1+1$, are all gradient flowlines of $h_-$. Similarly,   $\gamma^+_k\fc [0,t^+_k] \to \bS_{H_+}$, for $1\le k \le m_2$, and $\gamma^+_{m_2+1}\fc [0,+\infty) \to S^+\subset\bS_{H_+}$, are all gradient flowlines of $h_+$. They are required to satisfy:
                \begin{enumerate}
                    \item $q_- =\lim\limits_{t\to-\infty} \gamma^-_1(t)$, and $\gamma^-_1(0) = \lim\limits_{s\to-\infty} u^-_1(s,\cdot)$.
                    \item For all $2\le k \le m_1$, $\lim\limits_{s\to+\infty} u^-_{k-1}(s,\cdot) = \gamma^-_k(0)$ and $\gamma^-_k(t^-_k)=\lim\limits_{s\to-\infty} u^-_{k}(s,\cdot)$.
                    \item $\lim\limits_{s\to+\infty}u^-_{m_1}(s,\cdot)=\gamma^-_{m_1+1}(0)$ and $\gamma^-_{m_1+1}(t^-_{m_1+1})= \lim\limits_{s\to-\infty}w(s,\cdot).$
                    \item 
                    $\lim\limits_{s\to+\infty}w(s,\cdot)=\gamma^+_1(0)$, and $\gamma^+_1(t_1)=\lim\limits_{s\to-\infty}u^+_1(s,\cdot)$.
                    \item For all $2\le k \le m_2$, $\lim\limits_{s\to+\infty} u^+_{k-1}(s,\cdot) = \gamma^+_k(0)$ and $\gamma^+_k(t^+_k)=\lim\limits_{s\to-\infty} u^+_{k}(s,\cdot)$.\item$\lim\limits_{s\to+\infty} u_{m_2}(s,\cdot) = \gamma^+_{m_2+1}(0)$, and $\lim\limits_{t\to+\infty} \gamma^+_{m_2+1}(t) = q_+$.
                \end{enumerate}
            \end{enumerate}
           
            \item Let $m\in \N$. A \textbf{continuation flowline with $m$ cascades} from $q_-$ to $q_+$ is a continuation flowline of type $(m_1;m_2)$ from $q_-$ to $q_+$ for some $0\le m_1$, $0\le m_2$ such that $m_1 + m_2=m-1$.

       \end{enumerate}
   \end{defin}

\begin{rem}
\begin{itemize}
    \item Note that by the existence and uniqueness of solutions to ODEs, $\left(u^-_k\right)_{k=1}^{m_1}, w$ and $\left(u^+_k\right)_{k=1}^{m_2}$, from the last definition, determine the rest of the data; We chose to make them explicit for the geometric clarity and ease of dealing with Floer-Gromov limits.
    \item  In the definition of a continuation flowline with $m$ cascades, note that while we require that $m\ge1$, namely that we cannot have zero cascades in a continuation flowline, the continuation solution is allowed to be constant in the $\R$ variable.
\end{itemize}

\end{rem}

   Let us move to discuss continuation flowlines with cascades, which generalize the Floer flowlines with cascades.
   
   \begin{defin}

    Let $q_-\in\Crit(h_-)$ and $q_+\in\Crit(h_+)$.

     \begin{itemize}
    \item For every $m_1,m_2\in \Z_{\geq0}$, denote \textbf{the moduli space of continuation flowlines with cascades of type $(m_1;m_2)$ connecting $q_-$ to $q_+$} by $\cM^{cont}_{(m_1;m_2)}(q_-,q_+;\cH,\cJ)$.
        \item For every $m\in\N$. We denote \textbf{the moduli space of continuation flowlines with $m$ cascades connecting $q_-$ to $q_+$} by $\cM^{cont}_m(q_-,q_+;\cH,\cJ)$.
        \item  We denote the \textbf{space of all continuation flowlines with cascades connecting $q_-$ to $q_+$} by $\cM^{cont}(q_-,q_+;\cH,\cJ) := \bigcup_{m=1}^\infty \cM^{cont}_m(q_-,q_+;\cH,\cJ)$.
    \end{itemize}

   \end{defin}
   \begin{rem}
   \begin{itemize}
       \item  When $\cH$ and $\cJ$ are clear from the context we abbreviate to $\cM^{cont}(q_-,q_+)$, for every $m\in \N$ to $\cM^{cont}_m(q_-,q_+)$, and for every $m_1,m_2\in \Z_{\geq0}$ to $\cM^{cont}_{(m_1;m_2)}(q_-,q_+)$, respectively.
       \item Note that $\R^{m-1}$ acts freely on $$\cM^{cont}_m(q_-,q_+)=\bigcup_{ \substack{
                   m_1,m_2\geq0,\\
                    m_1+m_2=m-1
                }}\cM^{cont}_{(m_1;m_2)}(q_-,q_+)$$
                by additive shifts in the $\R$ coordinate of the Floer solutions.
   \end{itemize}

   \end{rem}

 \begin{prop}\label{prop: structure_and_dimension_of_moduli_space_of_continuation_w_cascascades}

     For every $t\in \R$ denote by $\varphi_{h_-,g_-}^t$ and $\varphi_{h_+,g_+}^t$ the gradient flows of $h_-$ and $h_+$ on $\bS_{H_-}$ and $\bS_{H_+}$, respectively, at time $t$. For every connected critical submanifold $S\subset\bS_{H_-}$, define $\varphi^-_{S}\fc S\times \R\to S$ by $\varphi^-(q,t)=\varphi_{h_-,g_-}^t(q)$ for every $(q,t)\in S\times \R$. Similarly, for every connected critical submanifold $S\subset\bS_{H_+}$, define $\varphi^+_{S}\fc S\times \R\to S$ by $\varphi^+(q,t)=\varphi_{h_+,g_+}^t(q)$ for every $(q,t)\in S\times \R$.

    Let $q_-\in\Crit(h_-)$ and $q_+\in\Crit(h_+)$ be two critical points lying on two critical submanifolds $S^-\subseteq \bS_{H_-}$ and $S^+\subseteq\bS_{H_+}$ respectively.
     Denote by $W^u(q^-)\subset S^-$ and $W^s(q^+)\subset S^+$ the unstable and stable manifolds of $q^-$ and $q^+$ respectively.
         
        \begin{itemize}
            \item Let $m_1,m_2\in \Z_{\geq0}$. The moduli space $\cM^{cont}_{(m_1;m_2)}(q_-,q_+)$ can be constructed as follows: 
            \begin{itemize}
               
                \item If $m_1=m_2=0$, then $\cM^{cont}_{(0;0)}(q_-,q_+)$ is given by the following fiber product:
                $$W^u(q_-)\times_{(\id,\ev^-)}\cM^{cont}(S^-,S^+)\times_{(\ev^+,\id)}W^s(q_+).$$
                
                \item If $m_1=0$ and $m_2\geq1$ then the moduli space  $\cM^{cont}_{(m_1;m_2)}(q_-,q_+)$
                is given by the union over all the critical submanifolds $S^1_+,\ldots,S^{m_2}_+\subset\bS_{H_+}$ of the fiber products of the following type:
                
            \end{itemize}
          \end{itemize}

        \begin{multline*}            W^u(q^-)\times_{\left(\id,\ev^-\right)}\cM^{cont}(S^-,S^1_+)\times_{\left(\ev^+,\,\varphi_{S^1_+}^+\circ(\ev^-,-\id)\right)}\left(\cM(S^1_+,S^2_+)\times [0,+\infty)\right)\times_{\left(\ev^+,\,\varphi_{S^2_+}^+\circ(\ev^-,-\id)\right)} \ldots \\\ldots \times_{\left(\ev^+,\,\varphi_{S^{m_2-1}_+}^+\circ(\ev^-,-\id)\right)} \left(\cM(S^{m_2-1}_+,S^{m_2}_+)\times [0,+\infty)\right)\\
        \times_{\left(\ev^+,\,\varphi_{S^{m_2}_+}^+\circ(\ev^-,-\id)\right)} \cM(S^{m_2}_+,S^+)\times_{\left(\ev^+,\id\right)}W^s(q_+),
        \end{multline*}
        
           \begin{itemize}
                   \item If $m_1\geq 1$ and $m_2=0$ then the moduli space  $\cM^{cont}_{(m_1;m_2)}(q_-,q_+)$
                is given by the union over all the critical submanifolds $S^2_-,\ldots,S^{m_1+1}_-\subset\bS_{H_-}$ of the fiber products of the following type:
           \end{itemize}

        \begin{multline*}            W^u(q^-)\times_{\left(\id,\ev^-\right)}\cM(S^-,S^2_1)\times_{\left(\ev^+,\,\varphi_{S^2_-}^-\circ(\ev^-,-\id)\right)}\left(\cM(S^2_-,S^3_-)\times [0,+\infty)\right)\times_{\left(\ev^+,\,\varphi_{S^3_-}^-\circ(\ev^-,-\id)\right)} \ldots \\\ldots \times_{\left(\ev^+,\,\varphi_{S^{m_1}_-}^-\circ(\ev^-,-\id)\right)} \left(\cM(S^{m_1}_-,S^{m_1+1}_-)\times [0,+\infty)\right)\\
        \times_{\left(\ev^+,\,\varphi_{S^{m_1+1}_-}^-\circ(\ev^-,-\id)\right)} \cM^{cont}(S^{m_1+1}_-,S^+)\times_{\left(\ev^+,\id\right)}W^s(q_+),
        \end{multline*}

    \begin{itemize}
        \item If $m_1,m_2\geq 1$ then the moduli space  $\cM^{cont}_{(m_1;m_2)}(q_-,q_+)$
                is given by the union over all the critical submanifolds $S^2_-,\ldots,S^{m_1+1}_-\subset\bS_{H_-}$ and $S^1_+,\ldots,S^{m_2}_+\subset\bS_{H_+}$ of the fiber products of the following type:
                
    \end{itemize}

        \begin{multline*}            W^u(q^-)\times_{\left(\id,\ev^-\right)}\cM(S^-,S^2_1)\times_{\left(\ev^+,\,\varphi_{S^2_-}^-\circ(\ev^-,-\id)\right)}\left(\cM(S^2_-,S^3_-)\times [0,+\infty)\right)\times_{\left(\ev^+,\,\varphi_{S^3_-}^-\circ(\ev^-,-\id)\right)} \ldots \\\ldots \times_{\left(\ev^+,\,\varphi_{S^{m_1}_-}^-\circ(\ev^-,-\id)\right)} \left(\cM(S^{m_1}_-,S^{m_1+1}_-)\times [0,+\infty)\right)\\
        \times_{\left(\ev^+,\,\varphi_{S^{m_1+1}_-}^-\circ(\ev^-,-\id)\right)} \left(\cM^{cont}(S^{m_1+1}_-,S^1_+)
        \times [0,+\infty)\right)
        \\
        \times_{\left(\ev^+,\,\varphi_{S^1_+}^+\circ(\ev^-,-\id)\right)}\left(\cM(S^1_+,S^2_+)
        \times [0,+\infty)\right)\times_{\left(\ev^+,\,\varphi_{S^2_+}^+\circ(\ev^-,-\id)\right)} \ldots \\\ldots \times_{\left(\ev^+,\,\varphi_{S^{m_2-1}_+}^+\circ(\ev^-,-\id)\right)} \left(\cM(S^{m_2-1}_+,S^{m_2}_+)\times [0,+\infty)\right)\\
        \times_{\left(\ev^+,\,\varphi_{S^{m_2}_+}^+\circ(\ev^-,-\id)\right)} \cM(S^{m_2}_+,S^+)\times_{\left(\ev^+,\id\right)}W^s(q_+),
        \end{multline*}

         \begin{itemize}
        \item  A generic choice of $\cJ$ and $\cH$ makes all the spaces in the product into manifolds and the evaluation maps submersions, hence $\cM^{cont}_m(q_-,q_+)$ is a smooth manifold with corners, as a disjoint union of fiber products of smooth manifolds with boundary.

    \item  Let $m\in \N$ and let $u\in \cM^{cont}_m(q_-,q_+)$. For every choice $\tau$ of families of compatible trivializations
        of $\det TM$ along the orbits in the connected components of $\bS_H$, the local dimension of $\cM_m(q_-,q_+)$ near $u$ is given by 
    
        $$
        \mu^\tau_{FMB}(q_+) - \mu^\tau_{FMB}(q_-) + 2\sum_{k=1}^m c_1^\tau(u_k) + m-1,$$
        where  $\left(u_k\right)_{k=1}^m$ are the cascades of $u$. \qed
        \end{itemize}
    
    \end{prop}

     \begin{defin}

     Let $q_-,q_+\in\Crit(h)$ be a pair of critical points of $h$.
     For every $m\in \Z$, $m\ge1$, we define \textbf{the moduli space of unparametrized Floer flowlines with $m$ cascades connecting $q_-$ to $q_+$}, as the quotient of $\cM_m(q_-,q_+;H,J)$, by the $\R^m$ action acting by additive shifts on the cascdes, and we denote it by $\widehat{\cM}_m(q_-,q_+;H,J)$, namely
     \[ \widehat{\cM}_m(q_-,q_+;H,J):= \cM_m(q_-,q_+;H,J) \mathbin{\big/} \R^m.\]
     For $m=0$ we set
      \[ \widehat{\cM}_0(q_-,q_+;H,J):= \cM_0(q_-,q_+;H,J) \mathbin{\big/} \R,\] where the $\R$ action acts on the single morse flowline by translations.
     \end{defin}
     \begin{rem}As before, we abbreviate to $\widehat{\cM}_m(q_-,q_+)$ when $H$ and $J$ are clear from the context.
     \end{rem}
     \begin{defin}\label{def: moduli of unparametrized Floer flowlines of m cascades}

          Let $q_-,q_+\in\Crit(h)$ be a pair of critical points of $h$. We define the \textbf{space of unparametrized Floer flowlines with at most $m$ cascades connecting $q_-$ to $q_+$} by $\widehat{\cM}_{\le m}(q_-,q_+) := \bigcup_{k=0}^m \widehat{\cM}_k(q_-,q_+)$.
     \end{defin}
     The space $\widehat{\cM}_{\le m}(q_-,q_+)$ can be given the structure of a smooth manifold, using the charts provided by gluing of Floer solutions. Namely, a representative,  $\left(\left(u_k\right)_{k=1}^{m^\prime},\left(t_k\right)_{k=2}^{m^\prime}\right)$ of an element lying in the boundary of $\cM_{m^\prime}(q_-,q_+)$ for some $m^\prime\leq  m $, has some of the parameters $t_k$ equal to $0$. On one hand, such element can be seen as a degeneration of a family of elements with all $t_k>0$, which provides a one sided neighborhood.
     On the other hand, it can be seen as a degeneration by breaking of a family of elements with $m^{\prime\prime}<m^\prime$ cascades, for which gluing provides a structure of a one sided neighborhood. Gluing these neighborhoods provide manifold charts. See \cite{Banyaga_Hurtubise_2013_Cascades}, for a detailed construction in the setting of Morse homology.

     \begin{prop} 
         
         Let $q_-,q_+\in\Crit(h)$ be a pair of critical points of $h$. For every $m\in \N$, for every choice, $\tau$, of families of compatible trivializations
        of $\det TM$ along the orbits in the connected components of $\bS_H$ and for every representative, $u$, of any equivalence class in $\widehat{\cM}_{\le m}(q_-,q_+)$ the dimension of $\widehat{\cM}_{\le m}(q_-,q_+)$ is given by 
        \[
        \dim \widehat{\cM}_{\le m}(q_-,q_+) = \mu^\tau_{FMB}(q_+) - \mu^\tau_{FMB}(q_-) + 2\sum_{k=1}^{m^\prime} c_1^\tau(u_k) -1,\]
        where $m^\prime$ is the number of cascades in $u$ and   $\left(u_k\right)_{k=1}^{m^\prime}$ are the cascades of $u$. \qed
     \end{prop}
     \begin{defin}\label{def: moduli of unparametrized Floer flowlines}

Let $q_-,q_+\in\Crit(h)$ be a pair of critical points of $h$. The 
         \textbf{space of all unparametrized Floer flowlines with cascades connecting $q_-$ to $q_+$} by $\widehat{\cM}(q_-,q_+) := \bigcup_{m=0}^\infty \widehat{\cM}_{\le m}(q_-,q_+)$.
   \end{defin}

    \begin{defin}\label{def:unparamterizedContWithCascades}

        Let $q_-\in\Crit(h_-)$ and $q_+\in\Crit(h_+)$ be a pair of critical points of $h_-$ and $h_+$, respectively.
    
     For every $m\in \Z$, $m\ge 1$ we define \textbf{the moduli space of unparametrized continuation flowlines with $m$ cascades connecting $q_-$ to $q_+$}, as the quotient of $\cM^{cont}_m(q_-,q_+;H,J)$, by the $\R^{m-1}$ action acting by additive shifts on the cascades which are Floer solutions, and we denote it by $\widehat{\cM}^{cont}_m(q_-,q_+;H,J)$, namely
     \[ \widehat{\cM}^{cont}_m(q_-,q_+;H,J):= \cM^{cont}_m(q_-,q_+;H,J) \mathbin{\big/} \R^{m-1}.\] 
     \end{defin}
     \begin{rem}As before, we abbreviate to $\widehat{\cM}^{cont}_m(q_-,q_+)$ when $\cH$ and $\cJ$ are clear from the context.
     \end{rem}
     \begin{defin}\label{def:unparamterizedContWithCascades-general}
          
          Let $q_-\in\Crit(h_-)$ and $q_+\in\Crit(h_+)$ be a pair of critical points of $h_-$ and $h_+$, respectively.           We define the \textbf{space of unparametrized continuation flowlines with at most $m$ cascades connecting $q_-$ to $q_+$} by $\widehat{\cM}^{cont}_{\le m}(q_-,q_+) := \bigcup_{k=1}^m \widehat{\cM}^{cont}_k(q_-,q_+)$.
     \end{defin}
     Similarly to the case of uparametrized Floer flowlines, the space $\widehat{\cM}^{cont}_{\le m}(q_-,q_+)$ can be given the structure of a smooth manifold, using the charts provided by gluing of Floer or continuation solutions.

     \begin{prop} 
    
    Let $q_-\in\Crit(h_-)$ and $q_+\in\Crit(h_+)$ be a pair of critical points of $h_-$ and $h_+$, respectively.
    
    For every $m\in \N$, for every choice, $\tau$, of families of compatible trivializations
    of $\det TM$ along the orbits in the connected components of $\bS_{H_-} \cup \bS_{H_+}$ and for every representative, $u$, of any equivalence class in $\widehat{\cM}^{cont}_{\le m}(q_-,q_+)$ the dimension of $\widehat{\cM}^{cont}_{\le m}(q_-,q_+)$ is given by 
        \[
        \dim \widehat{\cM}^{cont}_{\le m}(q_-,q_+) = \mu^\tau_{FMB}(q_+) - \mu^\tau_{FMB}(q_-) + 2\sum_{k=1}^{m^\prime} c_1^\tau(u_k),\]
        where $m^\prime$ is the number of cascades in $u$ and   $\left(u_k\right)_{k=1}^{m^\prime}$ are the cascades of $u$. \qed
     \end{prop}
     \begin{defin}

    Let $q_-\in\Crit(h_-)$ and $q_+\in\Crit(h_+)$ be a pair of critical points of $h_-$ and $h_+$, respectively.
         \textbf{space of all unparametrized continuation flowlines with cascades connecting $q_-$ to $q_+$} by $\widehat{\cM}^{cont}(q_-,q_+) := \bigcup_{m=0}^\infty \widehat{\cM}^{cont}_{\le m}(q_-,q_+)$.
   \end{defin}

    \begin{defin}
        Let $J$ be a time-independent almost complex structure on $M$. A \textbf{rooted bubble tree} consists of the following data \[\mathbf{T}=\left(\left(V,E\right), \left(B_v\right)_{v\in V}, \left(p^-_{e}\right)_{e\in E}, \left(p^+_{e}\right)_{e\in E},v_r\right ),\] such that:
        \begin{enumerate}
            \item $(V,E)$ is a rooted combinatorial tree, namely, an undirected graph without cycles, whose set of vertices is $V$ and its set of edges is $E$, with $v_r\in V$ a distinguished vertex called the root.
            \item For every vertex $v\in V$, $B_v\colon S^2 \to M$ is a $J$-holomorphic sphere.
            \item For every edge $e=(v_1,v_2)\in E$, connecting vertices $v_1$ and $v_2$, $p^+_e,p^-_e\in S^2$ are points, such that for the edge in the reverse direction $\bar{e}:=(v_2,v_1)$, they satisfy:
            \[p^-_e=p^+_{\bar{e}},\text{ and } p^+_e=p^-_{\bar{e}},\] 
            and the maps $B_{v_1},B_{v_2}$ satisfy
            \[B_{v_1}(p^-_e)=B_{v_2}(p^+_e).\]
            \end{enumerate}
       We denote by $r(\mathbf{T})$ the root of the rooted bubble tree $\mathbf{T}$, namely the vertex $v_r$, and by $B_{r(\mathbf{T})}$ the map $B_{v_r}$. For every $v\in V$, the holomorphic sphere $B_v$ is called a \textbf{bubble of} $\mathbf{T}$.
    \end{defin}

    \begin{defin}
    
    Let $S^- \subset \bS_{H_-}$ and $S^+ \subset \bS_{H_+}$ be corresponding critical submanifolds. A \textbf{bubbled continuation solution} connecting $S^-$ to $S^+$ is a tuple $\mathbf{u} = \left(u, (\mathbf{T}_k)_{k=1}^a, (p_k)_{k=1}^a\right)$ where $u \in \cM^{cont}(S^-, S^+)$ and either:
    \begin{itemize}
        \item $a = 0$ and the collections $(\mathbf{T}_k)$ and $(p_k)$ are empty; or
        \item $a \in \N$ and $p_k \in S^1 \times \R$ are points such that for all $1 \le k \le a$: 
        $$ u(p_k) \in \im B_{r(\mathbf{T}_k)}. $$
    \end{itemize}
    In the latter case, for every $1 \le k \le a$, the rooted bubble tree is holomorphic with respect to the almost complex structure $J_{p_k}$ on $M$ obtained by evaluating $\cJ$ at the point $p_k \in S^1 \times \R$.

Furthermore, we say that the bubbled continuation solution $\mathbf{u}$ converges at $-\infty$ to $\lim\limits_{s \to -\infty} u(s, \cdot)$, and we denote this limit by $\lim\limits_{s \to -\infty} \mathbf{u}(s, \cdot)$. Similarly, we say that $\mathbf{u}$ converges at $+\infty$ to $\lim\limits_{s \to +\infty} u(s, \cdot)$ and denote this limit by $\lim\limits_{s \to +\infty} \mathbf{u}(s, \cdot)$.

    The space of all bubbled continuation solutions connecting $S^-$ to $S^+$ is denoted by $\cM^{cont, bbl}(S^-, S^+; \cH, \cJ)$. 

    If $H_- = H_+$, $J_- = J_+$, and the homotopies $\cH, \cJ$ are constant (i.e., $H_s = H_-$ and $J_s = J_-$ for all $s \in \R$), then any bubbled continuation solution is called a \textbf{bubbled Floer solution} and its underlying continuation solution is called the \textbf{underlying Floer solution}. 

In this case, the space of all bubbled Floer solutions connecting $S^-$ to $S^+$ is denoted by $\cM^{bbl}(S^-, S^+; H, J)$, where $H=H_-$ and $J=J_-$.

\end{defin}

\begin{rem}
    As usual, when $\cH$ and $\cJ$ are clear from the context, we abbreviate these spaces to $\cM^{cont, bbl}(S^-, S^+)$ and $\cM^{bbl}(S^-, S^+)$, respectively.
\end{rem}

\begin{defin}

    A \textbf{broken bubbled continuation solution} for $(\cH, \cJ)$ is a tuple $\left((\mathbf{u}_k)_{k=1}^b, j\right)$, with $b, j \in \N$ and $1 \le j \le b$, such that for every $1 \le k \le b$:
    \begin{enumerate}
        \item if $1 \le k \le j-1$, then $\mathbf{u}_k$ is a bubbled Floer solution for $(H_-, J_-)$;
        \item if $k=j$, then $\mathbf{u}_k$ is a bubbled continuation solution for $(\cH, \cJ)$;
        \item if $j+1 \le k \le b$, then $\mathbf{u}_k$ is a bubbled Floer solution for $(H_+, J_+)$.
    \end{enumerate}
    Additionally, the bubbled solutions $(\mathbf{u}_k)_{k=1}^b$ satisfy:
    $$ \lim_{s \to \infty} \mathbf{u}_k(s, \cdot) = \lim_{s \to -\infty} \mathbf{u}_{k+1}(s, \cdot), $$
    for all $1 \le k \le b-1$.

Furthermore, we say that the broken bubbled continuation solution $\mathbf{u} = ((\mathbf{u}_k)_{k=1}^b, j)$ converges at $-\infty$ to $\lim\limits_{s \to -\infty} \mathbf{u}_1(s, \cdot)$, and we denote this limit by $\lim\limits_{s \to -\infty} \mathbf{u}(s, \cdot)$. Similarly, we say that $\mathbf{u}$ converges at $+\infty$ to $\lim\limits_{s \to +\infty} \mathbf{u}_b(s, \cdot)$ and denote this limit by $\lim\limits_{s \to +\infty} \mathbf{u}(s, \cdot)$. Moreover, let $S^- \subset \bS_{H_-}$ and $S^+ \subset \bS_{H_+}$ be critical submanifolds such that $\lim\limits_{s \to -\infty} \mathbf{u}(s, \cdot) \in S^-$ and $\lim\limits_{s \to +\infty} \mathbf{u}(s, \cdot) \in S^+$. We then say that $\mathbf{u}$ connects $S^-$ to $S^+$.

    If $H_- = H_+$, $J_- = J_+$, and the homotopies $\cH, \cJ$ are constant (i.e., $H_s = H_-$ and $J_s = J_-$ for all $s \in \R$), then any broken bubbled continuation solution is called a \textbf{broken bubbled Floer solution}.
\end{defin}

   \begin{defin}\label{def:bubbled floer flowline with cascades} 
   
   Let $S^-,S^+\subset \bS_H$ be critical submanifolds and let $q_-,q_+\in\Crit(h)$ be a pair of critical points of $h$, lying on $S^-$ and $S^+$ respectively.
       \begin{enumerate}
        \item A \textbf{bubbled Floer flowline with zero cascades} from $q_-$ to $q_+$, $\mathbf u = (\gamma)$,  is a Floer flowline with zero cascades, that is a gradient flowline of $h$, i.e. a smooth map $\gamma \fc \R \to \mathbf{S}$, such that 
           $$  \dot\gamma=\nabla_gh\circ \gamma,\qquad\text{and}\qquad
            \lim_{t\to\pm\infty}\gamma(t) = q_\pm.$$
            Note that in this case $S^- = S^+$.
        
        \item A \textbf{bubbled Floer flowline with} $\mathbf{1}$ \textbf{cascade} from $q_-$ to $q_+$ consists of data
            \[\mathbf{u}=\left( \left(S^-,S^+\right),\left(u\right),\left(\gamma_1,\gamma_2\right)\right),\]
            where
            \begin{enumerate}
                \item $u\in \cM^{bbl}\left(S^-,S^{+}\right)$, a bubbled Floer solution, non-constant in the $\R$ variable.
                
                \item $\gamma_1\fc (-\infty,0] \to \bS_H$ and $\gamma_{2}\fc [0,+\infty) \to \bS_H$, are all gradient flowlines of $h$, namely $\dot{\gamma}_1 = \nabla_gh\circ\gamma_1$ and $\dot{\gamma}_2 = \nabla_gh\circ\gamma_2$, such that:
                \begin{enumerate}
                    \item $q_- =\lim\limits_{t\to-\infty} \gamma_1(t)$, and $\gamma_1(0) = \lim\limits_{s\to-\infty} u(s,\cdot)$.
                    \item $\lim\limits_{s\to+\infty} u(s,\cdot) = \gamma_{2}(0)$, and $\lim\limits_{t\to+\infty} \gamma_{2}(t) = q_+$.
                \end{enumerate}
            \end{enumerate}

                    \item Let $m\in\N$ with $m\geq2$. A \textbf{bubbled Floer flowline with} $\mathbf{m}$ \textbf{cascades} from $q_-$ to $q_+$ consists of data
            \[\mathbf{u}=\left( \left(S^k\right)_{k=1}^{m+1},\left(u_k\right)_{k=1}^m,\left(t_k\right)_{k=2}^{m},\left(\gamma_k\right)_{k=1}^{m+1}\right),\]
            where
            \begin{enumerate}
                \item each $S^k$ is a connected component of $\bS_H$,
                \item $S^1 = S^-$ and $S^+ = S^{m+1}$.
                \item $u_k \in \cM^{bbl}\left(S^k,S^{k+1}\right)$, a bubbled Floer solution, non-constant in the $\R$ variable.
                \item $t_2,\ldots,t_m\in (0,+\infty)$ are positive numbers.
                \item $\gamma_1\fc (-\infty,0] \to \bS_H$, $\gamma_k\fc [0,t_k] \to \bS_H$, for $2\le k \le m$ and $\gamma_{m+1}\fc [0,+\infty) \to \bS_H$, are all gradient flowlines of $h$, namely $\dot{\gamma}_k = \nabla_gh\circ\gamma_k$ for $1\le k \le m+1$, such that:
                \begin{enumerate}
                    \item $q_- =\lim\limits_{t\to-\infty} \gamma_1(t)$, and $\gamma_1(0) = \lim\limits_{s\to-\infty} u_1(s,\cdot)$.
                    \item For all $2\le k \le m$, $\lim\limits_{s\to\infty} u_{k-1}(s,\cdot) = \gamma_k(0)$ and $\gamma_k(t_k)=\lim\limits_{s\to-\infty} u_{k}(s,\cdot)$.
                    \item $\lim\limits_{s\to+\infty} u_{m}(s,\cdot) = \gamma_{m+1}(0)$, and $\lim\limits_{t\to+\infty} \gamma_{m+1}(t) = q_+$.
                \end{enumerate}
            \end{enumerate}
            
       \end{enumerate}

Furthermore, if $\mathbf{u}$ is a bubbled Floer flowline cascades from $q_-$ to $q_+$ we say that $\mathbf{u}$ converges at $-\infty$ to $q_-$, and denote  $\lim\limits_{s \to -\infty} \mathbf{u}(s, \cdot)=q_-$. Similarly, we say that $\mathbf{u}$ converges at $+\infty$ to $q_+$ and denote $\lim\limits_{s \to +\infty} \mathbf{u}(s, \cdot)=q_+$. 
    
   \end{defin}

    \begin{defin}

    Let $q_-\in\Crit(h_-)$ and $q_+\in\Crit(h_+)$ be two critical points lying on two critical submanifolds $S^-\subseteq \bS_{H_-}$ and $S^+\subseteq\bS_{H_+}$ respectively.
       \begin{enumerate}

        \item A \textbf{bubbled continuation flowline with cascades of type $(0;0)$} from $q_-$ to $q_+$
         consists of data
            \[\left(\left( \left(S^-\right),\left(\gamma_1^-\right)\right),(w),\left( \left(S^+\right),\left(\gamma_1^+\right)\right)\right)\]
            where
            \begin{enumerate}
               
                \item $w \in \cM^{cont,bbl}\left(S^-,S^+;\cH,\cJ\right)$, a bubbled continuation solution.
                \item $\gamma^-_1\fc (-\infty,0] \to S^-\subset \bS_{H_-}$ is a  gradient flowline of $h_-$ and $\gamma^+_1\fc [0,+\infty) \to S^+\subset \bS_{H_+}$ is a  gradient flowline of $h_+$.  They are required to satisfy:
                \begin{enumerate}
                    \item $q_- =\lim\limits_{t\to-\infty} \gamma^-_1(t)$, and $\gamma^-_1(0) = \lim\limits_{s\to-\infty} w(s,\cdot)$.
                    \item 
                    $\lim\limits_{s\to+\infty}w(s,\cdot)=\gamma^+_1(0)$, and $\lim\limits_{t\to+\infty} \gamma^+_{1}(t) = q_+$.
                \end{enumerate}
            \end{enumerate}
The \textbf{ordered appearance} of this bubbled continuation flowline with cascades is $$((w_i)_{i=1}^1,(\gamma_i)_{i=1}^2),$$
where $w_1=w$, $\gamma_1=\gamma_1^-$ and $\gamma_2=\gamma_1^+$.

             \item Let $m_1\in \N$. A \textbf{buubled continuation flowline with cascades of type $(m_1;0)$} from $q_-$ to $q_+$
         consists of data
            \[\left(\left( \left(S_-^k\right)_{k=1}^{m_1+1},\left(u^-_k\right)_{k=1}^{m_1},\left(t^-_k\right)_{k=2}^{m_1+1},\left(\gamma^-_k\right)_{k=1}^{m_1+1}\right),(w),\left( \left(S^+\right),\left(\gamma^+_1\right)\right)\right)\]
            where
            \begin{enumerate}
                \item For every $1\leq k \leq m_1+1 $, $S_-^k$ is a connected component of $\bS_{H_-}$, moreover, $S_-^1 = S^-$.
                \item For every $1\leq k\leq m_1$, $u^-_k \in \cM^{bbl}\left(S_-^k,S_-^{k+1};H_-,J_-\right)$ is a bubbled Floer solution.
                \item $w \in \cM^{cont,bbl}\left(S^{m_1+1}_-,S^+;\cH,\cJ\right)$, a bubbled continuation solution.
                \item $\gamma^-_1\fc (-\infty,0] \to S^-\subset\bS_{H_-}$, $\gamma^-_k\fc [0,t^-_k] \to \bS_{H_-}$, for $2\le k \le m_1+1$, are all gradient flowlines of $h_-$. Similarly, $\gamma^+_1\fc [0,+\infty) \to S^+\subset\bS_{H_+}$, is a gradient flowline of $h_+$.
                They are required to satisfy:
                \begin{enumerate}
                    \item $q_- =\lim\limits_{t\to-\infty} \gamma^-_1(t)$, and $\gamma^-_1(0) = \lim\limits_{s\to-\infty} u^-_1(s,\cdot)$.
                    \item For all $2\le k \le m_1$, $\lim\limits_{s\to+\infty} u^-_{k-1}(s,\cdot) = \gamma^-_k(0)$ and $\gamma^-_k(t^-_k)=\lim\limits_{s\to-\infty} u^-_{k}(s,\cdot)$.
                    \item $\lim\limits_{s\to+\infty}u^-_{m_1}(s,\cdot)=\gamma^-_{m_1+1}(0)$ and $\gamma^-_{m_1+1}(t^-_{m_1+1})= \lim\limits_{s\to-\infty}w(s,\cdot).$
                    \item 
                    $\lim\limits_{s\to+\infty}w(s,\cdot)=\gamma^+_1(0)$, and $\lim\limits_{t\to+\infty} \gamma^+_{1}(t) = q_+$.
                \end{enumerate}
            \end{enumerate}

    The \textbf{ordered appearance} of this bubbled continuation flowline with cascades is $$\left((w_k)_{k=1}^{m_1+1},(t_k)_{k=2}^{m_1+1},(\gamma_k)_{k=1}^{m_1+2}\right)$$
where 
\begin{itemize}
    \item For every $1\leq k\leq m_1+1$, if $k=m_1+1$ then $w_k=w$, otherwise $w_k=u_k^-$;
    \item For every $2\leq k\leq m_1+1$ set $t_k=t_k^+$;
    \item For every $1\leq k\leq m_1+2$, if $k=m_1+2$ then $\gamma_k=\gamma_1^+$, otherwise $\gamma_k=\gamma_k^-$.
\end{itemize}

        \item Let $m_2\in \N$. A \textbf{bubbled continuation flowline with cascades of type $(0;m_2)$} from $q_-$ to $q_+$
         consists of data
            \[\left(\left( \left(S^-\right),\left(\gamma_1^-\right)\right),(w),\left( \left(S_+^k\right)_{k=1}^{m_2+1},\left(u^+_k\right)_{k=1}^{m_2},\left(t^+_k\right)_{k=1}^{m_2},\left(\gamma^+_k\right)_{k=1}^{m_2+1}\right)\right)\]
            where
            \begin{enumerate}
                \item For every $1\leq k\leq m_2+1$, $S_+^k$ is a connected component of $\bS_{H_+}$, moreover, $S^+=S_+^{m_2+1}$.
                \item For every $1\leq k\leq m_2$, $u^+_k \in \cM^{bbl}\left(S_+^k,S_+^{k+1};H_+,J_+\right)$ is a bubbled Floer solution.
                \item $w \in \cM^{cont.bbl}\left(S^{-},S_+^{1};\cH,\cJ\right)$, a bubbled continuation solution.
                \item $\gamma^-_1\fc (-\infty,0] \to S^-\subset \bS_{H_-}$ is a gradient flowline of $h_-$. Additionally, $\gamma^+_k\fc[0,t^+_k] \to \bS_{H_+}$, for $1\le k \le m_2$, and $\gamma^+_{m_2+1}\fc [0,+\infty) \to S^+\subset\bS_{H_+}$, are all gradient flowlines of $h_+$.  They are required to satisfy:
                \begin{enumerate}
                    \item $q_- =\lim\limits_{t\to-\infty} \gamma^-_1(t)$, and $\gamma^-_1(0) = \lim\limits_{s\to-\infty} w(s,\cdot)$.
                    \item 
                    $\lim\limits_{s\to+\infty}w(s,\cdot)=\gamma^+_1(0)$, and $\gamma^+_1(t_1)=\lim\limits_{s\to-\infty}u^+_1(s,\cdot)$.
                    \item For all $2\le k \le m_2$, $\lim\limits_{s\to+\infty} u^+_{k-1}(s,\cdot) = \gamma^+_k(0)$ and $\gamma^+_k(t^+_k)=\lim\limits_{s\to-\infty} u^+_{k}(s,\cdot)$.
                    
                    \item$\lim\limits_{s\to+\infty} u_{m_2}(s,\cdot) = \gamma^+_{m_2+1}(0)$, and $\lim\limits_{t\to+\infty} \gamma^+_{m_2+1}(t) = q_+$.
                \end{enumerate}
            \end{enumerate}

    The \textbf{ordered appearance} of this bubbled continuation flowline with cascades is $$\left((w_k)_{k=1}^{m_2+1},(t_k)_{k=2}^{m_2+1},(\gamma_k)_{k=1}^{m_2+2}\right)$$
where 
\begin{itemize}
    \item For every $1\leq k\leq m_2+1$, if $k=1$ then $w_k=w$, otherwise $w_k=u_{k-1}^+$;
    \item For every $2\leq k\leq m_2+1$ set $t_k=t_{k-1}^+$;
    \item For every $1\leq k\leq m_2+2$, if $k=1$ then $\gamma_k=\gamma_1^-$, otherwise $\gamma_k=\gamma_{k-1}^+$.
\end{itemize}

             \item Let $m_1,m_2\in \N$. A \textbf{bubbled continuation flowline with cascades of type $(m_1;m_2)$} from $q_-$ to $q_+$
         consists of data
            \[\left(\left( \left(S_-^k\right)_{k=1}^{m_1+1},\left(u^-_k\right)_{k=1}^{m_1},\left(t^-_k\right)_{k=2}^{m_1+1},\left(\gamma^-_k\right)_{k=1}^{m_1+1}\right),(w),\left( \left(S_+^k\right)_{k=1}^{m_2+1},\left(u^+_k\right)_{k=1}^{m_2},\left(t^+_k\right)_{k=1}^{m_2},\left(\gamma^+_k\right)_{k=1}^{m_2+1}  \right)\right)\]
            where
            \begin{enumerate}
                \item For every $1\leq k\leq m_1+1$,  $S_-^k$ is a connected component of $\bS_{H_-}$, moreover $S_-^1=S^-$.
                \item For every $1\leq k\leq m_2+1$, $S_+^k$ is a connected component of $\bS_{H_+}$, moreover,  $S^+=S_+^{m_2+1}$.
                \item For every $1\leq k\leq m_1$, $u^-_k \in \cM^{bbl}\left(S_-^k,S_-^{k+1};H_-,J_-\right)$ is a bubbled Floer solution.
                \item For every $1\leq k\leq m_2$,  $u^+_k \in \cM^{bbl}\left(S_+^k,S_+^{k+1};H_+,J_+\right)$ is bubbled Floer solution.
                \item $w \in \cM^{cont,bbl}\left(S^{m_1+1}_-,S_+^{1};\cH,\cJ\right)$, a bubbled continuation solution.
                \item $\gamma^-_1\fc (-\infty,0] \to S^-\subset\bS_{H_-}$, $\gamma^-_k\fc [0,t^-_k] \to \bS_{H_-}$, for $2\le k \le m_1+1$, are all gradient flowlines of $h_-$. Similarly,   $\gamma^+_k\fc [0,t^+_k] \to \bS_{H_+}$, for $1\le k \le m_2$, and $\gamma^+_{m_2+1}\fc [0,+\infty) \to S^+\subset\bS_{H_+}$, are all gradient flowlines of $h_+$. They are required to satisfy:
                \begin{enumerate}
                    \item $q_- =\lim\limits_{t\to-\infty} \gamma^-_1(t)$, and $\gamma^-_1(0) = \lim\limits_{s\to-\infty} u^-_1(s,\cdot)$.
                    \item For all $2\le k \le m_1$, $\lim\limits_{s\to+\infty} u^-_{k-1}(s,\cdot) = \gamma^-_k(0)$ and $\gamma^-_k(t^-_k)=\lim\limits_{s\to-\infty} u^-_{k}(s,\cdot)$.
                    \item $\lim\limits_{s\to+\infty}u^-_{m_1}(s,\cdot)=\gamma^-_{m_1+1}(0)$ and $\gamma^-_{m_1+1}(t^-_{m_1+1})= \lim\limits_{s\to-\infty}w(s,\cdot).$
                    \item 
                    $\lim\limits_{s\to+\infty}w(s,\cdot)=\gamma^+_1(0)$, and $\gamma^+_1(t_1)=\lim\limits_{s\to-\infty}u^+_1(s,\cdot)$.
                    \item For all $2\le k \le m_2$, $\lim\limits_{s\to+\infty} u^+_{k-1}(s,\cdot) = \gamma^+_k(0)$ and $\gamma^+_k(t^+_k)=\lim\limits_{s\to-\infty} u^+_{k}(s,\cdot)$.\item$\lim\limits_{s\to+\infty} u_{m_2}(s,\cdot) = \gamma^+_{m_2+1}(0)$, and $\lim\limits_{t\to+\infty} \gamma^+_{m_2+1}(t) = q_+$.
                \end{enumerate}
            \end{enumerate}
           
    The \textbf{ordered appearance} of this bubbled continuation flowline with cascades is $$\left((w_k)_{k=1}^{m_1+m_2+1},(t_k)_{k=2}^{m_1+m_2+1},(\gamma_k)_{k=1}^{m_1+m_2+2}\right)$$
where 
\begin{itemize}
    \item For every $1\leq k\leq m_1+m_2+1$, if $k\leq m_1$ then $w_k=u_k^-$, if $k=m_1+1$ then $w_k=w$, otherwise $w_k=u_{k-m_1-1}^+$;
    \item For every $2\leq k\leq m_1+m_2+1$, if $k\leq m_1+1$ set $t_k=t_{k}^-$, otherwise $t_k=t_{k-m_1-1}^+$;
    \item For every $1\leq k\leq m_1+m_2+2$, if $k\leq m_1+1$ then $\gamma_k=\gamma_k^-$, otherwise $\gamma_k=\gamma_{k-m_1-1}^+$.
\end{itemize}

            \item Let $m\in \N$. A \textbf{bubbled continuation flowline with $m$ cascades} from $q_-$ to $q_+$ is a continuation flowline of type $(m_1;m_2)$ from $q_-$ to $q_+$ for some $0\le m_1$, $0\le m_2$ such that $m_1 + m_2=m-1$.

       \end{enumerate}

Furthermore, if $\mathbf{u}$ is a bubbled continuation flowline cascades from $q_-$ to $q_+$ we say that $\mathbf{u}$ converges at $-\infty$ to $q_-$, and denote  $\lim\limits_{s \to -\infty} \mathbf{u}(s, \cdot)=q_-$. Similarly, we say that $\mathbf{u}$ converges at $+\infty$ to $q_+$ and denote $\lim\limits_{s \to +\infty} \mathbf{u}(s, \cdot)=q_+$. 
   \end{defin}

\begin{defin}

    A \textbf{broken bubbled continuation flowline} for $(\cH, \cJ)$ is a tuple $\left((\mathbf{u}_k)_{k=1}^b, j\right)$, with $b, j \in \N$ and $1 \le j \le b$, such that for every $1 \le k \le b$:
    \begin{enumerate}
        \item if $1 \le k \le j-1$, then $\mathbf{u}_k$ is a bubbled Floer flowline for $(H_-, J_-)$;
        \item if $k=j$, then $\mathbf{u}_k$ is a bubbled continuation flowline for $(\cH, \cJ)$;
        \item if $j+1 \le k \le b$, then $\mathbf{u}_k$ is a bubbled Floer flowline for $(H_+, J_+)$.
    \end{enumerate}
    Additionally, the bubbled flowlines $(\mathbf{u}_k)_{k=1}^b$ satisfy:
    $$ \lim_{s \to \infty} \mathbf{u}_k(s, \cdot) = \lim_{s \to -\infty} \mathbf{u}_{k+1}(s, \cdot), $$
    for all $1 \le k \le b-1$.

Furthermore, we say that the broken bubbled continuation flowline $\mathbf{u} = ((\mathbf{u}_k)_{k=1}^b, j)$ converges at $-\infty$ to $\lim\limits_{s \to -\infty} \mathbf{u}_1(s, \cdot)$, and we denote this limit by $\lim\limits_{s \to -\infty} \mathbf{u}(s, \cdot)$. Similarly, we say that $\mathbf{u}$ converges at $+\infty$ to $\lim\limits_{s \to +\infty} \mathbf{u}_b(s, \cdot)$ and denote this limit by $\lim\limits_{s \to +\infty} \mathbf{u}(s, \cdot)$. Moreover, let $S^- \subset \bS_{H_-}$ and $S^+ \subset \bS_{H_+}$ be critical submanifolds such that $\lim\limits_{s \to -\infty} \mathbf{u}(s, \cdot) \in S^-$ and $\lim\limits_{s \to +\infty} \mathbf{u}(s, \cdot) \in S^+$. We then say that $\mathbf{u}$ connects $S^-$ to $S^+$.

    If $H_- = H_+$, $J_- = J_+$, and the homotopies $\cH, \cJ$ are constant (i.e., $H_s = H_-$ and $J_s = J_-$ for all $s \in \R$), then any broken bubbled continuation flowline is called a \textbf{broken bubbled Floer flowline}.
\end{defin}

 Given a smooth manifold $X$, the natural action of $\R$ on itself induces an $\R$-action on the space of paths in $X$, denoted by $C^\infty(\R,X)$. For every $s_0\in \R$ and $\gamma\in C^\infty(\R,X)$, the action of $s_0$ on $\gamma$, denoted by $s_0\cdot \gamma$, is defined by$(s_0\cdot \gamma)(s) =s_0\cdot \gamma(s)= u(s+s_0)$ for all $s\in \R$. Similarly, the cylinder $\R\times S^1$ admits a natural free $\R$-action, given by $s_0\cdot x=(s+s_0,t)$ for every $s_0\in \R$ and $x=(s,t)\in \R\times S^1$. That action naturally induces an $\R$-action on the space of cylinders in $X$, denoted by $C^\infty(\R\times S^1,X)$. For every $s_0\in \R$ and $u\in C^\infty(\R\times S^1,X)$, the action of $s_0$ on $u$, denoted by $s_0\cdot u$, is defined by$(s_0\cdot u)(s,t) =s_0\cdot u(s,t)= u(s+s_0,t)$ for all $(s,t)\in \R\times S^1$.

Also, let us equip $\R\times S^1$ with the standard flat Riemannian metric, obtained as a pull back of the Euclidean metric on $\R^2$ under to quotient $\R\times \R\to \R\times (\R/\Z)$. Given a point $x\in \R\times S^1$ and a positive number $\epsilon>0$ we denote by $B(x,\epsilon)$ the open ball in $\R\times S^1$ of radius $\epsilon>0$, with respect to this metric, around a point $x$.

 \begin{defin} 
 
    Let $S^-,S^+\subseteq \bS_{H}$ be two critical submanifolds.
    
Let $(u_i)_{i\in \N}$ be a sequence of Floer solutions such that $u_i\in \cM(S^-,S^+;H,J)$. Furthermore, let $(s_i)_{i\in \N}$ be a sequence of real numbers and $\mathbf{w} = \left(w,\left(\mathbf{T}_j\right)_{j=1}^{a},\left(p_j\right)_{j=1}^{a}\right)\in\cM^{bbl}(S^-,S^+)$ be a bubbled Floer solution.

 For every $1\leq j\leq a$ denote by $(B_{j,\alpha})_{\alpha\in I_j}$ the collection of all the bubbles of $\mathbf{T}_j$. We say that $(s_i\cdot u_i)_{i\in \N}$ \textbf{$C^\infty_{loc}$ converges with bubbling to a $\mathbf{w}$}, if for every $\alpha\in I$ and $i\in \N$ there exists a point $x_i^\alpha=(s_i^\alpha,t_i^\alpha)\in \R\times S^1$ and a positive number $\epsilon_i^\alpha>0$ such that the following holds:
    \begin{enumerate}
        \item On $S^1 \times \R \setminus \left\lbrace p_1, \ldots, p_a\right\rbrace$ the sequence $( s_i\cdot u_i)_{i\in \N}$ convergence $C^\infty_{loc}$ to $w$.
        \item For every $1\leq j\leq a$ and $\alpha\in I_j$ we have $\lim\limits_{i\to\infty} x_i^\alpha=p_j$ and $\lim\limits_{i\to\infty} \epsilon_i^\alpha=0$;
     
        \item For every $1\leq j\leq a$ and $\alpha\in I_j$, we sequence $((s_i\cdot u_i)\vert^{}_{B(x_i^\alpha,\epsilon_i^\alpha)})_{i\in \N}$ convergence to the bubble $B_{j,\alpha}$ after appropriate rescaling and extension to the ideal point at infinity. 
        \item $\lim\limits_{i\to\infty} E(u_i) = E(w) + \sum\limits_{j=1}^a\sum\limits_{\alpha \in I_j}\omega(B_\alpha).$
    \end{enumerate}

    \end{defin}

 \begin{defin}

    Let $S^-\subseteq \bS_{H_-}$ and $S^+\subseteq\bS_{H_+}$ be two critical submanifolds. Let $(u_i)_{i\in \N}$ be a sequence of continuation solutions such that $u_i\in \cM^{cont}(S^-,S^+;\cH,\cJ)$. Furthermore, let $(s_i)_{i\in \N}$ be a sequence of real numbers and $\mathbf{w} = \left(w,\left(\mathbf{T}_j\right)_{j=1}^{a},\left(p_j\right)_{j=1}^{a}\right)$ be a bubbled solution such that one of the following hold: 
\begin{enumerate}
    \item $s_i=0$ for every $i\in \N$, and $\mathbf{w}\in \cM^{cont,bbl}(S^-,S^+;\cH,\cJ)$;
    \item $\lim\limits_{i\to+\infty}s_i=+\infty$ and $\mathbf{w}\in \cM^{bbl}(S,S^+;H_+,J_+)$, for some critical submanifold $S\subseteq\bS_{H_+}$ for $H_+$;
    \item $\lim\limits_{i\to+\infty}s_i=-\infty$ and $\mathbf{w}\in \cM^{bbl}(S^-,S;H_-,J_-)$, for some critical submanifold $S\subseteq\bS_{H_-}$ for $H_-$;
\end{enumerate}

 For every $1\leq j\leq a$ denote by $(B_{j,\alpha})_{\alpha\in I_j}$ the collection of all the bubbles of $\mathbf{T}_j$. We say that $(s_i\cdot u_i)_{i\in \N}$ \textbf{$C^\infty_{loc}$ converges with bubbling to a $\mathbf{w}$}, if for every $\alpha\in I$ and $i\in \N$ there exists a point $x_i^\alpha=(s_i^\alpha,t_i^\alpha)\in \R\times S^1$ and a positive number $\epsilon_i^\alpha>0$ such that the following holds:
    \begin{enumerate}
        \item On $S^1 \times \R \setminus \left\lbrace p_1, \ldots, p_a\right\rbrace$ the sequence $( s_i\cdot u_i)_{i\in \N}$ convergence $C^\infty_{loc}$ to $w$.
        \item For every $1\leq j\leq a$ and $\alpha\in I_j$ we have $\lim\limits_{i\to\infty} x_i^\alpha=p_j$ and $\lim\limits_{i\to\infty} \epsilon_i^\alpha=0$;
     
        \item For every $1\leq j\leq a$ and $\alpha\in I_j$, we sequence $((s_i\cdot u_i)\vert^{}_{B(x_i^\alpha,\epsilon_i^\alpha)})_{i\in \N}$ convergence to the bubble $B_{j,\alpha}$ after appropriate rescaling and extension to the ideal point at infinity. 
        \item $\lim\limits_{i\to\infty} E(u_i) = E(w) + \sum\limits_{j=1}^a\sum\limits_{\alpha \in I_j}\omega(B_\alpha).$
    \end{enumerate}

    \end{defin}

\begin{defin}

    Let $m\in \Z_{\geq0}$ and let $(\mathbf{u}_i)_{i\in \N}$ be a sequence of continuation flowlines with $m$ cascades. Let $\mathbf{w}=(\mathbf{w}_j)_{j=1}^b$ be a broken bubbled continuation flowline with cascades, where each $\mathbf{w}_j$ is a flowline with $m_j$ cascades, for some $m_j\in \Z_{\geq0}$. For every $1\leq j\leq b$, if $\mathbf{w}_j$ has zero cascades, then we write $\mathbf{w}_j=(\gamma^j)$ where $\gamma^j\fc \R\to \bS_{H_-}\cup \bS_{H_+}$ is a gradient flow line, otherwise we write $$\mathbf{w}_j= \left(\left(w^j_k\right)_{k=1}^{m_j}, \left(\tau^j_k\right)_{k=2}^{m_j},\left(\bar\gamma^j_k\right)_{k=1}^{m_j+1}\right).$$

    We say that $(\mathbf{u}_i)_{i\in \N}$ \textbf{converges to the broken bubbled continuation flowline with cascades} $\mathbf{w}$, if the following conditions hold:

    \begin{enumerate}
        \item In case when $\mathcal{H}$ and $\mathcal{J}$ are constant, and for every $1\leq j\leq b$, $\mathbf{w}_j$ has zero cascades, namely $\mathbf{w}_j = \left(\gamma^j\right)$, then for every $i\in \N$, $\mathbf{u}_i$ has zero cascades, namely $(u_i)_{i\in \N} = (\gamma_i)_{i\in \N}$ where  $\gamma_i\fc \R \to \mathbf{S}_H$ are  gradient flowlines, for every $i\in \N$. Here, we require that there exist $b$ sequences of real numbers $\left( s^j_i\right)_{i=1}^\infty$, with $1\le j \le b$, such that for all $1\le j \le b$ the sequence of flowlines $(s_i^j\cdot\gamma_i)_{i\in\N}$, $C^\infty_{loc}$ converges to $\gamma^j$.

        \item In case when $\mathcal{H}$ and $\mathcal{J}$ are nonconstant, or one of the pieces $\mathbf{w}_j$ has at least one cascade, then $m\ge1$ and we write $\mathbf{u}_i$ in its ordered appearance $\mathbf{u}_i=\left(\left(u^i_k\right)_{k=1}^m,\left(t^i_k\right)_{k=2}^m,\left(\gamma^i_k\right)_{k=1}^{m+1}\right)$ for every $i\in \N$. We require that there exist maps $(\Gamma^i)_{i\in \N}$ that encode for each piece in the broken flowline, $\mathbf{w}$, which sequence of cascades Floer-Gromov converges to a broken limit containing this piece. I.e. 
        for every $i\in \N$ there exists a surjective map
\[\Gamma^i\fc \big[1,\sum_{j=1}^b m_j\big]\cap \N \to \left[ 1, m\right]\cap \N,\]
            which is monotone non-decreasing, and also there exist real numbers $s^i_\lambda\in \R$, for each $\lambda\in \big[1,\sum_{j=1}^b m_j\big]\cap \N$, such that the sequence $(s^i_\lambda\cdot u^i_{\Gamma^i(\lambda)})_{i\in \N}$ converges $C^\infty_{loc}$ with bubbling to $w_\lambda$, where $w_\lambda$ denotes $w_k^j$ for $1\leq j\leq b$ and $1\leq k\leq m_j$ satisfying $\lambda = \sum_{p=1}^{j-1} m_p + k$ with $m_j\geq 1$.

            For every $\lambda\in \big[1,\sum_{j=1}^b m_j - 1\big]\cap \N$ set
            \begin{align*}
                &\tau_\lambda = \begin{cases}
                    \tau^j_{k+1} & \text{if } \lambda = \sum_{p=1}^j m_p + k \text{ for } 1\le k \leq m_j - 1, \,\text{where}\,m_j\geq2, \\
                    +\infty  & \text{if } \lambda = \sum_{p=1}^j m_p,\,\text{where}\,m_j\geq 1,
                \end{cases}\\
                \text{and } & \tau^i_\lambda = \begin{cases}
                    t^i_{\Gamma^i(\lambda)+1}, & \text{if } \lambda = \max \left\lbrace \lambda^\prime \in [1,\sum_{p=1}^b m_p - 1]\cap \N \,\mid\, \Gamma^i(\lambda^\prime) = \Gamma^i(\lambda) \right\rbrace \\
                    0,  & \text{otherwise},
                \end{cases}
            \end{align*}
            for every $i\in \N$. We require that $\lim\limits_{i\to\infty} \tau^i_\lambda = \tau_\lambda.$ This ensures that the finite gradient pieces inside each $\mathbf{w}_j$ are approached by gradient pieces of time-length converging to the correct one, and moreover, the breaks between each $\mathbf{w}_j$ are approached by gradient segments of time-length that converges to infinity.

        Moreover, for every $1\leq j\leq b$, if $\mathbf{w}_j=\left(\bar\gamma^j\right)$ is a flowline with zero cascades, that is $\bar\gamma^j\fc \R \to \mathbf{S}_{H_-}\cup \mathbf{S}_{H_+}$ is a gradient flowline, then we require that there exists a sequence of integers $(k_i)_{i\in \N}$ between $1$ and $m$ and a sequences of real numbers $\left( s_i\right)_{i=1}^\infty$ such that 
            $$\lim_{i\to\infty}s_i=-\infty,\qquad\text{and}\qquad\lim_{i\to+\infty} (s_i+t_{k_i}^i)=+\infty,$$
            and in addition, one of the following holds:

            \begin{itemize}
                \item If for every $i\in \N$, we denote by $\gamma_i\fc \R\to \bS_{H_-}\cup \bS_{H_+}$ the unique gradient line flow satisfying $\gamma_i(s_i)=\gamma_{k_i+1}^i(0)$, then the sequence $(\gamma_i)_{i\in \N}$ converges in $C^\infty_{loc}$ to $\bar{\gamma}_j$.
                \item If for every $i\in \N$, we denote by $\gamma_i\fc \R\to \bS_{H_-}\cup \bS_{H_+}$ the unique gradient line flow satisfying $\gamma_i(s_i)=\gamma_{k_i}^i(t_{k_i}^i)$, then the sequence $(\gamma_i)_{i\in \N}$ converges in $C^\infty_{loc}$ to $\bar{\gamma}_j$. Here we denote $t^i_1=0$ for every $i\in \N$.
            \end{itemize}

            \item $\lim\limits_{i\to \infty} E(\mathbf{u}_i) = \sum\limits_{j=1}^b E(\mathbf{w}_j).$

\end{enumerate}
\end{defin}

    \begin{thm}\label{thm:subsequence converging to a broken bubbled flowline with cascades}
         Let $\left(\mathbf{u}_i\right)_{i\in \N}$ sequence of Floer (continuation) flowlines with $m$ cascades. If there exists a number $E>0$ such that for all $i\in \N$, $E(\mathbf{u}_i)<E$, then $\left(\mathbf{u}_i\right)_{i\in \N}$ has a subsequence converging to a broken bubbled Floer (continuation) flowline with cascades.
    \end{thm}
    To provide an explanation for Theorem~\ref{thm:subsequence converging to a broken bubbled flowline with cascades} we recall convergence to broken solution in the Morse--Bott Floer context.
    
    \begin{defin} A sequence of Floer (continuation) solutions $u_i\in\cM(S^-,S^+)$ ( $u_i\in\cM^{cont}(S^-,S^+;H_s)$), is said to \textbf{Floer-Gromov converge to a broken bubbled Floer (continuation) solution} $\left(\left(\mathbf{w}_k\right)_{k=1}^b, (j)\right)$, where each $\mathbf{w}_k$ is given as $\mathbf{w}_k = \left(w^k,\left(\mathbf{T}^k_j\right)_{j=1}^{a^k},\left(p^k_j\right)_{j=1}^{a^k}\right)$, if the following holds:
    \begin{enumerate}
        \item $\lim\limits_{i\to\infty}E_{top}(u_i) = \sum\limits_{k=1}^b E_{top}(w^k) + \sum\limits_{\alpha\in I} \omega(B_\alpha)$, where $\left\lbrace B_\alpha\right\rbrace_{\alpha\in I}$ is the collection of all sphere bubbles in the broken bubbled solution, $I$ is some index set indexing it.
        \item For all $1\le k \le b$ there exists a sequence $s^k_i\in \R$ such that \[u_i( \cdot +s^k_i, \cdot) \xrightarrow[i\to\infty]{} \mathbf{w}_k,\]
        as a $C^\infty_{loc}$ convergence with bubbling.
    \end{enumerate}
    \end{defin}
    
    \begin{rem}\label{rem: FloerGromov convergence properties}
        Moreover, the following properties hold:
        \begin{enumerate}
            \item Let $\tau$ be a trivialization of $\det TM$ along all the $1$-periodic orbits appearing as asymptotes of $w^k$, $1\le k \le b$. Then for $i$ large enough, 
            \[c_1^\tau(u_i) = \sum\limits_{k=1}^b c_1^\tau(w^k) + \sum\limits_{\alpha\in I} c_1(B_\alpha).\]
            \item The images of $u_i$, converge in the Hausdorff distance to the union of the images of $w_k$, $1\le k \le b$ and the spheres $B_\alpha$, $\alpha\in I$.
        \end{enumerate}
    \end{rem}

    \begin{thm}\label{thm:floerGromovCompactness}
        Let $u_i$ be a sequence of a sequence of Floer (continuation) solutions $u_i\in\cM(S^-,S^+)$ ( $u_i\in\cM^{cont}(S^-,S^+;H_s)$), assume that there exists $C>0$ such that $E(u_i)<C$ for all $i$ large enough. Then after passing to a subsequence $i_r$, $u_{i_r}$ Floer-Gromov converges to a broken bubbled Floer (continuation) solution $\left(\left(\mathbf{w}_k\right)_{k=1}^b, (j)\right)$.   
    \end{thm}
    \begin{proof} (Sketch): The argument follows the standard proof in the literature (e.g. \cite[Theorem 9.1.7]{Audin_Damian_2014_Morse_and_Floer} and \cite[Corollary 3.4]{Salamon_1999_notes}), which is based on obtaining a uniform bound on the derivatives of $u_i$ in the complement of a finite set of ($i$-dependent) points in the domain, and applying the Arzela-Ascoli theorem to $u_i$, and its shifts in the $s$-coordinates by sequences of numbers $s^k_i$. The bubble trees are obtained by the standard local rescaling argument, see \cite[Chapter 4]{McDuff_Salamon_J_curves_2012}.
    Each sequence $s^k_i$ is responsible for the convergence to a different piece in the broken limit, and the only component needed to adapt the proof is the following claim.
    \begin{rem}
        A similar statement holds when $(u_i)_i$ is a sequence of solutions for a sequence of Hamiltonians and almost complex structures $(H_{s,i})_i$,  $(J_i)_i$, converging to a limiting pair $H_{s,\infty}$, $J_\infty$.
    \end{rem}
    \begin{thm}
    Let $u\colon S^1\times \R \to M$ be a solution of the Floer or continuation equation, with $E_{top}(u)<\infty$, then there exist critical submanifolds $\Lambda^-,\Lambda^+$ and $1$-periodic orbits $q^-\in \Lambda^-$ and $q^+ \in \Lambda^+$, such that
    \[\lim\limits_{s\to\pm\infty}u(s,\cdot) = q^{\pm}(\cdot),\qquad \lim\limits_{s\to\pm\infty}\left\vert\frac{d}{ds}u(s,\cdot)\right\vert = 0. \tag*{\qedsymbol}\]
    \end{thm}
    A similar theorem appears with proof in the context of Rabinowitz Floer homology, in \cite[Theorem 25]{Fauck_2016_thesis}.
    \end{proof}

    We can now explain the proof of Theorem \ref{thm:subsequence converging to a broken bubbled flowline with cascades}:
    \begin{proof}[Sketch of proof of Theorem \ref{thm:subsequence converging to a broken bubbled flowline with cascades}]
    Let us denote by $(u^n_k)^m_{k=1}$ the cascades appearing in $\mathbf{u}_n$
    Since $E(\mathbf{u}_n) < E$ and the energies of each cascade are nonnegative, we deduce that for all $1\le k \le m$ and for all $n$: $E(u^n_k) < E$. Thus we may apply Theorem \ref{thm:floerGromovCompactness} to obtain subsequences, still denoted $\left(u^n_k\right)_{n=1}^\infty$, each converging to a broken bubbled Floer or continuation solution, $\left(\left(\mathbf{w}^k_i\right)_{i=1}^{b_k}, (j)\right)$.
    
    Now consider the gradient flow segments $\gamma^n_k$ appearing in each $\mathbf{u}_n$. By the standard compactness arguments for gradient flow solutions, they converge either to finite time segments, or to a, possibly broken, infinite time gradient flow solution.

    All the resulting broken pieces can be compartmentalized into a broken bubbled flowline with cascades, where each piece begins at a gradient flowline of infinite time, follows along cascades and finite time gradient segments, and ends at the next gradient flowline of infinite time.
    \end{proof}

    \subsection{Energy bounds and bounds on the number of cascades}\label{ss: bounds for energy and cascades}
    We show that the Floer and continuation flowlines that participate in the count for the Floer differential operator, continuation morphisms and chain homotopies, all have uniformly bounded energy and a uniformly bounded number of cascades, allowing us later to apply compactness arguments for them.
    First we prove an energy bound. When interpreting the following proposition, take $N=1$ for  flowlines in the chain differential, $N=0$ for  flowlines in the chain continuation map and $N=-1$ for flowlines in the chain homotopy.

    \begin{prop}\label{prop:universal energy bound}
         Let $(M,\omega)$ be a closed monotone symplectic manifold. Let $H_-$ and $H_+$ be Hamiltonians on $(M,\omega)$
         satisfying the \textbf{MB} condition, such that each Hamiltonian is constant along each of its connected critical submanifolds, and let $\cH=(H_s)_{s\in \R}$ be a monotone homotopy of Hamiltonians satisfying $H_s = H_-$ for every $s\leq 0$, and $H_s = H_+$ for $s\geq 1$. Let $h_-$ and $h_+$ be two Morse functions on $\bS_{H_-}$ and $\bS_{H_+}$, respectively. Let $g_-$ and $g_+$ be two Riemannian metrics on $\bS_{H_-}$ and $\bS_{H_+}$, respectively. 

         For every integer $N\in \Z$ there exists a constant $E\geq0$ such that for every pair of almost complex structures $J_-$ and $J_+$ on $M$ compatible with $\omega$, possibly time-dependent, for every homotopy of compatible almost complex structures $\cJ=(J_s)_{s\in \R}$, satisfying $J_s = J_-$  for every $s\leq 0$ and $J_s = J_+$ for every $s\geq1$, for every critical points $q_-\in \Crit(h_-)$ and $q_+\in \Crit(h_+)$, for every $m\in \N$, for every a continuation flowline $\mathbf{u}$ with $m$ cascades $(u_k)_{k=1}^m$ and for every trivialization $\tau$ of $\det TM$ along the asymptotes of the cascades, if
         $$m+N = \mu^\tau_{FMB}\left(q^+;H_+\right) - \mu^\tau_{FMB}\left(q^-;H_-\right) + 2\sum_{k=1}^m c^\tau_1\left(u_k\right) + m,$$
    then the topological energy of $\mathbf{u}$ is bounded from above by $E$, that is,  
    $$E_{top}(\mathbf{u})\leq E.$$
   \end{prop}

    \begin{proof}
        Let $N \in \Z$. Since $H_-$ and $H_+$ satisfy the \textbf{MB} condition and $M$ is a closed symplectic manifold, it follows that $\bS_{H_-}$ and $\bS_{H_+}$ are compact. Therefore, the sets of critical points $\Crit(h_-)$ and $\Crit(h_+)$ of $h_-$ and $h_+$, respectively, are finite. For every $1$-periodic orbit $\alpha \in S_{H_-} \cup S_{H_+}$, let us choose a capping $\hat{\alpha}$, in a manner that is compatible with index and action along every connected component of $S_{H_-} \cup S_{H_+}$. This is done by picking an arbitrary orbit at each connected component, choosing a capping for it, and extending a capping to any other orbit along the component, but choosing paths from it to any other orbit, lifting them to cylinders, and attaching them to the chosen capping.
        Denote by $\tau_0$ the trivialization of $\det TM$ along all $1$-periodic orbits corresponding to critical points of either $h_-$ or $h_+$, induced by the chosen cappings.
        
        The constant $E$ that will be considered in this proof is defined as the sum of three constants. Each of these three constants is well-defined due to the fact that $\Crit(h_-)$ and $\Crit(h_+)$ are finite sets. The constants are the following:

\begin{itemize}
   
    \item The first constant $K_1$ is defined as follows:
    $$K_1 = \left| \max_{q \in \Crit(h_+)} \int_{S^1} H_+ \circ q(t) \, dt - \min_{p \in \Crit(h_-)} \int_{S^1} H_- \circ p(t) \, dt \right|.$$
    \item The second constant $K_2$ is defined as follows:
    $$K_2 = \left| \max_{q \in \Crit(h_+)} \omega(\hat{q}) - \min_{p \in \Crit(h_-)} \omega(\hat{p}) \right|.$$
     \item The third constant $K_3$ is defined as follows:
    $$K_3 = \frac{|\kappa|}{2} \left| N+\max_{p \in \Crit(h_-)} \mu^{\tau_0}_{FMB}(p; H_-) - \min_{q \in \Crit(h_+)} \mu^{\tau_0}_{FMB}(q; H_+) \right|.$$
    Here $\kappa$ stands for the monotonicity constant of $(M, \omega)$.
\end{itemize}

Consider $E = K_1 + K_2 + K_3$. Let $J_-$ and $J_+$ be almost complex structures on $M$ compatible with $\omega$, possibly time-dependent, and let $\cJ=(J_s)_{s\in \R}$ be a homotopy of compatible almost complex structures, satisfying $J_s = J_-$  for every $s\leq 0$ and $J_s = J_+$ for every $s\geq1$. Let $q_- \in \Crit(h_-)$ and $q_+ \in \Crit(h_+)$ be critical points, let $m \in \N$ be an integer, let $\mathbf{u}$ be a continuation flowline with $m$ cascades $(u_k)_{k=1}^m$, and let $\tau$ be a trivialization of $\det TM$ along the asymptotes of the cascades. Assume that
$$ m + N = \mu^\tau_{FMB}(q^+; H_+) - \mu^\tau_{FMB}(q^-; H_-) + 2 \sum_{k=1}^m c^\tau_1(u_k) + m . $$
 First we show that
\begin{equation}\label{eq: energy of flowline with cascades}
    E_{top}(\mathbf{u}) = \int_{S^1} H_+\circ q^+(t) \, dt - \int_{S^1} H_-\circ q^-(t) \, dt + \sum_{k=1}^m \omega(u_k).
\end{equation}

For every $1 \le k \le m$, denote by $q_k^-$ and $q_k^+$ the asymptotes of $u_k$ at $-\infty$ and $+\infty$, respectively. Additionally, for every $1 \le k \le m$, we define the Hamiltonians $H_k^-$ and $H_k^+$ to be $H_+$ or $H_-$ according to the following rules:
\begin{itemize}
    \item If $u_k$ is a Floer solution for $(H_-, J_-)$, then we denote $H_k^- = H_k^+ = H_-$.
    \item If $u_k$ is a Floer solution for $(H_+, J_+)$, then we denote $H_k^- = H_k^+ = H_+$.
    \item If $u_k$ is a continuation solution for $(\cH, \cJ)$, then we denote $H_k^- = H_-$ and $H_k^+ = H_+$.
\end{itemize}

Note that $H_1^-=H_-$ and $H_m^+=H_+$, and that the pairs $q_1^-$, $q^-$, and $q^+_m$, $q^+$, belong to the same critical submanifolds, each, respectively. Additionally, for every $1\leq k\leq m-1$ we have $H_k^+=H_{k+1}^-$ and the $1$-periodic orbits $q_k^+$ and $q_{k+1}^-$ belong to the same critical submanifold of $H_k^+=H_{k+1}^-$, and hence for every $t\in S^1$ we have
$$H_k^+\circ q_k^+(t)=H_{k+1}^-\circ q_{k+1}^-(t).$$
Recall that for every $1\leq k\leq m$ the topological energy of $u_k$ is given by 
$$E_{top}(u_k)=\int_{S^1}H_k^+\circ q_k^+(t)\,dt-\int_{S^1} H_k^-\circ q_k^-(t)\,dt+\omega(u_k).$$
Since the topological energy of $\mathbf{u}$ is the sum of the topological energies of $u_1, \ldots, u_m$, we can deduce Equation~\eqref{eq: energy of flowline with cascades} as follows:
\begin{align*}
    E_{top}(\mathbf{u})&=\sum_{k=1}^m E_{top}(u_k)\\
    &=\sum_{k=1}^m\left(\int_{S^1}H_k^+\circ q_k^+(t)\,dt-\int_{S^1} H_k^-\circ q_k^-(t)\,dt+\omega(u_k)\right)\\
    &=\int_{S^1}H_m^+\circ q_m^+(t)\,dt-\int_{S^1}H_1^-\circ q_1^-(t)\,dt+\sum_{k=1}^{m}\omega(u_k).
\end{align*}
as required, where the last equality follows since the sum of the integrals telescopes.

Now, since for every $1\leq k\leq m-1$ the orbits $q_k^+$ and $q_{k+1}^-$ belong to the same critical submanifols, we deduce by Claim \ref{claim: area of level set cylinders vanish} that $\omega(\widehat{q_k^+})=\omega(\widehat{q_{k+1}^-})$, where these cappings are the cappings that we chose previously. Thus, we obtain
\begin{align*}
    \sum_{k=1}^m \omega(u_k) &= \sum_{k=1}^m \omega\left(\widehat{q^-_k}\mathbin{\#} u_k \mathbin{\#} \left(\widehat{q^+_k}\right)^{-1}\right) + \sum_{k=1}^m \left( \omega\left(\widehat{q^+_k}\right)-\omega\left(\widehat{q^-_k}\right)\right)\\
    &= \sum_{k=1}^m \omega\left(\widehat{q^-_k}\mathbin{\#} u_k \mathbin{\#} \left(\widehat{q^+_k}\right)^{-1}\right) + \omega\left(\widehat{q_m^+}\right)-\omega\left(\widehat{q_1^-}\right)\\
        &= \sum_{k=1}^m \omega\left(\widehat{q^-_k}\mathbin{\#} u_k \mathbin{\#} \left(\widehat{q^+_k}\right)^{-1}\right) + \omega\left(\widehat{q_+}\right)-\omega\left(\widehat{q_-}\right),
\end{align*}

where the second equality follows since the second sum telescopes, and the third equality follows the fact that the pairs $q_1^-$, $q^-$, and $q^+_m$, $q^+$, belong to the same critical submanifolds.

Since $(M,\omega)$ is monotone with monotonicity constant $\kappa\in \R$, we deduce that for every continuous sphere $w\fc S^2\to M$ we have $\omega(w)=\kappa c_1(w)$. Therefore
$$\sum_{k=1}^m \omega(u_k) = \sum_{k=1}^m \kappa c_1\left(\widehat{q^-_k}\mathbin{\#} u_k \mathbin{\#} \left(\widehat{q^+_k}\right)^{-1}\right) + \omega\left(\widehat{q_+}\right)-\omega\left(\widehat{q_-}\right).$$

Additionally, since the trivialization $\tau_0$ is obtained from the capping we chose for the $1$-periodic orbits corresponding to the critical points $\Crit(h_-)\cup\Crit(h_+)$, we deduce that for every $1\leq k\leq m$ we have
$$c_1\left(\widehat{q^-_k}\mathbin{\#} u_k \mathbin{\#} \left(\widehat{q^+_k}\right)^{-1}\right)=c_1^{\tau_0}(u_k).$$
This proves that
 $$\sum_{k=1}^m \omega(u_k) = \kappa\sum_{k=1}^m  c_1^{\tau_0}\left( u_k\right) + \omega\left(\widehat{q_+}\right)-\omega\left(\widehat{q_-}\right).$$
Thus we deduce that
$$E_{top}(\mathbf{u})=\int_{S^1}H_m^+\circ q_m^+(t)\,dt-\int_{S^1}H_1^-\circ q_1^-(t)\,dt+\kappa\sum_{k=1}^m  c_1^{\tau_0}\left( u_k\right) + \omega\left(\widehat{q_+}\right)-\omega\left(\widehat{q_-}\right),$$
and in particular
$$E_{top}(\mathbf{u})\leq K_1+K_2+\kappa\sum_{k=1}^m  c_1^{\tau_0}\left( u_k\right).$$
The assumption
$$ m + N = \mu^\tau_{FMB}(q^+; H_+) - \mu^\tau_{FMB}(q^-; H_-) + 2 \sum_{k=1}^m c^\tau_1(u_k) + m . $$
and the fact that
\begin{align*}
  \mu^\tau_{FMB}(q^+; H_+) - \mu^\tau_{FMB}(q^-; H_-) &+ 2 \sum_{k=1}^m c^\tau_1(u_k)\\
  &=\mu^{\tau_0}_{FMB}(q^+; H_+) - \mu^{\tau_0}_{FMB}(q^-; H_-) + 2 \sum_{k=1}^m c^{\tau_0}_1(u_k),  
\end{align*}
as the expression does not depend on the chosen trivialization, imply that
$$\sum_{k=1}^m c^{\tau_0}_1(u_k)=\frac{1}{2}\left(N+\mu^{\tau_0}_{FMB}(q^-; H_-)-\mu^{\tau_0}_{FMB}(q^+; H_+)\right),$$
and hence $\kappa\sum_{k=1}^m  c_1^{\tau_0}\left( u_k\right)\leq K_3$. This proves that
$$E_{top}(\mathbf{u})\leq K_1+K_2+K_3=E,$$
which completes the proof.

\end{proof}

    Next we prove a lower bound for the energy of every non-constant Floer solution.
    \begin{prop}\label{prop:minimal energy quantum} 
   Let $(M, \omega)$ be a closed monotone symplectic manifold. Let $H$ be an autonomous Hamiltonian satisfying the \textbf{MB} condition, and let $h \fc \bS_H \to \R$ be a Morse function. Then there exists $\hbar > 0$ such that for every almost complex structure $J$ on $M$, any critical submanifolds $S^-, S^+ \subset \bS_H$, and every non-constant Floer solution $u \in \cM(S^-, S^+; H, J)$ connecting $S^-$ to $S^+$, the energy of $u$ satisfies
$$ E_{top}(u) \ge \hbar. $$
    \end{prop}
    \begin{proof}
        Denote by $\kappa$ the monotonicity constant of $(M, \omega)$. For every $1$-periodic orbit $x$ of $H$ that is a critical point of $h$, choose a capping $\hat{x}$. Since $\bS_H$ is compact, $\Crit(h)$ is finite. For this reason, the set
$$ A = \left\{ \int_{S^1} H \circ x(t) \, dt - \int_{S^1} H \circ y(t) \, dt + \omega(\hat{x}) - \omega(\hat{y}) : x, y \in \Crit(h) \right\} $$
is a finite subset of $\R$. Denote 
$$ B = (A + \kappa \Z) \cap (0, +\infty) = \{ a + \kappa m : a \in A, \, m \in \Z, \, a + \kappa m > 0 \}. $$

Consider $\hbar = \inf B$. Since $\kappa \Z$ is discrete and $A$ is finite, it follows that $B$ is discrete; in particular, $0$ is not an accumulation point of $B$. This implies that $\hbar > 0$.

Let $J$ be an almost complex structure on $M$, let $S^-, S^+ \subset \bS_H$ be critical submanifolds, and let $u \in \cM(S^-, S^+; H, J)$ be a non-constant Floer solution connecting $S^-$ to $S^+$. Let us prove that $E_{top}(u) \ge \hbar$.

Since $u$ is non-constant in the $\R$ variable, $E_{top}(u) > 0$. Denote by $q^- \in S^-$ and $q^+ \in S^+$ the orbits to which $u$ is asymptotic at $-\infty$ and $+\infty$, respectively, and by $\widehat{q^-}$ and $\widehat{q^+}$ their corresponding cappings. Then:
\begin{align*}
    E_{top}(u) &= \int_{S^1} H \circ q^+ \, dt - \int_{S^1} H \circ q^- \, dt + \omega(u) \\
    &= \int_{S^1} H \circ q^+ \, dt - \int_{S^1} H \circ q^- \, dt + \omega\left(\widehat{q^-} \# u \# (\widehat{q^+})^{-1}\right) - \left(\omega(\widehat{q^-}) - \omega(\widehat{q^+})\right) \\
    &= \left(\int_{S^1} H \circ q^+ \, dt - \int_{S^1} H \circ q^- \, dt\right) + \omega(\widehat{q^+}) - \omega(\widehat{q^-}) + \kappa c_1\left(\widehat{q^-} \# u \# (\widehat{q^+})^{-1}\right).
\end{align*}
In the last equality, we used the fact that $\widehat{q^-} \# u \# (\widehat{q^+})^{-1}$ is a sphere and $\kappa$ is the monotonicity constant of $(M, \omega)$. This equality shows that $E_{top}(u) \in B$, which proves that 
$$ E_{top}(u) \ge \inf B = \hbar, $$
as required.

    \end{proof}

    As a corollary of the previous two propositions, we obtain a bound on the number of cascades:
    
    \begin{coroll}\label{cor: bound on number of cascades}
Let $(M, \omega)$ be a closed monotone symplectic manifold. Let $H_-$ and $H_+$ be Hamiltonians on $(M, \omega)$ satisfying the \textbf{MB} condition, such that each Hamiltonian is constant along each of its critical submanifolds, and let $\cH = (H_s)_{s \in \R}$ be a monotone homotopy of Hamiltonians such that $H_s = H_-$ for every $s \le 0$ and $H_s = H_+$ for $s \ge 1$. Let $h_-$ and $h_+$ be Morse functions on $\bS_{H_-}$ and $\bS_{H_+}$, respectively.

For every integer $N \in \Z$, there exists a constant $\mathbf{m} \ge 0$ such that for any almost complex structures $J_-$ and $J_+$ compatible with $\omega$, possibly time-dependent, any homotopy of compatible almost complex structures $\cJ = (J_s)_{s \in \R}$ satisfying $J_s = J_-$ for $s \le 0$ and $J_s = J_+$ for $s \ge 1$, any Riemannian metrics $g_-$ and $g_+$ on $\bS_{H_-}$ and $\bS_{H_+}$, and any critical points $q_- \in \Crit(h_-)$ and $q_+ \in \Crit(h_+)$, if a continuation flowline $\mathbf{u}$ with $m$ cascades $(u_k)_{k=1}^m$ and a trivialization $\tau$ of $\det TM$ along its asymptotes satisfies
$$ m + N = \mu^\tau_{FMB}(q^+; H_+) - \mu^\tau_{FMB}(q^-; H_-) + 2 \sum_{k=1}^m c^\tau_1(u_k) + m, $$
then the number of cascades $m$ is bounded by $\mathbf{m}$:
$$ m \le \mathbf{m}. $$
    \end{coroll}
    \begin{proof}
      Let $N \in \Z$. Denote by $E$ and $\hbar$ the constants that are guaranteed by Proposition~\ref{prop:universal energy bound} and by Proposition~\ref{prop:minimal energy quantum}, respectively. 

Consider $\mathfrak{m} = 1 + \frac{E}{\hbar}$. Let $J_-$ and $J_+$ be almost complex structures on $M$ compatible with $\omega$, which are possibly time-dependent; let $\cJ = (J_s)_{s \in \R}$ be a homotopy of compatible almost complex structures satisfying $J_s = J_-$ for every $s \le 0$ and $J_s = J_+$ for every $s \ge 1$; let $g_-$ and $g_+$ be Riemannian metrics on $\bS_{H_-}$ and $\bS_{H_+}$, respectively; let $q_- \in \Crit(h_-)$ and $q_+ \in \Crit(h_+)$ be critical points; let $m \in \N$ be an integer; let $\mathbf{u}$ be a continuation flowline with $m$ cascades $(u_k)_{k=1}^m$; and let $\tau$ be a trivialization of $\det TM$ along the asymptotes of the cascades. Assume that
$$ m + N = \mu^\tau_{FMB}(q^+; H_+) - \mu^\tau_{FMB}(q^-; H_-) + 2 \sum_{k=1}^m c^\tau_1(u_k) + m. $$

We know that $E_{top}(\mathbf{u}) = \sum_{k=1}^m E_{top}(u_k)$. The topological energies of $u_1, \ldots, u_m$ are all non-negative, and except for at most one (in the case of a continuation flowline which is not Floer), they are all bounded from below by $\hbar$. Hence:
$$ E \ge E_{top}(\mathbf{u}) = \sum_{k=1}^m E_{top}(u_k) \ge (m-1) \hbar, $$
therefore $m \le \frac{E}{\hbar} + 1 = \mathfrak{m}$, as required.
    \end{proof}

    \subsection{Obstructions on solutions for perturbation data close to the ``special'' one}\label{ss: transversality}
    Let $\mathscr{J}^\infty(M,\omega)$ denote the space of all almost complex structure on $M$, compatible with $\omega$, equipped with the $C^\infty$ topology.

    Let $J_{std}$ denote the standard complex structure on $\CP^n$.

       \begin{thm}\label{tmh: transversality for ind=1} Let $0 < \Delta < 1$ and $\ell \in \Z_{\geq 0}$. Define $H_\ell \fc \CP^n \to \R$ by $H_\ell(z) = h(\Delta, \ell, \mu(z))$ for every $z \in \CP^n$. Denote by $\mathscr{N}_{\ell}\subseteq \mathscr{J}^\infty(\CP^n,\omega_{FS})$ the collection of all compatible almost complex structures $J$ that satisfy the following property:

For every pair of critical points $q_+$ and $q_-$ of the Morse functions associated with the critical submanifolds of $H_\ell$, and for every Floer flowline $u \in \cM(q_-, q_+; H_\ell, J)$ consisting of an $m$-cascade $(u_k)_{k=1}^m$, if 
\[\mu^\tau_{FMB}(q_+;H_\ell) - \mu^\tau_{FMB}(q_-;H_\ell) + 2\sum_{k=1}^m c_1^\tau(u_k) = 1,\]
then $(q_-, q_+)$ cannot be any of the following pairs:
\begin{itemize}
    \item $(\hat{x}^\ell_i, \hat{x}^\ell_j)$ for any $1 \leq i, j \leq \ell$.
    \item $(\check{x}^\ell_i, \check{x}^\ell_j)$ for any $0 \leq i, j \leq \ell + n$.
    \item $(\check{x}^\ell_i, \hat{x}^\ell_j)$ for any $0 \leq i \leq \ell + n$ and $1 \leq j \leq \ell$, except where $i=j$ and $n=1$, in which case $u$ does not intersect $D_\infty$.
    \item $(\hat{x}^\ell_i, \check{x}^\ell_j)$ for any $1 \leq i \leq \ell$ and $0 \leq j \leq \ell + n$, except where $j \in \{i-1, i+n\}$. Moreover, if $j=i-1$, then $u$ does not intersect $D_\infty$.
\end{itemize}

Then $\mathscr{N}_{\ell}$ contains an open neighborhood of the standard almost complex structure $J_{std}$.
    \end{thm}
    \begin{proof}
        Assume the contrary, then every neighborhood of $J_{std}$ contains an almost complex structure for which there exists a counterexample, namely a Floer flowline with cascades invalidating one of the conditions. By picking a countable neighborhood base of $J_{std}$, we thus obtain a sequence of counterexamples, namely, almost complex structures $J_n$, converging to $J_{std}$ in the $C^\infty$ topology, and a sequence of Floer flowlines with cascades, $(u^n)_{n=1}^\infty = \left(\left(u^n_k\right)_{k=1}^{m_n}\right)_{n=1}^\infty$, with $\left(\left(u^n_k\right)_{k=1}^{m_n}\right) \in \cM(q_n^-,q_n^+;H,J_n)$, with $\left(q_n^-,q_n^+\right)$ violating the list of forbidden values given in the theorem's formulation, and such that 
        \[\mu^\tau_{FMB}(q_n^+;H_\ell) - \mu^\tau_{FMB}(q_n^-;H_\ell) + 2\sum_{k=1}^{m_n} c_1^\tau(u^n_k) = 1.\]
        Since the list of values possible for $\left(q_n^-,q_n^+\right)$ is finite, there exists a constant subsequence of $\left(q_n^-,q_n^+\right)$, which, by abuse of notation, we still index by $n$, whose value we denote by $\left(q_n^-,q_n^+\right) = \left(q^-,q^+\right)$ .
        Next, by Corollary \ref{cor: bound on number of cascades}  the sequence $\left(m_n\right)_{n=1}^\infty$ is bounded, hence by passing to a subsequence (still indexed by $n$), we may assume $m_n = m$ is independent of $n$.

        By Proposition \ref{prop:universal energy bound} there exists $E$ such that $E_{top}\left(u^n\right) \le E$. By Theorem \ref{thm:subsequence converging to a broken bubbled flowline with cascades} there exists a subsequence, which, by abuse of notation we still denote by $\left(u^n\right)_{n=1}^\infty$, which Floer-Gromov converges to a broken bubbled Floer flowline with cascades $\left(\left(\mathbf{w}_k\right)_{k=1}^b, (j)\right)$, where each $\mathbf{w}_k$ has $m_k$ cascades. Let us enumerate the Floer solutions appearing as the cascades in all the pieces of  $\left(\mathbf{w}_k\right)_{k=1}^b$, in order of appearence, as $w_\lambda$ for $\lambda \in \left[1,a\right]\cap\N$, where $a=\sum_{k=1}^b m_k$, and index the collection of all the bubbles in all the trees, by $\left\lbrace B_\alpha \right\rbrace_{\alpha \in I}$ for some arbitrary indexing set $I$.
        Then, by applications of Remark \ref{rem: FloerGromov convergence properties} we obtain \[\mu^\tau_{FMB}(q^+;H_\ell) - \mu^\tau_{FMB}(q^-;H_\ell) + 2\sum_{k=1}^a c_1^\tau(w_k) + 2\sum_{\alpha\in A} c_1(B_\alpha) = 1. \] 
        Moreover, in the case that for all $n$,  $(u^n)_{n=1}^\infty = \left(\left(u^n_k\right)_{k=1}^{m_n}\right)_{n=1}^\infty$, are connecting $\hat{x}^\ell_i$ to $\check{x}^\ell_{i-1}$, and are intersecting $D_\infty$, it follows that also their limit $\left(\left(\mathbf{w}_k\right)_{k=1}^b, (j)\right)$ intersects $D_\infty$, possibly in an  asymptotic internal intersection. And, in the case that $n=1$ and for all $n$,  $(u^n)_{n=1}^\infty = \left(\left(u^n_k\right)_{k=1}^{m_n}\right)_{n=1}^\infty$, and are connecting $\check{x}^\ell_i$ to $\hat{x}^\ell_{i}$, are intersecting $D_\infty\cup D_1$, it follows that also their limit $\left(\left(\mathbf{w}_k\right)_{k=1}^b, (j)\right)$ intersects $D_\infty\cup D_1$, possibly in an  asymptotic internal intersection. Thus, $\left(\left(\mathbf{w}_k\right)_{k=1}^b, (j)\right)$ violates Theorem \ref{thm: diff_obst}, a contradiction.
        \end{proof}

        Let $0 < \Delta'\leq\Delta < 1$ and $\ell,\ell' \in \Z_{\geq 0}$ with $\ell\leq \ell'$. Define $H_\ell,H'_{\ell'} \fc \CP^n \to \R$ by $H_\ell(z) = h(\Delta, \ell, \mu(z))$ and $H'_{\ell'}(z) = h(\Delta', \ell', \mu(z))$ for every $z \in \CP^n$. Let $0 < \Delta'\leq\Delta < 1$ and $\ell,\ell' \in \Z_{\geq 0}$ with $\ell\leq \ell'$. Define $H_\ell,H'_{\ell'} \fc \CP^n \to \R$ by $H_\ell(z) = h(\Delta, \ell, \mu(z))$ and $H'_{\ell'}(z) = h(\Delta', \ell', \mu(z))$ for every $z \in \CP^n$. 
        Let us define the following space of monotone homotopies of Hamiltonians:
        \[\mathscr{H}\left(H_\ell,H'_{\ell'}\right) = \left\lbrace \mathcal{H}=\left(H_s\right)_{s\in\R} \,\mid\, \partial_s H_s\geq0,\; H_s = H_\ell \text{ for } s\le 0, \text{ and } H_s = H'_{\ell'} \text{ for } s\ge 1 \right\rbrace,\] equipped with the $C^\infty$ topology.
        Let $\mathcal{H}_0 \in \mathscr{H}\left(H_\ell,H'_{\ell'}\right)$ denote a homotopy as in Theorem \ref{thm: continuation_obst}.

        Moreover, let us denote by $\mathscr{J}_h^\infty(M,\omega)$ ($h$ stands for ``homotopy'') the space of all homotopies $\mathcal{J}=\left(J_s\right)_{s\in\R}$ of almost complex structures such that:
        \begin{enumerate}
            \item For all $s\in\R$, $J_s$ is an $\omega$ compatible almost complex structure on $M$.
            \item There exists an almost complex structure $J_-$ such that $J_s = J_-$ for $s\le 0$.   
            \item There exists an almost complex structure $J_+$ such that $J_s = J_+$ for $s\ge 1$.
        \end{enumerate}
        We equip the space $\mathscr{J}_h^\infty(M,\omega)$ with the $C^\infty$ topology, and we denote by $\mathcal{J}_{std}$ the, constant in $s$, homotopy $\left(J_{std}\right)_{s\in R}$.
    
  \begin{thm}\label{tmh: transversality for ind=0}
        There exists a $C^\infty$-neighborhood $\mathscr{N}_{\mathcal{H}_0}$ of $\mathcal{H}_0$, and a $C^\infty$-neighborhood $\mathscr{N}_{\mathcal{J}_{std}}$ of $\mathcal{J}_{std}$ such that for all $\mathcal{H} \in \mathscr{N}_{\mathcal{H}_0}$ and $\mathcal{J}\in \mathscr{N}_{\mathcal{J}_{std}}$ if $u \in \cM^{cont}\left(q^-,q^+;\mathcal{H},\mathcal{J}\right)$ is a continuation flowline with $m$-cascades $\left(u_k\right)_{k=1}^m$, such that 
        \[\mu^\tau_{FMB}(q_+;H'_{\ell'}) - \mu^\tau_{FMB}(q_-;H_\ell) + 2\sum_{k=1}^m c_1^\tau(u_k) = 0,\]
        then, as in Theorem \ref{thm: continuation_obst}, $(q^-,q^+)$ cannot be any of the following:
        \begin{itemize}
            \item $(\check{x}^\ell_i, \hat{y}^{\ell'}_j)$ for every $0 \leq i \leq \ell + n$ and $1 \leq j \leq \ell'$.
            \item $(\hat{x}^\ell_i, \check{y}^{\ell'}_j)$ for every $1 \leq i \leq \ell$ and $0 \leq j \leq \ell' + n$.
            \item $(\hat{x}^\ell_i, \hat{y}^{\ell'}_j)$ for every $1 \leq i \leq \ell$ and $1\leq j\leq \ell'$, except where $i=j$.
            \item $(\check{x}^\ell_i, \check{y}^{\ell'}_j)$  for every $0 \leq i \leq \ell+n$ and $0\leq j\leq \ell'+n$, except where $i=j$.
        \end{itemize}
    \end{thm}
    \begin{proof}
        Assume the contrary, then every pair of neighborhoods of $\mathcal{H}_0$ and $\mathcal{J}_{std}$ contain a homotopy of Hamiltonians, and a homotopy of almost complex structures for which there exists a counterexample to the theorem's statement, namely a continuation flowline with cascades invalidating one of the conditions. By picking countable neighborhood base of $\mathcal{H}_0$ and of $\mathcal{J}_{std}$, we thus obtain a sequence of counterexamples, namely, a sequence of homotopies of Hamiltonians $\mathcal{H}_n$ and a sequence of homotopies of almost complex structures $\mathcal{J}_n$, converging to $\mathcal{H}_0$ and $\mathcal{J}_{std}$ in the $C^\infty$ topology, respectively, and a sequence of continuation flowlines with cascades, $(u^n)_{n=1}^\infty = \left(\left(u^n_k\right)_{k=1}^{m_n},(j_n)\right)_{n=1}^\infty$, with $\left(\left(u^n_k\right)_{k=1}^{m_n},(j_n)\right) \in\cM^{cont}\left(q^-,q^+;\mathcal{H}_n,\mathcal{J}_n\right)$, with $\left(q_n^-,q_n^+\right)$ violating the list of forbidden values given in the theorem's formulation, and such that 
        \[\mu^\tau_{FMB}(q_n^+;H'_{\ell'}) - \mu^\tau_{FMB}(q_n^-;H_\ell) + 2\sum_{k=1}^{m_n} c_1^\tau(u^n_k) = 0.\]
        Since the list of possible values for $\left(q_n^-,q_n^+\right)$ is finite, there exists a constant subsequence of $\left(q_n^-,q_n^+\right)$, which, by abuse of notation, we still index by $n$, whose value we denote by $\left(q_n^-,q_n^+\right) = \left(q^-,q^+\right)$ .
        Next, by Corollary \ref{cor: bound on number of cascades}  the sequence $\left(m_n\right)_{n=1}^\infty$ is bounded, hence by passing to a subsequence (still indexed by $n$), we may assume $m_n = m$ is independent of $n$. Moreover the set of possible values for $j_n$ is now finite, and by passing to a subsequence we may assume $j_n=j$ is a constant subsequence.

        By Proposition \ref{prop:universal energy bound} there exists $E$ such that $E_{top}\left(u^n\right) \le E$. By Theorem \ref{thm:subsequence converging to a broken bubbled flowline with cascades} there exists a subsequence, which, by abuse of notation we still denote by $\left(u^n\right)_{n=1}^\infty$, which Floer-Gromov converges to a broken bubbled continuation solution $\left(\left(\mathbf{w}_k\right)_{k=1}^b,(j')\right)$, for the data $\left(\mathcal{H}_0,\mathcal{J}_{std}\right)$, where each $\mathbf{w}_k$ has $m_k$ cascades. Let us enumerate the Floer and continuation solutions appearing as the cascades in all the pieces of  $\left(\mathbf{w}_k\right)_{n=1}^\infty$, in order of appearance, as $w_\lambda$ for $\lambda \in \left[1,a\right]\cap\N$, where $a=\sum_{k=1}^b m_k$, and index the collection of all the bubbles in all the trees, by $\left\lbrace B_\alpha \right\rbrace_{\alpha \in I}$ for some arbitrary indexing set $I$.
        Then, by applications of Remark \ref{rem: FloerGromov convergence properties} we obtain \[\mu^\tau_{FMB}(q^+;H'_{\ell'}) - \mu^\tau_{FMB}(q^-;H_\ell) + 2\sum_{k=1}^a c_1^\tau(w_k) + 2\sum_{\alpha\in A} c_1(B_\alpha) = 0. \] 
        Thus, $\left(\left(\mathbf{w}_k\right)_{k=1}^b, (j')\right)$ violates Theorem \ref{thm: continuation_obst}, a contradiction.
        \end{proof}

    Let $\ell,\ell'\in \Z_{\geq0}$ and assume that $\ell\leq \ell'$. Let $\chi\fc \R\to \R$ be a smooth non-increasing function, satisfying, that there exists an $0<\ve<1 $ such that $\chi(s)=\Delta$ for every $s \in [0,\ve)$ and $\chi(s)=\Delta'$ for every $s\in(1-\ve,1]$. Similarly, let $\lambda\fc \R\to \R$ be a smooth non-decreasing function, satisfying $\lambda(s)=\ell$ for every $s \in [0,\ve)$ and $\lambda(s)=\ell'$ for every $s\in(1-\ve,1]$. 

    Now, define a $2$-cube of Floer data by $\mathcal{H}_0\fc [0,1]\times[0,1]\times \CP^n\to \R$ by $\mathcal{H}_0(s_0,s_1,z)=h(\chi(s_1),\lambda(s_0),\mu(z))$ for every $(s_0,s_1,z)\in \R\times \CP^n$. 
    We denote by $\mathcal{J}_{std}$ constant $2$-cube of almost complex stuctures, whose value is $J_{std}$. 
    Note the the $s_0$ edges of this cube, when taken as $1$-cubes, correspond to the continuation maps from $h(\Delta,\ell,\mu(\cdot))$ to $h(\Delta,\ell',\mu(\cdot))$ and from $h(\Delta',\ell,\mu(\cdot))$ to $h(\Delta',\ell',\mu(\cdot))$, and the $s_1$ edges, as $1$-cubes correspond to the continuation maps from $h(\Delta,\ell,\mu(\cdot))$ to $h(\Delta',\ell,\mu(\cdot))$ and from $h(\Delta,\ell',\mu(\cdot))$ to $h(\Delta',\ell',\mu(\cdot))$:
    
    \begin{figure}[htbp]
    \centering
    \includegraphics[width=0.33\linewidth]{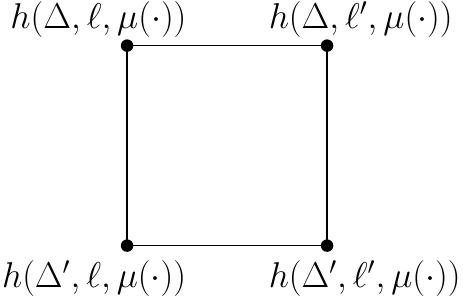}
    \caption{The $2$-cube and the Hamiltonians at its corners.}
    \label{fig:figh-2-cube}
    \end{figure}
    
    Let us define the following space of $2$-cubes of Hamiltonians, fix $\epsilon<\frac{1}{2}$:
    \[\mathscr{H}\left(\ell,\Delta,\ell',\Delta'\right) = \left\lbrace \mathcal{H}\text{ is a monotone $2$-cube} \,\mid\, \mathcal{H} = \mathcal{H}_0 \text{ in } \epsilon\text{-balls around the vertices}\right\rbrace\] equipped with the $C^\infty$ topology.

    Moreover, let us denote by $\mathscr{J}_2^\infty(M,\omega)$ (the $2$ denoting $2$-cubes), the space of all $2$-cubes of compatible almost complex structure, such that they are constant in $\epsilon$-balls around the vertices.

    \begin{thm} 
        There exists a $C^\infty$-neighborhood $\mathscr{N}_{\mathcal{H}_0}$ of $\mathcal{H}_0$ in $\mathscr{H}\left(\ell,\Delta,\ell',\Delta'\right)$, and a $C^\infty$-neighborhood $\mathscr{N}_{\mathcal{J}_{std}}$ of $\mathcal{J}_{std}$ in $\mathscr{J}_2^\infty(M,\omega)$ such that for all $\mathcal{H} \in \mathscr{N}_{\mathcal{H}_0}$ and $\mathcal{J}\in \mathscr{N}_{\mathcal{J}_{std}}$ if $(u,\gamma)$ is a pair of:
        \begin{enumerate}
            \item $\gamma$ a negative gradient trajectory on the $2$-cube $[0,1]\times[0,1]$, for the standard Morse function from Formula from Section~\ref{sss:Ham cubes FMB}, 
            \item and, $u$, a continuation flowline with $m$-cascades, $\left(u_k\right)_{k=1}^m$, connecting $q^-$ to $q^+$, with respect to the Floer data $\mathcal{H}\circ\gamma$ and $\mathcal{J}\circ\gamma$,
        \end{enumerate}   such that 
        \[\mu^\tau_{FMB}(q^+;h(\Delta',\ell',\mu(\cdot))) - \mu^\tau_{FMB}(q^-;h(\Delta,\ell,\mu(\cdot))) + 2\sum_{k=1}^m c_1^\tau(u_k) = -1,\]
        then, as in Theorem \ref{thm: htpy_obst}, $(q^-,q^+)$ cannot be any of the following:
        \begin{itemize}
            \item $(\check{x}^\ell_i, \check{y}^{\ell'}_j)$ for every $0 \leq i \leq \ell + n$ and $0 \leq j \leq \ell'+n$.
            \item $(\hat{x}^\ell_i, \hat{y}^{\ell'}_j)$ for every $1 \leq i \leq \ell$ and $1 \leq j \leq \ell' $.
            \item $(\check{x}^\ell_i, \hat{y}^{\ell'}_j)$ for every $0 \leq i \leq \ell+n$ and $1\leq j\leq \ell'$.
            \item $(\hat{x}^\ell_i, \check{y}^{\ell'}_j)$ for every
            $1\leq i\leq\ell$ and $0\leq j\leq\ell'+n$,
            except where $j=i+n-1$.
        \end{itemize}
    \end{thm}
     \begin{proof}
        Assume the contrary, then every pair of neighborhoods of $\mathcal{H}_0$ and $\mathcal{J}_{std}$ in $\mathscr{H}\left(\ell,\Delta,\ell',\Delta'\right)$ and $\mathscr{J}_2^\infty(M,\omega)$, that contains a $2$-cube of Hamiltonians, and a $2$-cube of almost complex structures for which there exists a counterexample to the theorem's statement, namely, a pair $(\gamma,u)$ with $u$ a continuation flowline  with cascades for the Floer data obtained by tha cube and $\gamma$, invalidating one of the conditions. By picking countable neighborhood base of $\mathcal{H}_0$ and of $\mathcal{J}_{std}$, we thus obtain a sequence of counterexamples, namely, a sequence of cubes of Hamiltonians $\mathcal{H}_k$ and a sequence of cubes of almost complex structures $\mathcal{J}_k$, converging to $\mathcal{H}_0$ and $\mathcal{J}_{std}$ in the $C^\infty$ topology, respectively, and a sequence of pairs $(\gamma_k,u_k)$, with $\gamma_k$ a negative gradient flowline on $[0,1]\times[0,1]$, and $u_k$ a continuation flowline with cascades, for the Floer data $\mathcal{H}_k\circ \gamma$, and $\mathcal{J}_k\circ \gamma$.

        Similarly to the proof of Theorems \ref{tmh: transversality for ind=1}  and \ref{tmh: transversality for ind=0}, we pass to a subsequence converging to a broken bubbled solution. Note that here $\gamma_k$ may converge to a broken negative gradient trajectory as well. Thus, we obtain either $(\gamma_\infty,u_\infty)$ satisfying condition (i) of Theorem \ref{thm: htpy_obst}, or $(\gamma^1_\infty,u_\infty)$ and $(\gamma^2_\infty,v_\infty)$ satisfying conditions (ii) and (iii) of Theorem \ref{thm: htpy_obst}, but violating the list of impossible pairs of asymptotes $(q^-,q^+)$, hence a contradiction.
    \end{proof}

\subsection{Tower}\label{ss:tower}

 Let $\Delta\in (0,1)$ and consider the acceleration data $(H_\ell)_{\ell\geq0}$ as we defined in Section~\ref{ss: acc. data}. For every compatible almost complex structure $J$ on $\CP^n$ we denote by $D_\ell(J)$ the collection of pairs $(q^-,q^+)$ of critical points of the associated Morse function on the critical submanifold with respect to $H_\ell$, for which the moduli space of Floer flowlines of index $1$ connecting $q^-$ to $q^+$ are empty, also, we denote by $D^{bbl}_\ell(J)$ the collection of pairs $(q^-,q^+)$ of critical points of the associated Morse function on the critical submanifold with respect to $H_\ell$, for which the moduli space of broken bubbled Floer flowlines of index $1$ connecting $q^-$ to $q^+$ are empty. Theorem~\ref{tmh: transversality for ind=1} says that for every $\ell\geq 0$ there is an open neighborhood $\mathscr{N}_\ell$ of $J_{std}$ such that for every $J\in \mathscr{N}_\ell$ we have $D^{bbl}_\ell(J_{std})\subset D_\ell(J)$. 

 Similarly, for every pair of homotopies of Floer data, namely, $\cH$ a homotopy of Hamiltonians from $H_\ell$ to $H_{\ell+1}$, and $\cJ$ a homotopy of compatible almost complex structures, denote by $C_\ell(\cH,\cJ)$ the collection of pairs of critical points of the associated Morse functions on the critical submanifold with respect to $H_\ell$ and $H_{\ell+1}$, for which the moduli space of continuation flowlines of index $0$ are empty, also, we denote denote by $C^{bbl}_\ell(\cH,\cJ)$ the collection of pairs of critical points of the associated Morse functions on the critical submanifold with respect to $H_\ell$ and $H_{\ell+1}$, for which the moduli space of broken bubbled continuation flowlines of index $0$ are empty. Theorem~\ref{tmh: transversality for ind=0} says that for every $\ell\geq 0$ there is a $C^\infty$-neighborhood $\mathscr{N}_{\mathcal{H}_\ell}$ of $\mathcal{H}_\ell$, with $\mathcal{H}_\ell$ a homotopy as in Theorem \ref{thm: continuation_obst}, and a $C^\infty$-neighborhood $\mathscr{N}_{\cH_\ell,\mathcal{J}_{std}}$ of $\mathcal{J}_{std}$ such that for every $\cH\in\mathscr{N}_{\mathcal{H}_\ell}$ and $\cJ\in \mathscr{N}_{\mathcal{J}_{std}}$ we have  $C^{bbl}_\ell(\cH,\cJ)\subset C_\ell(\cH,\cJ)$.

Our construction of a direct system of Floer complexes, for the sake of computing the relative symplectic cohomology of the ball, relies on constructing acceleration data which is both regular and satisfies the vanishing properties of matrix elements described in Theorems \ref{thm: diff_obst} and \ref{thm: htpy_obst}. In Proposition \ref{prop: finite tower} we construct it in stages, as a stabilizing limit of finite sequences, where in each step we extend the sequence of length $m$ by $1$, while possibly changing only the $m$'th step. hence there exists a stable limit.

Moreover, we make sure that all almost complex structures we construct remain equal to $J_{std}$ on a shell $A(r_2,r_1):=\overline{B_{r_2}-B_{r_1}}$ for some fixed $1>r_2>r_1>\Delta$. This is necessary to be able to apply the no escape lemma (Lemma \ref{lemma: no escape}).

Regularity for continuation maps is achieved by a generic perturbation of the homotopy, whereas regularity for the Floer differential is achieved by perturbing $J_{std}$ in the complement of $A(r_2,r_1)$. Since all the Floer complex generators lie in the complement of $A(r_2,r_1)$, all Floer solutions intersect the support of the perturbation in an open set, hence regularity is guaranteed. See \cite[Chapter 2]{Fauck_2016_thesis} for transversality results in Morse-Bott with cascades, in the Rabinowitz Floer homology, \cite[Proposition 3.5]{BO_Morse_Bott_2009}, for transversality in Floer homology with cascades for autonomous Hamiltonians, and \cite{FHS_1995_Transversality} for transvesality in the nondegenerate case. many of the ideas carry to the Morse-Bott case.

 \begin{prop}\label{prop: finite tower}
 
     Let $m\in \N$. Let $(\cH'_\ell,\cJ'_\ell)_{\ell=0}^m$ be a sequence of pairs of regular homotopies, satisfing the following:
     \begin{itemize}
         \item $\cH'_\ell\in\mathscr{N}_{\mathcal{H}_\ell}$, for every $0\leq \ell\leq m$.
         \item  $\cJ_\ell\in \mathscr{N}_{\cH_\ell,\mathcal{J}_{std}}$, for every $0\leq \ell\leq m$.
         \item There exists a sequence of compatible almost complex structures $(J_\ell)_{\ell=0}^{m+1}$ such that 
         \begin{itemize}
         \item $J_\ell\in \cN_\ell$ and the pair $(H_\ell,J_\ell)$ is regular, for every $0\leq \ell\leq m+1$.
         \item
             On the shell  $A(r_2,r_1)$, $J_\ell = J_{std}$.
             \item  $\cJ'_\ell$ connects $J_\ell$ to $J_{\ell+1}$, for every $0\leq \ell\leq m$.
         \end{itemize}
     \end{itemize}
     Then there exists a sequence $(\cH''_\ell,\cJ''_\ell)_{\ell=0}^{m+1}$ of pairs of regular homotopies, satisfying the following:
     \begin{itemize}
         \item For every $0\leq \ell\leq m-1$ we have $(\cH''_\ell,\cJ''_\ell)=(\cH'_\ell,\cJ'_\ell)$.
        \item $\cH''_\ell\in\mathscr{N}_{\mathcal{H}_\ell}$, for every $0\leq \ell\leq m+1$.
         \item  $\cJ''_\ell\in \mathscr{N}_{\cH_\ell,\mathcal{J}_{std}}$, for every $0\leq \ell\leq m+1$.
         \item There exists a sequence of compatible almost complex structures $(J'_\ell)_{\ell=0}^{m+2}$ such that 
         \begin{itemize}
         \item $J'_\ell\in \cN_\ell$ and the pair $(H_\ell,J'_\ell)$ is regular, for every $0\leq \ell\leq m+2$.
         \item
             On the shell  $A(r_2,r_1)$, $J_\ell = J_{std}$.
             \item  $\cJ''_\ell$ connects $J'_\ell$ to $J'_{\ell+1}$, for every $0\leq \ell\leq m+1$.
         \end{itemize}
        
     \end{itemize}
 \end{prop}
 \begin{proof}
     First, for every $0\leq \ell\leq m-1$, set $(\cH''_\ell,\cJ''_\ell)=(\cH'_\ell,\cJ'_\ell)$. Choose $\cH''_{m+1}\in \mathscr{N}_{\mathcal{H}_{m+1}}$ and $\cJ''_{m+1}\in \mathscr{N}_{\mathcal{H}_{m+1},\cJ_{std}}$ such that $(\cH''_{m+1},J''_{m+1})$ is regular. Also, if we denote by $J'_{m+1}$ and $J'_{m+2}$ the negative and positive ends of $\cJ''_{m+1}$, we may assume that $J'_{m+1}\in \mathscr{N}_{m+1}$, $J'_{m+2}\in \cN_{m+2}$, that on
             the shell  $A(r_2,r_1)$, $J'_{m+1} = J'_{m+2} = J_{std}$, and that the pairs $(H_{m+1},J'_{m+1})$, $(H_{m+2},J'_{m+2})$ are regular. Finally, we can choose $\cH''_m\in \mathscr{N}_{\cH_m}$ and $\cJ''_m\in \mathscr{N}_{\cH_m,\cJ_{std}}$ such that the pair $(\cH''_m,\cJ''_m)$ is regular and that $J'_m$ and $J'_{m+1}$ are the negative and positive ends of $\cJ''_{m}$. This completes the proof.

    \begin{figure}[H]
    \centering
    \includegraphics[width=0.7\linewidth]{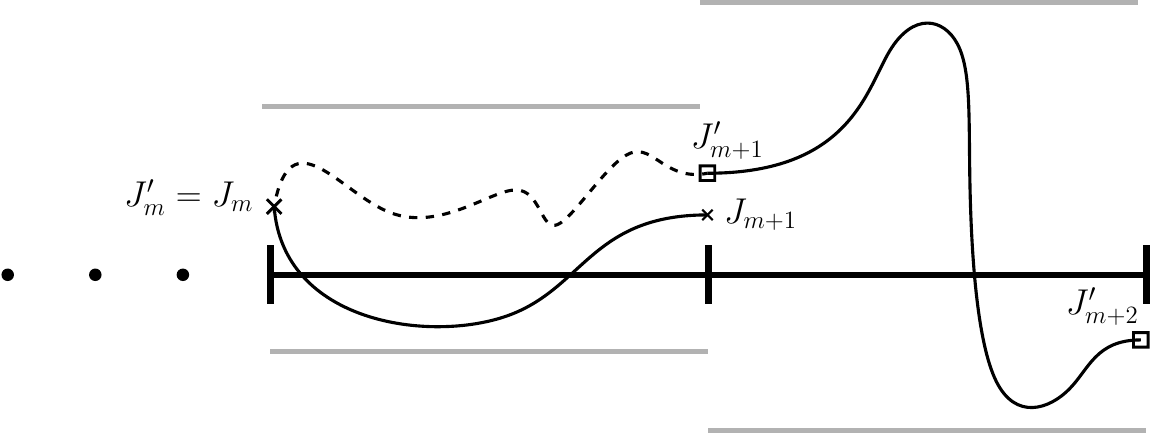}
    \caption{The recursion step is demonstrated. $J'_{m+1}$ is depicted to be chosen in $\mathcal{N}_{m+1}$ and the dashed line represents the homotopy constructed.}
    \label{fig:TowerRecursionStep}
    \end{figure}
 \end{proof}

\section{Proof of Proposition~\ref{prop: asymp_intersection}}\label{app: proof of Seidel's lemma}

The following proof of Proposition~\ref{prop: asymp_intersection} uses a classical result from Fourier theory; we will formulate and prove it for completeness:
\begin{lemma}\label{lemma: series of F coeff}
    Let $f\fc S^1\to \R$ be a differentiable function and assume that $f'$ is square-integrable. Then the series $\sum_{k\in \Z}|\hat{f}(k)|$ converges.
\end{lemma}
\begin{proof}
    Since $f'$ is square-integrable, by the Bessel inequality, the series $\sum_{k\in \Z}|\hat{f'}(k)|^2$ converges. Additionally, for every $k\in \Z$ we have $\hat{f'}(k)=2\pi i k \hat{f}(k)$. Thus, for every $k\in \Z\setminus\{0\}$ we have
    $$|\hat{f}(k)|=\frac{1}{2\pi}\cdot\frac{1}{|k|}\cdot|2\pi i k\hat{f}(k)|=\frac{1}{2\pi}\cdot\frac{1}{|k|}\cdot|\hat{f'}(k)|\leq\frac{1}{2}\left(\left(\frac{1}{2\pi}\cdot\frac{1}{|k|}\right)^2+|\hat{f'}(k)|^2\right).$$
    This shows that
    \begin{align*}
        \sum_{k\in \Z}|\hat{f}(k)|&=|\hat{f}(0)|+\sum_{k\in \Z\setminus\{0\}} |\hat{f}(k)|\\
        &\leq |\hat{f}(0)|+\sum_{k\in \Z\setminus\{0\}} \frac{1}{2}\left(\frac{1}{4\pi^2}\cdot\frac{1}{k^2}+|\hat{f'}(k)|^2\right)\\
        &= |\hat{f}(0)|+\frac{1}{4\pi^2}\cdot\frac{\pi^2}{6}+\frac{1}{2}\sum_{k\in \Z\setminus\{0\}} |\hat{f'}(k)|^2\\
         &\leq |\hat{f}(0)|+\frac{1}{24}+\frac{1}{2}\cdot\sum_{k\in \Z} |\hat{f'}(k)|^2,
    \end{align*}
    Thus we deduce that the series $\sum_{k\in \Z}|\hat{f}(k)|$ converges.
\end{proof}

\begin{proof}[Proof of Proposition~\ref{prop: asymp_intersection}]
    Write $u$ in coordinates $u(s,t)=(v(s,t),b(s,t))$ where $v\fc S\to W$ and $b\fc S\to B$ are smooth functions. Since $u$ is a Floer solution, it satisfies the Floer Equation~\eqref{eq: Floer eq}, which is given here as
    $$\frac{\partial}{\partial s}\begin{pmatrix}
        w\\ b
    \end{pmatrix}+\begin{pmatrix}
        J'&0\\0&i
    \end{pmatrix}\left(\frac{\partial}{\partial t}\begin{pmatrix}
        w \\ b
    \end{pmatrix}-\begin{pmatrix}
        V(w) \\ 2\pi \alpha i b
    \end{pmatrix}\right)=\begin{pmatrix}
        0 \\ 0
    \end{pmatrix}$$

   Thus, we obtain
    \begin{equation}\label{eq: PDE for asmptotic intersection}
        \frac{\partial b}{\partial s}+i\frac{\partial b}{\partial t}+2\pi\alpha b=0.
    \end{equation}
    Let us express $b$ using its Fourier series with respect to the $t$ variable; for every $k\in \Z$ there exists a smooth function $a_k\fc\R\to \C$ such that
    $$b(s,t)=\sum_{k\in \Z}a_k(s)e^{2\pi i kt},$$
    for every $(s,t)\in S$.
    Thus, for every $k\in \Z$, Equation~\eqref{eq: PDE for asmptotic intersection} implies the ordinary differential equation:
    $$a_k'-2\pi (k-\alpha) a_k=0.$$
    Solving this equation yields for every $k\in\Z$ there exists $c_k\in \C$ such that
    $$b(s,t)=\sum_{k\in \Z} c_k e^{2\pi(k-\alpha)s} e^{2\pi i k t},$$
    for every $(s,t)\in S$.

    Here we have two cases:
    \begin{enumerate}
        \item Assume $S=[R,+\infty)\times S^1$ for some $R\in \R$. Since $u$ converges to $p\times\{0\}$, we deduce that $\lim_{s\to+\infty} b(s,t)=0$. This implies that for every $k\in \Z$, if $k-\alpha\geq0$ then $c_k=0$. Given $k\in \Z$, since $\alpha\notin \Z$ we get that $k-\alpha<0$ if and only if $k\leq \lfloor\alpha\rfloor$. 

        Since $b$ is smooth, Lemma~\ref{lemma: series of F coeff} implies that for every $s\in[R,+\infty)$ the series $$\sum_{k\in \Z}|c_k|\cdot e^{2\pi(k-\alpha)s}$$
        converges.
        
        Let $k_0$ be the maximal $k\in \Z$ such that $c_k\neq0$. Denote $B=\sum_{k=-\infty}^{k_0-1} |c_k| e^{2\pi(k-\alpha)R} $. Then for every $(s,t)\in S$ we have
        \begin{align*}
            \left|b(s,t)-c_{k_0} e^{2\pi(k_0-\alpha)s} e^{2\pi i k_0 t}\right|&=\left|\sum_{k=-\infty}^{k_0-1} c_k e^{2\pi(k-\alpha)s} e^{2\pi i k t}\right|\\
            &\leq \sum_{k=-\infty}^{k_0-1} |c_k| e^{2\pi(k-\alpha)s} \\
            &= \sum_{k=-\infty}^{k_0-1} |c_k| e^{2\pi(k-\alpha)(s-R)}e^{2\pi(k-\alpha)R} \\
            &\leq \sum_{k=-\infty}^{k_0-1} |c_k| e^{2\pi(k_0-1-\alpha)(s-R)}e^{2\pi(k-\alpha)R} \\
            &=Be^{2\pi(k_0-1-\alpha)(s-R)}.
        \end{align*}
    
    Dividing this estimate by the absolute value of the leading term,
we obtain
\[
\sup_{t\in S^1}
\frac{\left|b(s,t)-c_{k_0}e^{2\pi(k_0-\alpha)s}
e^{2\pi i k_0t}\right|}
{|c_{k_0}|e^{2\pi(k_0-\alpha)s}}
\leq
\frac{B}{|c_{k_0}|e^{2\pi(k_0-\alpha)R}}
e^{-2\pi(s-R)}.
\]
The right-hand side tends to zero as $s\to+\infty$. Thus, for
$s$ sufficiently large, the straight-line homotopy between
$b(s,\cdot)$ and
$c_{k_0}e^{2\pi(k_0-\alpha)s}e^{2\pi i k_0(\cdot)}$
does not pass through zero. Consequently, for $s$ large enough 
\[
\wind(b(s,\cdot))
=\wind\left(
  \left(c_{k_0}e^{2\pi(k_0-\alpha)s}e^{2\pi i k_0t}\right)_{t\in[0,1]}
\right)
=k_0\leq\lfloor\alpha\rfloor.
\]
    
    Since the constant capping at $p$ has negative orientation, we deduce that for $s$ large enough we have
    $$u\cdot (W\times\{0\})=-\wind(b(s,\cdot))\geq-\lfloor\alpha\rfloor.$$

  \item Assume $S=(-\infty,R]\times S^1$ for some $R\in\R$. Since $u$ converges to $p\times\{0\}$, we deduce that $\lim_{s\to-\infty} b(s,t)=0$. This implies that for every $k\in \Z$, if $k-\alpha\leq0$ then $c_k=0$. Given $k\in \Z$, since $\alpha\notin \Z$ we get that $k-\alpha>0$ if and only if $k\geq 1+\lfloor\alpha\rfloor$. 

        Since $b$ is smooth, Lemma~\ref{lemma: series of F coeff} implies that for every $s\in(-\infty,R]$ the series $$\sum_{k\in \Z}|c_k|\cdot e^{2\pi(k-\alpha)s}$$
        converges.
        
        Let $k_0$ be the minimal $k\in \Z$ such that $c_k\neq0$. Denote $B=\sum_{k=k_0+1}^{\infty} |c_k| e^{2\pi(k-\alpha)R} $. Then for every $(s,t)\in S$ we have
        \begin{align*}
            \left|b(s,t)-c_{k_0} e^{2\pi(k_0-\alpha)s} e^{2\pi i k_0 t}\right|&=\left|\sum_{k=k_0+1}^{\infty} c_k e^{2\pi(k-\alpha)s} e^{2\pi i k t}\right|\\
            &\leq \sum_{k=k_0+1}^{\infty} |c_k| e^{2\pi(k-\alpha)s} \\
            &= \sum_{k=k_0+1}^{\infty} |c_k| e^{2\pi(k-\alpha)(s-R)}e^{2\pi(k-\alpha)R} \\
            &\leq \sum_{k=k_0+1}^{\infty} |c_k| e^{2\pi(k_0+1-\alpha)(s-R)}e^{2\pi(k-\alpha)R} \\
            &=Be^{2\pi(k_0+1-\alpha)(s-R)}.
        \end{align*}
    Dividing this estimate by the absolute value of the leading term,
    we obtain
    \[
    \sup_{t\in S^1}
    \frac{\left|b(s,t)-c_{k_0}e^{2\pi(k_0-\alpha)s}
    e^{2\pi i k_0t}\right|}
    {|c_{k_0}|e^{2\pi(k_0-\alpha)s}}
    \leq
    \frac{B}{|c_{k_0}|e^{2\pi(k_0-\alpha)R}}
    e^{2\pi(s-R)}.
    \]
    The right-hand side tends to zero as $s\to-\infty$. Thus, for
    $s$ sufficiently negative, the straight-line homotopy between
    $b(s,\cdot)$ and
    $c_{k_0}e^{2\pi(k_0-\alpha)s}e^{2\pi i k_0(\cdot)}$
    does not pass through zero. Consequently, for $s$ sufficiently negative 
    \[
    \wind(b(s,\cdot))
    =\wind\left(
      \left(c_{k_0}e^{2\pi(k_0-\alpha)s}e^{2\pi i k_0t}\right)_{t\in[0,1]}
    \right)
    =k_0\geq1+\lfloor\alpha\rfloor.
    \]
    
    Since the constant capping at $p$ has positive orientation, we deduce that for $s \ll 0$ we have
    $$u\cdot (W\times\{0\})=\wind(b(s,\cdot))\geq1+\lfloor\alpha\rfloor.$$

    \end{enumerate}
\end{proof}

\section{Complement of a ball}\label{app: complement of the ball}

\subsection{Acceleration data, complement of a ball}\label{ss: acc. data CP^n - ball}
Let $\Delta\in(0,1)$.

Let $g\fc (0,1)_\Delta\times[0,+\infty)_s\times [0,1]_r\to \R$ be a smooth function, satisfying the following:
\begin{itemize}
    \item  $\frac{\partial g}{\partial \Delta}\geq0$, $\frac{\partial g}{\partial s}\geq0$ and  $\frac{\partial^2}{\partial r^2}g\ge0$;
    \item For every $(\Delta,r)\in(0,1)\times[0,1]$ we have $\lim_{s\to\infty} g(\Delta,s,r)=\left\{\begin{array}{cc}
        +\infty, & r< \Delta, \\
        0, & r\geq\Delta.
    \end{array}\right.$
    
    \item For every $(\Delta,s)\in(0,1)\times[0,+\infty)$ and $r\in (\Delta,1]$ we have $g(\Delta,s,r)<0$, $\frac{\partial}{\partial r} g(\Delta,s,r)\in (-1,0)$ and $\frac{\partial^2}{\partial r^2} g(\Delta,s,r)=0$;
    
    \item There exists a smooth function $\delta\fc(0,1)_\Delta\times[0,+\infty)_s\to \R$ satisfying 
    \begin{itemize}
        \item $\frac{\partial}{\partial s}\delta<0$;
        \item For every $(\Delta,s)\in(0,1)\times[0,+\infty)$ we have $\delta(\Delta,s)=\Delta$ for $s=0$ and $0<\delta(\Delta,s)<\Delta$ otherwise;
        \item For every $\Delta\in(0,1)$ we have  $\lim_{s\to\infty}\delta(\Delta,s)=0$
    \end{itemize}
    such that for every $(\Delta,s)\in(0,1)\times[0,+\infty)$ and $r\in[0,\Delta-\delta(\Delta,s)]$ we have $\frac{\partial}{\partial r} g(\Delta,s,r)=-(s+\frac{1}{2})$, and for every $r\in(\Delta-\delta(\Delta,s),\Delta)$, we have $\frac{\partial^2}{\partial r^2}g>0$.
    
    \item For every $(\Delta,r)\in(0,1)\times [0,1]$ we have $g(\Delta,0,r)<0$.
\end{itemize}

Now, for every $\Delta\in(0,1)$ we use the acceleration datum $(G_i)_{i\in\ \N}$, consisting of autonomous Hamiltonians associated to the complement of an open ball of capacity $\Delta$ in $\CP^n$, given by 
\begin{equation}
    G_\ell(z)=g(\Delta,\ell,\mu(z)),
\end{equation}
for every $z\in \CP^n$ and $\ell\in\Z_{\geq0}.$ The function $g$ defines the homotopies between the Hamiltonians within the acceleration data, as well as the homotopies between acceleration data corresponding to two different capacities. This will be crucial for constructing the restriction maps used in the proof of Theorem~\ref{thm: SH of complement of ball in CPn}. 

Let us describe the generators of the Floer--Morse--Bott complexes for this acceleration datum. Given $\Delta\in(0,1)$ and $\ell\geq0$, the Hamiltonian $G_\ell$ has $\ell+2$ critical submanifolds, which come in three types of families.
\begin{enumerate}
    \item The global maximum at $0$,
    \item $\ell$ spheres, $S^{2n-1}$, consisting of the points where the slope is integral,
    \item and the global minima, forming a copy of $\CP^{n-1}$ at the divisor at infinity.
\end{enumerate}
By choosing a perfect Morse function on each of these critical submanifolds we get the following generators: The first family contributes a single generator which we denote by $\check{x}_{\ell+1}^\ell$; each sphere contributes a pair of generators corresponding to the minimum and maximum of a perfect Morse function, yielding in total $\check{x}_1^\ell,\hat{x}_1^\ell,\ldots,\check{x}_\ell^\ell,\hat{x}_\ell^\ell$, where the hats correspond to maxima and the checks to minima; lastly, the divisor at infinity contributes $n$ generators, which we denote by $\check{x}_{-(n-1)}^\ell,\ldots,\check{x}_0^\ell$, where $\check{x}_i^\ell$ has Morse index $-2i$.

\begin{rem}\label{rem: G_0}
    Note that the differential of the Floer complex of $G_0$ vanishes identically, due to index consideration, as all generators have even index.
\end{rem}

\subsubsection{The Robbin--Salamon index for orbits in $\CP^n\setminus D_\infty$}

According to Claim~\ref{claim: Int B(1)->CPn-D_infty}, the complement of $D_\infty=\CP^{n-1}$ in $\CP^n$ is symplectomorphic to the open ball $\Int B(1)\subset\C^n$ equipped with the standard symplectic form $\omega_0$.

The computations we present here are similar to those in \cite[Section 3.3]{Oancea_survey}, and they will be carried out with respect to the standard trivialization $\tau$ of the tangent bundle of $\Int B(1)\subset\C^n$.

By Remark~\ref{rem: dynamics of a radial Ham} the Hamiltonian vector field $G_\ell=g_\ell\circ \mu$ over $\Int B(1)$ is given by
$$X_{G_\ell}(z)=-2\pi  g'_\ell(\pi\|z\|^2)J_0z,$$ 
$z\in \Int B(1) $, and its Hamiltonian flow is given by
$$
\varphi_{G_\ell}^t(z)=e^{-2\pi g'_\ell(\pi\|z\|^2)J_0 t}z,
$$
for every $t\in\R$ and $z\in \Int B(1) $, where $J_0$ is the standard complex structure on $\R^{2n}$. Therefore, the linearization of the Hamiltonian flow is
$$
d_z\varphi_{G_\ell}^t(Y)
=
e^{-2\pi g'_\ell(\pi\|z\|^2)J_0t}Y
-
4\pi^2 t\,g''_\ell(\pi\|z\|^2)\langle z,Y\rangle\,
e^{-2\pi  g'_\ell(\mu(z))J_0t}J_0z,
$$
for every $t\in \R$, $z\in \Int B(1)$ and $Y\in T_z\R^{2n}$. In particular, along the constant orbit at $\{0\}$ we have $g''_\ell(0)\langle 0,Y\rangle=0$, and therefore
$$
d_0\varphi_{G_\ell}^t(Y)=e^{-2\pi  g'_\ell(0)J_0t}Y,
$$
for every $t\in\R$ and $Y\in T_0\R^{2n}$. Since $g'_\ell(0)=-(\ell+1/2)$, a direct computation shows that 
$$
\mu_{RS}^\tau(x;G_\ell)=n(2\ell+1),
$$
where $x\fc[0,1]\to \Int B(1)$ is the constant loop $x(t)=0$ for every $t\in[0,1]$.

According to our computation, a non-constant $1$-periodic orbit $x\fc [0,1]\to \Int B(1)$ of $G_\ell$ has the form $x(t)=e^{2\pi ktJ_0}z$ for every $t\in[0,1]$, for some $1\leq k\leq \ell$, and for some $z\in \Int B(1)$ satisfying $g'_\ell(\pi\|z\|^2)=-k$.

Take $1\leq k\leq \ell$ and denote
$$
S_k=\{\,z\in \Int B(1)\,|\, g'_\ell(\pi\|z\|^2)=-k\,\}.
$$
Thus $S_k$ is a sphere.  
Take $z\in S_k$, denote $\xi_z=T_z S_k\cap J_0 T_z S_k$, and note that $\xi_z$ is a symplectic linear subspace of $T_z\R^{2n}$ of codimension $2$.
Consider the symplectic trivialization of $T\R^{2n}$ given by
$$
T_w\R^{2n}
=
\R\langle (2\pi \|z\|^2)^{-1}z\rangle
\oplus
\R\langle  2\pi J_0z\rangle
\oplus
\xi_z,
$$
for every $w\in\R^{2n}$.

Note that for every $t\in[0,1]$ we have

\begin{align*}
    d_z\varphi_{G_\ell}^t\big((2\pi \|z\|^2)^{-1}z\big)&=e^{2\pi kJ_0t}\big((2\pi \|z\|^2)^{-1}z\big)\\
    &- 4\pi^2 t\,g''_\ell(\pi\|z\|^2)\langle z,\big((2\pi \|z\|^2)^{-1}z\big)\rangle\,
e^{2\pi  kJ_0t}J_0z\\
&=e^{2\pi kJ_0t}\left((2\pi \|z\|^2)^{-1}z
    -  t\,g''_\ell(\pi\|z\|^2)\,
2\pi J_0z\right),
\end{align*}

Additionally, for every $Y\in (\R\langle z\rangle)^\bot$ we have
$$
d_z\varphi_{G_\ell}^t(Y)=e^{2\pi k tJ_0}Y.
$$

Thus we can present $d_z\varphi_{G_\ell}^t$ in our trivialization as $[d_z\varphi_{G_\ell}^t]=\Psi(t)\cdot\chi(t)$, where $\Psi,\chi\fc[0,1]\to\Sp(2n)$ are given by $\Psi(t)=e^{2\pi k tJ_0}$ and
$$
\chi(t)=
\begin{pmatrix}
    \begin{pmatrix}
        1&0\\-t\,g''_\ell&1
    \end{pmatrix}&\\&\Id_{2n-2}
\end{pmatrix},$$
for every $t\in[0,1]$.

Let $K\fc[0,1]\times[0,1]\to\Sp(2n)$ be the homotopy that connects, with fixed endpoints, the path $\Psi(t)\chi(t)$ to the concatenation of $\Psi(t)$ and $\Psi(1)\chi(t)$, given by
$$
K(s,t)=
\begin{cases}
\Psi\!\left(\dfrac{2t}{s+1}\right)\,\chi(st), & t\le \dfrac{s+1}{2},\\[0.8em]
\Psi(1)\,\chi\!\left((s+2)t-(s+1)\right), & t\ge \dfrac{s+1}{2},
\end{cases}
$$
for every $(s,t)\in[0,1]\times[0,1]$, as in \cite{CFHW_1996_ApSH_II}.

Since the Robbin--Salamon index is homotopy invariant and additive under concatenation, we deduce that
$$
\mu_{RS}([d_z\varphi_{G_\ell}^t])
=
\mu_{RS}(\Psi)+\mu_{RS}(\Psi(1)\chi(\cdot)).
$$

Additionally, since $\Psi(t)=e^{2\pi k tJ_0}$ for every $t\in[0,1]$, and by Corollary~\ref{coroll: RS of rotation}, we know that $\mu_{RS}(\Psi)=2nk$.

As well as, since $\Psi(1)=\Id$, we get from the product, zero, and shear properties of the Robbin--Salamon index, together with Lemma~\ref{lemma: RS of T}, that
\begin{align*}
 \mu_{RS}(\Psi(1)\chi(\cdot))&=\mu_{RS}(\chi)  \\
 &=\mu_{RS}\left(\begin{pmatrix}
        1&0\\-t\,g''_\ell&1
    \end{pmatrix}_{t\in[0,1]}\right)+\mu_{RS}(\Id_{2n-2})\\
    &=-\mu_{RS}\left(\begin{pmatrix}
        1&-t\,g''_\ell\\0&1
    \end{pmatrix}_{t\in[0,1]}\right)+0\\
    &=-\frac{1}{2}\sign{g''_\ell(\|z\|^2)}\\
    &=-\frac{1}{2}.
\end{align*}

Therefore we obtain that
$$
\mu_{RS}^\tau(x;G_\ell)
=
\mu_{RS}([d_z\varphi_{G_\ell}^t])
=
2nk-\tfrac{1}{2}.
$$

\subsubsection{RS-index of orbits on $D_\infty = \CP^{n-1}$}

According to Claim~\ref{claim: mu in O(1)} and the definition of the Hamiltonian $G_\ell$, we can identify a neighborhood of $D_\infty \subset \CP^n$ with a neighborhood $U$ of the zero section of $\cO(1)$ over $D_\infty=\CP^{n-1}$. In these coordinates, $G_\ell$ takes the form
$$G_\ell(w,\lambda) = \frac{-\epsilon}{1+\|\lambda\|_h^2},$$
for some $\epsilon\in(0,1)$, and the Hamiltonian flow is given by
$$\varphi_{G_\ell}^t(w,\lambda) = (w, e^{-2\pi i\epsilon t}\lambda),$$
for every $(w,\lambda) \in U$ and $t \in \R$. 

Over the zero section, the tangent bundle $T\cO(1)|_{\CP^{n-1}}$ decomposes into the sum of the tangent bundle of the zero section and the vertical bundle:
$$T\cO(1)|_{\CP^{n-1}} = T\CP^{n-1} \oplus \cO(1)|_{\CP^{n-1}}.$$
As a conclusion from Remark~\ref{rem: dynamics of a radial Ham}, each point of $D_\infty \subset \CP^n$ is a constant $1$-periodic orbit for $G_\ell$. The linearization of the Hamiltonian flow along such a constant orbit $w_0 \in D_\infty$ can be expressed as
$$d_{w_0}\varphi_{G_\ell}^t(X,Y) = (X, e^{-2\pi i \epsilon t}Y),$$
for every $(X,Y) \in T_{w_0} \CP^{n-1} \oplus \cO(1)_{w_0}$.

By the product and zero axioms of the Robbin--Salamon index, the index of this path of symplectic matrices reduces to the index of the path restricted to the normal bundle:
$$\Psi(t) = D\varphi^t|_{\cO(1)_{w_0}} = e^{-2\pi\epsilon it} \cdot \id.$$
Note that $\frac{d}{dt}\Psi(t) = -2\pi\epsilon i e^{-2\pi\epsilon it} \cdot \id$. The crossing form is given by
$$\omega\left( Y, \dot{\Psi}(t) Y \right) = \left\langle Y, -2\pi\epsilon e^{-2\pi\epsilon it} Y \right\rangle.$$
A time $t \in [0,1]$ is a crossing if and only if $-\epsilon t \in \Z$, which occurs only at $t=0$.

Since the crossing form is negative definite at each crossing (acting on the $2$-dimensional real fiber), the signature is $-2$. Thus, the Robbin--Salamon index is
$$\mu_{RS}(\Psi) = \frac{1}{2}\text{sign}(\Psi(0)) + \sum_{j=1}^\ell \text{sign}(\Psi(t_j)) = \frac{1}{2}(-2) =-1.$$

\subsubsection{RS-index of orbits on $\CP^n\setminus\left(D_1\cup\ldots\cup D_n\cup D_\infty\right)$}

This section focuses on computing the Robbin--Salamon index for $1$-periodic orbits of $G_\ell$ that do not intersect the divisors $D_1,\ldots,D_n,D_\infty$. Recall that Claim~\ref{claim: CPn - all divisors} asserts that there is a  symplectomorphism $\Psi\fc (V\times \T^n,\Omega_0)\to (\CP^n\setminus (D_1\cup\ldots\cup D_n\cup D_\infty),\omega_{FS})$, where
$$ V = \left\{(x_1,\ldots,x_n) \in \R^n \, : \, \sum_{j=1}^n x_j < \frac{1}{\pi} \,\,\text{and}\,\, \forall \, j \in \{1,\ldots,n\} \,\, \text{we have} \,\, x_j > 0 \right\}, $$
also, $\T^n = (\R/\Z)^n$ is the $n$-torus equipped with angular coordinates $(\theta_1,\ldots,\theta_n)$, and $\Omega_0=\pi\sum_{j=1}^n dx_j\wedge d\theta_j$.

Since $\T^n$ is a Lie group, its tangent bundle admits a natural trivialization. This gives rise to a trivialization $\tau$ of $T(V\times \T^n)$, since $V$ is an open subset of $\R^n$. The Robbin--Salamon index in this section will be computed with respect to this trivialization.

Also, by Remark~\ref{rem: mu in CPn - all divisors}, the Hamiltonian $\Psi^*G_\ell$ is given by
$$ \Psi^*G_\ell(x,\theta_1,\ldots,\theta_n) = g_\ell(\pi\|x\|_1), $$
and the Hamiltonian vector field is given by
$$ X_{\Psi^*G_\ell}(x,\theta_1,\ldots,\theta_n) = -g'_\ell(\pi\|x\|_1) \sum_{j=1}^n \partial\theta_j, $$
for every $(x,\theta_1,\ldots,\theta_n) \in V \times \T^n$, where $\|x\|_1 = |x_1| + \cdots + |x_n|$ for every $x = (x_1,\ldots,x_n) \in \R^n$.

Therefore, the flow of $X_{\Psi^*G_\ell}$ is a translation in the $\theta$-variables:
$$ \varphi_{\Psi^*G_\ell}^t(x,\theta_1,\ldots,\theta_n) = (x, \theta_1 - tg'_\ell(\pi\|x\|_1), \ldots, \theta_n - tg'_\ell(\pi\|x\|_1)), $$
for every $(x,\theta_1,\ldots,\theta_n) \in V \times \T^n$ and $t \in \R$.

Thus, a point $(x_0,\theta_0) \in V \times \T^n$ is a fixed point of the time-$1$ map $\varphi_{\Psi^* G_\ell}^1$ if and only if $g'_\ell(\pi\|x_0\|_1) \in \Z$.

Let $(x_0,\theta_0) \in V \times \T^n$ be a fixed point of the time-$1$ map $\varphi_{\Psi^* G_\ell}^1$, and let $x$ be its $1$-periodic orbit. With respect to the trivialization $\tau$, the linearization  of $\{d_{(x_0,\theta_0)}\varphi_{\Psi^* G_\ell}^t\}_{t\in[0,1]}$ along $x$ is represented by the path $\Psi \fc [0,1] \to \Sp(2n)$ given by
$$
\Psi(t) = \begin{pmatrix}
I_n & 0 \\
- t a J_n & I_n
\end{pmatrix},
$$
for $t \in [0,1]$, where $J_n \in \R^{n \times n}$ is the all-ones matrix and $a = \pi g''_\ell(\pi\|x_0\|_1)$. Thus, by Lemma~\ref{lemma: RS of T} and the Shear property of the Robbin--Salamon index, we obtain
$$ \mu_{RS}(\Psi) = -\mu_{RS}(\Psi^T) = -\frac{1}{2}\sign(aJ_n) = -\frac{1}{2}\sign(a). $$
Here we used the fact that $na$ is the only non-zero eigenvalue of $aJ_n$ with multiplicity $1$; hence $\sign(aJ_n) = \sign(na) = \sign(a)$. According to the definition of $g_\ell$, we have $g''_\ell(\pi\|x_0\|_1) > 0$ (see Section~\ref{ss: acc. data}), and thus
$$ \mu_{RS}^\tau(x; \Psi^*G_\ell) = -\frac{1}{2}. $$

\subsection{Floer--Morse--Bott indices}\label{ss: computations of FMB for complement} 
Recall that the Floer--Morse--Bott index of a critical point $p$ of a Morse function $h$ on a critical submanifold $S$, with respect to an autonomous Hamiltonian $H$ satisfying that \textbf{MB} condition and a trivialization $\tau$ along $p$, is given by
$$ \mu_{FMB}^\tau(p;H) = \mu_{RS}^\tau(p;H) + \frac{1}{2}\dim M - \frac{1}{2}\dim S + \ind_h p $$
where $\ind_h p$ is the Morse index of $p$ as a critical point of $h$.

\begin{rem}
Note that the Floer--Morse--Bott index is defined globally modulo $2N$, where $N$ is the minimal first Chern number. We denote this index by $\mu_{FMB}$. This holds for every $1$-periodic orbit of a Hamiltonian satisfying the \textbf{MB} condition, where the index is computed using a trivialization obtained from cappings. The modulo $2N$ ambiguity arises because changing the capping by attaching a sphere results in the index being shifted by twice the first Chern class evaluated on that sphere.
\end{rem}

We now iterate over the orbits for which we computed the Robbin--Salamon index, and calculate the corresponding Floer--Morse--Bott indices of the $1$-periodic orbits of the Hamiltonian $G_\ell$ for $\ell \in \Z_{\geq 0}$ and $\Delta \in (0,1)$.

\begin{enumerate}
    \item For the trivialization $\tau$ over $\CP^n \setminus D_\infty$ discussed in Section~\ref{sss: CP^n - D_infty}, there are two types of $1$-periodic orbits for $G_\ell$: the maximum of $G_\ell$ and the non-constant periodic orbits in $\CP^n \setminus D_\infty$. Their Robbin--Salamon indices were computed in Section~\ref{sss:RS in Cn}:
    \begin{itemize}
        \item The maximum $x$ of $G_\ell$ is a constant and isolated $1$-periodic orbit; thus, the critical submanifold $S$ containing it is a singleton. For this reason, $\dim S=0$ and $\ind_h x = 0$, where $h \colon S \to \R$ is a Morse function. Thus,
        $$ \mu_{FMB}^\tau(x;G_\ell) = n(2\ell+1) + \frac{1}{2} \cdot 2n - \frac{1}{2} \cdot 0 + 0 = 2n(\ell+1). $$

        \item If $x$ is a non-constant orbit of $G_\ell$ lying in the $(2n-1)$-dimensional sphere 
        $$ S_k = \{z \in \CP^n \colon g'_\ell(\mu(z)) = -k\} $$
        for some $1 \leq k \leq \ell$, we showed that $\mu_{RS}^\tau(x;G_\ell) = 2nk - \frac{1}{2}$. Therefore,
        $$ \mu_{FMB}^\tau(x;G_\ell) = 2nk - \frac{1}{2} + \frac{1}{2} \cdot 2n - \frac{1}{2}(2n - 1) + \ind_h x = 2nk + \ind_h x, $$
        where $h$ is a perfect Morse function on the sphere $S_k$. Such a function has exactly two critical points (the minimum and the maximum) with Morse indices $0$ and $2n-1$, respectively. Thus, $\ind_h x \in \{0, 2n-1\}$.

        Thus, modulo $2(n+1)$, which is twice the minimal first Chern number we get that 
        $$\mu_{FMB}(\hat{x}_k^\ell;G_\ell)=2nk+2n-1=2(n+1)(k+1)-2k-3\equiv -2k-3 \mod 2(n+1),$$
        for every $1\leq k\leq \ell$, and 
        $$\mu_{FMB}(\check{x}_k^\ell;G_\ell)=2nk=2(n+1)k-2k\equiv -2k \mod 2(n+1),$$
        for every $0\leq k\leq \ell$.
    \end{itemize}

    \item For the trivialization $\tau$ over $\CP^n \setminus \{0\}$ discussed in Section~\ref{sss: CP^n - 0}, all orbits of $G_\ell$ considered are constant and lie in the divisor $D_\infty = \CP^{n-1}$. Their Robbin--Salamon indices were computed in Section~\ref{sss: RS on CPn-1}. If $h$ is a perfect Morse function on $\CP^{n-1}$ and $w$ is a critical point of $h$, we proved that $\mu_{RS}^\tau(w;G_\ell) = -1$. Therefore,
    $$ \mu_{FMB}^\tau(w;G_\ell) = - 1 + \frac{1}{2} \cdot 2n - \frac{1}{2}(2n - 2) + \ind_h w =  \ind_h w, $$
    where $\ind_h w \in \{0, 2, \dots, 2n-2\}$.

In particular, all the possibilities for $\mu^\tau_{FMB}(w;G_\ell)$ are even numbers, and hence the classes $\mu_{FMB}(w;G_\ell)$ modulo $2(n+1)$ are even.
    \item For the trivialization $\tau$ over $\CP^n \setminus (D_1 \cup \dots \cup D_n \cup D_\infty)$ discussed in Section~\ref{sss: CPn - all divisors}, the $1$-periodic orbits of $G_\ell$ lie on the spheres $S_1, \dots, S_\ell$ as before, but are disjoint from the divisors. According to Section~\ref{sss: RS on CPn minus divsors}, any such orbit $x$ has $\mu_{RS}^\tau(x;G_\ell) = -\frac{1}{2}$. Therefore,
    $$ \mu_{FMB}^\tau(x;G_\ell) = -\frac{1}{2} + \frac{1}{2} \cdot 2n - \frac{1}{2}(2n - 1) + \ind_h x = \ind_h x, $$
    where $h$ is a perfect Morse function on the $(2n-1)$-sphere, yielding $\ind_h x \in \{0, 2n-1\}$.
\end{enumerate}

\subsubsection{Asymptotic intersection number}

The second type of intersections with which we deal is not the intersections of Floer solutions themselves with a divisor, but rather the intersections of their asymptotics.

Let $\Sigma$ be a Riemann surface with a puncture $p$, and let $X$ be a finite dimensional smooth manifold containing a submanifold $Y$ of codimension $2$. Suppose $u\fc \Sigma \to X$ is a smooth map that is asymptotic at $p$ to a point $y \in Y$, and assume that $\im (u) \cap Y = \varnothing$. The \textbf{asymptotic intersection number} of $u$ with $Y$, denoted $u * Y$, is defined to be the intersection number of the capped map $u \# c$ with $Y$ at $y$, where $c$ is the constant capping map at $y$. The orientation of $c$ is taken to be positive if the puncture at $y$ is negative, and negative if the puncture is positive.

Let $R\in\R$ and let us denote by $S_+(R), S_-(R)$ the half-infinite cylinders: 
$$S_+(R)=[R,+\infty)\times S^1\qquad\text{and}\qquad S_-(R)=(-\infty,R]\times S^1.$$

Let $(M^{2n},\omega)$ be a symplectic manifold and let $J$ be a compatible almost complex structure on $(M,\omega)$. Assume that $(M,J)$ is biholomorphic to a product $(W,J')\times (B,i)$ where $(W,J')$ is an $(n-1)$-complex-dimensional complex manifold and $B$ is an open disk in $\C$ centered at the origin equipped with the standard complex structure. 

Let $H\fc M\to \R$ be a Hamiltonian and assume that its Hamiltonian vector field $X_H$ has the form
$$X_H(w,z)=V(w)+2\pi \alpha i z,$$
where $V$ is a smooth vector field on $W$ and $\alpha\in \R\setminus\Z$ is a non-integer real number.

The following result provides a slight generalization of a formula due to Seidel \cite[Eq. (7.22)]{Seidel_Fukaya_A_infty_Strcs_Assoc_to_Lef_Fib. III}.

\begin{prop}
    Let $p\in W$ be a zero of $V$, let $S$ be a half-infinite cylinder and let $u\fc S\to M$ be a Floer solution with respect to $H$ and $J$ that converges to $(p,0)$ and satisfies $\im u \cap (W\times\{0\})=\varnothing$.
    \begin{itemize}
        \item If $S=S_+(R)$, for some $R\in \R$, then $u* (W\times\{0\})\geq-\lfloor\alpha\rfloor$.
        \item If $S=S_-(R)$, for some $R\in \R$, then $u* (W\times\{0\})\geq 1+\lfloor\alpha\rfloor$.
    \end{itemize}
\end{prop}

For completeness, and to ensure the compatibility of our sign conventions, the proof Proposition~\ref{prop: asymp_intersection} is provided in Appendix~\ref{app: proof of Seidel's lemma}.

Let us present some examples that will be useful later:

\begin{exam}\label{exam: computations using Seidel's lemma for complement}
    Let $\Delta \in (0,1)$ and $\ell \in \Z_{\geq 0}$. Let $S$ be a half-infinite cylinder and let $u \fc S \to \CP^n$ be a Floer solution with respect to $H_\ell$ and the standard complex structure $J$ on $\CP^n$. Assume that $\im u \cap (D_1 \cup \dots \cup D_n \cup D_\infty) = \varnothing$.
    
    \begin{itemize}
        \item Assume that the asymptote of $u$ lies on $D_\infty = \CP^{n-1}$. According to Claim~\ref{claim: mu in O(1)} and the definition of the Hamiltonian $G_\ell$, we can identify a neighborhood of $D_\infty \subset \CP^n$ with a neighborhood $U$ of the zero section of the line bundle $\cO(1) \to \CP^{n-1}$. In these coordinates, $G_\ell$ takes the form
        $$G_\ell(w,\lambda) = \frac{-\epsilon}{1 + \|\lambda\|_h^2},$$
        for some $\epsilon\in(0,1)$, and the Hamiltonian vector field is given by
        $$X_{G_\ell}(w,\lambda) = -2\pi i \epsilon \lambda,$$
        for every $(w,\lambda) \in U$. 

        Since $\cO(1)$ is a holomorphic bundle, there is a neighborhood of the asymptote of $u$ which is biholomorphic to a product $W \times B$, where $W$ is an open set in $\CP^{n-1}$ and $B \subset \C$ is a small ball, both equipped with the standard complex structure. Moreover, we can choose this biholomorphism such that $X_{G_\ell}$ takes the form
        $$X_{G_\ell}(w,z) = -2\pi i \epsilon z,$$
        for every $(w,z) \in W \times B$.

        Thus, by Proposition~\ref{prop: asymp_intersection} and the fact that $\im u \cap D_\infty = \varnothing$, if $u$ has an asymptote at $+\infty$ on $D_\infty$, then
        $$u \cdot D_\infty = u \cdot (W \times \{0\}) \geq 1.$$
        Similarly, if $u$ has an asymptote at $-\infty$ on $D_\infty$, then
        $$u \cdot D_\infty = u \cdot (W \times \{0\}) \geq 0.$$

        \item Assume that the asymptote of $u$ is the constant loop $\check{x}_{\ell+1}^\ell$ located at $0 \in \Int B(1) \subset \CP^n$. Let $B_1, \dots, B_n \subset \C$ be small standard disks centered at the origin, such that $P = B_1 \times \dots \times B_n \subset B^{}_{\Delta} \subset \CP^n$. By definition, on $P$ the Hamiltonian $G_\ell$ takes the form
        $$G_\ell(z_1, \dots, z_n) = -\pi (\ell+1/2) (|z_1|^2 + \dots + |z_n|^2).$$
        Thus, the Hamiltonian vector field is given by
        $$X_{G_\ell}(z_1, \dots, z_n) = \sum_{j=1}^n (2\pi i (\ell+1/2) z_j).$$

        Let $1 \leq j \leq n$. By Proposition~\ref{prop: asymp_intersection} and the fact that $\im u \cap D_j = \varnothing$, if $u$ is asymptotic at $+\infty$ to $0$, then
        $$u \cdot D_j = u \cdot (B_j \times \{0\}) \geq -\lfloor (\ell+1/2) \rfloor=-\ell.$$
        Similarly, if $u$ is asymptotic at $-\infty$ to $0$, then
        $$u \cdot D_j = u \cdot (B_j \times \{0\}) \geq \lfloor (\ell+1/2) \rfloor + 1 = \ell+1. $$
    \end{itemize}
\end{exam}

\subsection{Obstructions}
  Theorems~\ref{thm: diff_obst_for_complement} and \ref{thm: continuation_obst_for_complement} below focus on the non-existence of certain Floer and continuation flowlines with cascades associated with our acceleration data. The main tool in their proofs is the positivity of intersection discussed in the previous section. To ensure that the conditions for positivity of intersection hold, we should verify that for a given broken bubbled flowline $u$ with bubbles $\{B_\alpha\}_{\alpha\in A}$ and $m$ cascades $u_1, \ldots, u_m$ having asymptotes $q_-^k$ and $q_+^k$ at $-\infty$ and $+\infty$ respectively, any of these asymptotes that is not constant does not intersect the divisors $D_1, \ldots, D_n$ for every $1 \leq k \leq m$. Given Hamiltonians $G_\ell$ and $G'_{\ell'}$ from our acceleration data associated with balls of capacity $\Delta$ and $\Delta'$ respectively, for some $\ell, \ell' \in \Z_{\geq 0}$, we note that any element of the unitary group $\U(n)$ induces a symplectomorphism of $\CP^n$ that preserves the divisor $D_\infty$, the Hamiltonians $G_\ell$ and $G'_{\ell'}$, and all the Floer--Morse--Bott indices of their $1$-periodic orbits. We thus deduce that there exists a unitary matrix $U$ such that the transformed divisors $U(D_1) \cup \dots \cup U(D_n)$ do not intersect any of the non-constant $1$-periodic orbits in the collection $\{q_-^1, q_+^1, \dots, q_-^m, q_+^m\}$, and do not contain the images of $u_1, \ldots, u_m$ or $\{B_\alpha\}_{\alpha \in A}$. Therefore, we can relabel $U(D_1), \ldots, U(D_n)$ as $D_1, \ldots, D_n$, allowing us to assume this configuration in the proof of each of these theorems.

Throughout this section and under the assumption above, we utilize two distinct trivializations of the tangent bundle of $\CP^n$ along the $1$-periodic orbits of $G_\ell$ and $G'_{\ell'}$, denoted by $\tau_B$ and $\tau_T$. 

The trivialization $\tau_B$ (where the subscript $B$ refers to the \textbf{ball}) is defined as follows: for the constant generators $$\check{x}_{\ell+1}^\ell, \check{x}_{0}^\ell, \ldots, \check{x}_{-(n-1)}^\ell,\check{y}_{\ell+1}^{\ell'}, \check{y}_{0}^{\ell'}, \ldots, \check{y}_{-(n-1)}^{\ell'},$$
the trivialization $\tau_B$ is obtained from the constant cappings; for the non-constant generators $$\check{x}_1^\ell, \hat{x}_1^\ell, \ldots, \check{x}_\ell^\ell, \hat{x}_\ell^\ell,\check{y}_1^{\ell'}, \hat{y}_1^{\ell'}, \ldots, \check{y}_{\ell'}^{\ell'}, \hat{y}_{\ell'}^{\ell'},$$ the trivialization $\tau_B$ is the one for $\CP^n \setminus D_\infty$ discussed in Section~\ref{sss: CP^n - D_infty}. 

The trivialization $\tau_T$ (where the subscript $T$ refers to the \textbf{torus}) is defined as follows: for the constant generators $$\check{x}_{\ell+1}^\ell, \check{x}_{0}^\ell, \ldots, \check{x}_{-(n-1)}^\ell,\check{y}_{\ell+1}^{\ell'}, \check{y}_{0}^{\ell'}, \ldots, \check{y}_{-(n-1)}^{\ell'},$$
the trivialization $\tau_T$ is obtained from the constant cappings; for the non-constant generators $$\check{x}_1^\ell, \hat{x}_1^\ell, \ldots, \check{x}_\ell^\ell, \hat{x}_\ell^\ell,\check{y}_1^{\ell'}, \hat{y}_1^{\ell'}, \ldots, \check{y}_{\ell'}^{\ell'}, \hat{y}_{\ell'}^{\ell'},$$
the trivialization $\tau_T$ is the one for $\CP^n \setminus (D_1 \cup \ldots \cup D_n \cup D_\infty)$ discussed in Section~\ref{sss: CPn - all divisors}.

The following proposition concerns cascades and bubbles contained in the divisor $D_\infty$. We prove that for our acceleration data, all of these components are constant. Fix $\ell,\ell'\in \Z_{\geq0}$ with $\ell\leq \ell'$ and $0<\Delta\leq \Delta'<1$. To use this result for concatenation of bubbled broken flowlines with cascades, we consider a finite sequence of Hamiltonians $(F_i)_{i=0}^N$ on $\CP^n$, given by $F_i=G_{a_i}^{\Delta_i}=g(\Delta_i,a_i,\mu(\cdot))$ for every $0\leq i\leq N$. Here, the indices form a finite sequence $\ell=a_0\leq \cdots\leq a_N=\ell'$ in $\Z_{\geq0}$, with corresponding parameters $\Delta=\Delta_0\leq\cdots\leq\Delta_N=\Delta'$. We then consider $N$ homotopies, one from $F_i$ to $F_{i+1}$ for each $i\in\{0,\ldots,N-1\}$, defined via the parameters of $g$, and we use $J_0$ for the almost complex structure. We denote $F_-=F_0$ and $F_+=F_N$.

\begin{prop}\label{prop: contained components complement}
Let $u$ be a concatenation of broken bubbled flowlines with
cascades with respect to the Hamiltonians $(F_i)_{i=0}^N$, with first negative asymptote at a generator $q_-$ for $F_-$ and last positive asymptote at a generator $q_+$ for
$F_+$. Denote the collection of the cascades of $u$ by $\cU$ and the collection of the
sphere bubbles by $\cB$. Let $\tau$ be any symplectic trivialization of $\det TM$ along the asymptotes of the cascades of $u$. If
    $$ \mu^\tau_{FMB}(q_+;F_+) - \mu^\tau_{FMB}(q_-;F_-) + 2\sum_{v\in\cU} c_1^\tau(v) + 2\sum_{B\in\cB} c_1(B) \;\in\;\{-1,0,1\}, $$
then every cascade or bubble of $u$ that is contained in $D_\infty$ is constant.
\end{prop}

\begin{proof} First, denote
\[\delta=\mu^\tau_{FMB}(q_+;F_+) - \mu^\tau_{FMB}(q_-;F_-) + 2\sum_{v\in\cU} c_1^\tau(v) + 2\sum_{B\in\cB} c_1(B), \]
and assume that $\delta\in\{-1,0,1\}$. We will show that every cascade or bubble of $u$ that is contained in $D_\infty$ is constant.

Second, note that $\delta$ is independent of the trivialization $\tau$, thus we can assume that $\tau=\tau_T$, where $\tau_T$ is defined on page~\pageref{p: def tau_T}.

We denote by $\cU_\infty$ and $\cB_\infty$ the subcollections of $\cU$ and $\cB$, respectively, whose
images are contained in $D_\infty$.
Also, set
\[
\begin{split}
C^{\mathrm{tot}}=\sum_{v\in\cU}c_1^{\tau_T}(v)
    +\sum_{B\in\cB}c_1(B),\qquad\qquad
C^{\mathrm{rem}}=\sum_{v\in\cU\setminus\cU_\infty}c_1^{\tau_T}(v)
+\sum_{B\in\cB\setminus\cB_\infty}c_1(B),
\end{split}
\]
and
\[
R:=C^{\mathrm{tot}}-C^{\mathrm{rem}}
=\sum_{v\in\cU_\infty}c^{\tau_T}_1(v)+\sum_{B\in\cB_\infty}c_1(B).
\]
Now, note that the Hamiltonian vector field vanishes on $D_\infty$. Thus, every cascade $v$
contained in $D_\infty$ extends to a $J_0$-holomorphic sphere $\bar v$, by removal of singularities.
Since for every $v\in \cU_\infty$, over the asymptotes of $v$, the trivialization $\tau$ is obtained from the constant cappings, thus we deduce that $c_1^\tau(v)=c_1(\bar v)$. Since the minimal first Chern number of $\CP^n$ is $n+1$, we deduce that $R\in(n+1)\Z$. Also, since for every $v\in \cU_\infty$, the sphere $\bar v$ is pseudo-holomorphic, we deduce together with monotonicity that $c_1(\bar v)=(n+1)\omega(\bar v)\geq0$. Similarly, for every $B\in \cB_\infty$ we have $c_1(B)\geq0$. Therefore, $R\in (n+1)\Z_{\geq0}$. We shall show that $R<n+1$, and thus deduce that $R=0$.

Let us separate into several cases:
\begin{itemize}
\item Assume $\cU_\infty=\cU$ or all the cascades of $u$ are constant at the origin. In these cases the first Chern numbers of all the elements of $\cU\cup\cB$ are nonnegative, thus we can write
$C^{\mathrm{tot}}=(n+1)d$ for some $d\geq0$. Also, in these cases $q_-,q_+$ are constant and are located either at the origin or contained in $D_\infty$, thus their Floer--Morse--Bott indices are
even, and hence $\delta=0$.

If $q_-,q_+$ are located at the origin, their Floer--Morse--Bott indices are $2n(\ell+1)$ and $2n(\ell'+1)$, respectively, and hence the fact $\delta=0$ implies that
\[
(n+1)d=n(\ell-\ell')\leq0,
\]
which implies that $d=0$.

Otherwise, assume that $q_-,q_+$ are on $D_\infty$, they have the form $q_-=\check{x}_i^\ell$ and $q_+=\check{y}_j^{\ell'}$ for some $-(n-1)\leq i,j\leq0$, and we computed that $\mu_{FMB}^{\tau}(q_-;F_-)=-2i$ and $\mu_{FMB}^{\tau}(q_+;F_+)=-2j$. Thus, since $\delta=0$, we get that
\[
(n+1)d=j-i\leq n-1,
\]
which implies that $d=0$.

Thus in both cases we have $d=0$, and in particular $R=0$. Since each term defining $R$ is nonnegative and is positive for a
non-constant component, it follows that all the components contained in $D_\infty$ are constant.

\item Assume that $\cU\setminus \cU_\infty\neq \varnothing$ and not all the cascades are constant and located at the origin.

As explained before, we can assume without loss of generality that the divisors $D_1,\ldots,D_n$ do not intersect any of the non-constant asymptotes of the cascades of $u$, and do not contain the images of its non-constant cascades or bubbles. We may also assume that these divisors do not intersect any of the constant asymptotes in $D_\infty$.

Fix $0\leq i\leq N$ and a generator $q$ of $F_i$. For now, for brevity, denote $a=a_i$. Let us denote by $\mathcal D(q)$ the collection of the divisors from $D_1,\ldots,D_n,D_\infty$ that intersect $q$. For every $D\in\mathcal D(q)$, write the normal Hamiltonian vector
field as $2\pi i\beta_D z$, as in Example~\ref{exam: computations using Seidel's lemma for complement}, and define
\[
e_-(q)=\sum_{D\in\mathcal D(q)}(1+\lfloor\beta_D\rfloor),
\qquad
e_+(q)=-\sum_{D\in\mathcal D(q)}\lfloor\beta_D\rfloor.
\]
By Proposition~\ref{prop: asymp_intersection}, these are lower bounds for
the sums of asymptotic intersections at a negative and a positive end,
respectively, of all the elements $v\in \cU\setminus\cU_\infty$ which are not constant at the origin. Note that $\beta_{D_\infty}\in(-1,0)$, and at the origin, for every $j\in\{1,\ldots,n\}$ we have $\beta_{D_j}=a+\frac12$. As usual, denote by $\check{x}_{-(n-1)}^{a},\ldots,\check{x}_{a+1}^{a},\hat{x}_1^{a},\ldots,\hat{x}_{a}^{a}$ the generators associated to $F_i$ with a suitable choice of Morse functions on the critical submanifolds. Writing $\mu_T(q)=\mu_{FMB}^{\tau_T}(q;F_i)$, we obtain:
\begin{center}
\small
\setlength{\tabcolsep}{3pt}
\renewcommand{\arraystretch}{1.3}
\begin{tabular}{@{}c@{\quad}c@{\quad}c@{\quad}c@{\quad}c@{\quad}c@{\quad}c@{}}
\toprule
$q$ & Location or range & $\mu_T(q)$ & $e_-$ & $e_+$
    & $\mu_T(q)-2e_-$ & $\mu_T(q)+2e_+$\\
\midrule
$\check{x}_{a+1}^{a}$ & $\text{origin}$ & $2n(a+1)$ & $n(a+1)$ & $-na$ & $0$ & $2n$\\
$\check{x}_{j}^{a}$ & $1\leq j\leq a$ & $0$ & $0$ & $0$ & $0$ & $0$\\
$\hat{x}_{j}^{a}$ & $1\leq j\leq a$ & $2n-1$ & $0$ & $0$ & $2n-1$ & $2n-1$\\
$\check{x}_{j}^{a}$ & $-(n-1)\leq j\leq0$ & $-2j$ & $0$ & $1$ & $-2j$ & $2-2j$\\
\bottomrule
\end{tabular}
\end{center}
In particular,
\[
0\leq\mu_T(q)-2e_-(q)\leq2n-1,
\qquad 0\leq\mu_T(q)+2e_+(q)\leq2n.\tag{$\ast$}
\]

Omit the cascades from $\cU_\infty$, the bubbles from $\cB_\infty$, and constant cascades at the origin. By the assumption, at least one cascade remains. We denote by $u_1,\ldots,u_m$ the remaining cascades in their order of appearance in $u$.

Recall that $F_i=G^{\Delta_i}_{a_i}$ for every $i\in\{0,\ldots,N\}$, and for every $0\leq i\leq j\leq N$ it holds that $a_i\leq a_j$.
Given any two consecutive remaining cascades, if the adjacent asymptotes lie in $D_\infty$, they contribute at least $1+0=1$ to the asymptotic intersection with $D_\infty$. If they lie at the origin, they contribute at least $-na_i+n(a_j+1)\geq n$, for some $0\leq i\leq j\leq N$, to the sum of the asymptotic intersections with $D_1,\ldots,D_n$. Otherwise, their asymptotic intersection contributions are zero. The same estimates apply to consecutive cascades between consecutive broken flowlines. Therefore, we have
\begin{equation*}
    u_k\mathord{\raisebox{-0.8ex}{$\big\vert$}}_{[0,+\infty)\times S^1} * D \,\,+\,\, u_{k+1}\mathord{\raisebox{-0.8ex}{$\big\vert$}}_{(-\infty,0]\times S^1} * D \;\ge\; 0,
\end{equation*}
for every $1\leq k\leq m-1$, where $D=D_1\cup\ldots\cup D_n\cup D_\infty$. Therefore, together with positivity of intersections, and Theorem~\ref{thm: c1 vs intersection number} we have
\begin{equation*}
    C^{\mathrm{rem}}\;\geq\;
    \left(u_1\mathord{\raisebox{-0.8ex}{$\big\vert$}}_{(-\infty,0]\times S^1} * D \,\,+\,\, u_m\mathord{\raisebox{-0.8ex}{$\big\vert$}}_{[0,\infty)\times S^1} * D \right).
\end{equation*}

Assume that $u_1$ is either a Floer solution for $F_i$ or a continuation trajectory from $F_i$ to $F_{i+1}$, and $u_m$ is either a Floer solution of $F_j$ or a continuation trajectory from $F_{j-1}$ to $F_j$, for some $0\leq i\leq j\leq N$.

Denote by $q$ the negative asymptote of $u_1$ and by $q'$ the positive asymptote of $u_m$. Any omitted cascades before $u_1$ or after $u_m$ are contained in $D_\infty$ or constant at the origin. Since $a_0\leq\ldots\leq a_N$, from the table above we deduce that
\[
e_-(q_-;F_-)\;\le\;e_-(q;F_i)\;\le\;
u_1\mathord{\raisebox{-0.8ex}{$\big\vert$}}_{(-\infty,0]\times S^1} * D,
\]
and
\[
e_+(q_+;F_+)\;\le\;e_+(q';F_j)\;\le\;
u_m\mathord{\raisebox{-0.8ex}{$\big\vert$}}_{[0,\infty)\times S^1} * D.
\]

Therefore, we get
\[
C^{\mathrm{rem}}\geq e_-(q_-;F_-)+e_+(q_+;F_+).
\]
Consequently, from $(\ast)$ we deduce that
\[
\begin{split}
\delta&=\mu_T(q_+) - \mu_T(q_-) + 2C^{\mathrm{tot}}\\
&=\mu_T(q_+) - \mu_T(q_-) + 2C^{\mathrm{rem}}+2R\\
&\geq\mu_T(q_+) - \mu_T(q_-) + 2\left(e_-(q_-;F_-)+e_+(q_+;F_+)\right)+2R\\
&=
\bigl(\mu_T(q_+)+2e_+(q_+;F_+)\bigr)
-\bigl(\mu_T(q_-)-2e_-(q_-;F_-)\bigr)+2R\\
&\geq-(2n-1)+2R.
\end{split}
\]
Thus
\[
2R\leq\delta+2n-1\leq1+2n-1=2n.
\]
Hence $0\leq R\leq n<n+1$, and since $R\in(n+1)\Z_{\geq0}$ we have $R=0$. Thus, as in the previous case, since each term defining $R$ is nonnegative and is positive for a
non-constant component, it follows that all the components contained in $D_\infty$ are constant.

\end{itemize}
\end{proof}

The first list of obstructions we present here focuses on the differential:

\begin{thm}\label{thm: diff_obst_for_complement}
    Let $0 < \Delta < 1$ and $\ell \in \Z_{\geq 0}$. Define $G_\ell \fc \CP^n \to \R$ by $G_\ell(z) = g(\Delta, \ell, \mu(z))$ for every $z \in \CP^n$. Let $q_+$ and $q_-$ be critical points of the Morse functions associated with the critical submanifolds with respect to $G_\ell$. 
    
    Suppose $u$ is a broken bubbled Floer flowline with cascades connecting $q_-$ to $q_+$, with cascades $u_1, \dots, u_m$ and a finite collection $\{B_\alpha\}_{\alpha\in A}$ of bubbles. Let $\tau$ be a trivialization of $\det TM$ along the asymptotics of $u_1,\ldots,u_m$. If 
    $$ \mu^\tau_{FMB}(q_+;G_\ell) - \mu^\tau_{FMB}(q_-;G_\ell) + 2\sum_{k=1}^m c_1^\tau(u_k) + 2\sum_{\alpha\in A} c_1(B_\alpha) = 1, $$
    then the pair $(q_-, q_+)$ cannot be any of the following:
    \begin{itemize}
        \item $(\hat{x}^\ell_i, \hat{x}^\ell_j)$ for every $1 \leq i, j \leq \ell$.
        \item $(\check{x}^\ell_i, \check{x}^\ell_j)$ for every $-(n-1) \leq i, j \leq \ell + 1$.
        \item $(\check{x}^\ell_i, \hat{x}^\ell_j)$ for every $-(n-1) \leq i \leq \ell + 1$ and $1 \leq j \leq \ell$, except where $i=j$ and $n=1$, in which case $u$ does not intersect $D_\infty$.
        \item $(\hat{x}^\ell_i, \check{x}^\ell_j)$ for every $1 \leq i \leq \ell$ and $-(n-1) \leq j \leq \ell + 1$, except where $j \in \{i+1, i-n\}$. Moreover, if $j=i+1$ then $u$ does not intersect $D_\infty$.
    \end{itemize}
\end{thm}

\begin{rem}\label{rem: independent of dim formula by tau_for_complement}
    The expression 
    $$ \mu^\tau_{FMB}(q_+;G_\ell) - \mu^\tau_{FMB}(q_-;G_\ell) + 2\sum_{i=1}^m c_1^\tau(u_i) $$
    is independent of $\tau$.
\end{rem}

This remark allows us to replace the trivialization $\tau$ from the formulation of Theorem~\ref{thm: diff_obst_for_complement} with $\tau_B$ or $\tau_T$ as needed during the proof. We will do the same in the proofs of Theorem~\ref{thm: continuation_obst_for_complement} as well.

\begin{proof}[Proof of Theorem \ref{thm: diff_obst_for_complement}] 
As explained before, we can assume without loss of generality that the divisors $D_1, \ldots, D_n$ do not intersect any of the non-constant asymptotes of $u_1, \ldots, u_m$, and do not contain the images of $u_1, \ldots, u_m$ and $\{B_\alpha\}_{\alpha \in A}$. Moreover, Proposition \ref{prop: contained components complement} guarantees that all cascades and bubbles contained in $D_\infty$ are constant and hence can be safely ignored in the proof.

We treat each case separately:
\begin{itemize}
    \item Let $1\leq i,j\leq \ell$ and assume that $(q_-,q_+)=(\hat{x}^\ell_i,\hat{x}^\ell_j)$.
    
    As computed in Section~\ref{ss: computations of FMB for complement}, both $\mu_{FMB}(\hat{x}_i^\ell;G_\ell)$ and $\mu_{FMB}(\hat{x}_j^\ell;G_\ell)$ are odd modulo $2(n+1)$. In particular, the difference 
    $$\mu_{FMB}(\hat{x}_j^\ell;G_\ell) - \mu_{FMB}(\hat{x}_i^\ell;G_\ell) $$
    is even, and hence not equal to $1$ modulo $2(n+1)$. Proposition~\ref{prop: traj and modular indices} then leads to a contradiction.

    \item Let $-(n-1)\leq i,j\leq \ell+1$ and assume that $(q_-,q_+)=(\check{x}^\ell_i,\check{x}^\ell_j)$.

    As computed in Section~\ref{ss: computations of FMB for complement}, both $\mu_{FMB}(\check{x}_i^\ell;G_\ell)$ and $\mu_{FMB}(\check{x}_j^\ell;G_\ell)$ are even modulo $2(n+1)$. In particular, the difference 
   $$\mu_{FMB}(\check{x}_j^\ell;G_\ell) - \mu_{FMB}(\check{x}_i^\ell;G_\ell) $$
    is even, and hence not equal to $1$ modulo $2(n+1)$. Proposition~\ref{prop: traj and modular indices} then leads to a contradiction.

    \item Let $-(n-1)\leq i\leq \ell+1$ and $1\leq j\leq \ell$, and assume that $(q_-,q_+)=(\check{x}^\ell_i,\hat{x}^\ell_j)$. 
    
    We should separate into two cases:
    \begin{enumerate}
    \item Assume $i=\ell+1$.  Let us use the trivialization $\tau_T$. 

With respect to this trivialization, we find that $\mu_{FMB}^{\tau_T}(\hat{x}_j^\ell;G_\ell) = 2n-1$ and $\mu_{FMB}^{\tau_T}(\check{x}_i^\ell;G_\ell) = 2n(\ell+1)$. Thus, by the positivity of intersections from Remark~\ref{rem: positivity of intersections for our acc. data} and Example~\ref{exam: computations using Seidel's lemma for complement}, we deduce that 
$$ \sum_{k=1}^m c_1^\tau(u_k) + \sum_{\alpha\in A} c_1(B_\alpha) =u\cdot (D_1\cup\ldots\cup D_n\cup D_\infty)\geq u\cdot (D_1\cup\ldots\cup D_n\cup D_\infty)\geq n(\ell+1). $$
It follows that
\begin{align*}
    1 &= \mu^\tau_{FMB}(q_+;G_\ell) - \mu^\tau_{FMB}(q_-;G_\ell) + 2\sum_{k=1}^m c_1^\tau(u_k) + 2\sum_{\alpha\in A} c_1(B_\alpha) \\
    &\geq 2n-1-2n(\ell+1)+2n(\ell+1)\\
    &=2n-1.
\end{align*}
For $n \in \N$, the only possibility for this inequality to hold is if $n=1$. Moreover, since the inequality is actually an equality in this case, we deduce from Theorem~\ref{thm: c1 vs intersection number} and that fact that $\PD(D_1\cup D_\infty)=c_1$ that 
$$ u \cdot (D_1 \cup D_\infty) = \sum_{k=1}^m c_1^\tau(u_k) + \sum_{\alpha\in A} c_1(B_\alpha) = \ell+1. $$
Since $u\cdot D_1\geq \ell+1$ and $u\cdot D_\infty\geq0$ we deduce that $u\cdot D_\infty=0$.

Next we move to the trivialization $\tau_B$. 
With respect to this trivialization, since $n=1$ we find that $\mu_{FMB}^\tau(\hat{x}_j^\ell;G_\ell) = 2j+2-1=2j+1$ and $\mu_{FMB}^\tau(\check{x}_i^\ell;G_\ell) = 2(\ell+1)$. Thus, by the positivity of intersections from Remark~\ref{rem: positivity of intersections for our acc. data}, we deduce:
$$ \sum_{k=1}^m c_1^\tau(u_k) + \sum_{\alpha\in A} c_1(B_\alpha) =2 u \cdot D_\infty=0. $$
It follows that
\begin{align*}
     1 &= \mu^\tau_{FMB}(q_+;G_\ell) - \mu^\tau_{FMB}(q_-;G_\ell) + 2\sum_{k=1}^m c_1^\tau(u_k) + 2\sum_{\alpha\in A} c_1(B_\alpha) \\
     &=  2j+1 -2(\ell+1),
\end{align*}
and hence $\ell+1= j\leq \ell$, and this is a contradiction.

    \item Assume that $1 \leq i \leq \ell$. Let us use the trivialization $\tau_T$. 

With respect to this trivialization, we find that $\mu_{FMB}^{\tau_T}(\hat{x}_j^\ell;G_\ell) = 2n-1$ and $\mu_{FMB}^{\tau_T}(\check{x}_i^\ell;G_\ell) = 0$. Thus, by the positivity of intersections from Remark~\ref{rem: positivity of intersections for our acc. data} and Example~\ref{exam: computations using Seidel's lemma}, we deduce that 
$$ \sum_{k=1}^m c_1^\tau(u_k) + \sum_{\alpha\in A} c_1(B_\alpha) \geq 0. $$
It follows that
$$ 1 = \mu^\tau_{FMB}(q_+;G_\ell) - \mu^\tau_{FMB}(q_-;G_\ell) + 2\sum_{k=1}^m c_1^\tau(u_k) + 2\sum_{\alpha\in A} c_1(B_\alpha) \geq 2n-1. $$
For $n \in \N$, the only possibility for this inequality to hold is if $n=1$. Moreover, since the inequality is actually an equality in this case, we deduce from Theorem~\ref{thm: c1 vs intersection number} and that fact that $\PD(D_1\cup D_\infty)=c_1$ that 
$$ u \cdot (D_1 \cup D_\infty) = \sum_{k=1}^m c_1^\tau(u_k) + \sum_{\alpha\in A} c_1(B_\alpha) = 0. $$
When $n=1$, the divisors $D_1$ and $D_\infty$ are simply the poles of $S^2 \cong \CP^1$. The positivity of intersections implies that $u$ does not intersect them; in particular, $u$ preserves the winding of the orbits $\check{x}_i^\ell$ and $\hat{x}_j^\ell$, which implies $i=j$.

\item Assume that $-(n-1) \leq i \leq 0$. Let us use the trivialization $\tau_B$. 

With respect to this trivialization, we find that $\mu_{FMB}^\tau(\hat{x}_j^\ell;G_\ell) = 2nj+2n-1$ and $\mu_{FMB}^\tau(\check{x}_i^\ell;G_\ell) = -2i$. Thus, by the positivity of intersections from Remark~\ref{rem: positivity of intersections for our acc. data} and Example~\ref{exam: computations using Seidel's lemma for complement}, we deduce:
$$ \sum_{k=1}^m c_1^\tau(u_k) + \sum_{\alpha\in A} c_1(B_\alpha) \geq c_1^\tau(u_1) \geq (n+1) (u_1 \cdot D_\infty) \geq (n+1)\cdot 0=0. $$
It follows that
\begin{align*}
     1 &= \mu^\tau_{FMB}(q_+;G_\ell) - \mu^\tau_{FMB}(q_-;G_\ell) + 2\sum_{k=1}^m c_1^\tau(u_k) + 2\sum_{\alpha\in A} c_1(B_\alpha) \\
     &\geq  2nj+2n-1 + 2i  \\
     &\geq 2n+2n-1-2(n-1) \\
     &= 2n+1 \\
     &\geq 3,
\end{align*}
which is a contradiction.
    \end{enumerate}

    \item Let $1\leq i\leq \ell$ and $-(n-1)\leq j\leq \ell+1$, and assume that $(q_-,q_+)=(\hat{x}_i^\ell,\check{x}_j^\ell)$.
    We should separate into two cases:
    
    \begin{enumerate}
    \item Assume $j=\ell+1$.  Let us first consider the trivialization $\tau_B$.

With respect to this trivialization, we find that $\mu_{FMB}^\tau(\check{x}_j^\ell;G_\ell) = 2n(\ell+1)$ and $\mu_{FMB}^\tau(\hat{x}_i^\ell;G_\ell) = 2ni+2n-1$. Thus Theorem~\ref{thm: c1 vs intersection number} implies that
$$ \sum_{k=1}^m c_1^\tau(u_k) + \sum_{\alpha\in A} c_1(B_\alpha) = (n+1) u \cdot D_\infty. $$
It follows that
\begin{align*}
1 &= \mu^\tau_{FMB}(q_+;G_\ell) - \mu^\tau_{FMB}(q_-;G_\ell) + 2\sum_{k=1}^m c_1^\tau(u_k) + 2\sum_{\alpha\in A} c_1(B_\alpha) \\
&= 2n(\ell+1) - 2ni - 2n + 1 + 2(n+1)u \cdot D_\infty \\
&= 2n(\ell-i) + 1 + 2(n+1)u \cdot D_\infty,
\end{align*}
yielding the relation
$$ u \cdot D_\infty = \frac{n(i-\ell)}{n+1}. $$
By the positivity of intersections (Remark~\ref{rem: positivity of intersections for our acc. data}), we have $u \cdot D_\infty \geq 0$, which implies that $i-\ell$ must be a multiple of $n+1$. Writing $i-\ell = k(n+1)$ for some $k \in \Z_{\geq 0}$, it follows that $u \cdot D_\infty = nk$.

Next we move to the trivialization $\tau_T$. 

With respect to this trivialization, we find that $\mu_{FMB}^\tau(\check{x}_j^\ell;G_\ell) = 2n(\ell+1)$ and $\mu_{FMB}^\tau(\hat{x}_i^\ell;G_\ell) = 2n-1$. Thus Theorem~\ref{thm: c1 vs intersection number}, Remark~\ref{rem: positivity of intersections for our acc. data} and Example~\ref{exam: computations using Seidel's lemma for complement} imply that
$$ \sum_{k=1}^m c_1^\tau(u_k) + \sum_{\alpha\in A} c_1(B_\alpha) = u \cdot (D_1 \cup \ldots \cup D_n \cup D_\infty) \geq -n\ell+nk. $$
It follows that
\begin{align*}
1 &= \mu^\tau_{FMB}(q_+;G_\ell) - \mu^\tau_{FMB}(q_-;G_\ell) + 2\sum_{k=1}^m c_1^\tau(u_k) + 2\sum_{\alpha\in A} c_1(B_\alpha) \\
&= 2n(\ell+1) - (2n-1) + 2u \cdot (D_1 \cup \ldots \cup D_n \cup D_\infty) \\
&\geq 2n\ell+1-2n\ell+2nk\\
&= 1+2nk,
\end{align*}
thus
$$ 0= u \cdot D_\infty =nk\leq 0, $$
and since $n\geq1$ we deduce that $k=0$. Therefore $i=\ell$, i.e. $j=i+1$. Additionally, $u\cdot D_\infty =k\cdot n=0$, and from the positivity of intersections, Remark~\ref{rem: positivity of intersections for our acc. data} we deduce that $\im u \cap D_\infty=\varnothing$.

        \item Assume that $1 \leq j \leq \ell$. Let us first consider the trivialization $\tau_B$.

With respect to this trivialization, we find that $\mu_{FMB}^\tau(\check{x}_j^\ell;G_\ell) = 2nj$ and $\mu_{FMB}^\tau(\hat{x}_i^\ell;G_\ell) = 2ni+2n-1$. Thus Theorem~\ref{thm: c1 vs intersection number} implies that
$$ \sum_{k=1}^m c_1^\tau(u_k) + \sum_{\alpha\in A} c_1(B_\alpha) = (n+1) u \cdot D_\infty. $$
It follows that
\begin{align*}
1 &= \mu^\tau_{FMB}(q_+;G_\ell) - \mu^\tau_{FMB}(q_-;G_\ell) + 2\sum_{k=1}^m c_1^\tau(u_k) + 2\sum_{\alpha\in A} c_1(B_\alpha) \\
&= 2nj - 2ni - 2n + 1 + 2(n+1)u \cdot D_\infty \\
&= 2n(j-i-1) + 1 + 2(n+1)u \cdot D_\infty,
\end{align*}
yielding the relation
$$ u \cdot D_\infty = \frac{n(i+1-j)}{n+1}. $$
By the positivity of intersections (Remark~\ref{rem: positivity of intersections for our acc. data}), we have $u \cdot D_\infty \geq 0$, which implies that $i+1-j$ must be a multiple of $n+1$. Writing $i+1-j = k(n+1)$ for some $k \in \Z_{\geq 0}$, it follows that $u \cdot D_\infty = nk$.

Next we move to the trivialization $\tau_T$. 

With respect to this trivialization, we find that $\mu_{FMB}^\tau(\check{x}_j^\ell;G_\ell) = 0$ and $\mu_{FMB}^\tau(\hat{x}_i^\ell;G_\ell) = 2n-1$. Thus Theorem~\ref{thm: c1 vs intersection number} and Remark~\ref{rem: positivity of intersections for our acc. data} imply that
$$ \sum_{k=1}^m c_1^\tau(u_k) + \sum_{\alpha\in A} c_1(B_\alpha) = u \cdot (D_1 \cup \ldots \cup D_n \cup D_\infty) \geq u \cdot D_\infty. $$
It follows that
\begin{align*}
1 &= \mu^\tau_{FMB}(q_+;G_\ell) - \mu^\tau_{FMB}(q_-;G_\ell) + 2\sum_{k=1}^m c_1^\tau(u_k) + 2\sum_{\alpha\in A} c_1(B_\alpha) \\
&= 0 - (2n-1) + 2u \cdot (D_1 \cup \ldots \cup D_n \cup D_\infty) \\
&\geq -2n + 1 + 2u \cdot D_\infty,
\end{align*}
and thus
$$ nk = u \cdot D_\infty \leq n. $$
Since $k \in \Z_{\geq 0}$ and $nk \leq n$, we deduce that $k=0$ or $k=1$. If $k=1$ we get that $i+1-j=n+1$, namely, $j=i-n$. Otherwise, $k=0$ and thus $i+1-j=0$ which means that $j=i+1$. Additionally, in this case $u\cdot D_\infty =k\cdot n=0$, and from the positivity of intersections, Remark~\ref{rem: positivity of intersections for our acc. data} we deduce that $\im u \cap D_\infty=\varnothing$.

\item Assume that $-(n-1) \leq j \leq 0$. Let us first consider the trivialization $\tau_B$.

With respect to this trivialization, we find that $\mu_{FMB}^\tau(\check{x}_j^\ell;G_\ell) = -2j$ and $\mu_{FMB}^\tau(\hat{x}_i^\ell;G_\ell) = 2ni+2n-1$. Thus Theorem~\ref{thm: c1 vs intersection number} implies that
$$ \sum_{k=1}^m c_1^\tau(u_k) + \sum_{\alpha\in A} c_1(B_\alpha) = (n+1) u \cdot D_\infty. $$
It follows that
\begin{align*}
1 &= \mu^\tau_{FMB}(q_+;G_\ell) - \mu^\tau_{FMB}(q_-;G_\ell) + 2\sum_{k=1}^m c_1^\tau(u_k) + 2\sum_{\alpha\in A} c_1(B_\alpha) \\
&= -2j - 2ni - 2n + 1 + 2(n+1)u \cdot D_\infty ,
\end{align*}
yielding the relation
$$ u \cdot D_\infty = \frac{j+n(i+1)}{n+1}. $$

Next we move to the trivialization $\tau_T$. 

With respect to this trivialization, we find that $\mu_{FMB}^\tau(\check{x}_j^\ell;G_\ell) = -2j$ and $\mu_{FMB}^\tau(\hat{x}_i^\ell;G_\ell) = 2n-1$. Thus Theorem~\ref{thm: c1 vs intersection number} and Remark~\ref{rem: positivity of intersections for our acc. data} imply that
$$ \sum_{k=1}^m c_1^\tau(u_k) + \sum_{\alpha\in A} c_1(B_\alpha) = u \cdot (D_1 \cup \ldots \cup D_n \cup D_\infty) \geq u \cdot D_\infty. $$
It follows that
\begin{align*}
1 &= \mu^\tau_{FMB}(q_+;G_\ell) - \mu^\tau_{FMB}(q_-;G_\ell) + 2\sum_{k=1}^m c_1^\tau(u_k) + 2\sum_{\alpha\in A} c_1(B_\alpha) \\
&= -2j - (2n-1) + 2u \cdot (D_1 \cup \ldots \cup D_n \cup D_\infty) \\
&\geq -2j-2n + 1 + 2u \cdot D_\infty\\
&= -2j-2n + 1 +2\cdot\frac{j+n(i+1)}{n+1},
\end{align*}
and thus $i \leq j+n$.

Since $u \cdot D_\infty = (j+n(i+1))/(n+1)$ is an integer, we deduce that
$$ j+n(i+1) \equiv 0 \pmod{n+1}. $$
Given $n \equiv -1 \pmod{n+1}$, we get $j-i-1 \equiv 0 \pmod{n+1}$. Consequently, there exists $k \in \Z$ such that $j-i-1 = k(n+1)$. Based on the ranges $1 \leq i \leq n+j$ and $-(n-1)
\leq j \leq 0$, it follows that
$$ -(n+1) \leq j-i-1 \leq -2, $$
which proves that $j-i-1 =-( n+1)$, and hence $j=i-n$.
    \end{enumerate}
    
\end{itemize}
\end{proof}

\begin{thm}\label{thm: continuation_obst_for_complement}
    Let $0 < \Delta\leq\Delta' < 1$ and $\ell,\ell' \in \Z_{\geq 0}$ with $\ell\leq \ell'$. Define $G_\ell,G'_{\ell'} \fc \CP^n \to \R$ by $G_\ell(z) = g(\Delta, \ell, \mu(z))$ and $G'_{\ell'}(z) = g(\Delta', \ell', \mu(z))$ for every $z \in \CP^n$. Let $q_-$ and $q_+$ be critical points of the Morse functions associated with the critical submanifolds with respect to $G_\ell$ and $G'_{\ell'}$, respectively. 
    
    Suppose $u$ is a broken bubbled continuation flowline with cascades connecting $q_-$ to $q_+$, with cascades $u_1, \dots, u_m$ and a finite collection $\{B_\alpha\}_{\alpha\in A}$ of bubbles. Let $\tau$ be a trivialization of $\det TM$ along the asymptotics of $u_1,\ldots,u_m$. If 
    $$ \mu^\tau_{FMB}(q_+;G'_{\ell'}) - \mu^\tau_{FMB}(q_-;G_\ell) + 2\sum_{k=1}^m c_1^\tau(u_k) + 2\sum_{\alpha\in A} c_1(B_\alpha) = 0, $$
    then the pair $(q_-, q_+)$ cannot be any of the following:
    \begin{itemize}
        \item $(\check{x}^\ell_i, \hat{y}^{\ell'}_j)$ for every $-(n-1)\leq i \leq \ell + 1$ and $1 \leq j \leq \ell'$.
        \item $(\hat{x}^\ell_i, \check{y}^{\ell'}_j)$ for every $1 \leq i \leq \ell$ and $-(n-1) \leq j \leq \ell' + 1$.
        \item $(\hat{x}^\ell_i, \hat{y}^{\ell'}_j)$ for every $1 \leq i \leq \ell$ and $1\leq j\leq \ell'$, except where $i=j$.
        \item $(\check{x}^\ell_i, \check{y}^{\ell'}_j)$  for every $-(n-1) \leq i \leq \ell+1$ and $-(n-1)\leq j\leq \ell'+1$, except where $i=j$.
    \end{itemize}
\end{thm}

\begin{proof}
    As explained before, we can assume without loss of generality that the divisors $D_1, \ldots, D_n$ do not intersect any of the non-constant asymptotes of $u_1, \ldots, u_m$, and do not contain the images of $u_1, \ldots, u_m$ and $\{B_\alpha\}_{\alpha \in A}$. Moreover, Proposition \ref{prop: contained components complement} guarantees that all cascades and bubbles contained in $D_\infty$ are constant and hence can be safely ignored in the proof.

We treat each case separately:
\begin{itemize}

    \item Let $-(n-1) \leq i \leq \ell + 1$ and $1 \leq j \leq \ell'$, and assume that $(q_-,q_+)=(\check{x}^\ell_i, \hat{y}^{\ell'}_j)$.
    
    As computed in Section~\ref{ss: computations of FMB for complement},  $\mu_{FMB}(\hat{y}^{\ell'}_j;G'_{\ell'})$ is odd and $\mu_{FMB}(\check{x}^\ell_i;G_\ell)$ is even modulo $2(n+1)$. In particular, the difference 
    $$\mu_{FMB}(\hat{y}^{\ell'}_j;G'_{\ell'}) - \mu_{FMB}(\check{x}^\ell_i;G_\ell) $$
    is odd, and hence not equal to $0$ modulo $2(n+1)$. Proposition~\ref{prop: traj and modular indices} then leads to a contradiction.

    \item Let $1 \leq i \leq \ell$ and $-(n-1) \leq j \leq \ell'+1$, and assume that $(q_-,q_+)=(\hat{x}^\ell_i, \check{y}^{\ell'}_j)$.
    
    As computed in Section~\ref{ss: computations of FMB for complement},  $\mu_{FMB}(\check{y}^{\ell'}_j;G'_{\ell'})$ is even and $\mu_{FMB}(\hat{x}^\ell_i;G_\ell)$ is odd modulo $2(n+1)$. In particular, the difference 
    $$\mu_{FMB}(\check{y}^{\ell'}_j;G'_{\ell'}) - \mu_{FMB}(\hat{x}^\ell_i;G_\ell) $$
    is odd, and hence not equal to $0$ modulo $2(n+1)$. Proposition~\ref{prop: traj and modular indices} then leads to a contradiction.

    \item Let $1\leq i\leq \ell$ and $1\leq j\leq \ell'$, and assume that $(q_-,q_+)=(\hat{x}^\ell_i,\hat{y}^{\ell'}_j)$.
    
    Let us use the trivialization $\tau_T$.

    With respect to this trivialization, we find that $\mu_{FMB}^\tau(\hat{y}_j^{\ell'};G'_{\ell'}) =\mu_{FMB}^\tau(\hat{x}_i^\ell;G_\ell) = 2n-1$. By the positivity of intersections from Remark~\ref{rem: positivity of intersections for our acc. data}, we deduce that 
$$ \sum_{k=1}^m c_1^\tau(u_k) + \sum_{\alpha\in A} c_1(B_\alpha) \geq 0. $$

It follows that
\begin{align*}
    0 &= \mu^\tau_{FMB}(q_+;G'_{\ell'}) - \mu^\tau_{FMB}(q_-;G_\ell) + 2\sum_{k=1}^m c_1^\tau(u_k) + 2\sum_{\alpha\in A} c_1(B_\alpha)\\
    &=2n-1 - (2n-1) + 2\sum_{k=1}^m c_1^\tau(u_k) + 2\sum_{\alpha\in A} c_1(B_\alpha)\\
    &= 2\sum_{k=1}^m c_1^\tau(u_k) + 2\sum_{\alpha\in A} c_1(B_\alpha)\\
    &\geq0,
\end{align*}
which yields, using Theorem~\ref{thm: c1 vs intersection number} and Remark~\ref{rem: positivity of intersections for our acc. data} that
   $$ 0\leq u\cdot D_\infty\leq u\cdot(D_1\cup\ldots\cup D_n\cup D_\infty)=\sum_{k=1}^m c_1^\tau(u_k) + \sum_{\alpha\in A} c_1(B_\alpha)\\
    =0,$$ 
    thus we deduce that $u\cdot D_\infty=0$.

     Next we move to the trivialization $\tau_B$. With respect to this trivialization, we find that $\mu_{FMB}^\tau(\hat{y}_j^{\ell'};G'_{\ell'}) = 2nj+2n-1$ and $\mu_{FMB}^\tau(\hat{x}_i^\ell;G_\ell) = 2ni+2n-1$. From Theorem~\ref{thm: c1 vs intersection number} we deduce that
$$ \sum_{k=1}^m c_1^\tau(u_k) + \sum_{\alpha\in A} c_1(B_\alpha)=(n+1) u\cdot D_\infty=0. $$
It follows that
\begin{align*}
     0 &= \mu^\tau_{FMB}(q_+;G'_{\ell'}) - \mu^\tau_{FMB}(q_-;G_\ell) + 2\sum_{k=1}^m c_1^\tau(u_k) + 2\sum_{\alpha\in A} c_1(B_\alpha) \\
     &= 2nj+2n-1-(2ni+2n-1) \\
     &=2n(j-i),
\end{align*}
and hence $i=j$.

    \item Let $-(n-1)\leq i\leq \ell+1$ and $-(n-1)\leq j\leq \ell'+1$, and assume that $(q_-,q_+)=(\check{x}^\ell_i,\check{y}^{\ell'}_j)$.

    We have nine cases to check:
    
    \begin{enumerate}
    \item Assume that $i=\ell+1$ and $1\leq j\leq \ell'$. Let us use the trivialization $\tau_T$. 

    With respect to this trivialization, we find that $\mu_{FMB}^\tau(\check{y}_j^{\ell'};G'_{\ell'}) = 0$ and $\mu_{FMB}^\tau(\check{x}_i^\ell;G_\ell) =2n(\ell+1)$. By Theorem~\ref{thm: c1 vs intersection number} and Example~\ref{exam: computations using Seidel's lemma for complement}, we deduce that 
$$ \sum_{k=1}^m c_1^\tau(u_k) + \sum_{\alpha\in A} c_1(B_\alpha) =u\cdot(D_1\cup\ldots\cup D_n\cup D_\infty)\geq u\cdot(D_1\cup\ldots\cup D_n)\geq n(\ell+1). $$

It follows that
\begin{align*}
    0 &= \mu^\tau_{FMB}(q_+;G'_{\ell'}) - \mu^\tau_{FMB}(q_-;G_\ell) + 2\sum_{k=1}^m c_1^\tau(u_k) + 2\sum_{\alpha\in A} c_1(B_\alpha)\\
    &=-2n(\ell+1)+ 2\sum_{k=1}^m c_1^\tau(u_k) + 2\sum_{\alpha\in A} c_1(B_\alpha)\\
    &\geq-2n(\ell+1)+2n(\ell+1)\\
    &=0,
\end{align*}
which yields, using Theorem~\ref{thm: c1 vs intersection number} and Remark~\ref{rem: positivity of intersections for our acc. data} that
   $$ u\cdot(D_1\cup\ldots\cup D_n)+u\cdot D_\infty=\sum_{k=1}^m c_1^\tau(u_k) + \sum_{\alpha\in A} c_1(B_\alpha)\\
    =n(\ell+1),$$ 
    thus since $u\cdot D_\infty\geq0$ and $u\cdot(D_1\cup \ldots\cup D_n)\geq n(\ell+1)$ we deduce that $u\cdot D_\infty=0$.

  Next we move to the trivialization $\tau_B$. With respect to this trivialization, we find that $\mu_{FMB}^\tau(\check{y}_j^{\ell'};G'_{\ell'}) = 2nj$ and $\mu_{FMB}^\tau(\check{x}_i^\ell;G_\ell) = 2n(\ell+1)$. From Theorem~\ref{thm: c1 vs intersection number} we deduce that
$$ \sum_{k=1}^m c_1^\tau(u_k) + \sum_{\alpha\in A} c_1(B_\alpha)=(n+1) u\cdot D_\infty=0. $$
It follows that
\begin{align*}
     0 &= \mu^\tau_{FMB}(q_+;G'_{\ell'}) - \mu^\tau_{FMB}(q_-;G_\ell) + 2\sum_{k=1}^m c_1^\tau(u_k) + 2\sum_{\alpha\in A} c_1(B_\alpha) \\
     &=2nj-2n(\ell+1),
\end{align*}
and hence $i=\ell+1=j$ as required.

 \item Assume that $i=\ell+1$ and $-(n-1)\leq j\leq 0$. Let us use the trivialization $\tau_T$. 

    With respect to this trivialization, we find that $\mu_{FMB}^\tau(\check{y}_j^{\ell'};G'_{\ell'}) =-2j$ and $\mu_{FMB}^\tau(\check{x}_i^\ell;G_\ell) = 2n(\ell+1)$. By Theorem~\ref{thm: c1 vs intersection number} and Example~\ref{exam: computations using Seidel's lemma for complement}, we deduce that 
$$ \sum_{k=1}^m c_1^\tau(u_k) + \sum_{\alpha\in A} c_1(B_\alpha) =u\cdot(D_1\cup\ldots\cup D_n\cup D_\infty)\geq u_1\cdot(D_1\cup\ldots\cup D_n)+u_m\cdot D_\infty\geq n(\ell+1)+1. $$

It follows that
\begin{align*}
    0 &= \mu^\tau_{FMB}(q_+;G'_{\ell'}) - \mu^\tau_{FMB}(q_-;G_\ell) + 2\sum_{k=1}^m c_1^\tau(u_k) + 2\sum_{\alpha\in A} c_1(B_\alpha)\\
    &= -2j-2n(\ell+1)+2\sum_{k=1}^m c_1^\tau(u_k) + 2\sum_{\alpha\in A} c_1(B_\alpha)\\
    &\geq-2j-2n(\ell+1)+2(n(\ell+1)+1)\\
     &\geq-2j+2.
\end{align*}
Thus we deduce that $j\geq 1$ which contradicts that assumption $0\geq j$.

        \item Assume that $1\leq i\leq \ell$ and $j=\ell'+1$. Let us use the trivialization $\tau_T$. 

    With respect to this trivialization, we find that $\mu_{FMB}^\tau(\check{y}_j^{\ell'};G'_{\ell'})=2n(\ell'+1)$ and $\mu_{FMB}^\tau(\check{x}_i^\ell;G_\ell) = 0$. By Theorem~\ref{thm: c1 vs intersection number} and Example~\ref{exam: computations using Seidel's lemma for complement}, we deduce that 
$$ \sum_{k=1}^m c_1^\tau(u_k) + \sum_{\alpha\in A} c_1(B_\alpha) =u\cdot(D_1\cup\ldots\cup D_n\cup D_\infty)\geq u_m\cdot(D_1\cup\ldots\cup D_n)\geq n\cdot(-\ell')=-n\ell'. $$

It follows that
\begin{align*}
    0 &= \mu^\tau_{FMB}(q_+;G'_{\ell'}) - \mu^\tau_{FMB}(q_-;G_\ell) + 2\sum_{k=1}^m c_1^\tau(u_k) + 2\sum_{\alpha\in A} c_1(B_\alpha)\\
    &=2n(\ell'+1)+ 2\sum_{k=1}^m c_1^\tau(u_k) + 2\sum_{\alpha\in A} c_1(B_\alpha)\\
    &\geq 2n(\ell'+1)-2n\ell'\\
    &=2n,
\end{align*}
and this is a contradiction.

    \item Assume that $i=\ell+1$ and $j=\ell'+1$.
        If $\ell=\ell'$, then $i=j$. Assume therefore that $\ell<\ell'$.
        Let us use the trivialization $\tau_B$.
        
        With respect to this trivialization, we find that
        $\mu_{FMB}^{\tau_B}(\check{x}_i^\ell;G_\ell)
        =2n(\ell+1)$ and
        $\mu_{FMB}^{\tau_B}(\check{y}_j^{\ell'};G'_{\ell'})
        =2n(\ell'+1)$.
        Since both asymptotes lie outside $D_\infty$,
        Theorem~\ref{thm: c1 vs intersection number} and
        the positivity of intersections from
        Remark~\ref{rem: positivity of intersections for our acc. data}
        imply that
        $$
        \sum_{k=1}^m c_1^{\tau_B}(u_k)
        +\sum_{\alpha\in A}c_1(B_\alpha)
        =(n+1)u\cdot D_\infty\geq0.
        $$
        It follows that
        \begin{align*}
        0
        &=\mu_{FMB}^{\tau_B}(q_+;G'_{\ell'})
          -\mu_{FMB}^{\tau_B}(q_-;G_\ell)
          +2\sum_{k=1}^m c_1^{\tau_B}(u_k)
          +2\sum_{\alpha\in A}c_1(B_\alpha)\\
        &=2n(\ell'-\ell)
          +2\sum_{k=1}^m c_1^{\tau_B}(u_k)
          +2\sum_{\alpha\in A}c_1(B_\alpha)\\
        &\geq2n(\ell'-\ell)>0,
        \end{align*}
        which is a contradiction.

  \item Assume that $1\leq i\leq \ell$ and $1\leq j\leq \ell'$. Let us use the trivialization $\tau_T$. 

    With respect to this trivialization, we find that $\mu_{FMB}^\tau(\check{y}_j^{\ell'};G'_{\ell'}) =\mu_{FMB}^\tau(\check{x}_i^\ell;G_\ell) = 0$. By the positivity of intersections from Remark~\ref{rem: positivity of intersections for our acc. data}, we deduce that 
$$ \sum_{k=1}^m c_1^\tau(u_k) + \sum_{\alpha\in A} c_1(B_\alpha) \geq 0. $$

It follows that
\begin{align*}
    0 &= \mu^\tau_{FMB}(q_+;G'_{\ell'}) - \mu^\tau_{FMB}(q_-;G_\ell) + 2\sum_{k=1}^m c_1^\tau(u_k) + 2\sum_{\alpha\in A} c_1(B_\alpha)\\
    &= 2\sum_{k=1}^m c_1^\tau(u_k) + 2\sum_{\alpha\in A} c_1(B_\alpha)\\
    &\geq0,
\end{align*}
which yields, using Theorem~\ref{thm: c1 vs intersection number} and Remark~\ref{rem: positivity of intersections for our acc. data} that
   $$ 0\leq u\cdot D_\infty\leq u\cdot(D_1\cup\ldots\cup D_n\cup D_\infty)=\sum_{k=1}^m c_1^\tau(u_k) + \sum_{\alpha\in A} c_1(B_\alpha)\\
    =0,$$ 
    thus we deduce that $u\cdot D_\infty=0$.

    Next we move to the trivialization $\tau_B$. 

With respect to this trivialization, we find that $\mu_{FMB}^\tau(\check{y}_j^{\ell'};G'_{\ell'}) = 2nj$ and $\mu_{FMB}^\tau(\check{x}_i^\ell;G_\ell) = 2ni$. From Theorem~\ref{thm: c1 vs intersection number} we deduce that
$$ \sum_{k=1}^m c_1^\tau(u_k) + \sum_{\alpha\in A} c_1(B_\alpha)=(n+1) u\cdot D_\infty=0. $$
It follows that
\begin{align*}
     0 &= \mu^\tau_{FMB}(q_+;G'_{\ell'}) - \mu^\tau_{FMB}(q_-;G_\ell) + 2\sum_{k=1}^m c_1^\tau(u_k) + 2\sum_{\alpha\in A} c_1(B_\alpha) \\
     &=2n(j-i),
\end{align*}
and hence $i=j$.

  \item Assume that $1\leq i\leq \ell$ and $-(n-1)\leq j\leq 0$. Let us use the trivialization $\tau_T$. 

    With respect to this trivialization, we find that $\mu_{FMB}^\tau(\check{y}_j^{\ell'};G'_{\ell'})=-2j$ and $\mu_{FMB}^\tau(\check{x}_i^\ell;G_\ell) = 0$. By Theorem~\ref{thm: c1 vs intersection number} and Example~\ref{exam: computations using Seidel's lemma}, we deduce that  
$$ \sum_{k=1}^m c_1^\tau(u_k) + \sum_{\alpha\in A} c_1(B_\alpha) \geq u_m\cdot D_\infty\geq 1. $$

It follows that
\begin{align*}
    0 &= \mu^\tau_{FMB}(q_+;G'_{\ell'}) - \mu^\tau_{FMB}(q_-;G_\ell) + 2\sum_{k=1}^m c_1^\tau(u_k) + 2\sum_{\alpha\in A} c_1(B_\alpha)\\
    &= -2j+2\sum_{k=1}^m c_1^\tau(u_k) + 2\sum_{\alpha\in A} c_1(B_\alpha)\\
    &\geq-2j+2\\
    &\geq2,
\end{align*}
which yields a contradiction.

    \item Assume that $-(n-1)\leq i\leq 0$ and $j=\ell'+1$. Let us use the trivialization $\tau_T$. 

    With respect to this trivialization, we find that $\mu_{FMB}^\tau(\check{y}_j^{\ell'};G'_{\ell'})=2n(\ell'+1)$ and $\mu_{FMB}^\tau(\check{x}_i^\ell;G_\ell) = -2i$. By Theorem~\ref{thm: c1 vs intersection number} and Example~\ref{exam: computations using Seidel's lemma for complement}, we deduce that 
    \begin{align*}
         \sum_{k=1}^m c_1^\tau(u_k) + \sum_{\alpha\in A} c_1(B_\alpha) &=u\cdot(D_1\cup\ldots\cup D_n\cup D_\infty)\\
         &\geq u_1\cdot D_\infty+u_m\cdot(D_1\cup\ldots\cup D_n)\\
         &\geq 0+n(-\ell')\\
         &=-n\ell'.
    \end{align*}

It follows that
\begin{align*}
    0 &= \mu^\tau_{FMB}(q_+;G'_{\ell'}) - \mu^\tau_{FMB}(q_-;G_\ell) + 2\sum_{k=1}^m c_1^\tau(u_k) + 2\sum_{\alpha\in A} c_1(B_\alpha)\\
    &=2n(\ell'+1)+2i+ 2\sum_{k=1}^m c_1^\tau(u_k) + 2\sum_{\alpha\in A} c_1(B_\alpha)\\
    &\geq 2n(\ell'+1)+2i-2n\ell'\\
    &\geq 2n+2i.
\end{align*}
thus we deduce that
$$-(n-1)\leq i\leq -n$$
and this is a contradiction.

  \item Assume that $-(n-1)\leq i\leq 0$ and $1\leq j\leq \ell'$. Let us use the trivialization $\tau_B$. 

    With respect to this trivialization, we find that $\mu_{FMB}^\tau(\check{y}_j^{\ell'};G'_{\ell'})=2nj$ and $\mu_{FMB}^\tau(\check{x}_i^\ell;G_\ell) = -2i$. By Theorem~\ref{thm: c1 vs intersection number} we know that
$$ \sum_{k=1}^m c_1^\tau(u_k) + \sum_{\alpha\in A} c_1(B_\alpha) = (n+1)u\cdot D_\infty. $$

It follows that
\begin{align*}
    0 &= \mu^\tau_{FMB}(q_+;G'_{\ell'}) - \mu^\tau_{FMB}(q_-;G_\ell) + 2\sum_{k=1}^m c_1^\tau(u_k) + 2\sum_{\alpha\in A} c_1(B_\alpha)\\
    &=2nj+2i+2(n+1)u\cdot D_\infty
\end{align*}
which yields that
$$ u\cdot D_\infty=-\frac{nj+i}{n+1},$$
is an integer.
Since $n\equiv -1 \pmod{n+1}$ and $nj+i\equiv0 \pmod{n+1}$ we deduce that 
$$-j+i\equiv0 \pmod{n+1},$$
thus there exists $r\in \Z$ such that $i=j+(n+1)r$. Thus we can write $u\cdot D_\infty=-(j+r)$.

Additionally, from Example~\ref{exam: computations using Seidel's lemma for complement} and Remark~\ref{rem: positivity of intersections for our acc. data} we deduce that
$$u\cdot D_\infty\geq u_1* D_\infty\geq 0,$$
but
\[ u\cdot D_\infty=-\frac{nj+i}{n+1}\le -\frac{n\cdot 1-(n-1)}{n+1}=-\frac{1}{n+1}<0,\]
a contradiciton.

\item Assume that $-(n-1)\leq i\leq 0$ and $-(n-1)\leq j\leq 0$. Let us use the trivialization $\tau_B$. 

    With respect to this trivialization, we find that $\mu_{FMB}^\tau(\check{y}_j^{\ell'};G'_{\ell'})=-2j$ and $\mu_{FMB}^\tau(\check{x}_i^\ell;G_\ell) = -2i$. By Theorem~\ref{thm: c1 vs intersection number} we know that
$$ \sum_{k=1}^m c_1^\tau(u_k) + \sum_{\alpha\in A} c_1(B_\alpha) = (n+1)u\cdot D_\infty. $$

It follows that
\begin{align*}
    0 &= \mu^\tau_{FMB}(q_+;G'_{\ell'}) - \mu^\tau_{FMB}(q_-;G_\ell) + 2\sum_{k=1}^m c_1^\tau(u_k) + 2\sum_{\alpha\in A} c_1(B_\alpha)\\
    &=2(i-j)+2(n+1)u\cdot D_\infty
\end{align*}
which yields that
$$ u\cdot D_\infty=\frac{j-i}{n+1}$$
is an integer and in particular $j-i\in (n+1)\Z$. But $|j-i|\le n-1 < n+1$ so we deduce $j-i=0$, so $i=j$.
    \end{enumerate}
\end{itemize}
\end{proof}

\subsection{Computations of Floer complexes}

    The purpose of this section is to compute the Floer complexes, with coefficients in $\Z$, of the Hamiltonians from the acceleration data presented in Section~\ref{ss: acc. data CP^n - ball}, and the continuation maps between them.

    To achieve transversality for the Floer differentials, we perturb
the almost complex structures while keeping the Hamiltonians and
the auxiliary Morse data fixed. For continuation maps, we perturb
the Hamiltonian homotopies and the associated families of almost
complex structures, keeping the endpoint Hamiltonians fixed and
the homotopies constant outside $[0,1]$. The Hamiltonian homotopies
are chosen to remain monotone. Whenever the no-escape lemma is
used, the perturbations are chosen to leave the Hamiltonians
unchanged and the almost complex structures equal to $J_{std}$
on the corresponding fixed shell outside the balls.

The obstruction statements of
Theorems~\ref{thm: diff_obst_for_complement}
and~\ref{thm: continuation_obst_for_complement}
persist under sufficiently small such perturbations. Analogous statements for the ball are explained in section \ref{ss: transversality}. \textbf{Note:} Throughout this subsection, references to
Theorems~\ref{thm: diff_obst_for_complement}
and~\ref{thm: continuation_obst_for_complement}
also refer to their extensions to the chosen regular
perturbations, as justified above.

Let $\Delta \in (0,1)$. For every $\ell \in \Z_{\geq 0}$, we will compute the Floer cochain complex of the Hamiltonian $G_\ell \fc \CP^n \to \R$ which is given by $G_\ell(z) = g(\Delta, \ell, \mu(z))$ for every $z \in \CP^n$.

Fix $\ell \in \Z_{\geq 0}$. For every $1 \leq i \leq \ell$, let $r_{\ell,i}^\Delta \in [0,1)$ denote the number satisfying $\frac{\partial}{\partial r}|_{r=r_{\ell,i}^\Delta} g(\Delta, \ell, r) = -i$. For every $0 \leq i \leq \ell$, let $g_{\ell,i}^\Delta = g(\Delta, \ell, r_{\ell,i}^\Delta)$, and let $g_{\ell,0}^\Delta = g(\Delta, \ell, 1)$. Additionally, for every $1 \leq i \leq \ell$, define the $(2n-1)$-spheres
$$S_{\ell,i} = \{z \in \CP^n : g'_\ell(\mu(z)) = -i\},$$
each equipped with a perfect Morse function. We also equip the divisor $D_\infty = \CP^{n-1}$ with a perfect Morse function. 

As before, let $\check{x}_{\ell+1}^\ell$ denote the constant orbit of $G_\ell$ at $0 \in \text{Int } B(1) \subset \CP^n$. For every $1 \leq i \leq \ell$, let $\check{x}^\ell_i$ and $\hat{x}^\ell_i$ denote the minimum and maximum points in $S_{\ell,i}$, respectively. Finally, for every $-(n-1) \leq i \leq 0$, let $\check{x}_i^\ell$ denote the critical point on $D_\infty$ with Morse index $-2i$.

\begin{thm}\label{thm: CF(G_l;Z)}
    For every $\ell\in \Z_{\geq0}$ there exist signs $A_{\ell,1},\ldots,A_{\ell,\ell},B_{\ell,1},\ldots,B_{\ell,\ell}\in \{-1,1\}$,

    such that for every $\ell\in \Z_{\geq0}$, the differential $d\fc CF(G_\ell;\Z)\to CF(G_\ell;\Z)$ satisfies
    \begin{itemize}
        \item  For every $1\leq i\leq \ell$ we have $d \hat{x}^\ell_i=A_{\ell,i}\check{x}^\ell_{i+1}+B_{\ell,i}\check{x}^\ell_{i-n}$.
     
    \item For every $-(n-1)\leq i\leq \ell+1$  we have $d\check{x}^\ell_i=0$.
    \end{itemize}
\end{thm}

Similar techniques allow us to compute the continuation maps between the Floer cochain complexes of the acceleration data of two different balls:

Let $0<\Delta\leq \Delta'<1$. Consider additionally the corresponding acceleration data $(G'_\ell)_{\ell\geq0}$ for the complement of the ball $B^{}_{\Delta'}=\mu^{-1}([0,\Delta'])$, given by $G'_\ell(z)=g(\Delta',\ell,\mu(z))$, for every $z\in \CP^n$ and $\ell\in \Z_{\geq0}$. For every $\ell\in \Z_{\geq0}$ denote by
$\check{y}_{-(n-1)}^\ell,\ldots,\check{y}_{\ell+1}^\ell,\hat{y}_1^\ell,\ldots,\hat{y}_\ell^\ell$ the generators for $CF(G'_\ell)$, as in Theorem~\ref{thm: CF(G_l;Z)}.

Let $\ell,\ell'\in \Z_{\geq0}$ and assume that $\ell\leq \ell'$. Let $\chi\fc \R\to \R$ be a smooth non-decreasing function, satisfying $\chi(s)=\Delta$ for every $s<0$ and $\chi(s)=\Delta'$ for every $s>1$. Similarly, let $\lambda\fc \R\to \R$ be a smooth non-decreasing function, satisfying $\lambda(s)=\ell$ for every $s<0$ and $\lambda(s)=\ell'$ for every $s>1$. 

Now, define a monotone homotopy $K\fc \R\times \CP^n\to \R$ from $G_\ell$ to $G'_{\ell'}$ by $K(s,z)=g(\chi(s),\lambda(s),\mu(z))$ for every $(s,z)\in \R\times \CP^n$. The homotopy $K$ defines a continuation map $\Phi\fc CF(G_\ell;\Z)\to CF(G'_{\ell'};\Z)$.

\begin{thm}\label{thm: continuation maps for G_l, over Z}
    For every $0<\Delta\leq\Delta'<1$ and $\ell,\ell'\in \Z_{\geq0}$ with $\ell\leq \ell'$ there exist signs
    $$\check{C}_{\ell,-(n-1)},\ldots,\check{C}_{\ell,\ell+1},\hat{C}_{\ell,1},\ldots,\hat{C}_{\ell,\ell}\in \{-1,1\}$$
    such that the continuation map $\Phi\fc CF(G_\ell;\Z)\to CF(G'_{\ell'};\Z)$ satisfies
    \begin{itemize}
        \item For every $-(n-1)\leq i\leq \ell+1$ we have $\Phi(\check{x}_i^\ell)=\check{C}_{i} \check{y}_i^{\ell'}$.

         \item For every $1\leq i\leq \ell$ we have $\Phi(\hat{x}_i^\ell)=\hat{C}_{i} \hat{y}_i^{\ell'}$.
    \end{itemize}
\end{thm}

\subsubsection{Computation of Floer complexes of admissible Hamiltonians}

Let $\ell \in \Z_{\geq 0}$, and define $\tilde{G}_\ell \fc \Int B(1) \to \R$ by 
$$ \tilde{G}_\ell(z) = G_\ell \circ \iota(z) = g(\Delta, \ell, \pi\|z\|^2), $$
for every $z \in \Int B(1)$, where $G_\ell$ is the $\ell$-th function in the acceleration data defined in Section~\ref{ss: acc. data CP^n - ball}. Moreover, the Hamiltonian $\tilde{G}_\ell$ is autonomous and satisfies the \textbf{MB} condition. Its critical submanifolds are $\{0\}$ and the spheres 
$$ S_j = \{z \in \Int B(1) : g_\ell'(\pi\|z\|^2) = -j\}, $$
for every $1 \leq j \leq \ell$, where $g_\ell \fc [0,1] \to \R$ is given by $g_\ell(r) = g(\Delta, \ell, r)$ for every $r \in [0,1]$.

We denote the generators of $CF(\tilde{G}_\ell; \Z)$ by $ \check{x}_1^\ell,\hat{x}_1^\ell \dots,  \check{x}_\ell^\ell,\hat{x}_\ell^\ell,,\check{x}_{\ell+1}^\ell$, where for every $1 \leq j \leq \ell$ the generator $\check{x}_j^\ell$ corresponds to the minimum point of a Morse function on $S_j$, and $\hat{x}_j^\ell$ to the maximum point of this function.

The main result of this section is the following proposition:

\begin{prop}\label{prop: CF of J-shaped in C^n, for CPn-Ball}
    There exist signs $A_{\ell,1}, \dots, A_{\ell,\ell} \in \{-1,1\}$ such that the differential $d \fc CF(\tilde{G}_\ell; \Z) \to CF(\tilde{G}_\ell; \Z)$ satisfies:
    \begin{itemize}
        \item $d\hat{x}_i^\ell = A_{\ell,i} \check{x}_{i+1}^\ell$ for every $1 \leq i \leq \ell$.
        \item $d\check{x}_i^\ell = 0$ for every $1 \leq i \leq \ell+1$.
    \end{itemize}
\end{prop}

\begin{proof}
    We begin by computing the Floer complex of an auxiliary Hamiltonian. Let $\ell \in \Z_{\geq 0}$. Since $\frac{\partial}{\partial r} g(\Delta,\ell,r)\in (-1,0)$, and $\frac{\partial^2}{\partial r^2} g(\Delta,\ell,r)=0$ for every $r\in(\Delta,1]$ we deduce that that $\frac{\partial}{\partial r} g(\Delta,\ell,\cdot)$ is constant on $(\Delta,1]$ and equals to $-\delta$ for some $\delta\in(0,1)$.

    and define $F_\ell \fc \Int B(1) \to \R$ by 
    $$ F_\ell(z) = -\delta R(z) =- \pi\delta\|z\|^2 $$
    for every $z \in \Int B(1)$. Then $F_\ell$ has only one $1$-periodic orbit $\gamma$, located at the origin. Since $\tilde{G}_\ell$ is $C^2$-small near this orbit, and this orbit is a maximum, we deduce that the Robbin-Salamon index of $\gamma$ with respect to the constant trivialization is $n$; hence,
    $$ \mu_{FMB}^\tau(\gamma;F_\ell) = \mu_{RS}^\tau(\gamma;F_\ell) + \frac{1}{2}\dim \R^{2n} - \frac{1}{2}\dim \{0\} + \ind_h \gamma = n + n = 2n, $$
    where $h$ is a Morse function on $\{0\}$, which must be constant. Thus, the Floer homology of $F_\ell$ is
    $$ HF^*(F_\ell;\Z) = \begin{cases}
        \Z \cdot [\gamma], & * = 2n, \\
        0, & * \neq 2n.
    \end{cases} $$

    Returning to the Hamiltonian $\tilde{G}_\ell$, let $\tau$ denote the trivialization of the $1$-periodic orbits of $\tilde{G}_\ell$ obtained from cappings. From Section~\ref{ss: computations of FMB for complement}, we know that $\mu_{FMB}^\tau(\check{x}_{\ell+1}^\ell)=2n(\ell+1)$, and for every $1 \leq i \leq \ell$:
    $$ \mu_{FMB}^\tau(\check{x}_i^\ell) = 2ni \quad \text{and} \quad \mu_{FMB}^\tau(\hat{x}_i^\ell) = 2ni + 2n - 1 = 2n(i+1) - 1. $$

    In the case $n=1$, since the differential increases the FMB-index by $+1$, we deduce that for every $1 \leq i \leq \ell$ there exists $C_i \in \Z$ such that $d\check{x}_i^\ell = C_i \hat{x}_i^\ell$. Fix $1\le i \le \ell$.
    Note that $\C^n$ is an exact symplectic manifold, hence action is well defined on all asymptotic orbits of the cascades. Since Floer solutions are action increasing, and since the action of  $\hat{x}^\ell_i$ equals the action of $\check{x}^\ell_i$, all cascades are of zero energy, hence constant. Therefore, there are no cascades in the flowline. Thus $C_i=0$ by the usual count of Morse differential in the Morse homology of $S^1$.

    Consequently, for any $n \in \N$, since the differential increases the FMB-index by $+1$, there exist $A_{\ell,1}, \dots, A_{\ell,\ell} \in \Z$ such that the differential $d \fc CF(\tilde{G}_\ell;\Z) \to CF(\tilde{G}_\ell;\Z)$ satisfies:
    \begin{itemize}
        \item $d\hat{x}_i^\ell = A_{\ell,i} \check{x}_{i+1}^\ell$ for every $1 \leq i \leq \ell$.
        \item $d\check{x}_i^\ell = 0$ for every $1 \leq i \leq \ell+1$.
    \end{itemize}
    To complete the proof, it remains to show that $A_{\ell,i} \in \{-1, 1\}$ for every $1 \leq i \leq \ell$.

    Note that $F_\ell$ and $\tilde{G}_\ell$ are admissible Hamiltonians with $F_\ell \leq \tilde{G}_\ell$, and they share the same slope at infinity. Thus, by Claim~\ref{claim: qi for Ham with the same slope}, we have an isomorphism $HF^*(\tilde{G}_\ell;\Z) \cong HF^*(F_\ell;\Z)$.

    Let $1 \leq i \leq \ell$. Since $HF^{2n(i+1)}(\tilde{G}_\ell;\Z) \cong HF^{2n(i+1)}(F_\ell;\Z) = 0$, we deduce that the restriction of the differential
    $$ d|_{CF^{2n(i+1)-1}(\tilde{G}_\ell;\Z)} \fc CF^{2n(i+1)-1}(\tilde{G}_\ell;\Z) \to CF^{2n(i+1)}(\tilde{G}_\ell;\Z) $$
    is surjective. Since $CF^{2n(i+1)-1}(\tilde{G}_\ell;\Z) = \Z \cdot \hat{x}_{i}^\ell$ and $CF^{2n(i+1)}(\tilde{G}_\ell;\Z) = \Z \cdot \check{x}_{i+1}^\ell$, it follows that $A_{\ell,i} \in \{-1, 1\}$ as required.
\end{proof}

    \subsection{Proof of Theorem~\ref{thm: CF(G_l;Z)} and Theorem~\ref{thm: continuation maps for G_l, over Z}}\label{ss: proof of computations of CF + continuations, for G_l}

Let us start with the following auxiliary lemma:
\begin{lemma}\label{lemma: generators of HF(G_l)}
    For every $\ell\in \Z_{\geq0}$ and each $j$ satisfying $-(n-1)\leq j\leq \ell+1$, we have
    $$HF^*(G_\ell;\Z)=\Z\cdot[\check{x}_j^\ell],$$
    where $*=\mu^\tau_{FMB}(\check{x}_{j}^\ell)\pmod{2(n+1)}$ and $\tau$ is the trivialization that is induced by the constant cappings of $\check{x}_j^\ell$ for $j\in\{\ell+1,0,-1,\ldots,-(n-1)\}$ and a capping of $\check{x}_j^\ell$ that does not intersect $D_\infty$ otherwise.
\end{lemma}

\begin{proof}
    Let $\ell\in \Z_{\geq0}$ and $-(n-1)\leq j\leq \ell+1$. We have two cases to check:
    \begin{itemize}
        \item Assume that $j\leq 0$. Note that $\mu_{FMB}^\tau(\check{x}_j^\ell)=-2j$. The Hamiltonian $G_0$ has $(n+1)$ $1$-periodic orbits, all of which are constant and have an even Floer--Morse--Bott index with respect to the constant cappings; therefore, the differential of $CF(G_0;\Z)$ is zero. Also, we know that $\check{x}_j^0$ is the only $1$-periodic orbit of $G_0$ with Floer--Morse--Bott index $-2j$, hence $HF^{-2j}(G_0;\Z)=\Z\cdot[\check{x}_j^0]$.

        Let $F \fc CF(G_0;\Z) \to CF(G_\ell;\Z)$ be the continuation map. By Theorem~\ref{thm: continuation_obst_for_complement}, there exists an integer $k\in \Z$ such that $F(\check{x}^0_j) = k \check{x}^{\ell}_{j}$. 

        Since 
        $$F_*([\check{x}_j^0]) = [k\check{x}_j^\ell] = k[\check{x}_j^\ell]$$
        and $[\check{x}_j^0]$ is a generator for $HF^{-2j}(G_0;\Z)$, we deduce that 
        $$F_*({HF^{-2j}(G_0;\Z)})=\Z\cdot k[\check{x}_j^\ell].$$ 
        On the other hand, since $\CP^n$ is a closed symplectic manifold, we know that the continuation map $F$ is a quasi-isomorphism; see \cite[Theorem 4]{Floer_1989_monotone_mfd}. Therefore, the induced map 
        $$F_*|_{HF^{-2j}(G_0;\Z)} \fc HF^{-2j}(G_0;\Z) \to HF^{-2j}(G_\ell;\Z)$$
        is an isomorphism, and hence
        $$HF^{-2j}(G_\ell;\Z) = \im F_*|_{HF^{-2j}(G_0;\Z)} = \Z\cdot k[\check{x}_j^\ell].$$
        Thus, $[\check{x}_j^\ell] \in HF^{-2j}(G_\ell;\Z) = \Z \cdot k[\check{x}_j^\ell]$, which implies that $k \in \{-1, 1\}$, and hence $HF^{-2j}(G_\ell;\Z) = \Z \cdot [\check{x}_j^\ell]$.

\item Assume that $j\geq 1$. Note that $\mu_{FMB}^\tau(\check{x}_j^\ell)=2nj$. Consider the Hamiltonian $G\fc \CP^n\to \R$ given by 
        $$G(z)=-(j-1/2)\mu(z),$$
        for every $z\in \CP^n$.
        The Hamiltonian $G$ has $(n+1)$ $1$-periodic orbits, all of which are constant and have an even Floer--Morse--Bott index with respect to the constant cappings; therefore, the differential of $CF(G;\Z)$ is zero. Denote by $x_j$ the $1$-periodic orbit of $G$ which is located at the origin of $\Int B(1)\subset \CP^n$. We know that $x_j$ is the only $1$-periodic orbit of $G$ with Floer--Morse--Bott index $2nj$, hence $HF^{2nj}(G;\Z)=\Z\cdot[x_j]$.

        Let $F \fc CF(G;\Z) \to CF(G_\ell;\Z)$ be the continuation map. Since $G$ coincides with $G_{j-1}$ in a neighborhood of $x_j$, up to the addition of a constant, by Theorem~\ref{thm: continuation_obst_for_complement}, there exists an integer $k\in \Z$ such that $F(x_j) = k \check{x}^{\ell}_{j}$. 

        Since 
        $$F_*([x_j]) = [k\check{x}_j^\ell] = k[\check{x}_j^\ell]$$
        and $[x_j]$ is a generator for $HF^{2nj}(G;\Z)$, we deduce that 
        $$F_*(HF^{2nj}(G;\Z))=\Z\cdot k[\check{x}_j^\ell].$$ 
        On the other hand, since $\CP^n$ is a closed symplectic manifold, we know that the continuation map $F$ is a quasi-isomorphism; see \cite[Theorem 4]{Floer_1989_monotone_mfd}. Therefore, the induced map 
        $$F_*|_{HF^{2nj}(G;\Z)} \fc HF^{2nj}(G;\Z) \to HF^{2nj}(G_\ell;\Z)$$
        is an isomorphism, and hence
        $$HF^{2nj}(G_\ell;\Z) = \im F_*|_{HF^{2nj}(G;\Z)} = \Z\cdot k[\check{x}_j^\ell].$$
        Thus, $[\check{x}_j^\ell] \in HF^{2nj}(G_\ell;\Z) = \Z \cdot k[\check{x}_j^\ell]$, which implies that $k \in \{-1, 1\}$, and hence $HF^{2nj}(G_\ell;\Z) = \Z \cdot [\check{x}_j^\ell]$.
    \end{itemize}
\end{proof}

  \begin{rem}\label{rem: conc from thm: diff_obst for G_l}
    Let $n\in \N$ and $\ell\in \Z_{\geq0}$. By Theorem~\ref{thm: diff_obst_for_complement}, there are integers
$$A_{\ell,1}, \dots, A_{\ell,\ell}, B_{\ell,1}, \dots, B_{\ell,\ell} \in \Z,$$
such that the differential $d \fc CF(G_\ell; \Z) \to CF(G_\ell; \Z)$ satisfies:
\begin{itemize}
    \item For every $1 \leq i \leq \ell$, we have $d \hat{x}^\ell_i = A_{\ell,i} \check{x}^\ell_{i+1} + B_{\ell,i} \check{x}^\ell_{i-n}$.
    \item For every $-(n-1) \leq i \leq \ell+1$, we have $d \check{x}^\ell_i = 0$.
\end{itemize}

Indeed this follows from Theorem~\ref{thm: diff_obst_for_complement}, 
except for the case where $n=1$. In that case, we have $d\check{x}_i^\ell=C_i\hat{x}_i^\ell$ for every $1\leq i\leq \ell$, where $C_1,\ldots,C_\ell\in \Z$. Fix $1 \le i \le \ell$.
By Theorem~\ref{thm: diff_obst_for_complement} every flowline with cascades connecting $\check{x}^\ell_i$ to $\hat{x}^\ell_i$ does not intersect $D_\infty$, and by the no escape lemma (Lemma \ref{lemma: no escape}), it is thus contained in a ball, and hence the action functional is well defined on all asymptotic orbits of the cascades. Since Floer solutions are action increasing, and since the action of  $\hat{x}^\ell_i$ equals the action of $\check{x}^\ell_i$, all cascades are of zero energy, hence constant. Therefore, there are no cascades in the flowline. Thus $C_i=0$ by the usual count of Morse differential in the Morse homology of $S^1$, and we deduce that $C_1=\cdots=C_\ell=0$.

\end{rem}
   \begin{proof}[Proof of Theorem~\ref{thm: CF(G_l;Z)}]

Let $\ell\in \Z_{\geq0}$, and let us compute the differential of $CF(G_\ell;\Z)$. By Remark~\ref{rem: conc from thm: diff_obst for G_l}, the differential $d \fc CF(G_\ell; \Z) \to CF(G_\ell; \Z)$ satisfies
\begin{itemize}
    \item $d \hat{x}^\ell_i = A_{\ell,i} \check{x}^\ell_{i+1} + B_{\ell,i} \check{x}^\ell_{i-n}$, for every $1 \leq i \leq \ell$, 
    \item $d \check{x}^\ell_i = 0$ for every $-(n-1) \leq i \leq \ell+1$,
\end{itemize}
for some $A_{\ell,1},\ldots,A_{\ell,\ell},B_{\ell,1},\ldots,B_{\ell,\ell}\in \Z$.

 By Theorem~\ref{thm: diff_obst_for_complement}, for every $1\leq i\leq \ell$ we know that if $u$ is a Floer bubbled flowline with cascades that connects $\hat{x}_i^{\ell}$ to $\check{x}_{i+1}^\ell$ then $u$ does not intersects $D_\infty$. Therefore the flowlines with cascades that contribute to the differential of $\hat{x}_i^\ell$  in $CF(G_\ell;\Z)$ are the same as those whose contribute the differential of $\hat{x}_i^\ell$ in $CF(\tilde{G}_\ell;\Z)$, where $\tilde{G}_\ell=G_\ell|_{\CP^n\setminus D_\infty}$. Thus by Proposition~\ref{prop: CF of J-shaped in C^n, for CPn-Ball} we deduce that $A_{\ell,1},\ldots, A_{\ell,\ell}\in \{-1,1\}$.

As a result, the Floer cochain complex $CF(G_\ell;\Z)$ decomposes into $n+1$ subcomplexes. For each $-(n-1) \leq j \leq 1$, the submodule 
$$C_{\ell,j}=CF^{-2j-1}(G_\ell;\Z)\oplus CF^{-2j}(G_\ell;\Z)$$
generated by
$$
\check{x}_{j}^\ell,\ \hat{x}_{j+n}^\ell,\ \check{x}_{j+n+1}^\ell,\ \hat{x}_{j+2n+1}^\ell,\ \ldots,\ 
\hat{x}_{j+k_{\ell,j}(n+1)-1}^\ell,\ \check{x}_{j+k_{\ell,j}(n+1)}^\ell,
$$
where $k_{\ell,j} = \left\lfloor \frac{\ell + 1 - j}{n+1} \right\rfloor$, is a subcomplex of $CF(G_\ell;\Z)$. Additionally, by Lemma~\ref{lemma: generators of HF(G_l)} we get  that for every $-(n-1)\leq j\leq 1$ we have 
$$H^*(C_{\ell,j})=\left\{\begin{array}{lc}
    HF^{-2j}(G_\ell;\Z), & *=-2j, \\
     0, & *\neq -2j, 
\end{array}\right.=\left\{\begin{array}{lc}
    \Z\cdot[\check{x}_{j+_{k_{\ell,j}(n+1)}}^\ell], & *=-2j, \\
     0, & *\neq -2j. 
\end{array}\right.$$

 Thus Lemma~\ref{lemma: zigzag complex}, applied for all the cochain subcomplexes $C_{\ell,0},\ldots,C_{\ell,n}$, implies that for every $1\leq i\leq \ell$ we have $B_{\ell,i}\in\{-1,1\}$, which completes the proof.

    \end{proof}

\begin{proof}[Proof of Theorem~\ref{thm: continuation maps for G_l, over Z}]
     As we mentioned in the proof of Theorem~\ref{thm: CF(G_l;Z)}, the Floer cochain complexes $CF(G_\ell;\Z)$, $CF(G'_{\ell'};\Z)$ decompose into $n+1$ subcomplexes. For each $j$ between $-(n-1)$ and $1$, the submodules 
$$C_{\ell,j}=CF^{-2j-1}(G_\ell;\Z)\oplus CF^{-2j}(G_\ell;\Z),\qquad\text{and}\qquad C'_{\ell',j}=CF^{-2j-1}(G'_{\ell'};\Z)\oplus CF^{-2j}(G'_{\ell'};\Z)$$
generated by

$$
\check{x}_{j}^\ell,\ \hat{x}_{j+n}^\ell,\ \check{x}_{j+n+1}^\ell,\ \hat{x}_{j+2n+1}^\ell,\ \ldots,\ 
\hat{x}_{j+k_{\ell,j}(n+1)-1}^\ell,\ \check{x}_{j+k_{\ell,j}(n+1)}^\ell,
$$
and

$$
\check{y}_{j}^{\ell'},\ \hat{y}_{j+n}^{\ell'},\ \check{y}_{j+n+1}^{\ell'},\ \hat{y}_{j+2n+1}^{\ell'},\ \ldots,\ 
\hat{y}_{j+k_{{\ell'},j}(n+1)-1}^{\ell'},\ \check{y}_{j+k_{{\ell'},j}(n+1)}^{\ell'},
$$
where $k_{\ell,j} = \left\lfloor \frac{\ell + 1 - j}{n+1} \right\rfloor$ and $k_{{\ell'},j} = \left\lfloor \frac{{\ell'} + 1 - j}{n+1} \right\rfloor$, are subcomplexes of $CF(G_\ell;\Z)$ and $CF(G'_{\ell'};\Z)$, respectively. Additionally, by Lemma~\ref{lemma: generators of HF(G_l)} we get that for every $-(n-1)\leq j\leq 1$ we have 
$$H^*(C_{\ell,j})=\left\{\begin{array}{lc}
    HF^{-2j}(G_\ell;\Z), & *=-2j, \\
     0, & *\neq -2j,
\end{array}\right.=\left\{\begin{array}{lc}
    \Z\cdot[\check{x}_j^\ell], & *=-2j, \\
     0, & *\neq -2j,
\end{array}\right.$$
and 
$$H^*(C'_{\ell',j})=\left\{\begin{array}{lc}
    HF^{2j}(G'_{\ell'};\Z), & *=-2j, \\
     0, & *\neq -2j,
\end{array}\right.=\left\{
\begin{array}{lc}
    \Z\cdot[\check{y}_j^{\ell'}], & *=-2j, \\
     0, & *\neq -2j.
\end{array}\right.$$
By Theorem~\ref{thm: continuation_obst_for_complement} there are $\check{C}_{-(n-1)},\ldots,\check{C}_{\ell+1},\hat{C}_{1},\ldots,\hat{C}_{\ell}\in \Z$ such that $\Phi\check{x}_i^\ell=\check{C}_i \check{y}_i^{\ell'}$ for every $-(n-1)\leq i\leq \ell+1$ and $\Phi\hat{x}_i^\ell=\hat{C}_i \hat{y}_i^{\ell'}$ for every $1\leq i\leq \ell$.

Since $\CP^n$ is a closed symplectic manifold, the continuation map $\Phi\fc CF(G_\ell;\Z)\to CF(G'_{\ell'};\Z)$ is a quasi-isomorphism. Thus Theorem~\ref{thm: CF(G_l;Z)} and Lemma~\ref{lemma: zigzag and morphisms}, which is applied for all the pairs $(C_{\ell,-(n-1)},C'_{\ell,-(n-1)}),\ldots,(C_{\ell,1},C'_{\ell,1})$ of subcomplexes the pair of cochain complexes $(CF(G_\ell;\Z) ,CF(G'_{\ell'};\Z))$, imply that 
$$|\check{C}_{-(n-1)}|=\cdots=|\check{C}_{\ell+1}|=|\hat{C}_1|=\cdots=|\hat{C}_\ell|=1,$$
as required.

\end{proof}

 \subsection{Energy of Floer flowlines with cascades for the complement of a ball}

Let $\Delta \in (0,1)$. We describe the complexes and connecting maps of the $1$-ray
$$CF(G_0) \to CF(G_1) \to CF(G_2) \to \cdots$$
where for every $\ell \in \Z_{\geq 0}$, the Hamiltonian $G_\ell \fc \CP^n \to \R$ is given by $G_\ell(z) = g(\Delta, \ell, \mu(z))$.

Theorem~\ref{thm: CF(G_l;Z)} and Theorem~\ref{thm: continuation maps for G_l, over Z} describe the Floer complexes of the Hamiltonians $(G_i)_i$ over $\Z$. To describe them with coefficients in the Novikov ring (in the sense of Section~\ref{sss: weighted_CF}), we must compute the energy of the Floer and continuation flowlines with cascades appearing in the differential and continuation maps.

Fix $\ell \in \Z_{\geq 0}$. For every $1 \leq i \leq \ell$, let $r_{\ell,i}^\Delta \in [0,1)$ denote the number satisfying $\frac{\partial}{\partial r}|_{r=r_{\ell,i}^\Delta} g(\Delta, \ell, r) = -i$. For every $0 \leq i \leq \ell$, let $g_{\ell,i}^\Delta = g(\Delta, \ell, r_{\ell,i}^\Delta)$, and let $g_{\ell,0}^\Delta = g(\Delta, \ell, 1)$. Additionally, for every $1 \leq i \leq \ell$, define the $(2n-1)$-spheres
$$S_{\ell,i} = \{z \in \CP^n : g'_\ell(\mu(z)) = -i\},$$
each equipped with a perfect Morse function. We also equip the divisor $D_\infty = \CP^{n-1}$ with a perfect Morse function. 

As before, let $\check{x}_{\ell+1}^\ell$ denote the constant orbit of $G_\ell$ at $0 \in \text{Int } B(1) \subset \CP^n$. For every $1 \leq i \leq \ell$, let $\check{x}^\ell_i$ and $\hat{x}^\ell_i$ denote the minimum and maximum points on $S_{\ell,i}$, respectively. Finally, for every $-(n-1) \leq i \leq 0$, let $\check{x}_i^\ell$ denote the critical point on $D_\infty$ with Morse index $-2i$.
\begin{thm}\label{thm: CF(G_l) + continuation maps}
    For every $\ell\in \Z_{\geq0}$ there exist signs 
    $$A_{\ell,1},\ldots,A_{\ell,\ell},B_{\ell,1},\ldots,B_{\ell,\ell},\check{C}_{\ell,-(n-1)},\ldots,\check{C}_{\ell,\ell+1},\hat{C}_{\ell,1},\ldots,\hat{C}_{\ell,\ell}\in \{-1,1\},$$
    such that for every $\ell\in \Z_{\geq0}$, the differential $d\fc CF(G_\ell)\to CF(G_\ell)$ satisfies
    \begin{itemize}
        \item  $$d \hat{x}^\ell_i=A_{\ell,i}T^{E(\hat{x}_{i}^\ell,\check{x}_{i+1}^\ell)}\check{x}^\ell_{i+1}+B_{\ell,i}T^{E(\hat{x}_{i}^\ell,\check{x}_{i-n}^\ell)}\check{x}^\ell_{i-n},$$
    where $$E(\hat{x}_{i}^\ell,\check{x}_{i+1}^\ell)=\left\{\begin{array}{ll}
        g_{\ell,i+1}^\Delta - g_{\ell,i}^\Delta - i( r_{\ell,i}^\Delta-r_{\ell,i+1}^\Delta) + r_{\ell,i+1}^\Delta, & i\leq \ell-1, \\
       g(\Delta,\ell,0) - g_{\ell,\ell}^\Delta - \ell r_{\ell,\ell}^\Delta, & i=\ell,
    \end{array}\right.$$ and $$E(\hat{x}_{i}^\ell,\check{x}_{i-n}^\ell)=\left\{\begin{array}{ll}
        g_{\ell,i-n}^\Delta - g_{\ell,i}^\Delta + i(r_{\ell,i-n}^\Delta - r_{\ell,i}^\Delta) + n(1 - r_{\ell,i-n}^\Delta), & i\geq n+1, \\
        g_{\ell,0}^\Delta - g_{\ell,i}^\Delta +i(1- r_{\ell,i}^\Delta), & i\leq n,
    \end{array}\right.$$
    for every $1\leq i\leq \ell$.
    
    \item $d\check{x}^\ell_i=0$ for every $-(n-1)\leq i\leq \ell+1$,
    \end{itemize}
   and the continuation map $\Phi_\ell\fc CF(G_\ell)\to CF(G_{\ell+1})$ satisfies 
   
    \begin{itemize}

    \item $$\Phi_\ell \hat{x}^\ell_i=\hat{C}_{\ell,i}T^{E(\hat{x}_i^\ell,\hat{x}_i^{\ell+1})}\hat{x}^{\ell+1}_{i},$$
     where 
    $$E(\hat{x}_i^\ell,\hat{x}_i^{\ell+1})=g_{\ell+1,i}^\Delta - g_{\ell,i}^\Delta - i( r_{\ell,i}^\Delta-r_{\ell+1,i}^\Delta ),$$
    for every $1\leq i\leq \ell$.

        \item  $$\Phi_\ell \check{x}^\ell_i=\check{C}_{\ell,i}T^{E(\check{x}_i^\ell,\check{x}_i^{\ell+1})}\check{x}^{\ell+1}_{i},$$
    where 
    $$E(\check{x}_i^\ell,\check{x}_i^{\ell+1})=\left\{\begin{array}{ll}
      g_{\ell+1,\ell+1}^\Delta+(\ell+1)r_{\ell+1,\ell+1}^\Delta-g(\Delta,\ell,0), & i= \ell+1, \\
      g_{\ell+1,i}^\Delta - g_{\ell,i}^\Delta - i( r_{\ell,i}^\Delta-r_{\ell+1,i}^\Delta),  & 1\leq i\leq \ell,\\
       g_{\ell+1,0}^\Delta - g_{\ell,0}^\Delta,& i\leq0,
    \end{array}\right.$$
    for every $-(n-1)\leq i\leq \ell+1$,

    \end{itemize}
\end{thm}

As described in Section~\ref{ss: energy}, the topological energy of a Floer or continuation flowline with cascades $u$ connecting the $1$-periodic orbits $p$ and $q$ is given by
$$ E_{top}(u) = \cA_H(q, \hat{q}) - \cA_H(p, \hat{p}) + \frac{1}{2(n+1)} \left( \dim_u \cM_m(p, q) - \mu_{FMB}^\tau(q) + \mu_{FMB}^\tau(p) - m + 1 \right), $$
where $\hat{p}$ and $\hat{q}$ are cappings for $p$ and $q$, respectively. The trivialization $\tau$ is induced by these cappings, and $\cA_H$ denotes the action functional associated with the Hamiltonian in the case of the differential, or the homotopy in the case of continuation maps.

\begin{proof}[Proof of Theorem~\ref{thm: CF(G_l) + continuation maps}]Let $\ell \in \Z_{\geq 0}$. We begin by computing the topological energy of flowlines with cascades that contribute to the differential. In this case, the dimension $\dim_u \cM_m(p, q)$ equals $m$ for every $p, q$, and $u$. Thus, the topological energy of such a flowline with cascades is given by
$$ E_{top}(u) = \cA_{G_\ell}(q, \hat{q}) - \cA_{G_\ell}(p, \hat{p}) + \frac{1}{2(n+1)} \left( \mu_{FMB}^\tau(p) - \mu_{FMB}^\tau(q) + 1 \right). $$
Note that for constant orbits, we can choose the cappings to be constant. Consequently, their action equals the value of the Hamiltonian itself. Thus,
$$ \cA_{G_\ell}(\check{x}_{\ell+1}^\ell) = g(\Delta, \ell, 0), \qquad \text{and} \qquad \cA_{G_\ell}(\check{x}_{-(n-1)}^\ell) = \cdots = \cA_{G_\ell}(\check{x}_{0}^\ell) = g(\Delta, \ell, 1)=g_{\ell,0}^\Delta. $$
For the orbits $\hat{x}_1^\ell, \check{x}_1^\ell, \ldots, \hat{x}_\ell^\ell, \check{x}_\ell^\ell$, we choose cappings contained in $\CP^n \setminus D_\infty$, which is an exact symplectic manifold symplectomorphic to $\Int B(1)$. Let $\lambda = \frac{1}{2} \sum_{i=1}^n (x_i dy_i - y_i dx_i)$ be a primitive of $\omega_0$ on $\Int B(1)$. Then, for every $1 \leq i \leq \ell$, we find that
$$ \cA_{G_\ell}(\hat{x}_i^\ell) = \cA_{G_\ell}(\check{x}_i^\ell) = \int_{S^1} G_\ell \circ \check{x}_i^\ell(t) \, dt + \int_{S^1} (\check{x}_i^\ell)^* \lambda = g_{\ell,i}^\Delta + i r_{\ell,i}^\Delta. $$
Additionally, the Floer--Morse--Bott indices for these $1$-periodic orbits with respect to the trivialization $\tau$ induced by these cappings are as follows:
\begin{itemize}
    \item $\mu_{FMB}^\tau(\check{x}_{\ell+1}^\ell) = 2n(\ell+1)$.
    \item For every $1 \leq i \leq \ell$, $\mu_{FMB}^\tau(\check{x}_i^\ell) = 2ni$.
    \item For every $1 \leq i \leq \ell$, $\mu_{FMB}^\tau(\hat{x}_i^\ell) = 2ni + 2n - 1$.
    \item For every $-(n-1) \leq i \leq 0$, $\mu_{FMB}^\tau(\check{x}_i^\ell) = -2i$.
\end{itemize}

The only Floer flowlines with cascades that contribute to the differential are those connecting $\hat{x}_i^\ell$ to $\check{x}_{i+1}^\ell$ and $\check{x}_{i-n}^\ell$ for $1 \leq i \leq \ell$. Thus, for every $1 \leq i \leq \ell-1$, we have
$$ \mu_{FMB}^\tau(\hat{x}_{i}^\ell) -\mu_{FMB}^\tau(\check{x}^\ell_{i+1}) =  (2ni + 2n - 1) -2n(i+1) = -1,$$
and hence
\begin{align*}
E(\hat{x}_i^\ell, \check{x}_{i+1}^\ell) &= \cA_{G_\ell}(\check{x}_{i+1}^\ell) - \cA_{G_\ell}(\hat{x}_i^\ell) \\
&= g_{\ell,i+1}^\Delta + (i+1)r_{\ell,i+1}^\Delta - g_{\ell,i}^\Delta - ir_{\ell,i}^\Delta \\
&= g_{\ell,i+1}^\Delta - g_{\ell,i}^\Delta - i( r_{\ell,i}^\Delta-r_{\ell,i+1}^\Delta) + r_{\ell,i+1}^\Delta.
\end{align*}
Additionally, for $i=\ell$ we have 
$$ \mu_{FMB}^\tau(\hat{x}^\ell_{\ell})-\mu_{FMB}^\tau(\check{x}_{\ell+1}^\ell)  =  (2n\ell + 2n - 1) -2n(\ell+1) = -1,$$
and hence
$$E(\hat{x}_\ell^\ell, \check{x}_{\ell+1}^\ell) = \cA_{G_\ell}(\check{x}_{\ell+1}^\ell) - \cA_{G_\ell}(\hat{x}_\ell^\ell)= g(\Delta,\ell,0) - g_{\ell,\ell}^\Delta - \ell r_{\ell,\ell}^\Delta .$$

Now consider the second case for $1 \leq i \leq \ell$. If $i \geq n+1$, we find
$$\mu_{FMB}^\tau(\hat{x}_{i}^\ell) - \mu_{FMB}^\tau(\check{x}^\ell_{i-n}) = 2ni + 2n - 1 - (2n(i-n)) = 2n(n+1) - 1,$$
which implies 
\begin{align*}
E(\hat{x}_i^\ell, \check{x}_{i-n}^\ell) &= \cA_{G_\ell}(\check{x}_{i-n}^\ell) - \cA_{G_\ell}(\hat{x}_i^\ell) + \frac{1}{2(n+1)} \left( \mu_{FMB}^\tau(\hat{x}_i^\ell) - \mu_{FMB}^\tau(\check{x}_{i-n}^\ell) + 1 \right) \\
&= g_{\ell,i-n}^\Delta + (i-n)r_{\ell,i-n}^\Delta - g_{\ell,i}^\Delta - ir_{\ell,i}^\Delta + n \\
&= g_{\ell,i-n}^\Delta - g_{\ell,i}^\Delta + i(r_{\ell,i-n}^\Delta - r_{\ell,i}^\Delta) + n(1 - r_{\ell,i-n}^\Delta).
\end{align*}

On the other hand, if $i \leq n$, we have
$$\mu_{FMB}^\tau(\hat{x}_{i}^\ell) - \mu_{FMB}^\tau(\check{x}^\ell_{i-n}) = 2ni + 2n - 1 -(- 2(i-n)) = 2i(n+1) - 1,$$
yielding 
\begin{align*}
E(\hat{x}_i^\ell, \check{x}_{i-n}^\ell) &= \cA_{G_\ell}(\check{x}_{i-n}^\ell) - \cA_{G_\ell}(\hat{x}_i^\ell) + \frac{1}{2(n+1)} \left( \mu_{FMB}^\tau(\hat{x}_i^\ell) - \mu_{FMB}^\tau(\check{x}_{i-n}^\ell) + 1 \right) \\
&= g_{\ell,0}^\Delta - g_{\ell,i}^\Delta - ir_{\ell,i}^\Delta + i\\
&= g_{\ell,0}^\Delta - g_{\ell,i}^\Delta +i(1- r_{\ell,i}^\Delta).
\end{align*}

Now we can move to computing the topological energy of flowlines with cascades that contribute to the continuation map $\Phi_\ell \fc CF(G_\ell) \to CF(G_{\ell+1})$. In this case, the dimension $\dim_u \cM_m(p, q) $ equals $m-1$ for every $p, q$, and $u$. Thus, the topological energy of such a flowline is given by
$$ E_{top}(u) = \cA_{G_{\ell+1}}(q) - \cA_{G_\ell}(p) + \frac{1}{2(n+1)} \left( \mu_{FMB}^\tau(p) - \mu_{FMB}^\tau(q) \right). $$
As before, for constant orbits of $G_{\ell+1}$, we can choose the cappings to be constant. Consequently, their action equals the value of the Hamiltonian itself. Thus, $\cA_{G_{\ell+1}}(\check{x}_{\ell+2}^{\ell+1}) = g(\Delta, \ell+1, 0) $ and
$$ \cA_{G_{\ell+1}}(\check{x}_{-(n-1)}^{\ell+1}) = \cdots = \cA_{G_{\ell+1}}(\check{x}_{0}^{\ell+1}) = g(\Delta, \ell+1, 1=g_{\ell+1,0}^\Delta). $$
For the orbits $\hat{x}_1^{\ell+1}, \check{x}_1^{\ell+1}, \ldots, \hat{x}_{\ell+1}^{\ell+1}, \check{x}_{\ell+1}^{\ell+1}$, we again choose cappings contained in $\CP^n \setminus D_\infty$. Thus, for every $1 \leq i \leq \ell+1$, we find
$$ \cA_{G_{\ell+1}}(\hat{x}_i^{\ell+1}) = \cA_{G_{\ell+1}}(\check{x}_i^{\ell+1}) = \int_{S^1} G_{\ell+1} \circ \check{x}_i^{\ell+1}(t) \, dt + \int_{S^1} (\check{x}_i^{\ell+1})^* \lambda = g_{\ell+1,i}^\Delta + i r_{\ell+1,i}^\Delta. $$
Additionally, the Floer--Morse--Bott indices for these $1$-periodic orbits with respect to the trivialization $\tau$ are as follows:
\begin{itemize}
    \item $\mu_{FMB}^\tau(\check{x}_{\ell+2}^{\ell+1}) = 2n(\ell+2)$.
    \item For $1 \leq i \leq \ell+1$, $\mu_{FMB}^\tau(\check{x}_i^{\ell+1}) = 2ni$.
    \item For $1 \leq i \leq \ell+1$, $\mu_{FMB}^\tau(\hat{x}_i^{\ell+1}) = 2ni + 2n - 1$.
    \item For $-(n-1) \leq i \leq 0$, $\mu_{FMB}^\tau(\check{x}_i^{\ell+1}) = -2i$.
\end{itemize}

Given $1 \leq i \leq \ell$, note that $\mu_{FMB}^\tau(\hat{x}_i^\ell) = \mu_{FMB}^\tau(\hat{x}_i^{\ell+1})$, and hence
\begin{align*}
    E(\hat{x}_i^\ell, \hat{x}_i^{\ell+1}) &= \cA_{G_{\ell+1}}(\hat{x}_i^{\ell+1}) - \cA_{G_\ell}(\hat{x}_i^\ell) \\
    &= g_{\ell+1,i}^\Delta + i r_{\ell+1,i}^\Delta - g_{\ell,i}^\Delta - i r_{\ell,i}^\Delta \\
    &= g_{\ell+1,i}^\Delta - g_{\ell,i}^\Delta - i( r_{\ell,i}^\Delta-r_{\ell+1,i}^\Delta ).
\end{align*}

Now, let $-(n-1) \leq i \leq \ell+1$.
\begin{itemize}
\item If $i=\ell+1$, then $\mu_{FMB}^\tau(\check{x}_{\ell+1}^\ell) = \mu_{FMB}^\tau(\check{x}_{\ell+1}^{\ell+1})$, yielding
$$E(\check{x}_{\ell+1}^\ell, \check{x}_{\ell+1}^{\ell+1}) = \cA_{G_{\ell+1}}(\check{x}_{\ell+1}^{\ell+1}) - \cA_{G_\ell}(\check{x}_{\ell+1}^\ell)= g_{\ell+1,\ell+1}^\Delta+(\ell+1)r_{\ell+1,\ell+1}^\Delta-g(\Delta,\ell,0).$$
    \item If $1\leq i \leq \ell$, then $\mu_{FMB}^\tau(\check{x}_i^\ell) = \mu_{FMB}^\tau(\check{x}_i^{\ell+1})$, yielding
    $$   E(\check{x}_i^\ell, \check{x}_i^{\ell+1}) = \cA_{G_{\ell+1}}(\check{x}_i^{\ell+1}) - \cA_{G_\ell}(\check{x}_i^\ell) = g_{\ell+1,i}^\Delta - g_{\ell,i}^\Delta - i( r_{\ell,i}^\Delta-r_{\ell+1,i}^\Delta).$$

    \item If $i \leq 0$, then $\mu_{FMB}^\tau(\check{x}_i^\ell) = \mu_{FMB}^\tau(\check{x}_i^{\ell+1})$, and hence
    $$ E(\check{x}_i^\ell, \check{x}_i^{\ell+1}) = \cA_{G_{\ell+1}}(\check{x}_i^{\ell+1}) - \cA_{G_\ell}(\check{x}_i^\ell) = g_{\ell+1,0}^\Delta - g_{\ell,0}^\Delta. $$
\end{itemize}
This completes the proof.

\end{proof}

\subsection{Direct limit of Floer complexes}

Throughout this section, we fix $\Delta \in (0,1)$. As usual, for every $\ell \in \Z_{\geq0}$, define $g_\ell \fc [0,1] \to \R$ by $g_\ell(r) = g(\Delta, \ell, r)$ for every $r \in [0,1]$ and focus on the Hamiltonians $G_\ell$ on $\CP^n$ given by 
$G_\ell(z) = g_\ell(\mu(z)) = g(\Delta, \ell, \mu(z))$, for every $z \in \CP^n$.

The purpose of this section is to compute the direct limit of the 1-ray of the cochain Floer complexes
$$CF(G_0) \to CF(G_1) \to CF(G_2) \to \cdots$$
where the connecting maps are the continuation maps. The description of the Floer complexes and the continuation maps can be found in Theorem~\ref{thm: CF(G_l) + continuation maps}.

Define a cochain complex $(C, d)$ as follows: the module $C$ is defined as
$$C = \Lambda_{>0}\check{x}_{-(n-1)}\oplus\cdots\oplus \Lambda_{>0}\check{x}_0\oplus \bigoplus_{j\in \N} \left( \Lambda_{>-j\Delta}\check{x}_j \oplus \Lambda_{>-j\Delta}\hat{x}_j \right),$$
and the differential $d$ satisfies $d(\lambda\check{x}_j) = 0$ for every $\lambda\in \Lambda_{>0}$ when $-(n-1)\leq j\leq 0$, and $\lambda \in \Lambda_{>-j\Delta}$ when $j\geq 1$, 
$$d(\lambda\hat{x}_j) = \lambda\check{x}_{j+1} + T^n\lambda\check{x}_{j-n}$$
for every $j \in \Z_{\geq n+1}$ and $\lambda \in \Lambda_{>-j\Delta}$, and also 
$$d(\lambda\hat{x}_j) = \lambda\check{x}_{j+1} + T^j\lambda\check{x}_{j-n}$$
for every $j \in \{1\ldots,n\}$ and $\lambda \in \Lambda_{>-j\Delta}$.

\begin{thm}\label{thm: drct lim of CF(G_l)}
    The direct limit of the 1-ray
    $$CF(G_0) \to CF(G_1) \to CF(G_2) \to \cdots$$
    is isomorphic to $(C, d)$. Additionally, the canonical maps $(f_\ell\fc CF(G_\ell)\to C)_{\ell\geq0}$ satisfy,     for every $\ell\in \Z_{\geq0}$, that:
     \begin{itemize}
         \item for every $-(n-1)\leq j \leq 0$ we have $f_\ell(\check{x}_j^\ell) = T^{-g_{\ell,0}^\Delta}\check{x}_j$,
        \item for every $1 \leq j \leq \ell$ we have $f_\ell(\check{x}_j^\ell) =  T^{-(g_{\ell,j}^\Delta+jr_{\ell,j}^\Delta)}\check{x}_j$, 
        \item for $j=\ell+1$ we have $f_\ell(\check{x}_{\ell+1}^\ell)=T^{-g(\Delta,\ell,0)}\check{x}_{\ell+1}$,
        \item for every $1 \leq j \leq \ell$ we have $f_\ell(\hat{x}_j^\ell) =  T^{-(g_{\ell,j}^\Delta+jr_{\ell,j}^\Delta)}\hat{x}_j$,
    \end{itemize}
\end{thm}

First, let us show that we can ignore the signs in the differentials and continuation maps associated with the Floer complexes described in Theorem~\ref{thm: CF(G_l) + continuation maps}. Specifically, we will prove that one can change the bases of the complexes $CF(G_\ell)$ so that all signs appearing in the formulas for these differentials and continuation maps are equal to $1$. This result is a consequence of Lemma~\ref{lemma: zigzag complex}.

\begin{prop}\label{prop: omission of signs for CF(G)}
    For every $\ell\in \Z_{\geq0}$ there is a basis
    $$\hat{x}_1^\ell,\ldots,\hat{x}_\ell^\ell,\check{x}_{-(n-1)}^\ell,\ldots,\check{x}_{\ell+1}^\ell$$ 
    for $CF(G_\ell)$ such that the differential $d\fc CF(G_\ell)\to CF(G_\ell)$ satisfies
    \begin{itemize}
        \item  $$d \hat{x}^\ell_i=T^{E(\hat{x}_{i}^\ell,\check{x}_{i+1}^\ell)}\check{x}^\ell_{i+1}+T^{E(\hat{x}_{i}^\ell,\check{x}_{i-n}^\ell)}\check{x}^\ell_{i-n},$$
    where $$E(\hat{x}_{i}^\ell,\check{x}_{i+1}^\ell)=\left\{\begin{array}{ll}
        g_{\ell,i+1}^\Delta - g_{\ell,i}^\Delta - i( r_{\ell,i}^\Delta-r_{\ell,i+1}^\Delta) + r_{\ell,i+1}^\Delta, & i\leq \ell-1, \\
       g(\Delta,\ell,0) - g_{\ell,\ell}^\Delta - \ell r_{\ell,\ell}^\Delta, & i=\ell,
    \end{array}\right.$$ and $$E(\hat{x}_{i}^\ell,\check{x}_{i-n}^\ell)=\left\{\begin{array}{ll}
        g_{\ell,i-n}^\Delta - g_{\ell,i}^\Delta + i(r_{\ell,i-n}^\Delta - r_{\ell,i}^\Delta) + n(1 - r_{\ell,i-n}^\Delta), & i\geq n+1, \\
        g_{\ell,0}^\Delta - g_{\ell,i}^\Delta +i(1- r_{\ell,i}^\Delta), & i\leq n,
    \end{array}\right.$$
    for every $1\leq i\leq \ell$.
    
    \item $d\check{x}^\ell_i=0$ for every $-(n-1)\leq i\leq \ell+1$,
    \end{itemize}
   and the continuation map $\Phi_\ell\fc CF(G_\ell)\to CF(G_{\ell+1})$ satisfies 
   
    \begin{itemize}

    \item $$\Phi_\ell \hat{x}^\ell_i=T^{E(\hat{x}_i^\ell,\hat{x}_i^{\ell+1})}\hat{x}^{\ell+1}_{i},$$
     where 
    $$E(\hat{x}_i^\ell,\hat{x}_i^{\ell+1})=g_{\ell+1,i}^\Delta - g_{\ell,i}^\Delta - i( r_{\ell,i}^\Delta-r_{\ell+1,i}^\Delta ),$$
    for every $1\leq i\leq \ell$.

        \item  $$\Phi_\ell \check{x}^\ell_i=T^{E(\check{x}_i^\ell,\check{x}_i^{\ell+1})}\check{x}^{\ell+1}_{i},$$
    where 
    $$E(\check{x}_i^\ell,\check{x}_i^{\ell+1})=\left\{\begin{array}{ll}
      g_{\ell+1,\ell+1}^\Delta+(\ell+1)r_{\ell+1,\ell+1}^\Delta-g(\Delta,\ell,0), & i= \ell+1, \\
      g_{\ell+1,i}^\Delta - g_{\ell,i}^\Delta - i( r_{\ell,i}^\Delta-r_{\ell+1,i}^\Delta),  & 1\leq i\leq \ell,\\
       g_{\ell+1,0}^\Delta - g_{\ell,0}^\Delta,& i\leq0,
    \end{array}\right.$$
    for every $-(n-1)\leq i\leq \ell+1$,

    \end{itemize}
\end{prop}

\begin{proof}

For simplicity, throughout this proof we omit the powers of the formal variable $T$ in the Novikov ring, as they do not affect the signs; thus, we translate the discussion to the case of cochain complexes with $\Z$-coefficients. In our proof, we will show that by changing the signs of our given basis for $CF(G_\ell;\Z)$ for every $\ell\in \Z_{\geq0}$, we can achieve the desired result of this proposition.

By Theorem~\ref{thm: CF(G_l) + continuation maps}, for every $\ell \in \Z_{\geq 0}$, the cochain complex $CF(G_\ell)$ decomposes into $n+1$ subcomplexes. For each $-(n-1) \leq j \leq 1$, denote by $C_{\ell,j}$ the subcomplex generated by
$$
\check{x}_{j}^\ell,\ \hat{x}_{j+n}^\ell,\ \check{x}_{j+n+1}^\ell,\ \hat{x}_{j+2n+1}^\ell,\ \ldots,\ 
\hat{x}_{j+k_{\ell,j}(n+1)-1}^\ell,\ \check{x}_{j+k_{\ell,j}(n+1)}^\ell,
$$
where $k_{\ell,j} = \left\lfloor \frac{\ell + 1 - j}{n+1} \right\rfloor$.

Fix $-(n-1) \leq j \leq 1$. Define $C_{-1,j}=0$, let $\Phi_{-1}\fc C_{-1,j}\to C_{0,j}$ be the zero map, and denote $k_{-1,j}=-1$. Let us prove by induction that for every $\ell\in \Z_{\geq0}$, we can change the signs of the basis elements of the subcomplex $C_{\ell,j}$ such that:
\begin{itemize}
    \item For every $1\leq i\leq k_{\ell,j}$ we have $d\hat{x}_{j+i(n+1)-1}^\ell=\check{x}_{j+(i-1)(n+1)}^\ell+\check{x}_{j+i(n+1)}^\ell$;
    \item For every $1\leq i\leq k_{\ell-1,j}$ we have $\Phi_{\ell-1}(\hat{x}_{j+i(n+1)-1}^{\ell-1})= \hat{x}_{j+i(n+1)-1}^{\ell}$;
    \item For every $0\leq i\leq k_{\ell-1,j}$ we have $\Phi_{\ell-1}(\check{x}_{j+i(n+1)}^{\ell-1})= \check{x}_{j+i(n+1)}^{\ell}$.
\end{itemize} 

\textbf{The base case of the induction:} For $\ell=0$, we have $k_{0,j}=0$, therefore all the conditions of the induction statement are satisfied in the trivial sense.
 
\textbf{The step case of the induction:} Let $\ell\in \N$. Assume that for every $0\leq \ell'\leq \ell$, we can change the signs of the basis elements of the subcomplex $C_{\ell',j}$ such that:
\begin{itemize}
   \item For every $1\leq i\leq k_{\ell',j}$ we have $d\hat{x}_{j+i(n+1)-1}^{\ell'}=\check{x}_{j+(i-1)(n+1)}^{\ell'}+\check{x}_{j+i(n+1)}^{\ell'}$;
    \item For every $1\leq i\leq k_{\ell'-1,j}$ we have $\Phi_{\ell'-1}(\hat{x}_{j+i(n+1)-1}^{\ell'-1})= \hat{x}_{j+i(n+1)-1}^{\ell'}$;
    \item For every $0\leq i\leq k_{\ell'-1,j}$ we have $\Phi_{\ell'-1}(\check{x}_{j+i(n+1)}^{\ell'-1})= \check{x}_{j+i(n+1)}^{\ell'}$.
\end{itemize} 
 
If $\ell\leq j+n-1$, then $k_{\ell,j}=0$, and therefore the first two conditions of the induction statement are satisfied in the trivial sense. By Theorem~\ref{thm: continuation maps for G_l, over Z} we know that there is $\epsilon\in\{-1,1\}$ such that $\Phi_{\ell-1}(\check{x}_j^{\ell-1})=\epsilon\check{x}_j^\ell$, thus replacing $\check{x}_j^\ell$ by $\epsilon\check{x}_j^\ell$ completes the proof in this case. 

If $\ell=j+n$, then the second condition is trivial. Moreover, we have
$$C_{\ell-1,j}=\Z\cdot \check{x}_j^{\ell-1}\qquad\text{and}\qquad C_{\ell,j}=\Z\cdot \check{x}_j^{\ell}\oplus \Z\cdot \hat{x}_{j+n}^{\ell}\oplus \Z\cdot \check{x}_{j+n+1}^{\ell},$$
and by Theorem~\ref{thm: CF(G_l;Z)} and Theorem~\ref{thm: continuation maps for G_l, over Z} there are $\alpha,\beta,\gamma\in\{-1,1\}$ such that $\Phi_{\ell-1}(\check{x}_j^{\ell-1})=\alpha\check{x}_j^\ell$ and $d\hat{x}_{j+n}^\ell=\beta\check{x}_j^\ell+\gamma\check{x}_{j+n+1}^\ell$. Then replace the basis $\check{x}_j^\ell, \hat{x}_{j+n}^\ell, \check{x}_{j+n+1}^\ell$ by the basis $\alpha\check{x}_j^\ell, \alpha\beta \hat{x}_{j+n}^\ell, \alpha\beta\gamma\check{x}_{j+n+1}^\ell$. This change of basis works since
$$d(\alpha\beta\hat{x}_{j+n}^\ell)=\alpha\beta^2\check{x}_j^\ell+\alpha\beta\gamma\check{x}_{j+n+1}^\ell=\alpha\check{x}_j^\ell+\alpha\beta\gamma\check{x}_{j+n+1}^\ell.$$
Thus, the first and the third conditions hold as well.

Otherwise, assume that $\ell\geq j+n+1$. This case is fully covered by Lemma~\ref{lemma: positive signs for zigzag and maps}, and this completes the induction step and the whole proof.

\end{proof}

Now we move to the proof of Theorem~\ref{thm: drct lim of CF(G_l)}.

\begin{proof}[Proof of Theorem~\ref{thm: drct lim of CF(G_l)}]
	
	Let us compute the direct limit of the 1-ray
    
    $$CF(G_0) \xrightarrow{\Phi_0} CF(G_1) \xrightarrow{\Phi_1} CF(G_2) \longrightarrow \cdots.$$
	
	Let $j \in \Z_{\geq -(n-1)}$. The sequence $(g_{\ell,0}^\Delta)_{\ell=0}^\infty$ is strictly increasing and its limit is $0$, and if $j\geq 1$, the sequence $(g_{\ell,j}^\Delta+jr_{\ell,j}^\Delta)_{\ell=j}^\infty$ is strictly increasing and its limit is $j\Delta$. By Lemma~\ref{lemma: drct lim of 1-ray}, if $j\leq 0$ the module $\Lambda_{>0}$ is a direct limit for the 1-ray
	$$\Lambda_{\geq 0}\check{x}_j^0 \xrightarrow{\Phi_0} \Lambda_{\geq 0}\check{x}_j^{1} \xrightarrow{\Phi_{1}} \Lambda_{\geq 0}\check{x}_j^{2} \xrightarrow{\Phi_{2}} \cdots,$$
    and if $j\geq 1$ the module $\Lambda_{>-j\Delta}$ is a direct limit for the 1-ray
	$$\Lambda_{\geq 0}\check{x}_j^j \xrightarrow{\Phi_j} \Lambda_{\geq 0}\check{x}_j^{j+1} \xrightarrow{\Phi_{j+1}} \Lambda_{\geq 0}\check{x}_j^{j+2} \xrightarrow{\Phi_{j+2}} \cdots.$$
	If $j \geq 1$, the same argument shows that $\Lambda_{>-j\Delta}$ is a direct limit for the 1-ray
	$$\Lambda_{\geq 0}\hat{x}_j^j \xrightarrow{\Phi_j} \Lambda_{\geq 0}\hat{x}_j^{j+1} \xrightarrow{\Phi_{j+1}} \Lambda_{\geq 0}\hat{x}_j^{j+2} \xrightarrow{\Phi_{j+2}} \cdots.$$
	Since direct limits and direct sums commute, we conclude that the $\Lambda_{\geq 0}$-module
	$$C = \Lambda_{>0}\check{x}_{-(n-1)}\oplus\cdots\oplus \Lambda_{>0}\check{x}_0\oplus \bigoplus_{j\in \N} \left( \Lambda_{>-j\Delta}\check{x}_j \oplus \Lambda_{>-j\Delta}\hat{x}_j \right),$$
	is a direct limit in the category of modules of the 1-ray
		$$CF(G_0) \xrightarrow{\Phi_0} CF(G_1)\xrightarrow{\Phi_{1}} CF(G_2)\xrightarrow{\Phi_{2}} \cdots.$$

	For every $\ell \in \Z_{\geq 0}$, let $f_\ell \fc CF(G_\ell) \to C$ denote the $\Lambda_{\geq 0}$-linear homomorphisms guaranteed by the universal property of the direct limit. By Lemma~\ref{lemma: drct lim of 1-ray}, we deduce that for every $\ell \in \Z_{\geq 0}$:
    \begin{itemize}
        \item for every $-(n-1)\leq j \leq 0$ we have $f_\ell(\check{x}_j^\ell) = T^{-g_{\ell,0}^\Delta}\check{x}_j$,
        \item for every $1 \leq j \leq \ell$ we have $f_\ell(\check{x}_j^\ell) =  T^{-(g_{\ell,j}^\Delta+jr_{\ell,j}^\Delta)}\check{x}_j$, 
        \item for $j=\ell+1$ we have $f_\ell(\check{x}_{\ell+1}^\ell)=T^{-g(\Delta,\ell,0)}\check{x}_{\ell+1}$,
        \item for every $1 \leq j \leq \ell$ we have $f_\ell(\hat{x}_j^\ell) =  T^{-(g_{\ell,j}^\Delta+jr_{\ell,j}^\Delta)}\hat{x}_j$.
    \end{itemize}
    Additionally, these maps satisfy $f_i = f_j \circ \Phi_i^j$ for every $j > i \geq 0$, where $\Phi_i^j = \Phi_{j-1} \circ \cdots \circ \Phi_i$.
	Lemma~\ref{lemma: drct lim in mod cat is drct lim in ch cat} implies that there is a differential $d$ on $C$ such that $(C, d)$ is a direct limit in the category of cochain complexes of the 1-ray, and $f_j \fc CF(G_\ell) \to (C, d)$ is a cochain map for every $j \in \Z_{\geq 0}$. 
	
	Now, let us compute the differential $d$:
	\begin{itemize}
		\item Let $-(n-1)\leq j\leq 0$ and $\lambda \in \Lambda_{>0}$. There exists $\ell \in \Z_{\geq0}$ such that $\lambda T^{g_{\ell,0}^\Delta} \in \Lambda_{\geq 0}$, which implies $f_\ell(T^{g_{\ell,0}^\Delta}\lambda \check{x}_j^\ell) = \lambda\check{x}_j$. Since $d_\ell\check{x}_j^\ell = 0$ and $f_\ell$ is a cochain map, we have:
		$$d(\lambda\check{x}_j) = d(f_\ell(T^{g_{\ell,0}^\Delta}\lambda \check{x}_j^\ell)) = f_\ell(d_\ell(T^{g_{\ell,j}^\Delta}\lambda \check{x}_j^\ell)) = 0.$$

        	\item Let $j \in \Z_{\geq 0}$ and $\lambda \in \Lambda_{>-j\Delta}$. There exists $\ell \geq j$ such that $\lambda T^{g_{\ell,j}^\Delta+jr_{\ell,j}^\Delta} \in \Lambda_{\geq 0}$, which implies $f_\ell(T^{g_{\ell,j}^\Delta+jr_{\ell,j}^\Delta}\lambda \check{x}_j^\ell) = \lambda\check{x}_j$. Since $d_\ell\check{x}_j^\ell = 0$ and $f_\ell$ is a cochain map, we have:
		$$d(\lambda\check{x}_j) = d(f_\ell(T^{g_{\ell,j}^\Delta+jr_{\ell,j}^\Delta}\lambda \check{x}_j^\ell)) = f_\ell(d_\ell(T^{g_{\ell,j}^\Delta+jr_{\ell,j}^\Delta}\lambda \check{x}_j^\ell)) = 0.$$

		\item Let $j \in \N$ and $\lambda \in \Lambda_{>-j\Delta}$. As before, there exists $\ell \geq j+1$ such that $\lambda T^{g_{\ell,j}^\Delta+jr_{\ell,j}^\Delta} \in \Lambda_{\geq 0}$, so $f_\ell(T^{g_{\ell,j}^\Delta+jr_{\ell,j}^\Delta}\lambda \hat{x}_j^\ell) = \lambda\hat{x}_j$. Since 
		$$d_\ell\hat{x}^\ell_j=T^{E(\hat{x}_{j}^\ell,\check{x}_{j+1}^\ell)}\check{x}^\ell_{j+1}+T^{E(\hat{x}_{j}^\ell,\check{x}_{j-n}^\ell)}\check{x}^\ell_{j-n},$$
    where $$E(\hat{x}_{j}^\ell,\check{x}_{j+1}^\ell)=g_{\ell,j+1}^\Delta - g_{\ell,j}^\Delta - j( r_{\ell,j}^\Delta-r_{\ell,j+1}^\Delta) + r_{\ell,j+1}^\Delta,$$ and $$E(\hat{x}_{j}^\ell,\check{x}_{j-n}^\ell)=\left\{\begin{array}{ll}
        g_{\ell,j-n}^\Delta - g_{\ell,j}^\Delta + j(r_{\ell,j-n}^\Delta - r_{\ell,j}^\Delta) + n(1 - r_{\ell,j-n}^\Delta), & j\geq n+1, \\
        g_{\ell,0}^\Delta - g_{\ell,j}^\Delta +j(1- r_{\ell,j}^\Delta), & j\leq n,
    \end{array}\right.$$
		and $f_\ell$ is a cochain map, we compute:
		\begin{align*}
			d(\lambda\hat{x}_j) &= d(f_\ell(\lambda T^{g_{\ell,j}^\Delta+jr_{\ell,j}^\Delta} \hat{x}_j^\ell)) \\
			&= \lambda T^{g_{\ell,j}^\Delta+jr_{\ell,j}^\Delta} f_\ell(d_\ell \hat{x}_j^\ell) \\
			&= \lambda T^{g_{\ell,j}^\Delta+jr_{\ell,j}^\Delta} f_\ell\left(T^{E(\hat{x}_{j}^\ell,\check{x}_{j+1}^\ell)}\check{x}^\ell_{j+1}+T^{E(\hat{x}_{j}^\ell,\check{x}_{j-n}^\ell)}\check{x}^\ell_{j-n}\right) 
		\end{align*}
        Since 
        $$ \lambda T^{g_{\ell,j}^\Delta+jr_{\ell,j}^\Delta} f_\ell\left(T^{E(\hat{x}_{j}^\ell,\check{x}_{j+1}^\ell)}\check{x}^\ell_{j+1}\right)= \lambda T^{g_{\ell,j}^\Delta+jr_{\ell,j}^\Delta} T^{-(g_{\ell,j+1}^\Delta+(j+1)r_{\ell,j+1}^\Delta)}T^{E(\hat{x}_{j}^\ell,\check{x}_{j+1}^\ell)}\check{x}_{j+1}$$
        
        and
        $$ g_{\ell,j}^\Delta+jr_{\ell,j}^\Delta-(g_{\ell,j+1}^\Delta+(j+1)r_{\ell,j+1}^\Delta)+E(\hat{x}_{j}^\ell,\check{x}_{j+1}^\ell)=0,$$
        we deduce that

    $$ \lambda T^{g_{\ell,j}^\Delta+jr_{\ell,j}^\Delta} f_\ell\left(T^{E(\hat{x}_{j}^\ell,\check{x}_{j+1}^\ell)}\check{x}^\ell_{j+1}\right)= \lambda \check{x}_{j+1}$$

Similarly, if $j\geq n+1$ since
$$\lambda T^{g_{\ell,j}^\Delta+jr_{\ell,j}^\Delta} f_\ell\left(T^{E(\hat{x}_{j}^\ell,\check{x}_{j-n}^\ell)}\check{x}^\ell_{j-n}\right)=\lambda T^{g_{\ell,j}^\Delta+jr_{\ell,j}^\Delta} T^{-(g_{\ell,j-n}^\Delta+(j-n)r_{\ell,j-n}^\Delta)}T^{E(\hat{x}_{j}^\ell,\check{x}_{j-n}^\ell)}\check{x}_{j-n}$$
      and  $$g_{\ell,j}^\Delta+jr_{\ell,j}^\Delta-(g_{\ell,j-n}^\Delta+(j-n)r_{\ell,j-n}^\Delta)+E(\hat{x}_{j}^\ell,\check{x}_{j-n}^\ell)=n,$$     
we deduce that
$$\lambda T^{g_{\ell,j}^\Delta+jr_{\ell,j}^\Delta} f_\ell\left(T^{E(\hat{x}_{j}^\ell,\check{x}_{j-n}^\ell)}\check{x}^\ell_{j-n}\right)=\lambda T^{n}\check{x}_{j-n}.$$
This shows that
$$d(\lambda\hat{x}_j)=\lambda \check{x}_{j+1}+\lambda T^{n}\check{x}_{j-n}.$$
On the other hand, if $j\leq n$ since
$$\lambda T^{g_{\ell,j}^\Delta+jr_{\ell,j}^\Delta} f_\ell\left(T^{E(\hat{x}_{j}^\ell,\check{x}_{j-n}^\ell)}\check{x}^\ell_{j-n}\right)=\lambda T^{g_{\ell,j}^\Delta+jr_{\ell,j}^\Delta} T^{-g_{\ell,0}^\Delta}T^{E(\hat{x}_{j}^\ell,\check{x}_{j-n}^\ell)}\check{x}_{j-n}$$
      and  $$g_{\ell,j}^\Delta+jr_{\ell,j}^\Delta-g_{\ell,0}^\Delta+E(\hat{x}_{j}^\ell,\check{x}_{j-n}^\ell)=j,$$     
we deduce that
$$\lambda T^{g_{\ell,j}^\Delta+jr_{\ell,j}^\Delta} f_\ell\left(T^{E(\hat{x}_{j}^\ell,\check{x}_{j-n}^\ell)}\check{x}^\ell_{j-n}\right)=\lambda T^{j}\check{x}_{j-n}.$$
This shows that
$$d(\lambda\hat{x}_j)=\lambda \check{x}_{j+1}+\lambda T^{j}\check{x}_{j-n}$$
		as required.
	\end{itemize}
\end{proof}

Let us introduce a new presentation of $(C, d)$. Define a cochain complex $(\tilde{C}, \tilde{d})$ by
$$\tilde{C} = \Lambda_{>0}\check{x}_{-(n-1)}\oplus\cdots\oplus \Lambda_{>0}\check{x}_0\oplus \bigoplus_{j\in \N} \left( \Lambda_{>0}\check{x}_j \oplus \Lambda_{>0}\hat{x}_j \right),$$
where the differential $\tilde{d}$ satisfies $\tilde{d}(\lambda\check{x}_j) = 0$ for every $j \in \Z_{\geq -(n-1)}$ and $\lambda \in \Lambda_{>0}$, 
$$\tilde{d}(\lambda\hat{x}_j) = T^\Delta \lambda\check{x}_{j+1} + T^{n(1-\Delta)} \lambda\check{x}_{j-n}$$
for every $j \in \Z_{\geq n+1}$ and $\lambda \in \Lambda_{>0}$, and 
$$\tilde{d}(\lambda\hat{x}_j) = T^\Delta \lambda\check{x}_{j+1} + T^{j(1-\Delta)} \lambda\check{x}_{j-n}$$
for every $j \in \{1,\ldots,n\}$ and $\lambda \in \Lambda_{>0}$. 

Now, define a map $K \fc \tilde{C} \to C$ such that for every $\lambda \in \Lambda_{>0}$:
\begin{itemize}
    \item $K(\lambda\check{x}_j) = \lambda\check{x}_j$ for every $-(n-1) \leq j \leq 0$;
    \item $K(\lambda\check{x}_j) = T^{-j\Delta}\lambda\check{x}_j$ for every $j \in \N$;
    \item $K(\lambda\hat{x}_j) = T^{-j\Delta}\lambda\hat{x}_j$ for every $j \in \N$.
\end{itemize}

\begin{claim}\label{claim: new presentation of modified complex}
    The map $K$ is a cochain isomorphism.
\end{claim}

\begin{proof}
    It is straightforward to verify that $K$ is an isomorphism of modules. We check that $K$ is a cochain map.
    For every $j \in \Z_{\geq -(n-1)}$ and $\lambda \in \Lambda_{>0}$, we have:
    $$K \circ \tilde{d}(\lambda\check{x}_j) = K(0) = 0 = d(K(\lambda\check{x}_j)).$$
    
    Next, we verify the differential for the $\hat{x}_j$ generators. For every $j \in \Z_{\geq n+1}$ and $\lambda \in \Lambda_{>0}$, we have:
    \begin{align*}
        K \circ \tilde{d}(\lambda\hat{x}_j) &= K(T^\Delta \lambda\check{x}_{j+1} + T^{n(1-\Delta)} \lambda\check{x}_{j-n}) \\
        &= T^{-(j+1)\Delta} T^\Delta \lambda\check{x}_{j+1} + T^{-(j-n)\Delta} T^{n(1-\Delta)} \lambda\check{x}_{j-n} \\
        &= T^{-j\Delta} \lambda\check{x}_{j+1} + T^{-j\Delta + n\Delta + n - n\Delta} \lambda\check{x}_{j-n} \\
        &= T^{-j\Delta} \lambda\check{x}_{j+1} + T^{-j\Delta} T^n \lambda\check{x}_{j-n} \\
        &= T^{-j\Delta} (\lambda\check{x}_{j+1} + T^n \lambda\check{x}_{j-n}) \\
        &= d(T^{-j\Delta} \lambda\hat{x}_j) \\
        &= d \circ K(\lambda\hat{x}_j).
    \end{align*}
    
    Finally, for every $j \in \{1,\ldots,n\}$ and $\lambda \in \Lambda_{>0}$, note that $j-n \leq 0$, so $K$ acts as the identity on the $\check{x}_{j-n}$ terms. Thus, we have:
    \begin{align*}
        K \circ \tilde{d}(\lambda\hat{x}_j) &= K(T^\Delta \lambda\check{x}_{j+1} + T^{j(1-\Delta)} \lambda\check{x}_{j-n}) \\
        &= T^{-(j+1)\Delta} T^\Delta \lambda\check{x}_{j+1} + T^{j(1-\Delta)} \lambda\check{x}_{j-n} \\
        &= T^{-j\Delta} \lambda\check{x}_{j+1} + T^{-j\Delta + j} \lambda\check{x}_{j-n} \\
        &= T^{-j\Delta} (\lambda\check{x}_{j+1} + T^j \lambda\check{x}_{j-n}) \\
        &= d(T^{-j\Delta} \lambda\hat{x}_j) \\
        &= d \circ K(\lambda\hat{x}_j).
    \end{align*}
    Therefore, $d \circ K = K \circ \tilde{d}$, and we conclude that $K$ is a cochain isomorphism.
\end{proof}

\subsection{Proof of Theorem~\ref{thm: SH of complement of ball in CPn}}
We start with discussing some Quotients and Schauder bases.

Given $\alpha \geq \gamma > 0$ and $\beta > 0$, denote by $\tilde{D} \colon \widehat{\Lambda_{\geq 0}^\infty} \to \widehat{\Lambda_{\geq 0}^\infty}$ the unique endomorphism of $\widehat{\Lambda_{\geq 0}^\infty}$ that satisfies
$$\tilde{D}(e_1) = T^\gamma e_1 + T^\beta e_2$$
and
$$\tilde{D}(e_i) = T^\alpha e_i + T^\beta e_{i+1}$$
for every $i \geq 2$. Denote by $D$ the restriction of $\tilde{D}$ to $\widehat{\Lambda_{>0}^\infty}$. Thus $D$ satisfies
$$D(\lambda e_1) = T^\gamma \lambda e_1 + T^\beta \lambda e_2$$
and
$$D(\lambda e_i) = T^\alpha \lambda e_i + T^\beta \lambda e_{i+1}$$
for every $i \geq 2$ and $\lambda \in \Lambda_{>0}$. Let us compute the cokernel of $D$ and the natural quotient map $q\fc \widehat{\Lambda_{\geq 0}^\infty}\to\coker D$ on it.

\begin{prop}\label{prop: coker of D modified} $D$ is injective. In addition:
    \begin{itemize}
        \item If $\alpha < \beta$, then $\coker D$ is isomorphic to $\Lambda_{(0,\gamma]} \oplus \Lambda_{(0,\alpha]}^\infty$. Moreover, the natural quotient map $q \fc \widehat{\Lambda_{>0}^\infty} \to \coker D$ satisfies
        $$q(\lambda \epsilon_1) = \overline{\lambda} \epsilon_1$$
        for every $\lambda \in \Lambda_{>0}$, where $\overline{\lambda}$ is the image of $\lambda$ under the quotient $\Lambda_{>0} \to \Lambda_{(0,\gamma]} = \Lambda_{>0}/\Lambda_{>\gamma}$, and
        $$q(\lambda \epsilon_i) = \overline{\lambda} \epsilon_i$$
        for every $i \geq 2$ and $\lambda \in \Lambda_{>0}$, where $\overline{\lambda}$ is the image of $\lambda$ under the quotient $\Lambda_{>0} \to \Lambda_{(0,\alpha]} = \Lambda_{>0}/\Lambda_{>\alpha}$, and $(\epsilon_i)_{i=1}^\infty$ is the Schauder basis of $\widehat{\Lambda_{\geq 0}^\infty}$ given by $\epsilon_1 = e_1 + T^{\beta-\gamma} e_2$ and $\epsilon_i = e_i + T^{\beta-\alpha} e_{i+1}$ for every $i \geq 2$.

        \item If $\gamma \geq \beta$, then $\coker D$ is isomorphic to $\Lambda_{>0} \oplus \Lambda_{(0,\beta]}^\infty$. Moreover, the natural quotient map $q \fc \widehat{\Lambda_{>0}^\infty} \to \coker D$ satisfies
        $$q(\lambda \epsilon_1) = \lambda \epsilon_1$$
        for every $\lambda \in \Lambda_{>0}$, and
        $$q(\lambda \epsilon_i) = \overline{\lambda} \epsilon_i$$
        for every $i \geq 2$ and $\lambda \in \Lambda_{>0}$, where $\overline{\lambda}$ is the image of $\lambda$ under the quotient $\Lambda_{>0} \to \Lambda_{(0,\beta]} = \Lambda_{>0}/\Lambda_{>\beta}$, and $(\epsilon_i)_{i=1}^\infty$ is the Schauder basis of $\widehat{\Lambda_{\geq 0}^\infty}$ given by $\epsilon_1 = e_1$, $\epsilon_2 = T^{\gamma-\beta} e_1 + e_2$, and $\epsilon_i = T^{\alpha-\beta} e_{i-1} + e_i$ for every $i \geq 3$.

        \item If $\alpha \geq \beta>\gamma$, then $\coker D$ is isomorphic to $\Lambda_{(0,\gamma]} \oplus \Lambda_{>0} \oplus \Lambda_{(0,\beta]}^\infty$. Moreover, the natural quotient map $q \fc \widehat{\Lambda_{>0}^\infty} \to \coker D$ satisfies
        $$q(\lambda \epsilon_1) = \overline{\lambda} \epsilon_1$$
        for every $\lambda \in \Lambda_{>0}$, where $\overline{\lambda}$ is the image of $\lambda$ under the quotient $\Lambda_{>0} \to \Lambda_{(0,\gamma]} = \Lambda_{>0}/\Lambda_{>\gamma}$,
        $$q(\lambda \epsilon_2) = \lambda \epsilon_2$$
        for every $\lambda \in \Lambda_{>0}$, and
        $$q(\lambda \epsilon_i) = \overline{\lambda} \epsilon_i$$
        for every $i \geq 3$ and $\lambda \in \Lambda_{>0}$, where $\overline{\lambda}$ is the image of $\lambda$ under the quotient $\Lambda_{>0} \to \Lambda_{(0,\beta]} = \Lambda_{>0}/\Lambda_{>\beta}$, and $(\epsilon_i)_{i=1}^\infty$ is the Schauder basis of $\widehat{\Lambda_{\geq 0}^\infty}$ given by $\epsilon_1 = e_1 + T^{\beta-\gamma} e_2$, $\epsilon_2 = e_2$, and $\epsilon_i = T^{\alpha-\beta} e_{i-1} + e_i$ for every $i \geq 3$.
    \end{itemize}
\end{prop}

\begin{proof}[Proof of Proposition~\ref{prop: coker of D modified}]
\begin{itemize}
    \item Assume that $\alpha < \beta$. Since $\gamma \leq \alpha$, this implies $\gamma < \beta$. We use the Schauder basis $(\epsilon_i)_{i \in \N}$ given by $\epsilon_1 = e_1 + T^{\beta-\gamma} e_2$ and $\epsilon_i = e_i + T^{\beta-\alpha} e_{i+1}$ for every $i \geq 2$. For every $\lambda \in \Lambda_{>0}$, since $\gamma\leq\alpha<\beta$, we have
    $$D(\lambda e_1) = T^\gamma \lambda e_1 + T^\beta \lambda e_2 = T^\gamma \lambda (e_1 + T^{\beta-\gamma} e_2) = T^\gamma \lambda \epsilon_1.$$
    For every $i \geq 2$ and $\lambda \in \Lambda_{>0}$, we have
    $$D(\lambda e_i) = T^\alpha \lambda e_i + T^\beta \lambda e_{i+1} = T^\alpha \lambda (e_i + T^{\beta-\alpha} e_{i+1}) = T^\alpha \lambda \epsilon_i.$$
    Therefore,
    $$\im D = \widehat{\Lambda_{>\gamma}\epsilon_1 \oplus \bigoplus_{i=2}^\infty \Lambda_{>\alpha}\epsilon_i},$$
    and hence it follows from Lemma~\ref{lemma: quotient of completes modules 1} that
    $$\coker D = \widehat{\Lambda_{>0}^\infty} / \im D = \widehat{\bigoplus_{i=1}^\infty \Lambda_{>0}\epsilon_i} \Bigg/ \widehat{\Lambda_{>\gamma}\epsilon_1 \oplus \bigoplus_{i=2}^\infty \Lambda_{>\alpha}\epsilon_i} = \Lambda_{(0,\gamma]}\epsilon_1 \oplus \bigoplus_{i=2}^\infty \Lambda_{(0,\alpha]}\epsilon_i = \Lambda_{(0,\gamma]} \oplus \Lambda_{(0,\alpha]}^\infty.$$
    Moreover, the quotient map $q \fc \widehat{\bigoplus_{i=1}^\infty \Lambda_{>0}\epsilon_i} \to \Lambda_{(0,\gamma]}\epsilon_1 \oplus \bigoplus_{i=2}^\infty \Lambda_{(0,\alpha]}\epsilon_i$ satisfies $q(\lambda \epsilon_1) = \overline{\lambda} \epsilon_1$, where $\overline{\lambda}$ is the image of $\lambda$ under the quotient $\Lambda_{>0} \to \Lambda_{(0,\gamma]} = \Lambda_{>0}/\Lambda_{>\gamma}$, and $q(\lambda \epsilon_i) = \overline{\lambda} \epsilon_i$ for every $i \geq 2$, where $\overline{\lambda}$ is the image of $\lambda$ under $\Lambda_{>0} \to \Lambda_{(0,\alpha]} = \Lambda_{>0}/\Lambda_{>\alpha}$.

    \item Assume that $\gamma \geq \beta$. We use the Schauder basis $(\epsilon_i)_{i \in \N}$ given by $\epsilon_1 = e_1$, $\epsilon_2 = T^{\gamma-\beta} e_1 + e_2$, and $\epsilon_i = T^{\alpha-\beta} e_{i-1} + e_i$ for every $i \geq 3$. For every $\lambda \in \Lambda_{>0}$, we have
    $$D(\lambda e_1) = T^\gamma \lambda e_1 + T^\beta \lambda e_2 = T^\beta \lambda (T^{\gamma-\beta} e_1 + e_2) = T^\beta \lambda \epsilon_2.$$
    For every $i \geq 2$ and $\lambda \in \Lambda_{>0}$, since $\alpha\geq\gamma\geq \beta$, we have
    $$D(\lambda e_i) = T^\alpha \lambda e_i + T^\beta \lambda e_{i+1} = T^\beta \lambda (T^{\alpha-\beta} e_i + e_{i+1}) = T^\beta \lambda \epsilon_{i+1}.$$
    Therefore,
    $$\im D = \widehat{\bigoplus_{i=2}^\infty \Lambda_{>\beta}\epsilon_i},$$
    and hence it follows from Lemma~\ref{lemma: quotient of completes modules 2} that
    $$\coker D = \widehat{\Lambda_{>0}^\infty} / \im D = \widehat{\bigoplus_{i=1}^\infty \Lambda_{>0}\epsilon_i} \Bigg/ \widehat{\bigoplus_{i=2}^\infty \Lambda_{>\beta}\epsilon_i} = \Lambda_{>0}\epsilon_1 \oplus \bigoplus_{i=2}^\infty \Lambda_{(0,\beta]}\epsilon_i = \Lambda_{>0} \oplus \Lambda_{(0,\beta]}^\infty.$$
    Moreover, the quotient map $q$ satisfies $q(\lambda \epsilon_1) = \lambda \epsilon_1$ for every $\lambda \in \Lambda_{>0}$, and $q(\lambda \epsilon_i) = \overline{\lambda} \epsilon_i$ for every $i \geq 2$, where $\overline{\lambda}$ is the image of $\lambda$ under the quotient $\Lambda_{>0} \to \Lambda_{(0,\beta]} = \Lambda_{>0}/\Lambda_{>\beta}$.

    \item Assume that $\alpha \geq \beta>\gamma$. We use the Schauder basis $(\epsilon_i)_{i \in \N}$ given by $\epsilon_1 = e_1 + T^{\beta-\gamma} e_2$, $\epsilon_2 = e_2$, and $\epsilon_i = T^{\alpha-\beta} e_{i-1} + e_i$ for every $i \geq 3$. For every $\lambda \in \Lambda_{>0}$, we have
    $$D(\lambda e_1) = T^\gamma \lambda e_1 + T^\beta \lambda e_2 = T^\gamma \lambda (e_1 + T^{\beta-\gamma} e_2) = T^\gamma \lambda \epsilon_1.$$
    For every $i \geq 2$ and $\lambda \in \Lambda_{>0}$, we have
    $$D(\lambda e_i) = T^\alpha \lambda e_i + T^\beta \lambda e_{i+1} = T^\beta \lambda (T^{\alpha-\beta} e_i + e_{i+1}) = T^\beta \lambda \epsilon_{i+1}.$$
    Therefore,
    $$\im D = \widehat{\Lambda_{>\gamma}\epsilon_1 \oplus \bigoplus_{i=3}^\infty \Lambda_{>\beta}\epsilon_i}.$$
    By Lemma~\ref{lemma: quotient of completes modules 2}, we obtain
    \begin{align*}
     \coker D &= \widehat{\Lambda_{>0}^\infty} / \im D  \\
     &= \widehat{\bigoplus_{i=1}^\infty \Lambda_{>0}\epsilon_i} \Bigg/ \widehat{\left(\Lambda_{>\gamma}\epsilon_1 \oplus 0 \cdot \epsilon_2 \oplus \bigoplus_{i=3}^\infty \Lambda_{>\beta}\epsilon_i\right)}\\
     &= \Lambda_{(0,\gamma]}\epsilon_1 \oplus \Lambda_{>0}\epsilon_2 \oplus \bigoplus_{i=3}^\infty \Lambda_{(0,\beta]}\epsilon_i\\
     &= \Lambda_{(0,\gamma]} \oplus \Lambda_{>0} \oplus \Lambda_{(0,\beta]}^\infty.
    \end{align*}
    Moreover, the quotient map $q$ satisfies $q(\lambda \epsilon_1) = \overline{\lambda} \epsilon_1$ where $\overline{\lambda}$ is the image of $\lambda$ under $\Lambda_{>0} \to \Lambda_{(0,\gamma]} = \Lambda_{>0}/\Lambda_{>\gamma}$, $q(\lambda \epsilon_2) = \lambda \epsilon_2$ for every $\lambda \in \Lambda_{>0}$, and $q(\lambda \epsilon_i) = \overline{\lambda} \epsilon_i$ for every $i \geq 3$, where $\overline{\lambda}$ is the image of $\lambda$ under the quotient $\Lambda_{>0} \to \Lambda_{(0,\beta]} = \Lambda_{>0}/\Lambda_{>\beta}$.
\end{itemize}
\end{proof}

We now apply Proposition~\ref{prop: coker of D modified} to a specific setting parameterized by $n \in \N$, $\Delta \in (0,1)$, and an integer $0 \leq i \leq n-1$. Let us set
$$\alpha = n(1-\Delta), \quad \gamma = (n-i)(1-\Delta), \quad \text{and} \quad \beta = \Delta.$$
Since $n - i \geq 1$ and $\Delta \in (0,1)$, we immediately have $\beta > 0$ and $\alpha \geq \gamma > 0$, satisfying the initial hypotheses. 

The inequalities comparing $\alpha$ and $\gamma$ with $\beta$ translate directly into thresholds for $\Delta$:
\begin{align*}
    \alpha < \beta &\iff n(1-\Delta) < \Delta \iff n < (n+1)\Delta \iff \Delta > \frac{n}{n+1}, \\
    \gamma < \beta &\iff (n-i)(1-\Delta) < \Delta \iff (n-i) < (n-i+1)\Delta \iff \Delta > \frac{n-i}{n-i+1}.
\end{align*}
Because $0\leq i \leq n-1$, we have $1\leq n-i\leq n$ and hence  $\frac{n-i}{n-i+1} \leq \frac{n}{n+1}$. This naturally partitions the interval $(0,1)$ into three regimes, yielding the following corollary.

\begin{coroll}\label{cor: coker of D parameterized}
    Let $D$ be the operator defined above with parameters $\alpha = n(1-\Delta)$, $\gamma = (n-i)(1-\Delta)$, and $\beta = \Delta$. The cokernel of $D$ depends on the value of $\Delta$ as follows:
    \begin{itemize}
        \item If $\Delta > \frac{n}{n+1}$, then $\alpha < \beta$, and $\coker D$ is isomorphic to 
        $$ \Lambda_{(0, (n-i)(1-\Delta)]} \oplus \Lambda_{(0, n(1-\Delta)]}^\infty. $$
        
        \item If $\Delta \leq \frac{n-i}{n-i+1}$, then $\gamma \geq \beta$ and $\alpha \geq \beta$, and $\coker D$ is isomorphic to 
        $$ \Lambda_{>0} \oplus \Lambda_{(0, \Delta]}^\infty. $$
        
        \item If $\frac{n-i}{n-i+1} < \Delta \leq \frac{n}{n+1}$, then $\gamma < \beta$ and $\alpha \geq \beta$, and $\coker D$ is isomorphic to 
        $$ \Lambda_{(0, (n-i)(1-\Delta)]} \oplus \Lambda_{>0} \oplus \Lambda_{(0, \Delta]}^\infty. $$
    \end{itemize}
\end{coroll}

\begin{proof}[Proof of Theorem~\ref{thm: SH of complement of ball in CPn}]
 The differential on each odd-degree component of the completed complex is identified with the operator $D$, while the differential on even-degree components is zero.
 By Proposition~\ref{prop: coker of D modified}, $D$ is injective, so the odd-degree cohomology equals $\ker D=0$. Consequently, the relative symplectic cohomology vanishes in all odd degrees: $SH^{(2k+1)}_{\CP^n}(\CP^n \setminus \Int B^{}_{\Delta}; \Lambda_{\geq 0}) = 0$.

Furthermore, computations of topological energy and computing direct limits show that the relative symplectic cohomology in degree $0$ is isomorphic to the relative symplectic cohomology in degree $2n$. 

For even degrees $* = 2i$ with $0 \leq i \leq n-1$, the Floer differential computing the relative symplectic cohomology is identified with the operator $D \colon \widehat{\Lambda_{>0}^\infty} \to \widehat{\Lambda_{>0}^\infty}$. The parameters of this operator depend on the degree index $i$ and the capacity $\Delta$ via the relations:
$$\alpha = n(1-\Delta), \quad \gamma = (n-i)(1-\Delta), \quad \text{and} \quad \beta = \Delta.$$
We have $SH^{(2i)}_{\CP^n}(\CP^n \setminus \Int B^{}_{\Delta}; \Lambda_{\geq 0}) \cong \coker D$. The structure of this cokernel is governed by Corollary~\ref{cor: coker of D parameterized}. We analyze the capacity intervals as follows:

\textbf{Case 1: $\Delta \leq \frac{1}{2}$.} \\
For any index $0 \leq i \leq n-1$, we have $\frac{n-i}{n-i+1} \geq \frac{1}{2} \geq \Delta$. Thus, $\Delta \leq \frac{n-i}{n-i+1}$ holds uniformly. Applying the corresponding regime of Corollary~\ref{cor: coker of D parameterized}, the cokernel is $\Lambda_{>0} \oplus \Lambda_{(0, \Delta]}^\infty$ for all even degrees $2i$ with $1 \leq i \leq n$. By the identification of degrees $0$ and $2n$ via the direct limit computation, the module at degree $0$ is identical. Thus, $SH^*$ is given by $\Lambda_{>0} \oplus \Lambda_{(0, \Delta]}^\infty$ for all even degrees.

\textbf{Case 2: $\frac{j}{j+1} < \Delta \leq \frac{j+1}{j+2}$ for an integer $1 \leq j \leq n-1$.} \\
We have two cases to consider.
\begin{itemize}
    \item For $* = 2i$ with $i \in \{n-j, \ldots, n-1\}$, we have $\frac{n-i}{n-i+1} \leq \frac{j}{j+1} < \Delta$. Furthermore, since $j \leq n-1$, we have $\Delta \leq \frac{j+1}{j+2} \leq \frac{n}{n+1}$. Combining these yields $\frac{n-i}{n-i+1} < \Delta \leq \frac{n}{n+1}$. By Corollary~\ref{cor: coker of D parameterized}, this yields the module $\Lambda_{(0, (n-i)(1-\Delta)]} \oplus \Lambda_{>0} \oplus \Lambda_{(0, \Delta]}^\infty$.
    
    \item For $* = 2i$ with $i \in \{0, \ldots, n-j-1\}$, we have $\frac{n-i}{n-i+1} \geq \frac{j+1}{j+2} \geq \Delta$. Thus, $\Delta \leq \frac{n-i}{n-i+1}$. By Corollary~\ref{cor: coker of D parameterized}, this yields the module $\Lambda_{>0} \oplus \Lambda_{(0, \Delta]}^\infty$.
    
\end{itemize}

\textbf{Case 3: $\Delta > \frac{n}{n+1}$.} \\
In this regime, for any $0\leq i \leq n-1$, we have $\frac{n-i}{n-i+1} \leq \frac{n}{n+1} < \Delta$. Therefore by Corollary~\ref{cor: coker of D parameterized}, this yields the module $\Lambda_{(0, (n-i)(1-\Delta)]} \oplus \Lambda_{(0, n(1-\Delta)]}^\infty$. This completes the proof.
\end{proof}

\subsection{Proof of Proposition~\ref{prop: res from CPn to the complement a ball}}\label{ss: res maps from CPn to complement a ball}

The proof of Proposition~\ref{prop: res from CPn to the complement a ball} is similar to and easier than that of Proposition~\ref{prop: res from CPn to a ball}. Let $0<\Delta<1$. Consider the corresponding acceleration data $(G_\ell)_{\ell\geq0}$ for the complement of the ball $B^{}_{\Delta}=\mu^{-1}([0,\Delta])$, given by $G_\ell(z)=g(\Delta,\ell,\mu(z))$ for every $z\in \CP^n$ and $\ell\in \Z_{\geq0}$. For every $\ell\in \Z_{\geq 0}$, we have $G_0\leq G_\ell$ and $G_0\leq 0$. Thus, since 
$$\lim\limits_{\ell\to \infty}G_\ell(z)=\left\{\begin{array}{ll}
   0, & z\notin B,\\
    +\infty, &  z\in B,
\end{array}\right.$$ 
for every $z\in \CP^n$, set $\varepsilon_0=1$, and there exists a non-increasing sequence $(\varepsilon_\ell)_{\ell\geq1}$ such that $\lim\limits_{\ell\to \infty}\varepsilon_\ell=0$ and for every $\ell\in \Z_{\geq1}$ we have $\varepsilon_\ell\cdot G_0\leq G_\ell$. For every $\ell\in \Z_{\geq1}$, we denote $G'_\ell=\varepsilon_\ell\cdot G_0$. Note that $(G'_\ell)_{\ell\geq1}$ is an acceleration datum for $\CP^n$ consisting of $C^2$-small Morse--Bott functions. For every $\ell\in \Z_{\geq0}$, the Floer complex $CF(G'_\ell)$ of $G'_\ell$ has $n+1$ generators, denoted $\check{z}_{-(n-1)}^\ell,\ldots,\check{z}_1^\ell$. The differential is $0$, and the continuation map $\Phi'_\ell\fc CF(G'_\ell)\to CF(G'_{\ell+1})$ satisfies 
\begin{equation}\label{eq: cont for G of CPn}
    \Phi'_\ell(\check{z}_i^\ell)=T^{G'_{\ell+1}(\check{z}_i^{\ell+1})-G'_{\ell}(\check{z}_i^{\ell})}\check{z}_i^{\ell+1}
\end{equation} for every $-(n-1)\leq i\leq 1$, see \cite[Page 599]{Varolgunes_2021_MV_and_relSH} for the Morse case.

For every $\ell\in \Z_{\geq0}$ denote by $\check{x}_{-(n-1)}^\ell,\ldots,\check{x}_{\ell+1}^\ell,\hat{x}_1^\ell,\ldots,\hat{x}_\ell^\ell$ the generators for $CF(G_\ell)$. By Proposition~\ref{prop: omission of signs for CF(G)}, for every $\ell\in \Z_{\geq0}$, the bases $\check{x}_{-(n-1)}^\ell,\ldots,\check{x}_{\ell+1}^\ell,\hat{x}_1^\ell,\ldots,\hat{x}_\ell^\ell$ for $CF(G_\ell)$ can be chosen such that the differentials and continuation maps in the $1$-ray
$$CF(G_0)\xrightarrow{\Phi_0}CF(G_1)\xrightarrow{\Phi_1}CF(G_2)\to\cdots$$

are represented by matrices with entries in $\{0, 1\}$, where we ignore the powers of the formal variable $T$ of the Novikov ring. We  refer to such bases as \textbf{normalized bases}.

Let $\chi\fc \R\to \R$ be a smooth non-decreasing function, satisfying $\chi(s)=0$ for every $s<0$ and $\chi(s)=1$ for every $s>1$.  Then for every $\ell\in \Z_{\geq0}$ one can define a monotone homotopy $K_\ell\fc \R\times \CP^n\to \R$ from $G'_\ell$ to $G_\ell$ by $K_\ell(s,z)
=(1-\chi(s))\varepsilon_\ell G_0(z)
+\chi(s)G_\ell(z)$ for every $(s,z)\in \R\times \CP^n$. Thus for every $\ell\in \Z_{\geq0}$ the homotopy $K_\ell$ defines a continuation map $\Psi_\ell\fc CF(G'_\ell)\to CF(G_\ell)$.
Similarly to Theorem~\ref{thm: CF(G_l) + continuation maps}, to describe the continuation maps with coefficients in the Novikov ring, in the sense of Section~\ref{sss: weighted_CF}, we must compute the energy of the Floer and continuation flowlines with cascades that appear in the differential and continuation maps.

\begin{prop}\label{prop: continuation maps for restrictions for CP^n->CP^n-B}
    For every $\ell\in \Z_{\geq0}$ there exist signs
    $$\check{D}_{-(n-1)},\ldots,\check{D}_{1}\in \{-1,1\}$$
    such that the continuation map $\Psi_\ell\fc CF(G'_\ell)\to CF(G_\ell)$ satisfies
    \begin{itemize}
    \item $\Phi_\ell\circ \Psi_\ell=\Psi_{\ell+1}\circ\Phi'_\ell$ where $\Phi_\ell\fc CF(G_\ell)\to CF(G_{\ell+1})$ and $\Phi'_\ell\fc CF(G'_\ell)\to CF(G'_{\ell+1})$ are the continuation maps.
   
        \item For every $-(n-1)\leq i\leq 1$ we have
        $$\Psi(\check{z}_i^\ell)=\check{D}_i T^{E(\check{z}_i^\ell,\check{x}_i^\ell)}\check{x}_i^\ell,$$
        where
        $$E(\check{z}_i^\ell,\check{x}_i^\ell)=\left\{\begin{array}{ll}
        g_{\ell,0}^\Delta-\varepsilon_\ell g_{0,0}^\Delta,& i\leq0,\\
            g(\Delta,\ell,0)-\varepsilon_\ell\cdot g(\Delta,0,0), & 1= i,\,\ell=0, \\
            g_{\ell,1}^\Delta+r_{\ell,1}^\Delta-\varepsilon_\ell\cdot g(\Delta,0,0), & 1=i\leq \ell. 
        \end{array}\right.$$

    \end{itemize}
\end{prop}
\begin{proof}
   Let $\ell \in \Z_{\geq 0}$. Since $G'_\ell$ is obtained from $G_0$ by multiplication by a small constant, the fact that $\Phi_\ell\circ \Psi_\ell=\Psi_{\ell+1}\circ\Phi'_\ell$ can be proven similar to Proposition~\ref{prop: commuting of cont maps}. From the same reason, Theorem~\ref{thm: continuation maps for G_l, over Z} implies that there exist $\check D_{-(n-1)},\ldots,\check D_1\in \{-1,1\}$ and $E(\check{z}_{-(n-1)}^\ell,\check{x}_{-(n-1)}^\ell),\ldots,E(\check{z}_1^\ell,\check{x}_1^\ell)>0$ such that
 $$\Psi(\check{z}_i^\ell)=\check{D}_i T^{E(\check{z}_i^\ell,\check{x}_i^\ell)}\check{x}_i^\ell,$$
 for every $-(n-1)\leq i\leq 1$.

   The rest of the proof is similar to that of Theorem~\ref{thm: CF(G_l) + continuation maps}. Since we are interested in continuation flowlines with cascades, the dimension $\dim_u \wh\cM_m(p, q)$ equals $0$ for every $p, q$, and $u$. Thus, the topological energy of such a flowline with cascaded is given by
$$ E_{top}(u) = \cA_{G_{\ell}}(q, \hat{q}) - \cA_{G'_\ell}(p, \hat{p}) + \frac{1}{2(n+1)} \left( \mu_{FMB}^\tau(p) - \mu_{FMB}^\tau(q) \right). $$

As in the proof of Theorem~\ref{thm: CF(G_l) + continuation maps}, for constant orbits of $G_\ell$ and $G'_\ell$, choose the cappings to be constant, thus their action equals the value of the Hamiltonian itself, that is
$$ \cA_{G_\ell}(\check{x}_{\ell+1}^\ell) = g(\Delta, \ell, 0), \qquad \text{and} \qquad \cA_{G_\ell}(\check{x}_{-(n-1)}^\ell) = \cdots = \cA_{G_\ell}(\check{x}_{0}^\ell) = g(\Delta, \ell, 1)=g_{\ell,0}^\Delta, $$
and 

$$ \cA_{G'_\ell}(\check{z}_{1}^\ell) = \varepsilon_\ell\cdot g(\Delta, 0, 0), \qquad \text{and} \qquad \cA_{G'_\ell}(\check{z}_{-(n-1)}^\ell) = \cdots = \cA_{G'_\ell}(\check{z}_{0}^\ell) = \varepsilon_\ell\cdot g(\Delta, 0, 1)=\varepsilon_\ell\cdot g_{0,0}^\Delta. $$
For the orbits $\hat{x}_1^\ell, \check{x}_1^\ell, \ldots, \hat{x}_\ell^\ell, \check{x}_\ell^\ell$, we choose cappings contained in $\CP^n \setminus D_\infty$, which is an exact symplectic manifold symplectomorphic to $\Int B(1)$. Let $\lambda = \frac{1}{2} \sum_{i=1}^n (x_i dy_i - y_i dx_i)$ be a primitive of $\omega_0$ on $\Int B(1)$. Then, for every $1 \leq i \leq \ell$, we find that
$$ \cA_{G_\ell}(\hat{x}_i^\ell) = \cA_{G_\ell}(\check{x}_i^\ell) = \int_{S^1} G_\ell \circ \check{x}_i^\ell(t) \, dt + \int_{S^1} (\check{x}_i^\ell)^* \lambda = g_{\ell,i}^\Delta + i r_{\ell,i}^\Delta. $$
Additionally, the Floer--Morse--Bott indices for these $1$-periodic orbits with respect to the trivialization $\tau$ induced by these cappings are as follows:
\begin{itemize}
    \item $\mu_{FMB}^\tau(\check{x}_{\ell+1}^\ell) = 2n(\ell+1)$.
    \item For every $1 \leq i \leq \ell$, $\mu_{FMB}^\tau(\check{x}_i^\ell) = 2ni$.
    \item For every $1 \leq i \leq \ell$, $\mu_{FMB}^\tau(\hat{x}_i^\ell) = 2ni + 2n - 1$.
    \item For every $-(n-1) \leq i \leq 0$, $\mu_{FMB}^\tau(\check{x}_i^\ell)=\mu_{FMB}^\tau(\check{z}_i^\ell) = -2i$.
    \item $\mu_{FMB}^\tau(\check{z}_{1}^\ell) = 2n$.
\end{itemize}

Now, let $-(n-1) \leq i \leq 1$.
\begin{itemize}
\item If $-(n-1)\leq i\leq 0$, then $E(\check{z}_i^\ell, \check{x}_i^{\ell}) = \cA_{G_{\ell}}(\check{x}_i^{\ell}) - \cA_{G'_\ell}(\check{z}_i^\ell) = g_{\ell,0}^\Delta-\varepsilon_\ell g_{0,0}^\Delta.$

    \item Assume that $i=1$.
    \begin{itemize}
        \item If $\ell\geq1$ then
        $$\mu_{FMB}^\tau(\check{z}_1^\ell)-\mu_{FMB}^\tau(\check{x}_1^\ell)=2n-2n=0,$$
        and hence
       $$E(\check{z}_1^\ell, \check{x}_1^{\ell})  = \cA_{G_{\ell}}(\check{x}_1^{\ell}) - \cA_{G'_\ell}(\check{z}_1^\ell)= g_{\ell,1}^\Delta+r_{\ell,1}^\Delta-\varepsilon_\ell\cdot g(\Delta,0,0)$$

        \item If $\ell=0$, then 
       $$E(\check{z}_1^\ell, \check{x}_1^{\ell})  = \cA_{G_{\ell}}(\check{x}_1^{\ell}) - \cA_{G'_\ell}(\check{z}_1^\ell)= g(\Delta,\ell,0)-\varepsilon_\ell\cdot g(\Delta,0,0).$$

    \end{itemize}
\end{itemize}
This completes the proof.

\end{proof}

The following result is an analogue of Proposition~\ref{prop: omission of signs for CF->CF'} and demonstrates that the signs from Theorem~\ref{prop: continuation maps for restrictions for CP^n->B} can be normalized to $1$ via a suitable change of bases.

\begin{prop}\label{prop: omission of signs for CF'->CF, complement of ball}
     For every $\ell\in \Z_{\geq0}$, there is a normalized basis, $\check{x}_{-(n-1)}^\ell,\ldots,\check{x}_{\ell+1}^\ell,\hat{x}_1^\ell,\ldots,\hat{x}_\ell^\ell $ 
for $CF(G_\ell)$ such that the continuation map $\Psi_\ell\fc CF(G'_\ell)\to CF(G_\ell)$ satisfies
        $$\Psi(\check{z}_i^\ell)=T^{E(\check{z}_i^\ell,\check{x}_i^\ell)}\check{x}_i^\ell,$$
        where
        $$E(\check{z}_i^\ell,\check{x}_i^\ell)=\left\{\begin{array}{ll}
        g_{\ell,0}^\Delta-\varepsilon_\ell g_{0,0}^\Delta,& i\leq0,\\
            g(\Delta,\ell,0)-\varepsilon_\ell\cdot g(\Delta,0,0), & 1= i,\,\ell=0, \\
            g_{\ell,1}^\Delta+r_{\ell,1}^\Delta-\varepsilon_\ell\cdot g(\Delta,0,0), & 1=i\leq \ell, 
        \end{array}\right.$$
        
       for every $-(n-1)\leq i\leq 1$.
\end{prop}

\begin{proof} Let $\ell\in \Z_{\geq0}$. The desired result is an immediate corollary of Proposition~\ref{prop: omission of signs for CF(G)} applied to $CF(G_\ell;\Z)$, combined with the fact that the differential on $CF(G'_\ell;\Z)$ is zero and that $CF(G'_\ell;\Z)$ has at most one generator in each degree. Thus, we are able to change the signs of the generators of $CF(G'_\ell;\Z)$ to deduce the desired result.
\end{proof}

A direct using of the combination of Lemma~\ref{lemma: drct lim of 1-ray} and Equation~\eqref{eq: cont for G of CPn} implies the direct limit of the $1$-ray
$$CF(G'_0) \to CF(G'_1) \to CF(G'_2) \to \cdots$$
is isomorphic to the cochain complex $(C',d')$ given by $$C'=\bigoplus_{j=-(n-1)}^1\Lambda_{>0}\check{z}_j$$ where the differential $d'$ is zero, and for every $\ell\in \Z_{\geq0}$ the canonical map $f'_\ell\fc CF(G'_\ell)\to C'$ satisfies $f'_\ell(\check{z}_j^\ell)=T^{-G'_\ell(\check{z}_j^\ell)}\check{z}_j$ for every $-(n-1)\leq j\leq 1$.

Consider the cochain complex $(C, d)$ defined as follows: the module $C$ is given by
$$C = \Lambda_{>0}\check{x}_{-(n-1)}\oplus\cdots\oplus \Lambda_{>0}\check{x}_0\oplus \bigoplus_{j\in \N} \left( \Lambda_{>-j\Delta}\check{x}_j \oplus \Lambda_{>-j\Delta}\hat{x}_j \right),$$

and the differential $d$ satisfies $d(\lambda\check{x}_j) = 0$ for every $\lambda\in \Lambda_{>0}$ when $-(n-1)\leq j\leq 0$, and $\lambda \in \Lambda_{>-j\Delta}$ when $j\geq 1$, 
$$d(\lambda\hat{x}_j) = \lambda\check{x}_{j+1} + T^n\lambda\check{x}_{j-n}$$
for every $j \in \Z_{\geq n+1}$ and $\lambda \in \Lambda_{>-j\Delta}$, and also 
$$d(\lambda\hat{x}_j) = \lambda\check{x}_{j+1} + T^j\lambda\check{x}_{j-n}$$
for every $j \in \{1\ldots,n\}$ and $\lambda \in \Lambda_{>-j\Delta}$.

By Theorem~\ref{thm: drct lim of CF(G_l)}, the direct limits of the $1$-ray
$$CF(G_0) \to CF(G_1) \to CF(G_2) \to \cdots$$
is isomorphic to $(C, d)$.

Let $\Psi \fc C' \to C$ denote the direct limit of the continuation maps $(\Psi_\ell \fc CF(G'_\ell) \to CF(G_\ell))_{\ell\geq 0}$. Note that $\Psi$ is well-defined because Proposition~\ref{prop: continuation maps for restrictions for CP^n->CP^n-B} asserts that $\Psi_{\ell+1} \circ \Phi'_\ell = \Phi_\ell \circ \Psi_\ell$, which ensures that the maps $(\Psi_\ell)_{\ell\geq 0}$ define a morphism of directed systems. The following proposition describes this map.

\begin{prop}
    The map 
    $$\Psi \fc \bigoplus_{j=-(n-1)}^1\Lambda_{>0}\check{z}_j\to\Lambda_{>0}\check{x}_{-(n-1)}\oplus\cdots\oplus \Lambda_{>0}\check{x}_0\oplus \bigoplus_{j\in \N} \left( \Lambda_{>-j\Delta}\check{x}_j \oplus \Lambda_{>-j\Delta}\hat{x}_j \right)$$
    satisfies $\Psi(\lambda \check{z}_j) = \lambda \check{x}_j$ for every $-(n-1)\leq j\leq 1$ and $\lambda\in \Lambda_{>0}$.
\end{prop}

\begin{proof}
By the universal property of the direct limit, for every $\ell \in \Z_{\geq 0}$, the following diagram commutes:
$$
\xymatrix@R=2pc@C=3pc{
CF(G'_\ell) \ar[r]^{f'_\ell} \ar[d]_{\Psi_\ell} & C' \ar[d]^{\Psi} \\
CF(G_\ell) \ar[r]_{f_\ell} & C
}
$$
where $f_\ell$ and $f'_\ell$ are the canonical maps of the direct limits $C$ and $C'$, respectively. 

Let $-(n-1)\leq j\leq 1$ and $\lambda \in \Lambda_{>0}$. From the properties of direct limits, for any $\ell \geq 1$, the generator $\lambda \check{z}_j \in C'$ can be represented as:
$$\lambda \check{z}_j = f'_\ell \left( \lambda T^{ G'_\ell(\check{z}_j^\ell)} \check{x}_j^\ell \right).$$
Using the commutativity of the diagram, we have:
$$\Psi(\lambda \check{z}_j) = \Psi \circ f'_\ell \left( \lambda T^{ G'_\ell(\check{z}_j^\ell)} \check{z}_j^\ell \right) = f_\ell \circ \Psi_\ell \left( \lambda T^{ G'_\ell(\check{z}_j^\ell)} \check{z}_j^\ell \right).$$
By Proposition~\ref{prop: omission of signs for CF'->CF, complement of ball}, we have $\Psi_\ell(\check{z}_j^\ell) = T^{E(\check{z}_j^\ell, \check{x}_j^\ell)} \check{x}_j^\ell$. Thus:
$$\Psi(\lambda \check{z}_j) = \lambda T^{G'_\ell(\check{z}_j^\ell)} f_\ell \left( T^{E(\check{z}_j^\ell, \check{x}_j^\ell)} \check{x}_j^\ell \right).$$
Applying the formula for $f_\ell \fc CF(G_\ell) \to C$ from Theorem~\ref{thm: drct lim of CF(G_l)}, we obtain:
$$\Psi(\lambda \check{z}_j) =\left\{\begin{array}{ll}
   \lambda T^{G'_\ell(\check{z}_j^\ell)+E(\check{z}_j^\ell, \check{x}_j^\ell)-g_{\ell,0}^\Delta} \check{x}_j^\ell ,  &  j\leq0,\\
     \lambda T^{G'_\ell(\check{z}_j^\ell)+E(\check{z}_j^\ell, \check{x}_j^\ell)-(g_{\ell,1}^\Delta+r_{\ell,1}^\Delta)} \check{x}_j^\ell, & j=1.
\end{array}\right.$$
If $j=1$, then $G'_\ell(\check{z}_1^\ell)=\varepsilon_\ell\cdot  g(\Delta,0,0)$ and $E(\check{z}_0^\ell, \check{x}_0^\ell)=g_{\ell,1}^\Delta+r_{\ell,1}^\Delta-\varepsilon_\ell\cdot g(\Delta,0,0)$ which imply that
    \begin{align*}
        G'_\ell(\check{z}_j^\ell)+E(\check{z}_j^\ell, \check{x}_j^\ell)-(g_{\ell,1}^\Delta+r_{\ell,1}^\Delta)&=\varepsilon_\ell\cdot g(\Delta,0,0)+g_{\ell,1}^\Delta+r_{\ell,1}^\Delta\\
        &-\varepsilon_\ell\cdot g(\Delta,0,0)-(g_{\ell,1}^\Delta+r_{\ell,1}^\Delta)\\
        &=0
    \end{align*}
   Otherwise, $j\leq 0$, and hence $G'_\ell(\check{z}_j^\ell)=\varepsilon_\ell\cdot  g_{0,0}^\Delta$ and $E(\check{z}_j^\ell, \check{x}_j^\ell)=g_{\ell,0}^\Delta-\varepsilon_\ell\cdot g_{0,0}^\Delta$ which imply that
$$G'_\ell(\check{z}_j^\ell)+E(\check{z}_j^\ell, \check{x}_j^\ell)-g_{\ell,0}^\Delta=\varepsilon_\ell\cdot  g_{0,0}^\Delta+g_{\ell,0}^\Delta-\varepsilon_\ell\cdot g_{0,0}^\Delta-g_{\ell,0}^\Delta=0$$
 Therefore, $\Psi(\lambda \check{z}_j) = \lambda \check{x}_j$. This completes the proof.
\end{proof}

As in the proof of Theorem~\ref{thm: drct lim of CF(G_l)}, define cochain complex $(\tilde{C},\tilde{d})$  over $\Lambda_{\geq 0}$ as follows:
$$\tilde{C}=\Lambda_{>0}\check{x}_{-(n-1)}\oplus\cdots\oplus \Lambda_{>0}\check{x}_0\oplus \bigoplus_{i=1}^\infty \left(\Lambda_{>0} \check{x}_i\oplus \Lambda_{>0} \hat{x}_i\right),$$
where the differential $\tilde{d}$ satisfies $\tilde{d}(\lambda\check{x}_j) = 0$ for every $j \in \Z_{\geq -(n-1)}$ and $\lambda \in \Lambda_{>0}$, 
$$\tilde{d}(\lambda\hat{x}_j) = T^\Delta \lambda\check{x}_{j+1} + T^{n(1-\Delta)} \lambda\check{x}_{j-n}$$
for every $j \in \Z_{\geq n+1}$ and $\lambda \in \Lambda_{>0}$, and 
$$\tilde{d}(\lambda\hat{x}_j) = T^\Delta \lambda\check{x}_{j+1} + T^{j(1-\Delta)} \lambda\check{x}_{j-n}$$
for every $j \in \{1,\ldots,n\}$ and $\lambda \in \Lambda_{>0}$.

By Claim~\ref{claim: new presentation of modified complex}, the cochain complex $(\tilde{C},\tilde{d})$ is isomorphic to $(C,d)$ via the isomorphism $K\fc \tilde{C}\to C$ given by 
\begin{itemize}
    \item $K(\lambda\check{x}_j) = \lambda\check{x}_j$ for every $-(n-1) \leq j \leq 0$ and $\lambda \in \Lambda_{>0}$;
    \item $K(\lambda\check{x}_j) = T^{-j\Delta}\lambda\check{x}_j$ for every $j \in \N$ and $\lambda \in \Lambda_{>0}$;
    \item $K(\lambda\hat{x}_j) = T^{-j\Delta}\lambda\hat{x}_j$ for every $j \in \N$ and $\lambda \in \Lambda_{>0}$.
\end{itemize}
Define the cochain map $\tilde{\Psi} \fc C' \to \tilde{C}$ by $\tilde{\Psi} = (K)^{-1} \circ \Psi$. A straightforward computation shows that $\tilde{\Psi}$ satisfies $$\tilde{\Psi}(\lambda \check{z}_j) = T^{\max(j,0)\cdot\Delta}\lambda \check{x}_j$$ for every $-(n-1) \leq j \leq 1$ and $\lambda \in \Lambda_{> 0}$.

\subsubsection*{Computation using Schauder bases}

Let $-(n-1) \le j \le 1$. Another important computation for us is the map $\overline{\Psi}\fc \Lambda_{>0}e \to \coker D$, which is induced from the map $\Psi\fc \Lambda_{>0}e \to \widehat{\Lambda_{>0}^\infty}$ satisfying
$$\Psi(\lambda e) = T^{\max(j,0)\cdot \beta} \lambda e_1,$$
for every $\lambda\in \Lambda_{>0}$, where $e$ is the unit of $\Lambda_{\geq0}$ and $e_1$ is the first element in the standard basis of $\widehat{\Lambda_{\geq0}^\infty}$. This induced map is defined via the following commutative diagram:

\begin{equation}\label{diag: alg res maps modified}
  \xymatrix@R=2pc@C=3pc{
    0 \ar[r] \ar[d] & \Lambda_{>0}e \ar[d]_{\Psi} \ar[r]^{\id} & \Lambda_{>0}e \ar[d]_{\overline{\Psi}} \\
    \widehat{\Lambda_{>0}^\infty} \ar[r]^{D} & \widehat{\Lambda_{>0}^\infty} \ar[r]^{q} & \coker D
  }
\end{equation}
where $q\fc \widehat{\Lambda_{>0}^\infty}\to \coker D$ is the quotient map described in Proposition~\ref{prop: coker of D modified}.

\begin{prop}\label{prop: alg restriction maps modified}
    The map $\overline{\Psi} \fc \Lambda_{>0}e \to \coker D$ can be explicitly described as follows:
    \begin{enumerate}
        \item If $\alpha < \beta$, then the map
        $$\overline{\Psi} \fc \Lambda_{>0}e \to \Lambda_{(0,\gamma]}\epsilon_1 \oplus \bigoplus_{i\geq 2}\Lambda_{(0,\alpha]}\epsilon_i$$
        is $0$ in case where $j=1$, otherwise it satisfies, for every $\lambda \in \Lambda_{>0}$,
        $$\overline{\Psi}(\lambda e) =  [\lambda]_\gamma \epsilon_1 - T^{\beta-\gamma} [\lambda]_\alpha \sum_{i=2}^\infty (-T^{\beta-\alpha})^{i-2} \epsilon_i ,$$
        where $\epsilon_1$ is a generator of $\Lambda_{[0,\gamma]}$ and $\epsilon_i$ is a generator of the $(i-1)$-th summand in the direct sum $\bigoplus_{i\geq 2}\Lambda_{[0,\alpha]} \subset \coker D$.

        \item If $\gamma \geq \beta$, then the map
        $$\overline{\Psi} \fc \Lambda_{>0}e \to \Lambda_{>0}\epsilon_1 \oplus \bigoplus_{i\geq 2}\Lambda_{(0,\beta]}\epsilon_i$$
        satisfies, for every $\lambda \in \Lambda_{>0}$,
        $$\overline{\Psi}(\lambda e) = T^{\max(j,0)\cdot \beta} \lambda \epsilon_1.$$

        \item If $\alpha \geq \beta > \gamma$, then the map
        $$\overline{\Psi} \fc \Lambda_{>0}e \to \Lambda_{(0,\gamma]}\epsilon_1 \oplus \Lambda_{>0}\epsilon_2 \oplus \bigoplus_{i\geq 3}\Lambda_{(0,\beta]}\epsilon_i$$
        satisfies, for every $\lambda \in \Lambda_{>0}$,
        $$\overline{\Psi}(\lambda e) = T^{\max(j,0)\cdot \beta} \left( [\lambda]_\gamma \epsilon_1 - T^{\beta-\gamma} \lambda \epsilon_2 \right).$$
    \end{enumerate}
\end{prop}

\begin{proof}
    The commutativity of the diagram in Equation~\eqref{diag: alg res maps modified} dictates that $\overline{\Psi} \circ \id = q \circ \Psi$. Therefore, for every $\lambda \in \Lambda_{>0}$, we have:
    $$\overline{\Psi}(\lambda e) = q(\Psi(\lambda e)) = q(T^{\max(j,0)\cdot \beta} \lambda e_1) = T^{\max(j,0)\cdot \beta} q(\lambda e_1).$$

    \begin{enumerate}
        \item Assume that $\alpha < \beta$. In this case, we use the Schauder basis $(\epsilon_i)_{i=1}^\infty$ given by $\epsilon_1 = e_1 + T^{\beta-\gamma} e_2$ and $\epsilon_i = e_i + T^{\beta-\alpha} e_{i+1}$ for every $i \geq 2$. The element $e_1$ can be expressed in terms of the Schauder basis $(\epsilon_i)_{i=1}^\infty$ as:
        $$e_1 = \epsilon_1 - T^{\beta-\gamma} \sum_{i=2}^\infty (-T^{\beta-\alpha})^{i-2} \epsilon_i.$$
        By Proposition~\ref{prop: coker of D modified}, the quotient map $q \fc \widehat{\Lambda_{>0}^\infty} \to \coker D$ satisfies $q(\lambda \epsilon_1) = [\lambda]_\gamma \epsilon_1$ and $q(\lambda \epsilon_i) = [\lambda]_\alpha \epsilon_i$ for every $i \geq 2$. Therefore,
        \begin{align*}
            \overline{\Psi}(\lambda e) &= T^{\max(j,0)\cdot \beta} q(\lambda e_1) \\
            &= T^{\max(j,0)\cdot \beta} q\left( \lambda \epsilon_1 - \lambda T^{\beta-\gamma} \sum_{i=2}^\infty (-T^{\beta-\alpha})^{i-2} \epsilon_i \right) \\
            &= T^{\max(j,0)\cdot \beta} \left( [\lambda]_\gamma \epsilon_1 - T^{\beta-\gamma} [\lambda]_\alpha \sum_{i=2}^\infty (-T^{\beta-\alpha})^{i-2} \epsilon_i \right),
        \end{align*}
        for every $\lambda \in \Lambda_{>0}$.
        Note that if $j=1$ then
        $$\overline{\Psi}(\lambda e)= [T^\beta \lambda]_\gamma \epsilon_1 - T^{\beta-\gamma} [T^\beta \lambda]_\alpha \sum_{i=2}^\infty (-T^{\beta-\alpha})^{i-2} \epsilon_i =0.$$
        \item Assume that $\gamma \geq \beta$. Here, the Schauder basis is given by $\epsilon_1 = e_1$, with subsequent elements $\epsilon_i$ defined in terms of $e_i$ and $e_{i-1}$. Since $\epsilon_1 = e_1$, and Proposition~\ref{prop: coker of D modified} states that $q(\lambda \epsilon_1) = \lambda \epsilon_1$ for every $\lambda \in \Lambda_{>0}$, we immediately get:
        $$\overline{\Psi}(\lambda e) = T^{\max(j,0)\cdot \beta} q(\lambda e_1) = T^{\max(j,0)\cdot \beta} q(\lambda \epsilon_1) = T^{\max(j,0)\cdot \beta} \lambda \epsilon_1,$$
        for every $\lambda \in \Lambda_{>0}$.

        \item Assume that $\alpha \geq \beta > \gamma$. The Schauder basis begins with $\epsilon_1 = e_1 + T^{\beta-\gamma} e_2$ and $\epsilon_2 = e_2$. Thus, we can easily express $e_1$ as $e_1 = \epsilon_1 - T^{\beta-\gamma} \epsilon_2$. By Proposition~\ref{prop: coker of D modified}, $q(\lambda \epsilon_1) = [\lambda]_\gamma \epsilon_1$ and $q(\lambda \epsilon_2) = \lambda \epsilon_2$. Therefore,
        \begin{align*}
            \overline{\Psi}(\lambda e) &= T^{\max(j,0)\cdot \beta} q(\lambda e_1) \\
            &= T^{\max(j,0)\cdot \beta} q(\lambda \epsilon_1 - \lambda T^{\beta-\gamma} \epsilon_2) \\
            &= T^{\max(j,0)\cdot \beta} \left( [\lambda]_\gamma \epsilon_1 - T^{\beta-\gamma} \lambda \epsilon_2 \right),
        \end{align*}
        for every $\lambda \in \Lambda_{>0}$.
    \end{enumerate}
    This completes the proof.
\end{proof}

\subsubsection*{End of the proof}

As seen in the proof of Theorem~\ref{thm: SH of complement of ball in CPn}, and as it is mentioned in Section~\ref{ss: relSH using MB}, we can compute the relative symplectic cohomology using acceleration datum of Hamiltonians that satisfy the \textbf{MB} condition, therefore we have $SH_{\CP^n}(\CP^n\setminus\Int B;\Lambda_{\geq0}) \cong H(\widehat{\tilde{C}})$ and $SH_{\CP^n}(\CP^n;\Lambda_{\geq0}) \cong H(\widehat{C'}).$ Our next objective is to compute the induced map $(\widehat{\Psi})_* \fc SH_{\CP^n}(\CP^n;\Lambda_{\geq0}) \to SH_{\CP^n}(\CP^n\setminus\Int B;\Lambda_{\geq0})$.

\begin{prop}\label{prop: drctlim restriction maps CPn->CPn-B}
    Let $-(n-1) \leq j \leq 1$. Let $e$ be the generator of $SH_{\CP^n}^{-2j}(\CP^n;\Lambda_{\geq0}) \cong \Lambda_{>0} e$. The morphism 
    $$H(\widehat{\Psi})_{-2j} \fc SH_{\CP^n}^{-2j}(\CP^n;\Lambda_{\geq0}) \to SH_{\CP^n}^{-2j}(\CP^n\setminus\Int B^{}_{\Delta};\Lambda_{\geq0})$$
    depends on the value of $\Delta$ as follows:
    
    \begin{itemize}
        \item If $\Delta \leq \frac{1}{2}$, then the induced map satisfies
        $$H(\widehat{\Psi})_{-2j}(\lambda e) = \begin{cases}
            T^{\Delta} \lambda e_1, & j = 1, \\
            \lambda e_1, & -(n-1) \leq j \leq 0.
        \end{cases}$$

        \item If $\frac{i}{i+1} < \Delta \leq \frac{i+1}{i+2}$ for some integer $1 \leq i \leq n-1$, then the induced map satisfies
        $$H(\widehat{\Psi})_{-2j}(\lambda e) = \begin{cases}
            T^{\Delta} \lambda e_1, & j = 1, \\
            \lambda e_1, & i-n < j \leq 0, \\
            [\lambda]_{(j+n)(1-\Delta)} e_1 - T^{\Delta-(j+n)(1-\Delta)} \lambda e_2, & -(n-1) \leq j \leq i-n.
        \end{cases}$$

        \item If $\Delta > \frac{n}{n+1}$, then the induced map satisfies the following
        \begin{itemize}
            \item If $j=1$ then $H(\widehat{\Psi})_{-2}(\lambda e)=0$ for every $\lambda\in\Lambda_{>0}$.
            \item If $j\leq 0$ then 
            $$H(\widehat{\Psi})_{-2j}(\lambda e)=[\lambda]_{(j+n)(1-\Delta)} e_1 - T^{\Delta-(j+n)(1-\Delta)} [\lambda]_{n(1-\Delta)} \sum\limits_{k=2}^\infty (-T^{\Delta-n(1-\Delta)})^{k-2} e_k,$$
            for every $\lambda\in\Lambda_{>0}$.
        \end{itemize}
        
    \end{itemize}  
\end{prop}

\begin{proof}
Consider the subcomplexes
$$C_j = \bigoplus_{m \geq 0} \left( \Lambda_{>0} \check{x}_{j+m(n+1)} \oplus \Lambda_{>0} \hat{x}_{j+m(n+1)+n} \right), \quad C'_j = \Lambda_{>0} \check{z}_j,$$
of $(\tilde{C}, \tilde{d})$ and $(C', d')$, respectively. We identify these components with copies of $\widehat{\Lambda_{>0}^\infty}$ via the maps:
$$\lambda \hat{x}_{j+m(n+1)+n} \mapsto \lambda e_{m+1}, \quad \lambda \check{x}_{j+m(n+1)} \mapsto \lambda e_{m+1},$$
for every $m \geq 0$ and $\lambda \in \Lambda_{>0}$, and replace the notation $\check{z}_j$ by $e'_1$. Under these identifications, the map $\tilde{\Psi}_{-2j}$ satisfies $\tilde{\Psi}_{-2j}(\lambda e'_1) = T^{\max(j,0)\cdot\Delta} \lambda e_1$.

The restriction of the differential $\tilde{d}$ to $C_j$ is identified with the operator $D \colon \widehat{\Lambda_{>0}^\infty} \to \widehat{\Lambda_{>0}^\infty}$. By evaluating the differential on $e_1$, the parameters of this operator depend on the degree index $j$ and the capacity $\Delta$ via the relations:
\begin{itemize}
    \item If $j=1$, then $\alpha = n(1-\Delta)$, $\gamma = n(1-\Delta)$, and $\beta = \Delta$.
    \item If $-(n-1) \leq j \leq 0$, then $\alpha = n(1-\Delta)$, $\gamma = (j+n)(1-\Delta)$, and $\beta = \Delta$.
\end{itemize}

By Proposition~\ref{prop: alg restriction maps modified}, the induced map $H(\widehat{\Psi})_{-2j}$ depends on the value of $\Delta$ as follows:

\begin{itemize}
    \item If $\Delta \leq \frac{1}{2}$, then $\gamma \geq \beta$ for all $-(n-1)\leq j\leq 1$. Therefore we get that  $H(\widehat{\Psi})_{-2j}(\lambda e) =T^{\max(j,0)\cdot \Delta} \lambda e_1$ for every $\lambda\in\Lambda_{>0}$.
    
    \item If $\frac{i}{i+1} < \Delta \leq \frac{i+1}{i+2}$ for an integer $1 \leq i \leq n-1$, the behavior bifurcates based on $j$:
    \begin{itemize}

    \item For $j=1$, we have $\gamma = n(1-\Delta)$. Since $i \leq n-1$, it follows that $n \geq i+1$. Thus, $\Delta \leq \frac{i+1}{i+2} \leq \frac{n}{n+1}$, which rearranges to $n(1-\Delta) \geq \Delta$. This satisfies $\gamma \geq \beta$, yielding $H(\widehat{\Psi})_{-2j}(\lambda e)=T^\Delta \lambda e_1$ for every $\lambda\in\Lambda_{>0}$.
    
    \item For $i-n < j \leq 0$, we have $j+n \geq i+1$. Therefore, the assumed upper bound $\Delta \leq \frac{i+1}{i+2}$ implies $\Delta \leq \frac{j+n}{j+n+1}$. This implies that $(j+n)(1-\Delta) \geq \Delta$, meaning $\gamma \geq \beta$. This yields that $H(\widehat{\Psi})_{-2j}(\lambda e)= \lambda e_1$ for every $\lambda\in\Lambda_{>0}$.
    
    \item For $-(n-1) \leq j \leq i-n$, we have $j+n \leq i$. The assumed lower bound $\Delta > \frac{i}{i+1}$ implies $\Delta > \frac{j+n}{j+n+1}$, which rearranges to $(j+n)(1-\Delta) < \Delta$, meaning $\gamma < \beta$. However, since $i \leq n-1$, the assumed upper bound guarantees $\Delta \leq \frac{i+1}{i+2} \leq \frac{n}{n+1}$, ensuring $\alpha = n(1-\Delta) \geq \Delta = \beta$. Thus, $\alpha \geq \beta > \gamma$, yielding $H(\widehat{\Psi})_{-2j}(\lambda e)=[\lambda]_{(j+n)(1-\Delta)} e_1 - T^{\Delta-(j+n)(1-\Delta)} \lambda e_2$ for every $\lambda\in\Lambda_{>0}$.

    \end{itemize}

    \item If $\Delta > \frac{n}{n+1}$, then $\alpha < \beta$. In this case, the formulas are immediately obtained by placing the parameters $\alpha, \beta,$ and $\gamma$ into the formula of Proposition~\ref{prop: alg restriction maps modified}, distinguishing between the cases $j=1$ and $j \leq 0$.
\end{itemize}
This completes the proof.
\end{proof}

\begin{proof}[Proof of Proposition~\ref{prop: res from CPn to the complement a ball}]
By Theorem~\ref{prop: continuation maps for restrictions for CP^n->CP^n-B}, we have the commutativity $\Phi_\ell \circ \Psi_\ell = \Psi_{\ell+1} \circ \Phi'_\ell$, where the continuation maps $\Phi_\ell, \Phi'_\ell, \Psi_\ell$, and $\Psi_{\ell+1}$ appear in the following diagram:

\begin{equation*}
\xymatrix@R=2pc@C=3pc{
    CF(G'_\ell) \ar[r]^{\Phi'_\ell} \ar[d]^{\Psi_\ell} & CF(G'_{\ell+1}) \ar[d]^{\Psi_{\ell+1}} \\ 
    CF(G_\ell) \ar[r]^{\Phi_\ell} & CF(G_{\ell+1})
}
\end{equation*}
for every $\ell \in \Z_{\geq 0}$. Since $\drctlim CF(G'_\ell)=\Lambda_{>0}^{n+1}$ is a complete module, we deduce that the canonical map 
$$H(\drctlim CF(G'_\ell)) \to H(\widehat{\drctlim} CF(G'_\ell))$$
is an isomorphism (even in the cochain level), and in particular, surjective. Thus, by Proposition~\ref{prop: algebraic preparation for computing restriction maps}, the following diagram commutes:

\begin{equation*}
\xymatrix@R=2pc@C=3pc{
    H(\widehat{\tel}\, CF(G'_\ell)) \ar[r] \ar[d]^{\cong} & H(\widehat{\tel}\, CF(G_\ell)) \ar[d]^{\cong} \\ 
    H(\widehat{\drctlim} CF(G'_\ell)) \ar[r] & H(\widehat{\drctlim} CF(G_\ell))
}
\end{equation*}

Furthermore, the vertical arrows in this diagram are isomorphisms. The top horizontal arrow represents the restriction map $\res \fc SH_{\CP^n}(\CP^n;\Lambda_{\geq 0}) \to SH_{\CP^n}(\CP^n\setminus \Int B^{}_{\Delta};\Lambda_{\geq 0})$. Consequently, the commutativity of the diagram and the fact that the vertical arrows are isomorphisms allow us to identify the restriction map with the bottom horizontal arrow. The explicit formula for this bottom arrow is provided in Proposition~\ref{prop: drctlim restriction maps CPn->CPn-B}, which completes the proof.
\end{proof}

	\bibliographystyle{alpha}

\begin{thebibliography}{99}
		
		
		
        
        
        \bibitem{AGV_2024_framed_E2_structures} Abouzaid, Mohammed, Yoel Groman, and Umut Varolg\"une\c s, {\it Framed E2 structures in Floer theory,} Advances in Mathematics, {\bf 450}: 7109755, (2024).


		\bibitem{Audin_Damian_2014_Morse_and_Floer} Audin, Mich\'ele, and  Damian, Mihai, {\it Morse theory and Floer homology,} London: Springer, Vol. \textbf{2}. (2014).

        \bibitem{Bai_Xu_2026_integral_Floer} Bai, Shaoyun, and Xu, Guangbo, {\it Integral Hamiltonian Floer theory: Foundations,} arXiv preprint arXiv:2606.15847 (2026).

        \bibitem{Banyaga_Hurtubise_2013_Cascades} Banyaga, Augustin, and Hurtubise, David E., {\it Cascades and perturbed Morse-Bott functions,} Algebr. Geom. Topol, {\bf 13(1)}: 237--275, (2013).


        
		

        \bibitem{Biran_2001_Lag_barriers_and_symp_emb}	Biran, Paul, {\it  Lagrangian barriers and symplectic embeddings,}  GAFA, Geom. funct. anal. \textbf{11}, 407--464 (2001).

        \bibitem{BES_2025_Maurer_Cartan_elements_in_SH}	 Borman, Matthew Strom, El Alami, Mohamed, and Nick Sheridan, {\it Maurer-Cartan elements in symplectic cohomology from compactifications,}  Advances in Mathematics \textbf{482}, 110633 (2025).


        \bibitem{BSV_2022_QH_deformamation of SH}	 Borman, Matthew Strom, Sheridan, Nick, and Varolg\"une\c s, Umut, {\it Quantum cohomology as a deformation of symplectic cohomology,} Symplectic Geometry: A Festschrift in Honour of Claude Viterbo's 60th Birthday. Cham: Springer International Publishing \textbf{482}, 1073--1149 (2022).        

        \bibitem{BO_Morse_Bott_2009} Bourgeois, Fr{\'e}d{\'e}ric and Oancea, Alexandru, {\it Symplectic homology, autonomous Hamiltonians, and Morse--Bott moduli spaces,} 
        Duke Math. J., {\bf 146 (1)}:  71 --174, (2009).

		
    \bibitem{CFH_1995_SH_II_general} Cieliebak, Kai, Floer, Andreas, and Helmut Hofer, {\it  Symplectic homology II: A general construction} Mathematische Zeitschrift, {\bf 218(1)}: 103--122, (1995).
		
		
		
		\bibitem{CFHW_1996_ApSH_II} Cieliebak, Kai, Floer, Andreas, Hofer,  Helmut, and Wysocki, Konrad, {\it Applications of symplectic homology II: Stability of the action spectrum,} Mathematische Zeitschrift, {\bf 223}: 27--45, (1996).



		
		
		
		\bibitem{Dickstein_2025_PhD} Dickstein, Adi, {\it {Generalized measures in symplectic topology},} Tel-Aviv University (2025).
		\bibitem{DGPZ_2024_Symp_top_and_IVMs} Dickstein, Adi, Ganor, Yaniv, Polterovich, Leonid and Zapolsky, Frol, {\it Symplectic topology and ideal-valued measures,} Selecta Mathematica, {\bf 30(5)}(2024)
		
		\bibitem{Eliashberg_Polterovich_2010_qs} Eliashberg, Yakov, and Leonid Polterovich, {\it Symplectic quasi-states on the quadric surface and Lagrangian submanifolds,}  arXiv preprint arXiv:1006.2501, (2010).
        
		\bibitem{EP_2009_rigid_subsets} Entov, Michael, and Polterovich, Leonid, {\it  Rigid subsets of symplectic manifolds,} Compositio Mathematica, {\bf 145(3)}: 773--826, (2009).

    \bibitem{Fauck_2016_thesis} Fauck, Alexander, {\it  Rabinowitz-Floer homology on Brieskorn manifolds,} arXiv preprint arXiv:1605.07892 (2016).
        
        \bibitem{Floer_1989_monotone_mfd} Floer, Andreas, {\it  Symplectic fixed points and holomorphic spheres,} Communications in Mathematical Physics, {\bf 120(4)}: 575--611, (1989).

        \bibitem{FH_1994_SH_I_open_sets_Cn} Floer, Andreas, and Helmut Hofer {\it  Symplectic homology I open sets in $\C^n$,} Mathematische Zeitschrift, {\bf 215(1)}: 37--88, (1994).


        \bibitem{FHW_1994_SH_I} Floer, Andreas, Hofer, Helmut, and Wysocki, Konrad {\it   Applications of symplectic homology I,} Mathematische Zeitschrift, {\bf 217(1)}: 577-606, (1994).



        \bibitem{FHS_1995_Transversality} Floer, Andreas, Hofer, Helmut, and Salamon, Dietmar {\it   Transversality in elliptic Morse theory for the symplectic action,} Mathematische Zeitschrift, {\bf 80(1)}: 251--292, (1995).





        \bibitem{Frauenfelder_2004_MB} Frauenfelder, Urs, {\it  The Arnold-Givental conjecture and moment Floer homology,} International Mathematics Research Notices, {\bf 2004(42)}: 2179--2269, (2004).
		
		\bibitem{Ganatra_Pomerleano_2021_log_PSS} Ganatra, Sheel and Pomerleano, Daniel, {\it A log PSS morphism with applications to Lagrangian embeddings,} Journal of Topology, {\bf 14(1)}: 291--368, (2021).

        \bibitem{Griffiths_Harris_1978}   Griffiths, Phillip, and Harris, Joseph, {\it Principles of algebraic geometry,} John Wiley \& Sons (1978).
		

		\bibitem{Groman_Varolgunes_2023_locality_of_SH} Groman, Yoel, and Umut Varolg\"une\c s, {\it Locality of relative symplectic cohomology for complete embeddings,}  Compositio Mathematica, {\bf 159(12)}: 2551--2637, (2023).
		
		
        \bibitem{Gromov_1985_J_curves} Gromov, Mikhail, {\it Pseudo holomorphic curves in symplectic manifolds,} Inventiones mathematicae, {\bf 82(2)}: 307--347, (1985).
        		
	\bibitem{Gutt_indexes} Gutt, Jean, {\it The Conley-Zehnder index for a path of symplectic matrices,} arXiv preprint arXiv:1201.3728, (2012).
		
		
	\bibitem{Jo_Kang_2026_Biran_relSH}
    Jo, Hong-kwon, and Kang, Jungsoo,
    {\it Biran decomposition and relative symplectic cohomology,}
    arXiv preprint arXiv:2610.01334 (2026).	
	
		
        \bibitem{Lerman_1995_symplectic_cuts} Lerman, Eugene, {\it  Symplectic cuts,}  Mathematical Research Letters, {\bf 2(3)}: 247--258, (1995).		
        
	
        \bibitem{MSV_2024_heavy_sets_and_SH} Mak, Cheuk Yu, Sun, Yuhan, and Varolg\"une\c s, Umut, {\it  A characterization of heaviness in terms of relative symplectic cohomology,}  Journal of Topology, {\bf 17(1)}: e12327, (2024).	

        	
        \bibitem{McDuff_1993_symplectic _blow_up} McDuff, Dusa, {\it  Remarks on the uniqueness of symplectic blowing up,}  Symplectic geometry, {\bf 192}: 157--167, (1993).		
		
        \bibitem{McDuff_Salamon_J_curves_2012} McDuff, Dusa, and Dietmar Salamon, {\it $ J $-holomorphic curves and symplectic topology,} American Mathematical Society Vol 52, (2012).

        \bibitem{McDuff_Salamon_intro_2016} McDuff, Dusa, and Dietmar Salamon, {\it Introduction to symplectic topology,} Oxford University Press, (2016).
		
		
        
	   \bibitem{Oancea_survey} Oancea, Alexandru, {\it A survey of Floer homology for manifolds with contact type boundary or symplectic homology,} arXiv preprint arXiv: math/0403377, (2004).

    
		\bibitem{Pardon_2016_alg_vir_counting} Pardon, John, {\it An algebraic approach to virtual fundamental cycles on moduli spaces of pseudo-holomorphic curves,} Geometry \& Topology, {\bf 20(2)}: 779--1034, (2016).
		
       \bibitem{Pardon_2019_contact} Pardon, John, {\it Contact homology and virtual fundamental cycles,} Journal of the American Mathematical Society, {\bf 32(3)}: 825--919, (2019).
		
		
		
		
			\bibitem{PRSZ_2020_persistence} Polterovich, Leonid and Rosen, Daniel and Samvelyan, Karina and Zhang, Jun, {\it Topological persistence in geometry and analysis,} American Mathematical Soc., \textbf{74} (2020).

       
		\bibitem{Ritter_TQFT}  Ritter, Alexander F., {\it Topological quantum field theory structure on symplectic cohomology,}  Journal of Topology, {\bf 6(2)}: 391--489, (2013).


        \bibitem{RS_index} Robbin, Joel and Salamon, Dietmar, {\it The Maslov index for paths,} Topology, {\bf 32(4)}: 827--844, (1993).
		
		
		\bibitem{Rotman_intro_Homo_alg} Rotman, Joseph J., {\it An introduction to homological algebra,} New York: Springer, \textbf{2} (2009).


    \bibitem{Salamon_1999_notes} Salamon, Dietmar, {\it Lectures on Floer homology,} Symplectic geometry and topology (Park City, UT, 1997), {\bf 7}. 143--229, (1999).

        

        \bibitem{Schmaschke_2016_thesis} Schm\"aschke, Felix, {\it Floer homology of Lagrangians in clean intersection,}  arXiv preprint arXiv:1606.05327, (2016).
        
		\bibitem{Seidel_Fukaya_A_infty_Strcs_Assoc_to_Lef_Fib. III} Seidel, Paul, {\it Fukaya $A_\infty$-Structures Associated to Lefschetz Fibrations. III,} Journal of Differential Geometry, {\bf 117 (3)}: 485--589, (2021).




        \bibitem{Sun_2024_ind-bdd_relSH} Sun, Yuhan, {\it Index-bounded relative symplectic cohomology,} Algebraic \& Geometric Topology, {\bf 24(9)}: 4799--4836, (2024).


         \bibitem{Sun_2024_heavy + relSH} Sun, Yuhan, {\it  Heavy sets and index bounded relative symplectic cohomology,} Journal of Fixed Point Theory and Applications, {\bf 26(2)}: 21, (2024).
		
        
        
		
		\bibitem{Tonkonog_Varolgunes_Super_rigidity_of_certain_skeleta_using_relative_symplectic_cohomology} Tonkonog, Dmitry and Varolg\"une\c s, Umut, {\it Super-rigidity of certain skeleta using relative symplectic cohomology,} Journal of Topology and Analysis, 1--49, (2021).


        \bibitem{Usher_2010_energy_capacity_inequality}  Usher, Michael, {\it {The sharp energy-capacity inequality},} Communications in Contemporary Mathematics, {\bf 12(03)}: 457--473, (2010).
		

                        
		

                        
		\bibitem{Varolgunes_2018_PhD} Varolg\"une\c s, Umut: {\it {Mayer-Vietoris property for relative symplectic cohomology},} Massachusetts Institute of Technology (2018).
		
		\bibitem{Varolgunes_2021_MV_and_relSH} Varolg\"une\c s, Umut, {\it {Mayer-Vietoris property for relative symplectic cohomology},} Geometry \& Topology, {\bf 25(2)}: 547--642, (2021).
		
		
		
		

        
        \bibitem{Wehrheim_Woodward_2-10_Quilted} Wehrheim, Katrin, and Chris T. Woodward, {\it {Quilted floer cohomology},}  Geometry \& Topology , {\bf 14(2)}: 833--902, (2010).
        
\bibitem{Weibel_1994_hom_alg} Weibel, Charles Alexander: {\it {An Introduction to Homological Algebra},} Rutgers University, Department of Mathematics (1994).
		


        \bibitem{Wendl_16_Lec_SFT} Wendl, Chris, {\it {Lectures on symplectic field theory,}}  arXiv preprint arXiv:1612.01009, (2016).
        
		
		
	\end{thebibliography}

\noindent
\begin{tabular}{l}
{\bf Adi Dickstein} \\
D\'epartement de Math\'ematiques et de Statistique\\
Universit\'e de Montr\'eal\\
CP 6128 succ Centre-Ville\\
Montr\'eal QC H3C 3J7, Canada\\
{\em E-mail:}  \texttt{adi.dickstein@gmail.com}
\end{tabular}

\medskip

    \noindent
    \begin{tabular}{l}
    {\bf Yaniv Ganor} \\
    School of Mathematical Sciences,\\
    Faculty of Sciences,\\
    HIT - Holon Institute of Technology \\
    
    Golomb St. 52, Holon 5810201, Israel\\
    {\em E-mail:}  \texttt{ganory@gmail.com}
    \end{tabular}

\end{document}